\documentclass[11pt,
oneside
]{book}

\usepackage{amsmath,amssymb,amsthm}

\usepackage{newtxtext, newtxmath}
\usepackage[margin=3cm]{geometry}

\usepackage{stmaryrd}
\usepackage{tikz}
\usetikzlibrary{shapes,backgrounds,matrix,arrows,automata,positioning,decorations}
\usepackage{tikz-cd}
\usepackage{enumitem}
\usepackage[hyphens]{url}

\usepackage[
giveninits=true,  
uniquename=init,  
backend=biber,    
citestyle=authoryear, 
maxbibnames=99,    
url=false,        
doi=false,
isbn=false,
]{biblatex}
\bibliography{book.bib} 

\theoremstyle{plain}
\newtheorem{theorem}{Theorem}[chapter]
\newtheorem{proposition}[theorem]{Proposition}
\newtheorem{lemma}[theorem]{Lemma}
\newtheorem{corollary}[theorem]{Corollary}

\theoremstyle{definition}
\newtheorem{definition}[theorem]{Definition}
\newtheorem{example}[theorem]{Example}
\newtheorem{remark}[theorem]{Remark}

\newtheorem{notationnum}[theorem]{Notation}
\newtheorem*{notation}{Notation}

\newcommand{\hint}[1]{\textit{Hint.} #1}
\newcommand{\ourexercises}{\subsection*{Exercises for Section~\thesection}}
\newtheorem{ourexercise}{Exercise}[section]

\DeclareMathOperator{\mysetminus}{--}
\renewcommand{\setminus}{\mysetminus}

\renewcommand{\theta}{\vartheta}
\renewcommand{\phi}{\varphi}

\newcommand{\bfK}{\mathbf{K}}

\newcommand{\N}{\mathbb{N}}
\newcommand{\bN}{\mathbb{N}}
\newcommand{\bR}{\mathbb{R}}
\newcommand{\bZ}{\mathbb{Z}}

\newcommand{\cA}{\mathcal{A}}
\newcommand{\cB}{\mathcal{B}}
\newcommand{\cC}{\mathcal{C}}
\newcommand{\cD}{\mathcal{D}}

\newcommand{\cF}{\mathcal{F}}
\newcommand{\cG}{\mathcal{G}}
\newcommand{\cH}{\mathcal{H}}
\newcommand{\cI}{\mathcal{I}}
\newcommand{\cJ}{\mathcal{J}}
\newcommand{\cK}{\mathcal{K}}
\newcommand{\cL}{\mathcal{L}}
\newcommand{\cM}{\mathcal{M}}
\newcommand{\cN}{\mathcal{N}}
\newcommand{\cP}{\mathcal{P}}
\newcommand{\cQ}{\mathcal{Q}}
\newcommand{\cR}{\mathcal{R}}
\newcommand{\cS}{\mathcal{S}}
\newcommand{\cT}{\mathcal{T}}
\newcommand{\cU}{\mathcal{U}}
\newcommand{\cV}{\mathcal{V}}

\newcommand{\btwo}{{\bf 2}}
\newcommand{\Kel}{\mathrm{K}} 
\newcommand{\ClpCon}{\mathrm{ClpCon}} 

\newcommand{\wa}{\widehat{a}}
\newcommand{\wb}{\widehat{b}}
\newcommand{\simb}{\widetilde{b}}%
\newcommand{\simg}{\widetilde{g}}

\newcommand{\subword}{\sqsubseteq}

\newcommand{\End}{\mathrm{End}}
\newcommand{\three}{\mathbf{3}}
\newcommand{\cat}[1]{\ensuremath{\mathbf{#1}}}
\newcommand{\BA}{\cat{BA}}
\newcommand{\Pos}{\cat{Pos}}
\newcommand{\DL}{\cat{DL}}
\newcommand{\CDFrame}{\cat{CDFrame}}
\newcommand{\Set}{\cat{Set}}

\newcommand{\POS}{\cat{Pos}}

\newcommand{\dcpo}{\cat{dcpo}}
\newcommand{\TopCat}{\cat{Top}}
\newcommand{\SOB}{\cat{Sober}}
\newcommand{\CONT}{\cat{Domain}}
\newcommand{\Alg}{\cat{Alg}}
\newcommand{\Spec}{\cat{Spec}}
\newcommand{\BoolSp}{\cat{BoolSp}}
\newcommand{\Frame}{\cat{Frame}}
\newcommand{\SpFrame}{\cat{SpFrame}}

\newcommand{\Priest}{\cat{Priestley}}

\newcommand{\jpp}{{\normalfont\textsc{jpp}}}

\newcommand{\Pt}{\mathrm{Pt}} 
\newcommand{\PT}{\mathrm{PT}} 
\newcommand{\PTplus}{\mathrm{PT}^+}
\newcommand{\St}{\mathrm{St}}

\newcommand{\JT}{\cJ\mathrm{triv}}
\newcommand{\TopLat}{\mathit{Top}}

\newcommand{\UB}{\normalfont\textsc{UB}} %
\newcommand{\MUB}{\normalfont\textsc{MUB}} %
\newcommand{\Up}{\mathcal{U}} %
\newcommand{\Upfin}{\mathcal{U}_{\mathrm{fin}}} %
\newcommand{\cPfin}{\mathcal{P}_{\mathrm{fin}}} %
\newcommand{\Down}{\mathcal{D}} %
\newcommand{\Downfin}{\mathcal{D}_{\mathrm{fin}}} %
\newcommand{\ClD}{\mathrm{ClpD}} %
\newcommand{\ClU}{\mathrm{ClpU}} %
\newcommand{\Clp}{\mathrm{Clp}} %
\newcommand{\KO}{\mathcal{KO}}%
\newcommand{\KS}{\mathcal{KS}}%

\newcommand{\op}{\mathrm{op}} %
\newcommand{\Filt}{\mathrm{Filt}} %
\newcommand{\Idl}{\mathrm{Idl}} %
\newcommand{\PrFilt}{\mathrm{PrFilt}} %
\newcommand{\PrIdl}{\mathrm{PrIdl}} %
\newcommand{\CompPrFilt}{\mathrm{CompPrFilt}} %

\DeclareMathOperator{\ob}{\mathrm{ob}}
\DeclareMathOperator{\mor}{\mathrm{mor}}
\DeclareMathOperator{\cod}{\mathrm{cod}}
\DeclareMathOperator{\dom}{\mathrm{dom}}
\newcommand{\im}{\mathrm{im}}

\newcommand{\epp}{\normalfont\textsc{EPP}}
\newcommand{\RegA}{\mathrm{Reg}_A} %
\newcommand{\Reg}{\mathrm{Reg}}

\newcommand{\Diam}{\Diamond}
\newcommand{\Hom}{\mathrm{Hom}}
\newcommand{\id}{\mathrm{id}}
\newcommand{\Sub}{\mathrm{Sub}}
\newcommand{\intr}{\mathrm{int}} 
\newcommand{\cl}{\mathrm{cl}}   

\newcommand{\onto}{\twoheadrightarrow}
\newcommand{\twoheaduparrow}{{\rotatebox[origin=c]{90}{$\twoheadrightarrow$}}}
\newcommand\twoheaddownarrow{\mathrel{\rotatebox[origin=c]{-90}{$\twoheadrightarrow$}}}
\newcommand{\To}{\Rightarrow}
\newcommand{\gen}[1]{\langle #1 \rangle}
\newcommand{\into}{\hookrightarrow}
\newcommand{\hto}{\to} %
\newcommand{\wayb}{\mathrel{<\!\!\!<}}
\newcommand{\wwb}{\mathrel{<\!\!<\!\!<}}
\newcommand{\threeheaddownarrow}{\mathrel{\,\raisebox{10pt}{\rotatebox{-90}{$\twoheadrightarrow\hskip-9pt\twoheadrightarrow$}}\!}}
\newcommand{\waydown}{{\twoheaddownarrow}} %
\newcommand{\wayup}{{\twoheaduparrow}}
\newcommand{\covers}{\mathrel{-\mkern-4mu}<}

\newcommand{\liff}{\leftrightarrow}
\newcommand{\sem}[1]{\llbracket {#1} \rrbracket} %
\newcommand{\semI}[1]{\llbracket {#1} \rrbracket_I} %
\newcommand{\semE}[1]{\llbracket {#1} \rrbracket_E} %
\newcommand{\wh}[1]{\widehat{#1}}
\newcommand{\iffdef}{\stackrel{\mathrm{def}}{\Longleftrightarrow}}

\newcommand{\pos}[1]{\draw[fill=black] #1 circle (2pt);}
\newcommand{\pob}[1]{\draw[fill=black] #1 circle (2pt); \draw #1 circle (4pt);}
\newcommand{\po}[1]{\draw[fill=black] #1 circle (2pt);}
\newcommand{\polab}[3]{\draw[fill=black] #1 circle (2pt); \node[#3=.1cm] at #1 {#2};}
\newcommand{\li}[2]{\draw [thick] #1 -- #2;}
\newcommand{\leftright}[4]{\begin{tikzpicture}[baseline]%
    \node[anchor=base] (A) {$#1$};%
    \node[anchor=base] (B) [right of=A, node distance=4em] {$#2$};%

    \draw[->] ([yshift=-2pt]A.east) to node[midway,below=-2pt] {\scriptsize $#3$} ([yshift=-2pt]B.west);%
    \draw[->] ([yshift=2pt]B.west) -- node[midway,above=-2pt] {\scriptsize $#4$} ([yshift=2pt]A.east);
  \end{tikzpicture}}

\newcommand{\initarr}[1]{
  \node (INIT#1) [left=0.5cm of #1] {};
  \path (INIT#1) edge (#1);
}

\usepackage{hyphenat}
\usepackage{imakeidx}

\newcommand{\emphind}[1]{\emph{#1}\index{#1}}
\makeindex 

\makeindex[name=notation,title={Notation},columns=1] 

\usepackage{fmtcount}

\newcounter{nlcounter}
\newcommand{\nl}[3]{%
\addtocounter{nlcounter}{1}%
\index[notation]{\alph{chapter}.\thenlcounter@{#1} \ : \ {#2}}
}

\newcommand{\notessec}{\section*{Notes for Chapter~\thechapter}}

\usepackage[placement=bottom,scale=1,opacity=1,vshift=5mm]{background}
\newsavebox{\prelimtxt}
\sbox{\prelimtxt}{%
    \parbox{\textwidth}{\footnotesize\textit{\rm This material will be published by Cambridge University Press as ``Topological Duality for Distributive Lattices: Theory and Applications'' by Mai Gehrke and Sam van Gool. This pre-publication is free to view and download for personal use only. Not for re-distribution, re-sale, or use in derivative works. \copyright Mai Gehrke and Sam van Gool, \today}}
    }
\backgroundsetup{contents=\usebox{\prelimtxt}}

\usepackage[pdfencoding=auto,
]{hyperref}
\hypersetup{
linkbordercolor={white},
colorlinks=true,
linkcolor=.,
anchorcolor=.,
citecolor=.,   
filecolor=.,
menucolor=.,
runcolor=.,
urlcolor=.,
}
\usepackage[hyphenbreaks]{breakurl} 

\begin{document}
\title{Topological Duality for Distributive Lattices: \\ Theory and Applications}
\author{Mai Gehrke and Sam van Gool}
\date{\today}
\hypersetup{pageanchor=false}
\begin{titlepage}
\maketitle
\end{titlepage}

\hypersetup{pageanchor=true}

\tableofcontents
\chapter*{Introduction}
\addcontentsline{toc}{chapter}{Introduction}

This book is a course on Stone-Priestley duality theory, with applications to
logic and the foundations of computer science. Our target audience includes
both graduate students and researchers in mathematics and computer science. The
main aim of the book is to equip the reader with the theoretical background
necessary for reading and understanding current research in duality and its
applications. We aim to be didactic rather than exhaustive, while we do give
technical details whenever they are necessary for understanding what the field
is about.

Distributive lattice structures are fundamental to logic, and thus appear
throughout mathematics and computer science. The reason for this is that the
notion of a distributive lattice is extremely basic: it captures a language
containing as its only primitives the logical operators `or', `and', `true' and
`false'. Distributive lattices are to the study of logic what rings and vector
spaces are to the study of classical algebra. Moreover, distributive lattices
also appear in, for example, ring theory and functional analysis.

A mathematical kernel that makes duality theory tick is the fact that the
structure of a lattice can be viewed in three equivalent ways. A distributive
lattice is all of the following:
\begin{enumerate}
\item a partially ordered set satisfying certain properties regarding upper and
    lower bounds of finite sets;
\item an algebraic structure with two idempotent monoid operations that
    interact well with each other;
\item a basis of open sets for a particular kind of topological or
    order-topological space.
\end{enumerate}
The first part of the book will define precisely the vague notions in this list
(`certain properties', `interact well', `a particular kind of'), and will prove
that these are indeed three equivalent ways of looking at distributive
lattices. The correspondence between algebraic and topological structure in the
last two items of the list can be cast in a precise categorical form, and is
then called a \emph{dual equivalence} or simply \emph{duality}. A duality
identifies an exciting, almost magical, and often highly useful intersection
point of algebra and topology.

Historically, Stone showed in the 1930s that distributive lattices are in a
duality with spectral spaces: a certain class of topological spaces with a
non-trivial specialization order, which are also the Zariski spectra of rings.
Stone's duality for distributive lattices is especially well-known in the  more
restrictive setting of Boolean algebras, obtained by adding an operator `not'
to the lattice signature, which satisfies the usual rules of logic: de Morgan's
laws and excluded middle. The restriction of Stone's duality to Boolean
algebras shows that they are in a duality with compact Hausdorff
zero-dimensional spaces. While the spaces associated to Boolean algebras are
better known than the slightly more general ones associated to distributive
lattices, the latter are vastly more versatile, having as continuous retracts,
among others, \emph{all} compact Hausdorff spaces, including connected spaces
such as the unit interval of the reals.

Nevertheless, Stone's duality for distributive lattices was for at least thirty
years seen by many as a lesser sibling of his duality for Boolean algebras, at
least partly due to the fact that the spaces that figure in it are not
Hausdorff, and the appropriate functions between the spaces are not all the
continuous ones. Priestley's seminal work in the 1970s lifted this obstacle, by
giving a first-class role to the specialization order that figures in Stone's
spectral spaces. Priestley reframed Stone's duality as one between distributive
lattices and certain \emph{partially ordered} topological spaces, now called
Priestley spaces. The first goal in this book is to build up the necessary
mathematics to prove Priestley's duality theorem, which we do in
Chapter~\ref{ch:priestley}; we also show there how it easily specializes to the
case of Boolean algebras. Building up to this first main result,
Chapters~\ref{ch:order}~and~\ref{chap:TopOrd} will teach the foundations of
order theory and topology that we rely on in the rest of the book. 

A unique feature of this book is that, in addition to developing general
duality theory for distributive lattices, we also show how it applies in a
number of areas within the foundations of computer science, namely, modal and
intuitionistic logics, domain theory and automata theory. The use of duality
theory in these areas brings to the forefront how much their underlying
mathematical theories have in common. It also prompts us to upgrade our
treatment of duality theory with various enhancements that are now commonly
used in the state-of-the-art research in the field. Most of these enhancements
make use of \emph{operators} on a distributive lattice: maps between lattices
that only preserve part of the lattice structure.

The simplest kind of operator is a map between lattices that respects the
structure of `and' and `true', but not necessarily `or' and `false'. If this
notion is understood as analogous to a linear mapping in linear algebra, then
it is natural to also consider more general binary, ternary, and $n$-ary
operators on lattices, which respect the structure of `and' and `true' in each
coordinate, as long as the remaining coordinates are fixed. The theory of
lattices with operators, and dualities for them, was developed in the second
half of the 20th century, roughly in two main chunks. First, in the 1950s, by
J\'onsson and Tarski, in the case of Boolean algebras, with immediate
applications to relation algebra, and the same theory was used heavily a little
later and very successfully for modal logic in the form of Kripke's semantics.
However, until the end of the 1980s, the duality theory for distributive
lattices with additional operations developed in the form of a great number of
isolated case-by-case studies. Starting with the work of Goldblatt, and also of
J\'onsson and this book’s first author, the general theory of distributive
lattices with additional operations came into a mature, more usable, form by
the 1990s. This theory is developed in Chapter~\ref{ch:methods}, which also
contains the first applications of duality theory, to free distributive
lattices, quotients and subspaces, implication-type operators, Heyting algebras
and Boolean envelopes.

In the development of the first four chapters of this book, we keep the use of
category theory to a minimum. In Chapter~\ref{ch:categories}, we then set the
results of the earlier chapters in the more abstract and general framework of
category theory. This development then also allows us to show how Priestley's
duality fits well in a more general framework for the interaction of topology
and order, which had been developed by Nachbin shortly before.  In
Chapter~\ref{chap:Omega-Pt}, we show how the various classes of topological
spaces with and without order, introduced by Stone, Priestley and others, all
relate to each other, and how they are in duality with distributive lattices
and their infinitary variant, frames.

Chapter~\ref{chap:DomThry} and \ref{ch:AutThry} contain two more modern applications of duality theory to theoretical computer science, namely to domain theory and to automata theory, respectively.
The domain theory that we develop in Chapter~\ref{chap:DomThry} is organized around three separate results: Hoffmann-Lawson duality; the characterization of those dcpos and domains, respectively, that fall under Stone duality; and Abramsky's celebrated 1991 Domain Theory in Logical Form paper. 

The duality-theoretic approach to automata theory that we develop in
Chapter~\ref{ch:AutThry} originates in  work linking profinite methods in
automata theory with duality theory \parencite{GGP2008}. It is organized around
a number of related results, namely: finite syntactic monoids can be seen as
dual spaces, and the ensuing effectivity of this powerful invariant for regular
languages; the free profinite monoid is the dual of the Boolean algebra of
regular languages expanded with residuation operations and, more generally,
topological algebras on Boolean spaces are duals of certain Boolean algebras
extended by residual operations. As an extended application example, we use
duality to give a profinite equational characterization for the class of
piecewise testable languages; and we end by discussing a characterization of
those profinite monoids for which the multiplication is open.

\subsection*{How to use this book}
This is a textbook on distributive lattices, spectral spaces and Stone and
Priestley dualities as they have developed and are applied in various areas at
the intersection of algebra, logic, and theoretical computer science. Our aim
is to get in a fairly full palette of duality tools as directly and quickly as
possible, then to illustrate and further elaborate these tools within the
setting of three emblematic applications: semantics of propositional logics,
domain theory in logical form, and the theory of profinite monoids for the
study of regular languages and automata. The text is based on lecture notes
from a 50-hour course in the \emph{Master Logique et Fondements de
l'Informatique} at Paris 7, which ran in the winters of 2013 and 2014. The fact
that it is based on notes from a course means that it reaches its goals while
staying as brief and to the point as possible. The other consequence of its
origin is that, while it is fully a mathematics course, the applications we aim
at are in theoretical computer science. The text has been expanded a bit beyond
what was actually said in the course, reaching research monograph level by the
very end of the last two chapters~\ref{chap:DomThry} and \ref{ch:AutThry}.
Nevertheless, we have focused on keeping the spirit of a lean and lively
textbook throughout, including only what we need for the applications, and
often deferring more advanced general theory to the application chapter where
it becomes useful and relevant.

While the original course on which the book is based covered the majority of
all the chapters of the book, there are several other options for its use. In
particular, a basic undergraduate course on lattices and duality could treat
just chapters~\ref{ch:order} through \ref{ch:priestley} and possibly selected
parts of \ref{ch:methods}, \ref{ch:categories}, and/or \ref{chap:Omega-Pt}. The
applications in the second part are fairly independent and can be included as
wanted, although the domain theory material in Chapter~\ref{chap:DomThry}
requires at least skeletal versions of Chapter~\ref{ch:categories}, and
Chapter~\ref{chap:Omega-Pt} in its entirety.

The first part, Chapters~\ref{ch:order} through \ref{chap:Omega-Pt}, is a
graduate level `crash course' in duality theory as it is practiced now.
Chapter~\ref{ch:order} introduces orders and lattices, and in particular the
distributive lattices that we will be concentrating on, as well as the finite
case of Stone duality, where topology is not yet needed.
Chapter~\ref{chap:TopOrd} introduces the topological side of the dualities. In
this chapter, we elaborate the interaction between order and topology, which is
so central to the study of spaces coming from algebraic structures. For this
purpose we have bent our philosophy of minimum content somewhat by introducing
the class of stably compact spaces and Nachbin's equivalent class of compact
ordered spaces. We believe that this setting provides the right level of
generality for understanding the connection between Stone's original duality
for spectral spaces and Priestley duality. The class of stably compact spaces,
being the closure of spectral spaces under continuous retracts, is also a more
robust setting than spectral spaces for a number of further applications that
we do not cover in this book, such as continuous domain theory and duality for
sheaf representations of algebras. The basic mathematical content of Priestley
duality is given in Chapter~\ref{ch:priestley}. Chapter~\ref{ch:methods}
introduces the most important general methods of modern duality theory: duality
for additional operations and sub-quotient duality, which then allows us to
immediately give first applications to propositional logics.
Chapter~\ref{ch:categories} then introduces categorical concepts such as
adjunctions, dualities, filtered colimits, and cofiltered limits, which play a
fundamental role in duality theory. This allows us to give a full categorical
account of Priestley duality by the end of the chapter.
Chapter~\ref{chap:Omega-Pt} treats the Omega-Point duality and Stone's original
duality for distributive lattices and makes the relationship between these
dualities and Priestley's version clear.

The duality theory developed in the first six chapters of the book is applied
to two different parts of theoretical computer science in the last two
chapters, which provide an entry into research-level material on these topics.
These two chapters are independent from each other, and have indeed
traditionally been somewhat separate in the literature, but our treatment here
shows how both topics in fact can be understood using the same
duality-theoretic techniques that we develop in the first part of the book.
When using this book for a course, a lecturer can freely choose material from
either or both of these chapters, according to interest.
Chapter~\ref{chap:DomThry}, on domain theory, contains a duality-theoretic
exposition of the solutions to domain equations, a classical result in the
semantics of programming languages. Chapter~\ref{ch:AutThry} develops a duality
theory for algebraic automata theory, and shows in particular how finite and
profinite monoids can be viewed as instances of the dual spaces of lattices
with operators that we study in this book.

We have given some bibliographic references throughout the text. We want to emphasize
here that these references are not in any way meant to give an exhaustive bibliography for
the vast amount of existing research in duality theory. They 
are rather intended as useful entry points into the research literature appropriate for
someone learning this material, who will then find many further references there. 
Furthermore, at the end of several chapters, we have added 
a number of small, additional notes giving technical pointers pertaining to specific topics discussed there --
again, we do not mean to imply exhaustivity.
When we introduce special or less standard notation, we use `Notation' blocks, which 
are occasionally numbered when we need  to refer back to them later.
The book ends with a listing of the most-used notations and an index of concepts.

In each chapter, all numbered items follow one and the same counter, with the exception
of exercises, which follow a separate numbering, indicating not only chapter but also section.
On the topic of exercises: this book contains many of them, varying greatly in difficulty.
In earlier chapters, many of the
exercises are routine verifications, but necessary to do for a learner who wants to master the material.
Especially in later chapters, there are exercises that could be viewed as small research
projects, although we refrain from stating open problems as exercises: for the less obvious
exercises, we have included hints, and references where available.

\subsection*{Comparison to existing literature and innovative aspects}
The first part of this book, Chapters~\ref{ch:order} through
\ref{chap:Omega-Pt}, covers quite classical material and may be compared to
existing textbooks. The closest are probably \emph{Distributive Lattices}
\parencite{BalDwi1974}, and \emph{Introduction to Lattices and Order}
\parencite{DavPri2002}.  Another classical gentle introduction to the field,
but focusing more on point-free topologies and frame theory than we do here, is
\emph{Topology via Logic} \parencite{Vic1989}.   Of these, \cite{BalDwi1974} is
probably the closest in spirit to our treatment, as it gets to the duality
quickly and then applies it. However, that book's applications to algebras of
propositional logic focus on varieties that are less central today. Davey and
Priestley's textbook has been very successful and has in particular managed to
attract a theoretical computer science readership to these topics. However, it
focuses more on the lattices and order \emph{per se} and the duality is covered
only as one of the final crowning chapters. Davey and Priestley's book is
therefore an excellent way in to ours and we recommend it as supplemental
reading in case students are needing additional details or to build up
mathematical maturity. The recent textbooks on Boolean algebras
\parencite{HaGi2008} and on Duality theories for Boolean algebras
\parencite{givant2014duality} are also relevant but are of course restricted to
the Boolean setting. There are also a number of classical references in lattice
theory by Grätzer, the most recent versions being \emph{Lattice Theory:
Foundation} \parencite{Gratzer11} and \emph{General Lattice Theory}
\parencite{Gratzer03}, which each contain material on duality theory and its
applications to lattice theory. 

Here we aim to get the dualities in place as soon as possible and then use
them. Where we differ the most from the existing books within this first part
is with our emphasis on the interaction between order and topology in
Chapter~\ref{chap:TopOrd}, and in placing Priestley duality within the wider
context of category theory (Chapter~\ref{ch:categories}) and Omega-point
duality (Chapter~\ref{chap:Omega-Pt}). Chapter~\ref{chap:TopOrd} provides a
textbook-level didactic account of the interaction between topology and order
culminating with the equivalence between Nachbin's compact ordered spaces and
stably compact spaces. In Chapter~\ref{ch:methods} we develop duality theory
methods for analyzing the structure of distributive lattices and operators on
them. All of these topics have become central in research in recent decades but
are so far difficult to access without delving in to the specialized
literature.

Perhaps the most important omission of this book is the theory of canonical
extensions, which was central to the already-mentioned work of J\'onsson and
others, in addition to duality. While this theory is very close to both of this
book's authors' hearts, and closely related to duality theory, this book is not
about that, and canonical extensions thus do not play a big role in this book,
at least not explicitly. Still, we will occasionally make references to
canonical extensions where appropriate. Along with and closely related to this
omission, we decided to take a point-set rather than a point-free approach to
the topics of this book. Point-free approaches focus on the algebraic side of
duality and thus avoid the point-set world of topology, which inherently
involves non-constructive principles. Duality is in a sense the justification
for the point-free approach since it makes the link between the algebraic and
the point-set worlds. In this book we remain fully anchored in the
set-theoretic approach to topology, in particular making use of the Axiom of
Choice as necessary. We do this as it is more easily accessible for a general
audience, and because our end applications in denotational semantics and
profinite algebras in automata theory are, in their currently practiced form,
focused on point-set topology. That being said, our focus on duality shows the
way and familiarizes the reader with the dual, point-free approach, thus making
them ready to embrace this approach. In this direction, one of our hopes with
this book is that it will entice some readers to learn about canonical
extensions and related point-free techniques. We believe the technique of
canonical extensions to be complementary to, and at least as important as,
duality, but so far less well-established in the literature. 

Many research monographs include similar material to that of the first part of
this book, but are not explicitly targeted at readers who are first learning
about the field, while this is a primary aim of our book. Classical such
monographs, closest in content to the first part of this book, are \emph{Stone
spaces} \parencite{Johnstone1986} and the book \emph{A compendium of continuous
lattices} \parencite{Getc80},  re-edited as \emph{Continuous lattices and
domains} \parencite{Getc2003}. More recently, the monograph \emph{Spectral
spaces} \parencite{DicSchTre2019} studies the same class of spaces as we do in
this book, but coming from a ring-theoretic perspective and emphasizing less
the order-theoretic aspects. The monograph \emph{Non-Hausdorff Topology and
Domain Theory} \parencite{Goubault2013} is close in spirit to our treatment in
Chapter~\ref{chap:TopOrd}, especially in its treatment of stably compact
spaces,  and also addresses a theoretical computer science audience. A
difference with our treatment here is that \parencite{Goubault2013} is focused
on non-Hausdorff topologies and therefore does not treat the (Hausdorff) patch
topology as central, as we do here. Related to our Chapter~\ref{chap:Omega-Pt}
is the monograph \emph{Frames and Locales} \parencite{PicPul2012} focused on
frames and point-free topology, and Chapter~\ref{chap:Omega-Pt} of this book
can be used as a preparation for jumping into that work.

The applications to domain theory and automata theory are treated  in
Chapters~\ref{chap:DomThry} and \ref{ch:AutThry}, respectively. These two
applications, and in particular the fact that we treat them in one place, as
applications of a common theory, are perhaps the most innovative and special
aspects of this book. Domain theory is the most celebrated application of
duality in theoretical computer science and our treatment is entirely new.
Automata theory is a relatively new application area for duality theory and has
never been presented in textbook format before. More importantly, both topics
are at the forefront of active research seeking to unify semantic methods with
more algorithmic topics in finite model theory. While previous treatments
remained focused on the point of view of domains/profinite algebra, with
duality theory staying peripheral, a shared innovative aspect of the
presentations of these topics in this book is that both are presented squarely
as applications of duality.

Finally, a completely original contribution of this book, which emerged during
its writing, precisely thanks to our treatment of the two topics as an
application of a common theory, is the fact that a notion we call
``\emph{preserving joins at primes}'' turns out to be central in both the
chapter on domain theory and in that on automata theory. This notion was
introduced in the context of automata theory and topological algebra in
\parencite{Geh16}; its application to domain theory is new to this book and
reflects a key insight of Abramsky's Domain Theory in Logical Form. We believe
this point to be an exciting new direction for future research in the field
that we hope some readers of this book will be inspired to take up.

\newpage
\section*{Acknowledgements}
We thank Cambridge University Press, in particular David Tranah and Anna
Scriven, for their trust and help with publishing this book.

We would also like to thank all our colleagues and students who have
inspired and encouraged us to complete this book. More specifically, many have
used and commented on early versions and their comments have been immensely
useful.

Samson Abramsky, Jorge Almeida, Achim Jung, Jean Goubault-Larrecq, Jean-Éric
Pin, and Benjamin Steinberg have all provided us with invaluable advice and
support on several of the later parts of the book; their reading of parts of
our manuscript at various stages and their very detailed and useful feedback
have been of great help. We would also like to thank Jim de Groot and Luca
Reggio, who both proofread the manuscript very carefully, and also Célia
Borlido and Victor Iwaniack who also made many useful comments. Jérémie Marques
has made vital mathematical contributions to the book, also acknowledged
specifically in Chapter~\ref{ch:AutThry}.

We are very grateful for the friendship, comments, and support that we have
received from those already mentioned above, as well as from many other
colleagues and co-authors with whom we have discussed and worked on several of
the topics discussed in this book, specifically Clemens Berger, Guram
Bezhanishvili, Nick Bezhanishvili, Mikołaj Bojańczyk, Thomas Colcombet, Dion
Coumans, Mirna Džamonja, Marcel Erné, 
Wesley Fussner, Silvio Ghilardi, Serge Grigorieff, John Harding, Tomáš
Jakl, André Joyal, Alexander Kurz, Vincenzo Marra, Paul-André Melliès, George
Metcalfe, Alessandra Palmigiano, Daniela Petri\c san, Hilary Priestley, 
Carlos Simpson, Howard Straubing, Yde Venema, and Glynn Winskel. We are
also very grateful to the many students who attended our lectures and worked
with the material in this book. This book is primarily written for them.

We dedicate this book to the memory of Bernhard Banaschewski, Bjarni
J{\'o}nsson, and Klaus Keimel.

Finally, we thank our families for their love and support.

\chapter{Order and lattices}\label{ch:order}
In this chapter we first introduce basic notions from order theory: preorders, partial orders, and lattices. We then zoom in on distributive lattices. In the finite case, we prove from first principles a duality theorem, which is a blueprint for the more advanced duality theorems that follow later in this text.

\section{Preorders, partial orders, suprema and infima}\label{sec:orders}
A binary relation $\preceq$ on a set $P$ is called
\begin{itemize}
  \item \emph{reflexive} if $p \preceq p$ for all $p \in P$,
  \item \emph{transitive} if $p \preceq q \preceq r$ implies $p \preceq r$ for all $p,q,r \in P$,
  \item \emph{anti-symmetric} if $p \preceq q$ and $q \preceq p$ imply $p = q$ for all $p, q \in P$,
  \item a \emph{preorder} if it is reflexive and transitive,
  \item a \emphind{partial order} if it is reflexive, transitive and anti-symmetric.
\end{itemize}
A \emph{preordered set}\index{preorder} is a tuple $(P,\preceq)$ with $\preceq$ a preorder on the set $P$. A \emphind{poset}\index{partial order} (short for partially ordered set) is a pair $(P,\leq)$ with $\leq$ a partial order on the set $P$.
Two elements $p$ and $q$ are \emphind{comparable} in a preorder $\preceq$ if at least one of $p \preceq q$ and $q \preceq p$ holds, and \emphind{incomparable} otherwise. The adjective `partial' in `partial order' refers to the fact that not all elements in a partial order are comparable.
A preorder is called \emph{total} or \emph{linear} if any two of its elements are comparable.
A \emphind{total order} or \emphind{linear order} or \emphind{chain} is a total preorder which is moreover anti-symmetric. A poset is called an \emphind{anti-chain} if no distinct elements are comparable.
The \emphind{strict part} of a partial order is the relation $<$ defined by $p < q$ if, and only if, $p \leq q$ and $p \neq q$. Notice that, if we have access to equality, then to specify a partial order $\leq$, it suffices to specify its strict part $<$, from which we can then define $p \leq q$ if, and only if, $p < q$ or $p = q$.
An \emphind{equivalence relation} is a preorder $\preceq$ which is moreover \emph{symmetric}, that is, $p \preceq q$ implies $q \preceq p$ for all $p, q \in P$. In this context, comparable elements are called \emphind{equivalent}.

\begin{example} \index{diamond poset} \index{Boolean algebra on 2 atoms}
 The  \emph{Hasse diagram} of the `diamond' poset $D = \{a,b,c,d\}$, with partial order ${\leq}$ whose strict part is $\{(a,b),(a,c),(a,d),(b,d),(c,d)\}$ is depicted in Figure~\ref{fig:diamond}.
  This partial order is not linear, because we have neither $b \leq c$ nor $c \leq b$.

  \begin{figure}
    \begin{center}
      \begin{tikzpicture}
        \polab{(0,0)}{$b$}{left}
        \polab{(1,-1)}{$a$}{below}
        \polab{(1,1)}{$d$}{above}
        \polab{(2,0)}{$c$}{right}
        \li{(0,0)}{(1,-1)}
        \li{(0,0)}{(1,1)}
        \li{(1,-1)}{(2,0)}
        \li{(1,1)}{(2,0)}
      \end{tikzpicture}
    \end{center}
    \vspace{-5mm}
    \caption{The `diamond' poset $(D,\leq)$}
    \label{fig:diamond}
  \end{figure}
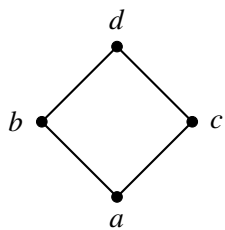

\end{example}
Notice that, in the above example, even though $a \leq d$, we did not draw an edge between $a$ and $d$ in the Hasse diagram. This is due to the fact that $a \leq d$ can be inferred by transitivity from the  order relations $a \leq b$ and $b \leq d$, which are depicted in the diagram. Thus, we only need to draw the `covering' relations in the diagram.

We now give the general definition of \emph{Hasse diagram}.
\begin{definition}\label{def: Hasse}
  For elements $p$ and $q$ of a poset $P$, we say that $q$ \emph{covers} $p$ if $p < q$ and there is no $r\in P$ such that $p < r < q$. We denote this relation by $p\covers q$.%
  \nl{$p \covers q$}{$q$ covers $p$, $q$ is an immediate successor of $p$}

  The elements of a poset are represented in the Hasse diagram as nodes.
  An edge is drawn from a node $p$ to a node $q$ whenever $q$
  covers $p$. In addition, in order not to have to indicate the direction of edges
  by arrows, the convention is that moving up along an edge in the diagram 
  corresponds to moving up in the order. Thus, in particular, points drawn at 
  the same height are incomparable. 
\end{definition}
All finite posets are represented by their Hasse diagrams, as
  well as some infinite ones. However, most infinite posets are not, such as for
  example the usual order on the unit interval, which has an empty covering relation.

There are several interesting classes of maps between preordered sets. (We use the word `map'\index{map} interchangably with `function' throughout this book.)
Let $(P,\preceq_P)$ and $(Q,\preceq_Q)$ be preordered sets and $f \colon P \to Q$ a function. The function $f$ is called
\begin{itemize}
  \item \emphind{order preserving} or \emphind{monotone} if $p \preceq_P p'$ implies $f(p) \preceq_Q f(p')$ for all $p, p' \in P$,
  \item \emphind{order reflecting} if $f(p) \preceq_Q f(p')$ implies $p \preceq p'$ for all $p, p' \in P$,
  \item an \emphind{order-embedding} if it is both order preserving and order reflecting,
  \item an \emphind{order-isomorphism} if it is order preserving and has an order preserving inverse.
\end{itemize}
Note that order-embeddings between posets are always injective, but not all
injective order preserving maps between posets are order-embeddings!
See Exercise~\ref{exe:orderemb}\ref{exe:injnotemb}. A function $f$ between
preordered sets is an order-isomorphism if, and only if, $f$ is a surjective
order-embedding, see Exercise~\ref{exe:orderemb}\ref{exe:surjembisiso}. If
$(P,\preceq_P)$ is a preordered set and $P'$ is a subset of $P$, then the
\emphind{inherited order} on $P'$ is the intersection of $\preceq_P$ with $P' \times
P'$, that is, it is such that the inclusion map $i \colon P' \hookrightarrow P$
is an order-embedding.

An elementary but important operation on preorders is that of `turning upside
down'. If $P$ is a preorder, we denote by $P^\mathrm{op}$\, the \emphind{opposite} of $P$, that is, the preorder with the same underlying set as $P$, but with preorder $\preceq'$ defined by $p \preceq' q$ if, and only if, $q \preceq p$, where $\preceq$ denotes the original preorder on $P$. A function $f \colon P \to Q$ is called \emphind{order-reversing} or \emphind{antitone} if it is order preserving as a function $f \colon P^{\mathrm{op}} \to Q$, that is, if for all $p, p' \in P$, if $p \preceq_P p'$ then $f(p') \preceq_Q f(p)$. An \emphind{anti-isomorphism} between $P$ and $Q$ is, by definition, an isomorphism between $P^{\mathrm{op}}$ and $Q$.
\nl{$P^\op$}{the opposite of a poset $P$}{}

\begin{notationnum}\label{not:function-composition}
Throughout this book, when $f \colon X \to Y$ and $g \colon Y \to Z$ are functions, we write $g \circ f$ for the \emphind{functional composition} of $f$ and $g$, to be read as `$g$ after $f$', that is, $g \circ f$ is the function $X \to Z$ defined by $(g \circ f)(x) = g(f(x))$ for all $x \in X$. We will sometimes omit the symbol $\circ$ and just write $gf$. In Notation~\ref{relation-notations}, we will introduce a slightly different notation for \emphind{relational composition}, as is common in the literature.
\nl{$g \circ f$}{composition of the functions $f$ and $g$, $g$ after $f$}{gcircf} 
\end{notationnum}

\begin{example} \label{exa:finitetotal}
  For any natural number $n$, the finite set $\mathbf{n} := \{0,1,\dots,n-1\}$
    is totally ordered by the usual ordering of natural numbers. 
\end{example}

\begin{example} \label{exa:totalorders}
  The sets of natural numbers $\mathbb{N}$, integers $\mathbb{Z}$, rational numbers $\mathbb{Q}$, and real numbers $\mathbb{R}$, with the usual orders, are total orders. 
\end{example}

\begin{example} \label{exa:preordernaturals}
  On the set of natural numbers $\mathbb{N}$, define a relation $\preceq$ by
  \[ p \preceq q \iff p = 0 \text{ or } (p \neq 0 \text{ and } q \neq 0).\]
  Note that $\preceq$ is a preorder, but not a partial order. We define the \emphind{poset reflection} of this preorder as follows (see Exercise~\ref{exe:reflection} for the general idea). Consider the quotient of $\mathbb{N}$ by the equivalence relation that identifies all non-zero numbers; denote this quotient by $P$, and equip it with the least preorder such that the quotient map $\mathbb{N} \to P$ is order preserving. Then $P$ is a poset, and any other order-preserving function from $\mathbb{N}$ to a poset $(Q, \leq)$ factors through it.
\end{example}

\begin{example}\label{exa:logicequivalence}
  Let $F$ be a set of formulas in some logic with a relation of derivability $\vdash$ between formulas of $F$. More concretely, $F$ can be the set of sentences in a first-order signature and $\vdash$ derivability with respect to some first-order theory. The relation $\vdash$ is rarely a partial order, as there are usually many syntactically different formulas which are mutually derivable in the logic. The poset reflection (see Exercise~\ref{exe:reflection}) consists of the $\vdash$-equivalence classes of formulas in $F$.
\end{example}

\begin{example} \label{exa:sequenceorders}
  Denote by $\mathbf{2}^*$ the set of finite sequences over the two-element set $\mathbf{2} = \{0,1\}$.
  The binary operation of \emphind{concatenation} is defined by juxtaposition of sequences. That is, 
  given sequences $p,r\in \mathbf{2}^*$ of length $n$ and $m$, respectively, $pr$ is the sequence of 
  length $n+m$ whose $i$th entry is the $i$th entry of $p$ if $i\leq n$ and is the $(i-n)$th entry of $r$ otherwise.
  \begin{enumerate}
    \item For $p, q \in \mathbf{2}^*$, define
          \[ p \leq_P q \iff \text{ there exists } r \in \mathbf{2}^* \text{ such that } pr = q.\]
          Note that $\leq_P$ is a partial order on $\mathbf{2}^*$ (see Exercise~\ref{exe:prefixorder}). The poset $(\mathbf{2}^*,\leq_P)$ is called the \emphind{full infinite binary tree}. %
          The partial order $\leq_P$ on $\mathbf{2}^*$ is called the \emphind{prefix order}.\label{exa:binaryprefix}
    \item For $p = (p_1,\dots,p_n)$ and $q = (q_1,\dots,q_m) \in \mathbf{2}^*$, define $p \leq_{\mathrm{lex}} q$ if, and only if, $p$ is a prefix of $q$ or, at the least index such that $p_i \neq q_i$, we have $p_i \leq q_i$ in $\mathbf{2}$.
          Note that $\leq_{\mathrm{lex}}$ is a total order on $\mathbf{2}^*$ (see Exercise~\ref{exe:lexico}). The partial order $\leq_{\mathrm{lex}}$ is called the \emph{lexicographic}\index{lexicographic order} or \emphind{dictionary order} on $\mathbf{2}^*$.\label{exa:binarylexico}
  \end{enumerate}
\end{example}

We define the fundamental notions of supremum and infimum.

\begin{definition}\label{def:supinf}
  Let $(P,\preceq)$ be a preorder. Let $S \subseteq P$.
  \begin{itemize}
    \item an element $s_0$ of $P$ is called a \emphind{lower bound} of $S$ if $s_0 \preceq s$ for all $s \in S$;
    \item an element $s_1$ of $P$ is called an \emphind{upper bound} of $S$ if $s \preceq s_1$ for all $s \in S$;
    \item a lower bound $s_0$ of $S$ is called an \emphind{infimum} or \emphind{greatest lower bound}of $S$ if, for any lower bound $s'$ of $S$, $s' \preceq s_0$;
    \item an upper bound $s_1$ of $S$ is called a \emphind{supremum} or \emphind{least upper bound} of $S$ if, for any upper bound $s'$ of $S$, $s_1 \preceq s'$.
  \end{itemize}
  In the special case where $S = \emptyset$, an element $s_0$ which is a supremum of $S$ is called a \emphind{bottom} or \emph{minimum} element of $P$, meaning that $s_0 \preceq s$ for all $s \in S$. Similarly, an element $s_1$ which is an infimum of $S$ is called a \emphind{top} or \emph{maximum} element of $P$.
\end{definition}
In a poset, any set has at most one infimum and at most one supremum
(see Exercise~\ref{exe:infsupunique}). If a unique infimum of a subset $S$ exists,
it is denoted by $\bigwedge S$ and is also known
as the \emphind{meet} of $S$. The supremum of $S$, if it exists uniquely, is
denoted by $\bigvee S$ and is known as the
\emphind{join}. 
In the case where $S =\{a,b\}$, we also write $a \wedge b$ and $a \vee b$, and if $S =
\{a_1,\dots,a_n\}$ we write $a_1 \vee \cdots \vee a_n$ and $a_1 \wedge \cdots
\wedge a_n$. The bottom element, if it exists, is denoted by $\bot$ or $0$, and the top element by $\top$ or $1$. If $S$ is a subset of a poset $P$, we denote the set of \emph{maximal} elements in $S$ by $\max(S)$\index{max}\index{min}; that is,
$$\max(S) := \{s \in S \ | \ \text{ for all } s' \in P, \text{ if } s \leq s' \text{ and } s' \in S, \text{ then } s' = s\}.$$
Similarly, the set of minimal elements in $S$ is denoted by $\min(S)$. In contrast to maximal elements of a set, the supremum of a set does not need to belong to the set itself. Note that non-empty subsets of a poset may not have any minimal or maximal elements; see the examples below. Also, postulating the existence of maximal or minimal elements in certain posets is related to choice principles; see our discussion of \emph{Zorn's Lemma}, Lemma~\ref{lem:zorn}, in the next chapter.

\begin{remark}\index{maximal vs. maximum}\index{supremum!vs. maximum}
  There are subtle but important differences between the words `maximum', `maximal' and `supremum'. An element is \emph{maximal} in a subset $S$ of a poset if there is no other element in $S$ that lies strictly above it, while it is a \emph{maximum} element in $S$ if all other elements of $S$ lie below it. Note that in a \emph{totally} ordered set, the concepts maximal and maximum are equivalent, but not in general. Finally, an important distinction between these two concepts and that of supremum is that for an element to be a supremum, it is \emph{not} needed that it lies in the set itself, while this is part of the definition for maximal and maximum elements, see Exercise~\ref{exe:umalsup}. For this reason, the supremum of a set S depends on the ambient poset, see Exercise~\ref{exe:counterexamples-complete}.
\end{remark} %

Infima and suprema may fail to exist. There are three different situations in which this can happen: a set can either have no lower (or upper) bounds at all, or its set of lower (or upper) bounds has incomparable maximal (or minimal) elements, or for infinite sets, the set of upper bounds may be non-empty but not have all such above a minimal upper bound (or a non-empty set of lower bounds not all witnessed by maximal lower bounds).

We illustrate the above ideas with three examples.
\begin{example}\label{exa:butterfly}
  In the poset $(P,\leq)$ whose Hasse diagram is depicted below, the set $S = \{a,b\}$ does not have an infimum, because $c$ and $d$ are incomparable maxim\emph{al} lower bounds of $S$, and hence neither is a maxim\emph{um} lower bound.
  \index{butterfly poset}

  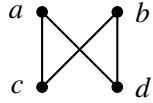
\begin{figure}[htp]
    \begin{center}
      \begin{tikzpicture}
        \polab{(0,1)}{$a$}{left}
        \polab{(1,1)}{$b$}{right}
        \polab{(0,0)}{$c$}{left}
        \polab{(1,0)}{$d$}{right}
        \li{(0,0)}{(0,1)}
        \li{(0,0)}{(1,1)}
        \li{(1,0)}{(0,1)}
        \li{(1,0)}{(1,1)}
      \end{tikzpicture}
    \end{center}
    \vspace{-5mm}
    \caption{The `butterfly' poset $(P,\leq)$}
    \label{fig:butterfly}
  \end{figure}
\end{example}

\begin{example}\label{exa:naturalsinfsup}
  In the set $\mathbb{N}$ of natural numbers with its usual total order, any subset has an infimum, which is in fact a minimum, but the only subsets having a supremum are the finite subsets. For any finite subset, the supremum is in fact a maximum.
\end{example}

\begin{example}\label{exa:rationalsinf}
  In the set $\mathbb{Q}$ of rational numbers with its usual total order, the subset $\{\frac{1}{n} \ | \ n \in \mathbb{N}_{\geq 1}\}$ has an infimum, $0$, but it does not have a minimum. Furthermore $\{q\in \mathbb{Q}\mid q\leq \sqrt{2}\}$ has upper bounds but no least one.
\end{example}

We end this section by introducing the important concepts of \emph{adjunction} and \emph{Galois connection} between preordered sets.\index{adjunction!between preorders}
\begin{definition}\label{dfn:poset-adjunction-def}
  Let $(P,\preceq_P)$ and $(Q,\preceq_Q)$ be preordered sets, and let $f \colon P \to Q$ and $g \colon Q \to P$ be functions. The pair $(f,g)$ is called an \emph{adjunction}, with $f$ the \emph{left} or \emph{lower adjoint} and $g$ the \emph{right} or \emph{upper adjoint}, provided that $f$ and $g$ are both order preserving, and for every $p \in P$ and $q \in Q$,\index{adjoint!between preorders}
  \[ f(p) \preceq_Q q \text{ if, and only if, } p \preceq_P g(q).\]
  An adjunction between $P^\op$ and $Q$ is called a \emphind{Galois connection} or \emphind{contravariant adjunction}.
\end{definition}

The notion of adjunction between (pre)orders is very important and will come into play in many places in this book. It is also a useful precursor to the concept of adjunction between categories that we will encounter later in Definition~\ref{dfn:adjunction}. Exercises~\ref{exe:adjunctions}~and~\ref{exe:adjointexistsiff} collect some important basic facts about adjunctions between preorders, which will be used throughout the book. 

\begin{example}
  Let $f \colon \mathbb{Z} \into \mathbb{Q}$ be the order-embedding which sends each integer $x \in \mathbb{Z}$ to itself, regarded as a rational number. The map $f$ has a right adjoint, $g$, which sends each rational $y \in \mathbb{Q}$ to its \emph{floor}, that is, $g(y)$ is the largest integer below $y$. The map $f$ also has a left adjoint, which sends a rational $y \in \mathbb{Q}$ to its \emph{ceiling}, that is, the smallest integer above $y$.
\end{example}

\emph{Contravariant} adjunctions occur particularly often in mathematics since, as we will see in the following example, they arise naturally any time we have a binary relation between two sets. 

\begin{example}\label{exa:galoisconnection}
  Fix a relation $R \subseteq X \times Y$ between two sets. For any $a \subseteq X$ and $b \subseteq Y$, define the sets $u(a) \subseteq Y$ and $\ell(b) \subseteq X$ by
  \begin{align*}
    u(a)    & := \{y \in Y \ \mid \ \text{for all } x \in a, x {R} y \}, \\
    \ell(b) & := \{ x \in X \ \mid \ \text{for all } y \in b, x {R} y \}
  \end{align*}
  The pair of functions $u \colon \mathcal{P}(X) \leftrightarrows \mathcal{P}(Y) \colon \ell$ is a Galois connection between the posets $(\mathcal{P}(X), \subseteq)$ and $(\mathcal{P}(Y), \subseteq)$, that is, for any $a \subseteq X$ and $b \subseteq Y$, we have $b \subseteq u(a)$ if, and only if, $a \subseteq \ell(b)$.
\end{example}

The name Galois connection refers to the work of Galois in the theory of rings and fields in which the central object of study is a Galois connection induced by the binary relation between subfields and automorphisms of a field given by the automorphism being the identity on the subfield. While this example is historically important, it is not so central to the topics of this book. So if you are not familiar with Galois Theory, we provide the following two classical examples, in order theory and logic, respectively.
First, in the special case when $R$ is a pre-order on a set $X$, $u(a)$ is the set of \emph{common upper bounds} for the elements of $a$, and $\ell(b)$ is the set of \emph{common lower bounds} for the elements of $b$.
Second, the Galois connection between theories and model classes studied in logic is also a special case of the Galois connection $(u,\ell)$, as follows. Suppose that $S$ is a set of structures, $F$ is a set of logical formulas, and suppose we are given a relation of `interpretation', $\models$ from $S$ to $F$, where, for $M \in S$ and $\phi \in F$, the relation $M \models \phi$ is read as ``$\phi$ holds in $M$''. Then in the  Galois connection of Example~\ref{exa:galoisconnection}, $u$ sends a set of models $a$ to its \emph{theory}, that is, the set of formulas that hold in every model of $a$, and $\ell$ sends a set of formulas $b$ to its \emph{class of models}, that is, the set of models in which every formula from $b$ holds.
Finally, not only are Galois connections obtained from a binary relation as in Example~\ref{exa:galoisconnection} omnipresent in mathematics, but, as we will see later, by topological duality theory, all Galois connections between distributive lattices are of this form; see in particular Proposition~\ref{prop:priestley-2-eq} and Exercise~\ref{exe:adjunctionpreservedbyduality}.

\ourexercises

\begin{ourexercise}\label{exe:hasse}
  Sketch the Hasse diagrams for the preorders described in Example~\ref{exa:finitetotal}, Example~\ref{exa:preordernaturals} and Example~\ref{exa:sequenceorders}(\ref{exa:binaryprefix}).
\end{ourexercise}

\begin{ourexercise}\label{exe:prefixorder}\index{free monoid}
  For any set $A$, let $A^*$ denote the set of finite sequences of elements of $A$. Prove that the relation $\leq_P$ on $A^*$ defined by
  \[ u \leq_P v \iffdef \text{ there exists } w \in A^* \text{ such that } uw = v\]
  is a partial order.

 \emph{Note.} This is a special case of the opposite of the so-called Green pre-order $\leq_{\mathcal{R}}$, which exists on any monoid.
\end{ourexercise}

\begin{ourexercise}\phantomsection\label{exe:lexico}\index{lexicographic order}
  Consider the relation $\leq_{\mathrm{lex}}$ on $\mathbf{2}^*$ defined in Example~\ref{exa:sequenceorders}\ref{exa:binarylexico}.
  \begin{enumerate}
    \item Prove that $\leq_{\mathrm{lex}}$ is a total order.
    \item Prove that, even though $\mathbf{2}^*$ is infinite, the total order $\leq_{\mathrm{lex}}$ is the transitive closure of its covering relation. That is, show that $p\leq_{\mathrm{lex}}q$ if, and only if, there are $r_0,\dots, r_n$ with $p\covers r_0\covers\dots r_n\covers q$.
  \end{enumerate}
\end{ourexercise}

\begin{ourexercise}\phantomsection\label{exe:orderemb}
  \begin{enumerate}
    \item Give an example of an injective order-preserving map between posets which is not an order-embedding.\label{exe:injnotemb}\index{order-preserving!and injective}
    \item Prove that a surjective order-embedding between posets is an order-isomorphism.\label{exe:surjembisiso}\index{order-embedding!surjective}
  \end{enumerate}
\end{ourexercise}

\begin{ourexercise}\label{exe:reflection}
  If $(P,\preceq)$ is a preordered set, define
  \[ p \equiv q \iff p \preceq q \text{ and } q \preceq p.\]
  \begin{enumerate}
    \item Prove that $\equiv$ is an equivalence relation on $P$.
    \item Prove that there is a well-defined \emph{smallest} partial order $\leq$ on the quotient set $P/{\equiv}$ such that the quotient map $f \colon P \to P/{\equiv}$ is order preserving.
    \item Prove that, for any order-preserving $g \colon P \to Q$ with $Q$ partially ordered, there exists a unique order-preserving $\overline{g} \colon P/{\equiv} \to Q$ such that $\overline{g} \circ f = g$.
  \end{enumerate}
  The partial order $P/{\equiv}$ defined in this exercise is called the \emphind{poset reflection} of the preorder $P$.
\end{ourexercise}

\begin{ourexercise}\label{exe:infsupunique}
  Let $(P,\preceq)$ be a preorder and $S \subseteq P$.
  \begin{enumerate}
    \item Prove that if $s_0$ and $s_0'$ are both infima of $S$, then $s_0 \preceq s_0'$ and $s_0' \preceq s_0$.
    \item Conclude that in a partial order, any set has at most one supremum and at most one infimum.
  \end{enumerate}
\end{ourexercise}

\begin{ourexercise}\label{exe:umalsup}
  Draw a graph with three nodes labeled `maximum', `maximal', and `supremum', and directed edges denoting that the existence of one implies the existence of the other. Do any more implications hold in finite posets? In totally ordered sets? In finite totally ordered sets?
\end{ourexercise}

\begin{ourexercise}\label{exe:adjunctions}
  Let $(P,\preceq_P)$ and $(Q,\preceq_Q)$ be preordered sets and $f \colon P \leftrightarrows Q \colon g$ a pair of order-preserving maps between them.
  \begin{enumerate}
    \item Prove that $(f,g)$ is an adjunction if, and only if, for every $p \in P$, $p \preceq_P gf(p)$, and for every $q \in Q$, $fg(q) \preceq_Q q$.
  \end{enumerate}
  For the rest of this exercise, assume that $(f,g)$ is an adjunction.
  \begin{enumerate}[resume]
    \item Prove that $fgf(p) \equiv f(p)$ and $gfg(q) \equiv g(q)$ for every $p \in P$ and $q \in Q$.
          \item\label{itm:equalities} Conclude that, in particular, if $P$ and $Q$ are posets, then $fgf = f$ and $gfg = g$.
    \item \label{itm:minimumimage} Prove that, if $P$ is a poset, then for any $p \in P$, $gf(p)$ is the minimum element above $p$ that lies in the image of $g$.
    \item Formulate and prove a similar statement to the previous item about $fg(q)$, for $q \in Q$.
    \item \label{itm:leftadjointpreservessup} Prove that, for any subset $S \subseteq P$, if the supremum of $S$ exists, then $f\big(\bigvee S\big)$ is the supremum of the direct image $f[S]$.
    \item Prove that, for any subset $T \subseteq Q$, if the infimum of $T$ exists, then $g\big(\bigwedge T\big)$ is the infimum of $g[T]$.
  \end{enumerate}
  In words, the last two items say that \emph{lower adjoints preserve existing suprema} and \emph{upper adjoints preserve existing infima}. In Exercise~\ref{exe:adjointexistsiff} of the next section we will see that a converse to this statement holds in the context of complete lattices.
\end{ourexercise}

\begin{ourexercise}\label{exe:adjunctioninjectivesurjective}
  Let $f \colon P \leftrightarrows Q \colon g$ be an adjunction between posets. Show that $f$ is injective if, and only if, $g$ is surjective, and that $f$ is surjective if, and only if, $g$ is injective. \hint{Part (\ref{itm:equalities}) of Exercise~\ref{exe:adjunctions} can be useful here.}
\end{ourexercise}

\section{Lattices}\label{sec:lattices}
A \emph{(bounded) lattice}\index{lattice!order-theoretic definition} is a partially ordered set $L$ in which every finite subset has a supremum and an infimum.  In fact, to be a lattice, it is sufficient that the empty set and all two-element sets have suprema and infima (see Exercise~\ref{exe:lattsuff}). A \emphind{complete lattice} $C$ is a partially ordered set in which every subset has a supremum and an infimum. In fact, for a partially ordered set to be a complete lattice, it is sufficient that every subset has a supremum (see Exercise~\ref{exe:complattsuff}).\footnotemark

An interesting equivalent definition of lattices is the following. A \emph{lattice}\index{lattice!algebraic definition} is a tuple $(L,\vee,\wedge,\bot,\top)$, where $\vee$ and $\wedge$ are binary operations on $L$ (that is, functions $L \times L \to L$), and $\bot$ and $\top$ are elements of $L$ such that the following axioms hold:
\begin{enumerate}
  \item the operations $\vee$ and $\wedge$ are \emphind{commutative}, that is, ${a \vee b = b \vee a}$ and ${a \wedge b = b \wedge a}$ for all $a, b \in L$;
  \item the operations $\vee$ and $\wedge$ are \emphind{associative}, that is, ${(a \vee b) \vee c = a \vee (b \vee c)}$ and ${(a \wedge b) \wedge c = a \wedge (b \wedge c)}$ for all $a, b, c \in L$;
  \item the operations $\vee$ and $\wedge$ are \emphind{idempotent}, that is, $a \vee a = a$ and $a \wedge a = a$ for all $a \in L$;
  \item the \emphind{absorption laws} $a \wedge (a \vee b) = a$ and $a \vee (a \wedge b) = a$ hold for all $a, b \in L$;
  \item the element $\bot$ is \emph{neutral for $\vee$} and the element $\top$ is \emph{neutral for $\wedge$}, that is, $\bot \vee a = a$ and $\top \wedge a = \top$ for all $a \in L$.%
\end{enumerate}
Given a lattice~$(L,\vee,\wedge,\bot,\top)$ according to this algebraic definition, define
\begin{equation}\label{eq:order-from-lattice}
  a \leq_L b \iff a \wedge b = a.
\end{equation}
Then $\leq_L$ defines a partial order on the set $L$ which makes $L$ into a lattice according to the order-theoretic definition, and the binary infimum and supremum are given by $\wedge$ and $\vee$, respectively. Conversely, given a lattice~${(L,\leq)}$ according to the order-theoretic definition, it is easy to check that the operations of binary join ($\vee$), binary meet ($\wedge$), and the elements $\top$ and $\bot$ make $L$ into a lattice according to the algebraic definition, and that $\leq$ is given by (\ref{eq:order-from-lattice}), which is then also equivalent to $a \vee b = b$. The somewhat tedious but instructive Exercise~\ref{exe:latticedefs} asks you to verify the claims made in this paragraph.\index{lattice!order from operations}

A \emphind{semilattice} is a structure $(L, \cdot, 1)$ where $\cdot$ is a
commutative, associative and idempotent binary operation, and $1$ is a neutral element for the
operation $\cdot$. The operation $\cdot$ can then either be seen as the
operation $\wedge$ for the partial order defined by $a \leq b$ if, and only if,
$a \cdot b = a$, or as the operation $\vee$ for the opposite partial order defined by $a
\leq b$ if, and only if, $a \cdot b = b$. When $(L, \vee, \wedge, \bot, \top)$
is a lattice, we call $(L,\vee,\bot)$ and $(L,\wedge,\top)$ the
\emphind{join-semilattice} and \emphind{meet-semilattice} \emph{reducts}\index{reduct!semilattice} of $L$,
respectively.

\begin{remark}
  Many authors use the word `lattice' for posets with all binary infima and suprema, and then `bounded lattices' are those that also have $\top$ and $\bot$. We require that all finite infima and suprema exist. Note that the non-empty finite infima and suprema are guaranteed to exist as soon as binary ones do (Exercise~\ref{exe:lattsuff}), while the empty infimum and supremum are just the $\top$ and $\bot$, respectively. In duality theory bounds are very convenient, if one does not have them, one should simply add them. Accordingly, we suppress the adjective `bounded' and use just `lattice' for the bounded ones and we will only specify when once in a while we have an unbounded lattice, sublattice, or lattice homomorphism.
\end{remark}
\subsection*{Homomorphisms, products, sublattices, quotients}
We briefly recall a few basic algebraic notions that we will need. For more information, including detailed proofs of these statements, we refer the reader to a textbook on universal algebra, such as for example \cite{BurSan2000,Wechler1992}.
A function $f \colon L \to M$ between lattices is called a \emphind{lattice homomorphism} if it preserves all the lattice operations; that is, $f(\bot_L) = \bot_M$, $f(\top_L) = \top_M$, and $f(a \vee_L b) = f(a) \vee_M f(b)$, $f(a \wedge_L b) = f(a) \wedge_M f(b)$ for all $a, b \in L$. 
Lattice homomorphisms are always order preserving, and injective lattice homomorphisms are always order-embeddings (see Exercise~\ref{exe:injlatthom}). We also call an injective homomorphism between lattices a \emphind{lattice embedding}.\index{homomorphism!between lattices} Similarly, bijective lattice homomorphisms are always order isomorphisms, and we call an bijective homomorphism between lattices a \emphind{lattice isomorphism}.\index{isomorphism!between lattices}

A simple induction (Exercise~\ref{exe:lattsuff}) shows that, if a function $f \colon L \to M$ between lattices preserves $\bot$ and $\vee$, then it preserves all finite joins; we say that such a function \emph{preserves finite joins}\index{finite-join-preserving}\index{semilattice homomorphism}, or also that it is a \emph{homomorphism for the join-semilattice reducts}. Note that the statement that $f$ preserves finite joins does not in general imply anything about the preservation of any other suprema that may exist in $L$. As a rule, whenever we write `$f$ preserves joins', even though we may sometimes omit the adjective `finite', we still only refer to preservation of finite suprema, as these are generally the only suprema that exist in a lattice. If we want to signal that $f$ preserves other suprema than just the finite ones, then we will always take care to explicitly say so. Symmetrically, a function $f \colon L \to M$ that preserves $\top$ and $\wedge$ is called \emphind{finite-meet-preserving}, and the same remarks about the adjective `finite' apply here.

The \emphind{Cartesian product}\index{product} of an indexed family $(L_i)_{i \in I}$ of lattices is the lattice structure on the product set $L := \prod_{i \in I} L_i$ given by pointwise operations; for example, $(\bot_L)_i = \bot_{L_i}$ and $(\top_L)_i = \top_{L_i}$ for every $i \in I$, and if $a = (a_i)_{i \in I}, b = (b_i)_{i \in I} \in L$ then $(a \vee b)_i = a_i \vee_{L_i} b_i$ and $(a \wedge b)_i = a_i \wedge_{L_i} b_i$. In this way, $L$ becomes a lattice, whose partial order is also the product order, that is, $a \leq_L b \iff a_i \leq_{L_i} b_i$ for every $i \in I$, and each projection map $\pi_i \colon L \onto L_i$ is a surjective homomorphism (see Exercise~\ref{exe:product-lattices}).
A \emph{sublattice}\index{sublattice} of a lattice $M$ is a subset $M'$ such that, $\bot$ and $\top$ are in $M'$ and for every $a, b \in M'$, both $a \vee b$ and $a \wedge b$ are in $M'$.
In this case, $M'$ is a (bounded) lattice in its own right, and the inclusion map $i \colon M' \into M$ is a lattice homomorphism. Also, the \emphind{direct image} or \emphind{forward image}\nl{$f[S]$}{forward image under a function $f$ of a subset $S$ of the domain}{fS}, denoted $f[L]$, of any lattice homomorphism $f \colon L \to M$ is a sublattice of the codomain, and if $f$ is moreover an order embedding, then the domain lattice $L$ is isomorphic to this direct image $f[L]$. If $L, L'$ are lattices and there exists a surjective homomorphism $f \colon L \onto L'$, then we say that $L'$ is a \emphind{homomorphic image} of $L$. An \emphind{unbounded sublattice} of a lattice $M$ is a subset $M'$ such that, for any $a, b \in M'$, both $a \vee b$ and $a \wedge b$ in $M'$. Note that an unbounded sublattice may or may not have a top and bottom element, and even if it does, these need not coincide with the top and bottom elements of $M$. What we call sublattice here is sometimes called \emph{bounded} sublattice, and the term sublattice then more generally refers to a subset closed under binary, but not necessarily empty joins and meets. We mostly consider sublattices with bounds in this book, unless noted explicitly otherwise.

A subset $L'$ of a lattice $L$ may fail to be a sublattice of $L$, even if it is a bounded lattice when equipped with the partial order inherited from $L$: the value of a join or meet may change when moving to a subset. In a similar vein, if $L$ is a complete lattice, then a sublattice $L'$ of $L$ may be a complete lattice in itself, while it is not necessarily the case that the supremum of any subset of $L'$ coincides with the supremum of the same subset in $L$. These are important distinctions; Exercise~\ref{exe:counterexamples-complete} suggests examples that show the difference. When $L$ is a complete lattice, we will reserve the term \emphind{complete sublattice} for a subset $L'$ of $L$ such that, for any subset $S$ of $L'$, the supremum of $S$ in $L$ and the infimum of $S$ in $L$ both also belong to $L'$. The above remarks then imply that a sublattice which is a complete lattice may fail to be a complete sublattice.

A \emph{congruence}\index{congruence!on a lattice}\index{quotient!of a lattice} on a lattice $L$ is an equivalence relation $\theta \subseteq L \times L$ such that, for any two pairs $(a,a') \in \theta$ and $(b,b') \in \theta$, the pairs $(a \vee b, a' \vee b')$ and $(a \wedge b, a' \wedge b')$ are both also in $\theta$. If $\theta$ is a congruence on a lattice $L$, then the quotient set $L/{\theta}$ carries a unique lattice structure which makes the quotient map $p \colon L \to L/{\theta}$ into a lattice homomorphism. Indeed, the operations $[a]_{\theta} \vee [b]_{\theta} := [a \vee b]_{\theta}$ and $[a]_{\theta} \wedge [b]_{\theta} := [a \wedge b]_{\theta}$ give a well-defined lattice structure on the set $L/{\theta}$, with bottom element $[\bot]_{\theta}$ and top element $[\top]_{\theta}$. If $f \colon L \to M$ is any lattice homomorphism, the \emphind{kernel} of $f$ is the equivalence relation defined by
\nl{$\ker f$}{kernel of a homomorphism $f$}{}
\[\ker f := \{(a,a') \in L \times L \ | \ f(a) = f(a')\},\]
which is a congruence on $L$ with the property that the quotient lattice
$L/{\ker f}$ is isomorphic to the direct image of $f$. In particular, if $f
\colon L \onto M$ is a surjective homomorphism, then the codomain $M$ is
isomorphic to $L/{\ker f}$. In other words, the homomorphic images of $L$ are,
up to isomorphism, all the quotients of $L$. These facts together are known as
the \emphind{first isomorphism theorem}\index{lattice!first isomorphism
theorem} for lattices: any lattice homomorphism $f \colon L \to M$ can be
factored as a surjective homomorphism followed by an embedding. Indeed, we have
that $f = e \circ p$, where $p \colon L \to L/{\ker f}$ is the quotient, and $e
\colon L/{\ker f} \to M$ is the embedding of the direct image of
$f$.\label{firstiso}

\subsection*{Distributivity}
A lattice $L$ is called \emph{distributive}\index{distributivity}\index{lattice!distributive} if
\begin{equation}\label{eq:dist1}
  \text{ for all } a, b, c \in L, \quad a \wedge (b \vee c) = (a \wedge b) \vee (a \wedge c),
\end{equation}
or, equivalently,
\begin{equation} \label{eq:dist2}
  \text{ for all } a, b, c \in L, \quad a \vee (b \wedge c) = (a \vee b) \wedge (a \vee c).
\end{equation}
The proof that (\ref{eq:dist1}) and (\ref{eq:dist2}) are indeed equivalent is left as Exercise~\ref{exe:disteq}. In particular, a lattice $L$ is distributive if, and only if, its opposite $L^\op$ is distributive; we say that distributivity is a \emphind{self-dual property}. Also, products, sublattices and homomorphic images of distributive lattices are again distributive; this is an immediate consequence of the fact that distributivity is described equationally, and can also be proved directly (see Exercise~\ref{exe:product-lattices} for the product case).

An easy inductive argument shows that in a distributive lattice $L$, we have:
\begin{align*}
  \text{ for any $a \in L$ and $F \subseteq L$ finite, } a \wedge \bigvee F = \bigvee_{b \in F} (a \wedge b),
\end{align*}
and, again equivalently,
\begin{align*}
  \text{ for any $a \in L$ and $F \subseteq L$ finite, } a \vee \bigwedge F = \bigwedge_{b \in F} (a \vee b).
\end{align*}

If $L$ is a \emph{complete} lattice, we say that the \emphind{Join Infinite Distributive law} (JID) holds in $L$ if
\begin{align} \label{eq:JID}
  \text{ for any $a \in L$ and $S \subseteq L$, } a \wedge \bigvee S = \bigvee_{b \in S} (a \wedge b).
\end{align}
A complete lattice in which the Join Infinite Distributive law holds is called a \emphind{frame}. We will encounter frames in Chapter~\ref{chap:Omega-Pt} and the chapters following it. Note that a distributive lattice may be complete, but fail to be a frame (see Exercise~\ref{exe:counterexamples-complete}.\ref{ite:compDLnotframe}).

A frame is the same thing as a \emphind{complete Heyting algebra}; we will encounter Heyting algebras in Section~\ref{sec:esakia-heyting} in Chapter~\ref{ch:methods}. Note however that the natural structure-preserving maps for frames and for complete Heyting algebras may differ: a \emphind{frame homomorphism} is required to preserve arbitrary joins and finite meets, while a \emphind{Heyting homomorphism} is required to preserve finite joins, finite meets, and the Heyting implication (see Exercise~\ref{exe:morphismsHAdifferent}).

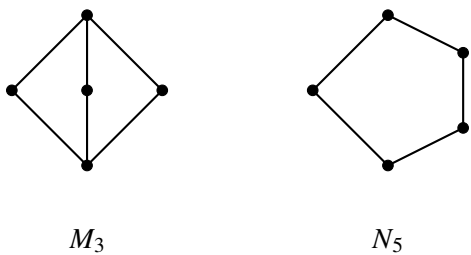
\begin{figure}[htp]
  \begin{center}
    \begin{tikzpicture}
      \po{(-2,0)}
      \po{(-3,1)}
      \po{(-2,1)}
      \po{(-1,1)}
      \po{(-2,2)}
      \draw[thick] (-2,0) -- (-3,1) -- (-2,2) -- (-1,1) -- (-2,0) -- (-2,2);
      \node at (-2,-1) {$M_3$};
      \po{(2,0)}
      \po{(1,1)}
      \po{(3,.5)}
      \po{(3,1.5)}
      \po{(2,2)}
      \node at (2,-1) {$N_5$};
      \draw[thick] (2,0) -- (1,1) -- (2,2) -- (3,1.5) -- (3,.5) -- (2,0);
    \end{tikzpicture}
    \vspace{-5mm}
  \end{center}
  \caption{The lattices $M_3$ and $N_5$.}
  \label{fig:nondist}
\end{figure}

There are two `minimal' counterexamples to distributivity,\index{distributivity!failure of}\index{M3@$M_3$}\index{N5@$N_5$} namely the non-distributive lattices $M_3$ and $N_5$, depicted in Figure~\ref{fig:nondist}. Indeed, the following proposition, which you will be asked to prove in Exercise~\ref{exe:distchar}, characterizes distributive lattices in terms of `forbidden substructures'.

\begin{proposition}\label{prop:distchar}
  Let $L$ be a lattice. Then $L$ is distributive if, and only if, $L$ does not contain an unbounded sublattice which is isomorphic to $M_3$ or $N_5$.
\end{proposition}
Note that Proposition~\ref{prop:distchar} does not require the existence of a sublattice (with top and bottom) isomorphic to $M_3$ or $N_5$. Stated more explicitly, `$L$ contains an unbounded sublattice isomorphic to $M_3$' means: there are five distinct elements $a,b,c,d,e$ in $L$ such that $a = b \wedge c = b \wedge d = c \wedge d$, $e = b \vee c = b \vee d = c \vee d$, but $a$ need not be the bottom element of $L$, and $e$ need not be the top element of $L$. We leave it to the reader to similarly spell out the definition of `$L$ contains an unbounded sublattice isomorphic to $N_5$'.%

\subsection*{Directed and filtering sets}
\label{def:directed} Let $P$ be a poset. A subset $D \subseteq P$ is called \emphind{directed} if it is non-empty, and for any $p, q \in D$, there exists $r \in D$ such that $r \geq p$ and $r \geq q$.
Equivalently, $D$ is directed if any finite subset of $D$ has an upper bound in $D$ (see Exercise~\ref{exe:directed}).  A \emphind{directed join} is the supremum of a directed set.
Order-dually, a subset $F \subseteq P$ is called \emphind{filtering} if it is non-empty, and for any $p, q \in F$, there exists $r \in F$ such that $r \leq p$ and $r \leq q$.
Again, equivalently, $F$ is filtering if any finite subset of $F$ has a lower bound in $F$.
We will sometimes use the word \emphind{up-directed} as a synonym for `directed', and \emphind{down-directed} as a synonym for `filtering'.

In lattice theory, directed and filtering subsets of a lattice often appear, and indeed, certain filtering sets called \emph{prime filters} are central to the duality theory in Chapter~\ref{ch:priestley} (also see Exercise~\ref{exe:idealdirected}).  Directed and filtering sets are also important in topology and domain theory, as we will see in Chapters~\ref{chap:TopOrd},~\ref{chap:Omega-Pt}, and \ref{chap:DomThry}. In particular, we will often encounter the notion of a \emphind{filtering collection} of subsets of a set $X$, which is just a filtering subset of the poset $(\mathcal{P}(X), \subseteq)$. More explicitly, a non-empty collection $\mathcal{F}$ of subsets of a set $X$ is \emph{filtering} if, for any $S, T \in \mathcal{F}$, there exists $R \in \mathcal{F}$ such that $R \subseteq S \cap T$.

A poset $P$ is \emph{directedly complete} \index{directedly complete poset} if any directed subset of $P$ has a supremum in $P$. Directedly complete posets are called \emphind{dcpo}'s, for short. Directed joins are what is needed to complete a lattice into a complete lattice (see Exercise~\ref{exe:directed}). A categorified version of this result will be proved in Proposition~\ref{prop:complete-category}.

\subsection*{Complements and Boolean algebras}
If $a$ is an element in a lattice, an element $b$ is called a \emph{complement} of $a$ if $a \wedge b = \bot$ and $a \vee b = \top$. A \emphind{Boolean algebra} is a distributive lattice in which every element has a complement. The complement of an element in a distributive lattice is unique if it exists, see Exercise~\ref{exe:distuncomp}.a. If $L$ is a Boolean algebra, we denote by $\neg a$ the unique complement of an element $a$.

More directly, a \emph{Boolean algebra} is a tuple $(B,\wedge,\vee,\bot,\top,\neg)$ such that $(B,\wedge,\vee,\bot,\top)$ is a distributive lattice, and for all $a \in B$, $a \wedge \neg a = \bot$ and $a \vee \neg a = \top$. The term `Boolean algebra' comes from another (the original) equational definition: a Boolean algebra is term-equivalent to a  commutative ring with unit in which all elements are idempotent; one only has to change which operations are viewed as basic (see Exercise~\ref{exe:BAeq}).

Regarding maps between Boolean algebras, for any lattice homomorphism $f \colon L \to M$, where $L$ and $M$ are Boolean algebras, the function $f$ must also preserve the operation $\neg$, that is, $f(\neg a) = \neg f(a)$ for all $a \in L$ (see Exercise~\ref{exe:distuncomp}.b). We thus have an unambiguously defined notion of \emph{homomorphism} between Boolean algebras.\index{homomorphism!between Boolean algebras}
However, a Boolean algebra in general has many sublattices that are not themselves Boolean algebras. We call a sublattice $L$ of a Boolean algebra $B$ a \emph{subalgebra}\index{subalgebra!of a Boolean algebra} if it is closed under the operation $\neg$.

As a subclass of distributive lattices, Boolean algebras take up a very special position: every distributive lattice has a ``minimal'' Boolean algebra sitting around it, called its \emphind{Boolean envelope} or \emphind{free Boolean extension}\footnote{In some literature, this object is also called the \emphind{Booleanization}, but this term has also been used with other meanings, so we will avoid it.}. In categorical terms, Boolean algebras form a \emph{full reflective subcategory} of distributive lattices; see Chapter~\ref{ch:categories}. We end this section by spelling out what this means.
\nl{$\bar{f}$}{unique extension of a function $f$ given by a universal property determined by the context}{}
\begin{definition}\label{def:booleanenvelope}
  Let $L$ be a distributive lattice. A Boolean algebra $B$, together with an injective homomorphism $e \colon L \into B$, is called a \emphind{Boolean envelope} of $L$ if, for every lattice homomorphism $h \colon L \to A$, with $A$ a Boolean algebra, there exists a unique homomorphism $\bar{h} \colon B \to A$ such that $\bar{h} \circ e = h$, , that is, such that the following diagram commutes:
  \begin{center}
    \begin{tikzpicture}
      \matrix (m) [matrix of math nodes,row sep=3em,column sep=3em,minimum width=3em]
      {
        B & A \\
        L &   \\
      };
      \path[-stealth]
      (m-1-1) edge node [above] {$\bar{h}$} (m-1-2)
      (m-2-1) edge node [left] {$e$} (m-1-1)
      (m-2-1) edge node [below] {$h$} (m-1-2);
    \end{tikzpicture}
  \end{center}
\end{definition}
In Proposition~\ref{prop:finite-boolenv} below, we will construct the Boolean envelope of a finite distributive lattice as a small application of the finite duality that we prove there. In fact, any distributive lattice has a Boolean envelope, and it is unique up to isomorphism; see Propositions~\ref{prop:boolenv-unique}~and~\ref{prop:boolenv} in Chapter~\ref{ch:priestley}. We will denote by $L^-$ the (unique up to isomorphism) Boolean envelope of a distributive lattice $L$.
\nl{$L^-$}{Boolean envelope of a distributive lattice $L$}{}
We will also see there (Corollary~\ref{cor:boolenv-hull}) that the Boolean envelope of $L$ has a further special property, namely that for any embedding $h \colon L \into A$ with $A$ a Boolean algebra, the Boolean subalgebra of $A$ generated by the image of $h$ is isomorphic to the Boolean envelope of $L$. That is, if $h$ is injective, then so is $\bar{h}$.

\ourexercises
\begin{ourexercise}\label{exe:lattsuff}
  Prove that, for a poset $(P,\leq)$ to be a lattice, it is sufficient that joins and meets exist for the empty set and for all two-element subsets. Also show that, if $f \colon L \to M$ is a function between lattices which preserves the empty join and joins of two-element sets, then $f$ preserves all finite joins.
\end{ourexercise}

\begin{ourexercise}\label{exe:complattsuff}
  Prove that, for a poset $(P,\leq)$ to be a complete lattice, it is sufficient that every subset has a supremum. By order-duality, it is also sufficient that every subset has an infimum.
\end{ourexercise}

\begin{ourexercise}
  Let $L$ be a lattice. Prove that $\bot$ is an \emphind{identity element} for $\vee$, that is, $a \vee \bot = a$ for all $a \in L$. By order-duality, $\top$ is an identity element for $\wedge$.
\end{ourexercise}

\begin{ourexercise}\phantomsection\label{exe:latticedefs}
  \begin{enumerate}
    \item Prove that the relation $\leq_L$ defined from a lattice according to the algebraic definition is a partial order in which all finite subsets have joins and meets. Further show that the binary and empty meets and joins in this partial order agrees with the operations of the original lattice.  \hint{Use the result of Exercise~\ref{exe:lattsuff}.}
    \item Prove that, if $(L,\leq)$ is a lattice according to the order-theoretic definition, then $(L,\vee,\wedge,\bot,\top)$ is a lattice according to the algebraic definition, where the operations denote the binary and empty joins and meets with respect to the partial order $\leq$. Further show that the partial order obtained from the binary operations gives back the original partial order. 
  \end{enumerate}
\end{ourexercise}

\begin{ourexercise}\label{exe:easy-counterexamples}
  Find at least one example of each of the following:
  \begin{enumerate}
    \item a partial order which is not a total order;
    \item a supremum which is not a maximum;
    \item a poset in which finite joins exist but which is not a lattice;
    \item a lattice which is not a complete lattice;
    \item a subset of a lattice which is not a sublattice, even though it is a lattice in the inherited partial order;
    \item an order-embedding between lattices that preserves binary joins, but not the bottom element.
  \end{enumerate}
\end{ourexercise}

\begin{ourexercise}\label{exe:counterexamples-complete}
  This exercise shows that there are some subtleties with the notion of completeness in lattices. Let $X$ be an infinite set. We will consider sublattices of the complete lattice $\cP(X)$, ordered by subset inclusion.
  \begin{enumerate}
    \item\label{ite:complete-sublattice-not-sub} 
    Consider the set $$L := \{ f \in \cP(X) \ \mid \ f\text{ is finite} \} \cup \{X\}.$$ Show that $L$ is a bounded sublattice of $\cP(X)$, and that $L$ is a complete lattice in the order inherited from $\cP(X)$, but also give an example of a subset $E$ of $L$ such that  the supremum of $E$ in $L$ is different from the supremum of $E$ in $\cP(X)$.

    \item \label{ite:compDLnotframe} Show that $L$ of the previous item is not a frame.
    \item A subset $c$ of $X$ is called \emphind{cofinite} if $X \setminus c$ is finite. Show that the set $$L' := \{c \in \cP(X)\mid c\text{ is cofinite}\} \cup \{\emptyset\}$$ is a bounded sublattice of $\cP(X)$, and a complete lattice in the order inherited from $\cP(X)$. 
    \item Let $K$ be a subposet of $\cP(X)$ containing both $L$ and $L'$. Show that if $K$  is a complete lattice in the order inherited from $\cP(X)$, then $K=\cP(X)$.
  \end{enumerate}
\end{ourexercise}

\begin{ourexercise}\label{exe:injlatthom}
  Let $f \colon L \to M$ be a function between lattices. Prove that
  \begin{enumerate}
    \item if $f$ is injective and preserves $\vee$ or $\wedge$, then $f$ is an order-embedding;
    \item if $f$ is bijective and preserves $\vee$ or $\wedge$, then $f$ is a lattice isomorphism.
  \end{enumerate}
\end{ourexercise}

\begin{ourexercise}\label{exe:disteq}
  Let $L$ be a lattice. Prove that (\ref{eq:dist1}) and (\ref{eq:dist2}) are equivalent. \hint{You need to use the absorption laws twice.}

  Note that this gives a proof that $L$ is distributive if, and only if, $L^\op$ is distributive; the property of distributivity is \emph{self-dual}.
\end{ourexercise}

\begin{ourexercise}\label{exe:product-lattices}
  Let $(L_i)_{i \in I}$ be an indexed family of lattices.
  \begin{enumerate}
    \item Prove that $\prod_{i \in I} L_i$, as defined in the text, is indeed a lattice.
    \item Prove that the partial order on $\prod_{i \in I}L_i$ is given by the pointwise product of the partial orders on the $L_i$.
    \item Prove that $\prod_{i\in I}L_i$ is distributive if, and only if, $L_i$ is distributive for every $i \in I$.
  \end{enumerate}
\end{ourexercise}

\begin{ourexercise}\label{exe:distchar}
  Prove Proposition~\ref{prop:distchar}.
\end{ourexercise}

\begin{ourexercise}\phantomsection\label{exe:distuncomp}
  \begin{enumerate}
    \item Prove that, in a distributive lattice, any element has at most one complement.
          \item\label{itm:lathomBA} Prove that any lattice homomorphism between Boolean algebras preserves the operation $\neg$.
    \item Prove that any function between Boolean algebras that preserves $\neg$, $\bot$ and $\vee$ is a lattice homomorphism.
  \end{enumerate}
\end{ourexercise}

\begin{ourexercise}\phantomsection\label{exe:directed}
  \begin{enumerate}
    \item   Let $P$ be a poset. Prove that, for any subset $D \subseteq P$, $D$ is directed if, and only if, any finite subset of $D$ has an upper bound in $D$. (Note that the empty set is always a subset of $D$, and that any element of $D$ will be an upper bound for it.)
    \item Let $L$ be a bounded lattice. Prove that $L$ is a complete lattice if, and only if, $L$ is a dcpo.
    \item Give an example of a dcpo that is not a lattice.
  \end{enumerate}
\end{ourexercise}

\begin{ourexercise}\label{exe:BAeq}
  Let $(B,+,\cdot,0,1)$ be a commutative ring with unit in which $a^2 = a$ for all $a \in B$. Define $a \leq b$ if, and only if, $a \cdot b = a$. Prove that $\leq$ is a distributive lattice order on $B$, and that
  every element of $B$ has a complement with respect to $\leq$. \hint{First show that $a + a = 0$ for all $a \in B$.}

  Conversely, let $(B,\wedge,\vee,\bot,\top,\neg)$ be a Boolean algebra. Define, for any $a, b \in B$, $a + b := (a \wedge \neg b) \vee (\neg a \wedge b)$, $a \cdot b := a \wedge b$, $0 := \bot$ and $1 := \top$. Prove that $(B,+,\cdot,0,1)$ is a commutative ring with unit in which $a^2 = a$ for all $a \in B$. 

  Finally, show that the composition of these two assignments in either order yields the identity.
  
  \emph{Note.} The operation $+$ defined in this exercise is known as \emphind{symmetric difference}.%
  \nl{$a + b$}{symmetric difference of two Boolean algebra elements $a$ and $b$}{plus}
\end{ourexercise}

\begin{ourexercise}\label{exe:adjointexistsiff}
  This exercise guides you through a proof of the `\emphind{adjoint functor theorem for complete lattices}'. Let $C, D$ be complete lattices and $f \colon C \to D$ a function.
  \begin{enumerate}
    \item Suppose that $f$ preserves all joins. For each $d \in D$, define $g(d) := \bigvee \{c \in C \ | \ f(c) \leq d\}$. Prove that $g$ is upper adjoint to $f$.
    \item Conclude from the previous item and Exercise~\ref{exe:adjunctions}.\ref{itm:leftadjointpreservessup} that a function $f$ between complete lattices possesses an upper adjoint if, and only if, $f$ preserves all joins.
    \item Conclude from the previous item, applied to $C^\op$ and $D^\op$, that a function $f$ between complete lattices possesses a lower adjoint if, and only if, $f$ preserves all meets.
  \end{enumerate}
\end{ourexercise}

\begin{ourexercise}\label{exe:adjointsfix}
  Let $f \colon C \leftrightarrows D \colon g$ be an adjunction between complete lattices, with $f$ left adjoint to $g$.
  \begin{enumerate}
    \item Show that the image of $g$ in $C$ is closed under all meets.
    \item Show that the image of $g$ is isomorphic to the image of $f$ in $D$.
    \item Give an example showing that the image of $g$, despite being a complete lattice in its own right, need not be a sublattice of $C$.
  \end{enumerate}
\end{ourexercise}

\begin{ourexercise}\phantomsection\label{exe:congruencelattice}
  \begin{enumerate}
    \item Prove that the collection of congruences on a lattice $L$ is a complete lattice under the inclusion order.
    \item Prove that this lattice of congruences is always distributive, even if $L$ is not.
  \end{enumerate}
\end{ourexercise}

\section{Duality for finite distributive lattices}\label{sec:finDLduality}
Lattices were introduced as abstract structures in the previous section. In this section we show that finite distributive lattices can be represented in a more concrete way, namely, as certain  sets of subsets closed under intersection and union equipped with the inclusion order. That is, the abstract algebraic operations of meet and join on finite distributive lattices are always, up to isomorphism, the set-theoretic operations of intersection and union.
This representation gives rise to our first example of a \emph{duality}. From it, we will also deduce a duality for finite Boolean algebras.

For any set $S$, we denote by $\mathcal{P}(S)$ the \emphind{power set} of $S$, that is, the collection of all subsets of $S$. The \emphind{inclusion order} of subsets gives a partial order on $\mathcal{P}(S)$, which is in fact a distributive lattice. (Indeed, $\mathcal{P}(S)$ is even a \emph{Boolean algebra}.) Any sublattice of $\mathcal{P}(S)$ is a distributive lattice, too. Conversely, any distributive lattice is a sublattice of a power set lattice, as we will see in Chapter 3. %
In this section, we will prove a stronger result for \emph{finite} distributive lattices (Proposition~\ref{prop:birkhoff}).

Let $(P,\preceq)$ be a preorder. An \emphind{up-set} is a subset $U \subseteq P$ such that whenever $p \in U$ and $p \preceq q$, we have $q \in U$. A \emphind{down-set} is a subset $D \subseteq P$ such that whenever $p \in D$ and $q \preceq p$, we have $q \in D$.
The collection of up-sets of a preorder is closed under arbitrary unions and intersections, and the preorder $\preceq$ thus specifies two complete sublattices of $\mathcal{P}(P)$, namely, the complete sublattice $\Up(P,\preceq)$ of \emph{up-sets} with respect to $\preceq$, and the complete sublattice $\Down(P, \preceq)$ of \emph{down-sets} with respect to $\preceq$, see Exercise~\ref{exe:updownsublat}.%
\nl{$\Down(P)$}{the complete lattice of down-sets of a poset $P$}{}
\nl{$\Up(P)$}{the complete lattice of up-sets of a poset $P$}{}

It follows that, for any subset $S$ of $P$,
there exists a smallest up-set containing $S$, which we denote by ${\uparrow} S$ and call the \emph{up-set generated by $S$}.
The set ${\uparrow} S$ contains those elements $p \in P$ for which there exists $s \in S$ such that $s \preceq p$, that is,
\[ {\uparrow} S = \{ p \in P \ \mid \ s \preceq p \text{ for some } s \in S \} \ . \]
\index{up-set!generated by a subset}\index{generated up-set}%
\index{down-set!generated by a subset}\index{generated down-set}\nl{${\downarrow} S$}{the down-set generated by a subset $S$ of a poset $P$}{}%
\nl{${\uparrow} S$}{the up-set generated by a subset $S$ of a poset $P$}{}
The \emph{down-set generated by $S$} is defined similarly, and denoted ${\downarrow} S$. %
In particular, for any $p \in P$, the up-set generated by $p$, which we denote by ${\uparrow} p$, %
is the set of elements above $p$, and the down-set generated by $p$, which we denote by ${\downarrow} p$, %
is the set of elements below $p$. 
The up-sets of the form ${\uparrow} p$ and the down-sets of the form ${\downarrow} p$ are called 
\emph{principal} up-sets and down-sets, respectively. 
By a \emphind{convex set} we mean a set that is an intersection of an up-set and a down-set (also see Exercise~\ref{ex:convex}). By a \emphind{finitely generated up-set} we mean an up-set of the form ${\uparrow} F$ for $F$ a finite (possibly empty) subset, or equivalently, a set that is a finite (possibly empty) union of principal up-sets. Similarly, a \emphind{finitely generated down-set} is a finite union of principal down-sets; these up-sets will be particularly important in Chapter~\ref{chap:DomThry}, also see Exercise~\ref{exe:fingendown}.
\nl{${\uparrow} p$}{the principal up-set generated by an element $p$ of a poset $P$}{}
\nl{${\downarrow} p$}{the principal down-set generated by an element $p$ of a poset $P$}{}
\index{principal up-set} \index{principal down-set} \index{up-set!principal} \index{down-set!principal}

 Notice that, if $U$ is an up-set, then its complement, $P \setminus U$, is a down-set, and vice versa; therefore, $\Up(P,\preceq)$ is isomorphic to $\Down(P,\preceq)^\op$, or, said otherwise, $\Up(P,\preceq)$ and $\Down(P,\preceq)$ are \emph{anti-isomorphic}.

\begin{notation}
  Throughout this book, when $U$ is a subset of a set $P$, we often use the notation $U^c$ instead of $P \setminus U$ for the \emphind{complement} of a subset $U$. Note that this abbreviated notation $U^c$ assumes that the `ambient' set $P$ is clear from the context.%
  \nl{$U^c$}{the complement of a subset $U$}{}
\end{notation}

Let $j$ be an element of a lattice $L$. Then $j$ is called (finitely)\footnote{We will omit the adjective `finitely'. On occasion we need  stronger versions of join-irreducibility, namely with respect to arbitrary subsets. In that case we will use the adjective \emph{completely}, see for example Exercise~\ref{exe:discrete-duality}. As a general rule throughout this book, when no further qualifying adjective is given for a concept involving joins and meets, we mean \emph{finite} joins and meets.} \emphind{join-irreducible} if, whenever $j = \bigvee S$ for a finite $S \subseteq L$, we have $j \in S$. Notice that $\bot$ is never a join-irreducible element, because $\bot = \bigvee \emptyset$. A useful equivalent definition of this concept is: an element $j$ of $L$ is join-irreducible if, and only if, $j \neq \bot$ and, whenever $j = a \vee b$ for $a, b \in L$, we have $j = a$ or $j = b$ (see Exercise~\ref{exe:joinirralt}).
We denote by $\cJ(L)$ the poset of join-irreducible elements of $L$, where the order is the restriction of the order on $L$.
Similarly, $m \in L$ is a (finitely) \emphind{meet-irreducible} element if $m = \bigwedge S$ implies $m \in S$ for any finite $S \subseteq L$, $\top$ is never meet-irreducible, and $\cM(L)$ denotes the poset of meet-irreducible elements of $L$. Note that meet-irreducibles of $L$ are exactly the same thing as join-irreducibles of $L^\op$.%
\nl{$\cJ(L)$}{the poset of (finitely) join-irreducible elements of a lattice $L$}{}
\nl{$\cM(L)$}{the poset of (finitely) meet-irreducible elements of a lattice $L$}{}

An important and useful fact about \emph{finite} lattices is that there are `enough' join-irreducibles to separate elements, in the following sense.\index{enough join-irreducibles}
\begin{lemma}\label{lem:finiteperfect}
  Let $L$ be a finite lattice. For any $a, b \in L$, if $a \nleq b$, then there exists $j \in \cJ(L)$ such that $j \leq a$ and $j \nleq b$.
\end{lemma}
\begin{proof}
  The set $T := ({\downarrow} a) \setminus ({\downarrow} b)$ of elements that are below $a$ but not below $b$ is non-empty, as it contains $a$. Since $L$ is finite, pick a minimal element $j$ of $T$. This element $j$ must be join irreducible. Indeed, suppose that $j = \bigvee S$ for some finite $S \subseteq L$. Then, since $\bigvee S \nleq b$, pick $c \in S$ such that $c \nleq b$. Since $c \leq j \leq a$, we have $c \in T$, so the minimality of $j$ implies that $j = c$.
\end{proof}
\nl{$\hat{a}$}{representation of an element $a$ of a finite distributive lattice as a down-set of its dual space}{}
For any finite lattice $L$, consider the function
\begin{align*}
  \widehat{(-)} \colon & L \to \Down(\cJ(L))                                         \\
                       & a \mapsto \widehat{a} := \{j \in \cJ(L) \ | \ j \leq_L a\},
\end{align*}
which sends every element of the lattice to the down-set of join-irreducibles below it. This function $\widehat{(-)}$ is obviously order preserving, and Lemma~\ref{lem:finiteperfect} says precisely that $\widehat{(-)}$ is an order-embedding. When is $\widehat{(-)}$ surjective, and hence an order isomorphism? We will give the answer in Proposition~\ref{prop:birkhoff}.

An element $j$ in a lattice $L$ is called \emphind{join-prime} if, for every finite $S \subseteq L$, $j \leq \bigvee S$ implies $j \leq a$ for some $a \in S$. Note that any join-prime element in a lattice is in particular join irreducible (see Exercise~\ref{exe:joinirrjoinp}). \emph{Meet-prime} elements are defined in the order-dual way: an element $m$ in a lattice $L$ is \emphind{meet-prime} if, for every finite $S \subseteq L$, $m \geq \bigwedge S$ implies $m \geq a$ for some $a \in S$. A meet-prime element of $L$ is again the same thing as a join-prime element of $L^\op$.
\begin{proposition}\label{prop:birkhoff}
  Let $L$ be a finite lattice. The following are equivalent:
  \begin{enumerate}[label=(\roman*)]
    \item the lattice $L$ is distributive;
    \item every join-irreducible element of $L$ is join-prime;
    \item the function $\widehat{(-)}$ is surjective, and thus an isomorphism.
  \end{enumerate}
\end{proposition}
\begin{proof}
  (i) $\Rightarrow$ (ii). Let $j$ be join irreducible. If $j \leq \bigvee S$ for some finite $S$, then
  \[ j = j \wedge \big( \bigvee S \big) = \bigvee_{a \in S} (j \wedge a),\]
  where we use the distributive law in the last step. Since $j$ is join irreducible, $j = j \wedge a$ for some $a \in S$, which means that $j \leq a$.

  (ii) $\Rightarrow$ (iii).  If $D \in \Down(\cJ(L))$ and $j \in \cJ(L)$, then $j \leq \bigvee D$ if, and only if, there exists $a \in D$ such that $j \leq a$, which in turn is equivalent to $j \in D$ because $D$ is a down-set. Thus, $\widehat{\bigvee D} = D$, showing that $\widehat{(-)}$ is surjective.

  Finally, (iii) $\Rightarrow$ (i) is clear, because the lattice $\Down(\cJ(L))$ is distributive, and distributivity is preserved by isomorphism.
\end{proof}
Note that the proof of implication (i) $\Rightarrow$ (ii) in Proposition~\ref{prop:birkhoff} did not use the assumption that $L$ is finite. Therefore, this proposition also implies that in any distributive lattice $L$, \emph{join-prime} and \emph{join-irreducible} are synonymous. The order dual of Proposition~\ref{prop:birkhoff} can be used to characterize the property of distributivity (which is self-dual) in terms of \emph{meet-prime} and \emph{meet-irreducible} elements. Moreover, for a finite distributive lattice $L$, the posets $\cJ(L)$ and $\cM(L)$ are isomorphic (see Exercise~\ref{exe:kappa}).

If $L$ is a finite distributive lattice, then we call $\cJ(L)$ the \emphind{dual poset} of $L$. If $P$ is a finite poset, we call $\Down(P)$ the \emphind{dual distributive lattice} of $P$. With this terminology, Proposition~\ref{prop:birkhoff} implies that any finite distributive lattice is isomorphic to its double dual, that is,
\[ L \cong \Down(\cJ(L)) \ .\]
We get a similar `double dual' result if we start from a finite poset $P$.

\begin{proposition}\label{prop:finite-poset-double-dual}
  Let $P$ be a finite poset and $D \in \Down(P)$. Then $D$ is join irreducible in $\Down(P)$ if, and only if, $D$ is a principal down-set. In particular, $\cJ(\Down(P))$ is a poset isomorphic to $P$.
\end{proposition}

\begin{proof}
  Suppose $D$ is join irreducible. Since $P$ is finite, we have 
  \[ D = \bigcup_{p \in D} {\downarrow} p. \]
Since $D$ is join irreducible, we can pick $p \in D$ such that $D = {\downarrow} p$, as required.

  Conversely, if $D = {\downarrow} p$ for some $p \in P$, then $D$ is non-empty. Also, if $D = A_1 \cup A_2$ for some down-sets $A_1, A_2$ of $P$, then we have $p \in A_i$ for some $i \in \{1,2\}$. Since $A_i$ is a down-set, we then obtain $D = {\downarrow} p \subseteq A_i \subseteq D$, so $D = A_i$.
\end{proof}

We thus see that any finite distributive lattice can be represented isomorphically as the lattice of down-sets of a finite poset (Proposition~\ref{prop:birkhoff}). Moreover, for a given finite distributive lattice $L$, the finite poset $P$ for which $L \cong \Down(P)$ is unique up to isomorphism: indeed, Proposition~\ref{prop:finite-poset-double-dual} implies that if $L \cong \Down(P)$ for some finite poset $P$, then $\cJ(L) \cong \cJ(\Down(P)) \cong P$. 

To turn this representation result into a \emph{duality}, we now consider morphisms. If $f \colon P \to Q$ is an order-preserving function between preorders, then the \emphind{inverse image} map
\begin{align*}
  \Down(f) \colon & \Down(Q) \to \Down(P) \\
                  & D \mapsto f^{-1}(D)
\end{align*}
\nl{$f^{-1}(S)$}{inverse image of a subset $S$ of the codomain under a function $f$}{}
is a lattice homomorphism.  We will prove in Proposition~\ref{prop:finDLmorphisms} below that, in the special case where $P$ and $Q$ are finite posets, every lattice homomorphism $\Down(Q) \to \Down(P)$ arises in this way. Before we can do so, we need a general lemma about adjunctions between finite lattices.
\begin{lemma}\label{lem:loweradj-preserves-joinprime}
  Let $g \colon D \leftrightarrows E \colon h$ be an adjunction between finite lattices, and further suppose that $h$ preserves joins. Then, for any join-prime element $j$ in $D$, the element $g(j)$ is join prime in $E$.
\end{lemma}
\begin{proof}
  Let $S \subseteq E$ be any subset such that $g(j) \leq \bigvee S$. Then, since $h$ is upper adjoint to $g$, $j \leq h(\bigvee S) = \bigvee h[S]$, where we use that $h$ preserves joins. Since $j$ is join prime, pick $s \in S$ such that $j \leq h(s)$. Since $g$ is lower adjoint to $h$, $g(j) \leq s$. Thus, $g(j)$ is join prime.
\end{proof}
\begin{proposition}\label{prop:finDLmorphisms}
  Let $P$ and $Q$ be finite posets. For any lattice homomorphism $h \colon \Down(Q) \to \Down(P)$, there exists a unique order-preserving $f \colon P \to Q$ such that $h = \Down(f)$.
\end{proposition}
\begin{proof}
  Since $h$ preserves all meets, it has a lower adjoint, $g$, by Exercise~\ref{exe:adjointexistsiff}(c). By Lemma~\ref{lem:loweradj-preserves-joinprime}, since $h$ also preserves all joins, $g$ sends join-prime elements to join-prime elements. Now, if $p \in P$, then ${\downarrow}p$ is join prime by Proposition~\ref{prop:finite-poset-double-dual}, and thus $g({\downarrow}p)$ is join prime. By the other direction of Proposition~\ref{prop:finite-poset-double-dual}, pick the unique $f(p) \in Q$ such that $g({\downarrow}p) = {\downarrow}f(p)$.
  Notice that the function $f \colon P \to Q$ thus defined is order preserving, because $g$ is order preserving, so $p \leq p'$ implies $f(p) \in {\downarrow} f(p) = g({\downarrow} p) \subseteq g({\downarrow} p') = {\downarrow} f(p')$.

  Moreover, for any $E \in \Down(Q)$ and $p \in P$, we have, using the adjunction and the definition of $f$, that
  \[p \in h(E) \iff {\downarrow} p \subseteq h(E) \iff g({\downarrow}p) \subseteq E \iff {\downarrow} f(p) \subseteq E \iff f(p) \in E,\]
  so that $h(E) = f^{-1}(E)$, as required.
  The uniqueness of $f$ is left as Exercise~\ref{exe:completeuniquenesspart}.
\end{proof}
Summing up, we have associated to every finite distributive lattice $L$ a finite poset $\cJ(L)$, that we called the dual poset of $L$, and, conversely, to every finite poset $P$, a finite distributive lattice $\Down(P)$ that we called the dual distributive lattice of $P$. We have proved that:
\begin{enumerate}
  \item[(1)] every finite distributive lattice is isomorphic to its double dual (Proposition~\ref{prop:birkhoff})
  \item[(2)] every finite poset is isomorphic to its double dual (Proposition~\ref{prop:finite-poset-double-dual}),
  \item[(3)] homomorphisms between finite distributive lattices are in one-to-one correspondence with order-preserving functions between their dual posets (Proposition~\ref{prop:finDLmorphisms}).
\end{enumerate}
The reversal of direction of arrows when moving to `the other side', that is, a function from $P$ to $Q$ gives a function from $\Down(Q)$ to $\Down(P)$, is what makes the correspondence in (3) `\emph{dual}'. What all this means, in practice, is that finite distributive lattices with homomorphisms between them are \emph{essentially the same thing} as finite posets with order-preserving functions between them.

In fancier terms, we have proved the following theorem.  \index{duality!for finite distributive lattices} \index{Birkhoff duality} \index{duality!Birkhoff}
\begin{theorem}\label{thm:birkhoffduality}
  The functors $\Down$ and $\cJ$ constitute a duality between the category $\DL_f$ of finite distributive lattices with homomorphisms and the category $\Pos_f$ of finite posets with order-preserving functions.
\end{theorem}
We will define the precise meaning of the terms (`category', `functor', `duality') used in this theorem in Section~\ref{sec:external} in Chapter~\ref{ch:categories}, but the reader who is not yet familiar with these terms can rest assured that the mathematical content of the theorem consists precisely of items (1), (2), and (3) above. For a precise explanation of why what we have proved here shows that the functors form a duality, see Example~\ref{exa:cat-duality-examples} on p.~\pageref{exa:cat-duality-examples}.

An easy generalization of Birkhoff duality, which does not involve topology, is what we call \emphind{discrete duality} for distributive lattices. This generalization starts from the observation that for any (not necessarily finite) poset $P$, one may still recover $P$ from the distributive lattice $\Down(P)$, now as the poset of \emph{completely} join-prime elements. This yields a dual equivalence between the category $\Pos$ of posets with order-preserving functions and a category $\DL^+$ of complete lattices that are join generated by their completely join-prime elements; see Exercise~\ref{exe:discrete-duality} below, and also Example~\ref{exa:adjunctions}.\ref{itm:DLplus} and Theorem~\ref{thrm:Raney53} later in this book.

\subsection*{Duality for finite Boolean algebras}
We end this section by describing how Theorem~\ref{thm:birkhoffduality} specializes to finite Boolean algebras.
An \emphind{atom} of a lattice $L$ is a minimal non-bottom element, that is, an element $j \in L$ such that $j \neq \bot$ and $\bot \leq a \leq j$ implies $a = \bot$ or $a = j$ for any $a \in L$.
\begin{proposition}\label{prop:birkhoffBA}
  Let $L$ be a finite distributive lattice. The following are equivalent:
  \begin{enumerate}
    \item[(i)] the distributive lattice $L$ is a Boolean algebra;
    \item[(ii)] every join-irreducible element of $L$ is an atom;
    \item[(iii)] the order on $\cJ(L)$ is trivial, that is, distinct elements are incomparable;
    \item[(iv)] the distributive lattice $L$ is isomorphic to $\mathcal{P}(\cJ(L))$.
  \end{enumerate}
\end{proposition}
\begin{proof}
  (i) $\Rightarrow$ (ii). Let $j \in \cJ(L)$. Suppose that $\bot \leq a \leq j$. Since $a \vee \neg a = \top$, we have $j \leq a \vee \neg a$. Since $L$ is distributive, by Proposition~\ref{prop:birkhoff} $j$ is join prime, so $j \leq a$ or $j \leq \neg a$. If $j \leq a$, then $a = j$, and we are done. If $j \leq \neg a$, then $a \leq j \leq \neg a$, so $a = a \wedge \neg a = \bot$.

  (ii) $\Rightarrow$ (iii). Clear from the definition of atom.

  (iii) $\Rightarrow$ (iv). By Proposition~\ref{prop:birkhoff}, $L$ is isomorphic to $\Down(\cJ(L))$. By (iii), any subset of $\cJ(L)$ is a down-set.

  (iv) $\Rightarrow$ (i). Clear because $\mathcal{P}(\cJ(L))$ is a Boolean algebra.
\end{proof}
Proposition~\ref{prop:birkhoffBA} shows in particular that every finite Boolean algebra $L$ is of the form $\mathcal{P}(S)$, where $S$ is the finite set of atoms (= join-irreducibles) of $L$.

Note also from (iii) in Proposition~\ref{prop:birkhoffBA} that, if $L$ is a finite Boolean algebra and $M$ is a finite distributive lattice, then \emph{any} function $f \colon \cJ(L) \to \cJ(M)$ is order preserving. Thus, the lattice homomorphisms $M \to L$ are exactly the inverse images of functions $\cJ(L) \to \cJ(M)$. In particular, if both $L$ and $M$ are finite Boolean algebras, then homomorphisms from $L$ to $M$ correspond to functions $\cJ(M) \to \cJ(L)$. We conclude:\index{duality!for finite Boolean algebras} \index{Birkhoff duality} \index{Stone duality!finite case}
\begin{theorem}\label{thm:birkhoffdualityBA}
  The functors $\mathcal{P}$ and $\cJ$ constitute a duality between the category $\BA_f$ of finite Boolean algebras with homomorphisms and the category $\Set_f$ of finite sets with functions.
\end{theorem}

As promised at the end of the previous section, we end this first chapter by using the dualities to give a simple concrete description of the \emphind{Boolean envelope} of a finite distributive lattice, defined in the previous section.
\begin{proposition}\label{prop:finite-boolenv}
  Let $L$ be a finite distributive lattice. Then the finite Boolean algebra $\mathcal{P}(\mathcal{J}(L))$, with the embedding $\widehat{(-)} \colon L \to \cP(\cJ(L))$, is a Boolean envelope of $L$.
\end{proposition}
\begin{proof}
  Let $h \colon L \to A$ be a lattice homomorphism, with $A$ a Boolean algebra. By the results in this section, we may assume, up to isomorphism, that $L = \mathcal{D}(P)$ for a poset $P$, $A = \mathcal{P}(X)$ for a set $X$, and $h = f^{-1}$ for a function $f \colon X \to P$. The function $\bar{h}$, defined by sending any $u \in \mathcal{P}(P)$ to $f^{-1}(u)$, is clearly a homomorphism extending $h$. The uniqueness is left as Exercise~\ref{ex:unique-finite-boolean-envelope}; or see the more general proof of Proposition~\ref{prop:boolenv} in Chapter~\ref{ch:priestley}.
\end{proof}

This intimate connection between distributive lattices and Boolean algebras, and the various dualities for them, will be made more precise in Section~\ref{sec:boolenv-duality} in Chapter~\ref{ch:priestley}, and further in Section~\ref{sec:StoneSpaces} in Chapter~\ref{chap:Omega-Pt}.

As in the case of distributive lattices, the duality between finite Boolean algebras and finite sets fairly easily extends to a duality between sets and Boolean algebras that are complete and \emphind{atomic}, that is, for any $a \neq \bot$, there is an atom $j \leq a$. This is called \emphind{discrete duality} between sets and complete and atomic Boolean algebras (see Exercise~\ref{ex:CABASet}).

\ourexercises
\begin{ourexercise}\label{exe:updownsublat}
  Let $P$ be a poset. Prove that $\Down(P)$ is a complete sublattice of $\cP(P)$. Use order-duality to deduce that $\Up(P)$ is also a complete sublattice of $\cP(P)$.
\end{ourexercise}

\begin{ourexercise}\label{ex:convex}
  Prove that a subset $C$ of $P$ is convex if, and only if, for any $p, q \in C$, if $p \leq r \leq q$, then $r \in C$.
\end{ourexercise}

\begin{ourexercise}\label{exe:joinirralt}
  Prove that an element $j$ in a lattice $L$ is join irreducible if, and only if, $j \neq \bot$ and for any $x,y \in L$, if $j = x \vee y$, then $j = x$ or $j = y$.
\end{ourexercise}

\begin{ourexercise}\phantomsection\label{exe:perfect}
  \begin{enumerate}
    \item Formulate a lemma which says that any finite lattice has enough meet-irreducibles, and prove it (you may use Lemma~\ref{lem:finiteperfect}).
    \item Give an example of a lattice which does not have enough join-irreducibles.
    \item Give an example of a complete lattice which has enough join-irreducibles, but not enough meet-irreducibles.
    \item Can you find an example as in (c) with the additional requirement that the complete lattice is a \emph{frame}? That it is a \emph{Boolean algebra}? If not, formulate and prove a result about this.
  \end{enumerate}
\end{ourexercise}

\begin{ourexercise}\phantomsection\label{exe:joinirrjoinp}
  \begin{enumerate}
    \item Prove that any join-prime element in a lattice is join irreducible.
    \item Give an example of a lattice and a join-irreducible element in it that is not join prime.
  \end{enumerate}
\end{ourexercise}

\begin{ourexercise}\label{exe:completeuniquenesspart}
  Complete the uniqueness part of the proof of Proposition~\ref{prop:finDLmorphisms} by showing that, if $f, f' \colon P \to Q$ are order preserving and $\Down(f) = \Down(f')$, then $f = f'$.%
\end{ourexercise}

\begin{ourexercise}\label{exe:kappa}
  Let $L$ be a finite distributive lattice.
  \begin{enumerate}
    \item Prove that, for any $j \in \cJ(L)$, the set
          \[ L \setminus ({\uparrow} j) = \{a \in L \ | \ j \nleq a\}\]
          has a maximum. We denote this maximum by $\kappa(j)$.
    \item Prove that, for any $j \in \cJ(L)$ and $a \in L$, $j \nleq a$ if, and only if, $a \leq \kappa(j)$.
    \item Prove that $\kappa(j)$ is meet irreducible for any $j \in \cJ(L)$.
    \item Prove that $\kappa \colon \cJ(L) \to \cM(L)$ is an order-isomorphism.
  \end{enumerate}
\end{ourexercise}

\begin{ourexercise}\label{exe:fingendown}
  This exercise describes a variant of Proposition~\ref{prop:finite-poset-double-dual} that works for posets that are not necessarily finite, and gives a universal property for the construction. 
  Also compare this exercise with Exercise~\ref{exe:discrete-duality} below.
  Let $P$ be any poset. 
  \begin{enumerate}
    \item Prove that a down-set $D$ of $P$ is finitely generated if, and only if, $\max(D)$ is finite and $D = {\downarrow}\max(D)$.
  \end{enumerate}
  Denote by $\Downfin(P)$ the poset of finitely generated down-sets of $P$, ordered by inclusion.
  \begin{enumerate}[resume]
    \item Show that $\Downfin(P)$ is a join-subsemilattice of $\Down(P)$.
    \item Prove that the function $P \to \cJ(\Downfin(P))$, defined by sending $p \in P$ to ${\downarrow} p$, is an order isomorphism between $P$ and $\cJ(\Downfin(P))$.
    \item Show that, for any join-semilattice $(L, \vee, \bot)$ and any order-preserving function $f \colon P \to L$, there exists a unique join-preserving function $\bar{f} \colon \Downfin(P) \to L$ such that $\bar{f}({\downarrow} p) = f(p)$ for every $p \in P$.
  \end{enumerate}
  The last item shows that $\Downfin(P)$ is the \emphind{free join-semilattice} over the poset $P$; we will revisit this construction in Example~\ref{exa:adjunctions}.\ref{itm:free-join-semi} in a more general context.
\end{ourexercise}

\begin{ourexercise}\label{ex:unique-finite-boolean-envelope}
  Prove that the function $\bar{h}$ defined in the proof of Proposition~\ref{prop:finite-boolenv} is unique. \hint{It suffices (why?) to prove that the sublattices $\mathcal{D}(P)$ and $\mathcal{U}(P)$ together generate $\mathcal{P}(P)$. For this, note first that any singleton $\{p\}$ can be obtained as the intersection of the down-set ${\downarrow} p$ and the up-set ${\uparrow} p$; then use that any $u \in \mathcal{P}(P)$ can be written as a (finite) union of singletons.}
\end{ourexercise}

\begin{ourexercise}\label{exe:discrete-duality}
  The final two exercises of this first chapter outline a \emphind{discrete duality} for distributive lattices, and show how it specializes to Boolean algebras. Your solution can essentially follow the same scheme as the proofs for \emph{finite} duality that we gave in this chapter, so it makes for a useful exercise to check that you have understood those.

  We call an element $j$ of a complete lattice $L$ \emphind{completely join-prime} if, for any subset $S$ of $L$, $j \leq \bigvee S$ implies $j \leq s$ for some $s \in S$.  We denote by $\cJ^\infty(L)$ the set of completely join-prime elements of $L$.  We say $L$ is \emph{generated by its completely join-prime elements}\index{down-set lattice}\index{completely join-prime!elements generate} if for every $a \in L$ there exists $ S\subseteq\cJ^\infty(L)$ such that $a = \bigvee S$. Note that in this case we in fact must have $a = \bigvee ({\downarrow}a \cap \cJ^\infty(L))$ for every $a \in L$.

  Throughout the exercise, $P$ and $Q$ denote arbitrary posets.
  \begin{enumerate}
    \item Recall from Exercise~\ref{exe:updownsublat} that, for any poset $P$, $\Down(P)$ is a complete sublattice of $\cP(P)$. Show that, for any $p \in P$, ${\downarrow} p$ is completely join prime. Conclude that $\Down(P)$ is generated by its completely join-prime elements.
    \item Prove that a complete lattice is generated by its completely join-prime elements if, and only if, for any $a, b \in L$, if $a \nleq b$, then there exists $j \in \cJ^\infty(L)$ such that $j \leq a$ and $j \nleq b$.
    \item Let $L$ be a complete lattice generated by its completely join-prime elements. Prove that the function $\widehat{(-)}$, defined by
          \[\widehat{a} := \{j \in \cJ^\infty(L) \ \mid \ j \leq_L a\},\]
          is an order isomorphism between $L$ and $\Down(\cJ^\infty(L))$. Conclude in particular that $L$ is completely distributive.
    \item Prove that, for any order-preserving function $f \colon P \to Q$, the inverse image map $f^{-1} \colon \Down(Q) \to \Down(P)$ is a complete lattice homomorphism.
    \item Let $h \colon M \to L$ be a complete lattice homomorphism between complete  lattices, and let $g$ be its lower adjoint, which exists by Exercise~\ref{exe:adjointexistsiff}. Prove that, if $p$ is completely join prime in $L$, then $g(p)$ is completely join prime in $M$.
    \item Prove that for any complete lattice homomorphism $h \colon \Down(Q) \to \Down(P)$, there exists a unique order-preserving $f \colon P \to Q$ such that $h = f^{-1}$.
    \item Show that, if $L$ is a complete lattice generated by its completely join-primes, and $M$ is a complete sublattice of $L$, then $M$ is also generated by its completely join-primes.
    \item \label{ite:realinterval-nocjp} Show that the real unit interval $[0,1]$ has no completely join-prime elements, but that it is the image of a complete lattice homomorphism with domain $\Down([0,1])$.
  \end{enumerate}
  \emph{Note.} This exercise shows that a lattice $L$ is isomorphic to one of the form $\Down(P)$, for $P$ a poset, if, and only if, $L$ is complete and generated by its completely join-prime elements. An alternative characterization is that $L$ is 
  complete, completely distributive, and generated by its completely join-irreducible elements. See Section~\ref{sec:dom} and Exercise~\ref{exer:compldist-jirr}.
\end{ourexercise}

\begin{ourexercise}\phantomsection\label{ex:CABASet}
  \begin{enumerate}
    \item Prove that, in a complete Boolean algebra, an element is completely join prime if, and only if, it is an atom.
    \item Conclude, using Exercise~\ref{exe:discrete-duality} that every complete and atomic Boolean algebra is isomorphic to one of the form $\cP(S)$.
    \item In contrast with Exercise~\ref{exe:discrete-duality}.\ref{ite:realinterval-nocjp}, show that if $h \colon M \onto L$ is a surjective complete lattice homomorphism between complete Boolean algebras and $M$ is atomic, then so is $L$. \hint{Show that if $a$ is an atom of $M$, then either $h(a)=0$ or $h(a)$ is an atom of $L$.}
  \end{enumerate}
\end{ourexercise}

\notessec
We recommend \cite{DavPri2002} as supplementary reading for a more detailed introduction to order and lattice theory. In addition, the classic \cite{BalDwi1974} remains a good resource for the basics of lattice theory, although it is somewhat outdated when it comes to the more advanced theory.

The main results in Section~\ref{sec:finDLduality} are essentially due to \cite{Bir1933} and are also consequences of the more general results of \cite{Sto1937/38}, which we will present later in this book.
There exist similar, but more involved results for (finite) lattices that are not necessarily distributive. In this text, however, we will limit ourselves to distributive lattices. More information and references on duality for general lattices can be found, for example, in \cite{DavPri2002} and in the introduction of our paper \cite{GehGoo2014}; also see \cite{MosJip14,MosJip14b}.
\chapter{Topology and order}\label{chap:TopOrd}

In this chapter we present some material at the interface of topology and order theory. We begin by recalling basic material on topology and order, leading up to the culminating result of the chapter, Theorem~\ref{thrm:COSpace-StabCompSp}, which provides an equivalence between certain compact Hausdorff topological spaces equipped with orders, which were first introduced by Nachbin, and certain non-Hausdorff spaces known as stably compact spaces. Nachbin's thesis of the same title as this chapter was first published in Portuguese in 1950 and later translated to English \parencite{Nachbin64}. It is a nice text and we recommend it as supplemental reading.

\section{Topological spaces}

We expect readers to be familiar with the basic definitions and notions of topology. Nevertheless, we give them here in order to fix notation and nomenclature. The exact placement of definitions can be located using the index. 

A \emphind{topological space} is a pair $(X,\tau)$ where $X$ is a non-empty set and $\tau$ is a bounded sublattice of $\cP(X)$ which is closed under arbitrary unions. The elements of $\tau$ are called \emphind{open sets} while their complements are said to be \emphind{closed sets}. We will often simply write $X$ for a topological space (if the collection of opens is clear). The collection of opens will then be denoted by $\Omega(X)$. The collection of closed subsets of $X$ is denoted $\cC(X)$. A subset $K \subseteq X$ is called \emphind{clopen} if it is both closed and open, that is, if both $K$ and $K^c$ are in $\tau$. The collection of clopen subsets of a space is denoted $\Clp(X)$.
\nl{$\Omega(X)$}{the collection of open subsets of a topological space $X$}{}
\nl{$\cC(X)$}{the collection of closed subsets of a topological space $X$}{}
\nl{$\Clp(X)$}{the collection of clopen subsets of a topological space $X$}{}

A function $f\colon X\to Y$ between topological spaces is \emphind{continuous} provided the inverse image function $f^{-1}\colon\cP(Y)\to\cP(X)$ takes opens of $Y$ to opens of $X$. That is, $f$ is continuous if, and only if, there is a restriction of $f^{-1}$ to $\Omega(Y)$ which makes the following diagram commute:
\[
\begin{tikzpicture}
\node(PY) at (0,1.5) {$\cP(Y)$};
\node(PX) at (3,1.5) {$\cP(X)$};
\node(OY) at (0,0) {$\Omega(Y)$};
\node(OX) at (3,0) {$\Omega(X)$};

\draw[->,>=stealth'] (PY)to node[above] {$f^{-1}$} (PX);
\draw[->,dashed,>=stealth'] (OY) to node[below] {} (OX);
\draw[right hook->,>=stealth'] (OY) to (PY);
\draw[right hook->,>=stealth'] (OX) to (PX);

\end{tikzpicture}
\]
\begin{example}
For any set $X$, the set $\tau=\cP(X)$ is a topology on $X$.  This topology is known as the \emphind{discrete topology} on $X$. Notice that any function from a discrete space to any topological space is continuous.
\end{example}

\begin{example}
For any set $X$, the set $\tau=\{\emptyset,X\}$  is a topology on $X$.  This topology is known as the \emphind{indiscrete topology} on $X$. Notice that any function from any topological space  to an indiscrete space are continuous.
\end{example}

\begin{example}\label{exp:real}
Let $\bR$ be the set of real numbers. The usual topology on $\bR$ consists of those sets $U\subseteq\bR$ with the property that, for each $x\in U$, there exists $\varepsilon>0$ such that the interval
\[
(x-\varepsilon,x+\varepsilon)=\{y\in\bR\mid x-\varepsilon<y<x+\varepsilon\}
\]
is entirely contained in $U$. It is not hard to see that a function $f\colon\bR\to\bR$ is continuous with respect to this topology if, and only if, it satisfies the usual epsilon-delta definition of continuity (see Exercise~\ref{exer:real}).
\end{example}

A \emphind{subspace} of a topological space $(X,\tau)$ is given by a subset
$Y\subseteq X$ and is equipped with the \emphind{subspace topology}, defined as 
\[
\tau\upharpoonright Y=\{U\cap Y\mid U\in\tau\}.
\]

A continuous function $f\colon X\to Y$ between topological spaces is said to be an \emphind{open mapping} provided $f[U]=\{f(x)\mid x\in U\}$ is open in $Y$ for any open $U\subseteq X$. Similarly,  $f\colon X\to Y$ is said to be a \emphind{closed mapping} provided $f[C]=\{f(x)\mid x\in C\}$ is closed in $Y$ for any closed $C\subseteq X$. Further, $f$ is said to be an \emphind{embedding} provided it is injective and $f^{-1}\colon\im(f)\to X$ is also continuous, where we view $f[X]=\im(f)$ as a topological space in the subspace topology. Finally, $f$ is a \emphind{homeomorphism} provided it is a bijection and both $f$ and $f^{-1}$ are continuous.

Since open sets are closed under unions, for any subset $S$ of a topological space $(X, \tau)$, there is a largest open set, $\intr(S)$, that is contained in $S$. The set $\intr(S)$ is called the \emphind{interior} of $S$; in a formula, 
\nl{$\intr(S)$}{interior of a subset $S$ of a topological space}{}
\[\intr(S) = \bigcup\{U \subseteq X \ | \ U \subseteq S \text{ and } U \in \tau\}.\] 
We note that $\intr$ is the upper adjoint to the inclusion map $\Omega(X) \into \mathcal{P}(X)$ (see Exercise~\ref{exer:int-adjoint}). If $x \in X$ is a point and $S$ is any subset of $X$, then $S$ is a \emphind{neighborhood} of $x$ if $x \in \intr(S)$, that is, if there exists an open set $U \subseteq S$ such that $x \in U$.  The interior of a set $S$ is sometimes denoted by $S^\circ$. 
\nl{$S^\circ$}{interior of a subset $S$ of a topological space}{}
Symmetrically, since closed sets are closed under arbitrary intersections, any subset $S$ has a \emphind{closure}, $\cl(S)$, which is defined as the smallest closed set containing $S$; in a formula, 
\[\cl(S) = \bigcap \{C\subseteq X \ | \ S \subseteq C\text{ and } C\text{ is } \tau\text{-closed}\}.\] 
\nl{$\cl(S)$}{closure of a subset $S$ of a topological space}{}
\nl{$\overline{S}$}{closure of a subset $S$ of a topological space}{}
The map $\cl$ is the lower adjoint to the inclusion map $\cC(X) \into \mathcal{P}(X)$. The closure of a set $S$ is sometimes denoted by $\overline{S}$. A subset $S$ of $X$ is called \emphind{dense} if $\cl(S) = X$. 

Notice that arbitrary intersections of topologies on a fixed set $X$ are again topologies. Accordingly, for any collection $\cS$ of subsets of $X$, there is a least topology containing $\cS$. We call this the \emphind{topology generated by $\cS$}. 
When $\tau$ is a topology on a set $X$, a collection $\cS$ of subsets of $X$ is called a \emph{subbase} for the topology $\tau$ if $\tau = \gen{\cS}$. Note that $\gen{\cS}$ is the collection of subsets of $X$ that can be written as arbitrary unions of finite intersections of elements of $\cS$.
We say that a subbase $\cS$ for a topology $\tau$ is a \emphind{base} if, for every $x\in X$ and $U\in\tau$ with $x\in U$, there is $V\in\cS$ with $x\in V\subseteq U$. Notice that, in this case, finite intersections are not needed to generate $\tau$, that is, every open of $\tau$ is simply a union of elements from $\cS$. An equivalent way of saying this is that a subbase $\cB$ is a base for the topology generated by $\cB$ if and only if, for every $x\in X$ 
\[
\cB_x=\{V\in\cB\mid x\in V\}
\]
is filtering. In particular, the closure of a subbase $\cS$ under finite intersections is always a base, but bases are more general than that; see Exercise~\ref{exer:base}.

Let $(X_i)_{i \in I}$ be a collection of topological spaces indexed by a set $I$. The \emphind{product space} $\prod_{i\in I} X_i$ is the Cartesian product of the sets $X_i$, equipped with the topology generated by the subbase consisting of the sets
\[
\pi_i^{-1}(V), \text{ where }i\in I \text{ and } V\subseteq X_i \text{ is open in  }X_i,
\]
where $\pi_i\colon \prod_{j\in I} X_j\to X_i$ is the projection onto the $i$th coordinate for $i\in I$.
Let $X$ be a topological space. We recall the five main \emphind{separation axioms} that may hold for $X$.
\begin{itemize}
\item $X$ is $T_0$ (or Kolmogorov) provided, for all $x,y\in X$ with $x\neq y$ there is an open $U\subseteq X$ which contains exactly one of $x$ and $y$;
\item $X$ is $T_1$ (or Fr\'echet) provided, for all $x,y\in X$ with $x\neq y$ there is an open $U\subseteq X$ with
$x\in U$ and $y\not\in U$;
\item $X$ is $T_2$ (or Hausdorff) provided, for all $x,y\in X$ with $x\neq y$ there are opens $U,V\subseteq X$ with
$x\in U$ and $y\in V$ and $U\cap V=\emptyset$;
\item $X$ is $T_3$ (or regular) provided $X$ is $T_1$ and, for all $x\in X$ and closed $C\subseteq X$ with $x\not\in C$ there are opens $U,V\subseteq X$ with $x\in U$ and $C\subseteq V$ and $U\cap V=\emptyset$;
\item $X$ is $T_4$ (or normal) provided $X$ is $T_1$ and, for all closed $C,D\subseteq X$ such that $C \cap D = \emptyset$, there are opens $U,V\subseteq X$ with $C\subseteq U$ and $D\subseteq V$ and $U\cap V=\emptyset$.
\end{itemize}
Hausdorff spaces may also be characterized as those spaces satisfying a \emph{closed graph theorem}, that is, a space $X$ is Hausdorff if, and only if, the graph of any continuous function $Y \to X$ is closed (see Exercise~\ref{exer:Hausdorff}).

\subsection*{Compactness}

Let $S\subseteq X$ where $X$ is a topological space. An \emphind{open cover} $\cU$ of $S$ is a collection of open sets $\cU\subseteq\Omega(X)$ such that $S\subseteq \bigcup\cU$. A subset $K\subseteq X$ is \emphind{compact} provided every open cover $\cU$ of $K$ contains a finite subcover, that is, a finite subcover $\cU'\subseteq\cU$ which is also a cover of $K$. An equivalent definition of compactness using closed sets instead of open sets is the following. For a collection $\cA \subseteq \cC(X)$ of closed sets of $X$, say that $\cA$ has the \emphind{finite intersection property} with respect to $K$ if for every finite subcollection $\cA'$ of $\cD$, $\big(\bigcap \cA'\big) \cap K \neq \emptyset$. Then $K$ is compact if, and only if, for every collection $\cA \subseteq \cC(X)$ that has  the finite intersection property with respect to $K$, we have $\big(\bigcap \cA\big) \cap K \neq \emptyset$. In particular, the space $X$ itself is compact if every open cover of $X$ contains a finite subcover, or equivalently, if every collection of closed sets with the finite intersection property (with respect to $X$) has a non-empty intersection (see Exercise~\ref{exer:covers} for further equivalent definitions of compactness for a space).

\begin{example}\label{exp:compact}
The subspace $[0,1]=\{x\in \bR\mid 0\leq x\leq 1\}$ of $\bR$ with the usual topology is compact (see Exercise~\ref{exer:Icompact}).
\end{example}

A topological space $X$ is \emphind{locally compact} provided that, for each $x\in X$ and each $U\in\Omega(X)$ with $x\in U$, there are $V\in\Omega(X)$ and $K\subseteq X$ compact such that
\[
x\in V\subseteq K\subseteq U.
\]
Note that, if $X$ is Hausdorff, then compactness implies local compactness, while this is not the case in general (see Exercise~\ref{exer:loccompact}).

\begin{proposition}\label{prop:proj-along-comp}
Let $X$ and $Y$ be topological spaces and $\pi_Y\colon X\times Y \to Y$ the projection onto the second coordinate. If $X$ is compact, then $\pi_Y$ is a closed mapping.
\end{proposition}

\begin{proof}
Let $C\subseteq X\times Y$ be closed and suppose $y\not\in\pi_Y[C]$. That is, for each $x\in X$, we have $(x,y)\not\in C$. Thus, as $C$ is closed,  for each $x\in X$, there are $U_x$ open in $X$ and $V_x$ open in $Y$ with $(x,y)\in U_x\times V_x$ and
\[
C\cap (U_x\times V_x)=\emptyset.
\]
Since $x\in U_x$ for each $x\in X$, the collection $\{U_x\mid x\in X\}$ is an open cover of $X$. Also, as $X$ is compact, there is a finite subset $M\subseteq X$ such that $\{U_x\mid x\in M\}$ covers $X$. Now
letting $V=\bigcap\{V_x\mid x\in M\}$ we have $y\in V$, $V\subseteq Y$ open, and $V\cap\pi_Y[C]=\emptyset$.
\end{proof}

We now recall a result, which requires a non-constructive principle, and is often useful for proving compactness.

\begin{theorem}[Alexander Subbase Theorem]\label{thm:alexander}
Let $X$ be a topological space and $\cS$ a subbase for the topology on $X$. If every cover $\cU\subseteq\cS$ of $X$ has a finite subcover, then $X$ is compact.
\end{theorem}

Any proof of the Alexander Subbase Theorem must necessarily use a non-constructive principle, and we here use Zorn's Lemma, which will also be crucially used in Chapter~\ref{ch:priestley} in the proof of the Stone prime filter-ideal theorem (Theorem~\ref{thm:DPF}). Note that while the latter has historically been called a ``Lemma'', we treat it here as a postulate, as is common, since it is equivalent in Zermelo-Fraenkel set theory to the Axiom of Choice.  We have made the choice to freely use non-constructive principles like Zorn's Lemma in this book, although they could sometimes have been avoided.  For more information about equivalences between various choice principles, see for example the books~\cite{Jec1973, How1998}.

\begin{lemma}[Zorn's Lemma]\index{Zorn's Lemma}\index{Axiom of Choice}\label{lem:zorn}
Let $S$ be a non-empty partially ordered set such that if $C \subseteq S$ is totally ordered, then there exists an upper bound $c$ of $C$ in $S$. Then $S$ has a maximal element, that is, there exists $s \in S$ such that for any $s' \in S$, if $s' \geq s$, then $s' = s$.
\end{lemma}
\begin{proof}[Proof that Zorn's Lemma implies Alexander Subbase Theorem.] We ask the reader to fill in the gaps of this proof sketch in Exercise~\ref{exer:alexander}. We prove the contrapositive statement. Suppose that $X$ is not compact. Zorn's Lemma guarantees that there exists a \emph{maximal} open cover $\cC$ of $X$ which does not have a finite subcover. The subcollection $\cC \cap \cS$ of $\cC$ can be shown to still be a cover of $X$, using the maximality of $\cC$. Now $\cC \cap \cS$ is a cover by elements from the subbase, which can not have a finite subcover.
\end{proof}

An earlier, equivalent, variation of Zorn's Lemma, that we will sometimes use, is the following.
\begin{lemma}[Hausdorff maximality principle]\index{Hausdorff maximality principle}\label{lem:hausdorff-maximality} 
  For any totally ordered subset $C$ of a partially ordered set $P$, there is a totally ordered subset $C' \supseteq C$ which is maximal among totally ordered subsets of $P$, with respect to subset inclusion.
\end{lemma}

\ourexercises

\begin{ourexercise}\label{exer:real}
Consider $\bR$ equipped with the collection of subsets $U\subseteq\bR$ with the property that, for each $x\in U$, there exists $\varepsilon>0$ such that the interval $(x-\varepsilon,x+\varepsilon)$ is entirely contained in $U$.
\begin{enumerate}
\item Show that $(\bR,\tau)$ is a topological space;
\item Show that the collection of intervals $(r,s)$, where both $r$ and $s$ are rational, forms a base for $\bR$.
\item Show that $f\colon\bR\to\bR$ is continuous if, and only if, for every $x\in\bR$ and every $\varepsilon>0$, there exists $\delta>0$ such that, for all $y\in\bR$ with $|x-y|<\delta$ we have $|f(x)-f(y)|<\varepsilon$.
\end{enumerate}
\end{ourexercise}

\begin{ourexercise}\label{exer:open-closed-map}
Give examples of continuous maps which are:
\begin{enumerate}
\item neither open nor closed,
\item open but not closed,
\item closed but not open.
\end{enumerate}
\end{ourexercise}

\begin{ourexercise}\label{exer:int-adjoint}
Let $X$ be a topological space. Prove that the interior map $\intr \colon \mathcal{P}(X) \to \Omega(X)$ is upper adjoint to the inclusion map $\iota \colon \Omega(X) \into \mathcal{P}(X)$, and that the closure map $\cl \colon \mathcal{P}(X) \to \cC(X)$ is  lower adjoint to the inclusion $\iota' \colon \cC(X) \into \mathcal{P}(X)$.
\end{ourexercise}

\begin{ourexercise}\label{exer:emb-homeo-subsp}
Let $X$ and $Y$ be topological spaces, and $f\colon X\to Y$ a continuous injection.
\begin{enumerate}
\item Show that $f$ need not be an embedding;
\item Show that $f$ is an embedding if, and only if, $f$ co-restricted to $\im(f)$ is a homeomorphism. Here, the \emphind{co-restriction} of $f\colon X\to Y$ to a subspace $S\subseteq Y$, where $\im(f)\subseteq S$, is the function with domain $X$, codomain $S$, and the same action as $f$.
\item Show that if the continuous map $f$ is a bijection, then it is a homeomorphism if, and only if, it is open and if, and only if, it is closed.
\item Show that if $X$ is a subset of $Y$ and $f \colon X \into Y$ is the inclusion map, then $X$ is a subspace of $Y$ if, and only if, $f$ is an embedding.
\end{enumerate}
\end{ourexercise}

\begin{ourexercise}\label{exer:base}
Let $(X,\tau)$ be a topological space, and $\cB$ and $\cS$ be, respectively, a base and a subbase for $\tau$.
\begin{enumerate}
\item Show that ${\langle}\cB{\rangle}=\{\bigcup\cC\mid \cC\subseteq\cB\}$.
\item \label{itm:base-from-subbase} Let $\cT$ be the closure of $\cS$ under finite intersections. Show that $\cT$ is a base for $\tau$.
\item Give an example of a base for a topology which is not closed under binary intersections.
\item Let $Y$ be a topological space. Show that $f\colon Y\to X$ is continuous if, and only if, $f^{-1}(U)$ is open in $Y$ for each $U\in\cS$.
\end{enumerate}
\end{ourexercise}

\begin{ourexercise}\label{exer:covers} \
  Let $X$ be a topological space. Show that all of the following statements are equivalent:
  \begin{enumerate}[label=(\roman*)]
    \item for all $\mathcal{R} \subseteq \Omega(X)$, if $\bigcup \mathcal{R} = X$, then there exists a finite $\mathcal{R}' \subseteq \mathcal{R}$ such that $\bigcup \mathcal{R}' = X$;
    \item for all directed $\mathcal{D} \subseteq \Omega(X)$, if $\bigcup \mathcal{D} = X$, then $X \in \mathcal{D}$;
    \item for all $\mathcal{S} \subseteq \cC(X)$, if $\bigcap \mathcal{S} = \emptyset$, then there exists a finite $\mathcal{S}' \subseteq \mathcal{S}$ such that $\bigcap \mathcal{S}' = \emptyset$;
    \item \label{itm:fin-int-prop} for all filtering $\mathcal{F} \subseteq \cC(X)$, if $\bigcup \mathcal{F} = \emptyset$ then $\emptyset \in \mathcal{F}$;
    \item for all $\mathcal{R} \subseteq \Omega(X)$, $\mathcal{S} \subseteq \cC(X)$, if $\bigcap \mathcal{S} \subseteq \bigcup \mathcal{R}$, then there exist finite $\mathcal{R}' \subseteq \mathcal{R}$ and finite $\mathcal{S}' \subseteq \mathcal{S}$ such that $\bigcap \mathcal{S}' \subseteq \bigcup \mathcal{R}'$;
    \item for all directed $\mathcal{D} \subseteq \Omega(X)$ and filtering $\mathcal{F} \subseteq \cC(X)$, if $\bigcap \mathcal{F} \subseteq \mathcal{D}$, then there exist $F \in \mathcal{F}$ and $U \in \mathcal{D}$ such that $F \subseteq U$.
  \end{enumerate}
  \end{ourexercise}


\begin{ourexercise}\label{exer:comp}
Suppose $(X,\tau)$ is a compact topological space.
\begin{enumerate}
\item Let $C\subseteq X$ be closed. Show that $C$ is compact.
\item Find a compact space $X$ with a compact subset $K$ which is not closed.
\item Show that if $X$ is Hausdorff and $K\subseteq X$ is compact, then $K$ is closed.
\end{enumerate}
\end{ourexercise}

\begin{ourexercise}\label{exer:coherence} 
  \begin{enumerate}
  \item Show that the finite union of compact sets is compact;
  \item Give an example to show that the intersection of two compact sets need not be compact.
  \end{enumerate}
  \end{ourexercise}


\begin{ourexercise}\label{exe:filter-compact}
  Let $X$ be a compact topological space, and let $\cF$ be a filter of clopen subsets of $X$. Prove that, for any clopen set $K$ in $X$,
  \[ K \in \cF \iff \bigcap \cF \subseteq K \ . \] 
  \hint{Use the characterization of compactness given in Exercise~\ref{exer:covers}.}
\end{ourexercise}

\begin{ourexercise}\label{exer:alexander}
  This exercise is based on \cite[Exercise~3.12.2, p.~221]{Engelking89} and asks you to fill in the details of the proof of the Alexander Subbase Theorem. 
  Let $X$ be a topological space and $\cS$ a subbase for the topology on $X$.
  \begin{enumerate}
  \item Prove that, in the following sub-poset of $(\mathcal{P}(\Omega(X)), \subseteq)$, 
    \[\mathbf{C} := \{ \cC \subseteq \Omega(X) \ : \ \cC \text{ is a cover of } X \text{ and } \cC \text{ has no finite subcover}\},\] any totally ordered subset $(\cC_i)_{i \in I}$ of $\mathbf{C}$ has an upper bound in $\mathbf{C}$. \hint{Show that $\bigcup_{i \in I} \cC_i$ is in $\mathbf C$.}

    Suppose $X$ is not compact. Then by Zorn's Lemma we can pick a maximal element $\cM$ of the poset $\mathbf{C}$.
    \item Prove that, for any open sets $U, V$, if $V \in \cM$ and $U \subseteq V$, then $U \in \cM$.

    \item Prove that, for any finite number of open sets $U_1, \dots U_n$, if $U_i \not\in \cM$ for every $1 \leq i \leq n$, then $\bigcap_{i=1}^n U_i \not\in \cM$. \hint{Use the maximality of $\cM$ to get finite subcovers $\mathcal{F}_i$ of $\cM \cup \{U_i\}$ for every $i$, and show that 
    $\mathcal{F} := \bigcup_{i=1}^n \mathcal{F}_i \cup\{\bigcap_{i=1}^n U_i\}$ is then a finite subcover of $\cM \cup \{\bigcap_{i=1}^n U_i\}$.}
      \item Conclude that $\cM \cap \cS$ is a cover of $X$ that does not have a finite subcover.
  \end{enumerate}
\end{ourexercise}

\begin{ourexercise}\label{exer:sep}
Show that $T_4$ implies $T_3$, which implies $T_2$, which implies $T_1$, which implies $T_0$. Further, show that all these implications are strict.
\end{ourexercise}

\begin{ourexercise}\phantomsection\label{exer:Hausdorff}
\begin{enumerate}
  \item
Show that a topological space $X$ is Hausdorff if, and only if, the diagonal relation
\[
\Delta_X=\{(x,x)\mid x\in X\}
\]
is closed in the product topology on $X\times X$.

\item Let $f \colon Y \to X$ be a continuous function and suppose that $X$ is Hausdorff. Prove that the \emphind{graph} of $f$, that is, the relation 
  \[ \{(y,f(y)) \ \mid \ y \in Y\}\]
  is closed in the product topology on $Y \times X$.
\item Conclude that a topological space $X$ is Hausdorff if, and only if, the graph of any continuous function $Y \to X$ is closed.
\end{enumerate}
\end{ourexercise}

\begin{ourexercise}\label{exer:easy-product}
  Let $X$ and $Y$ be topological spaces. 
  \begin{enumerate}
    \item Prove that, for any $y_0 \in Y$, the function $i \colon X \to Y \times X$ defined by $i(x) := (x,y_0)$ is continuous.
    \item Prove that the projection maps $\pi_X \colon X \times Y \to X$ and $\pi_Y \colon X \times Y \to Y$ are open.
  \end{enumerate}
\end{ourexercise}

\begin{ourexercise}\label{exer:Icompact}
Show that the subspace $[0,1]=\{x\in \bR\mid 0\leq x\leq 1\}$ of $\bR$ with the usual topology is compact.
\end{ourexercise}

\begin{ourexercise}\label{exer:compHaus-normal}
Show that any compact Hausdorff space is normal. \hint{Show first that it is regular.}
\end{ourexercise}

\begin{ourexercise}\label{exer:Tychonoff}
(Tychonoff's Theorem) Show that the product of compact spaces is again compact. \hint{A proof can be found in any standard reference on topology, for example \cite[Theorem 3.2.4, p.~138]{Engelking89}.}
\end{ourexercise}

\begin{ourexercise}\phantomsection\label{exer:loccompact}
\begin{enumerate}
	\item Show that any compact Hausdorff space is locally compact. \hint{Use that compact Hausdorff spaces are regular, as proved in Exercise~\ref{exer:compHaus-normal}.}
	\item Find a topological space which is compact but not locally compact.
\end{enumerate}
\end{ourexercise}

\begin{ourexercise}[Quotient space]\label{exer:quotspace}
Let $X$ be a topological space and $\equiv$ an equivalence relation on $X$. The
\emphind{quotient space} of $X$ by $\equiv$ is the space based on $X/{\equiv}$
whose open sets are those $U\subseteq X/{\equiv}$ such that
\[
q^{-1}(U)=\bigcup\{[x]_\equiv\mid [x]_\equiv\in U\}=\{x\in X\mid [x]_\equiv\in U\}
\]
is open in $X$, where $q\colon X\to X/{\equiv}, x\mapsto [x]_\equiv$ is the canonical quotient map.
\begin{enumerate}
	\item Show that the topology on $X/{\equiv}$ is the finest topology on $X/{\equiv}$ making $q\colon X\to X/{\equiv}$ continuous.
	\item Show that $X/{\equiv}$ is a $T_1$ space if, and only if, every equivalence class of $\equiv$ is closed in $X$.
	\item Show that if $X/{\equiv}$ is a Hausdorff space, then $\equiv$ is necessarily a closed as a subset of the product space $X\times X$. 	
\item Show that if the quotient map is open, then $X/{\equiv}$ is Hausdorff if, and only if, $\equiv$ is a closed in $X\times X$. 	
	\item Show that if $f\colon X\twoheadrightarrow Y$ is a continuous surjection, then $f$ factors through the canonical quotient map $q\colon X\to X/\ker(f)$ by a unique continuous bijection $\tilde{f}\colon X/\ker(f)\to Y$.
	\item Give an example in which $\tilde{f}$ is not a homeomorphism.
	\item Show that if $f$ is open or closed then $\tilde{f}$ is a homeomorphism. But show by giving an example that this condition is not necessary.
\end{enumerate}
\end{ourexercise}

\section{Topology and order}\label{sec:TopOrd}

Let $(X,\tau)$ be a topological space. The \emphind{specialization order}, of $\tau$ is the binary relation $\leq_\tau$ on $X$ defined by
\[
    x\leq_\tau y \iff \text{for every }  U\in\Omega(X),\ \,\text{if } x\in U \text{
    then } y\in U\ .
\]
\nl{$\leq_\tau$}{specialization order associated with a topology $\tau$}{}
This relation $\leq_\tau$ is clearly reflexive and transitive and thus, for any topological space, the specialization order is a preorder on $X$. It is not hard to see that it is a partial order if, and only if, $X$ is $T_0$. Moreover, for any $y \in X$, the principal down-set, ${\downarrow}y$, of $y$ in the specialization order is the closure of the singleton set $\{y\}$. In particular, $T_1$ spaces can be characterized as those spaces having a trivial specialization order. The reader is asked to prove these statements in Exercise~\ref{exer:specorder}. 

 In analysis and algebraic topology, the spaces studied are almost always Hausdorff, so in these fields the interaction with order theory is minimal. However, in applications of topology to algebra and logic, almost all spaces are $T_0$ but not $T_1$, so in these areas, as in theoretical computer science, the interaction of topology and order plays an important role. 

A subset of a topological space $X$ is said to be \emphind{saturated} provided it is an up-set in the specialization order, or equivalently, provided it is an intersection of opens. Note that a subset $K\subseteq X$ is compact if, and only if, its saturation ${\uparrow} K$ is compact. We denote by $\KS(X)$ the collection of \emphind{compact-saturated} subsets of $X$ (that is, subsets of $X$ that are both compact and saturated). As we will see later on, beyond the Hausdorff setting, but in the presence of compactness, $\KS(X)$ is in many aspects the right generalization of the closed subsets. The following fact is often useful and illustrates the consequence of compactness in terms of the specialization order. Recall that, for a subset $S$ of a poset, $\min S$ is the (possibly empty) set of minimal points of $S$. 
\nl{$\KS(X)$}{the collection of compact-saturated subsets of a topological space $X$}{}

\begin{proposition}\label{prop:compact-implies minpoints}
Let $X$ be a $T_0$ space and $K\subseteq X$ compact, then $K\subseteq {\uparrow}\min(K)$, where we consider $X$ in its specialization order. In particular, if $K$ is compact-saturated then $K= {\uparrow}\min(K)$. 
\end{proposition}

\begin{proof}
We first show that that if $D\subseteq K$ is down-directed in $(X,\leq_\tau)$, then there is a lower bound of $D$ in $K$. Let $D$ be a down-directed set in $(X,\leq_\tau)$, and suppose that $K$ contains no lower bounds of $D$. 
Then the directed collection of open subsets $\{({\downarrow }x)^c\mid x\in D\}$ is an open cover $K$. Therefore, by compactness, it follows that there is $x\in D$ with $K\subseteq ({\downarrow }x)^c$. In particular, $x\not\in K$ and thus $D\not\subseteq K$.

Now let $x\in K$ and, by the Hausdorff Maximality Principle (Lemma~\ref{lem:hausdorff-maximality}), let $C$ be a maximal chain in $K$ containing $x$. Then, by the above argument, $C$ has a lower bound $x'\in K$. Now by maximality of $C$, it follows that $x'\in\min(K)$ and thus $K\subseteq {\uparrow}\min(K)$.
\end{proof}

\subsection*{The lattice of all topologies on a set}

Let $X$ be a set. Note that the collection  
\[
\TopLat(X):=\{\tau\in\cP(\cP(X))\mid \tau \text{ is a topology}\}.
\]
is closed under arbitrary intersections and thus (see Exercise~\ref{exe:complattsuff} in Chapter~\ref{ch:order}) it is a complete lattice in the inclusion order. Infima are given by intersections, while suprema are given by the topologies generated by unions. The least topology on $X$ is the indiscrete topology, while the largest is the discrete topology. We will often make use of the binary \emphind{join of topologies} on a given set $X$.

The interaction of compactness and the Hausdorff separation axiom is illuminated by looking at $\TopLat(X)$ as a complete lattice. Indeed, by inspecting the definitions (see Exercise~\ref{ex:down-and-up-in-Top}), note that
\[
T_2(X):=\{\tau\in\TopLat(X)\mid \tau \text{ is } T_2\}
\]
is an up-set in $\TopLat(X)$, while
\[
T_{\it Comp}(X):=\{\tau\in\TopLat(X)\mid (X,\tau) \text{ is compact}\}
\]
is a down-set in $\TopLat(X)$. It follows that the set of compact-Hausdorff topologies on a set $X$ is a convex\index{convex} subset of $\TopLat(X)$. The following very useful result tells us that it is in fact an anti-chain.

\begin{proposition}\label{prop:compHaus-incomp}
Let $X$ be a set and $\sigma$ and $\tau$ topologies on $X$ with $\sigma\subseteq\tau$. If $\sigma$ is Hausdorff and $\tau$ is compact, then $\sigma=\tau$.
\end{proposition}

\begin{proof}
Since $T_2(X)$ is an up-set and $\TopLat_{\it Comp}(X)$ is a down-set, the hypotheses on $\sigma$ and $\tau$ imply that both are simultaneously compact and Hausdorff. 
Note that in any $T_1$ space, and thus in particular in a Hausdorff space, any set is saturated.
Moreover, in compact-Hausdorff spaces, being closed is equivalent to being compact (and saturated), see Exercise~\ref{exer:compHaus}. Also, any set which is compact in a bigger topology remains so in the smaller topology.
Thus, we have the following sequence of (bi)implications for any subset $U\subseteq X$:
\begin{align*}
U\in\tau  &\iff U^c\in \cC(X,\tau)\\
               & \iff U^c\in \KS(X,\tau) \\
               & \implies U^c\in \KS(X,\sigma)\\
               & \iff U^c\in \cC(X,\sigma)\iff U\in\sigma.\qedhere
\end{align*}
\end{proof}

\subsection*{Order-topologies}
From topology we get order, but it is also possible to go the other way. Especially in computer science applications where second-order structure such as a topology is difficult to motivate, topologies induced by orders play an important role; see also the applications to domain theory in Chapter~\ref{chap:DomThry}.

As we have seen, for any topological space, the closures of points are equal to their principal down-sets for the specialization order. Thus, if a topological space $X$ has specialization order $\leq$, then \emph{at least} each set of the form $({\downarrow} x)^c$, for $x \in X$, must be open. We now proceed in the converse direction. Let $(P,\leq)$ be a partially ordered set. We define several topologies on $P$ for which the specialization order coincides with $\leq$.
\begin{itemize}
\item The \emphind{upper topology} on $P$, $\iota^\uparrow(P)$, is defined as the least topology in which ${\downarrow} p$ is closed for every $p\in P$. That is, the upper topology is given by
\[
\iota^\uparrow(P)=\langle({\downarrow} p)^c\mid p\in P\rangle.
\]
\nl{$\iota^\uparrow(P)$}{the upper topology on a poset $P$}{}
\item The  \emphind{Scott topology} on $P$ consists of those up-sets which are inaccessible by directed suprema. That is, an up-set $U\subseteq P$ is Scott open if, and only if, $\bigvee D\in U$ implies $U\cap D\neq\emptyset$ for all directed subsets $D\subseteq P$. We denote the Scott topology on $P$ by $\sigma(P)$.
\nl{$\sigma(P)$}{the Scott topology on a poset $P$}{}

\item The \emphind{Alexandrov topology} on $P$ is the largest topology on $P$ yielding $\leq$ as its specialization order. That is,
\[
\alpha(P)=\{ U\subseteq P\mid U\text{ is an up-set}\}.
\]
\nl{$\alpha(P)$}{the Alexandrov topology on a poset $P$}{}
\end{itemize}
The upper, the Scott, and the Alexandrov topologies all have the original order $\leq$ as their specialization order. In fact, if we denote by $\TopLat(P,\leq)$ the complete lattice of topologies on $P$ yielding $\leq$ as their specialization order, it is not hard to see that this is precisely the closed interval $[\iota^\uparrow(P),\alpha(P)]$ in $\TopLat(P)$.

Clearly, there are order-dual definitions for each of these topologies, which have the reverse of $\leq$ as their specialization order. For example, the \emphind{lower topology} on $P$ is defined by 
\[
\iota^\downarrow(P)=\langle({\uparrow} p)^c\mid p\in P\rangle.
\] The \emphind{dual Alexandrov topology} has all down-sets as open sets. One can also consider the order-dual of the Scott topology but this is not so common, as the motivation for having closed sets which are stable under directed joins comes from a model of computing in which a computation is considered as the directed join of all its partial computations or finite approximations; more on this in Chapter~\ref{chap:DomThry}.

Using the above `one-sided' topologies as building blocks, we now also define a number of `two-sided' topologies on a partially ordered set $(P,\leq)$.  These two-sided topologies are all $T_1$ and thus have trivial specialization order.
\begin{itemize}
\item The \emphind{interval topology} on $P$ is the join of the upper and lower topologies. That is,
\[
\iota(P):=\iota^\uparrow(P)\vee\iota^\downarrow(P).
\]
The usual topology on the reals is in fact the interval topology given by the usual order on the reals.
\item The \emphind{Lawson topology} on $P$ is the join of the Scott and the lower topologies. That is,
\[
\lambda(P):=\sigma(P)\vee\iota^{\downarrow}(P).
\]
\item Note that the join of the Alexandrov and dual Alexandrov topologies is the \emphind{discrete topology} on $P$.
\end{itemize}

\ourexercises

\begin{ourexercise}
  \label{ex:down-and-up-in-Top}
  Prove that $T_2(X)$ is an up-set in $\TopLat(X)$, and that $T_{Comp}(X)$ is a down-set in $\TopLat(X)$.
\end{ourexercise}
\begin{ourexercise}\label{exer:closedmap}
Let $X$ be a compact space, $Y$ a Hausdorff space, and $f\colon X\to Y$ be a continuous map.
\begin{enumerate}
\item Show that the map $f$ is closed. \hint{Use Proposition~\ref{prop:compHaus-incomp}.}
\item Show that if $f$ is a bijection, then it is a homeomorphism.
\item Show that  the  co-restriction of $f$ to its image $X\twoheadrightarrow\im(f), x\mapsto f(x)$ is  a quotient map. In particular, as soon as $f$ is surjective, it is a quotient map (see Exercise~\ref{exer:quotspace}).
\end{enumerate}
\end{ourexercise}

\begin{ourexercise}\label{exer:specorder}
Let $X$ be a topological space.
\begin{enumerate}
\item Show that the specialization order on $X$ is a preorder;
\item Show that $X$ is $T_0$ if, and only if, the specialization order on $X$ is a partial order;
\item Show that $X$ is $T_1$ if, and only if, the specialization order on $X$ is trivial. That is, $x\leq y$ if, and only if, $x=y$;
\item Show that $x\leq y$ in the specialization order if, and only if, $x\in\overline{\{y\}}$. That is, $\overline{\{y\}}={\downarrow}y$;
\item Show that a subset $S\subseteq X$ is an intersection of open sets if, and only if, it is an up-set in the specialization order.
\end{enumerate}
\end{ourexercise}

\begin{ourexercise}\label{exer:interval-ordertop}
Let $X$ be a set and $\leq$ an order on $X$. Show that a topology $\tau$ on $X$ has $\leq$ as its specialization order if, and only if,
\[
\iota^\uparrow(X,\leq)\subseteq\tau\subseteq\alpha(X,\leq).
\]
\end{ourexercise}

\begin{ourexercise}\label{exer:cont-implies-op}
Show that if a function $f\colon X\to Y$ between topological spaces is continuous, then it is order-preserving with respect to the specialization orders on $X$ and $Y$. Give an example to show that the converse is false.
\end{ourexercise}

\begin{ourexercise}\label{exer:TOPfin}
Show that if $P$ is a finite partially ordered set then $\iota^\uparrow(P)=\alpha(P)$. Conclude that, in a finite $T_0$ space, any up-set is open. Further show that, for any two partially ordered sets $P$ and $Q$, a map $F\colon P\to Q$ is order-preserving if, and only if, it is continuous with respect to the Alexandrov topologies on $P$ and $Q$. 

\emph{Note.} Using terminology that we will introduce in Definition~\ref{dfn:iso-between-cats}, this exercise shows that the category of finite partially ordered sets is isomorphic to the category of finite $T_0$ topological spaces. For further details, see Example~\ref{exa:topfiniso}.
\end{ourexercise}

\begin{ourexercise}\label{exer:compHaus}
Let $X$ be a compact Hausdorff space and $K\subseteq X$. Show that $K$ is closed if, and only if, it is compact if, and only if, it is compact-saturated.
\end{ourexercise}

\section{Compact ordered spaces}\label{sec:comp-ord-sp}
In this section we show that there is an isomorphism between certain compact spaces equipped with an order, first introduced by Nachbin, and certain $T_0$ spaces known as stably compact spaces. These spaces provide a well-behaved generalization of compact Hausdorff spaces and also contain the topological spaces dual to distributive lattices, which we call \emph{spectral} spaces, and which are the main object of study of this book, together with \emph{Priestley} spaces, their order-topological counterpart.

\begin{definition}\label{dfn:comp-ord-sp}
An \emphind{ordered space} is a triple $(X,\tau,\leq)$ such that
\begin{itemize}
\item $(X,\tau)$ is a topological space;
\item $(X,\leq)$ is a partially ordered set;
\item $\leq\ \subseteq X\times X$ is closed in the product topology.
\end{itemize}
An ordered space is said to be a \emphind{compact ordered space} provided the underlying topological space is compact.
\end{definition}
\begin{example}\label{exa:kord-finite}
  The set $\{0,1\}$, equipped with the discrete topology and the usual order, is a compact ordered space. More generally, any finite poset, equipped with the discrete topology, is a compact ordered space.
\end{example}

A \emph{morphism}\index{ordered space!morphism}\index{morphism between ordered spaces} from an ordered space $(X, \tau_X, \leq_X)$ to an ordered space $(Y, \tau_Y, \leq_Y)$ is a function $f \colon X \to Y$ that is both continuous as a map from the space $(X, \tau_X)$ to $(Y,\tau_Y)$ and order preserving as a map from the poset $(X, \leq_X)$ to $(Y, \leq_X)$. An \emphind{order-homeomorphism} between ordered spaces is a morphism that is both a homeomorphism and an order-isomorphism; this is sometimes also called \emph{isomorphism}\index{isomorphism between ordered spaces}. For an equivalent definition of order-homeomorphism, see Exercise~\ref{exer:orderhomeo}.

\begin{proposition}
Let $X$ be an ordered space. Then the underlying topological space is Hausdorff.
\end{proposition}

\begin{proof}
Since $X$ is an ordered space, $\leq$  is closed in $X\times X$ equipped with the product topology. Thus $\geq$ is also closed in $X\times X$ with the product topology and it follows that $\Delta_X=\,\leq\cap\geq$ is closed in $X\times X$ with the product topology. But this is equivalent to $X$ being Hausdorff (see Exercise~\ref{exer:Hausdorff}).
\end{proof}

The following proposition is an important technical tool in the study of compact ordered spaces.

\begin{proposition}\label{prop:cos-closed}
Let $X$ be a compact ordered space and $C\subseteq X$ a closed subset of $X$. Then ${\uparrow}C$ and ${\downarrow}C$ are also closed. In particular, ${\uparrow}x$ and ${\downarrow}x$ are closed for all $x\in X$.
\end{proposition}

\begin{proof}
If $C\subseteq X$ is closed in $X$, then $C\times X$ is closed in $X\times X$ equipped with the product topology. Now, as $X$ is an ordered space it follows that $\leq$ is closed and thus $(C\times X)\,\cap\leq$ is closed in $X\times X$. Consider the set 
\[
\pi_2[(C\times X)\,\cap\leq]={\uparrow}C,
\]
where $\pi_2\colon X\times X\to X$ is the projection on the second coordinate. By 
Proposition~\ref{prop:proj-along-comp} it follows that it is closed in $X$. Projecting $(C\times X)\,\cap\geq$ on the second coordinate shows that ${\downarrow}C$ is closed. Finally, as $X$ is Hausdorff, it is in particular $T_1$ and thus the singletons $x$ are all closed. It follows that ${\uparrow}x$ and ${\downarrow}x$ are closed for all $x\in X$.
\end{proof}

We can now derive the following very useful \emphind{order-separation property} for compact ordered spaces.

\begin{proposition}\label{prop:cos-ordsep}
Let $X$ be a compact ordered space. For all $x,y\in X$, if $x\nleq y$, then there are disjoint sets $U,V\subseteq X$ with $U$ an open up-set containing $x$ and $V$ an open down-set containing $y$.
\end{proposition}

\begin{proof}
Let $x,y\in X$ with $x\nleq y$. Then ${\uparrow}x$ and ${\downarrow} y$ are disjoint. Also, by Proposition~\ref{prop:cos-closed}, the sets ${\uparrow}x$ and ${\downarrow} y$ are closed. Now, since $X$ is compact ordered, it is compact Hausdorff and therefore also normal (see Exercise~\ref{exer:compHaus-normal}).  Thus there are open disjoint sets $U,V\subseteq X$ with ${\uparrow}x\subseteq U$ and ${\downarrow} y\subseteq V$. Finally, let $U'=({\downarrow}U^c)^c$ and $V'=({\uparrow}V^c)^c$, then one may verify that $U'$ is an open up-set, $V'$ is an open down-set, and we have
\[
x\in U' \ \text{ and }\ y\in V' \ \text{ and }\ U'\cap V'=\emptyset.\qedhere
\]
\end{proof}
To any compact ordered space, we now associate two $T_0$ spaces. These spaces are not $T_1$ as long as the order on $X$ is non-trivial. To be specific, if $(X,\tau,\leq)$ is a compact ordered space, then we define
\[
\tau^{\uparrow}=\tau\cap \Up(X,\leq)
\]
\nl{$\tau^{\uparrow}$}{the open up-set topology associated with an ordered space $(X,\tau,\leq)$}
and
\[
\tau^{\downarrow}=\tau\cap \Down(X,\leq).
\]
\nl{$\tau^{\downarrow}$}{the open down-set topology associated with an ordered space $(X,\tau,\leq)$}
In other words, $\tau^{\uparrow}$ is the intersection of the topology $\tau$ and the Alexandrov topology on $(X,\leq)$ and $\tau^{\downarrow}$ is the intersection of the topology $\tau$ and the dual Alexandrov topology on $(X,\leq)$. Accordingly, $\tau^{\uparrow}$  and $\tau^{\downarrow}$ are indeed topologies on $X$. We will often denote the topological space underlying the original ordered space $(X,\tau, \leq)$ simply by $X$, the space $(X,\tau^{\uparrow})$ by $X^{\uparrow}$, and the space $(X,\tau^{\downarrow})$ by $X^{\downarrow}$.
\nl{$X^{\uparrow}$}{the topological space with open up-set topology associated with an ordered space $X$}{}
\nl{$X^{\downarrow}$}{the topological space with open down-set topology associated with an ordered space $X$}{}

Note that if $(X,\tau,\leq)$ is a compact ordered space, then so is $(X,\tau,\geq)$. Thus any property of the spaces $X^{\uparrow}$ and their relation to $\leq$ implies that the order-dual property is true for the spaces $X^{\downarrow}$ and we will not always state both.

\begin{proposition}\label{prop:compord-specord}
Let $(X,\tau,\leq)$ be a compact ordered space. Then the specialization order of $X^{\uparrow}$ is $\leq$ and in particular $X^{\uparrow}$ is a $T_0$ space.
\end{proposition}

 \begin{proof}
 For each $x\in X$, ${\downarrow}x$ is closed in $(X,\tau)$ and it is a down-set. Therefore ${\downarrow}x$ is closed in $X^{\uparrow}$ and it follows that $\iota^{\uparrow}(X,\leq)\subseteq \tau^{\uparrow}$. Also, clearly $\tau^{\uparrow}\subseteq \Up(X,\leq)=\alpha(X,\leq)$ and thus the specialization order of $X^{\uparrow}$ is $\leq$ (see Exercise~\ref{exer:interval-ordertop}).
 \end{proof}

 A crucial fact, given in the following proposition, which will enable us to come back to a compact ordered space $X$ from $X^{\uparrow}$, is the fact that $X^{\uparrow}$ and $X^{\downarrow}$ are inter-definable by purely topological means without using the data of the original compact ordered space.

\begin{proposition}\label{prop:Fdown=Kup}
Let $(X,\tau,\leq)$ be a compact ordered space. Then
\[
V\in\tau^{\downarrow} \quad\iff\quad V^c\in\KS(X^{\uparrow}).
\]
\end{proposition}

 \begin{proof}
Note that $\cC(X^{\downarrow})=\cC(X,\tau)\cap\,\Up(X,\leq)$. Also, by definition, $\KS(X^{\uparrow})$ consists of those subsets of $X$ that are both compact with respect to $\tau^{\uparrow}$ and belong to $\Up(X,\leq)$. Thus we need to show that if $S\in\Up(X,\leq)$ then $S$ is closed relative to $\tau$ if, and only if, it is compact relative to $\tau^{\uparrow}$.

Let $S\in\Up(X,\leq)$. If $S$ is  closed relative to $\tau$, then $S$ is compact relative to $\tau$ (see Exercise~\ref{exer:comp}). But then it is also compact relative to the smaller topology  $\tau^{\uparrow}$ as required. For the converse, suppose now that $S$ is compact relative to $\tau^{\uparrow}$ and let $y\not\in S$. For each $x\in S$, since $x\nleq y$, by Proposition~\ref{prop:cos-ordsep}, there are disjoint sets $U_x,V_x\subseteq X$ with $U_x$ an open up-set containing $x$ and $V_x$ an open down-set containing $y$. It follows that $\{U_x\}_{x\in S}$ is an open cover of $S$ relative to $\tau^{\uparrow}$. Thus by compactness, there is a finite subset $F\subseteq S$ such that $\{U_x\}_{x\in F}$ covers $S$. Let
\[
V=\bigcap\{V_x\mid x\in F\}
\]
then $V$ is disjoint from the union of $\{U_x\}_{x\in F}$ and thus from $S$. Also, $V$ is open relative to $\tau^{\downarrow}$ and $y\in V$. That is, we have shown that $S$ is closed relative to $\tau^{\downarrow}$.
 \end{proof}

We are now ready to introduce a class of (unordered) topological spaces called \emph{stably compact spaces}. Stably compact spaces have a fairly complex definition but, as we will see, they are in fact none other than those spaces which occur as $X^{\uparrow}$ for $X$ a compact ordered space.

Before we give the definition (Definition~\ref{dfn:stab-comp-sp}), we need to identify two more properties of spaces, both related to the interaction of compactness and intersection. First, we call a compact space \emphind{coherent} provided the intersection of any two compact-saturated subsets is again compact. A space is called \emphind{well-filtered} provided for any filtering collection $\cF$ of compact-saturated sets and any open $U$ we have
 \[\bigcap\cF\subseteq U \quad \implies \quad \text{ there exists } K\in\cF \text{ such that } K\subseteq U.\]
One can show that if $X$ is well-filtered, then the collection of compact-saturated subsets of $X$ is closed under filtering intersections (see Exercise~\ref{exer:well-filtered}).  Notice that if a space is both coherent and well-filtered then the collection of compact-saturated sets is actually closed under arbitrary intersections. Also notice that compact Hausdorff spaces have both these properties since the compact-saturated sets are just the closed sets.

\begin{definition}\label{dfn:stab-comp-sp}
 A \emphind{stably compact space} is a topological space which is $T_0$, compact, locally compact, coherent, and well-filtered.
\end{definition}
 
 \begin{example}\label{ex:Sierpinski}
 Any finite $T_0$ space is stably compact. A particularly important stably compact space is the \emphind{Sierpinski space}
 \[
 \mathbb{S} =(\{0,1\},\{\emptyset,\{1\},\{0,1\}\}).
 \]
 Note that the topology on $\mathbb{S}$ is equal to $\tau^{\uparrow}$, where $(\{0,1\},\tau,\leq)$ is the compact ordered space described in Example~\ref{exa:kord-finite}.
 \end{example}
 
 Using Proposition~\ref{prop:Fdown=Kup} it is not hard to see that $X^\uparrow$ is stably compact whenever $X$ is a compact ordered space. We will now show that this is in fact one direction of a one-to-one correspondence between compact ordered spaces and stably compact spaces. To this end we need the notion of the co-compact dual of a topology.
 
 Let $\tau$ be a topology on a set $X$. The \emphind{co-compact dual} of $\tau$, denoted $\tau^\partial$, is the topology generated by the complements of compact-saturated subsets of $(X,\tau)$. That is,
 \nl{$\tau^\partial$}{the co-compact dual of a topology $\tau$}{}
 \[
 \tau^\partial:=\langle \{K^c\mid K\in\KS(X,\tau)\}\rangle.
 \]
 The compact-saturated sets are always closed under finite unions, so the collection of their complements is closed under finite intersection and is thus a base for $\tau^\partial$. In the case of a stably compact space, the compact-saturated sets are also closed under arbitrary intersections (see Exercise~\ref{exer:well-filtered}), so the collection of their complements is already a topology and we have
 \[
 \tau^\partial=\{ K^c\mid K\in\KS(X,\tau)\}.
 \]
 Further we may define the \emphind{patch topology} obtained from $\tau$ to be
\nl{$\tau^p$}{the patch topology obtained from a topology $\tau$}{}
 \[
 \tau^p=\tau\vee\tau^\partial.
 \]
 We can now show how to get back the topology of a compact ordered space $X$ from the topology of $X^\uparrow$.
 
 \begin{proposition}\label{prop:compordaspatch}
 Let $(X,\tau,\leq)$ be a compact ordered space. Then $(\tau^\uparrow)^\partial=\tau^\downarrow$ and
 $(\tau^\uparrow)^p=\tau$.
 \end{proposition}
 
 \begin{proof}
 The first equality is a just a restatement of Proposition~\ref{prop:Fdown=Kup} in terms of the co-compact dual topology. Once we have this, we may observe that Proposition~\ref{prop:cos-ordsep} tells us, among other things, that $\tau^\uparrow \vee (\tau^\uparrow)^\partial = \tau^\uparrow \vee \tau^\downarrow$ is a Hausdorff topology which is contained in $\tau$. But, by Proposition~\ref{prop:compHaus-incomp}, if a Hausdorff topology is below a compact topology, then in fact they are equal, so $(\tau^\uparrow)^p=\tau$ as desired.
 \end{proof}

 \begin{theorem}\label{thrm:COSpace-StabCompSp}
 The assignments
 \[
 (X,\tau,\leq)\ \mapsto \ (X,\tau^\uparrow)
 \]
 and
 \[
 \qquad(X,\rho)\quad \mapsto \ (X,\rho^p,\leq_\rho)
 \]
 establish a one-to-one correspondence between compact ordered spaces and stably compact spaces.%
 \end{theorem}

 \begin{proof}
 It is left as Exercise~\ref{exer:Xupstabcomp} to show that if $(X,\tau,\leq)$ is a compact ordered space, then $X^\uparrow$ is stably compact. Here, we show that if $(X,\rho)$ is a stably compact space then $(X,\rho^p,\leq_\rho)$ is compact ordered.
 
 As a first step, we show that $\leq_\rho$ is closed relative to the product topology induced by the patch topology. Let $x,y\in X$ with $x\nleq_\rho y$. Then, by definition of the specialization order, there is an open set $U\in\rho$ with $x\in U$ and $y\not\in U$. By local compactness of $\rho$ there exist $V\in\rho$ and $K\subseteq X$ which is compact such that $x\in V\subseteq K\subseteq U$. Since ${\uparrow}K$ is also compact and $K\subseteq {\uparrow}K\subseteq U$ since $U$ is an up-set in the specialization order, we may assume, without loss of generality, that $K\in\KS(X,\rho)$. It follows that both $V$ and $K^c$ are in $\rho^p$, that $x\in V$, $y\in K^c$, and, because $V$ is an up-set and $V\cap K^c=\emptyset$, we have
 \[
 (V\times K^c)\;\cap \leq_{\rho} \; =\emptyset.
 \]
 Thus we have shown that an arbitrary element $(x,y)$ of the complement of $\leq_\rho$ lies in a basic open of the product topology which is disjoint from $\leq_\rho$ as required.
 
 Next we show that $(X,\rho^p)$ is compact. By the Alexander Subbase Theorem, it suffices to show that covers by subbasic opens have finite subcovers. We use the subbase for $\rho^p$ given by
 \[
 \cS=\{K^c\mid K\in\KS (X,\rho)\}\cup\rho.
 \]
 Now let $\cC\subseteq\cS$ be a cover of $X$. Define
 \[
 \cC_\KS=\{K\in\KS (X,\rho)\mid K^c\in\cC\}
 \]
 and
 \[
 \cC_\rho=\cC\cap\rho.
 \]
 Then the fact that $\cC$ covers $X$ implies that
 \[
 \bigcap\cC_\KS\subseteq\bigcup\cC_\rho.
 \]
 Notice that $U=\bigcup\cC_\rho$ is open in $(X,\rho)$. Let $\cC_{\KS}'$ denote the closure of $\cC_{\KS}$ under finite intersections, that is, 
 \[
 \cC'_\KS :=\big\{\bigcap \cL\mid \cL\subseteq \cC_\KS \text{ is finite} \big\}
 \]
 Note that $\cC'_\KS$ is, by compactness and coherence, a collection of compact-saturated subsets of $(X,\rho)$, which is moreover filtering. Thus, by well-filteredness, it follows that there is a finite $ \cL\subseteq \cC_\KS$ with
 \[
 \bigcap\cL\subseteq\bigcup\cC_\rho.
 \]
 Now, this means that $\cC_\rho$ is an open cover of the compact-saturated set $\bigcap\cL$. By compactness there is a finite subcover $\cC'_\rho$ of $\bigcap\cL$. This in turn is equivalent to saying that
 \[
 \cC'=\{K^c\mid K\in \cL\}\cup\cC'_\rho
 \]
 is a finite subcover of the original cover $\cC$ as required.
 
 Thus it follows that $(X,\rho^p,\leq_\rho)$ is indeed a compact ordered space whenever $(X,\rho)$ is a stably compact space.

 Further, the combination of Proposition~\ref{prop:compord-specord} and Proposition~\ref{prop:compordaspatch} implies that the composition of the two assignments gives the identity on compact ordered spaces. It remains to show that the reverse composition yields the identity of stably compact spaces. 
  To this end, let $(X,\rho)$ be a stably compact space. Since $\rho\subseteq\rho^p$ and $\rho\subseteq\Up(X,\leq_\rho)$, it follows that $\rho\subseteq\rho^p\cap\,\Up(X,\leq_\rho)$. The fact that $\rho^p\cap\,\Up(X,\leq_\rho)\subseteq\rho$ follows from Exercise~\ref{exer:SCSsubbase}.
 \end{proof}
 
 While the correspondence of the above theorem provides an isomorphism for objects, the natural classes of maps for compact ordered spaces and stably compact spaces are \emph{not} the same. A natural notion of structure-preserving map for stably compact spaces is that of a continuous function. For compact ordered spaces the natural notion of structure preserving map is that of a function which is simultaneously continuous and order preserving.
It is not hard to see that every continuous and order-preserving map between compact ordered spaces is continuous for the corresponding stably compact spaces. However, the converse is not true in general. In fact, the continuous and order-preserving maps between compact ordered spaces correspond to the so-called proper maps between stably compact spaces (see Exercise~\ref{exer:propermaps}). We will show in Example~\ref{exa:stcomp-kord-iso} of Chapter~\ref{ch:categories} that this correspondence formally yields an isomorphism of categories.

We finish this section by recording a useful `translation' between various properties of subsets of a compact ordered space and its corresponding stably compact space. The proof is left as an instructive exercise in applying Theorem~\ref{thrm:COSpace-StabCompSp}.
\begin{proposition}\label{prop:with-without-order-scs}
Let $(X,\tau,\leq)$ be a compact ordered space, with $(X,\tau^\downarrow)$ and $(X,\tau^\uparrow)$ the stably compact spaces of $\tau$-open up-sets and $\tau$-open down-sets, respectively. For any subset $S$ of $X$,
\begin{enumerate}
  \item $S$ is saturated in $(X,\tau^{\uparrow})$ if, and only if, $S$ is an up-set in $(X,\tau,\leq)$ if, and only if, the complement of $S$ is saturated in $(X,\tau^{\downarrow})$;
  \item $S$ is closed in $(X,\tau^{\uparrow})$ if, and only if, $S$ is a closed down-set in $(X,\tau,\leq)$ if, and only if, $S$ is compact and saturated in $(X,\tau^{\downarrow})$;
  \item \label{itm:koisclup} $S$ is compact and open in $(X,\tau^{\uparrow})$ if, and only if, $S$ is a clopen up-set in $(X,\tau,\leq)$ if, and only if, the complement of $S$ is compact and open in $(X,\tau^{\downarrow})$.
\end{enumerate}
\end{proposition}

\ourexercises

\begin{ourexercise}\label{exe:compordaltdef}
Let $(X,\tau,\leq)$ be a triple such that $(X,\tau)$ is a topological space and $(X,\leq)$ is a poset. Prove that the following are equivalent:
\begin{enumerate}
\item[(i)] The order $\leq$ is closed in $(X,\tau) \times (X,\tau)$.
\item[(ii)] for every $x, y \in X$, if $x \nleq y$, then there exist open subsets $U$ and $V$ of $(X,\tau)$ such that $x \in U$, $y \in V$, and ${\uparrow}U \cap {\downarrow}V = \emptyset$.
\end{enumerate}
\end{ourexercise}

\begin{ourexercise}\label{exer:orderhomeo}
Let $f \colon X \to Y$ be a morphism between ordered spaces. Prove that $f$ is an order-homeomorphism if, and only if, there exists a morphism $g \colon Y \to X$ such that $g \circ f = \id_X$ and $f \circ g = \id_Y$.
\end{ourexercise}

\begin{ourexercise}\label{exer:ordnorm}
Ordered variants of regularity and normality also hold for compact ordered spaces. In particular, prove that if $C$ and $D$ are closed subspaces of a compact ordered space with
\[
{\uparrow} C\cap{\downarrow}D=\emptyset
\]
then there are disjoint sets $U,V\subseteq X$ with $U$ an open up-set containing $C$ and $V$ an open down-set containing $D$.
\end{ourexercise}

\begin{ourexercise}\label{ex:PXtopology}
  Let $X$ be a set.
  \begin{enumerate}
  \item Show that the set $2^X$ of functions  from $X$ to $2$ is bijective with $\cP(X)$.
  \item Show that the product topology on $2^X$, where $2$ is equipped with the Sierpinski topology, translates to the topology on $\cP(X)$ given by the subbase $\{\eta_x\mid x\in X\}$, where $\eta_x=\{S\subseteq X\mid x\in S\}$, and by the base $\{\eta_F\mid F\subseteq X\text{ is finite}\}$, where $\eta_F=\{S\subseteq X\mid F\subseteq S\}$.
  \item Show that the product topology on $2^X$, where $2$ is equipped with the discrete topology, translates to the topology on $\cP(X)$ given by the subbase $\{\eta_x,\mu_x\mid x\in X\}$, where $\mu_x=\{S\subseteq X\mid x\not\in S\}$, and by the base $\{\eta_F\cap\mu_G\mid F,G\subseteq X\text{ are finite}\}$, where $\mu_G=\{S\subseteq X\mid G\cap S=\emptyset\}$.
  \end{enumerate}
  \end{ourexercise}
  
\begin{ourexercise}\label{exer:unitint}
Show that $\mathbb{R}$ with its usual order and topology is an ordered space. Show that the unit interval $[0,1]$ is a compact ordered space.
\end{ourexercise}

\begin{ourexercise}\label{exer:prodCOS}
Show that $\{0,1\}$, equipped with the discrete topology and the usual order, is a compact ordered space. 
Also show that compact ordered spaces are closed under arbitrary Cartesian product, where the product is equipped with the product topology and the coordinate-wise order. Conclude in particular that for any set $X$, the space $\{0,1\}^X$ is a compact ordered space.
\end{ourexercise}

\begin{ourexercise}\label{exer:COSadjunction}
Let $(X,\tau,\leq)$ be a compact ordered space.
\begin{enumerate}
\item Show that  the inclusion $\Up(X,\leq)\hookrightarrow\cP(X)$ has a lower adjoint given by ${\uparrow}(\ )\colon\cP(X)\twoheadrightarrow\Up(X,\leq)$ and an upper adjoint given by $S\mapsto ({\downarrow}S^c)^c$.
\item Conclude that Proposition~\ref{prop:Fdown=Kup} implies that the inclusion $\KS(X^{\uparrow})\hookrightarrow\KS(X)=\cC(X)$ has a lower adjoint given by $C\mapsto{\uparrow}C$ and that the inclusion $\tau^{\uparrow}\hookrightarrow \tau$ has an upper adjoint given by $U\mapsto({\downarrow}U^c)^c$.
\end{enumerate}
\textit{Note.} We thank J\'er\'emie Marqu\`es for suggesting this reformulation of Proposition~\ref{prop:Fdown=Kup} in terms of upper and lower adjoints on opens and compact-saturated sets.
\end{ourexercise}

\begin{ourexercise}\phantomsection\label{exer:well-filtered}
  \begin{enumerate}
  \item Show that if $X$ is a well-filtered space then the intersection of any filtering collection of compact-saturated sets is again compact-saturated.
  \item Show that, if a collection $\mathcal{F}$ of subsets of $\mathcal{P}(X)$ is filtering and closed under finite intersections (that is, for any finite $\mathcal{S} \subseteq \mathcal{F}$, $\bigcap \mathcal{S} \in \mathcal{F}$), then $\mathcal{F}$ is closed under arbitrary intersections.
  \item Conclude that if $X$ is well-filtered, compact, and coherent, then any intersection of compact-saturated sets is again compact-saturated.
  \end{enumerate}
\end{ourexercise}

\begin{ourexercise}\label{exer:compHaus-stabcomp} 
 Show that if $X$ is a compact Hausdorff space, then $X$ is stably compact.
\end{ourexercise}

\begin{ourexercise}\label{exer:with-without-order}
Deduce Proposition~\ref{prop:with-without-order-scs} from Proposition~\ref{prop:compordaspatch} and Theorem~\ref{thrm:COSpace-StabCompSp}.
\end{ourexercise}

\begin{ourexercise}\label{exer:Xupstabcomp} \
 Show that if $X$ is a compact ordered space, then $X^\uparrow$ is stably compact.
\end{ourexercise}

\begin{ourexercise}
Let $(X,\tau)$ be a $T_0$ space. Show that $\leq_{\tau^\partial}\ =\ \geq_\tau$.
\end{ourexercise}

\begin{ourexercise}\label{exer:SCSsubbase} \
Let $(X,\tau)$ be a stably compact space, $\cB\subseteq\tau$, and $\cK\subseteq\KS(X,\tau)$ such that, for all $x, y \in X$, if $x \nleq y$, then there exist $U \in \cB$ and $K \in \cK$ such that $x \in U$ and $y \not\in K$.
Show that $\mathcal B$ is a subbase for $\tau$.
\end{ourexercise}

\begin{ourexercise}\label{exer:SCSretracts}
Show that continuous retracts of stably compact spaces are again stably compact. Here a \emphind{retract} of a topological space $X$ is a continuous function $f\colon X\to X$ with $f \circ f=f$.
\end{ourexercise}

\begin{ourexercise}\label{exer:propermaps}
Let $f\colon X\to Y$ be a continuous function between topological spaces. We call $f$ \emphind{proper}\index{proper map}\index{continuous!proper} provided the following two properties hold:
\begin{enumerate}
\item  ${\downarrow}f(C)$ is closed whenever $C\subseteq X$ is closed;
\item $f^{-1}[K]$ is compact for any $K\subseteq Y$ which is compact-saturated.
\end{enumerate}
Now let $X$ and $Y$ be compact ordered spaces and denote the corresponding pair of stably compact spaces by $X^{\uparrow}, X^{\downarrow}$ and $Y^{\uparrow}, Y^{\downarrow}$, respectively. Further let  $f\colon X\to Y$ be a function between the underlying sets. Show that the following conditions are equivalent:
\begin{enumerate}[label=(\roman*)]
\item the function $f$, viewed as a map between compact ordered spaces, is continuous and order preserving;
\item the function $f$, viewed as a map between the stably compact spaces $X^{\uparrow}$ and $Y^{\uparrow}$, is proper;
\item the function $f$, viewed as a map between the stably compact spaces $X^{\downarrow}$ and $Y^{\downarrow}$, is  proper.
 \end{enumerate}
\end{ourexercise}

\notessec
For readers who need further introduction to topology and for general topology references beyond what we have included here, we recommend  a classical book on General Topology such as \cite{Engelking89}. For a comprehensive reference on topology and order we recommend \cite{Getc2003}.

The correspondence between stably compact spaces and compact ordered spaces in the form given in Theorem~\ref{thrm:COSpace-StabCompSp} originates with the first edition of the Compendium on Continuous Lattices \parencite{Getc80}. More focused sources presenting the correspondence and the relation with the co-compact dual of a topology are \cite{Jung2004,Lawson2011}.
\chapter{Priestley duality}\label{ch:priestley}

In this chapter, we show how to extend the duality for finite distributive
lattices given in Chapter~\ref{ch:order} to \emph{all} distributive lattices.
The two key ideas, due to Stone, are to generalize the join/meet-prime elements of the
finite case to \emph{prime filters/ideals}
(Section~\ref{sec:primefiltersideals}), and to introduce \emph{topology} on the
structure dual to a distributive lattice (Section~\ref{sec:topologize}), which
leads us to a dual equivalence or \emph{duality}. The main technical tool is
Stone's Prime Filter-Ideal Theorem, Theorem~\ref{thm:DPF}. In this chapter we
elaborate a modern variant of Stone's original duality, \emph{Priestley
duality}. The precise connection between Priestley's duality and Stone's
original duality will be made in Theorem~\ref{thm:Stone-isom-Priestley} in
Chapter~\ref{chap:Omega-Pt}. After treating distributive lattice duality, we
show in the final Section~\ref{sec:boolenv-duality} how the more widely known
duality for Boolean algebras follows as an easy consequence. Throughout the
chapter, in-between introducing the general concepts and proving results about
them, we show how to compute dual spaces of several distributive lattices, as
`running examples'. We encourage the reader to work through the examples and
accompanying exercises in detail, as we believe this is crucial for developing
an intuition for dual spaces.

\section{Prime filters and ideals}\label{sec:primefiltersideals}
Recall from Section~\ref{sec:finDLduality} that every finite distributive
lattice can be represented by the poset of its join-irreducible elements. The
following example of an infinite distributive lattice shows why
join-irreducibles cannot be used in general.
\begin{example}\label{exa:nojoinirr}
\nl{$P \oplus Q$}{ordered sum of two posets $P$ and $Q$, with $P$ below $Q$}{}
Consider the lattice $L = \{\bot\} \oplus (\mathbb{N}^\op)^2$, depicted in
Figure~\ref{fig:nojoinirr}. Every non-bottom element $(m,n)$ in this lattice is the join
of the two elements $(m+1,n)$ and $(m,n+1)$ strictly below it. Therefore, there
are no join-irreducibles in $L$, $\cJ(L) = \emptyset$.
\begin{figure}[htp]
\begin{center}
\begin{tikzpicture}
\node[draw,circle,inner sep=1pt,fill,label=right:{$(0,0)$}] (1) at (0,3) {};
\node[draw,circle,inner sep=1pt,fill] (a) at (-.5,2.5) {};
\node[draw,circle,inner sep=1pt,fill] (b) at (.5,2.5) {};
\node[draw,circle,inner sep=1pt,fill] (c) at (-1,2) {};
\node[draw,circle,inner sep=1pt,fill] (d) at (0,2) {};
\node[draw,circle,inner sep=1pt,fill] (e) at (1,2) {};
\node[] (x) at (-1.5,1.5) {};
\node[draw,circle,inner sep=1pt,fill] (f) at (-.5,1.5) {};
\node[draw,circle,inner sep=1pt,fill] (g) at (.5,1.5) {};
\node[] (w) at (1.5,1.5) {};
\node[] (y) at (-1,1) {};
\node[draw,circle,inner sep=1pt,fill] (h) at (0,1) {};
\node[] (v) at (1,1) {};
\node[] (z) at (-.5,.5) {};
\node[] (u) at (.5,.5) {};
\node[draw,circle,inner sep=1pt,fill,label=right:{$\bot$}] (0) at (0,-.5) {};
\draw[thin] (1) -- (a);
\draw[thin] (1) -- (b);
\draw[thin] (a) -- (c);
\draw[thin] (a) -- (d);
\draw[thin] (b) -- (d);
\draw[thin] (b) -- (e);
\draw[thin] (c) -- (f);
\draw[thin] (d) -- (f);
\draw[thin] (d) -- (g);
\draw[thin] (e) -- (g);
\draw[thin] (f) -- (h);
\draw[thin] (g) -- (h);
\draw[thin] (f) -- (h);
\draw[thin,dotted] (c) -- (x);
\draw[thin,dotted] (f) -- (y);
\draw[thin,dotted] (h) -- (z);
\draw[thin,dotted] (h) -- (u);
\draw[thin,dotted] (g) -- (v);
\draw[thin,dotted] (e) -- (w);
\end{tikzpicture}
\end{center}
\label{fig:nojoinirr}
\caption{A distributive lattice with no join-irreducible elements}
\end{figure}
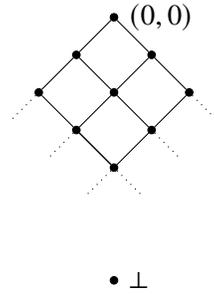
\end{example}
To find the correct notion that should replace `join-irreducible' in the case
of infinite distributive lattices, we need to change our perspective from
specific \emph{elements} of the lattice to specific \emph{subsets} of the
lattice. If $L$ is a distributive lattice and $j \in \cJ(L)$, then $j$ can be
uniquely represented by the collection $F_j$ of elements greater than or equal
to $j$. Elements of $F_j$ can be thought of as `approximations' of the element
$j$, which grow in precision as one moves downward in the set $F_j$; indeed,
$j$ is the infimum of this set $F_j$. In the case of infinite lattices, while
join-irreducibles $j$ themselves may fail to exist (Example~\ref{exa:nojoinirr}
above), these `approximating sets' $F$ will still exist. The formal notion of
`approximating set', for a join-irreducible, is that of a \emph{prime filter},
defined as follows.
\begin{definition}\index{filter}\index{filter!prime}\index{prime filter}
Let $L$ be a distributive lattice. A subset $F$ of $L$ is called a
\emph{filter} if it is non-empty, an up-set, and for any $a, b \in L$, if $a
\in F$ and $b \in F$, then $a \wedge b \in F$. A filter is called \emph{proper}
if $F \neq L$, or, equivalently, $\bot \not\in F$. A filter $F$ is called
\emph{prime} provided that $F$ is proper, and, for any $a, b \in L$, if $a \vee
b \in F$, then $a \in F$ or $b \in F$.
\end{definition}
One way to think of filters in a lattice is that a filter represents an `idealized element', in much the same way as Noether's ideals in rings. From this point of view, a non-principal filter stands for a meet that does not exist in the lattice. Note that if $S$ is a subset of a lattice and $S \subseteq S'$, if infima of both sets exist, then $\bigwedge S' \leq \bigwedge S$. Thus, since taking a meet over a larger set yields a smaller element, it is natural to postulate the \emph{reverse} inclusion on filters, as we will do here. The dual notion of \emph{ideal} in a lattice is introduced in Definition~\ref{def:ideal} below; there, the order of (non-reversed) subset inclusion is the natural one, since ideals in a lattice stand for an idealized join, and taking joins over larger sets yield larger elements.

With the point of view that filters stand for idealized meets, \emph{prime} filters stand for idealized join-irreducible elements. More precisely, in finite lattices, prime filters are in one-to-one correspondence with join-prime elements: to any join-prime element $j$, one may associate the prime filter $F_j := {\uparrow}j$ of elements greater than or equal to $j$ (see Exercise~\ref{exe:joinprime-primefilt}). Notice that, also here, $j \leq j'$ if, and only if, $F_{j'} \subseteq F_j$. Thus, this correspondence is an \emph{isomorphism of posets} if one equips the set of prime filters with the partial order of \emph{reverse} inclusion. 
All this motivates the following definition of a \emph{partial order on the set of filters}.

\begin{definition}\label{def:filterorder}\index{filters!partial order on}
If $F$ and $F'$ are filters in a lattice $L$, we say that $F$ \emph{is below} $F'$, $F \leq F'$ if, and only if, $F'$ is a subset of $F$. We denote by $\Filt(L)$ the poset of filters of $L$, and by $\PrFilt(L)$ the poset of prime filters of $L$.
\nl{$\Filt(L)$}{the poset of filters of a lattice $L$}{}
\nl{$\PrFilt(L)$}{the poset of prime filters of a lattice $L$}{}
\end{definition}

The following examples of infinite distributive lattices and their posets of prime filters will be used as running examples in this chapter.
\begin{example}\label{exa:Nplus1squared}
  Write $\bN \oplus 1$ for the total order on the set $\bN \cup \{\omega\}$ which extends $\bN$ by setting $n \leq \omega$ for all $n \in \bN \cup \{\omega\}$. All prime filters of $\bN \oplus 1$ are principal; they are the sets of the form $F_n := {\uparrow} n$ for $n \in (\bN \oplus 1) \setminus \{0\}$. Their order is the same as that of $\bN \oplus 1$; this is a highly exceptional case of a distributive lattice that is order-isomorphic to its poset of prime filters, via the isomorphism that sends $\omega$ to $F_\omega$ and $n \in \bN$ to $F_{n+1}$. %

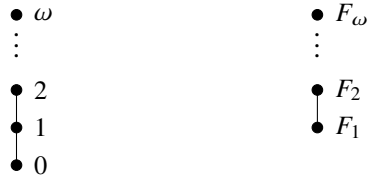
\begin{figure}[htp]
  \begin{center}
  \begin{tikzpicture}
  \polab{(0,0)}{\small $0$}{right}
  \polab{(0,.5)}{\small $1$}{right}
  \polab{(0,1)}{\small $2$}{right}
  \node at (0,1.7) {$\vdots$};
  \polab{(0,2)}{\small $\omega$}{right}
  \draw (0,0) -- (0,1);

  \polab{(4,.5)}{\small $F_{1}$}{right}
  \polab{(4,1)}{\small $F_{2}$}{right}
  \node at (4,1.7) {$\vdots$};
  \polab{(4,2)}{\small $F_\omega$}{right}
  \draw (4,0.5) -- (4,1);
  \end{tikzpicture}
\end{center}
  \caption{The distributive lattice $\mathbb{N} \oplus 1$ and its poset of prime filters.}
  \label{fig:nplusone}
  \end{figure}

 Now consider the Cartesian product $(\bN \oplus 1)^2 := (\bN \oplus 1) \times (\bN \oplus 1)$. Since $\bN \oplus 1$ is a distributive lattice, so is $(\bN \oplus 1)^2$. The prime filters of $(\bN \oplus 1)^2$ are the principal filters of the form $F_{(n,m)} = {\uparrow} (n, m)$, where $n = 0$ or $m = 0$ and not both. Indeed, for any $(n,m) \in (\bN \oplus 1)^2$, $(n,m) \leq (n,0) \vee (0,m)$, so if ${\uparrow} (n,m)$ is prime, then we must have $m = 0$ or $n = 0$. The partial order on the prime filters is the order of the disjoint union, that is, $(n,m) \leq (n',m')$ when either $n \leq n'$ and $m = m' = 0$, or $n = n' = 0$ and $m \leq m'$.
  \end{example}
  The fact that the prime filters of a Cartesian product split as a disjoint union is not a coincidence (see Exercise~\ref{exe:sumproductdualposet}).

\begin{example}\label{exa:NplusNop}
Consider the totally ordered set $\mathbb{N} \oplus \mathbb{N}^\op$, the ordered sum of $\mathbb{N}$ and its opposite; that is, the underlying set is $\{ (n, i)\mid n \in \mathbb{N}, i \in \{0,1\}\}$, and the order is given by $(n,i) \leq (m,j)$ if one of the following holds:
\begin{itemize}
  \item $i < j$, or
  \item $i = j = 0$ and $n \leq m$, or
  \item $i = j = 1$ and $m \leq n$.
\end{itemize}
In this distributive lattice, every $(n,i) \neq (0,0)$ is join irreducible, so $F_{(n,i)} = {\uparrow}(n,i)$ is a prime filter. Moreover, there is one non-principal prime filter, $F_{\omega} := \{(n, 1)\mid n \in \mathbb{N}\}$. The order $\leq$ on prime filters is inherited from the order on the lattice, and $F_{\omega} \leq F_{(n,i)}$ iff $i = 1$; see Figure~\ref{fig:nplusnop}.

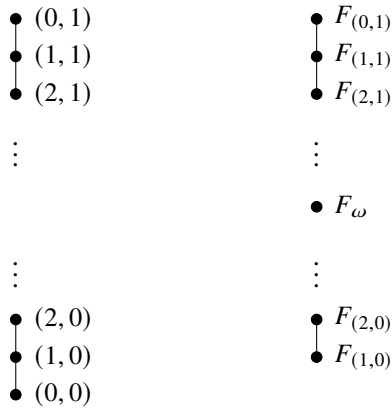
\begin{figure}[htp]
  \begin{center}
  \begin{tikzpicture}
  \polab{(0,0)}{\small $(0,0)$}{right}
  \polab{(0,.5)}{\small $(1,0)$}{right}
  \polab{(0,1)}{\small $(2,0)$}{right}
  \node at (0,1.7) {$\vdots$};
  \polab{(0,5)}{\small $(0,1)$}{right}
  \polab{(0,4.5)}{\small $(1,1)$}{right}
  \polab{(0,4)}{\small $(2,1)$}{right}
  \node at (0,3.3) {$\vdots$};
  \draw (0,0) -- (0,1);
  \draw (0,5) -- (0,4);

  \polab{(4,.5)}{\small $F_{(1,0)}$}{right}
  \polab{(4,1)}{\small $F_{(2,0)}$}{right}
  \node at (4,1.7) {$\vdots$};
  \polab{(4,2.5)}{\small $F_\omega$}{right}
  \polab{(4,5)}{\small $F_{(0,1)}$}{right}
  \polab{(4,4.5)}{\small $F_{(1,1)}$}{right}
  \polab{(4,4)}{\small $F_{(2,1)}$}{right}
  \node at (4,3.3) {$\vdots$};
  \draw (4,0.5) -- (4,1);
  \draw (4,5) -- (4,4);
  \end{tikzpicture}
\end{center}
  \caption{The distributive lattice $\mathbb{N} \oplus \mathbb{N}^\op$ and its poset of prime filters.}
  \label{fig:nplusnop}
  \end{figure}
\end{example}

The following example shows explicitly that prime filters in distributive lattices generalize prime numbers. They are in fact closely related to prime ideals of rings (see Exercise~\ref{exe:primes-ring}).
\begin{example}\label{exa:prime-divisibility}
Consider the set $\bN$, equipped with the partial order $|$ of \emphind{divisibility}, defined by $n \mid m$ iff there exists $q \in \mathbb{N}$ such that $m = qn$; in particular, $n \mid 0$ for all $n \in \mathbb{N}$. The partial order $(\bN, |)$ is a bounded distributive lattice: it has $1$ as its bottom, $0$ as its top, and for any $m, n \in \mathbb{N}$, $m \wedge n$ is the greatest common divisor of $m$ and $n$, and $m \vee n$ is the least common multiple of $m$ and $n$.

The prime filters of $(\bN, |)$ are $F_0 := \{0\}$, and the sets $F_{p^k} := \{n \in \mathbb{N} \ : \ p^k \mid n\}$, for every prime number $p$ and $k \geq 1$. It is easy to verify that the $F_{p^k}$ are indeed prime filters (see Exercise~\ref{exe:prime-divisibility}); we show that every prime filter is of this form. Suppose that $F$ is a prime filter of $(\bN, |)$ and $F \neq F_0$. Let $m$ be the minimal non-zero element of $F$. We show first that $F = {\uparrow} m$. Clearly, ${\uparrow} m \subseteq F$. For the converse inclusion, let $n \in F$ be arbitrary, $n > 0$. Then $n \wedge m \in F$ since $F$ is a filter. The greatest common divisor of $n$ and $m$ is non-zero and $\leq m$, and therefore equal to $m$, by minimality of $m$. We now show that $m = p^k$ for some prime $p$ and $k \geq 1$. First note that $m > 1$ since $F$ is proper. Pick a prime divisor $p$ of $m$ and pick $k$ maximal such that $p^k \mid m$. Then $m = p^k \vee \frac{m}{p^k}$, so, since $F$ is a prime filter, we must have either $p^k \in F$ or $\frac{m}{p^k} \in F$. But $\frac{m}{p^k} < m$, so it cannot be an element of $F$ by minimality of $m$. Thus, $p^k \in F$, and $0 < p^k \leq m$, so we get $m = p^k$, again by minimality of $m$. 

The partial order on the set of prime filters is given by, for any $p, q$ prime and $k, \ell \geq 1$, $F_{p^k} < F_0$, and $F_{p^k} \leq F_{q^\ell}$ if, and only if, $p = q$ and $k \leq \ell$; see Figure~\ref{fig:prime-naturalnumbers} below.
\begin{figure}[htp]
\begin{center}
\begin{tikzpicture}
  \polab{(-2,1)}{$F_4$}{left}
  \polab{(-2,0)}{$F_2$}{below}
  \draw (-2,0) -- (-2,1.5);
  \node at (-2, 2) {$\vdots$};

  \polab{(0,0)}{$F_3$}{below}
  \polab{(0,1)}{$F_9$}{left}
  \draw (0,0) -- (0,1.5);
  \node at (0, 2) {$\vdots$};
  
  \polab{(2,0)}{$F_5$}{below}
  \polab{(2,1)}{$F_{25}$}{left}
  \draw (2,0) -- (2,1.5);
  \node at (2, 2) {$\vdots$};

  \polab{(0,3)}{$F_0$}{right}
  \node at (4,1) {$\dots$};
\end{tikzpicture}
\end{center}
\caption{The poset of prime filters of $(\bN, |)$.}
\label{fig:prime-naturalnumbers}
\end{figure}
\end{example}

\begin{example}\label{exa:rationalintervals-primes}
Consider the set $L$ of those subsets of the real unit interval $[0,1]$ that can be written as a finite union of open rational intervals, that is, as a finite union of sets of the form $(p,q)$, $[0,q)$ or $(p,1]$, with $p, q \in \mathbb{Q} \cap [0,1]$. Note that, since a finite intersection of open rational intervals is again an open rational interval, $L$ is a sublattice of $\mathcal{P}([0,1])$ under the inclusion ordering, and therefore it is in particular a distributive lattice. 
There are three kinds of prime filters in the distributive lattice $L$:
\begin{itemize}
  \item for every $q \in [0,1]$, $F_q := \{U \in L \ \mid \ q \in U\}$;
  \item for every $q \in [0,1) \cap \mathbb{Q}$, $F_q^+ := \{U \in L \ \mid \ \text{for some } q < q', (q, q') \subseteq U\}$;
  \item for every $q \in (0,1] \cap \mathbb{Q}$, $F_q^- := \{U \in L \ \mid \ \text{for some } q' < q, (q', q) \subseteq U\}$.
\end{itemize}
One may prove (see Exercise~\ref{exe:rationalintervals}) that every prime filter of $L$ is of one of the above forms, that each $F_q$ is maximal in the partial order, while $F_q^+$ and $F_q^-$ are strictly below $F_q$. The partial order is depicted in Figure~\ref{fig:cherries}.
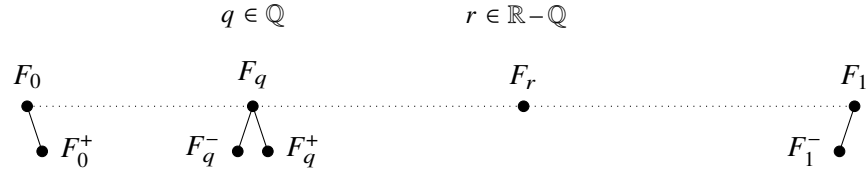
\begin{figure}
  \begin{center}
    \begin{tikzpicture}
      \node at (0,1.2) {\small $q \in \mathbb{Q}$};
      \node at (3.5, 1.2) {\small $r \in \mathbb{R} \setminus \mathbb{Q}$};
      \draw[dotted] (-3,0) -- (8,0);
      \polab{(-3,0)}{$F_0$}{above}
      \polab{(-2.8,-.6)}{$F_0^+$}{right}
      \polab{(7.8,-.6)}{$F_1^-$}{left}
      \polab{(8,0)}{$F_1$}{above}
      \polab{(0,0)}{$F_q$}{above}
      \polab{(-.2,-.6)}{$F_q^-$}{left}
      \polab{(.2,-.6)}{$F_q^+$}{right}
      \polab{(3.6, 0)}{$F_r$}{above}

      \draw (-2.8,-.6) -- (-3,0);
      \draw (7.8,-.6) -- (8,0);
      \draw (-.2,-.6) -- (0,0) -- (.2,-.6);
    \end{tikzpicture}
  \end{center}
  \caption{The poset of prime filters of the rational intervals in $[0,1]$}
  \label{fig:cherries}
\end{figure}
\end{example}

In the rest of this section, we will substantiate the claim that any distributive lattice contains `enough' prime filters. For this purpose, and also for other applications to follow, it will be convenient to introduce the order-dual notion to prime filters: prime ideals.
\begin{definition}\label{def:ideal}\index{ideal}\index{ideal!prime}\index{prime ideal}
Let $L$ be a distributive lattice. A subset $I$ of $L$ is called an \emph{ideal} if it is non-empty, a down-set, and for any $a, b \in L$, if $a \in I$ and $b \in I$, then $a \vee b \in I$. An ideal $I$ is called \emph{proper} if $I \neq L$, or, equivalently, $\top \not\in I$. An ideal $I$ is called \emph{prime} provided that $I$ is proper, and, for any $a, b \in L$, if $a \wedge b \in I$, then $a \in I$ or $b \in I$.
\end{definition}
Notice that a subset $I \subseteq L$ is a (prime) ideal if, and only if, $I$ is a (prime) filter in $L^\op$. Thus, it follows by order-duality from Exercise~\ref{exe:joinprime-primefilt} that, for a finite lattice $L$, the prime ideals in $L$ are in a one-to-one correspondence with meet-prime elements of $L$. This correspondence associates, to any meet-prime element $m$ of $L$, the prime ideal $I_m := {\downarrow}m$ of elements less than or equal to $m$. Notice that $m \leq m'$ if, and only if, $I_m \subseteq I_{m'}$. This motivates the following definition of a \emph{partial order on the set of ideals}. If $I$ and $I'$ are ideals in a lattice $L$, we say that $I$ \emph{is below} $I'$, $I \leq I'$, if, and only if $I$ is a subset of $I'$. We denote by $\Idl(L)$ the poset of ideals of $L$, and by $\PrIdl(L)$ the poset of prime ideals of $L$.
\nl{$\Idl(L)$}{the poset of ideals of a lattice $L$}{}
\nl{$\PrIdl(L)$}{the poset of prime ideals of a lattice $L$}{}

Just as join-prime and meet-prime elements are complementary notions (see Exercise~\ref{exe:kappa} in Chapter~\ref{ch:order}), prime filters and prime ideals are complementary, in the following (literal) sense.
\begin{lemma}\label{lem:filtidlcomp}\index{prime filter!characteristic function of} \index{prime ideal!characteristic function of}\index{prime filter!complement} \index{prime ideal!complement}
Let $L$ be a lattice and $F \subseteq L$. The following are equivalent:
\begin{enumerate}
\item[(i)] The set $F$ is a prime filter;
\item[(ii)] The set $I := L \setminus F$ is a prime ideal;
\item[(iii)] The characteristic function $\chi_F \colon L \to 2$, which sends an element $a$ of $L$ to $1$ if $a \in F$, and to $0$ if $a \not\in F$, is a lattice homomorphism.
\end{enumerate}
\end{lemma}
Exercise~\ref{exe:filtidlcomp} asks you to prove Lemma~\ref{lem:filtidlcomp}.
We note also that under the equivalence given in this lemma, the partial order on prime filters, which is by definition the \emph{reverse} subset inclusion order (Definition~\ref{def:filterorder}), transfers to the usual subset inclusion order on prime ideals. Note also that, for the characteristic functions corresponding to prime filters $F$ and $F'$, we have that $F' \leq F$ if, and only if, $F \subseteq F'$ if, and only if, $\chi_{F} \leq \chi_{F'}$ in the pointwise ordering on functions.

The main result about prime filters and prime ideals, which we will prove now, is that there are `enough' of them in any distributive lattice, in the following sense.
\begin{theorem}[Stone's Prime Filter-Ideal Theorem]\label{thm:DPF}\index{Stone's Prime Filter-Ideal theorem}
Let $L$ be a distributive lattice. If $F$ is a filter in $L$ and $I$ is an ideal in $L$ such that $F \cap I = \emptyset$, then there exists a prime filter $G$ in $L$ such that $F \subseteq G$ and $G \cap I = \emptyset$.
\end{theorem}
\begin{proof}
Consider the partially ordered set
\[\mathcal{S} := \{G \in \Filt(L) \ | \ F \subseteq G \text{ and } G \cap I = \emptyset\},\]
ordered by \emph{inclusion}. This set $\mathcal{S}$ is non-empty, because it contains $F$, and if $\mathcal{C} \subseteq \mathcal{S}$ is a chain in $\mathcal{S}$, then $\bigcup_{G \in \mathcal{C}} G$ is a filter (see Exercise~\ref{exe:unionoffilters}), and it belongs to $\mathcal{S}$. By Zorn's Lemma (Lemma~\ref{lem:zorn}), pick a maximal element $G$ in $\mathcal{S}$. We prove that $G$ is prime. By Lemma~\ref{lem:filtidlcomp}, we may show equivalently that $J := L \setminus G$ is a prime ideal. Since $G$ is a filter disjoint from $I$, we immediately see that $J$ is a down-set containing $I$, and hence in particular non-empty, and if $j_1 \wedge j_2 \in J$ then $j_1 \in J$ or $j_2 \in J$. It remains to show that $J$ contains the join of any two of its elements. Let $a_1, a_2 \in J$.
We use the filters generated by $G \cup \{a_1\}$ and $G \cup \{a_2\}$ (see Exercise~\ref{exe:generatefilter}). For $i = 1,2$, the set
\[ G_i := \gen{G \cup \{a_i\}}_{\mathrm{filt}} = \{b \in L \ | \ \text{there exists } g \in G \text{ such that } g \wedge a_i \leq b\}\]
is a filter that strictly contains $G$. By the maximality of $G$ in $\mathcal{S}$, $G_i$ must intersect $I$ non-trivially. Therefore, for $i = 1,2$, pick $g_i \in G$ such that $g_i \wedge a_i \in I$.

Define $g := g_1 \wedge g_2$. Then, using distributivity,
\[ g \wedge (a_1 \vee a_2) = (g \wedge a_1) \vee (g \wedge a_2) \leq (g_1 \wedge a_1) \vee (g_2 \wedge a_2).\]
Since $g_i \wedge a_i \in I$ for $i = 1,2$, it follows from this that $g \wedge (a_1 \vee a_2) \in I$. Therefore, since $J$ contains $I$, we have $g \wedge (a_1 \vee a_2) \in J$. Since $G$ is a filter, we have $g \in J$ or $a_1 \vee a_2 \in J$. However, $g \in J$ is impossible since $g_1 \in G$ and $g_2 \in G$. Thus, $a_1 \vee a_2 \in J$, as required.
\end{proof}

We make note of the fact that the proof of Theorem~\ref{thm:DPF} relies on Zorn's Lemma, which we already encountered in the context of the Alexander subbase theorem in Chapter~\ref{chap:TopOrd}. In fact, in Zermelo-Fraenkel set theory without choice, the statement of Theorem~\ref{thm:DPF} is strictly weaker than the axiom of choice. It is equivalent to both the ultrafilter theorem for Boolean algebras and to the Alexander subbase theorem (see, for example, \cite{How1998}).

From Theorem~\ref{thm:DPF}, we obtain the following \emph{representation theorem for distributive lattices}, due to \cite{Sto1937/38}.
\begin{theorem}[Stone representation for distributive lattices]\label{thm:DLrep}
Let $L$ be a lattice. The function
\begin{align*}
\widehat{(-)} \colon &L \to \Down(\PrFilt(L)) \\
&a \mapsto \widehat{a} := \{F \in \PrFilt(L) \ | \ a \in F\}
\end{align*}
is a well-defined lattice homomorphism. Moreover $\widehat{(-)}$ is injective if, and only if, $L$ is distributive. In particular, any distributive lattice $L$ embeds into the lattice of down-sets of the poset $\PrFilt(L)$.
\end{theorem}
\begin{proof}
We leave it as Exercise~\ref{exe:hatlatticehom} to prove that $\widehat{(-)}$ is a well-defined lattice homomorphism. If it is moreover injective, then $L$ is distributive, as $L$ is then isomorphic to a sublattice of a distributive lattice. Conversely, suppose that $L$ is distributive and let $a, b \in L$ be such that $a \nleq b$. Then the filter $F := {\uparrow}a$ generated by $a$ is disjoint from the ideal $I := {\downarrow}b$ generated by $b$. By Theorem~\ref{thm:DPF}, pick a prime filter $G$ containing $F$ and disjoint from $I$. Then $G \in \widehat{a} \setminus \widehat{b}$, so $\widehat{a} \not\subseteq \widehat{b}$.
\end{proof}
\nl{$\hat{a}$}{representation of an element $a$ of a distributive lattice as a down-set of its dual space}{}
Note that we use the same notation, $\widehat{(-)}$, here as the notation we used in Section~\ref{sec:finDLduality} for the function from a finite lattice $L$ to the lattice of down-sets of join-irreducible elements of $L$. Indeed, if $L$ is a finite lattice, then for any $a \in L$,  the set $\widehat{a}$ defined here in Theorem~\ref{thm:DLrep} and the set $\widehat{a}$ in Section~\ref{sec:finDLduality} correspond to each other, under the correspondence between prime filters and join-prime elements that holds in finite lattices (see Exercise~\ref{exe:joinprime-primefilt}).

Theorem~\ref{thm:DLrep} is less satisfactory than the representation theorem for \emph{finite} distributive lattices, because it does not give an isomorphism, but only a lattice embedding of a distributive lattice into a lattice of down-sets. In order to get an isomorphism, and then a full duality, we will introduce a topology on the set of prime filters in the next section.

\ourexercises

\begin{ourexercise}\label{exe:finitemeets}
Let $L$ be a lattice and $F \subseteq L$. Prove that the following are equivalent.\index{filter!characteristic function of}
\begin{enumerate}
\item[(i)] The set $F$ is a filter.
\item[(ii)] For every finite $S \subseteq L$, we have $\bigwedge S \in F$ if, and only if, $S \subseteq F$.
\item[(iii)] The characteristic function $\chi_F \colon L \to 2$ of $F$ preserves finite meets.
\end{enumerate}
\end{ourexercise}

\begin{ourexercise}\label{exe:idealdirected}
  Let $L$ be a lattice.
  \begin{enumerate}
  \item Prove that a subset $I$ of $L$ is an ideal if, and only if, $I$ is a directed down-set.
  \item Conclude that $F \subseteq L$ is a filter if, and only if, $F$ is a filtering up-set.
    \end{enumerate}
\end{ourexercise}
\begin{ourexercise}\label{exe:joinprime-primefilt}
Let $L$ be a lattice.
\begin{enumerate}
\item Prove that, for any $a \in L$, the set $F_a := {\uparrow}a = \{b \in L \ | \ a \leq b\}$ is a filter.
\item Prove that, for any join-prime element $j$ of $L$, the set $F_j = {\uparrow}j$ is a prime filter.
\item Prove that, for any prime filter $F$, if $\bigwedge F$ exists and belongs to $F$, then $j_F := \bigwedge F$ is join prime.
\item Prove that, if $L$ is finite, then the assignments $j \mapsto F_j$ and $F \mapsto j_F$ constitute a well-defined bijection between the poset of join-prime elements of $L$ and the poset of prime filters of $L$, ordered by reverse inclusion.
\end{enumerate}
\end{ourexercise}

\begin{ourexercise}
  Prove that every prime filter of the lattice $\bN \oplus \bN^\op$ of Example~\ref{exa:NplusNop} is either principal or equal to $F_\omega$.
\end{ourexercise}

\begin{ourexercise}
  Prove that every prime filter of the lattice $(\bN \oplus 1)^2$ of Example~\ref{exa:Nplus1squared} is principal, and thus, by the argument given in that example, of the stated form.
\end{ourexercise}
  \begin{ourexercise}\label{exe:rationalunitinterval}
    Consider the total order $\mathbb{Q} \cap [0,1]$. Show that there are three types of prime filters in this distributive lattice: 
    \begin{itemize}
      \item for every $q \in \mathbb{Q} \cap (0,1]$, the principal prime filter $F_q = {\uparrow} q$,
      \item for every $q \in \mathbb{Q} \cap [0,1)$, the prime filter $G_q := {\uparrow} q \setminus \{q\}$,  
      \item for every irrational $r \in (0,1)$, the prime filter $H_r := \{ q \in \mathbb{Q} \cap [0,1] \mid r < q \}$. 
    \end{itemize}
    Make a diagram of the partial order on these prime filters.
    \end{ourexercise}

\begin{ourexercise}      \label{exe:rationalintervals}
    Show that the poset of prime filters of the lattice of rational intervals from Example~\ref{exa:rationalintervals-primes} is as stated. You may proceed as follows.
    \begin{enumerate}
      \item Show that, for every $q \in [0,1]$, $F_q$ is a prime filter which is minimal in the inclusion ordering (and thus maximal in the poset of prime filters).
      \item Show that, for every prime filter $F$, there exists a unique $q \in [0,1]$ such that $F_q \subseteq F$. \hint{For existence, consider the intersection of the closures of the $U \in F$. For uniqueness, use that any two distinct points in $[0,1]$ have a rational point in between.}
      \item Show that, for every irrational $r \in (0,1)$, if $F$ is a prime filter and $F_r \subseteq F$, then $F = F_r$.
      \item Show that, for every rational $q \in (0,1)$, $F_q^-$ and $F_q^+$ are maximal prime filters.
      \item Using the previous item, show that, if $F$ is a prime filter and $F_q \subsetneq F$ for some $q \in (0,1)$, then $F = F_q^-$ or $F = F_q^+$.
    \end{enumerate}
  \end{ourexercise}
    
  \begin{ourexercise} \label{exe:prime-divisibility}
  Prove that, for any prime number $p$ and $k\geq 1$, the set $\{ n \in \mathbb{N} \ : \ p^k \mid n\}$ is a prime filter in $(\bN, |)$.
  \end{ourexercise}

  \begin{ourexercise}  \label{exe:primes-ring} (This exercise requires familiarity with basic ring theory; see the article \cite{Ban96} and the recent monograph \cite{DicSchTre2019} for more about the link between commutative algebra and lattice theory.)

  Let $R$ be a commutative ring with unit. Recall that a \emph{ring ideal} $I$ of $R$ is a subgroup of $R$ such that, for every $i \in I$ and $r \in R$, $ri$ is in $I$. For any subset $S \subseteq R$, there is a smallest ring ideal containing $S$, denoted here $\gen{S}_R$, and consisting of all elements that can be written as $\sum_{i=1}^n r_i s_i$ for some $r_1, \dots, r_n \in R$ and $s_1,\dots,s_n \in S$. A ring ideal $I$ is called \emph{prime} if $I \neq R$ and, whenever $r, s \in R$ such that $rs \in I$, either $r \in I$ or $s \in I$. A ring ideal $I$ is called \emph{radical} if, for every $r \in R$, if $r^n \in I$ for some $n \geq 1$, then $r \in I$. (In this exercise, we use the expressions `(prime) ring ideal' and `(prime) lattice ideal' to avoid confusion between the notions. By `radical ideal' we always mean radical ring ideal.)

    \begin{enumerate}
      \item Prove that, for any set $S$, the set
      \[ \sqrt{S} := \{ r \in R \mid r^n \in \gen{S}_R \text{ for some } n \geq 1 \}\]
      is the smallest radical ideal containing $S$, called the \emph{radical ideal generated by $S$}.
      \item Prove that the collection $\mathrm{RIdl}_{fg}(R)$ of finitely generated radical ideals of $R$ is a distributive lattice under the inclusion order.
      \item Prove that, if $P$ is a prime ring ideal of $R$, then 
      \[ \phi(P) := \{ I \in \mathrm{RIdl}_{fg}(R) \mid I \subseteq P \} \]
      is a prime lattice ideal of $\mathrm{RIdl}_{fg}(R)$.
      \item Prove that $\phi$ is an order-isomorphism between the prime ring ideals of $R$ and the prime lattice ideals of $\mathrm{RIdl}_{fg}(R)$.
    \end{enumerate}
  \end{ourexercise}

  \begin{ourexercise}  \label{exe:sumproductdualposet} 
Let $L$ and $M$ be distributive lattices. Prove that the prime filters of $L
\times M$ are exactly the filters that are either of the form $F \times
M$ for some prime filter $F$ of $L$, or of the form $L \times G$ for some
prime filter $G$ of $M$. Conclude that the poset of prime filters of $L \times
M$ is the disjoint union of the posets of prime filters of $L$ and of $M$.
\end{ourexercise}

\begin{ourexercise}\label{exe:filtidlcomp}
Prove Lemma~\ref{lem:filtidlcomp}.
\end{ourexercise}

\begin{ourexercise}\label{exe:unionoffilters}
Let $\mathcal{D} \subseteq \Filt(L)$ be a \emph{directed} collection of filters in the \emph{inclusion} order, that is, $\mathcal{D}$ is non-empty, and if $G_1, G_2 \in \mathcal{D}$, then there exists $G_0 \in \mathcal{D}$ such that $G_1 \cup G_2 \subseteq G_0$. Prove that $\bigcup_{G \in \mathcal{D}} G$ is a filter. Conclude in particular that the union of a non-empty chain of filters is a filter.
\end{ourexercise}

\begin{ourexercise}\label{exe:generatefilter}
Let $L$ be a lattice.\index{filter!generated}
\nl{$\gen{S}_{\mathrm{filt}}$}{filter generated by a subset $S$ of a lattice $L$}{}
\begin{enumerate}
\item Prove that any intersection of a collection of filters is a filter. Conclude that, in particular, for any $S \subseteq L$, there exists a smallest filter containing $S$.  We refer to this filter as \emph{the filter generated by $S$} and denote it by $\gen{S}_{\mathrm{filt}}$.
\item Prove that, for any $S \subseteq L$,
\[\gen{S}_{\mathrm{filt}} = \{b \in L \ | \ \text{there exists finite } T \subseteq S \text{ such that } \bigwedge T \leq b\}.\]
\item Prove that, if $G$ is a filter and $a \in L$, then
\[ \gen{G \cup \{a\}}_{\mathrm{filt}} = \{b \in L \ | \ \text{there exists } g \in G \text{ such that } g \wedge a \leq b\}.\]
Conclude that, if $I$ is a down-set which intersects $\gen{G \cup \{a\}}_{\mathrm{filt}}$ non-trivially, then there exists $g \in G$ such that $g \wedge a \in I$.
\item Formulate and prove analogous results for the \emph{ideal generated} by a
	set $S$, notation $\gen{S}_{\mathrm{idl}}$.\index{ideal!generated}
\nl{$\gen{S}_{\mathrm{idl}}$}{ideal generated by a subset $S$ of a lattice $L$}{}
\end{enumerate}
\end{ourexercise}

\begin{ourexercise}\label{exe:hatlatticehom}
Prove that, for any lattice $L$, the function $\widehat{(-)}$ defined in Theorem~\ref{thm:DLrep} is a well-defined lattice homomorphism.
\end{ourexercise}

\begin{ourexercise}\label{exe:finitesubsetsprimes}
Let $L := \{a \subseteq \mathbb{N} \ | \ a \text{ is finite or } a = \mathbb{N}\}$, the distributive lattice of finite subsets of $\mathbb{N}$ and $\mathbb{N}$ itself.
\begin{enumerate}
\item Show that, for each $n \in \mathbb{N}$, the set $x_n := \{a \in L \ | \ n \in a\}$ is a prime filter of $L$.
\item Show that $x_{\infty} := \{\mathbb{N}\}$ is a prime filter of $L$.
\item Show that $\PrFilt(L) = \{x_n \ | \ n \in \mathbb{N}\} \cup \{x_\infty\}$.
\item Explicitly describe the partial order $\leq$ on $\PrFilt(L)$.
\end{enumerate}
\end{ourexercise}

\begin{ourexercise}\label{exe:fincofsubsetsprimes}
A subset $a$ of a set $X$ is called \emphind{co-finite} (in $X$) if $X \setminus a$ is finite. Let $M := \{a \subseteq \mathbb{N} \ | \ a \text{ is finite or co-finite}\}$, the distributive lattice of finite or co-finite subsets of $\mathbb{N}$, and let $L$ be the lattice defined in the previous exercise, Exercise~\ref{exe:finitesubsetsprimes}.
\begin{enumerate}
\item Prove that, for any lattice homomorphism $h \colon L \to 2$,
    there exists a unique lattice homomorphism $h' \colon M \to 2$
    that extends $h$. \hint{Use
    Exercise~\ref{exe:distuncomp}.\ref{itm:lathomBA}.}
\item Write down an explicit bijection between $\PrFilt(L)$ and $\PrFilt(M)$. \hint{Use item (a) and Lemma~\ref{lem:filtidlcomp}.}
\item Describe the partial order $\leq$ on $\PrFilt(M)$. Which of the two directions of the bijection in the previous item is order preserving?
\end{enumerate}
\end{ourexercise}
\begin{ourexercise}\label{exe:intersection-primefilters}
Deduce from Stone's prime filter-ideal theorem that any filter $F$ in a distributive lattice $L$ is equal to the intersection of all prime filters $G$ that contain $F$.
\end{ourexercise}

\section{Priestley duality}\label{sec:topologize}
In this section, we introduce a topology on the poset $\PrFilt(L)$ that makes it into a special kind of compact ordered space, called \emph{Priestley space}, after H.~A. Priestley, who introduced this topology in \cite{Pri1970}. We will then show that the poset $\PrFilt(L)$, equipped with this topology, completely represents the distributive lattice $L$. Finally, we will also prove that homomorphisms between distributive lattices correspond to continuous order-preserving functions between their dual spaces.

\subsection*{The Priestley dual space of a lattice}
\begin{definition}\label{def:Priestleydual}
Let $L$ be a distributive lattice. The \emph{Priestley topology}, $\tau^p$, on the set $\PrFilt(L)$ is the topology generated by the subbase
\[ \mathcal{S} := \{\widehat{a} \ | \ a \in L\} \cup \{(\widehat{b})^c \ | \ b \in L\}.\]
\end{definition}

This definition makes $\PrFilt(L)$ into a special kind of compact ordered space, as we prove now.
\begin{proposition}\label{prop:priestleycompactord}
Let $L$ be a distributive lattice.
\begin{enumerate}
\item[1.] The topology $\tau^p$ on $\PrFilt(L)$ is compact.
\item[2.] For any $F, G \in \PrFilt(L)$, if $F \nleq G$, then there exists a clopen down-set $U$ in $\PrFilt(L)$ such that $G \in U$ and $F \not\in U$.
\end{enumerate}
In particular, $(\PrFilt(L),\tau^p,\leq)$ is a compact ordered space.
\end{proposition}
\begin{proof}
  1. By the Alexander Subbase Theorem (Theorem~\ref{thm:alexander}), it suffices to prove that any cover $\mathcal{C} \subseteq \mathcal{S}$ of $\PrFilt(L)$ by open sets from the subbase has a finite subcover. Let $S, T \subseteq L$ be such that
\[ \PrFilt(L) = \bigcup_{a \in S} \widehat{a} \cup \bigcup_{b \in T} (\widehat{b})^c.\]
Rephrasing slightly, this equality says that every prime filter $G$ of $L$ that contains the set $T$ must intersect the set $S$ non-trivially. In particular, writing $F$ for the filter generated by $T$ and $I$ for the ideal generated by $S$, every prime filter $G$ of $L$ that contains $F$ must intersect $I$ non-trivially. By Theorem~\ref{thm:DPF}, the intersection of $F$ and $I$ is non-empty, so pick $c \in F \cap I$. By Exercise~\ref{exe:generatefilter}, pick finite $T' \subseteq T$ and finite $S' \subseteq S$ such that $\bigwedge T' \leq c$ and $c \leq \bigvee S'$. In particular, $\bigwedge T' \leq \bigvee S'$, and hence, since $\widehat{(-)}$ is a lattice homomorphism (Theorem~\ref{thm:DLrep}), we have $\bigcap_{b \in T'} \widehat{b} \subseteq \bigcup_{a \in S'} \widehat{a}$. Rewriting this subset inclusion, we conclude that
\[ \PrFilt(L) = \bigcup_{a \in S'} \widehat{a} \cup \bigcup_{b \in T'} (\widehat{b})^c.\]

2. Let $F,G \in \PrFilt(L)$ be such that $F \nleq G$. By definition, there exists $a \in L$ such that $a \in G$ and $a \not\in F$. Thus, $U := \widehat{a}$ is a clopen down-set that contains $G$ but not $F$.

By item 2 and Exercise~\ref{exe:compordaltdef}, $(\PrFilt(L),\leq,\tau^p)$ is in particular a compact ordered space.
\end{proof}

Proposition~\ref{prop:priestleycompactord} gives rise to the following definitions.
\begin{definition}
An ordered topological space $(X,\tau,\leq)$ is called \emph{totally order-disconnected} (TOD) if, for every $x, y \in X$, whenever $x \nleq y$, there exists a clopen down-set $U$ in $X$ such that $y \in U$ and $x \not\in U$.
A compact ordered space which is totally order-disconnected is called a \emph{Priestley space}, sometimes also a CTOD space.
\index{totally order-disconnected}
\index{Priestley space}
\index{TOD}
\index{CTOD}
\end{definition}
The TOD property yields the following stronger version of the order normality property of compact ordered spaces (Exercise~\ref{exer:ordnorm}) for Priestley spaces.\index{order normality}\index{Priestley space!is order-normal}
\begin{proposition}\label{prop:order-normality-Priestley}
    Let $C$ and $D$ be closed sets in a Priestley space $X$. If ${\downarrow} C
    \cap {\uparrow} D = \emptyset$, then there exists a clopen down-set $K
    \subseteq X$ such that $C \subseteq K$ and $K \cap D = \emptyset$.
\end{proposition}
A proof of this proposition is outlined in
Exercise~\ref{exe:closure-normality-Priestley}.

We now examine some Priestley topologies in our running examples.
\begin{example}\label{exa:Nplus1squared-top}
Recall from Example~\ref{exa:Nplus1squared} the lattice $\bN \oplus 1$ and its  poset of prime filters $X := \{F_1, F_2, \dots\} \cup \{F_\omega\}$, see Figure~\ref{fig:nplusone}. For any $n \in \mathbb{N}$, $\widehat{n}$ is the finite subset $\{ F_i \mid 1 \leq i \leq n\}$, and $\widehat{\omega} = X$. It follows that, for every $1 \leq n < \omega$, the singleton $\{F_n\}$ is equal to $\widehat{n} \setminus \widehat{n-1}$, and is therefore clopen in the Priestley topology. Therefore, every subset of $X \setminus \{F_\omega\}$ is open. If, on the other hand $S \subseteq X$ and $F_\omega \in S$, then $S$ is open if, and only if, $S$ is co-finite (see Exercise~\ref{exe:fincofsubsetsprimes}). This topology is homeomorphic to the \emphind{one-point compactification} of the discrete countable space $\{F_1, F_2, \dots\}$, also see Exercise~\ref{exe:fincofsubsetsspace} and Example~\ref{exa:fincofsubsetsdual} below.

As shown in Example~\ref{exa:Nplus1squared}, the lattice $(\bN \oplus 1)^2$ has
as its dual poset the disjoint union of $X$ with itself. For $(n,m) \in (\bN
\oplus 1)^2$, the clopen down-set $\widehat{(n,m)}$ is equal to
\[ \{F_{(n',0)} \mid n' \leq n\} \cup \{F_{(0,m')} \mid m' \leq
m\}.\] 
It follows that the topology is also the disjoint union of the topologies on the two copies of $X$, that is, a subset $U \subseteq X + X$ is open if, and only if, its intersection with both the `left' and `right' copy of $X$ is open (see Exercise~\ref{exe:sumproduct}).
\end{example}

\begin{example}\label{exa:NplusNop-top}
  We continue Example~\ref{exa:NplusNop} of the lattice $\bN \oplus \bN^\op$, whose dual poset was depicted in Figure~\ref{fig:nplusnop}. For $n \in \mathbb{N}$, we have that
  \[ \widehat{(n,0)} = \{ F_{(k,0)} \mid k \leq n\}, \] and 
  \[ \widehat{(n,1)} = \{ F_{(k,1)} \mid k \geq n\} \cup \{ F_\omega \} \cup \{ F_{(m,0)} \mid m \in \mathbb{N}\}.\] 
  It follows that, for $n \geq 1$,
  \[\{F_{(n,0)}\} = \widehat{(n,0)} \setminus \widehat{(n+1,0)},\] and, for $n \geq 0$,
  \[\{F_{(n,1)}\} = \widehat{(n,1)} \setminus \widehat{(n-1,1)},\] so that again all the singleton sets except for $\{F_{\omega}\}$ are open in the Priestley topology. As in the previous example, one may again prove that the topology is homeomorphic to the one-point compactification of the discrete countable space $\{F_{(k,i)} \mid k \in \mathbb{N}, i \in \{0,1\}\}$. The only thing distinguishing this space from the space $X$ in the previous example is the partial order. %
\end{example}

\begin{example}\label{exa:primenumbers-topology}
  We continue Example~\ref{exa:prime-divisibility} of the lattice of divisibility $(\mathbb{N}, |)$, whose dual poset was depicted in Figure~\ref{fig:prime-naturalnumbers}. Clearly, $\widehat{0}$ is the entire space, and for any $n \geq 1$, we have 
  \[ \widehat{n} = \{ F_{p^k} \ : \ p^k \mid n, p \text{ prime}, k \geq 1 \}.\]
  In particular, for any prime $p$ and $k \geq 1$, we have that
  \[ \{F_{p^k}\} = \widehat{p^k} \setminus \widehat{p^{k-1}},\]
  so all singletons except $\{F_0\}$ are clopen.  Since $\widehat{n}$ is finite for every $n \geq 1$, it follows again that a clopen set $K$ in the Priestley topology contains $F_0$ if, and only if, it is co-finite, and the topology on the dual poset is again the one-point-compactification of a discrete countable space.
\end{example}

\begin{example}\label{exa:rationalintervals-topology}
We continue Example~\ref{exa:rationalintervals-primes} of the lattice of finite unions of rational intervals, whose dual poset was depicted in Figure~\ref{fig:cherries}. As a useful notation for the rest of this example (also see Notation~\ref{not:neutral} below), for every $r \in [0,1]$, we denote the prime filter $F_r$ defined in Example~\ref{exa:rationalintervals-primes} simply by $r$, and we write $r^+$ and $r^-$ for the prime filters that were denoted by $F_r^+$ and $F_r^-$, respectively, whenever they exist. Also, for every $r \in [0,1]$, let us write $X_r$ for the set of prime filters `located' at $r$, that is,
\begin{align*}
  X_r := \left\{ \begin{array}{ll}
    \{r\} &\text{ if $r$ is irrational,} \\
    \{r, r^+, r^-\} &\text{ if } r \in (0,1) \cap \mathbb{Q}, \\
    \{0, 0^+\} &\text{ if } r = 0, \\
    \{1, 1^-\} &\text{ if } r = 1.
  \end{array} \right.
\end{align*}
Now, towards describing the Priestley topology on the poset $X = \bigcup_{r \in [0,1]} X_r$, note that, for any $q \in \mathbb{Q} \cap [0,1]$,
\begin{align*}
  \widehat{[0,q)} &= \bigcup_{r < q} X_r \cup \{ q^- \},\\
  \widehat{(q,1]} &= \{q^+\} \cup \bigcup_{q < r} X_r,
\end{align*}
and, hence,
\begin{align*}
  \widehat{[0,q)}^c &= \{q\} \cup \widehat{(q,1]},\\
  \widehat{(q,1]}^c &= \widehat{[0,q)} \cup \{q\}.
\end{align*}
It follows that, in the Priestley topology, for any $q \in \mathbb{Q} \cap [0,1]$, the singleton $\{q\}$ is clopen, since it is equal to $\widehat{[0,q)}^c \cap \widehat{(q,1]}^c$. In fact, the collection of clopen sets in the Priestley topology is generated as a sublattice of $\mathcal{P}(X)$ by the sets $\widehat{[0,q)}$, $\widehat{(q,1]}$, and $\{q\}$, for $q \in \mathbb{Q} \cap [0,1]$. To see this, note that the complement of any of the generating sets is equal to a finite union of other generating sets: for $\widehat{[0,q)}$ and $\widehat{(q,1]}$ this was already shown above, and for $\{q\}$, note that
\[ \{q\}^c = \widehat{[0,q)} \cup \widehat{(q,1]}.\]
Now note also that the sublattice generated by the sets of these three forms is equal to the collection of finite unions of sets $\widehat{(p,q)}$, a finite set of rationals and $\widehat{[0,a)}$ and $\widehat{(b,1]}$. For an open subset $U$ in the Euclidean topology on $[0,1]$, write
\[\widetilde{U} := \bigcup \{ \widehat{I} \mid I \text{ a rational open interval in } [0,1] \text{ such that } I \subseteq U\},\] 
a Priestley open down-set of $X$. It then follows that any open set of the
Priestley topology can be written as a union $\widetilde{U} \cup V$, where $U$ is an open subset of $[0,1]$ in the Euclidean topology, and $V$ is a subset of $\mathbb{Q} \cap [0,1]$. It is instructive to give a direct proof that this is indeed a compact topology on $X$; such a proof is outlined in Exercise~\ref{exe:compact-rational-intervals}.
\end{example}

In Definition~\ref{def:Priestleydual} above we consider a topology on the set of prime filters of a lattice $L$. By Lemma~\ref{lem:filtidlcomp}, this set is in bijection with the set of prime ideals and with the set of homomorphisms of the lattice into $\btwo$. As a consequence, this topology may be translated to any of these three sets. Indeed, in the literature, depending on the time period, or the applications to be treated, the dual space of a distributive lattice may be based on any of these three sets of points. As we will see later in this book, for some applications, it is useful to be able to switch flexibly between these choices. This is the motivation for working with what we will call the \emphind{neutral dual space}, which has a `neutral' underlying set of points that comes equipped with named bijections to the prime filters, the prime ideals, and the homomorphisms into $\btwo$ that are connected as in Lemma~\ref{lem:filtidlcomp}.

\begin{notationnum}[The neutral dual space]\label{not:neutral}
\nl{$X_L$}{the (neutral) dual space of a distributive lattice $L$}{}
\nl{$F_x$}{prime filter corresponding to a point $x$ of a (neutral) dual space}{}
\nl{$I_x$}{prime ideal corresponding to a point $x$ of a (neutral) dual space}{}
\nl{$h_x$}{homomorphism into $2$ corresponding to a point $x$ of a (neutral) dual space}{}
Let $L$ be a distributive lattice. We fix a set $X_L$ of the same cardinality as $\PrFilt(L)$, $\PrIdl(L)$ and $\Hom_\DL(L,2)$, and we fix three bijections:
\begin{align*}
F_{(-)} &\colon X_L \to \PrFilt(L),  \\
I_{(-)} &\colon X_L \to \PrIdl(L), \text{ and} \\
h_{(-)} &\colon X_L \to \Hom_\DL(L,2)
\end{align*}
so that, for all $x\in X_L$ and for all $a\in L$, we have
\[
a\in F_x \quad\iff\quad a\not\in I_x \quad\iff\quad h_x(a)=\top.
\]
Thus, to any element $x \in X_L$ there corresponds a unique prime filter $F_x$ of $L$, a unique prime ideal $I_x$ of $L$, which is the set complement of $F_x$ in $L$, and a unique homomorphism $h_x \colon L \to 2$, which is the characteristic function of $F_x$. Elements of the set $X_L$ are called \emphind{points} of the distributive lattice $L$, and are usually denoted by lowercase letters towards the end of the alphabet.

We define a \emph{partial order} $\leq$ on $X_L$  by
\begin{align*}
x \leq y &\text{ if, and only if, } F_x \supseteq F_y, \\
		&\text{ if, and only if, } I_x \subseteq I_y, \\
		&\text{ if, and only if, } h_x \geq h_y \text{ pointwise}.
\end{align*}

For any $a \in L$, the set $\widehat{a}$ is defined by
\begin{align*}
\widehat{a} &= \{ x \in X \ | \ a \in F_x \} \\
			&= \{ x \in X \ | \ a \not\in I_x \} \\
			&= \{ x \in X \ | \ h_x(a) = \top\}.
\end{align*}
The \emph{Priestley topology}, $\tau^p$, on $X_L$ is the topology generated by the subbase $\{\widehat{a} \ | \ a \in L\} \cup \{\widehat{b}^c \ | \ b \in L\}$.
\end{notationnum}

\begin{definition}\label{def:Priestleydualspace}
Let $L$ be a distributive lattice. The tuple $(X_L, \leq, \tau^p)$ is called the \emph{Priestley dual space} of $L$.
\end{definition}
\begin{remark}[On the definition of the order on the dual space]\label{rem:order-choice}
Our definition of the order on the Priestley dual space is aligned with the
point of view on prime filters as idealized join-prime elements of the lattice;
also see our remarks preceding Definition~\ref{def:filterorder} above. However,
in some literature, including the original paper~\parencite{Pri1970}, the
ordering on the Priestley dual space of a lattice is the opposite of the one we
give here. While this is ultimately a matter of convention, it is also a
potential source of confusion, that the reader should be aware of, and we now
explain in some detail why the two choices naturally exist.

Indeed, since for any distributive lattice $L$, the reverse lattice $L^\op$ is also a distributive lattice, and the poset of prime filters of $L^\op$ is the opposite poset of the prime filters of $L$, both $(X_L, \leq, \tau^p)$ and $(X_L, \geq, \tau^p)$ are Priestley spaces. To build a dual equivalence between distributive lattices and ordered topological spaces, one needs to associate one of these two spaces to $L$ and the other to $L^\op$, but there is a free choice as to which space is associated to which lattice. With the choice we make here, $L$ is isomorphic to the clopen down-sets of its Priestley dual space; with the other choice, $L$ will be isomorphic to the clopen up-sets. 

While our choice is more natural relative to our view on prime filters and in comparison with the case of finite distributive lattices, where the dual is a subset of the lattice and is thus already equipped with a partial order, there are settings where putting the opposite order on the dual space can be a more natural choice. This is in particular the case in the more categorical point of view on the points of the dual space as the set of homomorphisms into $2$, and also relative to the convention for the specialization order in topology, where open sets are up-sets.

In this book, we mostly focus on the choice given in the above definition, and we notify the reader when we need to diverge from this choice; this will in particular be the case in Chapters~\ref{chap:Omega-Pt} and \ref{chap:DomThry}, in order to make the link with existing literature on frames and domains, where the other choice is more prevalent; see also Remark~\ref{rem:order-yoga}.
\end{remark}

\subsection*{Priestley representation}
We can use the Priestley topology to exactly identify which down-sets lie in the image of the map $\widehat{(-)}$ from the representation theorem for distributive lattices, Theorem~\ref{thm:DLrep}.
\begin{proposition}\label{prop:hatimage}
Let $L$ be a distributive lattice. The image of $\widehat{(-)} \colon L \to \Down(X_L)$ consists exactly of those down-sets in $X_L$ that are clopen in the Priestley topology.
\end{proposition}
\begin{proof}
It is clear from the definition of the Priestley topology that $\widehat{a}$ is clopen for every $a \in L$. Conversely, let $A$ be an arbitrary clopen down-set of $X_L$. For every $y \in A$ and $x \not\in A$, since $A$ is a down-set, we have that $x \not\leq y$, so, by definition of $\leq$, pick an element $b_{x,y} \in L$ such that $b_{x,y} \in F_y$ and $b_{x,y} \not\in F_x$. In particular, for every $y \in A$, the collection $\{(\widehat{b_{x,y}})^c \ | \ x \not\in A\}$ is a cover of $A^c$. Since $A$ is open, its complement, $A^c$, is closed, and hence compact. For every $y \in A$, pick a finite subcover $\{(\widehat{b_{x_i,y}})^c \ | \ i = 1, \dots, n\}$ of $A^c$. Then $\bigcap_{i=1}^n \widehat{b_{x_i,y}} \subseteq A$, so, defining $a_y := \bigwedge_{i=1}^n b_{x_i,y}$, we have $y \in \widehat{a_y} \subseteq A$. Doing this for every $y \in A$, we obtain a cover $\{\widehat{a_y} \ | \ y \in A\}$ of $A$, which is closed, and hence compact. Pick a finite subcover $\{\widehat{a_{y_i}} \ | \ j = 1,\dots,m \}$ of $A$. Defining $a := \bigvee_{j=1}^n a_{y_j}$, we have $A = \widehat{a}$, as required.
\end{proof}

Since $\widehat{(-)}$ was already shown to be a lattice embedding in Theorem~\ref{thm:DLrep}, we obtain the following representation result for distributive lattices.
\begin{corollary}\label{cor:doubledualDL}
Any distributive lattice is isomorphic to the lattice of clopen down-sets of its Priestley dual space.
\end{corollary}
For a Priestley space $X$, we denote by $\ClD(X)$ its lattice of clopen down-sets. Proposition~\ref{prop:hatimage} may then be rephrased as: for any distributive lattice $L$, the map $\widehat{(-)}$ is an isomorphism from $L$ to $\ClD(X)$. This is a generalization of Proposition~\ref{prop:birkhoff} to arbitrary distributive lattices. Indeed, when $L$ is a finite distributive lattice then its Priestley dual space is order-homeomorphic to the finite poset $\cJ(L)$ equipped with the discrete topology, so the clopen down-sets are just the down-sets, and Corollary~\ref{cor:doubledualDL} then gives $L \cong \Down(\cJ(L))$.

Analogously to Proposition~\ref{prop:finite-poset-double-dual} for finite posets, we will now also prove that every Priestley space is order-homeomorphic to its double dual. A basic idea in topology\footnote{This same idea will again be central, in the setting of sober spaces, in Chapter~\ref{chap:Omega-Pt}, see Definition~\ref{def:neighborhood-filter}.} is that for any $T_0$ topological space $(X,\rho)$, the set $\mathcal{N}(x)$ of open neighborhoods of a point $x$ can be used to uniquely pinpoint $x$ among all the points in $X$. We now apply this idea in the setting of Section~\ref{sec:comp-ord-sp}, where $\rho = \tau^{\downarrow}$, the topology of open down-sets of some ordered topological space $(X,\tau,\leq)$. We will show that, for a Priestley space, it actually suffices to consider the intersection of $\mathcal{N}(x)$ with the \emph{clopen} down-sets of $(X,\tau,\leq)$. A crucial step for proving this is the following.
An important connection between Priestley spaces and prime filters of clopen down-sets, that we will use below, is the following.
\begin{lemma}\label{lem:intersection-prime-principal}
  Let $(X, \tau, \leq)$ be a Priestley space. For any prime filter $\mathcal{F}$ of the lattice $\ClD(X)$ of clopen down-sets of $X$, the intersection $\bigcap \mathcal{F}$ is a principal down-set in $X$, that is, there exists a point $x \in X$ such that $\bigcap \mathcal{F} = {\downarrow}x$.
\end{lemma}
\begin{proof}
  Let $\cF$ be a prime filter of $\ClD(X)$, and write $\cI$ for the complement of $\cF$, which is a prime ideal of $\ClD(X)$. We first show that the collection of clopen sets
  \[ \cC := \cF \cup \{L^c \ \mid \ L \in \cI\} \] 
  has the finite intersection property.
  Indeed, if for any finite subsets $\cS \subseteq \cF$ and $\cT \subseteq \cI$, we would have $\bigcap \cS$ disjoint from $\bigcap \{ L^c \ \mid \ L \in \cT \}$, then this would give $\bigcap \cS \subseteq \bigcup \cT$, which is clearly impossible since $\cF$ and $\cI$ are disjoint. 
  Since the space $X$ is compact, we can therefore (see Exercise~\ref{exer:covers}.\ref{itm:fin-int-prop}) pick a point $x \in \bigcap \cC$.
  Then in particular $x \in \bigcap \cF$, and 
  since $\bigcap \cF$ is a down-set, we have 
  ${\downarrow} x \subseteq \bigcap \cF$. 
  It remains to prove that $\bigcap \cF \subseteq {\downarrow} x$.
  Let $y \in \bigcap \cF$. To show that $y \leq x$, by (the contrapositive of) total order disconnectedness, it suffices to prove that, for any $L \in \ClD(X)$, if $x \in L$ then $y \in L$.
  Let $L \in \ClD(X)$. If $x \in L$, then, since 
  $x \in \bigcap \{ L^c \ \mid \ L \in \cI\}$, 
  we have that $L \not\in \cI$. It follows that $L \in \cF$, and thus $y \in L$. 
\end{proof}

We now get the following `double dual' result for Priestley spaces.
\begin{proposition}\label{prop:priestley-double-dual}
  Let $(X,\tau,\leq)$ be an ordered topological space. Write $L$ for the lattice of clopen down-sets of $X$, and, for any $x \in X$,
  define
  \[ \mathcal{B}(x) := \{ K \in L \ | \ x \in K \} \ .\]
  We have the following properties.
  \begin{enumerate}
    \item The set $\mathcal{B}(x)$ is a prime filter of $L$ for every $x \in X$.
    \item The function $\beta \colon X \to X_L$, defined by sending a point $x \in X$ to the unique point $\beta(x) \in X_L$ with $F_{\beta(x)} = \mathcal{B}(x)$, is a continuous order-preserving function, satisfying $\beta^{-1}(\widehat{K}) = K$ for any $K \in L$.
    \item \label{itm:beta-order-homeo} The function $\beta$ is an order homeomorphism if, and only if, $(X,\tau,\leq)$ is a Priestley space.
  \end{enumerate}
\end{proposition}
\begin{proof}
  a. Notice that, $X_L\in\mathcal{B}(x)$ and for any $K_1, K_2 \in L$, we have $K_1 \cap K_2 \in \mathcal{B}(x)$ if, and only if, both $K_1$ and $K_2$ are in $\mathcal{B}(x)$, so $\mathcal{B}(x)$ is a filter. Clearly $\emptyset \not\in \mathcal{B}(x)$, and if $K_1 \cup K_2 \in \mathcal{B}(x)$, then $x \in K_1 \cup K_2$, so $x \in K_i$ for some $i \in \{1,2\}$. This shows that $\mathcal{B}(x)$ is prime.

  b. Note first that, for any $K \in L$ and $x \in X$ we have
  \[ \beta(x) \in \widehat{K} \iff K \in \mathcal{B}(x) \iff x \in K\ . \]
  In other words, $\beta^{-1}(\widehat{K}) = K$. Thus, by the definition of the Priestley topology on $X_L$, $\beta$ is a continuous function. 
  We now show that $\beta$ is order preserving. Let $x, y \in X$ with $x \leq y$. Then, if $K$ is a clopen down-set that contains $y$, then it also contains $x$. Thus, $\mathcal{B}(y) \subseteq \mathcal{B}(x)$. In light of Notation~\ref{not:neutral}, this means that $\beta(x) \leq \beta(y)$. 
  
  c. By Proposition~\ref{prop:priestleycompactord}, $X_L$ is a Priestley space for any distributive lattice $L$, so if $\beta$ is an order homeomorphism then $(X,\tau,\leq)$ is also a Priestley space. For the other implication, assume $(X,\tau,\leq)$ is a Priestley space. The total order-disconnected axiom says precisely that, for any $x, y \in X$, if $x \nleq y$, then $\mathcal{B}(y) \not\subseteq \mathcal{B}(x)$, so that $\beta(x) \nleq \beta(y)$. Since we saw in the previous item that $\beta$ is always continuous and order preserving, it follows that $\beta$ is a continuous order embedding. Now, by the previous item, for any $K \in L$ we have $\beta[K] = \widehat{K}$, so $\beta$ is open.
  It remains to prove that $\beta$ is surjective. Let $\chi \in X_L$ be arbitrary. By Lemma~\ref{lem:intersection-prime-principal}, pick $x \in X$ such that $\bigcap F_{\chi} = {\downarrow} x$. To finish the proof, we show that $\beta(x) = \chi$, that is, $\cB(x) = F_{\chi}$. Indeed, for any $K \in L$, we have 
  \[ K \in \cB(x) \iff x \in K \iff \bigcap F_{\chi} \subseteq K \iff K \in F_{\chi} ,\]
  where the last equivalence follows from the compactness of $X$, using Exercise~\ref{exer:covers}. 
\end{proof}

In the next subsection, we show how this object correspondence between distributive lattices and Priestley spaces extends to homomorphisms between lattices.%

\subsection*{Priestley duality for homomorphisms}
Let $L$ and $M$ be distributive lattices, and denote by $X$ and $Y$ the Priestley dual spaces of $L$ and $M$, respectively. First, let $f \colon X \to Y$ be an order-preserving function. We have a complete homomorphism $f^{-1} \colon \Down(Y) \to \Down(X)$, which sends any down-set $D$ of $Y$ to the down-set $f^{-1}(D)$ of $X$. Now assume that $f$ is moreover continuous. This means that the inverse image under $f$ of any clopen down-set in $Y$ is a clopen down-set in $X$. In other words, by Proposition~\ref{prop:hatimage}, $f^{-1}$ restricts correctly to a function between the images of the representation maps $\widehat{(-)}$ for $M$ and $L$. For any $b \in M$, let $h_f(b)$ denote the unique element of $L$ such that $f^{-1}(\widehat{b}) = \widehat{h_f(b)}$; that is, $h_f \colon M \to L$ is the restriction of $f^{-1}$ to $M$ and $L$, identifying them with their image under $\widehat{(-)}$ in the down-set lattices. In a diagram:
\begin{center}
\begin{tikzpicture}
  \matrix (m) [matrix of math nodes,row sep=3em,column sep=3em,minimum width=2em]
  {
     \Down(Y) & \Down(X) \\
     M & L \\};
  \path[-stealth]
    (m-1-1) edge node [above] {$f^{-1}$} (m-1-2)
    (m-2-1) edge node [above] {$h_f$} (m-2-2)
    (m-2-1) edge node [left] {$\widehat{(-)}$} (m-1-1)
    (m-2-2) edge node [right] {$\widehat{(-)}$} (m-1-2);
\end{tikzpicture}
\end{center}
Note that $h_f$ is a lattice homomorphism, because it is essentially the restriction of the complete lattice homomorphism $f^{-1}$. Thus we have shown that for any continuous and order-preserving $f \colon X \to Y$, there is a homomorphism $h_f \colon M \to L$ making the above diagram commute. In the following proof we will use the description of $h_f $ in terms of homomophisms into $2$. For this purpose, recall that, for any point $x \in X$, $h_x$ denotes the homomorphism $L \to 2$ naturally associated to $x$. Now note that, for any $x \in X$, $b \in M$,
\[ h_{f(x)}(b) = \top \iff f(x) \in \widehat{b} \iff x \in \widehat{h_f(b)} \iff h_x(h_f(b)) = \top,\]
that is, the homomorphism $h_{f(x)}$ associated to the point $f(x) \in Y$ is equal to the composition $h_x \circ h_f$.

Crucially, every lattice homomorphism from $M$ to $L$ arises as $h_f$ for some continuous order-preserving $f$, as we will now prove.\index{homomorphism!dual function}

\begin{proposition}\label{prop:priestleymorphisms}
Let $L$ and $M$ be distributive lattices with Priestley dual spaces $X$ and $Y$, respectively. For any lattice homomorphism $h \colon M \to L$, there exists a unique continuous order-preserving function $f \colon X \to Y$ such that $h = h_f$.
\end{proposition}
\begin{proof}
Let $h \colon M \to L$ be a lattice homomorphism.
Define $f \colon X \to Y$ to be the function which associates to each element $x \in X$ the element $f(x) \in Y$ corresponding to the homomorphism $h_x \circ h$. By the above remarks, this is the only function $f$ such that $h = h_f$. In particular, for any $a \in M$, we have $f^{-1}(\widehat{a}) = \widehat{h(a)}.$
Hence, $f$ is continuous, by the definition of the Priestley topologies on $X$ and $Y$. Moreover, if $x \leq x'$, then $h_x \geq h_{x'}$ pointwise, so that also $h_{f(x)} = h_x \circ h \geq h_{x'} \circ h = h_{f(x')}$ pointwise; thus, $f(x) \leq f(x')$.
\end{proof}

Note the similarity of Proposition~\ref{prop:priestleymorphisms} to
Proposition~\ref{prop:finDLmorphisms} in the finite Birkhoff duality. Indeed,
Proposition~\ref{prop:finDLmorphisms} is a special case of
Proposition~\ref{prop:priestleymorphisms}, if one takes into account the
correspondence between prime filters and join-irreducibles in the finite case
(see Exercise~\ref{exe:joinprime-primefilt}).

Exercise~\ref{exe:morphisms-concrete} gives an alternative concrete
way of defining the unique map $f$ for which $h = h_f$, by working with prime
filters or prime ideals. You are asked to prove there that, in the situation of
Proposition~\ref{prop:priestleymorphisms}, for any $x \in X$, the prime filter
corresponding to $f(x)$ is the inverse image under $h$ of the prime filter
corresponding to $x$; that is, $F_{f(x)} = h^{-1}(F_x)$. Similarly, $I_{f(x)} =
h^{-1}(I_x)$.

The next proposition highlights two important features of the morphism
correspondence given in Proposition~\ref{prop:priestleymorphisms}. 
\begin{proposition}\label{prop:surj-inj}
    Let $L$ and $M$ be distributive lattices with Priestley dual spaces $X$ and
    $Y$, respectively, and let $h \colon M \to L$ be a lattice homomorphism
    with dual continuous order-preserving function $f \colon X \to Y$.
    \begin{enumerate}
        \item The function $h$ is injective if, and only if, $f$ is surjective.
        \item The function $h$ is surjective if, and only if, $f$ is an order-embedding.
    \end{enumerate}
\end{proposition}
\begin{proof}
The first item is left as Exercise~\ref{exe:Priestley-injsurj}. For the second
item, suppose first that $f$ is an order-embedding. Let $a \in L$. Consider the
subsets $C := f[\widehat{a}]$ and $D := f[\widehat{a}^c]$ of $Y$, which are
both closed because $f$ is a closed map (Exercise~\ref{exer:closedmap}). We
show first that ${\downarrow} C \cap {\uparrow} D = \emptyset$. If there would exist a
point $y$ in this intersection, then there would exist $x \in \widehat{a}$ and
$x' \in \widehat{a}^c$ such that $f(x') \leq y \leq f(x)$. Since $f$ is an
order-embedding, we then get $x' \leq x$, which is impossible because
$\widehat{a}$ is a down-set. Thus, ${\downarrow} C \cap {\uparrow} D =
\emptyset$, and Proposition~\ref{prop:order-normality-Priestley} gives an
element $b \in M$ such that $C \subseteq \widehat{b}$ and $\widehat{b} \cap D
= \emptyset$. We show now that $h(b) = a$. On the one hand, $C \subseteq \widehat{b}$ gives
that $\widehat{a} \subseteq f^{-1}(\widehat{b}) = \widehat{h(b)}$, so $a \leq
h(b)$. On the other hand, if $x \in \widehat{h(b)}$, then $f(x) \in
\widehat{b}$, so $f(x) \not\in D$. In particular, $x \not\in \widehat{a}^c$, so
$x \in \widehat{a}$. Thus, $\widehat{h(b)} \subseteq \widehat{a}$, so that
$h(b) \leq a$, as required.

Conversely, suppose that $h$ is surjective. Let $x, x' \in X$ be such that
$x \nleq x'$. Pick $a \in L$ such that $x' \in \widehat{a}$ and $x \not\in
\widehat{a}$. Since $h$ is surjective, pick $b \in M$ such that $h(b) = a$.
Then $x' \in \widehat{h(b)} = f^{-1}(\widehat{b})$, so $f(x') \in \widehat{b}$,
but by a similar argument, $f(x) \not\in \widehat{b}$. We conclude that $f(x)
\nleq f(x')$ because $\widehat{b}$ is a down-set.
\end{proof}

Note that by Exercise~\ref{exer:closedmap} any continuous function between 
Priestley spaces is a closed mapping, and thus by Exercise~\ref{exer:quotspace} 
any surjective continuous function between Priestley spaces is a quotient map of 
the underlying topological spaces. Thus item (a) of the above proposition says
we have a one-to-one correspondence between sublattices of a distributive lattice
and Priestley quotient spaces, that is, the quotient spaces of the underlying 
topological space equipped with any Priestley order which makes the quotient map
order preserving. And item (b) tells us that we have a one-to-one correspondence 
between quotient lattices (and thus lattice congruences) and Priestley subspaces of 
the dual Priestley space.

We will refine this result in Section~\ref{sec:quotients-and-subs}, which gives
concrete methods for computing quotients and subs on either side of the
duality. Here, we just show how the two items in  Proposition~\ref{prop:surj-inj}, taken 
together, in fact can be used to understand Priestley duality for arbitrary lattice homomorphisms. 
Recall (p.~\pageref{firstiso}) that the first isomorphism theorem for lattices says that any
homomorphism $h \colon M \to L$ can be factored as a quotient map $p
\colon M \onto M/{\theta}$ followed by a lattice embedding $i \colon
M/{\theta} \into L$, where $\theta = \ker(h)$ is the kernel of $h$. Writing $Z$
for the Priestley dual space of the lattice $M/{\theta}$, 
Proposition~\ref{prop:surj-inj} shows that the quotient map $p$ is dual to an embedding
$j \colon Z \into Y$, and that the lattice embedding $i \colon M/{\theta}
\into L$ is dual to a quotient of Priestley spaces $q \colon X \onto Z$.
These two maps give a factorization of the continuous order-preserving function
$f \colon X \to Y$ dual to $h$. In summary, the following two triangles are
dual to each other.
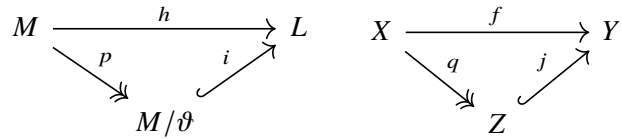
\begin{center}
\begin{tikzcd}
    M \arrow[rr,"h"] \arrow[rd, "p", two heads] & & L \\
     & M/{\theta} \arrow[ru,"i", hook] &
\end{tikzcd}
\quad
\begin{tikzcd}
    X \arrow[rr,"f"] \arrow[rd, "q", two heads] & & Y \\
     & Z \arrow[ru,"j",hook] &
\end{tikzcd}

\end{center}
We give two examples of homomorphisms and their duals, using the examples of Priestley spaces that we developed in this chapter.
\begin{example}
  Recall the lattices $\bN \oplus 1$ and $\bN \oplus \bN^\op$, whose Priestley dual spaces where described in Examples~\ref{exa:Nplus1squared-top}~and~\ref{exa:NplusNop-top}.
  Consider the (injective) lattice homomorphism $h \colon \bN \oplus 1 \to \bN \oplus \bN^\op$ which sends each $n \in \bN$ to $(n,0)$, and $\omega$ to $(0,1)$. %
  For  a prime filter $F_x$ in the dual space of $\bN \oplus \bN^\op$, the prime filter $F_{f(x)}$ is the inverse image under $h$ of $F_x$ (see Exercise~\ref{exe:morphisms-concrete}). Thus, the continuous order-preserving function $f$ from the dual space of $\bN \oplus \bN^\op$ to the dual space of $\bN \oplus 1$ sends every prime filter $F_{(n,0)}$ to $F_n$, for $n \geq 1$, and sends $F_\omega$, as well as any prime filter above it, to the prime filter $F_\omega$ of the dual space of $\bN \oplus 1$. 
\end{example}

\begin{example}\label{exa:primenumbers-maps}
Consider the lattice $(\bN, |)$, whose Priestley dual space was described in Example~\ref{exa:primenumbers-topology}. A natural number $n$ is called \emph{square-free} if it is a product of distinct primes, that is, for any prime $p$ and $k \geq 1$, if $p^k \mid n$ then $k = 1$. Note that the subset $S$ of square-free or zero numbers is a sublattice of $(\bN, |)$. A prime filter of $S$ is either $\{0\}$, or of the form $F_p \cap S$ for some prime $p$. Denote by $h$ the inclusion homomorphism $S \to \bN$. Note that $F_{p^k} \cap S = \{0\}$ for any $k \geq 2$. Therefore, the Priestley dual of $h$ is the function $f$ defined as follows: $f$ sends $F_0$ to $F_0$, for each prime number $p$, $f$ sends the prime filter $F_{p}$ of $\bN$ to the prime filter $F_p \cap S$ of $S$, and for every $k \geq 2$, $f$ sends the prime filter $F_{p^k}$ of $\bN$ to $F_0$; see Figure~\ref{fig:primenumbers-maps}.
\begin{figure}[htp]
  \begin{center}
  \begin{tikzpicture}
    \polab{(-2,1)}{$F_4$}{left}
    \polab{(-2,0)}{$F_2$}{below}
    \draw (-2,0) -- (-2,1.5);
    \node at (-2, 2) {$\vdots$};
  
    \polab{(0,0)}{$F_3$}{below}
    \polab{(0,1)}{$F_9$}{left}
    \draw (0,0) -- (0,1.5);
    \node at (0, 2) {$\vdots$};
    
    \polab{(2,0)}{$F_5$}{below}
    \polab{(2,1)}{$F_{25}$}{left}
    \draw (2,0) -- (2,1.5);
    \node at (2, 2) {$\vdots$};
  
    \polab{(0,3)}{$F_0$}{left}
    \node at (2.5,0.5) {$\dots$};

    \node at (3.5,1.5) {$\onto$};

    \polab{(5,0)}{$F_2 \cap S$}{below}
    \polab{(6.5,0)}{$F_3 \cap S$}{below}
    \polab{(8,0)}{$F_5 \cap S$}{below}
    \polab{(6.5,3)}{$F_0$}{right}
    \node at (9,0) {$\dots$};

    \draw[dotted,-stealth] (-2,0) to[bend right=25] (5,0);
    \draw[dotted,-stealth] (-2,1) to[bend left=20] (6.5,3);
    \draw[dotted,-stealth] (0,0) to[bend right=25] (6.5,0);
    \draw[dotted,-stealth] (0,1) to[bend left=20] (6.5,3);
    \draw[dotted,-stealth] (2,0) to[bend right=25] (8,0);
    \draw[dotted,-stealth] (2,1) to[bend left=20] (6.5,3);
    \draw[dotted,-stealth] (0,3) to[bend left=25] (6.5,3);
    \draw (5,0)--(6.5,3)--(6.5,0);
    \draw (8,0)--(6.5,3);
  \end{tikzpicture}
  \end{center}
  \caption{The dual of the homomorphism $h \colon S \into (\bN, |)$.}
  \label{fig:primenumbers-maps}
  \end{figure}

Also consider the function $g \colon \bN \to S$, which sends any $n \geq 0$ to its largest square-free divisor, and $0$ to $0$; concretely, if $n = p_1^{k_1} \cdots p_m^{k_m}$ for distinct primes $p_1, \dots, p_m$ and $k_1, \dots, k_m \geq 1$, then $g(n) := p_1 \cdots p_m$. Note that $g$ is upper adjoint to $h$: for any $s \in S$ and $n \in \bN$, we have that $h(s) \mid n$ if, and only if, $s \mid g(n)$. In particular, $g$ is meet-preserving. Note that $g$ is in fact a homomorphism: clearly $g(1) = 1$, and for any $n, n' \in \bN$, a prime $p$ divides $g(n \vee n')$ if, and only if, $p \mid n$ or $p \mid n'$, if, and only if, $p$ divides $g(n) \vee g(n')$; hence, $g(n \vee n') = g(n) \vee g(n')$. The dual of $g$ is the function which sends each prime filter $F_p \cap S$ of $S$ to the prime filter $F_p$ of $\bN$, and $\{0\}$ to $\{0\}$. The functions dual to $h$ and $g$ are also an adjoint pair between the Priestley posets of prime filters; we will come back to this example in Exercise~\ref{exe:prime-bifinite}, after we discuss `bifinite domains' in Section~\ref{sec:bifinite}.
\end{example}

As for the finite case in Chapter~\ref{ch:order}, we now state Priestley's duality theorem without having formally defined what a dual equivalence is. We refer to Chapter~\ref{ch:categories} for the precise definitions of the terms that are used, and the reader will then be able to deduce this theorem for themselves in Exercise~\ref{exe:priestley-is-duality}, also see the details given in Theorem~\ref{thm:priestley-duality-long}. For now, it suffices that the reader understands Corollary~\ref{cor:doubledualDL} and Proposition~\ref{prop:priestleymorphisms}, which, we will see, form the content of the theorem.
 \index{Priestley duality}\index{duality!Priestley}
\begin{theorem}\label{thm:priestleyduality}
The category $\DL$ of bounded distributive lattices with homomorphisms is dually equivalent to the category $\cat{Priestley}$ of Priestley spaces with continuous order-preserving functions.
\end{theorem}

\ourexercises

\begin{ourexercise}\label{exe:finitesubsetsspace}
  Let $L$ be the lattice of finite subsets of $\mathbb{N}$ and $\mathbb{N}$ itself, as defined in Exercise~\ref{exe:finitesubsetsprimes}. Recall from Exercise~\ref{exe:finitesubsetsprimes} the description of $\PrFilt(L)$ as $\{x_n \ | \ n \in \mathbb{N}\} \cup \{x_\infty\}$. Consider the Priestley topology on $\PrFilt(L)$.
  \begin{enumerate}
  \item Prove that, for every $n \in \mathbb{N}$, $\{x_n\}$ is clopen.
  \item Prove that a subset $K \subseteq \PrFilt(L)$ is clopen if, and only if, either $K$ is finite and does not contain $x_\infty$, or $K$ is co-finite and contains $x_\infty$.
  \end{enumerate}
  \end{ourexercise}

\begin{ourexercise}\label{exe:fincofsubsetsspace}
  Let $M$ be the lattice of finite or co-finite subsets of $\mathbb{N}$, as defined in Exercise~\ref{exe:fincofsubsetsprimes}, and let $L$ be the lattice of finite subsets of $\mathbb{N}$ and $\mathbb{N}$ itself. Recall from Exercise~\ref{exe:fincofsubsetsprimes} that there is a bijection between $\PrFilt(L)$ and $\PrFilt(M)$, which is order preserving in only one direction.
  \begin{enumerate}
  \item Prove that the bijection from Exercise~\ref{exe:fincofsubsetsprimes} is a homeomorphism between $\PrFilt(L)$ and $\PrFilt(M)$ in their respective Priestley topologies.
  \item Write down explicitly the isomorphism $\widehat{(-)}$ between $M$ and the clopen down-sets of $\PrFilt(M)$, and show that it restricts correctly to an isomorphism between $L$ and the clopen down-sets of $\PrFilt(L)$.
  \end{enumerate}
\end{ourexercise}

\begin{ourexercise} \label{exe:sumproduct}
Let $L$ and $M$ be distributive lattices with Priestley dual spaces $X$ and $Y$, respectively. Prove that the Priestley dual space of the Cartesian product $L \times M$ is homeomorphic to the disjoint sum of the Priestley spaces $X$ and $Y$. \hint{Recall Exercise~\ref{exe:sumproductdualposet} for the poset part.}
\end{ourexercise}
\begin{ourexercise} \label{exe:closedsubofPriestley}
Let $X$ be a Priestley space, and let $Y \subseteq X$ be a 
subset, which we consider as an ordered topological space equipped with 
the subspace topology and partial order inherited from $X$. Prove that $Y$ is a
Priestley space if, and only if, $Y$ is closed as a subset of $X$. \hint{Use the result of Exercise~\ref{exer:compHaus}.}
\end{ourexercise}

\begin{ourexercise}\label{exe:Priestley-binary-product}
  Let $X$ and $Y$ be Priestley spaces, and let $X \times Y$ be the ordered topological space
  obtained by equipping the set $X \times Y$ with the product topology and the point-wise order.
  \begin{enumerate}
    \item Prove that for any clopen down-sets $K \subseteq X$ and $K' \subseteq Y$, the product $K \times K'$ is a clopen down-set.
    \item  Prove that for any clopen down-set $L \subseteq X \times Y$, there exists a finite collection $(K_1,K_1'), \dots, (K_n,K_n')$ of pairs where, for each $1 \leq i \leq n$, $K_i \in \ClD(X)$ and $K_i' \in \ClD(Y)$, and 
  \[ L = \bigcup_{i=1}^n K_i \times K_i' \ . \] 
  \end{enumerate}
\end{ourexercise}

\begin{ourexercise}\label{exe:Priestley-product-space}
  Prove that the Cartesian product of any family of Priestley spaces, defined by
  equipping the product set with the product topology and the point-wise partial order, is again a Priestley space. 
\end{ourexercise}

\begin{ourexercise}\label{exe:closure-normality-Priestley}
Let $X$ be a Priestley space, $L$ the lattice dual to $X$, and let $C$ be a closed subset of $X$.
\begin{enumerate}
    \item Prove that, for any $x \not\in {\downarrow} C$, there exists a clopen
        down-set $K \subseteq X$ such that $x \in K$ and $K \cap {\downarrow}C
        = \emptyset$.
    \item \label{itm:down-closure} Conclude that 
        \[ {\downarrow} C = \bigcap \{ \widehat{a} \ \mid \ a \in L, C
        \subseteq \widehat{a} \}\ . \]
    \item \label{itm:up-closure} Prove that
        \[ {\uparrow} C = \bigcap \{ \widehat{b}^c \ \mid \ b \in L, C
        \subseteq \widehat{b}^c \}\ . \] 
    \item Show that, if $C, D$ are closed subsets of $X$ such that
        ${\downarrow} C \cap {\uparrow} D = \emptyset$, then there exists a
        clopen down-set $K \subseteq X$ such that $C \subseteq K$ and $K \cap D
        = \emptyset$.
\end{enumerate}
\end{ourexercise}
\begin{ourexercise}\label{exe:Priestley-filters}
  Let $L$ be a distributive lattice and $X$ its dual Priestley space. For any subset $S$ of $L$, define 
  \[ C(S) := \bigcap \{ \widehat{a} \mid a \in S \}\ , \] 
  and for any subset $A$ of $X$, define 
  \[ F(A) := \{ a \in L \mid A \subseteq \wh{a} \}\ .\]  
  \begin{enumerate}
    \item Show that the assignments $S \mapsto C(S)$ and $A \mapsto F(A)$ are a Galois connection between $\cP(L)$ and $\cP(X)$, that is, for any $S \subseteq L$ and $A \subseteq X$, 
    \[ A \subseteq C(S) \iff S \subseteq F(A) \ .\]
    \hint{You can either prove this directly, or show that $C$ and $F$ arise as in Example~\ref{exa:galoisconnection}.}
    \item Prove that, for any $S \in \cP(L)$, $C(S)$ is a closed down-set of $X$, and that every closed down-set of $X$ arises in this way. \hint{Use Exercise~\ref{exe:closure-normality-Priestley}.\ref{itm:down-closure}.}
    \item Prove that, for any $A \in \cP(X)$, $F(A)$ is a filter of $L$, and that every filter of $L$ arises in this way.
    \item Conclude that the lattice of filters of $L$ is anti-isomorphic to the lattice of closed down-sets of $X$.
    \item Prove that under this anti-isomorphism, prime filters of $L$ correspond to principal down-sets of $X$. Indeed, for any prime filter $\mu$ of $L$, show that the closed down-set $C(\mu)$ is equal to ${\downarrow}x$, where $x$ is the unique point of $X$ for which $F_x = \mu$.
  \end{enumerate}
\end{ourexercise}

\begin{ourexercise} \label{exe:dualofNsquaredlattice}
Describe the Priestley dual space of the distributive lattice given in Example~\ref{exa:nojoinirr}.
\end{ourexercise}
\begin{ourexercise} \label{exe:compact-rational-intervals}
Prove directly that the topology on the poset $X$ from Example~\ref{exa:rationalintervals-topology} is compact. 

You may proceed as follows. Let $\mathcal{C}$ be a cover of $X$ that consists of sets of the form $\widehat{[0,q)}, \widehat{(q,1]}$, and $\{q\}$, for $q \in \mathbb{Q} \cap [0,1]$, that is is, 
\[\mathcal{C} = \{ \widehat{[0,a)} : a \in L \} \cup \{ \widehat{(b,1]} : b \in R\} \cup \{ \{p \} : p \in P \},\]
for some subsets $L, R, P \subseteq \mathbb{Q} \cap [0,1]$.
\begin{enumerate}
\item Explain why it is enough to prove that any such cover $\mathcal{C}$ has a finite subcover. \hint{Find the appropriate Lemma in Chapter~\ref{chap:TopOrd}.}
\item Prove that $\sup L \geq \inf R$.
\item Prove that there exist $a \in L$ and $b \in R$ such that $a \geq b$. \hint{Towards a contradiction, if this were false, consider the prime filters around $s := \sup L, $ which will then equal $\inf R$ using the previous item.}
\item Conclude that $\mathcal{C}$ has a subcover consisting of at most three sets.
\end{enumerate}
\end{ourexercise}

\begin{ourexercise}\label{exe:morphisms-concrete}
Let $h \colon M \to L$ be a homomorphism between distributive lattices. Let $X$ and $Y$ be the Priestley dual spaces of $L$ and $M$, respectively, and $f \colon X \to Y$ the continuous order-preserving function dual to $h \colon M \to L$.  Prove that, for any $x \in X$, the prime filter $F_{f(x)}$ is equal to $h^{-1}(F_x)$, and the prime ideal $I_{f(x)}$ is equal to $h^{-1}(I_x)$.
\end{ourexercise}
\begin{ourexercise}\label{exe:Priestley-injsurj}
  Let $h \colon M \to L$ be a homomorphism between distributive lattices and $f \colon X \to Y$ the dual continuous order-preserving function. Prove that $h$ is injective if, and only if, $f$ is surjective, but that $f$ may be injective without $h$ being surjective.
  {\it Note.} We will see the correspondence between sublattices and quotient spaces in more detail in Section~\ref{sec:quotients-and-subs}.  
\end{ourexercise}

\section{Boolean envelopes and Boolean duality}\label{sec:boolenv-duality}
We now show how the original duality for Boolean algebras \parencite{Stone1936} is a special case of Priestley duality. In Chapter~\ref{chap:Omega-Pt}, we will see how, in turn, Priestley duality can also be derived and derives from Stone's original duality for distributive lattices \parencite{Sto1937/38}.

Recall the definition of \emphind{Boolean envelope} given in Definition~\ref{def:booleanenvelope}. Towards deriving Stone duality for Boolean algebras from Priestley duality, we now prove the existence and uniqueness of Boolean envelopes. First, uniqueness is straightforward, and in fact a consequence of the more general fact that ``free objects are unique up to isomorphism''. Compare for example Proposition~\ref{prop:free-unique} in Chapter~\ref{ch:methods} on the free distributive lattice, and Example~\ref{exa:adjunctions} in Chapter~\ref{ch:categories}, which puts free constructions in the general categorical setting of left adjoints to a forgetful functor.
\begin{proposition}\label{prop:boolenv-unique}
Let $L$ be a distributive lattice. If $e \colon L \to B$ and $e' \colon L \to B'$ are Boolean envelopes of $L$, then there is a unique isomorphism $\phi \colon B \to B'$ such that $\phi \circ e = e'$.
\end{proposition}
\begin{proof}
Let $\phi$ be the unique homomorphism $B \to B'$ such that $\phi \circ e = e'$, by the existence part of the defining property of the Boolean envelope of $e$, and similarly $\psi \colon B' \to B$ such that $\psi \circ e' = e$. Since $\psi \circ \phi \circ e = \psi \circ e' = e$, we must have $\psi \circ \phi = \id_B$ by the uniqueness part of the defining property of the Boolean envelope $e$.  Similarly, $\phi \circ \psi = \id_{B'}$.
\end{proof}
Note in particular that it follows from Proposition~\ref{prop:boolenv-unique} that the Boolean envelope of a Boolean algebra $L$ is simply $L$ itself. 
For the construction of the Boolean envelope of an arbitrary distributive lattice, we will use the following lemma, which will also be useful later.

\begin{lemma}\label{lem:Priestleybase}
Let $(X, \tau, \leq)$ be the Priestley space dual to a distributive lattice $L$. The collection
\[ \mathcal{B} := \{ \widehat{a} \setminus \widehat{b} \ | \ a, b \in L\}\]
is a base for the topology $\tau$, and any clopen set of $\tau$ is a finite union of sets in $\mathcal{B}$.
\end{lemma}

\begin{proof}
By the definition of $\tau$ and Exercise~\ref{exer:base}.\ref{itm:base-from-subbase}, the collection of finite intersections of sets of the form $\widehat{a}$ and $\widehat{b}^c$ is a base for $\tau$. This collection is equal to $\mathcal{B}$: if $P, N \subseteq L$ are finite sets, then 
\[ \bigcap_{a \in P} \widehat{a} \cap \bigcap_{b \in N} \widehat{b}^c = \widehat{\bigwedge P} \setminus \widehat{\bigvee N}.\]
Any open set is therefore equal to a (possibly infinite) union of sets in $\mathcal{B}$. If the open set is moreover closed, then it is compact, so that this union can be taken to be finite.
\end{proof}

\begin{proposition}\label{prop:boolenv}
Let $L$ be a distributive lattice. The Boolean algebra $L^-$ of clopen subsets of the Priestley dual space of $L$ together with the embedding $\widehat{(-)} \colon L \to L^-$ is the Boolean envelope of $L$. 
\end{proposition}

\begin{proof} 
We will prove that the Boolean algebra $L^-$, with the embedding $\widehat{(-)}$, has the universal property required of the Boolean envelope (Definition~\ref{def:booleanenvelope}). Let $h \colon L \to A$ be a homomorphism, with $A$ a Boolean algebra. We will construct the unique homomorphism $L^- \to A$ extending $h$ using duality. Write $Y$ for the dual Priestley space of $L$, and $X$ for the dual Priestley space of $A$. By Proposition~\ref{prop:priestleymorphisms}, pick the unique continuous order-preserving function $f \colon X \to Y$ such that $h = h_f$. We will now use this $f$ to define $\bar{h}$. Indeed, let $b \subseteq Y$ be any clopen set. By continuity of $f$, the set $f^{-1}(b)$ is clopen in $X$; define $\bar{h}(b)$ to be the unique element of $A$ such that $\widehat{\bar{h}(b)} = f^{-1}(b)$. Note that $\bar{h}$ is a homomorphism of Boolean algebras, because the function $f^{-1}$ is. Also, for any $a \in L$, we have $\widehat{\bar{h}(\widehat{a})} = f^{-1}(\widehat{a}) = \widehat{h(a)}$, showing that $\bar{h} \circ \widehat{(-)} = h$. To see that $\bar{h}$ is the unique such homomorphism, suppose that $g \colon L^- \to A$ is any homomorphism with $g \circ \widehat{(-)} = h$. Note that the set $E := \{ b \in L^- \mid g(b) = \bar{h}(b)\}$ is a Boolean subalgebra of $L^-$, which contains the image of $L$ under $\widehat{(-)}$. By Lemma~\ref{lem:Priestleybase}, any clopen set of $Y$ is a Boolean combination of sets of the form $\widehat{a}$ for $a \in L$; thus, $E = L^-$, so that $g = \bar{h}$.
\end{proof}

We note an interesting corollary to the proof of Proposition~\ref{prop:boolenv}.

\begin{corollary}\label{cor:boolenv-hull}
  Let $L$ be a distributive lattice and $e \colon L \to A$ a lattice embedding with $A$ a Boolean algebra. Then the Boolean algebra generated by the image of $e$ is isomorphic to $L^-$.
\end{corollary}

\begin{proof}
  Denote by $B$ the Boolean algebra generated by $\im(e)$. By Proposition~\ref{prop:boolenv}, let $\bar{e} \colon L^- \to A$ be the unique Boolean algebra homomorphism such that $\bar{e}(\widehat{a}) = e(a)$ for every $a \in L$. Since $\im(\bar{e})$ is a Boolean algebra containing $\im(e)$, it also contains $B$. By the proof of Proposition~\ref{prop:boolenv}, $\bar{e}$ is dual to the continuous function $f$ dual to $L \to A$. Since $e$ is injective, $f$ is surjective (see Exercise~\ref{exe:Priestley-injsurj}). Thus, $\bar{e}$ is injective too.
\end{proof}

Using our characterization of the Boolean envelope in Proposition~\ref{prop:boolenv}, we now in particular deduce the following characterization of Boolean algebras, from which we will then further obtain Stone duality for Boolean algebras.

\begin{proposition}[\cite{Nachbin47}]\label{prop:boolean-trivial-order}
A distributive lattice $L$ is a Boolean algebra if, and only if, the order on the dual Priestley space of $L$ is trivial.
\end{proposition}
\begin{proof}
Note that a distributive lattice $L$ is a Boolean algebra if, and only if, the embedding of $L$ into $L^-$ is an isomorphism. By Proposition~\ref{prop:boolenv}, this happens if, and only if, every clopen subset of the dual Priestley space $X$ of $L$ is a down-set. The latter is clearly the case if the order on $X$ is trivial, which shows the right-to-left direction. For the left-to-right direction, suppose that the order on $X$ is not trivial: pick $x, y \in X$ with $x \leq y$ and $y \nleq x$. Since $X$ is a Priestley space, pick a clopen down-set $K$ of $X$ with $x \in K$ and $y \not\in K$. The set $X \setminus K$ is clopen, but it contains $y$ and not $x$, so it is not a down-set.
\end{proof}

Note that Proposition~\ref{prop:boolean-trivial-order} implies in particular that the prime filters of a Boolean algebra form an anti-chain; in other words, every prime filter is maximal. This is the generalized version of the fact (Proposition~\ref{prop:birkhoffBA}) that a finite distributive lattice is a Boolean algebra exactly when every join-irreducible element is an atom. 

Stone duality for Boolean algebras may -- anachronistically, see the notes for this chapter -- be viewed as Priestley duality `without the order'. We end this section by giving the definitions in detail.

\begin{definition}\index{ultrafilter}\index{filter!maximal}
  A proper filter $F$ in a lattice $L$ is called a \emph{maximal filter} if it is inclusion-maximal among the proper filters of $L$, that is, for any proper filter $F'$ such that $F \subseteq F'$, we have $F = F'$.
\end{definition}

We note here an alternative characterization of maximal filters in a Boolean algebra: a proper filter $F$ in a Boolean algebra $B$ is prime if, and only if, for every $a \in B$, either $a$ or $\neg a$ is in $F$. A Filter satisfying the latter condition is usually called an \emphind{ultrafilter} in the literature. Exercise~\ref{exe:ultra-altchar} asks you to prove that a filter is indeed prime iff it is  maximal iff it is an ultrafilter. 

\begin{definition}\index{Boolean space}\index{Boolean algebra!dual space of}
A \emph{Boolean space}\footnote{Boolean spaces are also known as `Stone spaces'
in the literature~\parencite{Johnstone1986}, but in other references, the
name `Stone space' has been used to refer to the spaces that play a role in
Stone's more general duality for distributive lattices, and that we call
spectral spaces, see Chapter~\ref{chap:Omega-Pt}.} is a topological space
$(X,\tau)$ which is compact and \emph{totally disconnected}, that is, for any two
points $x, y \in X$, if $x \neq y$, then there exists a clopen subset $K
\subseteq X$ such that $x \in K$ and $y \not\in K$. 
\end{definition}
A space is called \emph{zero-dimensional} if its clopen subsets form a base. An
alternative characterization of Boolean spaces is that they are exactly the
topological spaces which are compact, Hausdorff and zero-dimensional (see
Exercise~\ref{exe:zerodimensional}).
\begin{remark}\label{rem:boolean-priestley}
Note that a topological space $(X,\tau)$ is Boolean if, and only if,
$(X,\tau,=)$ is a Priestley space, see Exercise~\ref{exe:Booleanspace}. Also, a
Priestley space is in particular a compact ordered space based on a Boolean
space. One may ask whether all such ordered spaces satisfy the TOD property. This is
not the case, as was shown in \cite{Stralka80}. In Exercise~\ref{exe:Stralka}
you are asked to work through the example given there.
\end{remark}

The \emph{dual space} of a Boolean algebra $L$ is defined as the dual space of $L$, viewed as a distributive lattice, under Priestley duality. By Proposition~\ref{prop:boolean-trivial-order}, the partial order on the space will be trivial in this case, and, because $(\widehat{a})^c = \widehat{\neg a}$ for any $a \in L$, the definition of the topology can also be slightly simplified, as follows.

\begin{definition}
Let $L$ be a Boolean algebra. The \emph{dual space} of $L$ is the set $X_L$ of ultrafilters of $L$, equipped with the topology generated by the sets
\[ \widehat{a} := \{x \in X_L \ | \ a \in F_x\} \text{ for } a \in L.\]
If $X$ is a Boolean space, its \emph{dual algebra} is the Boolean algebra of clopen subsets of $X$.
\end{definition}

We now use Proposition~\ref{prop:boolean-trivial-order} to deduce three
results from our account of Priestley duality in Section~\ref{sec:topologize}. 
The first follows immediately from
Corollary~\ref{cor:doubledualDL}:
\begin{corollary}\label{cor:doubledualBA}
Let $L$ be a Boolean algebra. Then $\widehat{(-)} \colon L \to \mathcal{P}(X_L)$ is an embedding of Boolean algebras, whose image consists of the clopen subsets of $X_L$. In particular, any Boolean algebra is isomorphic to the algebra of clopen subsets of its dual space.
\end{corollary}
Second, applying Proposition~\ref{prop:priestley-double-dual}.\ref{itm:beta-order-homeo},
in the special case where $X$ is a Boolean space, we get
\begin{corollary}\label{cor:doubledualBS}
Let $X$ be a Boolean space. Then the function sending $x \in X$ to the
point $\beta(x) \in X_{\Clp(X)}$ associated to the ultrafilter 
\[ \mathcal{B}(x) := \{ K \in \Clp(X) \ | \ x \in K \} \]
is a well-defined homeomorphism from $X$ to $X_{\Clp(X)}$.
\end{corollary}

Finally, for morphisms, one specializes the account from
Section~\ref{sec:topologize} by dropping `order-preserving'. If $X$ and $Y$ are
the dual spaces of Boolean algebras $L$ and $M$, and $f \colon X \to Y$ is a
continuous function, then $f^{-1} \colon \mathcal{P}(Y) \to \mathcal{P}(X)$
restricts correctly to the images of $M$ and $L$ under $\widehat{(-)}$. The
homomorphism $h_f \colon M \to L$ is defined by the condition that
$\widehat{h_f(a)} = f^{-1}(\widehat{a})$, for every $a \in M$.

\begin{corollary}
Let $L$ and $M$ be Boolean algebras with dual spaces $X$ and $Y$, respectively.
For any Boolean algebra homomorphism $h \colon M \to L$, there is a unique
continuous function $f \colon X \to Y$ such that $h = h_f$. 
\end{corollary}

In summary, just as we saw in Section~\ref{sec:topologize} for distributive lattices
and in Section~\ref{sec:finDLduality} for finite distributive lattices and
finite Boolean algebras, we have a duality between Boolean algebras and Boolean
spaces. Again, the categorical terminology `dually equivalent' will be defined
precisely in Definition~\ref{dfn:equivalence}. 
\begin{theorem}\label{thm:bool-stone-duality}
The category $\BA$ of Boolean algebras with homomorphisms is dually equivalent
to the category $\BoolSp$ of Boolean spaces with continuous functions.
\end{theorem}

We end this section by giving two well-known examples of Boolean spaces and their dual algebras, which will be used later in this book.
\begin{example}[Stone-{\v C}ech compactification of a discrete space]\label{exa:StoneCechcompactification}
  Let $X$ be an infinite set, and consider the Boolean algebra $\cP(X)$ of all subsets of $S$. The dual space $Y$ of $\cP(X)$ is, up to homeomorphism, the so-called \emphind{Stone-{\v C}ech compactification} of the topological space obtained by endowing the set $X$ with the discrete topology.%
  \nl{$\beta X$}{Stone-{\v C}ech compactification of a discrete space $X$}{beta}

  We now explain what this statement means. First of all, $Y$ is a \emphind{compactification} of the discrete space $X$, that is, $Y$ is a compact Hausdorff space, and $X$, viewed as a discrete space, sits inside $Y$ as a dense subspace. %
  To see that $X$ sits densely in the dual space $Y$ of $\cP(X)$, notice first that, for each $x\in X$, the principal filter $F_x={\uparrow}\{x\}$ is prime or, in other words, an ultrafilter; indeed, these are exactly the principal ultrafilters of $\cP(X)$, see Exercise~\ref{exe:StoneCech-principal}. We denote these points of $Y$ by the corresponding elements of $X$ and, since distinct points of $X$ give distinct principal filters, we thus consider $X$ as a subset of $Y$. Each point of $X$ is isolated in $Y$, because the clopen corresponding to the element $\{x\}$ of $\cP(X)$ contains only $x$. Furthermore, any non-empty basic clopen of $Y$ is of the form
  \[
  \widehat{T}=\{y\in Y \mid T\in F_y\}
  \]
  for some non-empty $T\subseteq X$. Since $T$ is non-empty, there is $x\in T$, and thus $T\in F_x$, or equivalently, $x\in \widehat{T}$, which shows that $X$ is a dense subspace of $Y$.
  
  The fact that $Y$ is the Stone-{\v C}ech compactification of $X$ further means that any other compactification of $X$ is a quotient of $Y$. In Exercise~\ref{exe:StoneCech} you are asked to show this for compactifications of $X$ that fall within Stone duality, that is, compactifications $Z$ of $X$ that are themselves Boolean spaces.
  
  The Stone-{\v C}ech compactification is usually denoted by $\beta$ and for this reason, we will henceforth denote the dual space of $\cP(X)$ by $\beta X$. As it is common in the literature to view the points of $\beta X$ as the ultrafilters of $\cP(X)$, we will also do that here. 
  In the context of our ``neutral dual space'' notation (Notation~\ref{not:neutral}), this means that we take our bijection $F_{(\_)}$ to be the identity function in this specific case. 
  The \emphind{remainder} of $\beta X$ is defined as the subspace ${}^* X := \beta X\setminus X$. The points of ${}^* X$ are the non-principal prime filters, also known as the \emph{free ultrafilters}\index{free ultrafilter}\index{ultrafilter!free}.
\nl{${}^*X$}{remainder of the Stone-{\v C}ech compactification of the (discrete) space $X$}{star}
\end{example}
  
  \begin{example}\label{exa:fincofsubsetsdual}
  Let $X$ be an infinite set, and consider the Boolean algebra $M$ of all subsets of $X$ which are either finite or co-finite as introduced in Exercise~\ref{exe:fincofsubsetsprimes}. It is easy to see, as above, that each element of $X$ gives a distinct principal ultrafilter of $M$ and that these points are isolated in the dual space $Y$ of $M$. So, again, $Y$ is a compactification of $X$. However, in this case there is just \emph{one} non-principal ultrafilter. To see this, let $F$ be an ultrafilter of $M$. If $F$ contains a finite set $S=\{x_1,\dots, x_n\}$, then 
  \[
  \{x_1\}\cup\ \dotsm \ \cup\{x_n\}\in F
  \] 
  and thus $\{x_i\}\in F$ for some $i$ and $F$ is principal. Thus, if $F$ is non-principal, then $F$ does not contain any finite sets. Further, for any finite set $S$, we have
  \[
  S\cup X\setminus S=X\in F
  \]
  and, because $F$ is prime, it follows that $X\setminus S\in F$. That is, $F$ must consist of all the co-finite sets and nothing else.
  
  The dual space of $M$ is known as the \emphind{one-point compactification} of $X$. We will denote it by $X_\infty=X\cup\{\infty\}$. 
  Its basic clopens are the finite subsets of $X$ and the co-finite subsets of $X_\infty$ that contain $\infty$. Thus the opens are all the subsets of $X$ and the co-finite subsets of $X_\infty$ that contain $\infty$.
  \end{example}

\ourexercises

\begin{ourexercise}\label{exe:ultra-altchar}
Let $F$ be a proper filter in a Boolean algebra $B$. Prove that $F$ is prime (and hence maximal) if, and only if, for every $a \in B$, either $a$ or $\neg a$ is in $F$. \hint{The left-to-right direction follows easily from $a \vee \neg a = \top$. For the right-to-left direction, consider the proof of (1) $\Rightarrow$ (2) in Proposition~\ref{prop:birkhoffBA}.}
\end{ourexercise}

\begin{ourexercise}\label{exe:Booleanspace}
Prove that a topological space $(X,\tau)$ is Boolean if, and only if, $(X,\tau,=)$ is a Priestley space.
\end{ourexercise}

\begin{ourexercise}\label{exe:intersection-ultrafilter}
Let $B$ be a Boolean algebra with dual space $(X,\tau)$. Deduce from Exercise~\ref{exe:Priestley-filters} that filters of $B$ are in a bijection with closed subsets of $X$, and that, in particular, if $\mu$ is an ultrafilter of $B$, then $\bigcap \{\wh{a} \mid a \in \mu\}$ is equal to the singleton $\{x\}$, where $x$ is the unique point in $X$ such that $\mu = F_x$.
\end{ourexercise}
\begin{ourexercise}\label{exe:zerodimensional}
Prove that a topological space $(X,\tau)$ is Boolean if, and only if, $(X,\tau)$ is compact, Hausdorff, and zero-dimensional (that is, the clopen subsets form a base).
\end{ourexercise}
\begin{ourexercise}\label{exe:compactcons-is-boolenv}
Prove that the Boolean envelope of a distributive lattice $L$ is isomorphic to the \emph{center} of the congruence lattice of $L$, that is, the lattice (which is always a Boolean algebra) of those congruences on $L$ that have a complement in the complete lattice of congruences on $L$.
\end{ourexercise}
\begin{ourexercise}\label{exe:StoneCech-principal}
  Show that the principal ultrafilters of the Boolean algebra $\cP(X)$ are exactly the 
  sets of the form
  \[ {\uparrow} \{x\} = \{T \in \cP(X) \ \mid \ x \in T \}, \] 
  for some $x \in X$.
\end{ourexercise}
  
  \begin{ourexercise}\label{exe:StoneCech} Let $X$ be a set, $\beta X$ the dual space of $\cP(X)$, and $i\colon X\to\beta X$ the inclusion (see Example~\ref{exa:StoneCechcompactification}). Further, let $B$ be a Boolean algebra with dual space $Z$ and let $e\colon X\to Z$ be any function from $X$ to the set underlying $Z$. Show that there is a unique continuous function $f\colon\beta X\to Z$ so that the following diagram commutes.
   
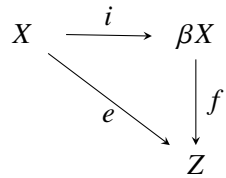
\begin{center}
    \begin{tikzpicture}
    \matrix (m) [matrix of math nodes,row sep=3em,column sep=3em,minimum width=3em]
    {
       X & \beta X \\
         & Z\\};
    \path[-stealth]
      (m-1-1) edge node [above] {$i$} (m-1-2)
      (m-1-1) edge node [below] {$e$} (m-2-2)
      (m-1-2) edge node [right] {$f$} (m-2-2);
    \end{tikzpicture}
      \end{center}
\emph{Note.} A generalization of this result allows one to define an \emph{ordered compactification} of a partially ordered set $P$, by considering the Priestley dual space of the lattice of down-sets of $P$.
  \end{ourexercise}

 \begin{ourexercise}\label{exe:beta-Stonemap}
  Again, let $X$ be a set, $\beta X$ the Stone-\v Cech compactification of $X$ or, equivalently, dual space of $\cP(X)$, and $i\colon X\to\beta X$ the canonical  inclusion map. 
\begin{enumerate}
\item Let $S\subseteq X$ and denote by $e\colon S\to X$ the corresponding inclusion. By dualizing the surjective $e^{-1}\colon\cP(X)\to \cP(S)$, verify that we get an embedding $\beta e\colon\beta S\to \beta X$ extending $e$. In other words, $\beta S$ may be seen as a closed subspace of $\beta X$. 
\item Show that, as a closed subspace of $\beta X$, $\beta S$ consists of those
    ultrafilters of $\beta X$ that are up-sets of their intersection with
    $\cP(X)$. Conclude that this is precisely the clopen subspace $\widehat{S}$ of $\beta X$.
\item Show that the Stone map 
\[
\widehat{(\ )}\colon\cP(X)\to \Clp(\beta X)
\] 
is given by $\widehat{T}=\beta T=\overline{T}$, where $\beta T$ is viewed as a subspace of $\beta X$ as above, and $\overline{(\ )}$ stands for topological closure in $\beta X$. 
\end{enumerate}
 \end{ourexercise}

\begin{ourexercise}[Cantor space]\label{exe:CantorSpace}
The \emphind{Cantor space} $\mathcal C$ can be constructed as a
subspace of the unit interval by repeated removal of middle thirds, as follows.
We define a family of closed subspaces $\{\mathcal C_n\}_{n=0}^\infty$ using
the following operation on real intervals. For $a,b\in\bR$ with $a<b$ define 
\[
F([a,b])=[a,(2a+b)/3]\cup[(a+2b)/3,b],
\]
which removes the open middle third of the interval, thus producing a union of two closed sub-intervals.
Now consider the following subspaces of the unit interval $I=[0,1]$ of the real line $\bR$:
\[
{\mathcal C}_{0}=I=I_{0,1}\quad\text{and}\quad {\mathcal C}_{n}=I_{n,1}\ \cup
\dots\ \cup I_{n,2^n} \ ,
\]
where, for each $n \geq 0$, the sequence of intervals $I_{n+1,1}, \dots,
I_{n+1,2^{n+1}}$ is recursively defined by 
\[
F(I_{n,k})= I_{n+1,2k-1}\cup I_{n+1,2k} \text{ for } 1\leq k\leq 2^n.
\]
Finally, let 
\[
\mathcal C=\bigcap_{n=0}^\infty \mathcal C_n\ .
\]
Let $X$ be the product space $2^{\bN}$, where $2=\{0,1\}$ is the two element
discrete space. 
\begin{enumerate}
\item Show that $X$ is homeomorphic to the Cantor space $\mathcal C$.
\item Let $\leq$ denote the lexicographic order on $X$. Show that it agrees with the order inherited from the unit intervals.
\item Show that $(X,\leq)$ is a Priestley space. Give a description of the dual lattice.
\item \label{itm:cantor-cover} Show that the covering relation associated with the order on $X$ is given by $x\covers y$ if, and only if, there is $n\in\bN$ with $x_i=y_i$ for $i<n$, $x_n<y_n$ (and thus, necessarily, $x_n=0$ and $y_n=1$), while $x_j=1$ and $y_j=0$ for all $j>n$. Further, describe this covering relation relative to the subspace $\mathcal C$ of $\bR$.
\item     Show that, in the case of $(X,\leq)$, we have that $(X,\preceq)$,
    where $${\preceq}={\covers}\cup\Delta$$ with $\Delta=\{(x,x)\mid x\in X\}$, is
    also a partially ordered set. Give an example to show that, in general, the reflexive closure of a cover relation associated to a partial order need not be a partial order.

\end{enumerate}
\end{ourexercise}

\begin{ourexercise}\label{exe:Stralka}
Let $X$ be the product space $2^{\bN}$ equipped with $\preceq$, the reflexive closure of the cover relation $\covers$ associated with the lexicographical order $\leq$ on $X$ obtained from the usual order on $2$. For more detail see Exercise~\ref{exe:CantorSpace}. 
\begin{enumerate}
\item Show that $\preceq$ is a closed relation on $X$. To this end you may proceed as follows: Suppose $x,y\in X$ with $y\not\preceq x$. Then $x\neq y$ and it is not the case that $y\covers x$. Since $x\neq y$ and $\leq$ is a total order on $X$ we have either $x<y$ or $y<x$.
    \begin{enumerate}[label=\arabic*.]
\item Suppose $x<y$. Show that there are clopen sets $U,V$ with $x\in U$, $y\in V$, and $(V\times U)\,\cap\leq\ =\emptyset$. Since $\preceq\ \subseteq\ \leq$, we conclude that $(V\times U)\,\cap\preceq\ =\emptyset$.
\item Suppose $y<x$. Since $x$ is not a cover of $y$ by assumption, show that there exist $z\in X$ with $y<z<x$, and $z$ is not part of any covering pair, that is, $z$ does not terminate in an infinite tuple of $0$'s or of $1$'s (see  Exercise~\ref{exe:CantorSpace}.\ref{itm:cantor-cover}). Further, show that with $U=(z,1]=\{w\in X\mid z<w\}$ and $V=[0,z)=\{w\in X\mid w<z\}$, we have that $x\in U$, $y\in V$, and $(V\times U)\,\cap\preceq\ =\emptyset$.
\end{enumerate}
\item Show that the compact ordered space $(X,\preceq)$ is not a Priestley space. Here we need to show that $(X,\preceq)$ does not have the TOD property. Note that a basic clopen in $X$ is of the form
\[
U=\big(\bigcap_{i\in F}\pi_i^{-1}(0)\big)\cap\big(\bigcap_{j\in G}\pi_j^{-1}(1)\big),
\]
where $F$ and $G$ are disjoint finite subsets of $\bN$.
\begin{enumerate}[label=\arabic*.]
\item Suppose $F\neq\emptyset$ and let $n=\max F$. Show that there is $x\in U$ with $x_i=1$ for all $i>n$. Conclude that $x$ has a cover $y$ that is not in $U$, thus showing that $U$ is not an up-set in $(X,\preceq)$.
\item Suppose that $F=\emptyset$ and show that $U$ is an up-set of $(X,\preceq)$ if, and only if, $G$ is an initial segment of $\bN$.
\item Show that a clopen subset of $X$ is a down-set with respect to $\preceq$ if, and only if, it is an down-set with respect to $\leq$ and thus conclude that $(X,\preceq)$ is not a Priestley space.
\item Find points $x$ and $y$ in $X$ showing that $(X,\preceq)$ does not have the TOD property.
\end{enumerate}
\end{enumerate}
\end{ourexercise}

\notessec
The influential American
mathematician Marshall H. Stone introduced Stone duality for both 
Boolean algebras
and bounded distributive lattices in a series of papers published in the 1930s
\parencite{Stone1934, Stone1935, Stone1936, Stone1937BA, Sto1937/38, Sto1938BA}.
Priestley duality \parencite{Pri1970} recasts Stone duality for bounded distributive lattices 
in terms of Nachbin's theory of ordered spaces \parencite{Nachbin64}, and identifies
an isomorphic dual category to the one given by \cite{Sto1937/38}. See also results in this direction by \cite{Ne59}. The connection between Stone's original duality and Priestley's variant has been explained in several places in the literature, for example~\cite{Cor1975,Fle2000}.

While Stone's representation theorem and duality for Boolean algebras \parencite{Stone1936} are widely known, his results for distributive lattices, which include both a representation theorem and a duality, have recently received less attention in the literature. One reason for this is the fact that the spaces in Stone's duality for distributive lattices \parencite{Sto1937/38} are non-Hausdorff. Priestley's work used a much nicer class of spaces, at the expense of equipping the space with a partial order. In this book, we decided to introduce Priestley's duality first, as it nicely exhibits the interplay between order and topology that is central to the field. We make the connection with Stone's original duality in Chapter~\ref{chap:Omega-Pt}, see Theorem~\ref{thm:Stone-isom-Priestley}.

While the main results of this chapter are well-known, the notation of a `neutral Priestley dual space' is more recent and originates with the papers \cite{GePr1, GePr2}. This notation is particularly useful for studying duality for additional operations which are order preserving in some coordinates and order reversing in others, as we will begin to do in the next chapter, as well as for studying operations that preserve or reverse both join and meet in each coordinate, which was the subject of the papers \cite{GePr1,GePr2}.
\chapter{Duality methods}\label{ch:methods}

In this chapter, we build on the Priestley duality of the previous chapter to
develop methods for analyzing distributive lattices and the morphisms between
them, with applications to logical systems associated to them. We first show in
Section~\ref{sec:free-description} how Priestley duality allows us to easily
compute \emph{free} distributive lattices and Boolean algebras, corresponding
to normal forms for propositional logic formulas. In
Section~\ref{sec:quotients-and-subs}, we establish correspondences between
sublattices of a lattice and quotients of its dual space, and also between
quotients of a lattice and closed subspaces of its dual space. These
Galois-type correspondences are often instrumental for computing dual spaces of
lattices in concrete cases; we will give a few basic examples already in this
chapter and we will make full use of these correspondences in the application
Chapters~\ref{chap:DomThry} and \ref{ch:AutThry}.
Sections~\ref{sec:unaryopduality} and \ref{sec:generalopduality} show how
duality theory can treat classes of functions that are more general than the
lattice homomorphisms that appear in Priestley duality: operators. We first
treat unary operators (that is, functions preserving only finite meets or finite
joins, but not both) in Section~\ref{sec:unaryopduality}, and then in
Section~\ref{sec:generalopduality} we focus on a particularly relevant case of
operators of arity 2, namely \emph{operators of implication type}. This case is
sufficiently complex to show the flavor of the fully general case, while
avoiding notational difficulties when dealing with operators of general arity.
In Sections~\ref{sec:kripke} and \ref{sec:esakia-heyting}, we give two
classical applications of operator duality, namely to deduce Kripke models for
modal logic, and to obtain Esakia duality for Heyting algebras, the algebraic
structures for intuitionistic propositional logic. Duality for operators will
also be applied in both Chapters~\ref{chap:DomThry} and \ref{ch:AutThry}; in
particular, Section~\ref{sec:generalopduality} contains some forward pointers
to where this theory is relevant in both of those chapters. The final section,
Section~\ref{sec:altchains}, contains an application of \emph{discrete} duality
which characterizes the Boolean subalgebra generated by a sublattice. This
result, which is of a more combinatorial flavor, will be used in
Chapter~\ref{ch:AutThry}. 

The methods we develop in this chapter can also be viewed in a more abstract
categorical form, and we will do so for some of them in the next chapter,
Chapter~\ref{ch:categories}. It is useful to first see them ``in action'': this
chapter shows the use of individual instruments, while the next chapter shows
how they all fit together in an ensemble.

\section{Free distributive lattices} \label{sec:free-description}
In this section, we apply Priestley duality to give concrete descriptions of
\emph{free} distributive lattices over sets. We also deduce a description of
the free Boolean algebra over a set as a corollary. We will point out some
connections to propositional logic.

While our discussion here focuses
on free distributive lattices, a large part of the development in this section
is an instance of a much more general \emph{universal algebra}, developed by
Birkhoff; a standard reference for this is the textbook \cite{BurSan2000}. This
section does not require knowledge of universal algebra, but we provide some
pointers for the interested reader. In Example~\ref{exa:adjunctions}, we will
relate the definition of free objects to the categorical notion of adjunction,
to be introduced in Chapter~\ref{ch:categories}.

The construction of free algebras is an important tool in universal algebra.
The reason for this is Birkhoff's Variety Theorem, which links axiomatization
by equations with the model theoretic constructions of quotients, subalgebras,
and Cartesian products, see \cite[Theorem~11.9]{BurSan2000}. The technical crux
of this theorem is the fact that if a class of algebras is axiomatized by
equations, then it contains free algebras over any set
\parencite[Theorem~10.12]{BurSan2000}, and any other algebra of the class is a
quotient of a free algebra \parencite[Corollary~10.11]{BurSan2000}. While this is
interesting, free algebras are also often notoriously difficult to understand
and transferring information to quotient algebras may also be challenging. We
will see here that, in the specific case of distributive lattices,
free distributive lattices are quite simple to understand \emph{dually},
as their dual spaces are Cartesian products of the two-element Priestley space.
Combined with the
methods of Section~\ref{sec:quotients-and-subs} below, which will allow us to
describe quotients of distributive lattices dually as closed subspaces, the free
algebras from Birkhoff's Variety Theorem provide a powerful tool for
distributive-lattice-based algebras. We will see this in action, for example,
in the applications in Section~\ref{sec:DTLF}.

We now define what it means for a distributive lattice to be \emph{free}. This
definition will look familiar if you have previously encountered, for example,
free groups or free semigroups. It is also similar to the definition of Boolean
envelope, Definition~\ref{def:booleanenvelope}. The proper general setting for
these definitions will be developed in Section~\ref{sec:external} in the
subsection on adjunctions, see Example~\ref{exa:adjunctions}.\ref{itm:freeDL}
in particular.
\begin{definition}
  Let $V$ be a set, $F$ a distributive lattice, and $e \colon V \to F$ a
  function. Then $F$ is said to be \emph{free}\index{free distributive lattice}
  over $V$ via $e$ provided that, for every distributive lattice $L$ and every
  function $f \colon V \to L$, there exists a unique lattice homomorphism
  $\bar{f} \colon F \to L$ such that $\bar{f} \circ e = f$, that is, such that the
  following diagram commutes: 
  \begin{center}
    \begin{tikzpicture}
      \matrix (m) [matrix of math nodes,row sep=3em,column sep=3em,minimum width=3em]
      {
        F & L \\
        V &   \\
      };
      \path[-stealth]
      (m-1-1) edge node [above] {$\bar{f}$} (m-1-2)
      (m-2-1) edge node [left] {$e$} (m-1-1)
      (m-2-1) edge node [below] {$f$} (m-1-2);
    \end{tikzpicture}
  \end{center}
\end{definition}
The property expressed by the diagram is called a \emph{universal property} of
$(F,e)$ with respect to distributive lattices, and a function $e \colon V \to
F$ as in this definition is called a \emph{universal arrow}. Given a universal
arrow $e \colon V \to F$ and $f \colon V \to L$ with $L$ a distributive
lattice, we call the homomorphism $\bar{f} \colon F \to L$ the \emph{unique
extension} of $f$ \emph{along} $e$.

Towards the general construction of free distributive lattices, we make two
basic observations about its definition, which in particular justify speaking
of `the' free distributive lattice over a set $V$.
\begin{proposition}\label{prop:free-unique}
  Let $V$ be a set. A free distributive lattice over $V$ is unique up to isomorphism.
  That is, if $e \colon V \to F$ and $e' \colon V \to F'$ are two universal arrows, then there exists a unique lattice isomorphism $\phi \colon F \to F'$ such that $\phi \circ e = e'$.
\end{proposition}
\begin{proof}
  Suppose $e \colon V \to F$ and $e' \colon V \to F'$ are universal arrows. Denote by $\phi \colon F \to F'$ the unique extension of $e'$ along $e$, and by $\psi \colon F' \to F$ the unique extension of $e$ along $e'$. Then the composite function $\psi \circ \phi \colon F \to F$ is an extension of $e$ along $e$, that is, it is a homomorphism with the property that $\psi \circ \phi \circ e = e$, since, by definition, $\phi \circ e = e'$ and $\psi \circ e' = e$. But $\id_F \colon F \to F$ is also an extension of $e$ along $e$, so by uniqueness, we have $\psi \circ \phi = \id_F$. By symmetry, $\phi \circ \psi = \id_{F'}$. Thus, $\phi$ is a lattice isomorphism, and $\phi \circ e = e'$ by construction, as required.
\end{proof}

\begin{notation}
\nl{$F_{\DL}(V)$}{free distributive lattice over a set $V$}{}
  For $V$ a set, we denote by $F_{\DL}(V)$ the free distributive lattice over $V$, and by $e \colon V \to F_{\DL}(V)$ the accompanying universal arrow, if they exist. By Proposition~\ref{prop:free-unique}, this notation is well-defined, if we consider isomorphic lattices as the same. Next, we give an algebraic and a dual construction showing the existence of $F_{\DL}(V)$ and $e \colon V \to F_{\DL}(V)$ for any set $V$.%
\end{notation}

\begin{proposition}\label{prop:free-generated}
  Suppose $F_{\DL}(V)$ is the free distributive lattice over $V$, and $e \colon V \to F_{\DL}(V)$ the accompanying universal arrow. Then $F_{\DL}(V)$ is generated by $e[V]$.
\end{proposition}
\begin{proof}
  Denote by $L$ the sublattice of $F_{\DL}(V)$ generated by $e[V]$. Denote by $p \colon F_{\DL}(V) \to L$ the unique extension along $e$ of the function $e' \colon V \to L$, defined as the co-restriction to $L$ of the function $e$. Write $i \colon L \to F_{\DL}(V)$ for the inclusion homomorphism; the various morphisms are depicted in the following diagram.
  
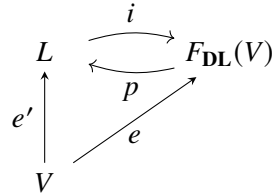
\begin{center}
    \begin{tikzpicture}
      \matrix (m) [matrix of math nodes,row sep=3em,column sep=3em,minimum width=3em]
      {
        L & F_{\mathbf{DL}}(V) \\
        V &                              \\
      };
      \path[-stealth]
      (m-2-1) edge node [left] {$e'$} (m-1-1)
      (m-2-1) edge node [below] {$e$} (m-1-2);
      \draw[->] (m-1-1) to [bend left=15] node [above] {$i$} (m-1-2);
      \draw[->] (m-1-2) to [bend left=15] node [below] {$p$} (m-1-1);
    \end{tikzpicture}
  \end{center}
  Note, similarly to the proof of Proposition~\ref{prop:free-unique}, that $i \circ p$ is an extension of $e$ along $e$, and therefore $i \circ p = \id_{F_{\DL}(V)}$. Hence, $i$ is surjective, proving that $L = F_{\DL}(V)$.
\end{proof}

Note that we have not yet established that a free distributive lattice and accompanying universal arrow actually exist for every set $V$. It is possible to give very general arguments for the existence, because distributive lattices have a finitary axiomatization, see for example \cite[Section~II.10]{BurSan2000}. Here we give two different concrete constructions of the free distributive lattice over $V$, one coming from the general algebraic considerations, the other using Priestley duality, and thus specific to distributive lattices. The advantage of the second construction is that it gives a concrete representation of the free distributive lattice as clopen down-sets of a particular Priestley space. This representation is often useful in applications; we will in particular make use of it in Chapters~\ref{chap:Omega-Pt}~and~\ref{chap:DomThry}.

For the algebraic construction, let $V$ be a set. Since the free distributive
lattice is generated by the image under the universal arrow, it makes sense
that we can build it by making all possible well-formed expressions over $V$
and then take a quotient, as we will do now. 
A \emphind{lattice term}\index{term} with variables
in $V$ is a well-formed expression built from $V$ using the operation symbols
$\vee$, $\wedge$, $\top$, and $\bot$. For example, $a$, $(a \vee b) \wedge a$,
and $\bot \vee a$ are examples of lattice terms with variables in $\{a,b\}$;
while these three terms are distinct syntactic objects, they clearly should be
considered `equivalent' from the perspective of (distributive) lattices. In
fact, the universal property tells us that the free distributive lattice, if it
exists, should have every $V$-generated distributive lattice as a quotient. We
will now define an equivalence relation $\equiv$ on the set of terms which
makes this idea precise.

We write $T(V)$ for the set of all lattice terms with variables in $V$.
Note that, given a function $f \colon V \to L$ with $L$ any structure equipped
with operations $\vee$, $\wedge$, $\top$ and $\bot$, we may inductively define
an \emph{interpretation function} $\tilde{f} \colon T(V) \to L$, by
$\tilde{f}(v) := f(v)$ for $v \in V$, $\tilde{f}(\top) := \top$,
$\tilde{f}(\bot) := \bot$, $\tilde{f}(t \vee t') := \tilde{f}(t) \vee
\tilde{f}(t')$ and $\tilde{f}(t \wedge t') := \tilde{f}(t) \wedge
\tilde{f}(t')$ for any terms $t, t'$. Now consider the equivalence relation
$\equiv$ on $T(V)$, defined by $t \equiv t'$ if, and only if, 
for every function  $f \colon V \to L$,  with
$L$ a distributive lattice, we have
$\tilde{f}(t) = \tilde{f}(t')$. The crucial insight is now that the set
$T(V)/{\equiv}$ of equivalence classes of lattice terms naturally admits the
structure of a distributive lattice, as follows. First note that the
equivalence relation $\equiv$ is congruential on $T(V)$, that is, if $a \equiv a'$
and $b \equiv b'$, then $a \wedge b \equiv a' \wedge b'$ and $a \vee b \equiv
a' \vee b'$, as is easily verified using the fact that $\tilde{f}$ in the
definition of $\equiv$ preserves the operations $\vee$ and $\wedge$ (see
Exercise~\ref{exe:freeDL-algebraic}). Therefore, we have well-defined
operations on $T(V)/{\equiv}$ given by \[ \top := [\top]_\equiv, \bot :=
[\bot]_\equiv, [a]_{\equiv} \vee [b]_{\equiv} := [a \vee b]_{\equiv}, \text{
and } [a]_{\equiv} \wedge [b]_{\equiv} := [a \wedge b]_{\equiv},\]
for any $a, b \in T(V)$.
To see that $(T(V)/{\equiv}, \top, \bot, \vee, \wedge)$ is a distributive
lattice under these operations, one may verify that the defining equations from
Section~\ref{sec:lattices} hold in $T(V)/{\equiv}$ (see
Exercise~\ref{exe:freeDL-algebraic}). Finally define $e\colon V\to
T(V)/{\equiv}$ by $e(v)=[v]_{\equiv}$.
\begin{proposition}\label{prop:freeDL-algebraic}
  The distributive lattice $T(V)/{\equiv}$ is free over $V$ via the universal arrow $e$.
\end{proposition}
\begin{proof}
  Let $f \colon V \to L$ be a function to a distributive lattice $L$. The  interpretation function $\tilde{f} \colon T(V) \to L$ has the property that $\tilde{f} \circ e = f$, and, for any $t \equiv t'$, we have $\tilde{f}(t) = \tilde{f}(t')$, by definition of $\equiv$. Therefore the function $\bar{f} \colon T(V) \to L$, defined by $\bar{f}([t]_{\equiv}) := \tilde{f}(t)$ is well-defined. Observe that the function $\bar{f}$ is a homomorphism by definition of the operations on $T(V)/{\equiv}$, and it is an extension of $V$. We leave uniqueness of $\bar{f}$ as an exercise (see Exercise~\ref{exe:freeDL-algebraic}).
\end{proof}
\begin{remark}\label{rem:LT-algebra}
  The above construction is essentially that of the \emphind{Lindenbaum-Tarski algebra} for a propositional logic without negation over variables $V$. A completely analogous development, replacing $\DL$ by $\BA$ throughout, gives a construction of the free Boolean algebra over a set of variables; also see Corollary~\ref{cor:freeBA-duality} below.
\end{remark}

Before taking on the general duality-theoretic construction, let's look at a small example.
\begin{example}\label{exa:twogen-logic}
  Let $V = \{p,q\}$. We get lattice terms $\bot,\top,p,q,p\wedge q$ and $p\vee q$, and it is not too difficult to see that, in any distributive lattice, the interpretation of all other terms must be equal to one of these. That is, $F_{\DL}(\{p,q\})$ must be the lattice depicted on the left in Figure~\ref{fig:freedltwo} below. There, the four join-primes of $F_{\DL}(\{p,q\})$ are identified as the circled nodes.
  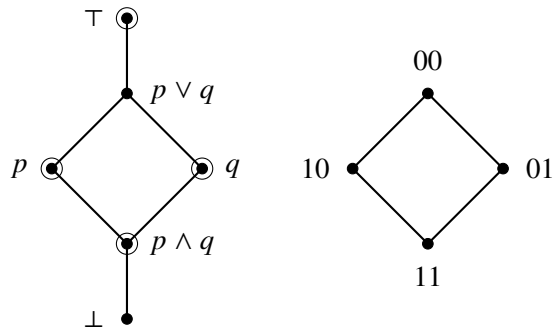
\begin{figure}[htp]
    \begin{center}
      \begin{tikzpicture}[node distance=15pt, label distance=1pt]
        \pos{(0,-1)} %
        \pob{(0,0)} %
        \pob{(-1,1)} %
        \pob{(1,1)} %
        \pos{(0,2)} %
        \pob{(0,3)} %

        \draw[thick] (0,-1) -- (0,0) -- (-1,1) -- (0,2) -- (0,3);
        \draw[thick] (0,0) -- (1,1) -- (0,2);

        \node[label=left:{$\bot$}] at (0,-1) {};
        \node[label=right:{$p \wedge q$}] at (0,0) {};
        \node[label=left:{$p$}] at (-1,1)  {};
        \node[label=right:{$q$}] at (1,1) {};
        \node[label=right:{$p \vee q$}] at (0,2) {};
        \node[label=left:{$\top$}] at (0,3) {};

        \po{(4,0)} %
        \po{(3,1)} %
        \po{(5,1)} %
        \po{(4,2)} %
        \draw[thick] (4,0) -- (3,1) -- (4,2) -- (5,1) -- (4,0);

        \node[label=below:{11}] at (4,0) {};
        \node[label=left:{10}] at (3,1) {};
        \node[label=right:{01}] at (5,1) {};
        \node[label=above:{00}] at (4,2) {};
      \end{tikzpicture}
    \end{center}
    \caption{The free distributive lattice on two generators, and its dual poset.}
    \label{fig:freedltwo}
  \end{figure}

  In the figure on the right, we depict the dual poset as $(2,\preceq)^{\{p,q\}}$, that is, the set of functions from $\{p,q\}$ to $2$, ordered pointwise, where the order on $2$ is given by $1\preceq 0$.  Each function $f \colon \{p, q\} \to 2$ is denoted as a pair $f(p)f(q)$, for example, $10$ is the function sending $p$ to $1$ and $q$ to $0$. 
  
  The idea behind this representation of the dual poset of $F_{\DL}(\{p,q\})$ is the following. As we know from Definition~\ref{def:Priestleydualspace}, we may choose to see the dual poset of $F_{\DL}(\{p,q\})$ as $\Hom_\DL(F_{\DL}(\{p,q\}),\btwo)$ in the reverse point-wise order. Since $F_{\DL}(\{p,q\})$ is generated by $p$ and $q$, each homomorphism $h \colon F_{\DL}(\{p,q\}) \to \btwo$ is totally determined by its restriction to $\{p,q\}$. Also, by the universal property, \emph{every} function from $\{p,q\}$ to $2$ is the restriction of such a homomorphism.

For example, the element $p$ in the lattice corresponds to the down-set  $\{10, 11\}$ of the  dual poset while the element $\bot$ corresponds to the empty down-set, and the element $p \vee q$ corresponds to the down-set $\{10, 01, 11\}$.
  Thus, we see that the concrete incarnation of $F_{\DL}(\{p,q\})$ obtained via duality is given by the universal arrow
  \[
    e\colon \{p,q\}\to \Down\big((2,\preceq)^{\{p,q\}}\big)\!, \ p\mapsto {\downarrow} \chi_{p},
  \]
  where $\chi_{p}$ is the characteristic function of the singleton $\{p\}$ and the order on $2$ is given by $1\preceq 0$.

  From a logic perspective, the elements of the dual poset correspond to lines in a truth table for propositional logic on two variables, and the subset associated to a formula $\phi$ is the set of lines in the truth table at which $\phi$ is true. This is precisely as conceived by logicians going back to the fundamental work by \cite{Boole1847}.

  We here obtained the universal arrow $e$  as a map to the \emph{down}-sets of $2^{\{p,q\}}$ in the point-wise order relative to the ``upside down'' order on $2$, in which $1\preceq 0$. Note that this map may alternatively be described using \emph{up}-sets and the poset $\btwo$, which is $\{0,1\}$ with the usual order, in which $0\leq 1$, as the map
  \[
    \{p,q\}\to \Up\big(\btwo^{\{p,q\}}\big)\!,\  p\mapsto {\uparrow} \chi_{p}.
  \]
  These are two alternate descriptions of one and the same concrete incarnation of the free distributive lattice over the set $\{p,q\}$. Our reason for sticking with the ``upside down'' version of $2$ is that it morally is not the lattice $\btwo$, but the poset dual to $\bf 3$, the three element chain, which is a free 
  distributive lattice on a single generator, see Remark~\ref{rem:freeDL-duality-altproof} and Exercise~\ref{exe:freeiscoproduct}.
\end{example}
We now proceed to make the corresponding argument for an arbitrary set $V$. %
For a set $V$, let $2^V$ denote the set of functions from $V$ to $2$, which may
be viewed as the $|V|$-fold product of the two-element set $2$. 
We will show in Proposition~\ref{prop:freeDL-duality} that the dual space of
$F_{\DL}(V)$ is order-homeomorphic to the following \emph{ordered generalized
Cantor space}\index{Cantor space!ordered}. %

\begin{definition}\label{dfn:ordered-general-cantor}
  Let $V$ be a set. We define the topology $\pi$ on $2^V$ to be the product topology, where $2$ carries the discrete topology. We define the partial order $\preceq$ on $2^V$, for $x, y \in 2^V$, by
  \[ f \preceq g \iff \text{ for all } v \in V, \text{ if } g(v) = 1 \text{ then } f(v) = 1.\]
  The ordered topological space $(2^V, \pi, \preceq)$ is called the \emph{ordered generalized Cantor space}\index{Cantor space!ordered} over $V$. The topological space $(2^V, \pi)$ is called \emphind{generalized Cantor space}, and $(2^\mathbb{N}, \pi)$ is the (classical) \emphind{Cantor space} (see Exercise~\ref{exe:CantorSpace}).
\end{definition}

Note that, by definition of the product topology, 
the topology $\pi$ on $2^V$ has as a subbase the sets of the form
\[ \sem{v \mapsto a} := \{ f \in 2^V \mid f(v) = a \},\]
where $v$ ranges over the elements of $V$ and $a = 0,1$.

Note also that $\preceq$ is the pointwise order on $2^V$ with respect to the order on $2$ in which $1 \preceq 0$, as in Example~\ref{exa:twogen-logic} above; see Remark~\ref{rem:free-dl-as-upsets} below for an alternative description using up-sets and the order $\btwo$, in which $0 \leq 1$.%

\begin{proposition}\label{prop:freeDL-duality}
  Let $V$ be a set. The Priestley dual space of the free distributive lattice over $V$ is order-homeomorphic to the ordered generalized Cantor space over $V$, and the function $s \colon V \to \ClD(2^V)$, which sends $v \in V$ to the clopen down-set $s(v) := \sem{v \mapsto 1}$, is a universal arrow.
\end{proposition}

Our proof here uses the algebraic fact, already proved above, that there exists a free distributive lattice over $V$. It is possible to give a proof ``from scratch'' that does not use this fact (see Exercise~\ref{exe:freeDL-dual-altproof}).
\begin{proof}
  Let $e \colon V \to F$ be a free distributive lattice over $V$ and denote by
  $(X,\tau,\leq)$ the Priestley dual space of $F$. We will first exhibit an
  order-homeomorphism $\phi$ between $X$ and the ordered generalized Cantor space over
  $V$. Note that, for any point $x \in X$, the corresponding homomorphism $h_x
  \colon F \to \btwo$ must be of the form $\bar{f_x}$ for some function $f_x
  \colon V \to 2$: indeed, $h_x \colon F \to \btwo$ is equal to $\bar{f_x}$ for
  $f_x := h_x \circ e$. Let us write $\phi$ for the surjective function from
  $2^V$ to $X$ that sends $f \colon V \to 2$ to the element $x \in X$ with $h_x
  = \bar{f}$.%

  We now show that $\phi$ is an order-embedding, that is, that for any functions $f, g \colon V \to 2$, we have $f \preceq g$ if, and only if, $\phi(f) \leq \phi(g)$ in $X$. First, if $\phi(f) \leq \phi(g)$ in $X$, then $\bar{f} \geq \bar{g}$ pointwise, so in particular $f = \bar{f} \circ e \geq \bar{g} \circ e = g$ pointwise, which means $f \preceq g$ by definition.
  For the other direction, note that, if $f \preceq g$, then the set $\{u \in F \mid \bar{f}(u) \geq \bar{g}(u)\}$ is a sublattice of $F$ containing $V$, and must therefore be equal to $F$ by Proposition~\ref{prop:free-generated}. It follows that $\phi$ is an order-isomorphism between the posets $(2^V, \preceq)$ and $(X,\leq)$.

  We now show that $\phi$ is a homeomorphism from $(2^V,\pi)$ to $(X,\tau)$.
  Note first that, by Proposition~\ref{prop:free-generated}, the set $e[V]$
  generates $F$. It follows that the image of $e[V]$ under the isomorphism
  $\widehat{(-)} \colon F \to \ClD(X)$ generates $\ClD(X)$. Thus, by definition
  of the Priestley topology $\tau$ on $X$, to establish continuity of $\phi$,
  it suffices to show that $\phi^{-1}(\widehat{e(v)})$ is clopen for
  every $v \in V$. Let $v \in V$ and $f \in 2^V$. We then have that $\phi(f)
  \in \widehat{e(v)}$ if, and only if, $\bar{f}(e(v)) = 1$, if, and only if,
  $f(v) = 1$. Thus, for any $v \in V$, we have 
  \begin{equation}\label{eq:phi-to-cylinder}
    \phi^{-1}(\widehat{e(v)}) = \sem{v \mapsto 1},
  \end{equation}
  which is a clopen set in the product topology $\pi$. Since the sets $\sem{v \mapsto 1}$ and their complements form a subbase for the clopen sets of the topology $\pi$, and $\phi$ is a bijection, we also immediately obtain from the equality (\ref{eq:phi-to-cylinder}) that $\phi$ is an open map. Finally, (\ref{eq:phi-to-cylinder}) shows that the diagram
  \begin{center}
    \begin{tikzpicture}
      \node (F) at (0,2) {$F$};
      \node (CLDX) at (2,2) {$\ClD(X)$};
      \node (V) at (0,0) {$V$};
      \node (CLD2V) at (2,0) {$\ClD(2^V)$};

      \draw[->] (V) to node[left] {$e$} (F);
      \draw[->] (F) to node[above] {$\widehat{(-)}$} (CLDX);
      \draw[->] (CLDX) to node[right] {$\phi^{-1}$} (CLD2V);
      \draw[->] (V) to node[above] {$s$} (CLD2V);
    \end{tikzpicture}
  \end{center}
  commutes, and therefore the function $s$ is a universal arrow, since the arrow $e$ is universal, $\phi^{-1} \circ \widehat{(-)}$ is an isomorphism, and isomorphisms preserve universal arrows (see Exercise~\ref{exe:universal-preserved-by-iso}).
\end{proof}

\begin{remark}\label{rem:free-dl-as-upsets}
  The alternative explanation of the free distributive lattice as up-sets instead of down-sets from Example~\ref{exa:twogen-logic} carries over to the general case. In particular, note that the function $s$ from Proposition~\ref{prop:freeDL-duality} may also be described as the function $s \colon V \to \ClU(\btwo^V)$ which sends $v$ to $\sem{v \mapsto 1}$, where now $\btwo$ is the order on $2$ in which $0 \leq 1$.

  A different kind of order symmetry in this context is the fact that the free distributive lattice is anti-isomorphic to itself. In terms of the concrete representation of the free distributive lattice as clopen down-sets of the ordered generalized Cantor space, we may consider the function $t \colon V \to \ClU(2^V)^\op$ that sends $v \in V$ to $t(v) := \sem{v\mapsto 0}$. Its unique extension $\bar{t}$ along $s$ is an isomorphism from $\ClD(2^V)$ to $\ClU(2^V)^\op$, which may be described concretely by sending a clopen down-set $D$ to the clopen up-set $2^V \setminus D$.
\end{remark}

\begin{remark}\label{rem:freeDL-duality-altproof}
  We outline an alternative, more abstract proof of Proposition~\ref{prop:freeDL-duality}, by using some forward references to the categorical language that we will develop in Chapter~\ref{ch:categories} 
  (see also Exercise~\ref{exe:freeiscoproduct}).  First, using the categorical fact that `left adjoints preserve colimits' (see Exercise~\ref{exe:left-ad-pres-colimits}), the free distributive lattice over a set $V$ must be the $V$-fold coproduct of the free distributive lattice over a singleton set. Now, the free distributive lattice on one generator $p$ is easily seen to be the three-element chain $\{ \bot < p < \top \}$, whose dual is $\{p, \top\}$. Using the fact that the category of Priestley spaces contains arbitrary products of finite posets (see Example~\ref{exa:priestley-as-profinite posets}), we may now recover Proposition~\ref{prop:freeDL-duality} directly from the categorical fact that a dual equivalence sends any coproduct of distributive lattices to a corresponding product of Priestley spaces.
\end{remark}

The \emphind{free Boolean algebra} over a set $V$ can be defined in an entirely analogous way to the free distributive lattice: it is a Boolean algebra $A$, together with a function $e \colon V \to A$, such that for every function $f \colon V \to B$ with $B$ a Boolean algebra, there exists a unique homomorphism $\bar{f} \colon A \to B$ such that $\bar{f} \circ e = f$. The proofs of Proposition~\ref{prop:free-unique}~and~\ref{prop:free-generated} can be carried out in the same way for the free Boolean algebra, as well as the algebraic construction of the free Boolean algebra: this is now the Lindenbaum-Tarski algebra of classical propositional logic. Alternatively, we may combine the free distributive lattice with the Boolean envelope\index{Boolean envelope} construction (see Section~\ref{sec:boolenv-duality}) to obtain the free Boolean algebra.
\begin{lemma}\label{lem:freeDLplusboolenv}
  Let $V$ be a set. The Boolean envelope of the free distributive lattice over $V$ is the free Boolean algebra over $V$.
\end{lemma}
\begin{proof}
  Let $i \colon V \to F_{\DL}(V)$ be the free distributive lattice over $V$, and let
  $j \colon F_{\DL}(V) \to F_{\DL}(V)^-$ be the Boolean envelope of $F_{\DL}(V)$. We claim that the
  composite $e := j \circ i \colon V \to F_{\DL}(V)^-$ has the required universal
  property for the free Boolean algebra. Indeed, for any function $f \colon V \to
  B$ with $B$ a Boolean algebra, there is first, by definition of the free
  distributive lattice, a unique lattice homomorphism $\bar{f} \colon F_{\DL}(V) \to B$
  with $\bar{f} \circ i = f$, and then, by definition of the Boolean envelope, a
  unique homomorphism $\bar{f}^- \colon F_{\DL}(V)^- \to B$ such that $\bar{f}^- \circ
  j = \bar{f}$. We are in the situation of the following diagram:
  \begin{center}
\begin{tikzpicture}
      \node (FV) at (0,3) {$F_{\DL}(V)$};
      \node (FVb) at (3,3) {$F_{\DL}(V)^-$};
      \node (V) at (0,0) {$V$};
      \node (B) at (3,0) {$B$};

      \draw[->] (V) to node[left] {$i$} (FV);
      \draw[->] (FV) to node[above] {$j$} (FVb);
      \draw[->] (V) to node[left,xshift=-5mm,yshift=-4mm] {$e$} (FVb);
      \draw[->,dashed] (FVb) to node[right] {$\bar{f}^-$} (B);
      \draw[->,dashed] (FV) to node[above,yshift=4mm,xshift=-3mm] {$\bar{f}$} (B);
      \draw[->] (V) to node[above] {$f$} (B);
\end{tikzpicture}
\end{center}
  and we get
  $$\bar{f}^- \circ e = \bar{f}^- \circ j \circ i = \bar{f} \circ i = f\ . $$ 
  The uniqueness of $\bar{f}^-$ is clear from the uniqueness parts
  of the universal properties of the free distributive lattice and the Boolean
  envelope. 
\end{proof}
\begin{corollary}\label{cor:freeBA-duality}
  Let $V$ be a set. The Boolean algebra of clopen sets of the generalized Cantor space, $(2^V, \pi)$, is the free Boolean algebra over $V$ via the function $s$ which sends any $v \in V$ to the clopen set $\sem{v \mapsto 1}$.
\end{corollary}
\begin{proof}
  Combine Proposition~\ref{prop:boolenv}, Proposition~\ref{prop:freeDL-duality}, and Lemma~\ref{lem:freeDLplusboolenv}.
\end{proof}
In case $V$ is countably infinite, the free Boolean algebra over $V$ is
countably infinite, and it is \emphind{atomless}, that is, it does not have any
atoms. In fact, one may use model-theoretic techniques to prove that this is, 
up to isomorphism, the \emph{unique} countably infinite and atomless Boolean algebra,
and it is a central structure in several parts of logic, see for example
\cite[Cor.~5.16]{BAhandbook}. Its dual is the
classical Cantor space $2^\bN$ (also see Exercise~\ref{exe:CantorSpace}).

\ourexercises

\begin{ourexercise}\label{exe:universal-preserved-by-iso}
Let $V$ be a set, and suppose that $e \colon V \to F$, $e' \colon V \to F'$ are functions to distributive lattices $F$ and $F'$, and $\phi \colon F' \to F$ is a lattice isomorphism such that $\phi \circ e' = e$. Prove that $e$ is a universal arrow if, and only if, $e'$ is a universal arrow.
\end{ourexercise}
\begin{ourexercise}\label{exe:freeBA-on-2}
How many elements does the free Boolean algebra on two generators have? Draw its Hasse diagram, labeling each of its elements with a corresponding Boolean algebra term.
\end{ourexercise}
\begin{ourexercise}
  \label{exe:locallyfinite}
  Let $L$ be a distributive lattice and suppose that $G \subseteq L$ is a set of generators for $L$. Prove that there is a surjective homomorphism $F_{\DL}(G) \to L$. Conclude that any \emphind{finitely generated distributive lattice} is finite. Formulate and prove the analogous result for Boolean algebras.

    \emph{Note.} This exercise is an instance of the general universal
    algebraic fact that, in a variety of algebras, all finitely generated
    algebras are finite if, and only if, the free finitely generated algebras are finite. Such varieties are known as \emphind{locally finite} varieties, and the results in this section show in particular that both the varieties of distributive lattices and of Boolean algebras are locally finite.
\end{ourexercise}

\begin{ourexercise}\label{exe:freeDL-algebraic}
  This exercise asks you to supply some details of the proof of Proposition~\ref{prop:freeDL-algebraic}.
  \begin{enumerate}
    \item Show that $\equiv$ is congruential on $T(V)$.
    \item Show that $T(V)/{\equiv}$, with operations defined as in the paragraph before Proposition~\ref{prop:freeDL-algebraic}, satisfies all the defining equations of a distributive lattice. \hint{This is almost immediate from the definition of $\equiv$.
    \item Prove in detail that the function $\bar{f}$ defined in the proof of Proposition~\ref{prop:freeDL-algebraic} is a homomorphism and that $\bar{f} \circ e = f$.}
    \item Prove that if $g \colon T(V)/{\equiv} \to L$ is a homomorphism and $g \circ e = f$, then $g = \bar{f}$.
  \end{enumerate}
\end{ourexercise}

\begin{ourexercise}\label{exe:Vto2}
  Let $V$ be a set and $L$ a distributive lattice. Show that the poset $L^V$ of all functions from $V$ to $L$ with the pointwise order is order isomorphic to the poset of all homomorphisms from $F_{\DL}(V)$ to $L$ also with the pointwise order.
\end{ourexercise}

\begin{ourexercise}\label{exe:freeDL-dual-altproof}
  This exercise outlines a proof of Proposition~\ref{prop:freeDL-duality} that does not rely on generalities from universal algebra, and gives a more concrete construction of the extension homomorphisms. %
  Let $V$ be a set, and $(2^V, \preceq, \pi)$ the ordered generalized Cantor space over $V$. %

  \begin{enumerate}
    \item Prove that, for each $v \in V$, the set $s(v) := \sem{v \mapsto 1}$ is a clopen down-set in $(2^V, \preceq, \pi)$.
    \item Let $K$ be a clopen down-set of $(2^V, \preceq, \pi)$. Prove, using the compactness of $2^V$ and the definition of the order, that there exists a finite collection $\mathcal{C}$ of finite subsets of $V$ such that $K = \bigcup_{C \in \mathcal{C}} \bigcap_{v \in C} s(v)$.
    \item Give a direct proof that, for any function $h \colon V \to L$, with $L$ a distributive lattice, there exists a unique homomorphism $\bar{h} \colon \ClD(2^V) \to L$ such that $\bar{h} \circ s = f$.
  \end{enumerate}
\end{ourexercise}

\section{Quotients and subs}\label{sec:quotients-and-subs}
In many applications of Priestley duality, including some in the later chapters
of this book, we use the dual space of a lattice to study its quotient lattices
and sublattices, or we use the dual lattice of a space to study its quotient
spaces and subspaces. Here, \emph{closed} subspaces of a Priestley space will
play an important role, because they are exactly the subspaces that are
themselves Priestley spaces in the inherited order and topology (see
Exercise~\ref{exe:closedsubofPriestley}). In this section, we exhibit two
Galois connections, between relations and subsets, which in particular allow us
to prove that (i) the quotients of a distributive lattice are in an
order-reversing bijection with the closed subspaces of its Priestley dual
space; and (ii) the sublattices of a distributive lattice are in an
order-reversing bijection with the Priestley quotients of its Priestley dual
space. While (i) and (ii) can also be deduced in a more abstract way from
Theorem~\ref{thm:priestleyduality} using category theory (see
Theorem~\ref{thm:subquotient-duality-categorically}), the Galois connections
that we develop in this chapter do a bit more. In concrete applications, these
Galois connections are what allow us to define so-called \emphind{equations} on
either side of the duality.

Let $L$ be a distributive lattice with Priestley dual space $X$. We will first
show how to associate with any binary relation $R$ on
$L$ a closed subspace $\sem{R}$ of $X$.   
Here, a pair that belongs to such a binary relation on a lattice 
is thought of as an equation,
and it acts as a constraint defining a subset on the spatial side of the duality.
This motivates the following notation and definition.

\begin{notationnum}\label{notation:lattice-equation}
When $L$ is a lattice, we will denote, in this subsection, pairs of elements 
by the notation $a \approx b$, instead of $(a,b)$.
\nl{$a \approx b$}{pair of lattice elements, seen as equation on the dual space}{approx}
\end{notationnum}

\begin{definition}\label{def:lattice-quotient-duality}
Let $L$ be a distributive lattice with  dual Priestley space $X$.
We define the binary relation ${\models} \subseteq X \times L^2$ by, 
for any pair of elements $a \approx b \in L^2$ and $x \in X$,
\nl{$x \models a \approx b$}{satisfaction relation between an element $x$ of a dual space and a pair $(a,b)$ of elements of a lattice, viewed as an equation}{}
\[ x \models a \approx b \iffdef x \in \widehat{a} \text{ if, and only if } x \in \widehat{b}.\]
For any binary relation $R \subseteq L^2$, we define 
\[ x \models R \iffdef \text{for every } a \approx b \text{ in } R,\  x \models a \approx b. \] 
Moreover, for any $a \approx b \in L^2$, we define a clopen subset
$\sem{a \approx b}$ of  $X$ by 
\[\sem{a \approx b} :=\{x \in X \ | \ x \models a \approx b \} = (\widehat{a} \cap
\widehat{b}) \cup (\widehat{a}^c \cap \widehat{b}^c),\]
\nl{$\sem{a \approx b}$}{space associated with an equation between lattice elements $a$ and $b$}{}
and we extend this assignment to any binary relation $R$ on $L$, by defining
\[    \sem{R} := \{x \in X \ | \ x \models R\} = \bigcap_{a \approx b \in R} \sem{a \approx b}. \] 
\nl{$\sem{R}$}{space associated with a binary relation $R$ on a lattice}
Conversely, for any element $x \in X$, we define a binary relation $\theta(x)$ on $L$ by
\nl{$\theta(x)$}{congruence on a lattice associated with an element $x$ of its dual space}{}
\[\theta(x) := \ker(h_x) = (F_x \times F_x) \cup (I_x \times I_x) = \{a \approx
b \ | \ x \models a \approx b \}.\]
We extend this assignment to subsets $S \subseteq X$ by defining
\nl{$\theta(S)$}{congruence on a lattice associated with a subset $S$ of its dual space}{}
\[\theta(S) := \bigcap_{x \in S} \theta(x) = \{a \approx b \ | \ a, b \in L \text{ such that } \widehat{a} \cap S = \widehat{b} \cap S\}.\]
\end{definition}
Note that the above definitions of $\sem{-}$ and $\theta$ are a special case of
the functions $u$ and $\ell$ introduced in Example~\ref{exa:galoisconnection},
in the case of the relation ${\models} \subseteq X \times L^2$.
 
\begin{proposition}\label{prop:quotientlattice-subspace}
  Let $L$ be a distributive lattice with dual Priestley space $X$.  The two functions $\sem{-} \colon \mathcal{P}(L^2) \leftrightarrows \mathcal{P}(X) \colon \theta$ form a Galois connection, whose fixed points on the left are the lattice congruences on $L$, and whose fixed points on the right are the closed subspaces of $X$.
\end{proposition}

\begin{proof}
  The functions indeed form a Galois connection since they are the contravariant adjoint pair
  between the powersets of $L^2$ and $X$, formed as in Example~\ref{exa:galoisconnection}, in the
  case of the relation ${\models}$.

  For the statement about fixed points, we need to show that the image of
  $\sem{-}$ consists of the closed subsets of $X$, and that the image of
  $\theta$ consists of the congruences on $L$. First, since each $\sem{a
  \approx b}$ is clopen, $\sem{R}$ is closed for every binary relation $R$.
  Conversely, let $C \subseteq X$ be any closed set. We show that $C =
  \sem{\theta(C)}$.  By adjunction, we have $C \subseteq \sem{\theta(C)}$. We
  prove the other inclusion by contraposition. Let $x \not\in C$ be arbitrary.
  Using the base for the Priestley topology given in
  Lemma~\ref{lem:Priestleybase}, pick $a, b \in L$ such that $x \in
  \widehat{a} \setminus \widehat{b}$ and $C$ is disjoint from $\widehat{a}
  \setminus \widehat{b}$. Note that the latter implies that the pair $a \approx a \wedge
  b$ is in $\theta(C)$, since for any $y \in C$, we have $y \in \widehat{a}$ if,
  and only if, $y \in \widehat{a} \cap \widehat{b}$. However, $x \in
  \widehat{a}$ but $x \not\in \widehat{a \wedge b}$, so $x \not\in
  \sem{\theta(C)}$, as required.

  We now show that the image of $\theta$ consists of the congruences on $L$.
Note that $\theta(S)$ is a congruence for any $S$, as it is an intersection of
congruences (see Exercise~\ref{exe:congruencelattice}). 
Let $\phi$ be a congruence on $L$. We show that $\phi =
\theta(\sem{\phi})$. The left-to-right inclusion holds by adjunction. For the
other direction, we reason by contraposition, and assume $(a,b) \not\in \phi$.
Denote by $p \colon L \onto L/{\phi}$ the lattice quotient by the congruence
$\phi$. Then $p(a) \neq p(b)$, and we assume without loss of generality that
$p(a) \nleq p(b)$. By Theorem~\ref{thm:DPF} applied to $L/{\phi}$, the filter
${\uparrow}p(a)$ and the ideal ${\downarrow}p(b)$, we can pick a prime filter $G$ in
$L/{\phi}$ containing $p(a)$ but not $p(b)$. Then the set 
\[ F_x := p^{-1}(G), \] 
the inverse image of $G$ under the
homomorphism $p$, is a prime filter of $L$ which contains
$a$ but not $b$. Moreover, for any $(c,d) \in \phi$, we have $c \in F_x$ if,
and only if, $p(c) \in G$, if, and only if, $d \in F_x$, since $p(c) = p(d)$ by
assumption. Thus, $x \in \sem{\phi}$, and we conclude that $(a,b) \not\in
\theta(\sem{\phi})$. \end{proof}
\begin{remark}
  An interesting feature of Proposition~\ref{prop:quotientlattice-subspace} is
that we can start, on either side of the duality, with an unstructured set. Given an  
arbitrary subset $R\subseteq L\times L$ or $S\subseteq X$, by applying the
Galois connection, we always end up with a structured subset, that is,
$\sem{R}$ is not just a subset of $X$, but it is a closed, hence Priestley, subspace of $X$, 
and $\theta(S)$ is not just a subset of $L\times L$ but it is a congruence of
$L$. 
\end{remark}

From Proposition~\ref{prop:quotientlattice-subspace}, we deduce the following theorem, which shows how to explicitly compute the subspace dual to a lattice quotient generated by some equations.
\nl{$\gen{R}$}{the lattice congruence generated by $R$, for $R$ a binary relation on a lattice $L$}{}
\begin{theorem}\label{thrm:quotientlattice-subspace}
  Let $L$ be a distributive lattice with dual Priestley space $X$, and let $R$ be a binary relation on $L$. Then the lattice congruence generated by $R$, $\gen{R}$, is equal to $\theta(\sem{R})$, and the Priestley dual space of $L/\gen{R}$ is order-homeomorphic to the closed subspace $\sem{R}$ of $X$.
\end{theorem}
\begin{proof}
  In general, for any adjunction $f \colon P \leftrightarrows Q \colon g$, for any $p \in P$, $gf(p)$ is the minimum of $\mathrm{im}(g) \cap {\uparrow}p$ (see Exercise~\ref{exe:adjunctions}.\ref{itm:minimumimage} in Chapter~\ref{ch:order}). In particular, using that $\mathrm{im}(\theta)$ consists of the congruences on $L$, $\theta(\sem{R})$ is the smallest congruence containing $R$. Let $Y$ denote the Priestley dual space of the quotient lattice $L / \gen{R}$. The dual of the quotient map $p \colon L \onto L / \gen{R}$ is the continuous order-preserving function $i \colon Y \to X$ which can be defined by $h_{i(y)} := h_y \circ p$, for every $y \in Y$. Note that $i$ is an order embedding: if $y \nleq y'$ in $Y$, then we can pick $b \in L / \gen{R}$ such that $y \in \widehat{b}$ and $y' \not\in \widehat{b}$, and since $p$ is surjective, we can pick $a \in L$ such that $b = p(a)$. Then $i(y) \in \widehat{a}$ and $i(y') \not\in \widehat{a}$, so $i(y) \nleq i(y')$. Finally, for any $x \in X$, we have that 
  \[ x \in \sem{\gen{R}} \text{ if, and only if, } \gen{R} \subseteq \theta(\{x\}) = \ker(h_x), \] 
  using the Galois connection between $\sem{-}$ and $\theta$. The latter holds if, and only if, there exists $y \in Y$ such that $h_x = h_y \circ p$. Therefore, the image of $i$ is equal to $\sem{\gen{R}} = \sem{\theta(\sem{R})} = \sem{R}$, using the general fact that $gfg = g$ for any adjunction $(f,g)$, see Exercise~\ref{exe:adjunctions}.\ref{itm:equalities}. 
  Thus, $i$ is a homeomorphism between $Y$ and the closed subspace $\sem{R}$, as required.
\end{proof}

A second theorem mirrors the previous one, but on the space side.
\begin{theorem}\label{thm:closure-via-duality}
  Let $X$ be a Priestley space with dual distributive lattice $L$, and let $S$ be a subset of $X$. Then the closure of $S$, $\overline{S}$, is equal to $\sem{\theta(S)}$, and the lattice dual to $\overline{S}$ is isomorphic to the quotient $L/\theta(S)$.
\end{theorem}
\begin{proof}
  Left as Exercise~\ref{exe:closure-via-duality}.
\end{proof}

We give a few examples of the duality between lattice quotients and closed subspaces. You are asked to verify some details in Exercises~\ref{exe:quot-subspace}~and~\ref{exe:BAcong}.

\begin{example}\label{exa:quot-subspace1}
  Let $B$ be the eight element Boolean algebra; denote its three atoms by $a$, $b$ and $c$, and
  its dual space by
  $X=\{x,y,z\}$, where $F_x={\uparrow}a$, $F_y={\uparrow}b$, and $F_z={\uparrow}c$. For the equation $b\approx a \vee b$, we get $\sem{b \approx a \vee b} = \{y,z\}$, and the generated congruence $\gen{b \approx a \vee b} = \theta(\sem{b \approx a \vee b})$ on $B$ is depicted in Figure~\ref{fig:quot-subspace1}.
\end{example}
  \begin{figure}[htp]
    \begin{center}
      \begin{tikzpicture}[node distance=15pt, label distance=1pt]
        \po{(0,-1)}

        \po{(-1,0)}
        \po{(1,0)}
        \po{(0,0)}

        \po{(-1,1)}
        \po{(0,1)}
        \po{(1,1)}

        \po{(0,2)}

        \draw[thick] (0,-1) -- (-1,0) -- (-1,1) -- (0,2) -- (0,1) -- (-1,0);
        \draw[thick] (0, -1) -- (0, 0) -- (1,1) -- ((1, 0) -- (0,-1);
        \draw[thick] (1,0) -- (0,1) -- (0, 2) -- (1,1);
        \draw[thick] (0,0) -- (-1,1);

        \node[label=left:{$a$}] at (-1,0) {};
        \node[label=left:{$a \vee b$}] at (-1,1) {};
        \node[label=right:{$b$}] at (0,0)  {};
        \node[label=right:{$c$}] at (1,0) {};

        \draw[thick, dashed, rotate around={315:(-.5,-.5)}] (-.5,-.5) ellipse (28pt and 8pt);
        \draw[thick, dashed, rotate around={315:(-.5,.5)}] (-.5,.5) ellipse (28pt and 8pt);
        \draw[thick, dashed, rotate around={315:(.5,.5)}] (.5,.5) ellipse (28pt and 8pt);
        \draw[thick, dashed, rotate around={315:(.5,1.5)}] (.5,1.5) ellipse (28pt and 8pt);
      \end{tikzpicture}
    \end{center}
    \caption{A congruence on the eight element Boolean algebra}
    \label{fig:quot-subspace1}
  \end{figure}
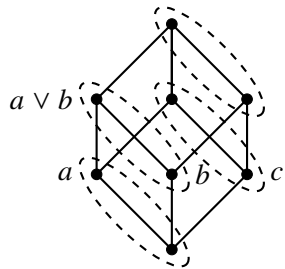
\begin{example}\label{exa:quot-subspace2}
  Let $X$ be the poset depicted on the left in Figure~\ref{fig:quot-subspace2} and let $S=\{x,y\}$.
  The dual lattice $\Down(X)$ and the congruence on $\Down(X)$ corresponding to the subspace 
  $S\subseteq X$ are as depicted on the right in the same figure. Furthermore, the subspace $S$ is 
  equal to $\sem{a\approx b}$, where $a$ is the down-set ${\downarrow} z$ and $b$ is the down-set ${\downarrow} x = \{x\}$. Note that the quotient of $\Down(X)$ by the congruence $\gen{a \approx b}$ is the three element chain, which is indeed the lattice dual to the induced poset on the subset $S \subseteq X$.
\end{example}
  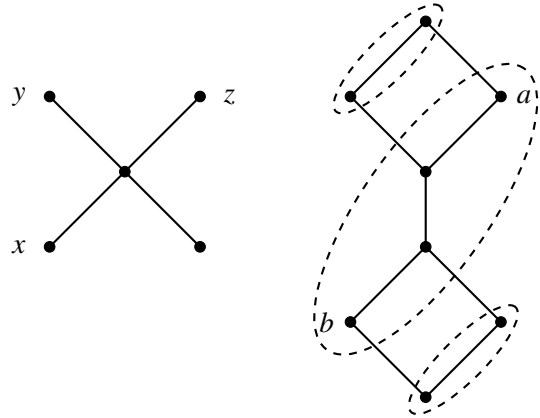
\begin{figure}[htp]
    \begin{center}
      \begin{tikzpicture}[node distance=15pt, label distance=1pt]
        \po{(-1,-1)}
        \po{(1,-1)}
        \po{(0,0)}
        \po{(-1,1)}
        \po{(1,1)}

        \draw[thick] (-1,-1)--(0,0)--(-1,1);
        \draw[thick] (1,-1)--(0,0)--(1,1);
        \node[label=left:{$x$}] at (-1,-1) {};
        \node[label=left:{$y$}] at (-1,1) {};
        \node[label=right:{$z$}] at (1,1)  {};
        \node[label=right:{$a$}] at (4.9,1) {};
        \node[label=left:{$b$}] at (3.1,-2) {};
        \po{(4,0)}
        \po{(3,1)}
        \po{(5,1)}
        \po{(4,2)}

        \po{(4,-1)}
        \po{(3,-2)}
        \po{(5,-2)}
        \po{(4,-3)}

        \draw[thick] (4,0)--(3,1)--(4,2)--(5,1)--(4,0);
        \draw[thick] (4,0)--(4,-1)--(3,-2)--(4,-3)--(5,-2)--(4,-1);

        \draw[thick, dashed, rotate around={45:(3.5,1.5)}] (3.5,1.5) ellipse (28pt and 8pt);
        \draw[thick, dashed, rotate around={45:(4.5,-2.5)}] (4.5,-2.5) ellipse (28pt and 8pt);
        \draw[thick, dashed, rotate around={55:(4,-0.5)}] (4,-0.5) ellipse (65pt and 24pt);

      \end{tikzpicture}
    \end{center}
    \caption{A poset, its down-set lattice, and the congruence corresponding to the subset $\{x,y\}$.  Here, $a$ corresponds to the down-set ${\downarrow}z$ and $b$ corresponds to the down-set $\{x\}$.}
    \label{fig:quot-subspace2}
  \end{figure}
\begin{example} \label{exa:remainderofbetaX}
  Recall the Stone-{\v C}ech compactification $\beta X$ of a set $X$, introduced in  Example~\ref{exa:StoneCechcompactification}. As remarked in Example~\ref{exa:StoneCechcompactification}, every point of $X$ is isolated in $\beta X$, so that $X$ is an \emph{open} subspace of $\beta X$. Recall also from that example that we identify points of $X$ with ultrafilters of $\cP(X)$. Note that the set ${}^*X := \beta X\setminus X$, known as the \emphind{remainder} of the Stone-{\v C}ech compactification $\beta X$, is closed. Elements of ${}^*X$ are known as \emph{free ultrafilters}\index{free ultrafilter}\index{ultrafilter!free}. Since ${}^*X$ is a closed subspace of the dual space of $\mathcal{P}(X)$, it corresponds to a quotient $B$ of the Boolean algebra $\mathcal{P}(X)$.
  We will show now that this quotient algebra $B$ is given by the following set of equations on $\cP(X)$: 
  \begin{equation} \label{eq:remainder-equations}
    \{ \{x\} \approx \emptyset \ \mid \ x \in X \} .
\end{equation}
  Indeed, if $\mu$ is a free ultrafilter and $x\in X$, then $\emptyset\not\in \mu$ and $\{x\}\not\in\mu$, so $\mu\models \{x\}\approx \emptyset$. On the other hand, if $\mu$ is principal, then there is $x\in X$ with $\{x\}\in\mu$, see Exercise~\ref{exe:StoneCech-principal}. Thus $\emptyset\not\in \mu$ and $\{x\}\in\mu$, so $\mu\not\models \{x\}\approx \emptyset$. We have shown that an ultrafilter $\mu$ of $\cP(X)$ satisfies all the equations in (\ref{eq:remainder-equations}) if, and only if, $\mu \in {}^* X$.

Write $\theta$ for the congruence on $\cP(X)$ generated by the set of equations (\ref{eq:remainder-equations}). Then $\theta$ identifies two elements $S$ and $T$ of $\cP(X)$ if, and only if, the symmetric difference $S+T$ of $S$ and $T$ is finite. Indeed, for any ultrafilter $\mu$ of $\cP(X)$, we have $\mu\models S\approx T$ if, and only if, $S+T\not\in\mu$ (Exercise~\ref{exe:BAcong}). Thus, $S$ and $T$ will be identified by $\theta$ if, and only if, $S+T$ does not belong to any free ultrafilter, which happens if, and only if, $S+T$ is finite.

The Boolean algebra quotient $B = \cP(X)/{\theta}$, which we identified in this example as the dual algebra of the closed subspace ${}^* X$ of $\beta X$, is often denoted $\mathcal{P}(X)/\mathrm{fin}$ in the literature, and is well studied in set theory and general topology.
\end{example}

\subsection*{Sublattices and Priestley quotient spaces}
A very similar story to the one above can be told for sublattices and quotients of Priestley spaces that are themselves Priestley spaces. In fact, the ensuing notion of Priestley space (in)equations is a very important tool in the theory of automata and regular languages, as we will see in Chapter~\ref{ch:AutThry}. We will now give the necessary definitions and statements of the relevant theorems for a general duality between sublattices and Priestley quotient spaces. We also note that, in fact, the lattice-quotient--closed-subspace duality from the previous subsection and the sublattice--Priestley-quotient-space duality are, in a sense, dual to each other.

We first introduce the notion of \emph{quotient space} in the context of Priestley spaces. This is an instance of a more general notion of quotient of ordered topological spaces, but we only need it in this setting.
\begin{definition}\label{def:compatible}
  A preorder $\preceq$ on a Priestley space $(X, \leq_X, \tau)$ is \emph{compatible} if ${\leq_X} \subseteq {\preceq}$ and, for any $x, y \in X$, if $x \not\preceq y$, then there exists a $\tau$-clopen $\preceq$-down-set $K$ in $X$ such that $y \in K$ and $x \not\in K$.\index{compatible!preorder}
\end{definition}
Compatible preorders on $X$ give an intrinsic description of those Priestley
spaces $Y$ that are quotients of $X$, that is, for which there exists a continuous
order-preserving surjective map $X \onto Y$. We just indicate here what this
means, and ask you to fill in the details in the exercises, see Exercises~\ref{exe:compatible-from-quotient}~and~\ref{exe:compatible-to-quotient}. 
First, for any $p
\colon X \onto Y$, the preorder $\preceq_p$ on $X$ defined by $x \preceq_p x'$ if, and only
if, $p(x) \leq_Y p(x')$ is compatible. 
Conversely, if $\preceq$ is a compatible
preorder on $X$, denote by $Y$ the poset reflection of the preordered set
$(X,\preceq)$, as defined in Exercise~\ref{exe:reflection}. 
That is, as a set,  $Y$ is $X/{\equiv}$ where ${\equiv} :=
{\preceq \cap \succeq}$. Let us write $q \colon X \onto Y$ for the quotient map.
Then the partial order $\leq_Y$ on the poset reflection is defined, for any $y =
q(x)$ and $y' = q(x')$ in $Y$, by $y \leq_Y y'$ if, and only if, $x \preceq x'$. The \emph{quotient topology} on $Y$ is defined by $\tau_Y := \{U \subseteq Y \ \mid \ q^{-1}(U) \in \tau_X\}$. Then the ordered topological space $(Y, \leq_Y, \tau_Y)$ is a Priestley space, and the function $q$ is a continuous order-preserving map. Moreover, for any Priestley space $(Z, \leq_Z, \tau_Z)$, any continuous order-preserving map $f \colon (X,\preceq,\tau_X) \to (Z, \leq_Z, \tau_Z)$ factors uniquely through the map $q$. We will say that $(Y, \leq_Y, \tau_Y)$ is the \emph{quotient of the Priestley space $X$} by the compatible preorder $\preceq$ and will denote it by $X/{\preceq}$.\index{quotient!of a Priestley space}\index{compatible preorder!quotient by}

\begin{example}\label{exa:Priestley order as a compatible order}
  Consider a Priestley space $(X, \tau,\leq)$ with dual lattice $L$. Notice (see Exercise~\ref{exe:Booleanspace}) that then $(X,\tau,=)$ is also a Priestley space. In fact, the dual lattice of this Priestley space is the Boolean algebra $L^-$ of all clopen subsets of $X$, also known as the Boolean envelope of $L$, see Corollary~\ref{cor:doubledualBA}. Of course, $L$ is a sublattice of $L^-$. The dual of the inclusion of $L$ into $L^-$ is simply the identity map on $X$, viewed as a continuous order-preserving map $\id_X\colon (X,\tau,=)\to(X, \tau,\leq)$. The corresponding compatible preorder is the partial order $\leq$ of $(X, \tau,\leq)$.

  More generally, any injective homomorphism $i \colon L \into B$, with $B$ a Boolean algebra, factors as the composition of $e \colon L \into L^-$ and $\bar{i} \colon L^- \into B$. Dually, denoting by $(Y,\pi)$ the Boolean space dual to $B$, this gives a quotient map of Boolean spaces $f \colon (Y,\pi) \to (X,\tau)$, followed by $\id_X \colon (X, \tau, =) \to (X,\tau, \leq)$. This shows that sublattices of Boolean algebras can be understood dually by a Boolean equivalence relation together with a Priestley order on the quotient.
\end{example}

As before, the correspondence between compatible preorders and sublattices allows us to view pairs of elements $(x,y)$ from $X$ as constraints on $L$ yielding sublattices. However, as Example~\ref{exa:Priestley order as a compatible order} already shows, equating elements of the space is not a fine enough notion to witness all sublattices. In fact, as we will see, spatial equations witness Boolean subalgebras, whereas we will need to think of pairs $(x,y)$ as spatial \emph{inequations} in order to be able to witness all sublattices. For this reason, we introduce the following notation.

\begin{notationnum}\label{notation:space-inequation}
Let $X$ be a set. We denote by $I(X)$ the set of all formal inequations $x\preceq x'$ where $x,x'\in X$. That is,
\[
I(X)=\{x\preceq y\mid x,y\in X\}.
\]
Given a set 
of inequations $T\subseteq I(X)$, we call $R_T := \{(x,y)\mid x\preceq y\in T\}$ the \emph{binary relation corresponding to} $T$. 
\end{notationnum}
Notice that, up to isomorphism, $I(X)$ is just $X^2$, and subsets $T\subseteq I(X)$ are just binary 
relations on $X$. However, when studying quotients of spaces, 
we need to be a bit more careful with our notation than in the analogous setting of Notation~\ref{notation:lattice-equation} for quotients of lattices, 
because we here need to carefully distinguish
whether we interpret a pair $(x,y)$ of points in a Priestley space 
as an \emph{inequation} or as an \emph{equation}. 
We begin with the first interpretation, for the second
interpretation, see Notation~\ref{notation:space-equation} below.

\begin{definition}\label{dfn:spatial-inequations}
Let $L$ be a distributive lattice with dual Priestley space $X$. We define the binary relation $\models\ \subseteq L\times I(X)$ by, for any element $a \in L$ and $x,y\in X$,
\nl{$a \models x \preceq y$}{satisfaction relation between an element $a$ of a lattice and a pair $(x,y)$ of elements of the dual space viewed as an inequation}{}
\[
a\ \models\ x\preceq y \ \iffdef\ y\in\widehat{a} \text{ implies }x\in \widehat{a}.
\]
When $a\ \models\ x\preceq y$, we will say that $a$ \emphind{satisfies the (spatial) inequation} $x\preceq y$. Similarly, for a subset $A\subseteq L$, we write $A\ \models\ x\preceq y$ and say that $A$ satisfies $x\preceq y$ provided every $a\in A$ satisfies $x\preceq y$.

This relation allows us to define a map $\preceq_{-}\colon\cP(L)\rightarrow \cP(X^2)$,  by, for $a\in L$ and $x,y\in X$,
\nl{$\preceq_a$}{preorder on the dual space associated with an element $a$ of a lattice}{}
\[ 
x \preceq_a y \iffdef a\ \models\ x\preceq y,
\]
and for any subset $A \subseteq L$, ${\preceq_{A}} := \bigcap_{a \in A} {\preceq_a}$, that is,
\nl{$\preceq_A$}{preorder on the dual space associated with a subset $A$ of a lattice}{}
\[
 x \preceq_A y \iffdef \ A\ \models\ x\preceq y.
 \]
We obtain a map $\semI{-}\colon\cP(X^2)\rightarrow\cP(L)$,  by, for $x,y\in X$
  \[
   \semI{(x,y)} := I_y \cup F_x = \{a \in L \ \mid \ a\ \models\ x\preceq y \},\]
\nl{$\semI{(x,y)}$}{sublattice associated with a pair of elements of the dual space viewed as an inequation}
and for any  $R\subseteq X^2$, we define the subset $\semI{R}$ of $L$ by
  \[
    \semI{R} := \bigcap_{(x,y) \in R} \semI{(x,y)} = \{a \in L \, \mid \, \text{ for all } (x,y) \in R, \ a\ \models\ x\preceq y \}.
  \]
\nl{$\semI{R}$}{sublattice associated with a binary relation on the dual space viewed as a set of inequations}
In order to be able to talk about inequations and sets of inequations directly, we will also use,
for any $(x,y) \in X^2$, the notation $\sem{x\preceq y}$
for $\semI{(x,y)}$, and, for any $T \subseteq I(X)$, the notation $\sem{T}$ for the set $\semI{R_T}$, where $R_T$ is the binary relation corresponding to $T$. 
\end{definition}

\begin{proposition}\label{prop:sublattice-quotientspace}
  Let $L$ be a distributive lattice with dual Priestley space $X$.  The two functions $\preceq_{-} \colon \mathcal{P}(L) \leftrightarrows \mathcal{P}(X^2)\colon \semI{-}$ form a Galois connection, whose fixed points on the left are the sublattices of $L$, and whose fixed points on the right are the compatible preorders on $X$.
\end{proposition}
\begin{proof} 
The functions indeed form a Galois connection since they are the contravariant adjoint pair between the powersets of $L$ and $X^2$, formed as in Example~\ref{exa:galoisconnection}, in the case of the relation $\models$.

Given $x, y\in X$, it is easy to verify that the set $\sem{x\preceq y}$ is a sublattice of $L$ (see Exercise~\ref{exe:sublattice-quotient-adjunction}.a), and thus so is $\semI{R}$, being an intersection of sublattices, for any $R\subseteq X^2$.
Now let $K$ be a sublattice of $L$. We want to show $K=\semI{\preceq_K}$. By adjunction, we have $K\subseteq\semI{\preceq_K}$.  For the converse inclusion, let $a\in\semI{\preceq_K}$ be arbitrary. 
Note that,  for any $(x,y) \in \widehat{a}^c \times \widehat{a}$, we have $a \not\models x \preceq y$, so $x \not\preceq_K y$, which means that we can pick $a_{(x,y)}\in K$ with $x\in\widehat{a_{(x,y)}}^c$ and $y\in\widehat{a_{(x,y)}}$. Now, for any fixed $x\in\widehat{a}^c$, we get 
\[
\widehat{a}\subseteq\bigcup\{\widehat{a_{(x,y)}}\mid y\in\widehat{a} \}
\]
Applying compactness, the union of a finite subcover yields an element $b_{x}\in K$ with $a\leq b_{x}$ and $x\not\in\widehat{b}_x$. Now, $\widehat{a}=\bigcap \{\widehat{b}_{x}\mid x\in\widehat{a}^c\}$ and, again by compactness, there is a finite set $F\subseteq \widehat{a}^c$ with 
\[
\widehat{a}=\bigcap_{x\in F}\widehat{b}_{x}. 
\]
That is, $a=\bigwedge_{x\in F} b_{x},$ which is an element of $K$.

Given $a\in L$, the relation ${\preceq_a} = (\widehat{a} \times X) \cup (X \times \widehat{a}^c)$ is clearly a preorder which contains $\leq_X$. Also, it is compatible since, if $x\not\preceq_a x'$, then $(x,x')\in \widehat{a}^c\times \widehat{a}$. That is, $x'\in\widehat{a}$ while $x\not\in\widehat{a}$ and the set $\widehat{a}$ is a $\preceq_a$-down-set
(see Exercise~\ref{exe:sublattice-quotient-adjunction}.b). Further, for $A\subseteq L$, $\preceq_A$ is again a compatible preorder since these are closed under arbitrary intersections (see Exercise~\ref{exe:sublattice-quotient-adjunction}.c).

Now, to see that all compatible preorders are fixed points of the Galois collection, let $\preceq$ be a compatible preorder on $X$. Again, we have 
$\preceq\ \subseteq \ \preceq_{\semI{\preceq}}$ by adjunction.
For the reverse inclusion, suppose $x\not\preceq y$. Since $\preceq$ is compatible, we have a clopen $V\subseteq X$ with $y\in V$, $x\not\in V$, and $V$ is a $\preceq$-down-set. Since ${\leq_X} \subseteq {\preceq}$, $V$ is also a $\leq_X$-down-set, and thus, since $V$ is also clopen, we can pick $a\in L$ with $V=\widehat{a}$. The fact that $\widehat{a}$ is a $\preceq$-down-set precisely says that  $a\in\semI{\preceq}$. Also, the fact that $y\in \widehat{a}$ and $x\not\in\widehat{a}$ means that $x\not\preceq_a y$ and thus $(x,y)\not\in{\preceq_{\semI{\preceq}}}$, as required. 
\end{proof}

In the same way as for quotient--subspace duality, we now get the following duality between sublattices and preorders on its dual space; compare Theorems~\ref{thrm:quotientlattice-subspace}~and~\ref{thm:closure-via-duality}.

\begin{theorem}\parencite{Sch02}\label{thrm:sublattices-inequations}
  Let $L$ be a distributive lattice with dual Priestley space $X$.
  \begin{enumerate}
    \item For any subset $A$ of $L$, the sublattice generated by $A$ is equal to $\semI{\preceq_A}$, and the Priestley dual space of this sublattice is order-homeomorphic to the Priestley quotient $X/{\preceq_A }$.
    \item For any binary relation $R$ on $X$, the smallest compatible preorder containing $R$ is equal to $\preceq_{\semI{R}}$, and the lattice dual to the Priestley quotient by this preorder is isomorphic to $\semI{R}$.
  \end{enumerate}
\end{theorem}

Theorem~\ref{thrm:sublattices-inequations} allows one to characterize sublattices of a lattice $L$ by sets of inequations over the dual space of $L$. As we will see in Section~\ref{sec:EilReittheory}, this is a method used in automata theory to prove the decidability of the membership problem for classes of regular languages. Apart from characterizing sublattices of a lattice $L$, Theorem~\ref{thrm:sublattices-inequations} can also be used to \emph{separate} sublattices, as identified in the following corollary.

\begin{corollary}
 Let $L$ be a distributive lattice with dual Priestley space $X$ and let $K_1$ and $K_2$ be two sublattices of $L$. Then $K_2\not\subseteq K_1$ if and only if there exists an inequation $x\preceq y$ over $X$ and $a\in K_2$ so that $a\not\models x\preceq y$ but $K_1\models x\preceq y$.
\end{corollary}

We now spell out the specialization of the subalgebra--quotient-space duality
to the Boolean case. Recall from Proposition~\ref{prop:boolean-trivial-order}that a Priestley space corresponds to a Boolean
algebra if, and only if, the order of the space is trivial. 
It follows that a compatible preorder $\preceq$ 
yields a Boolean quotient if, and only if, $\preceq$ is an equivalence relation or, equivalently,
if it is symmetric. We call symmetric compatible
preorders \emph{compatible equivalence 
relations},\index{compatible!equivalence relation} and, 
in situations where we are interested in Boolean subalgebras, we use the stronger notion 
of spatial \emph{equation} that we define now. 

\begin{notationnum}\label{notation:space-equation}
 For any set $X$, we denote by $E(X)$ the set of all formal equations $x\equiv y$ where $x,y\in X$. That is,
\[
E(X)=\{x\equiv y\mid x,y\in X\}.
\]
As before, when $T \subseteq E(X)$ is a set of formal equations, we denote by $R_T \subseteq X^2$ the corresponding
binary relation on $X$.
\end{notationnum}
Analogously to what we did above, we define a notion of satisfaction for formal equations on a space.
\begin{definition}\label{dfn:spatial-equations}
 Let $B$ be a Boolean algebra and $X$ its Priestley dual space. We define the binary relation $\models\ \subseteq B\times E(X)$ by, for any element $a \in B$ and $x,y\in X$ 

\nl{$a \models x \equiv y$}{satisfaction relation between an element $a$ of a Boolean algebra and a pair $(x,y)$ of elements of the dual space, viewed as an equation}{}
\[
a\ \models\ x\equiv y \ \iff\ x\in\widehat{a} \text{ if, and only if, }y\in \widehat{a}.
\]
When $a\ \models\ x\equiv y$, we will say that $a$ \emphind{satisfies the (spatial) equation} $x\equiv y$. Similarly, for a subset $A\subseteq B$, we write $A\ \models\ x\equiv y$ and say that $A$ satisfies $x\equiv y$ provided every $a\in A$ satisfies $x\equiv y$.
\nl{$A \models x \equiv y$}{satisfaction relation between a subset $A$ of a Boolean algebra and a pair $(x,y)$ of elements of the dual space, viewed as an equation}{}

This relation $\models$ allows us to define a map $\equiv_{-}\colon\cP(B)\rightarrow \cP(X^2)$,  by, for $a\in B$ and $x,y\in X$,
\nl{$\equiv_a$}{equivalence relation on the dual space associated with an element $a$ of a Boolean algebra}{}
\[ 
x \equiv_a y \iff a\ \models\ x\equiv y,
\]
and for any subset $A \subseteq B$, ${\equiv_{A}} := \bigcap_{a \in A} {\equiv_a}$, that is,
\[
 x \equiv_A y \iff \ A\ \models\ x\equiv y.
 \]
 \nl{$\equiv_A$}{equivalence relation on the dual space associated with a subset $A$ of a Boolean algebra}
Similarly we obtain a map $\semE{-}\colon\cP(X^2)\rightarrow\cP(B)$,  by, for $x,y\in X$
\nl{$\semE{(x,y)}$}{Boolean subalgebra associated with a pair of elements of the dual space viewed as an equation}{}
  \[
   \semE{(x,y)} := (I_x\cap I_y) \cup (F_x\cap F_y) = \{a \in B \ \mid  x \in \widehat{a}, \text{ if, and only if, } y \in \widehat{a} \},\]
and for any  $R\subseteq X^2$, we define the subset $\semE{R}$ of $B$ by
\nl{$\semE{R}$}{Boolean subalgebra associated with a binary relation $R$ on the dual space viewed as a set of  equations}{}
  \[
    \semE{R} := \bigcap_{(x,y) \in R} \semE{(x,y)} = \{a \in B \, \mid \, \text{ for all } (x,y) \in R, \ a\ \models\ x\equiv y \}.
  \]
In order to be able to talk about equations directly, we will also use the notation 
$\sem{x\equiv y}$ for $\semE{(x,y)}$ and, for any subset $T \subseteq E(X)$, we write
$\sem{T}$ for $\semE{R_T}$, where $R_T$ is the binary relation corresponding to $T$.
\end{definition}
 
 \begin{proposition}\label{prop:subBalg-quotientspace}
  Let $B$ be a Boolean algebra with dual Priestley space $X$.  The two functions 
  $\equiv_{-} \colon \mathcal{P}(B) \leftrightarrows \mathcal{P}(X^2)\colon \semE{-}$ 
  form a Galois connection, whose fixed points on the left are the Boolean subalgebras 
  of $B$, and whose fixed points on the right are the compatible equivalence relations 
  on $X$.
\end{proposition}
\begin{proof} 
See Exercise~\ref{exer:subBalg-quotientspace}.
\end{proof}

\begin{corollary}\label{cor:subalgebras-equations}
Let $B$ be a Boolean algebra with dual space $X$. 
\begin{enumerate}
	\item For any subset $A$ of $B$, the Boolean subalgebra generated by $A$ is equal to
	$\semE{\equiv_A}$, and the dual space of this subalgebra is homeomorphic to the 
	quotient space $X/{\equiv_A}$.
\item For any binary relation $R$ on $X$, the smallest compatible equivalence
	relation containing $R$ is equal to $\equiv_{\semE{R}}$, 
	and the Boolean algebra dual to the quotient by this equivalence relation is isomorphic
	to $\semE{R}$.
\end{enumerate}
\end{corollary}
\begin{proof}
Deduced from Proposition~\ref{prop:subBalg-quotientspace} in a similar way to Theorems~\ref{thrm:quotientlattice-subspace}~and~\ref{thm:closure-via-duality}.
\end{proof}

Note that, in part (a) of Corollary~\ref{cor:subalgebras-equations}, we could have used $\semI{-}$ instead of $\semE{-}$ since, when applying $\semI{-}$ to a compatible \emph{equivalence} relation, 
the corresponding lattice is a Boolean algebra. 
Further, $\semI{\equiv_A}$ and $\semE{\equiv_A}$ are also equal to $\semE{\preceq_A}$. 
Similarly, in part (b) of Corollary~\ref{cor:subalgebras-equations}, since the lattice $\semI{R}$ generates the Boolean algebra $\semE{R}$, the same compatible equivalence relation corresponds to both of these subsets of $B$. Finally,  $\equiv_{\semI{R}}$ and $\equiv_{\semE{R}}$ are both equal to $\preceq_{\semE{R}}$, since an inequality is satisfied by a Boolean subalgebra if, and only if, its opposite is also.

\begin{remark}\label{rem:inequations}
  One problem with compatible preorders, which has sometimes hampered their 
  successful application, is that it is difficult to understand intrinsically in a Priestley 
  space how to obtain the compatible preorder generated by a binary relation on 
  the space. The existence of the Galois connection of 
  Proposition~\ref{prop:sublattice-quotientspace} frees us from this problem and 
  allows us to specify sublattices from arbitrary sets $E$ of inequality constraints 
  of the form $x\preceq y$.

  This is similar to the fact that, in logic, we do not need to identify the full theory 
  of a class of structures that we are interested in, since we may be able to capture 
  it by a much smaller set of axioms.
\end{remark}

As mentioned above, we will see dual space (in)equations in action in 
Chapter~\ref{ch:AutThry}. Below, we will already give two elementary examples.

\begin{example}[Sublattices of the free distributive lattice]\label{exa:sub-of-free}
  Coming back to the example of the free distributive lattice on two generators, 
  Example~\ref{exa:twogen-logic}, consider the sublattices $L_1 = \{\bot, p, \top\}$ 
  and $L_2 = \{\bot, p \vee q, \top\}$. One may prove from the definitions that the 
  corresponding quotients of $2^{\{p,q\}}$ are given by the equivalence relations 
  $\equiv_1$ and $\equiv_2$, where $\equiv_1$ has classes $\{00, 01\}$ and 
  $\{10,11\}$, while $\equiv_2$ has classes $\{00\}$ and $\{01,10,11\}$. We leave 
  it as Exercise~\ref{exe:subs-of-freedl-two} to classify the other sublattices of this 
  lattice.
\end{example}

\begin{example}[Equations for a subalgebra]\label{exa:equationssubalg}
  In this example we will see how we can use extraneous structure on a dual
  space and lattice to identify a smaller set of equations for a subalgebra of
  an infinite Boolean algebra; this technique is further exploited in
  Section~\ref{sec:EilReittheory}.

  Consider the set $\bZ$ of integers, and denote by $\bZ^+$ the subset of positive 
  integers and by $\bZ^-$ the subset of negative integers; so 
  $\bZ = \bZ^- \cup \{0\} \cup \bZ^+$. Let $M$ be the Boolean subalgebra of 
  $\cP(\bZ)$ consisting of all those subsets $S$ of $\bZ$ such that \emph{both} 
  $S\cap\bZ^+$ is either finite or co-finite, \emph{and} $S\cap\bZ^-$ is either 
  finite or co-finite. One may then show (see Exercise~\ref{exe:equationssubalg}) 
  that the dual space of $M$ is the `two-point compactification of $\bZ$'
  \[
    \bZ_{-\infty}^{+\infty} := \bZ\cup\{-\infty,+\infty\},
  \]
  which topologically is the disjoint union of the one-point compactification 
  $\bZ^+\cup\{+\infty\}$ of $\bZ^+$ with the discrete topology, the one-point 
  compactification $\bZ^-\cup\{-\infty\}$ of $\bZ^-$ 
  with the discrete topology, and the one point space $\{0\}$. (For the one-point compactification, see 
  Example~\ref{exa:fincofsubsetsdual} and Exercise~\ref{exe:fincofsubsetsspace}.) 
  Note that, since $\bZ$ is the disjoint union of $\bZ^+$, $\bZ^-$, and $\{0\}$, every free ultrafilter of $\cP(\bZ)$ contains exactly one of $\bZ^+$ or $\bZ^-$, because a free ultrafilter clearly cannot contain $\{0\}$. The dual of the inclusion 
  $M\hookrightarrow\cP(\bZ)$ is the surjective function
  \begin{align*}
    \beta(\bZ) & \onto \bZ_{-\infty}^{+\infty}, \\
    \mu        & \mapsto\left\{
    \begin{array}{lll}
      k       & \text{ if } & \{k\}\in\mu \text{ where } k \in \bZ, \\
      +\infty & \text{ if } & \mu \text{ free and } \bZ^+\in\mu,                             \\
      -\infty & \text{ if } & \mu \text{ free and } \bZ^-\in\mu.               %
    \end{array}
    \right.
  \end{align*}
  Thus, the compatible preorder on $\beta\bZ$ corresponding to the subalgebra $M$ 
  of $\cP(\bZ)$ is the equivalence relation in which each $k\in \bZ$ is only related to 
  itself, and two free ultrafilters $\mu$ and $\nu$ are related provided they either both 
  contain $\bZ^+$, or both contain $\bZ^-$.
  That is, the remainder is split into two uncountable equivalence classes and each 
  free ultrafilter is related to uncountably many other free ultrafilters.

  By contrast, we will now show that, by using the successor structure on $\bZ$, the 
  subalgebra $M$ can be `axiomatized' by a much `thinner' set of equations. 
  The successor function $\bZ \to \bZ$, which sends $k \in \bZ$ to $k + 1$, has as its discrete dual
  the complete homomorphism 
  $\cP(\bZ) \to \cP(\bZ)$ which sends a subset $S \in \cP(\bZ)$ to the set 
  \[ S - 1 := \{ k \in \bZ \ \mid \ k + 1 \in S \} = \{ s - 1 \ \mid \ s \in S \}\ .\] 
  For 
  $\mu \in \beta\bZ$, write
  \[
    \mu+1:=\{S\in\cP(\bZ)\mid S-1\in\mu\}=\{S+1\mid S\in\mu\} .
  \]
  This is a well-defined function $\beta \bZ \to \beta \bZ$, as it is the Priestley dual of the homomorphism $S\mapsto S-1$ on $\cP(\bZ)$; it is also the unique continuous extension of
  the successor function $\bZ \to \bZ$, when we view the codomain as a subset of $\beta \bZ$.
  To describe the equational basis for the sublattice $M$ of $\cP(\bZ)$, 
  consider the set of equations $\mu + 1 \approx \mu$, as $\mu$ ranges over ${}^*\bZ$, 
  where we recall that ${}^*\bZ := \beta\bZ \setminus \bZ$, the \emph{remainder} 
  of $\beta\bZ$, see Example~\ref{exa:remainderofbetaX}.
  We will show that the sublattice $M$ of $\cP(Z)$ contains exactly those $S \in \cP(\bZ)$ that 
  satisfy all of these equations,  that is, we will prove that  
 \begin{equation}\label{eq:complete-eqs-M} 
    M \quad = \quad \sem{\,\mu+1\approx\mu\mid \mu\in{}^*\bZ\,}.
 \end{equation}
  To this end, note first that, for $S \in \cP(\bZ)$, $S$ satisfies $\mu+1 \approx \mu$ if, 
  and only if, both $\mu$ and $\mu + 1$ contain $S$, or neither $\mu$ nor $\mu + 1$ contains $S$.
  Now, for the left-to-right inclusion of (\ref{eq:complete-eqs-M}), let $S \in M$ and let $\mu$ be a free ultrafilter of $\cP(\bZ)$. We show that $S\models 
  \mu+1\approx\mu$. Since $\mu$ is prime and $\bZ^+\cup (\bZ \setminus \bZ^+)=\bZ\in\mu$, 
  it follows that either $\bZ^+\in\mu$ or $\bZ \setminus \bZ^+ \in\mu$. We treat the case 
  $\bZ^+\in\mu$ and leave the other as an exercise. Since $S \in M$, we have that $S \cap \bZ^+$ is either finite or co-finite. If $S\cap\bZ^+$ is finite then, as 
  $\mu$ is free, $S\not\in\mu$. Also $S\cap\bZ^+$ finite implies that $(S-1)\cap\bZ^+$ is 
  finite and thus $S-1\not\in\mu$. So $S\models \mu+1\approx\mu$. If on the other hand 
  $S\cap\bZ^+$ is co-finite, then, as $(\bZ^+\setminus S) \cup (S\cap\bZ^+)=\bZ^+\in\mu$,
  and $\mu$ is free, it follows that $S\cap\bZ^+\in\mu$. Furthermore, $S\cap\bZ^+$ co-finite 
  implies that $(S-1)\cap\bZ^+$ is also co-finite and by the same argument we have 
  $S-1\in\mu$ so that $S\models \mu+1\approx\mu$. 
  For the right-to-left inclusion of (\ref{eq:complete-eqs-M}), we reason contrapositively, and suppose $S\not\in M$. Then $S\cap\bZ^+$ is neither finite
  nor co-finite, or $S\cap\bZ^-$ is neither finite nor co-finite. Again, we
  treat the first case and leave the second as an exercise. If $S\cap\bZ^+$ is
  neither finite nor co-finite it follows that there is an infinite set $T\subseteq\bZ^+$ 
  such that, for each $k\in T$
  \[
    k\not\in S \quad\text{ but }\quad k+1\in S.
  \]
  By the Prime Filter-Ideal Theorem~\ref{thm:DPF}, here applied to the Boolean 
  algebra $\cP(\bZ)$, pick an ultrafilter  $\mu$ of $\cP(\bZ)$ which contains the 
  filter ${\uparrow}T$ and is disjoint from the ideal $I$ consisting of all finite subsets 
  of $\cP(\bZ)$. Since $\mu$ is disjoint from $I$, it is free, and since 
  ${\uparrow}T\subseteq\mu$ we have $T\in\mu$. Now as $S\cap T=\emptyset$ 
  it follows that $S\not\in\mu$ and since $T+1\subseteq S$, or equivalently, 
  $T\subseteq S-1$, it follows that $S-1\in\mu$. That is, we have exhibited a free 
  ultrafilter $\mu$ such that  $S\not\models \mu+1\approx\mu$ and this completes 
  the proof of (\ref{eq:complete-eqs-M}).
\end{example}

We note that Example~\ref{exa:equationssubalg} is related to a well-known 
language from descriptive complexity theory. The monoid of integers under 
addition is the so-called syntactic monoid of the language called `majority', 
consisting of all bitstrings with a majority of $1$'s; see Exercise~\ref{exe:majority}. 

In Chapter~\ref{ch:AutThry} we will also need the duality between complete Boolean 
subalgebras of a powerset algebra and equivalence relations on the underlying set. 
This is a dicrete version of the duality between spatial equations and Boolean subalgebras 
stated in Proposition~\ref{prop:subBalg-quotientspace}. 
Like all other correspondences in this section, it is also possible to treat this duality
via the Galois connection obtained from a relation 
${\models}\subseteq\cP(S)\times\cP(S\times S)$ defined by $a\models s\equiv s'$ if, and 
only if, $s\in a\ \iff\ s'\in a$, but we here give an alternative, direct proof, in Theorem~\ref{thrm:BA-sub-discrete} below.

Recall that a \emph{complete Boolean subalgebra}\index{subalgebra!of a complete
Boolean algebra} of a complete Boolean algebra $B$ is a Boolean subalgebra $A$
such that, for any $S \subseteq A$, the supremum and the infimum of $S$ in $B$ 
are in $A$. Note that, using negation, it suffices to assume this for suprema \emph{or}
for all infima. Also note that for $A$ to be a complete Boolean subalgebra it is 
\emph{not} sufficient to 
merely require that $A$ be a Boolean algebra which is complete in its order
(see Exercise~\ref{exe:counterexamples-complete}.\ref{ite:complete-sublattice-not-sub}).
In Theorem~\ref{thrm:BA-sub-discrete} below, we will consider the poset $\mathrm{Sub}_c(B)$ of
complete Boolean subalgebras ordered by inclusion, and the poset $\mathrm{Eq}(S)$
of equivalence relations on a set $S$, also ordered by inclusion. Both are 
complete lattices. 
For an equivalence relation ${\equiv}$ on a set $S$, we will call a subset $u$
of $S$ \emphind{invariant} provided $u\models s\equiv s'$ for all $s, s' \in S$ and
we write $\sem{\equiv}$ for the collection of invariant subsets of $S$. 
\nl{$\sem{\equiv}$}{Boolean subalgebra of $\cP(S)$ associated with an equivalence relation $\equiv$ on a set $S$}{}
\begin{theorem}\label{thrm:BA-sub-discrete}
  Let $S$ be a set and let $B := \cP(S)$ be its Boolean algebra of subsets. For any 
  equivalence relation ${\equiv}$ on $S$, the set $\sem{\equiv}$ is a complete Boolean 
  subalgebra of $B$. Moreover, the assignment ${\equiv} \mapsto \sem{\equiv}$ is an 
  anti-isomorphism between $\mathrm{Eq}(S)$ and $\mathrm{Sub}_c(B)$.
\end{theorem}
\begin{proof}
  Let $f \colon S \to S/{\equiv}$ be the quotient function associated to ${\equiv}$. Note 
  that, essentially by definition, $\sem{\equiv}$ is exactly the image of the dual complete 
  homomorphism $f^{-1} \colon \mathcal{P}(S/{\equiv})\into B$, which is injective and the 
  upper adjoint of the forward image function $f[-]\colon\cP(S)\to\cP(S/{\equiv})$, see
  Exercise~\ref{exe:discrete-duality}. 
  Conversely, given a complete Boolean subalgebra $A$ of $B$, the inclusion 
  $i\colon A\hookrightarrow B$ is an injective complete homomorphism and thus, it has
  a lower adjoint $i'\colon B\to A$ which is surjective. As in the finite case, see the proof
  of  Lemma~\ref{lem:loweradj-preserves-joinprime}, $i'$ sends completely join irreducibles
  to completely join irredicibles, and thus atoms to atoms. That is, $i'$ restricts to a function
  $f\colon B\to At(A)$. Note that since $i'$ preserves arbitrary joins, $i'$ is, up to isomorphism, 
  forward image under $f$, so $i'$ is surjective if, and only if, $f$ is, and thus $f$ is the quotient  
  function associated to a uniquely determined equivalence relation on $S$. Clearly these two   
  assignments are order reversing and inverse to each other.
\end{proof}

\ourexercises

\begin{ourexercise}\label{exe:closure-via-duality}
  Prove Theorem~\ref{thm:closure-via-duality}.
\end{ourexercise}

\begin{ourexercise}
  Consider the unique homomorphism $h \colon F_{\DL}(p,q) \to F_{\BA}(p)$ that sends $p$ to $p$ and $q$ to $\neg p$. Show that $h$ is surjective, and compute the dual injective function.
\end{ourexercise}

\begin{ourexercise}\label{exe:quot-subspace}
  Prove the statements made in Examples~\ref{exa:quot-subspace1} and~\ref{exa:quot-subspace2}.
\end{ourexercise}

\begin{ourexercise}
\label{exe:BAcong}
Let $\theta$ be a congruence on a Boolean algebra $B$.
\begin{enumerate}
  \item Prove that, for any $a, b \in B$, $a \mathrel{\theta} b$ if, and only if, $a + b \mathrel{\theta} 0$. Here, $a + b$ denotes the \emphind{symmetric difference} of $a$ and $b$, that is, $a + b := (a \wedge \neg b) \vee (\neg a \wedge b)$.
  \item Prove that, for any ultrafilter $\mu$ of $B$, we have $\mu \models a \approx b$ if, and only if, $a + b \not\in \mu$.
  \item Conclude that, for any binary relation $R$ on $B$, we have $x \in \sem{R}$ if, and only if, $x \not\in \widehat{a + b}$ for all $(a,b) \in R$.
\end{enumerate}
\end{ourexercise}

\begin{ourexercise}\label{exe:compatible-from-quotient}
  Let $p \colon X \to Y$ be a continuous order-preserving map between Priestley spaces. Prove that the relation ${\preceq} := p^{-1}({\leq_Y})$ on $X$ is a compatible preorder.
\end{ourexercise}

\begin{ourexercise}\label{exe:compatible-to-quotient}
  Let ${(X,\tau_X,\leq_X)}$ be a Priestley space and let 
  $\preceq$ be a compatible preorder on $X$. Define ${\equiv} := {\preceq \cap \succeq}$, let $Y := X/{\equiv}$, and denote by $q \colon X \onto Y$ the quotient map.
  \begin{enumerate}
    \item Show that, for any $x, x' \in X$, $x \leq_X x'$ implies $x \preceq x'$.
    \item Prove that the relation $\leq_Y$ on $Y$ defined, for $y = q(x)$ and
	    $y' = q(x')$ in $Y$, by $y \leq_Y y'$ if, and only if, $x \preceq x'$, is a well-defined partial order.
    \item Prove that, with the quotient topology $\tau_Y := \{U \subseteq Y \ \mid \ q^{-1}(U) \in \tau_X\}$, $(Y,\tau_Y,\leq_Y)$ is a Priestley space.
    \item Prove that $q \colon X \to Y$ is continuous and order preserving.
    \item Prove that, for any Priestley space $Z$ and any continuous ${f \colon X \to Z}$ such that $x \preceq x'$ implies $f(x) \leq_Z f(x')$, there exists a unique continuous order-preserving $\bar{f} \colon Y \to Z$ such that $f = \bar{f} \circ q$.
  \end{enumerate}
\end{ourexercise}

\begin{ourexercise}
  Prove the assertion made in Example~\ref{exa:Priestley order as a compatible order}.
\end{ourexercise}

\begin{ourexercise}\label{exe:sublattice-quotient-adjunction}
  Let $L$ be a distributive lattice and $X$ its Priestley dual space.
  \begin{enumerate}
  \item For $x,x'\in X$, show that $\sem{x\preceq x'}$ is a sublattice of $L$.
  \item For each $a\in L$, show that $\preceq_a$ is a compatible preorder on $X$, 
  that the equivalence classes of the corresponding equivalence relation are 
  $\widehat{a}$ and $\widehat{a}^c$, and that $\widehat{a}\leq\widehat{a}^c$ is 
  the only non-trivial relation in the quotient order. Conclude that $\semI{\preceq_a}$ 
  is the three element sublattice $\{0,a,1\}$ of $L$.
  \item Show that compatible preorders on $X$ are closed under arbitrary 
  intersections. 
  \end{enumerate}
\end{ourexercise}

\begin{ourexercise}\phantomsection\label{exer:subBalg-quotientspace}
\begin{enumerate}
\item Verify that the functions in Proposition~\ref{prop:subBalg-quotientspace} constitute the 
Galois connection given by the relation $a\models x\equiv y$.
\item Verify that, for any $R\subseteq X^2$, $\semE{R}$ is closed under complementation 
and is thus a Boolean subalgebra of $B$.
\item Prove that for any subset $A\subseteq B$, ${\equiv_A}={\preceq_{A'}}$, where 
$A'=A\cup\{\neg a\mid a\in A\}$ and conclude that it is a compatible equivalence relation on $X$.
\item Using the above and Proposition~\ref{prop:sublattice-quotientspace} complete the proof of 
Proposition~\ref{prop:subBalg-quotientspace} and its Corollary~\ref{cor:subalgebras-equations}.
\item State and prove results analagous to Proposition~\ref{prop:subBalg-quotientspace} and Corollary~\ref{cor:subalgebras-equations} for the discrete duality between complete and atomic 
Boolean algebras and sets.
\end{enumerate}
\end{ourexercise}

\begin{ourexercise}\label{exe:finitesubs-of-BA}
 Let $\theta$ be an equivalence relation on a Boolean space $X$. Prove that the following 
 conditions are equivalent:
\begin{enumerate}[label=(\roman*)]
\item $\theta$ is clopen as a subset of $X \times X$; 
\item $[x]_\theta$ is clopen for every $x \in X$;
\item $\theta$ has finitely many classes and the corresponding quotient map $X\to X/\theta$ 
is continuous (where $X/\theta$ is a finite discrete space).
\end{enumerate}
\emph{Hint.} (i)$\implies$(ii): Show that both $[x]_\theta = \pi_2[\theta \cap (\{x\} \times X)]$ 
and $\pi_2[\theta^c \cap (\{x\} \times X)]$ are closed. (ii)$\implies$(iii): Use compactness of $X$.
(iii)$\implies$(i): Consider the product map $X \times X\to X/\theta \times X/\theta$ and the inverse image
of the diagonal in $X/\theta \times X/\theta$.
\end{ourexercise}

\begin{ourexercise}\label{exe:subs-of-freedl-two}
  Building on Example~\ref{exa:sub-of-free}, identify the sublattices of the free distributive lattice 
  on two generators and the corresponding compatible quasi-orders on $2^{\{p,q\}}$.
\end{ourexercise}

\begin{ourexercise}\label{exe:equationssubalg}
  Let $M$ be the Boolean subalgebra of $\cP(\bZ)$ consisting of all those $S\subseteq\bZ$ 
  such that both $S\cap\bZ^+$ is either finite or co-finite and $S\cap\bZ^-$ is either finite or 
  co-finite, see Example~\ref{exa:equationssubalg}.
  \begin{enumerate}
    \item Show that $M\cong M^-\times\cP(\{0\})\times M^+$, where $M^-$ and $M^+$ are 
    the Boolean algebras of all finite or co-finite subsets of $\bZ^-$ and $\bZ^+$, respectively.
    \item Show that the dual space of $M$ is the topological sum (that is, disjoint union) of the 
    one-point compactification $\bZ^-\cup\{-\infty\}$ of $\bZ^-$, the one-point 
    compactification $\bZ^+\cup\{+\infty\}$ of $\bZ^+$, and the one-point space $\{0\}$.
  \end{enumerate}
\end{ourexercise}

\begin{ourexercise}\label{exe:quotbeta}
  Let $X$ be a set and $\beta X$ the dual space of $\cP(X)$.
  \begin{enumerate}
    \item Show that there is a one-to-one correspondence between each of
          \begin{enumerate}
            \item The Boolean subalgebras $\cB$ of $\cP(X)$;
            \item The continuous surjections $f\colon\beta X\to Y$ with $Y$ a Boolean space;
            \item The set functions $h\colon X\to Y$ with dense image.
          \end{enumerate}
    \item In particular show that for each $L\subseteq X$, we have
          \[
            \widehat{L}^{\beta X}=\overline{L}^{\beta X}
          \]
          and, for a subalgebra $\cB$ and corresponding continuous surjection $f$ and set 
          function with dense image $h$, we have
          \[
            L\in\cB\ \iff\ f[\widehat{L}^{\beta X}]\text{ is open in }Y\ \iff\ \overline{h[L]}^Y\text{  is open in }Y,
          \]
          and in this case we have
          \[
            \widehat{L}^{Y}=f[\widehat{L}^{\beta X}]=\overline{h[L]}^Y,
          \]
          where $\widehat{(\ )}$ is the Stone map, $\overline{(\ )}$ is topological closure, 
          and the decorations refer to the ambient space in question.
  \end{enumerate}
\end{ourexercise}

\section{Unary operators}\label{sec:unaryopduality}
In all the dualities discussed in this book so far, the morphisms on the algebraic side 
have been the homomorphisms, that is, the maps preserving all of the lattice structure. 
In this section we will relax this condition and study maps between distributive lattices 
that only preserve finite meets, but not necessarily finite joins. Such functions are also 
known as \emph{unary normal multiplicative operators} in the literature. The dualities
developed in this section originate with the seminal works \cite{JonTar1951, JonTar1952}; 
also see the notes at the end of this chapter.
We show that there is still a dual equivalence 
of categories if we generalize, on the space side, from continuous order-preserving 
functions to certain relations that are compatible with the order and topology of the 
Priestley space. We begin by introducing some necessary notation and defining what this means precisely. 

\begin{notationnum} \label{relation-notations} Here and in what follows, we use some common 
notations for composition, forward and inverse image for binary relations: let 
$R \subseteq X \times Y$ and $S \subseteq Y \times Z$ be binary relations. We often use 
\emph{infix notation}, writing $x{R}y$ for $(x,y) \in R$. The 
\emph{composition}\index{relational composition} $R \cdot S$ is the relation from $X$ to 
$Z$ defined by 
    \[ R \cdot S := \{(x,z) \ | \ \text{ there exists } y \in Y \text{ such
    that } x{R}y{S}z\}\ . \] 
\nl{$R \cdot S$}{relational composition of two relations $R$ and $S$, with $R$ before $S$}{cdot}
    Note that this left-to-right notation for relational composition differs from the right-to-left notation $\circ$ that we use for \emph{functional} composition (Notation~\ref{not:function-composition}). That is, if $f \colon X \to Y$ and $g \colon Y \to Z$ are functions with graphs $R_f \subseteq X \times Y$ and $R_g \subseteq Y \times Z$, respectively, then the relational composition $R_f \cdot R_g$ is the graph of the functional composition $g \circ f$.

    The \emph{converse} of $R$ is the relation $R^{-1} := \{(y,x) \in Y \times X : xRy\}$. 
    For any subset $U$ of $X$, we write $R[U]$ for the \emphind{relational direct image}, 
    that is, 
    \[R[U] := \{y \in Y \ | \text{ there exists } u \in U \text{ such that } u{R}y
    \}\ .\] 
The same set is sometimes denoted $\exists_R[U]$.
\nl{$R[U]$}{direct image under a binary relation $R$ of a subset $U$ of the domain}{squarebracket}
\nl{$\exists_R[U]$}{direct image under a binary relation $R$ of a subset $U$ of the domain}{}
For singleton subsets $\{x\}$ of $X$, we write $R[x]$ instead of $R[\{x\}]$.
Converse relations allow us, in particular, to define \emphind{relational
inverse image} $R^{-1}[V]$ for any $V \subseteq Y$, by taking the direct image
of the converse relation, that is,
\[R^{-1}[V] = \{x \in X \ | \ \text{ there exists } v \in V \text{ such that } x{R}v
\} \ . \] 
 Finally, we also use the \emphind{relational universal image}, $\forall_R$, defined, 
 for any $U \subseteq X$, by
  \begin{align}\label{eq:relational-univ}
      \forall_R[U] &:= \{y \in Y \ | \  \text{ for all } x \in X, \text { if }
    x{R}y, \text{ then } x \in U\},
  \end{align}
\nl{$\forall_R[U]$}{relational universal image under a binary relation $R$ of a subset $U$ of the domain}
and we note that
  \[  \forall_R[U] = \{y \in Y \ | \ R^{-1}[y] \subseteq U \} \ . \] 
  For later use, we note a convenient formula for switching between direct and universal relational image:
  \begin{equation}\label{eq:univ-direct-image}
    \forall_R[U] = Y \setminus R[X \setminus U], \quad \text{ for any } U \subseteq X.
  \end{equation}
\end{notationnum}
The import of the operation $\forall_R$ introduced in
(\ref{eq:relational-univ}) stems from the fact that if $R\subseteq X\times Y$
is any relation between sets, then the function $R^{-1}[-] \colon \mathcal{P}(Y) \to
\mathcal{P}(X)$ preserves arbitrary joins, and the function $\forall_R[-]
\colon \mathcal{P}(X) \to \mathcal{P}(Y)$ is its upper adjoint (see
Exercise~\ref{exe:box-duality-viacanext}.\ref{itm:forallisupperadjoint}).
Furthermore, we will now define the property of \emph{upward order-compatibility}
for a relation $R$ between Priestley spaces,
which turns out to be equivalent to the fact that the operations $R^{-1}$ and $\forall_R$ 
restrict correctly to the
sublattices $\Down(X)$ and $\Down(Y)$ (see
Exercise~\ref{exe:box-duality-viacanext}.\ref{itm:uppercompatible}). As we will
see in the proof of Proposition~\ref{prop:relation-dual-to-box} below, these
basic order-theoretic facts, combined with two topological continuity properties for 
the relation $R$, are fundamental to the duality for unary operators.

Note that the notion of `compatibility' for an arbitrary relation between Priestley spaces, 
which we introduce now, is distinct from the notions of compatibility 
for preorders and equivalence relations that we saw in the previous section. 

\begin{definition}\label{dfn:compatiblerelation}
  Let $X$ and $Y$ be Priestley spaces and let $R \subseteq X \times Y$ be a 
  relation. We say that $R$ is:
  \begin{itemize}
    \item \emphind{upward order-compatible} if ${\geq}\cdot R\cdot{\geq}\subseteq R$, 
    that is, for any $x, x' \in X$ and $y, y' \in Y$, whenever $x' \geq x {R} y \geq y'$, 
    we have $x' {R} y'$;
    \item \emphind{upper Priestley continuous} if, for every clopen up-set $K\subseteq Y$, 
    the set $R^{-1}[K]$ is clopen;
    \item \emphind{point-closed} if, for every $x \in X$, the set $R[x]$ is closed;
    \item \emphind{upward Priestley compatible} if $R$ is upward
        order-compatible, upper Priestley continuous, and
        point-closed.\index{Priestley compatible
        relation!upward}\index{compatible!relation on a Priestley space}
  \end{itemize}
\end{definition}

The first aim of this section is to prove (Proposition~\ref{prop:relation-dual-to-box}) 
that the finite-meet-preserving functions between two 
distributive lattices are in one-to-one correspondence with the upward Priestley 
compatible relations between their respective dual spaces. This correspondence will 
generalize the correspondence between homomorphisms and continuous order-preserving 
functions of Priestley duality (see Exercise~\ref{exe:functionalcompatible}). It also yields a 
new duality theorem, as we will remark at the end of this section and prove in the next 
chapter once we have the appropriate categorical terminology in place. We will also show 
at the end of this section how to obtain an analogous correspondence between 
finite-join-preserving functions and downward Priestley compatible relations.

\subsection*{The finite case}
To motivate our proof of the general case (Proposition~\ref{prop:relation-dual-to-box}), 
let us proceed as
we did in Chapters~\ref{ch:order}~and~\ref{ch:priestley} and first examine the
finite case. Let $L$ and $M$ be \emph{finite} distributive lattices. In
Chapter~\ref{ch:order}, we saw that every homomorphism $h \colon M \to L$
uniquely arises as $f^{-1}$ for some order-preserving $f \colon \cJ(L) \to
\cJ(M)$. This function $f$ was obtained as the restriction to $\cJ(L)$ of the lower adjoint of $h$. 
Crucially, we showed in Lemma~\ref{lem:loweradj-preserves-joinprime} that $f$ sends any element of $\cJ(L)$ to $\cJ(M)$. This way of defining a function dual to $h$ first uses that $h$ preserves finite \emph{meets} for the existence of the lower adjoint $f$, and then uses that $h$ preserves
finite \emph{joins} to show that $f$ sends join-primes to join-primes.

Now, when we consider functions between finite distributive lattices that only 
preserve meets, but not necessarily joins, the lower adjoint still exists, 
but it may no longer restrict correctly to join-prime elements. 
Instead, for a meet-preserving function 
$h\colon M\to L$ with lower adjoint $f\colon L\to M$, recall from 
Proposition~\ref{prop:birkhoff} that every element $a \in L$ is a finite join of join-prime 
elements. Therefore, since $f \colon L \to M$ preserves finite joins and $\cJ(L)$ 
join-generates $L$, the function $f$ is uniquely determined by its restriction to join-prime 
elements: for any $a \in L$, we have $f(a) = \bigvee f[{\downarrow} a \cap \cJ(L)]$. 
Moreover, since $\cJ(M)$ join-generates $M$, in order to dually encode the function $f$, 
and therefore $h$, it suffices to know the value $f(p)$ for each $p \in \cJ(L)$. To this end, 
we define
\begin{equation} \label{eq:dualrelation}
  R := \{(p,q) \in \cJ(L) \times \cJ(M) \ | \ q \leq f(p)\},
\end{equation}
so that $f(p)=\bigvee R[p]$ for every $p \in \cJ(L)$.
Note that we can express $R$ in terms of $h$ using the adjointness:
\begin{align*}
  q \leq f(p)\  & \iff\ \forall b\in M \ ( f(p)\leq b \implies q\leq b )                    \\
                & \iff\ \forall b\in M \ ( p\leq h(b) \implies q\leq b )                    \\
                & \iff \ \forall b\in M \ ( p\in \widehat{h(b)} \implies q\in\widehat{b} ),
\end{align*}
where we recall that $\widehat{(-)}$ denotes the lattice isomorphism between a finite 
distributive lattice and the down-set lattice of its poset of join-prime elements (Proposition~\ref{prop:birkhoff}).

We call the binary relation $R \subseteq \cJ(L) \times \cJ(M)$ defined by 
(\ref{eq:dualrelation}) \emph{the relation dual to $h$}. This relation $R$ is 
upward order-compatible, see Exercise~\ref{exe:finiterelationduality}. Moreover, if we view
the finite posets $\cJ(L)$ and $\cJ(M)$ as Priestley spaces by equipping them with
the discrete topology, then any relation between them 
is trivially upper Priestley continuous and point-closed.
The original 
meet-preserving function $h \colon M \to L$ can be recovered from $R$, using the following equality, which holds for any $b \in M$:
\begin{equation}\label{eq:boxfromrelation}
  \widehat{h(b)} = \forall_{R^{-1}}[\widehat{b}],
\end{equation}
where we recall that $\forall_{R}$ is the universal image defined in 
(\ref{eq:relational-univ}).
Writing out the definitions, (\ref{eq:boxfromrelation}) expresses the fact that, 
for any $q \in \cJ(M)$,
\[ q \leq h(b) \iff \text{ for all } p \in \cJ(L), \text{ if } p {R} q, \text{ then } p \leq b,\]
which can be proved using the lower adjoint $f$ of $h$ and the definition of 
$R$ (see Exercise~\ref{exe:finiterelationduality}).

This concludes our informal description, in the case of finite distributive lattices, 
of the duality between meet-preserving functions and upward Priestley compatible relations. 
Summing up, we have represented the meet-preserving function $h$ by an 
upward order-compatible relation 
$R$ between the dual posets $\cJ(L)$ and $\cJ(M)$, from which $h$ can be recovered 
as the upper adjoint of the unique join-preserving function $f$ which is defined for 
$p \in \cJ(L)$ by $f(p) := \bigvee R[p]$. As an instructive exercise (see 
Exercise~\ref{exe:finiterelationduality}), we invite you to verify the claims that were left 
unproved here, although they are also direct consequences of  the general duality 
theorem that we prove below (see Exercise~\ref{exe:box-duality-viacanext} for details 
about the relationship). A reader familiar with modal logic may have recognized in the equation 
(\ref{eq:boxfromrelation}) the definition of the $\Box$ (``box'') operator associated 
to a Kripke relation $R$; more on this in Section~\ref{sec:kripke} below.

\subsection*{The general case}
In the remainder of this section, we generalize the ideas outlined above to 
finite-meet-preserving functions $h$ between arbitrary distributive lattices that are 
not necessarily finite. As in Chapter~\ref{ch:priestley}, join-primes have to be 
replaced by points of the dual space, but the underlying ideas are the same as in 
the finite case.

\begin{proposition}\label{prop:relation-dual-to-box}
  Let $L$ and $M$ be distributive lattices with Priestley dual spaces $X$ and $Y$, 
  respectively. For any finite-meet-preserving $h \colon M \to L$, there exists a unique 
  upward Priestley compatible relation, $R \subseteq X \times Y$, such that,
  \begin{equation}\label{eq:h-is-box-of-relation}
    \text{ for any } b \in M, \ \widehat{h(b)} = \forall_{R^{-1}}[\widehat{b}].
  \end{equation}
  This relation $R$ may be defined explicitly, for $x \in X$, by
  \begin{equation}\label{eq:relation-image-of-point}
    R[x] := \bigcap \{ \widehat{b} \mid b \in M, \; x \in \widehat{h(b)} \},
  \end{equation}
  or, equivalently,
  \begin{equation}\label{eq:relation-dual-to-box-def}
    R := \{(x,y) \in X \times Y \ \mid \ \text{for all } b \in M, \text{ if } h(b) \in F_x, \text{ then } b \in F_y \}.
  \end{equation}
\end{proposition}

\begin{proof}
  Let $R \subseteq X \times Y$ denote the relation defined in (\ref{eq:relation-image-of-point}) 
  and (\ref{eq:relation-dual-to-box-def}). We will now establish three properties, which suffice to conclude: 
  \begin{enumerate}
    \item[(1)] $R$ satisfies 
  (\ref{eq:h-is-box-of-relation}); 
    \item[(2)] $R$ is upward Priestley compatible; and 
    \item[(3)] $R$ is the unique relation with properties (1) and (2).
  \end{enumerate}

  For (1), unfolding the definitions, we need to prove that, for any $x \in X$ and $b \in M$,
  \begin{equation} \label{eq:to-prove-h-is-box}
    x \in \widehat{h(b)} \iff \forall y \in Y, \text{ if } x R y, \text{ then } y \in \widehat{b}.
  \end{equation}
  The left-to-right direction is clear by definition of $R$. For the converse, we reason 
  contrapositively. Suppose that $x \not\in \widehat{h(b)}$; this means that the 
  homomorphism $h_x \colon L \to \btwo$ sends $h(b)$ to $0$. Since $h$ is 
  finite-meet-preserving, the composite function $k := h_x \circ h \colon M \to \btwo$ 
  is finite-meet-preserving. Thus, the set $F := k^{-1}(1)$ is a filter which does not 
  contain $b$. By the prime filter theorem (Theorem~\ref{thm:DPF}), pick a prime filter 
  $F_y$ containing $F$ and still not containing $b$. The fact that $F \subseteq F_y$ is 
  easily seen to be equivalent to $x{R}y$, while $y \not\in \widehat{b}$, establishing 
  that the right-hand-side of (\ref{eq:to-prove-h-is-box}) fails, as required.

  For (2), we first show that (1) already yields that $R$ is upper Priestley
  continuous. Indeed, any clopen up-set $K \subseteq Y$ is equal to $Y
  \setminus \widehat{b}$ for some $b \in M$, so that $R^{-1}[K] = X \setminus
  \forall_{R^{-1}}[\widehat{b}]$, by the formula for switching between
  universal and direct relational image (\ref{eq:univ-direct-image}). This
  means that $R^{-1}[K]$ is clopen, since by (1) it is equal to the complement
  of the clopen set $\widehat{h(b)}$. Also, $R$ is point-closed, since 
  (\ref{eq:relation-image-of-point})
  shows that,
  for any $x \in X$, $R[x]$ is a closed down-set in $Y$. 
  To finish the proof of (2), note that upward
  order-compatibility follows easily from the definitions, or also by an
  application of Exercise~\ref{exe:ordercompatible}.

  We now prove (3). Indeed, we will prove the following stronger fact, 
  namely, that for any upward Priestley compatible relations $R, S 
  \subseteq X \times Y$, we have
  \begin{equation} \label{eq:2-eq-box-relations}
    S \subseteq R \iff \text{ for every } b \in M, \ \forall_{R^{-1}}[\widehat{b}] \subseteq \forall_{S^{-1}}[\widehat{b}].
  \end{equation}
  Note that the uniqueness (3) follows from (\ref{eq:2-eq-box-relations}), 
  for if $S$ is any upward Priestley compatible relation satisfying 
  (\ref{eq:h-is-box-of-relation}), then $\forall_{S^{-1}}[\widehat{b}] = 
  \widehat{h(b)} = \forall_{R^{-1}}[\widehat{b}]$ for every $b \in M$, so $S = R$. 
  The left-to-right direction of (\ref{eq:2-eq-box-relations}) is immediate from the 
  definition of the universal image. For the converse, we will reason by 
  contraposition and use the fact that $R$ is point-closed. Suppose that 
  $S \not\subseteq R$; pick $x \in X$ and $y \in Y$ such that $y \in S[x]$ and 
  $y \not\in R[x]$. Since $R[x]$ is closed, and also a down-set by order-compatibility, 
  there exists $b \in M$ such that $R[x] \subseteq \widehat{b}$ and 
  $y \not\in \widehat{b}$, as follows from Proposition~\ref{prop:order-normality-Priestley}. 
  It now follows from 
  the definition of universal image that $x \in \forall_{R^{-1}}[\widehat{b}]$, but 
  $x \not\in \forall_{S^{-1}}[\widehat{b}]$ since $x{S}y$ but $y \not\in \widehat{b}$. 
  This concludes the proof of (\ref{eq:2-eq-box-relations}).
\end{proof}

\begin{definition}\label{dfn:relation-dual-to-box}
  The relation $R$ defined in Proposition~\ref{prop:relation-dual-to-box} is called 
  the \emph{dual relation} of the finite-meet-preserving function $h$.
\end{definition}
In  Exercise~\ref{exe:every-compatible-relation-comes-from} you are asked to 
use the uniqueness part of Proposition~\ref{prop:relation-dual-to-box} to derive
that, if $X$ and $Y$ are Priestley spaces,
then every upward Priestley compatible relation $R \subseteq X \times Y$ 
is the dual relation of a unique  
finite-meet-preserving function.

Proposition~\ref{prop:relation-dual-to-box} is the crucial new ingredient for the 
following extension of the Priestley duality theorem (Theorem~\ref{thm:priestleyduality}) 
to a larger collection of morphisms, namely all finite-meet-preserving functions.
\begin{theorem}\label{thm:unaryboxduality}
  The category of distributive lattices with finite\hyp{}meet\hyp{}preserving functions 
  is dually equivalent to the category of Priestley spaces with upward Priestley compatible 
  relations.
\end{theorem}
The proof uses as its crucial ingredient Proposition~\ref{prop:relation-dual-to-box} above, 
combined with some techniques from category theory, and will be given in 
Section~\ref{sec:duality-categorically}, p.~\pageref{unaryboxproof}. 

To finish this section, we draw a few further corollaries from 
Proposition~\ref{prop:relation-dual-to-box}. First, we show how it specializes to the Boolean case.

\begin{definition}\label{dfn:compatible-boolean}
  Let $X$ and $Y$ be Boolean spaces. A relation $R \subseteq X \times Y$ is
  called \emph{Boolean compatible}\index{compatible!relation on a
  Boolean space}\index{Boolean compatible relation} if it is point-closed and 
  continuous, that is, for any clopen $K \subseteq Y$, the set $R^{-1}[K]$ is clopen.
\end{definition}
Since the Priestley order on a Boolean space is trivial, note that a relation between 
Boolean spaces is upward Priestley compatible if, and only if, it is compatible according 
to Definition~\ref{dfn:compatible-boolean}. Combining this observation with Propositions~\ref{prop:relation-dual-to-box}~and~\ref{prop:boolean-trivial-order} 
allows us to deduce the following.
\begin{corollary}
  Let $A$ and $B$ be Boolean algebras with dual spaces $X$ and $Y$, respectively. 
  Finite-meet-preserving functions $B \to A$ are in a one-to-one correspondence with 
  Boolean compatible relations $R \subseteq X \times Y$.
\end{corollary}

For easy reference and future use, we also record the order-duals of 
Proposition~\ref{prop:relation-dual-to-box} and Theorem~\ref{thm:unarydiamondduality}, 
and the accompanying definitions.
\begin{definition}\label{dfn:downward-compatible}
  Let $X$ and $Y$ be Priestley spaces and let $R \subseteq X \times Y$ be a relation. 
  We say that $R$ is:
  \begin{itemize}
    \item \emph{downward order-compatible} if ${\leq} \cdot R \cdot {\leq} \subseteq R$, 
    that is, for any $x, x' \in X$ and $y, y' \in Y$, whenever $x' \leq x {R} y \leq y'$, we have 
    $x' {R} y'$;
    \item \emph{lower Priestley continuous} if, for every clopen down-set $K \subseteq Y$, 
    the set $R^{-1}[K]$ is clopen;
    \item \emph{downward Priestley compatible} if $R$ is downward order-compatible, 
    lower Priestley continuous, and point-closed.\index{Priestley compatible 
    relation!downward}\index{downward Priestley compatible}\index{compatible!relation on 
    a Priestley space}
  \end{itemize}
\end{definition}
Let $X$ and $Y$ be Priestley spaces and let $R \subseteq X \times Y$ be a downward Priestley
compatible relation. We define the following function:
\begin{align}\label{eq:diamondR}
  \exists_{R^{-1}}    \colon \ClD(Y) & \to \ClD(X) \nonumber \\
  b & \mapsto R^{-1}[b],
\end{align}
which we note is the same as $X \setminus \forall_{R^{-1}}[Y \setminus b]$ by 
(\ref{eq:univ-direct-image}). Conversely, any finite-join-preserving function 
$h \colon \ClD(Y) \to \ClD(X)$ is equal to $\exists_{R^{-1}}$ for a unique downward 
Priestley compatible relation $R_h \subseteq X \times Y$, which can be defined 
explicitly by
\begin{equation}\label{eq:Rdiamond}
  R_h := \{(x,y) \in X \times Y \ \mid \ \text{ for every } b \in \ClD(Y), \text{ if } y \in b, 
  \text{ then } x \in h(b)\}.
\end{equation}

\begin{proposition}\label{prop:relation-dual-to-diamond}
  Let $L$ and $M$ be distributive lattices with Priestley dual spaces $X$ and $Y$, 
  respectively. The assignments $R \mapsto \exists_{R^{-1}}$ (\ref{eq:diamondR}) and 
  $h\mapsto R_h$ (\ref{eq:Rdiamond}) form a bijection between finite-join-preserving 
  functions from $M$ to $L$ and downward Priestley compatible relations from $X$ to $Y$.
\end{proposition}
Proposition~\ref{prop:relation-dual-to-diamond} has essentially the same proof as 
Proposition~\ref{prop:relation-dual-to-box}. Instead of re-doing the entire proof, one 
may also appeal to order-duality to \emph{deduce} this proposition 
from Proposition~\ref{prop:relation-dual-to-box}, see 
Theorem~\ref{thm:unarydiamondduality} in Chapter~\ref{ch:categories} and the 
remarks following it.

We finish by examining two more special cases that may help elucidate the connection 
between the relations dual to finite-join- and finite-meet-preserving functions; in both of 
these cases, we examine a special setting where the two notions of dual relation 
interact with each other.

First, if $f \colon L \leftrightarrows M \colon g$ is an \emphind{adjoint pair} between distributive lattices, then the left adjoint $f$ is finite-join-preserving and the right adjoint $g$ is finite-meet-preserving (see Exercise~\ref{exe:adjunctions}). Denote by $X$ and $Y$ the Priestley dual spaces of $L$ and $M$, respectively. By the results of this section, $f$ has a dual downward Priestley compatible relation $R_f \subseteq Y \times X$ and $g$ has a dual upward Priestley compatible relation $R_g \subseteq X \times Y$. The two relations are closely related: $R_f$ is the relational \emphind{converse} of $R_g$, see Exercise~\ref{exe:converse-adjunction}.\index{relational converse!and adjunction}

Second, if $h \colon M \to L$ is a \emph{homomorphism} between the distributive lattices $L$ and $M$, then $h$ has \emph{both} an upward Priestley compatible dual relation $R_h$, because it preserves finite meets, \emph{and} a downward Priestley compatible dual relation $S_h$, because it preserves finite joins. The intersection of $R_h$ and $S_h$ can be seen to be a \emph{functional} relation, that is, for every $x \in X$, there is a unique $y \in Y$ such that $(x,y) \in R_h \cap S_h$. This relation $R_h \cap S_h$ is in fact equal to (the graph of) the continuous order-preserving function $f \colon X \to Y$ dual to $h$, as it was defined in Section~\ref{sec:topologize}; see Exercise~\ref{exe:functional-relation}.

\ourexercises

\begin{ourexercise}\label{exe:finiterelationduality}
  Let $h \colon M \to L$ be a meet-preserving function between finite distributive lattices and let $R$ be the relation defined in (\ref{eq:dualrelation}).
  \begin{enumerate}
    \item Prove that $R$ is upward order-compatible.
    \item Prove equation (\ref{eq:boxfromrelation}).
    \item Show directly (that is, without referring to Proposition~\ref{prop:relation-dual-to-box}) that the assignment $h \mapsto R$ is a bijection between meet-preserving functions from $M$ to $L$ and upward order-compatible relations from $\cJ(L)$ to $\cJ(M)$.
  \end{enumerate}
\end{ourexercise}

\begin{ourexercise}\label{exe:ordercompatible}
  Let $X$ and $Y$ be posets and $R \subseteq X \times Y$ a relation. Prove that  $R$ is upward order-compatible if, and only if, $R^{-1}$ is downward order-compatible if, and only if, for any subsets $S \subseteq X$ and $T \subseteq Y$, $R[S]$ is a down-set and $R^{-1}[T]$ is an up-set.
\end{ourexercise}

\begin{ourexercise}\label{exe:box-duality-viacanext}
  This exercise shows in more detail how the proof of Proposition~\ref{prop:relation-dual-to-box} generalizes the proof sketch for the finite case given in the beginning of the section, and how one could naturally arrive at the compatibility conditions of Definition~\ref{dfn:compatiblerelation}. It also establishes a link between duality for finite-meet-preserving maps and the canonical extension of such maps.

  Let $X$ and $Y$ be posets and let $R \subseteq X \times Y$ be a relation.
  \begin{enumerate}
      \item \label{itm:forallisupperadjoint} Prove that the pair of functions $R[-] \colon \cP(X) \leftrightarrows \cP(Y) \colon \forall_{R^{-1}}[-]$ is an adjoint pair.
      \item \label{itm:uppercompatible} Using Exercise~\ref{exe:ordercompatible}, show that it then follows that this adjunction restricts to a well-defined adjunction  $R[-] \colon \Down(X) \leftrightarrows \Down(Y) \colon \forall_{R^{-1}}[-]$ if, and only if, $R$ is upward order-compatible.
  \end{enumerate}

  Further assume that $X$ and $Y$ are Priestley spaces dual to distributive lattices $L$ and $M$, respectively.
  \begin{enumerate}
    \item[c.]
      Prove that $R$ is upward Priestley compatible if, and only if, the function $\forall_{R^{-1}}[-] \colon \Down(Y) \to \Down(X)$ factors through the embeddings $L \into \Down(X)$ and $M \into \Down(Y)$, that is, if there exists a function $h \colon M \to L$ such that $\widehat{h(b)} = \forall_{R^{-1}}[\widehat{b}]$ for all $b \in M$. Also show that such a function $h$, if it exists, must preserve finite meets.
    \item[d.] Show that, even if $R$ is upward Priestley compatible, the function $R[-]$ does \emph{not} necessarily factor through $L \into \Down(X)$ and $M \into \Down(Y)$.
    \item[e.] Prove that for any finite-meet-preserving function $h \colon M \to L$, there exists a unique completely meet-preserving function $h^\delta \colon \Down(Y) \to \Down(X)$, such that $h^\delta(\widehat{b}) = \widehat{h(b)}$ for any $b \in M$; show that the lower adjoint of $h^\delta$ is equal to the relational direct image function $R[-] \colon \Down(X) \to \Down(Y)$, where $R$ is the relation dual to $h$.

        {\it Note.} This item is straightforward to prove if one combines the earlier two items of this exercise with the results from this section. A more difficult exercise is to prove this last item without using the results from this section; see for example \cite{GehJon1994}. This can then be used to give an alternative proof of Theorem~\ref{thm:unaryboxduality}.
    \item[f.] Explain why the proof for the finite case, outlined earlier in this section, is a special case of the previous item.
  \end{enumerate}
\end{ourexercise}

\begin{ourexercise}\label{exe:functionalcompatible}
  This exercise shows that the correspondence between upward Priestley compatible relations and finite-meet-preserving functions generalizes Priestley duality for homomorphisms given in Chapter~\ref{ch:priestley}. Let $X$ and $Y$ be dual Priestley spaces of distributive lattices $L$ and $M$, respectively.
  \begin{enumerate}
    \item Prove that, if $f \colon X \to Y$ is a continuous order-preserving function, then $R_f := \{(x,y) \in X \times Y \mid f(x) \geq y \}$ is an upward Priestley compatible relation, which moreover has the property that $R[x]$ has a maximum for every $x \in X$.
    \item Prove that, if $R \subseteq X \times Y$ is an upward Priestley compatible relation and $R[x]$ has a maximum for every $x \in X$, then $f_R \colon X \to Y$ defined by $f_R(x) := \max R[x]$ is a continuous order-preserving function.
    \item Prove that a finite-meet-preserving function $f \colon M \to L$ is a homomorphism if, and only if, the dual relation $R$ is such that $R[x]$ has a maximum for every $x \in X$.
  \end{enumerate}
\end{ourexercise}

\begin{ourexercise}\label{exe:every-compatible-relation-comes-from}
  Let $X$ and $Y$ be Priestley spaces with dual lattices $L$ and $M$, respectively, and let $R$ be an upward Priestley compatible relation from $X$ to $Y$. Show that there exists a unique meet-preserving function $h \colon M \to L$ such that (\ref{eq:h-is-box-of-relation}) holds.
  Conclude that this $h$ is the unique meet-preserving function $M \to L$ for which $R$ is the dual relation.
\end{ourexercise}

\begin{ourexercise}\label{exe:converse-adjunction}
  This exercise shows that ``the dual relations of adjoint pairs are converse to each other''. \index{adjoint pair!dual relation}\index{relational converse!and adjunction}
  Let $f \colon L \leftrightarrows M \colon g$ be an adjoint pair between distributive lattices, and let $X$ and $Y$ be the Priestley dual spaces of $L$ and $M$, respectively. Since $f$ is a left adjoint, it preserves finite joins; let $R_f \subseteq Y \times X$ denote the downward Priestley compatible relation dual to $f$. Similarly, let $R_g \subseteq X \times Y$ denote the upward Priestley compatible relation dual to $g$. Prove that $(x,y) \in R_g$ if, and only if, $(y,x) \in R_f$. That is, $R_f$ and $R_g$ are converse relations. \hint{The first items of Exercise~\ref{exe:box-duality-viacanext} can also be useful here.}
\end{ourexercise}

\begin{ourexercise}\label{exe:functional-relation}\index{homomorphism!dual relations}\index{homomorphism!dual function}
  Let $h \colon M \to L$ be a homomorphism between distributive lattices and let $X$ and $Y$ be the Priestley dual spaces of $L$ and $M$, respectively. Denote by $R_h \subseteq X \times Y$ the upward Priestley compatible relation dual to $h$, viewed as a finite-meet-preserving function, and by $S_h \subseteq X \times Y$ the downward Priestley compatible relation dual to $h$, viewed as a finite-join-preserving function. Prove that the intersection $R_h \cap S_h$ of the two relations is a functional relation, and that this is the continuous order-preserving function dual to $h$, as defined in Section~\ref{sec:topologize}.
\end{ourexercise}

\section{Modal algebras and Kripke completeness}\label{sec:kripke}
In this section we show how the duality results of the previous section relate
to Kripke's possible world semantics for classical modal logic. We will see in Theorem~\ref{thm:Kcomplete} that the dual space of the free modal algebra gives a canonical model for the basic normal modal logic $\bfK$. In this view, Kripke semantics for modal logic is obtained from the Stone-J\'onsson-Tarski duality for operators that we developed in Section~\ref{sec:unaryopduality} by `forgetting the topology'.
We begin by giving a minimal introduction to modal logic, limiting ourselves to the parts that are needed for understanding the connection to duality theory; for much more material on modal logic and duality we refer to classic textbooks in the field, such as \cite{BRV2001, ChaZak1997}.

The basic normal modal logic, $\bfK$, is an extension of classical
propositional logic by a \emph{necessity operator}, $\Box$.\index{necessity
operator}\index{modal operator} Formally,
\emph{modal formulas}\index{modal formula}\index{modal logic!classical} 
are terms built from propositional variables and the
constants $\top$, $\bot$, using the binary operations $\vee$, $\wedge$, 
and unary operations $\neg$ and $\Box$; we use the common notational convention
that unary operators bind more strongly than binary ones, that is, the notation
$\Box p \to \neg q$ denotes the formula $(\Box p) \to (\neg q)$, which is
different from the formula $\Box (p \to \neg q)$. We also consider the
\emph{possibility operator} $\Diam$ as an abbreviation, defined by de Morgan
duality as $\Diam \phi := \neg \Box \neg \phi$, and implication may be defined
by $\phi \to \psi := \neg \phi \vee \psi$. 
Note that, if one considers modal operators in the setting of positive logic
(that is, distributive lattices or Heyting algebras) rather than of classical logic 
(that is, Boolean algebras), the interplay between $\Box$, $\Diam$, and the
lattice operations is more subtle, and they are no longer interdefinable by 
formulas. In this section, for simplicity we only treat the classical, Boolean 
algebra-based case, although it is possible to develop dualities for modal algebras 
based on distributive lattices and Heyting algebras~\parencite{Gol1989, CelJan99, 
Bez99}.  We will consider the distributive case in Example~\ref{exa:modal-alg-functor} 
and Exercise~\ref{exer:boxfunctor} in the next chapter, and we also revisit duality
for operators on distributive lattices in Section~\ref{sec:funcspace}, where it will 
be related to a function space construction; see in particular
Definition~\ref{def:uppervietoris}.\index{modal logic!positive}

\nl{$\vdash_{\bfK}$}{derivability in normal modal logic ${\bfK}$}{}
A notion of \emph{derivability} between modal formulas, $\vdash_{\bfK}$, may be defined by extending a Hilbert-style proof calculus for classical logic with one additional axiom, $\Box(p \to q) \to (\Box p \to \Box q)$, and one additional rule: from $\vdash_{\bfK} \phi$, infer $\vdash_{\bfK} \Box \phi$. We will not need to enter into details of the proof calculus here; see for example \cite[Section~3.6]{ChaZak1997} for more details. Crucially for us is the characteristic property of $\bfK$ that, for any formulas $\phi$ and $\psi$, $\vdash_{\bfK} \Box(\phi \wedge \psi) \liff (\Box \phi \wedge \Box \psi)$, and $\vdash_{\bfK} \Box \top \liff \top$. In other words, the operation $\phi \mapsto \Box \phi$ yields a finite-meet-preserving function on the set of $\vdash_{\bfK}$-equivalence classes of modal formulas. This leads to the following definition.
\begin{definition}\label{dfn:modal-algebra}
  A \emph{modal algebra} is a pair $(B, \Box)$, where $B$ is a Boolean algebra, and $\Box \colon B \to B$ is a finite-meet-preserving endofunction on $B$.
  A \emph{homomorphism} from a modal algebra $(B, \Box_B)$ to a modal algebra $(A, \Box_A)$ is a homomorphism $h \colon B \to A$ such that, for every $b \in B$, $h(\Box_B b) = \Box_A h(b)$.
\end{definition}
Note that, if $(B,\Box)$ is a modal algebra and $V$ is a set of variables, and $f \colon V \to B$ is any function, then any modal formula $\phi$ with propositional variables drawn from $V$ has a uniquely defined \emph{interpretation} $\bar{f}(\phi)$ in $B$, which may be defined inductively by setting $\bar{f}(v) := f(v)$ for $v \in V$, $\bar{f}(\top) := \top$, $\bar{f}(\Box \phi) := \Box \bar{f}(\phi)$, $\bar{f}(\phi \vee \psi) := \bar{f}(\phi) \vee \bar{f}(\psi)$, $\bar{f}(\neg \phi) := \neg \bar{f}(\phi)$, etc. In algebraic terms, $\bar{f}$ is the unique extension of $f$ to a homomorphism from the term algebra over $V$ to $B$.

It may be proved that the set of $\bfK$-derivable equivalence classes of modal
formulas in a fixed set of variables $V$ is, up to isomorphism, the \emph{free
modal algebra} over $V$. The proof is a straightforward adaptation to the
modal case of the construction of the Lindenbaum-Tarski algebra for
propositional logic, which we saw when we constructed the free distributive
lattice in Section~\ref{sec:free-description}, see Remark~\ref{rem:LT-algebra};
for all the details in the modal setting, see for example 
\cite[Section~7.5]{ChaZak1997} or \cite[Section~5.2]{BRV2001}. From this fact, one
obtains the following theorem, which provides the foundation for algebraic
modal logic, where it is known as the ``algebraic completeness theorem for
$\bfK$''.

\begin{theorem}\label{thm:K-algebraic-completeness}
  Let $\phi$ and $\psi$ be modal formulas with propositional variables among
  $x_1, \dots, x_n$. Then $\vdash_{\bfK} \phi \liff \psi$ if, and only if, for
  every function $f \colon \{x_1,\dots,x_n\} \to B$, with $B$ a modal algebra,
  we have $\bar{f}(\phi) = \bar{f}(\psi)$. In particular, $\vdash_{\bfK} \phi$
  if, and only if, $\bar{f}(\phi) = \top$ for every interpretation $f$.
\end{theorem}
\begin{proof}
  See, for example, \cite[Theorems~7.43 and 7.44]{ChaZak1997} or \cite[Theorem~5.27]{BRV2001}.
\end{proof}

We now show how this algebraic completeness theorem for $\bfK$ may be combined
with the duality in Theorem~\ref{thm:unaryboxduality} of the previous section to
obtain easy proofs of Kripke completeness theorems for $\bfK$. We first give
the definition of Kripke semantics. 
\begin{definition}
  A (discrete) \emphind{Kripke frame} is a pair $(X, R)$, with $X$ a set and $R
  \subseteq X \times X$ a relation. A \emphind{Kripke Boolean space} is a
  triple $(X, \tau, R)$, with $(X,\tau)$ a Boolean space and $R \subseteq X
  \times X$ a compatible relation on $X$.

For a set of propositional variables $V$, a \emphind{valuation} to a Kripke
frame $(X,R)$ is a function $c \colon X \to 2^V$. When $(X,\tau,R)$ is a Kripke
Boolean space, a valuation $c \colon X \to 2^V$ is called
\emph{admissible}\index{admissible valuation}\index{valuation!admissible, for Kripke Boolean models} if,
for every $v \in V$, the set
\[ \{ x \in X \ \mid \ c(x)(v) = 1 \} \]
is clopen.
A (discrete) \emphind{Kripke model} is a tuple $(X,R,c)$ where $c$ is a valuation, and a
\emphind{Kripke Boolean model} is a tuple $(X,\tau,R,c)$ where $(X,\tau,R)$ is a
Kripke Boolean space and $c$ is an admissible valuation.
\end{definition}
Kripke models can be used to define a concrete, set-based semantics for modal
logic, generalizing truth tables for propositional logic, as follows. Fix a set
of variables $V$. For any Kripke model $(X,R,c)$, 
\nl{$\models$}{forcing relation on Kripke models}{}
a \emph{forcing} or \emph{truth} relation, $\models$ is
defined between points of $X$ and modal formulas with variables in $V$, as
follows. For any $x \in X$ and $v \in V$, define $x \models v$ if, and only if,
$c(x)(v) = 1$, and also define $x \models \top$ to always hold. Then, extend
$\models$ inductively by defining, for any modal formulas $\phi$, $\psi$, that 
$x \models \phi \vee \psi$ if, and only if, $x \models \phi$ or $x \models
\psi$; $x \models \neg \phi$ if, and only if, it is not the case that $x \models \phi$; and finally,  
\begin{equation}\label{eq:box-semantics}
  x \models \Box \phi \stackrel{\mathrm{def}}{\iff} \text{ for every } y \in X, \text{ if } x {R} y \text{ then } y \models \phi.
\end{equation}
If $x \models \phi$, then we say that $\phi$ \emph{holds} or \emph{is true} at
$x \in X$. This definition in particular applies to a Kripke Boolean model
$(X,\tau,R,c)$, where the further condition that the valuation $c$ must be 
admissible implies that the \emphind{truth set} of $\phi$,
\nl{$\sem{\phi}$}{truth set of a modal formula $\phi$ in a Kripke model}{}
\[\sem{\phi} := \{x \in X \ \mid \ x \models \phi\},\] 
is clopen, for every modal formula $\phi$; also see 
Proposition~\ref{prop:kripke-semantics} and Exercise~\ref{exe:clopenvaluation} 
below. If $\phi$ is a modal formula with variables in $V$ and $(X,R,\tau,c)$ is a 
Kripke Boolean model, then $\phi$ is \emph{valid} on the model 
\index{valid!formula on a Kripke Boolean model}\index{valid!formula on a Kripke 
model} provided that $\phi$ is true at every $x \in X$. The formula $\phi$ is
\emph{valid} on a Kripke Boolean space if it is valid for every admissible valuation. 
On the other hand,  $\phi$ is valid on a (discrete) Kripke frame $(X,R)$ it is valid 
under every valuation. Clearly, if $(X,R,\tau)$ is a Kripke Boolean space and $\phi$ 
is a formula that is valid on the underlying discrete Kripke frame $(X,R)$,  then 
$\phi$ is in particular valid on the Kripke Boolean space $(X,R,\tau)$. The converse 
is not true in general; for an example (see Exercise~\ref{exe:GL}). 

\begin{example}\label{exa:kripkeframe}
  Consider the Kripke frame $(\mathbb{N}, <)$. Examples of valid formulas on this 
  Kripke frame are $\neg \Box \bot$, since for every $n \in \mathbb{N}$ there exists 
  $m \in \mathbb{N}$ with $n < m$, and also $\Box \Box p \to \Box p$, since $<$ 
  is a transitive relation. We leave it as an instructive exercise for the reader unfamiliar 
  with modal logic to check in detail that these formulas are indeed valid. An example 
  of a formula that is not valid on this Kripke frame is $\Box p \to p$: for the valuation 
  $c$ which sends $0\in\mathbb{N}$ to $(p\mapsto 0)$ and any other $n\in\mathbb{N}$ 
  to $(p\mapsto 1)$, we have that $0 \models \Box p$ but $0 \not\models p$, so that 
  $0 \not\models \Box p \to p$. Note that, under this valuation, $\Box p \to p$ does 
  happen to be true in all $n \in \mathbb{N} \setminus \{0\}$. In fact, one may prove 
  as another instructive exercise that the formula $\Box p \to p$ is valid on a Kripke 
  frame if, and only if, the relation of the frame is reflexive 
  (see Exercise~\ref{exe:kripke-examples}).
\end{example}
A first connection between Kripke frames and modal algebras is given by 
Corollary~\ref{cor:modal-object-duality} below, a consequence of the results in the 
previous section.
\begin{definition}
  For a Kripke Boolean space $(X,\tau,R)$, define its \emph{dual modal algebra} to 
  be the pair $(B,\forall_{R^{-1}})$, where $B$ is the Boolean algebra dual to $(X,\tau)$.
\end{definition}
The following is now a straightforward application of Proposition~\ref{prop:relation-dual-to-box}.
\begin{corollary}\label{cor:modal-object-duality}
  Every modal algebra $(B, \Box)$ is isomorphic to the dual modal algebra of a
  Kripke Boolean space $(X,\tau,R)$, which is unique up to isomorphism.
\end{corollary}
\begin{proof}
  By Stone duality, there is an up to isomorphism unique Boolean space
  $(X,\tau)$ dual to $B$. By Proposition~\ref{prop:relation-dual-to-box}, and
  the fact that `upward Priestley compatible' means `compatible'
  (Definition~\ref{dfn:compatible-boolean}) in the Boolean case, there is a
  unique compatible relation $R$ on $X$ such that $\Box = \forall_{R^{-1}}$.
\end{proof}

We now use Corollary~\ref{cor:modal-object-duality} to reformulate Kripke's 
semantics for modal logic in the language of modal algebras. Note first that 
functions $c \colon X \to 2^V$ are in a bijection with functions 
$f \colon V \to \cP(X)$, via the ``currying'' map which sends any function 
$c \colon X \to 2^V$ to the function $f_c : V \to \cP(X)$ defined by
\[ f_c(v) := \{ x \in X \ | \ c(x)(v) = 1\}.\]
The relationship between Kripke semantics and homomorphisms of modal 
algebras is now explained by the following proposition. Recall that $\bar{f}$ 
denotes the unique interpretation of modal formulas into a modal algebra 
extending a given interpretation $f$ of the variables. Also note that 
$(\cP(X), \forall_{R^{-1}})$ is a modal algebra for any Kripke frame $(X,R)$.

\begin{proposition}\label{prop:kripke-semantics}
  Let $(X,R)$ be a Kripke frame and $c \colon X \to 2^V$ a valuation. For any 
  modal formula $\phi$ and $x \in X$, $x \models \phi$ if, and only if, 
  $x \in \overline{f_c}(\phi)$.
\end{proposition}
\begin{proof}
    Let us denote by $h$ the function that sends any modal formula $\phi$ 
    to its truth set, that is, for any modal formula $\phi$, 
  \[ h(\phi) := \{ x \in X \mid x \models \phi\}.\] 
  The claim may then be rephrased as saying that 
  $$h(\phi) = \overline{f_c}(\phi)$$ 
  for all modal formulas $\phi$. To prove this, we
  show first that $h$ respects the interpretation in $(\mathcal{P}(X),
  \forall_{R^{-1}})$ of each of the operations used to
  construct modal formulas: for the
  Boolean operations, this is straightforward, and for any modal formula
  $\phi$, we have 
  \[ h(\Box \phi) = \{ x \in X \mid \forall y \in X, \text{ if } x R y \text{ then } y \in h(\phi)\} = \forall_{R^{-1}} (h(\phi)).\]
 To finish the proof, note that $h|_V = f_c$, by definition of $\models$ on
 variables. Thus, $h$ extends $f_c$, and must therefore be equal to $\overline{f_c}$, as required. 
\end{proof}

We may now combine Corollary~\ref{cor:modal-object-duality} and 
Proposition~\ref{prop:kripke-semantics} to give a strong completeness 
theorem for the modal logic $\mathbf{K}$. A set of formulas $\Gamma$ 
is called $\mathbf{K}$-\emph{consistent} if for every finite subset $F$ 
of $\Gamma$, the formula $\neg\big(\bigwedge_{\phi \in F} \phi\big)$ 
is not provable in $\mathbf{K}$. A set of formulas $\Gamma$ is called 
\emph{satisfiable} (in a Kripke model) if there exists a Kripke model 
and a point $x$ in it such that all the formulas in $\Gamma$ hold in $x$.

\begin{theorem}\label{thm:Kcomplete}
  Any $\mathbf{K}$-consistent set of formulas is satisfiable.
\end{theorem}
\begin{proof}
  Let $\Gamma$ be a $\mathbf{K}$-consistent set of formulas and let $V$ be the
  set of propositional variables that occur in $\Gamma$. Let $A$ denote the free
  modal algebra on $V$. For any formula $\phi$, we denote by $[\phi]$ the
  corresponding element of $A$, that is, its equivalence class up to
  $\mathbf{K}$-provability. Write $(X,\tau, R)$ for the Kripke Boolean space
  dual to $A$, using the duality of Corollary~\ref{cor:modal-object-duality}. 
  Denote by $T$ the filter in $A$ generated by $\{[\phi] \mid \phi \in \Gamma\}$. 
  Note that the consistency assumption means exactly that the bottom element 
  of $A$ is not in $T$: indeed, if we would have $[\bot] \in T$, then there would 
  have to exist a finite set $F$ of formulas in $\Gamma$ such that 
  $\bigwedge_{\phi \in F}[\phi] = [\bot]$, which would mean exactly that 
  $\neg\big(\bigwedge_{\phi \in F} \phi\big)$ is provable in $\mathbf{K}$. 
  Therefore, using compactness of the space $(X, \tau, R)$, pick a point $x$ 
  that is in the intersection of the sets $\widehat{[\phi]}$, for $\phi$ in $T$. 
  Consider the canonical valuation $c\colon X \to 2^V$ such that $f_c(v) = 
  \widehat{[v]}$ for every $v \in V$; explicitly, $c(x)(v) = 1$ if, and only if, 
  $x \in \widehat{[v]}$. Then, since both $[\phi] \mapsto \widehat{[\phi]}$ and 
  $\phi \mapsto \overline{f_c}(\phi)$ are homomorphisms from $A$ to the dual 
  modal algebra of $(X,\tau,R)$ with the same value on $V$, they must be equal. 
  In particular, for any $\phi \in \Gamma \subseteq T$, we have 
  $x \in \widehat{[\phi]} = \overline{f_c}(\phi)$, so that $x \models \phi$ by 
  Proposition~\ref{prop:kripke-semantics}.
\end{proof}

The Kripke model that we construct in the proof does not depend on the
particular choice of the set $\Gamma$, but only on the set $V$ of propositional
variables that occur in $\Gamma$. The Kripke model is then obtained as the dual
space of the free modal algebra on $V$. This Kripke model is known in the modal
logic literature as the \emphind{canonical model} for the logic $\mathbf{K}$,
on the set of variables $V$. In that literature, the points of this model are
often presented as ``maximal consistent sets''.  Maximal consistent sets are
ultrafilters of the free modal algebra in disguise (see Exercise~\ref{exe:mcs}).  
The Boolean algebra with operator $(\cP(X), \forall_{R^{-1}})$ constructed here 
is known as the \emphind{canonical extension} of the modal algebra $A$, also 
see the notes at the end of this chapter.

The proof technique for completeness given in Theorem~\ref{thm:Kcomplete} is
not limited to $\mathbf{K}$. In particular, if $S$ is any set of modal
formulas, which we think of as axioms, then a new modal logic $\mathbf{K} + S$
may be defined as the set of formulas that are derivable using the rules from
$\mathbf{K}$, but allowing in addition an appeal to any substitution instance
of the formulas in $S$ without proof. Some famous examples include $\mathbf{K4}
:= \mathbf{K} + \{ \Box p \to \Box \Box p\}$ and $\mathbf{S4} := \mathbf{K} +  \{
\Box p \to \Box \Box p, \Box p \to p\}$. 

A further interesting example of a modal axiom is \emphind{Löb's axiom}
\[ \lambda := \Box (\Box p \to p) \to \Box p\ ,  \]
which, when added to $\mathbf{K}$, yields the \emphind{Gödel-Löb logic}
$\mathbf{GL}$. Gödel's name is attached to the logic because of its
relevance in \emph{provability logic}, where `$\Box p$' is interpreted as the
provability of a proposition, encoded in some theory of arithmetic. In
duality-theoretic terms, the axiom is interesting because it has a non-trivial
behavior with respect to the topology on its Kripke Boolean spaces, and a
proper understanding of the logic $\mathbf{GL}$ requires using this topology
(see Exercise~\ref{exe:GL} below for more information). In other words, the axiom $\mathbf{GL}$ is \emph{not} preserved by canonical extension.

In a slightly different direction, \emph{finite model properties}, that is,
completeness with respect to a class of \emph{finite} models, are often
desirable in modal logic, and may also be obtained using duality methods. We do
not discuss this further here, but refer to the already cited modal logic
textbooks \parencite{BRV2001,ChaZak1997} for more information.

\subsection*{Modal algebra homomorphisms and bounded morphisms}
To end our exploration of duality for modal algebras in this section, we
briefly discuss duality for homomorphisms between modal algebras and the
corresponding notion of \emphind{bounded morphism} between the dual Kripke
Boolean spaces. 

Recall that we defined a homomorphism from a modal algebra $(B,\Box_B)$ 
to a modal algebra $(A,\Box_A)$ in Definition~\ref{dfn:modal-algebra} as a 
Boolean algebra homomorphism that preserves the box operation. The dual 
of such a homomorphism should clearly be a continuous function from the 
dual space $X_A$ to the dual space $X_B$ which satisfies an additional 
property with respect to the respective Kripke relations $R_A$ and $R_B$. 
The following definition and proposition show what this property is.

\begin{definition}\label{dfn:boundedmorphism}
  Let $(X, R)$ and $(Y, S)$ be Kripke frames. A function $f \colon X \to Y$ 
  is called a \emphind{bounded morphism} if it satisfies the following two properties:
  \begin{enumerate}
    \item for every $x, x' \in X$, if $x R x'$, then $f(x) S f(x')$,
    \item for every $x \in X$, $y \in Y$, if $f(x) S y$, then there exists $x' \in X$ such that $x R x'$ and $f(x') = y$.
  \end{enumerate}
\end{definition}
The first condition in Definition~\ref{dfn:boundedmorphism} is sometimes 
called the `forth' condition, and the second condition the `back' condition. 
Bounded morphisms are also called \emphind{back-and-forth morphisms} 
in the modal logic literature. We will see in the proof of 
Proposition~\ref{prop:bounded-morphisms} below that the two conditions 
correspond to two subset inclusions.
We need the following lemma, which, in a slightly more general context, 
has been referred to as ``Esakia's Lemma'' in the literature.

\begin{lemma}\label{lem:esakia-for-boxes}
  Let $R \subseteq X \times X$ be a relation on a Boolean space $X$ 
  such that $R[x]$ is closed for every $x \in X$. Let $f \colon X \to Y$ be 
  a continuous function from $X$ to a Boolean space $Y$. For every 
  $y \in Y$,
  \[ R^{-1}\big[f^{-1}(y)\big] = \bigcap \{ R^{-1}\big[f^{-1}(K)\big] \ \mid \ K \subseteq Y \text{ clopen and } y \in K \}.\]
\end{lemma}
\begin{proof}
  The left-to-right inclusion is obvious. For the other inclusion, suppose that $x \in X$ 
  is not in the left hand side. 
  Then $f[R[x]]$ is a closed set that does not contain $y$, so, using Proposition~\ref{prop:order-normality-Priestley}, there exists a clopen set $K \subseteq Y$ which contains $y$ and is disjoint from $f[R[x]]$. The latter means that $x \not\in R^{-1}\big[f^{-1}(K)\big]$, as required.
\end{proof}

\begin{proposition}\label{prop:bounded-morphisms}
  Let $(A, \Box_A)$ and $(B, \Box_B)$ be modal algebras with dual Kripke 
  Boolean spaces $(X, \tau_X, R)$ and $(Y, \tau_Y, S)$, respectively. Let 
  $h \colon B \to A$ be a homomorphism of the underlying Boolean algebras 
  and let $f \colon X \to Y$ the dual continuous function. The following are 
  equivalent:
  \begin{enumerate}
    \item[(i)] the function $h$ is a homomorphism of modal algebras,
    \item[(ii)] the function $f$ is a bounded morphism.
  \end{enumerate}
\end{proposition}
\begin{proof}
  Recall from Proposition~\ref{prop:relation-dual-to-box} that, for any $a \in A$ and $b \in B$,
  $\widehat{\Box_A a} = \forall_{R^{-1}}[\widehat{a}]$, and
  $\widehat{\Box_B b} = \forall_{S^{-1}}[\widehat{b}]$.
  Also, for any $b \in B$, $\widehat{h(b)} = f^{-1}(\widehat{b})$, by Stone duality for homomorphisms.
  Thus, (i) is equivalent to:
  \begin{equation}\label{eqn:modal-alg-hom-dually}
    \text{for every } b \in B, f^{-1}(\forall_{S^{-1}}[\widehat{b}]) = \forall_{R^{-1}}[f^{-1}(\widehat{b})].
  \end{equation}
  Using the fact that, for any binary relation $R$ on a set $X$,
  $\forall_{R^{-1}}[U] = X \setminus R^{-1}[X \setminus U]$
  (see equation (\ref{eq:univ-direct-image})), and making a change of variable $b := \neg b$, 
  (\ref{eqn:modal-alg-hom-dually}) is equivalent to: 
  \begin{equation}\label{eqn:bdd-morphism-before-canext}
    \text{ for every } b \in B, f^{-1}(S^{-1}[\widehat{b}]) = R^{-1}[f^{-1}(\widehat{b})].
  \end{equation}
  We will now show that the condition in (\ref{eqn:bdd-morphism-before-canext}) is equivalent to:
  \begin{equation}\label{eqn:bdd-morphism-sets}
    \text{ for every } y \in Y, f^{-1}(S^{-1}[y]) = R^{-1}[f^{-1}(y)].
  \end{equation}
  The condition in (\ref{eqn:bdd-morphism-sets}) is easily seen, by unraveling the 
  definitions, to be equivalent to (ii), so this will conclude the proof of the proposition.

  Now, to prove the equivalence of (\ref{eqn:bdd-morphism-before-canext}) and
  (\ref{eqn:bdd-morphism-sets}), note that (\ref{eqn:bdd-morphism-sets})
  implies (\ref{eqn:bdd-morphism-before-canext}), since inverse images preserve
  unions. Thus, assuming (\ref{eqn:bdd-morphism-sets}), we have that for
  \emph{any} subset $U$ of $Y$,
  \[ f^{-1}(S^{-1}[U]) = \bigcup_{y \in U} f^{-1}(S^{-1}[y]) = \bigcup_{y \in U} R^{-1}[f^{-1}(y)] = R^{-1}[f^{-1}(U)].\]
  For the other direction, we use Lemma~\ref{lem:esakia-for-boxes} twice.
  Applying this lemma first to the relation $S$ and the identity function on
  $Y$, we get $S^{-1}[y]  = \bigcap_{b \in F_y} S^{-1}[\widehat{b}]$. Now
  assume (\ref{eqn:bdd-morphism-before-canext}) holds. Then, for any $y \in Y$,
  we have 
\begin{align*}
    f^{-1}(S^{-1}[y]) & = f^{-1}\big(\bigcap_{b \in F_y} S^{-1}[\widehat{b}] \big)             \\
                      & = \bigcap_{b \in F_y} f^{-1}(S^{-1}[\widehat{b}])                      \\
                      & = \bigcap_{b \in F_y} R^{-1}[f^{-1}(\widehat{b})] = R^{-1}[f^{-1}(y)],
  \end{align*}
  where we have used Lemma~\ref{lem:esakia-for-boxes} again for the last equality.
\end{proof}
\begin{remark}
  The proof of Proposition~\ref{prop:bounded-morphisms} in fact shows that $h$ 
  is a homomorphism of modal algebras if, and only if, for \emph{every} subset 
  $U \subseteq Y$, $f^{-1}(S^{-1}[U]) = R^{-1}[f^{-1}(U)]$. In other words, the property 
  of being a homomorphism of modal algebras \emph{lifts} from $h$ to the complete 
  homomorphism $f^{-1}$ from $\mathcal{P}(Y)$ to $\mathcal{P}(X)$. That is, using 
  terminology that we do not introduce further in this book, the proposition shows that 
  the property of `preserving the box operation' is a \emphind{canonical} property.
\end{remark}

\ourexercises

\begin{ourexercise}\label{exe:clopenvaluation}
    Let $(X,\tau,R,c)$ be a Kripke Boolean model. Prove that, for any modal
    formula $\phi$, the set 
    \[ \{x \in X \ \mid \ x \models \phi\} \] 
    is clopen. \hint{Use the proof of
    Proposition~\ref{prop:kripke-semantics}.}
\end{ourexercise}
\begin{ourexercise}\label{exe:kripke-examples}
  Prove the claims made in Example~\ref{exa:kripkeframe}.
\end{ourexercise}

\begin{ourexercise}\label{exe:mcs}
  A \emphind{maximal consistent set} of modal formulas with respect to the
  modal logic $\mathbf{K}$ is a $\mathbf{K}$-consistent set of formulas that is not
  properly contained in any other $\mathbf{K}$-consistent set. Fix a set of
variables $V$. Prove that the set of ultrafilters of the free modal algebra
over $V$ is in a bijection with the set $\mathrm{MCS}(V)$ of maximal consistent
sets of modal formulas whose variables lie in $V$.  
Further, by the results in this chapter, the set of
ultrafilters comes equipped with a topology and a compatible relation $R$;
describe the corresponding topology and the relation on the set
$\mathrm{MCS}(V)$. \hint{The bijection is `almost' the identity function.}
\end{ourexercise}

\begin{ourexercise}\label{exe:GL}
  This exercise concerns Löb's axiom, $\lambda := \Box (\Box p \to p) \to \Box p$. We write $\mathbf{GL}$ for the logic $\mathbf{K} + \{\lambda\}$. In particular, the last parts of this exercise outline a proof that this axiom is not \emph{canonical}: there exists a Kripke Boolean space on which the axiom is valid (that is, with respect to admissible valuations), while it is not valid on the underlying Kripke frame (that is, with respect to all valuations).
  \begin{enumerate}
    \item Prove that $\lambda$ is equivalent to $\Diam q \to \Diam(q \wedge \neg \Diam q)$, where $q := \neg p$.
    \item Show that, if $\lambda$ is valid on a Kripke frame (that is, true under any valuation), then the Kripke frame must be transitive and irreflexive.
    \item Let $(X, \tau, R)$ be a Kripke Boolean space for which $R$ is
        transitive and irreflexive. Prove
        that $\lambda$ is valid on $(X, \tau, R)$ if, and only if, for every
        clopen set $K$, $y \in K$, and $x \in X$ such that $x {R} y$, there
        exists $y_0 \in K$ such that $x {R} y_0$ and $y_0 {R} z$ implies $z
        \not\in K$, for every $z \in X$.  
    \item Consider the Boolean space $X$ which is the one-point compactification of the countable set $\{n^-, n^+ : n \in \mathbb{N}\}$; write $\infty$ for the additional point `at infinity'. Define the relation $R$ on $X$ to be the smallest transitive relation satisfying the following conditions: for any $n, m \in \mathbb{N}$, $n^+ R \infty$, $\infty R m^-$, $n^+ R (n + 1)^+$ and $(m + 1)^- R m^-$. Prove that $R$ is a compatible relation. \hint{Drawing a picture may help: the relation $R$ makes the elements $m^-$ into a descending chain, the elements $n^+$ into an ascending chain, and puts $\infty$ in the middle between the two chains, with the $m^-$ chain on top and the $n^+$ chain below. Note the similarity with Figure~\ref{fig:nplusnop}, the dual of the lattice $\mathbb{N} \oplus \mathbb{N}^\op$.}
    \item Show that $\lambda$ is valid on the Kripke Boolean space $(X, \tau, R)$ defined in the previous item.
    \item Prove that $\lambda$ is not valid on the discrete \emph{Kripke frame} $(X,R)$
        underlying the space of the previous items; that is, find a
        (non-admissible!) valuation $c$ and a point $x$ in $X$ such that $x
        \not\models \lambda$ in the discrete Kripke model $(X,R, c)$.
  \end{enumerate}
\end{ourexercise}

\section{Operators of implication type}\label{sec:generalopduality}
In this section, we study duality for \emph{binary} operators between
distributive lattices. Such operators generalize the finite-join-preserving and
finite-meet-preserving functions for which we developed a duality in
Section~\ref{sec:unaryopduality}. We focus in this section on binary operators
that are `of implication type', in a sense to be made precise below. The theory
that we develop in this section holds more generally, and can be developed for
any `order type', and for operations of any arity, also larger than $2$. We
restrict ourselves here to binary operators of implication type, because these
are the types of operators that we will use in the applications in
Chapters~\ref{chap:DomThry} and \ref{ch:AutThry}, and because the theory for
these operators is sufficiently general to see what goes on in the general case.

\begin{definition}\label{dfn:implicationtype} Let $D$, $E$, and $F$ be
  distributive lattices. A function $h \colon D \times E \to F$ is called an
  \emphind{implication-type operator}\index{operator!of implication
  type}\index{implication-type!operator}, if, for any $d, d' \in D$, $e, e' \in E$,
  the following properties hold: 
  \begin{align} \label{eq:join-to-meet-in-first} h(d \vee d', e) &= h(d, e) \wedge
               h(d', e) \quad \text{ and } \quad h(\bot, e) = \top, \\
               \label{eq:meet-to-meet-in-second} h(d, e \wedge e') &= h(d, e) \wedge h(d, e') \quad \text{ and } \quad h(d, \top) = \top.
    \end{align} 
\end{definition} 
The two properties stated in (\ref{eq:join-to-meet-in-first}) together are sometimes referred to as `$h$ \emph{sends finite joins to meets} in the first coordinate', and the two properties in (\ref{eq:meet-to-meet-in-second}) as `$h$ \emph{sends finite meets to meets} or \emph{preserves meets} in the second coordinate'. As mentioned already, variations on this definition are possible. First, when the lattices are complete, one can require that $h$ sends \emph{all} joins to meets in the first coordinate and preserves \emph{all} meets in the second coordinate; this is then called a \emphind{complete implication-type operator}. Second, we will sometimes consider monotone versions of this definition. We will then say that $h$ is an \emphind{operator of diamond type} if it sends finite joins to joins in each of its coordinates, and a \emphind{complete operator of diamond type} if it does so for all joins (also see Exercises~\ref{exe:residualscomplete}~and~\ref{ex:ToR-is-implication}).\footnote{In the literature, the operators that we call `of diamond type' are often called `normal (additive) operators', and then operators `of box type' (that is, preserving finite meets in each coordinate) are called `dual normal operators', or `multiplicative operators'. A certain amount of bookkeeping is then required to deal with operators that mix the two, so that, for instance, what we call `implication-type operator' would be called a $(\partial,1)$-dual normal operator in the terminology of, for example, \cite{GePr1}. While variations are possible, the methods developed here are limited to operations for which some order dualization of selected input coordinates yields either an operator of box type or of diamond type. For example, the methods developed here do not apply to binary operations that only preserve join in the first coordinate, and only meets in the second.}

The terminology `implication-type' stems from the fact that, in many logical theories there is an implication operation and it is typically of implication-type as defined here. However, most commonly, implication operations in logic satisfy a stronger property, namely that of being the \emphind{residual} of the conjunction operation of the logic. Residuation is a higher arity variant of the very important concept of adjunction in order theory. While residuals may be defined relative to functions of any arity, for simplicity of notation, we only give the definition in the binary case. This notion will play an important role later in the book, in particular in Chapter~\ref{ch:AutThry}.

\begin{definition}\label{dfn:residual}
Let $P_1,P_2$, and $Q$ be preordered sets and let $h \colon P_1 \times P_2 \to Q$ be a binary order-preserving function. Suppose that for each $p_1 \in P_1$, the function $h_{p_1} \colon P_2 \to Q$, defined by $h_{p_1}(p_2) := h(p_1,p_2)$ has an upper adjoint $r_{p_1}\colon Q\to P_2$. Then we obtain a binary function
\[
r\colon P_1\times Q \to P_2, (p_1,q)\mapsto r_{p_1}(q).
\]
By the definition of upper adjoint, $r$ is uniquely determined by the property that for all $p_1\in P_1$, $p_2\in P_2$ and $q\in Q$ we have 
\[
h(p_1,p_2)\leq q\ \iff\ p_2\leq r(p_1,q).
\]
We call the function $r$ the \emphind{right residual} of $h$, and we may similarly define a \emphind{left residual} $l$ of $h$ by inverting the roles of $P_1$ and $P_2$. It is uniquely determined by the property that for all $p_1\in P_1$, $p_2\in P_2$ and $q\in Q$ we have 
\[
h(p_1,p_2)\leq q\ \iff\ p_1\leq l(q,p_2).
\]
If both left and right residuals exist, then $h$ is called a \emphind{residuated operation}, and $r$ and $l$ are called its (upper) residuals.\footnote{A more complete name for $r$ would be \emph{upper right residual}, and for $l$ \emph{upper left residual}. Note that the use of `left/right' in the definition of residual refers to the `variable' coordinate of $h$ when considering the existence of an upper adjoint; for example, the right residual $r$ is the upper adjoint when the right input coordinate of $h$ is variable, and the left input coordinate is fixed. This is a different use of `left/right' than in the definition of left/right adjoint, and to avoid confusion, in this context, the lower/upper terminology is preferred for referring to adjoints. Since we mostly consider \emph{upper} residuals in this book, we will usually omit the adjective `upper'.}
\end{definition}

From the join and meet preservation properties of adjoints, it follows that a residuated operation preserves all existing joins in either coordinate, and, for example, the right residual sends existing joins in its first coordinate and existing meets in its second to meets. As a consequence, a residuated binary operation between lattices is an operator, and its right residual is of implication-type, as is the operation obtained by switching the order of the coordinates of the left residual (see  Exercise~\ref{exe:residuals}). 

Note that Definition~\ref{dfn:implicationtype} can be equivalently formulated by saying that, for
every $d_0 \in D$, the unary operation $e \mapsto h(d_0, e)$ is
finite-meet-preserving when viewed as a map from $E$ to $F$, and for every $e_0
  \in E$, the operation $d \mapsto h(d, e_0)$ is finite-meet-preserving when
viewed as a map from $D^\op$ to $F$. Implication-type operators can thus be understood as \emphind{bilinear maps} from the meet-semilattice underlying $D^\op \times E$ to the meet-semilattice underlying $F$. Note
also that, for a function $h \colon D \times E \to F$, being of implication type
is a very different property from preserving finite meets as a map from the
Cartesian product lattice $D^\op \times E$ to the lattice $F$! 
(see
Exercise~\ref{exe:difference-finite-meet-binary}). Pursuing the above analogy, a
bilinear map of abelian groups is very different from a homomorphism from the
direct product group.

\begin{example}\label{exa:BoolandHeyt} For any Boolean algebra $B$, the function
  $h \colon B \times B \to B$ defined by $h(a,b) := \neg a \vee b$ is an operator
  of implication type. This follows from the distributive law and the fact that
  $\neg$ is a homomorphism from $B^\op$ to $B$. Note that this is a rather special
  example, because this operation $h$ is in fact uniquely definable from the order
  structure of the Boolean algebra, being the right residual of the meet operation. 
  More generally, \emphind{Heyting algebras} are distributive lattices for which the 
  meet operation admits a residual (see Exercise~\ref{exe:residuals}). We will discuss 
  this example in depth in Section~\ref{sec:esakia-heyting}.   
  \end{example}

\begin{example}\label{exa:MV} Consider the lattice $L =
    [0,1]$, the real unit interval with the usual ordering, and define the operation $h \colon L \times L \to L$ by $h(a,b) := \min(1-a+b,1)$. The operation $h$ is
  an operator of implication type, and is known in the literature as the
  \emph{{\L}ukaciewicz implication}. More generally, if $A$ is any lattice-ordered
  abelian group with identity element $e$ and containing a so-called \emphind{strong unit} $u$ then an operator of implication type $h$ may be defined
  on the unit interval $[e, u]$ by $h(a,b) := (u - a + b) \wedge u$. These
  kinds of implication-type operators play a central role in the study of multi-valued
  logic and \emphind{MV-algebras}, we have omitted the
  precise definitions here; see for example \cite{CDM, Mun2011}. \end{example}

\begin{example}\label{exa:residualfromdot}
  If $\cdot \colon X \times X \to X$ is
  a binary operation on a set $X$, then it induces an implication-type operator
  $\backslash$ on the power set $\mathcal{P}(X)$, defined by
  \[ u \backslash v := \{y \in X \mid \forall x \in u, \quad x \cdot y \in v\}.\]
  This kind of operator
  of implication type plays an important role in applications of duality theory to
  automata theory, and is actually part of a residuated family of operations on $\mathcal{P}(X)$, further see Section~\ref{sec:freeprofmonoid}. Also see
  Exercise~\ref{exe:discrete-dual-res-family} below for a generalization of this example
  to ternary relations on a set $X$, which one may think of as multi-valued binary
  operations.
\end{example}

Later in this book, we will also encounter operators of implication type when we construct the lattice dual to a \emph{function space} construction; see Section~\ref{sec:funcspace}.

\subsection*{Duality for operators of implication type}
As we did with finite-meet-preserving functions in
Section~\ref{sec:unaryopduality}, we will associate to every operator of
implication type a relation between the dual spaces. The underlying ideas are
very similar to that unary case: an operator $h \colon D \times E \to F$ should
be determined by its action on points of the dual space. Indeed, for a fixed $d
  \in D$, the operation $h_d \colon E \to F$ defined by $h_d(e) := h(d,e)$ is a
finite-meet-preserving function and thus has an associated dual relation $R_d
  \subseteq X_F \times X_E$, by the results of Section~\ref{sec:unaryopduality}.
We will show in this section how this family of relations $(R_d)_{d \in D}$ can
in fact be described by a single ternary relation $R \subseteq X_D \times X_F \times X_E$. The following additional notations for ternary relations will be useful to this end.

\begin{notation}\label{ternary-relation-notation} 
Let $R \subseteq X \times Y
\times Z$ be a ternary relation between sets $X$, $Y$ and $Z$. Let $U \subseteq
X$. We use the following notation for the direct image of $U$ under the relation
  $R$: \[ R[U,\_,\_] := \{(y,z) \in Y \times Z \ | \ \text{ there exists } x
    \in U \text{ such that } R(x,y,z)\}.\] Analogously, let $V \subseteq Y$ and $W
    \subseteq Z$. Then the notations $R[\_,V,\_]$ and $R[\_,\_,W]$ are defined in
  the same way. Similarly, we define the direct image \[ R[U, V, \_] := \{z \in Z \ | \ \text { there exists } (x,y) \in U \times V \text{ such that }
    R(x,y,z)\}.\] In particular, when $U$ is a singleton set $\{x\}$, we write
  $R[x,\_,\_]$ instead of $R[\{x\},\_,\_]$, and, similarly, when moreover $V = \{y\}$, we write $R[x,y,\_]$ for the set $R[\{x\},\{y\},\_]$.
\nl{$R[U,\_,\_]$}{direct image of a subset $U$ under a ternary relation $R$}{squarebracketblanks}
\nl{$R[U,V,\_]$}{direct image of subsets $U$ and $V$ under a ternary relation $R$}{squarebracketblank}
\end{notation}

In analogy with Section~\ref{sec:unaryopduality}, we now identify necessary and sufficient conditions for a relation $R$ to be a dual relation of an operator of implication type.
\begin{definition}\label{dfn:compatible-implication}
  Let $X$, $Y$, and $Z$ be Priestley spaces. A relation $R \subseteq X \times Y
  \times Z$ is \emph{compatible} \index{compatible!relation of implication
  type}\index{implication type!compatible relation of} (of implication type)\footnote{We again use the term `compatible' here, now for a ternary relation. Whenever we need to distinguish
  it from the other types of compatibility considered earlier in this chapter, we will call it `compatible of implication type', but we will sometimes just say `compatible' in order not to make the terminologytoo heavy.} if it satisfies the following properties:
  \begin{itemize}
    \item for any $x, x' \in X$, $y, y' \in Y$, $z, z' \in Z$, if $x' \geq x$, $y' \geq y$, and $z' \leq z$, and $R(x,y,z)$, then $R(x',y',z')$;
    \item for any clopen down-set $U$ of $X$ and any clopen up-set $V$ of $Z$, the set $R[U,\_,V]$ is clopen;
    \item for every $y \in Y$, the set $R[\_,y,\_]$ is closed.
  \end{itemize}
\end{definition}
We now prove that a compatible relation of implication type between Priestley
spaces one-to-one corresponds to an implication-type operator on the dual
distributive lattices.  Analogously to the \emph{universal image} defined for
binary relations, which gave a finite-meet-preserving function for every binary
relation, we will now associate an implication-type operator to a ternary relation.

Let $X$, $Y$, and $Z$ be Priestley spaces. For any relation $R \subseteq X \times Y \times Z$, $U \subseteq X$ and $V \subseteq Z$, define the subset $U \To_R V$ of $Y$ as:
\begin{align}\label{eq:imp-of-relation}
  U \To_R V & := \{y \in Y \mid \text{for all } x \in U, z \in Z, \text{ if } R(x,y,z) \text{ then } z \in V \}.
\end{align}
\nl{$\To_R$}{operator of implication type associated with a ternary relation $R$}{}

\begin{lemma}\label{lem:ToR-is-implication}
  If $R$ is a compatible relation of implication type, then $\To_R$ defines an implication-type operator from $\ClD(X) \times \ClD(Z)$ to $\ClD(Y)$.
\end{lemma}
\begin{proof}
  The set on the right-hand-side of (\ref{eq:imp-of-relation}) is equal to $Y \setminus R[U, \_, Z \setminus V]$, and is therefore a clopen down-set, by the compatibility conditions on $R$. The required equations for an implication-type operator now follow directly from the fact that the direct image operation $(U, V) \mapsto R[U,\_,V]$ preserves joins in each coordinate (see Exercise~\ref{ex:ToR-is-implication} for the details).
\end{proof}
We now show how to recover the compatible relation from an implication-type operator.
\begin{definition} \label{dfn:dual-relation-of-implication} Let $D$, $E$, and $F$ be distributive lattices with Priestley
  dual spaces $X_D$, $X_E$, and $X_F$, respectively.  Let $h \colon D \times E \to
    F$ be an operator of implication type. We define the \emphind{dual relation}
  $R_h \subseteq X_D \times X_F \times X_E$ as follows, for any $(x, y, z) \in X_D
    \times X_F \times X_E$:
  \begin{align}\label{eq:dual-ternary-rel}
    R_h(x,y,z) \stackrel{\mathrm{def}}{\iff} & \text{ for every } d
    \in D, e \in E,                                                                                                                  \\
                                             & \text{ if } d \in F_x \text{ and } h(d, e) \in F_y, \text{ then } e \in F_z.\nonumber
  \end{align}
\end{definition}
This definition can be given a logical interpretation: if we think of prime filters as complete theories, then $(x,y,z)$ is in the relation associated to an implication operation $\To$ if the triple $(x,y,z)$ ``respects'' a sort of modus ponens rule for the implication: when a proposition $\phi$ is in the theory of $x$ and $\phi \To \psi$ is in the theory of $y$, then $\psi$ must be in the theory of $z$.
\begin{proposition}\label{prop:implication-operator-objects}
  Let $h \colon D \times E \to F$ be an implication-type operator. The dual relation $R_h \subseteq X_D \times X_F \times X_E$ is a compatible relation of implication type, and it is the unique compatible relation $R$ such that, for any $d \in D$, $e \in E$, $\widehat{h(d,e)} = \widehat{d} \To_R \widehat{e}$.
\end{proposition}
\begin{proof}
  The proof follows the same general scheme as the proof of Proposition~\ref{prop:relation-dual-to-box} in Section~\ref{sec:unaryopduality}, but requires a bit more work. It is immediate from the definition of $R_h$ that it satisfies the order compatibility condition. We now prove that, for any $d \in D$, $e \in E$, we have
  \begin{equation}\label{eq:imp-dual-to-prove}
    \widehat{h(d,e)} = \widehat{d} \To_{R_h} \widehat{e}.
  \end{equation}
  For the left-to-right inclusion, suppose that $y \in \widehat{h(d,e)}$. To prove that $y \in \widehat{d} \To_R \widehat{e}$, let $x \in \widehat{d}$ and $z \in X_E$ such that $R(x,y,z)$. Then, by definition of $R_h$, since $y \in \widehat{h(d,e)}$, we have $z \in \widehat{e}$. Thus, $y \in \widehat{d} \To_{R_h} \widehat{e}$.

  For the right-to-left inclusion, we reason contrapositively. Suppose that $y \not\in \widehat{h(d,e)}$. Consider the subset of $D$ defined by
  \[ I := \{d' \in D \mid h(d',e) \in F_y\}.\]
  Using the fact that $h$ turns joins into meets in the first coordinate, we note that $I$ is an ideal in $D$, and by assumption it does not contain $d$. Using the Prime Filter Theorem~\ref{thm:DPF}, pick $x \in X_D$ such that $x \in \widehat{d}$ and $F_x$ is disjoint from the ideal $I$. Now consider the subset of $E$ defined by
  \[ G := \{e' \in E \mid \text{there exists } d' \in F_x \text{ such that } h(d',e') \in F_y\}.\]
  Note that $G$ does not contain $e$: if $d' \in F_x$, then $d' \not\in I$, so $h(d',e) \not\in F_y$. We now show that $G$ is a filter in $E$. First, $\top_E \in G$, because $h(\top_D, \top_E) = \top_F \in F_y$, and $\top_D \in F_x$. Also, $G$ is an up-set: if $e_1 \in G$ and $e_2 \geq e_1$, then $h(d', e_2) \geq h(d', e_1)$ for any $d'$, so if $h(d', e_1) \in F_y$ then also $h(d', e_2) \in F_y$. Finally, we show $G$ is closed under binary meets. Suppose that $e_1, e_2 \in G$. Pick $d_1, d_2 \in F_x$ such that $h(d_i, e_i) \in F_y$ for $i = 1, 2$. Then $d := d_1 \wedge d_2 \in F_x$, and
  \[ h(d, e_1 \wedge e_2) = h(d,e_1) \wedge h(d,e_2) \geq h(d_1,e_1) \wedge h(d_2,e_2) \in F_y,\]
  where we use in the last inequality that $h$ is order reversing in the first coordinate.
  Using the Prime Filter Theorem~\ref{thm:DPF}, pick $z \in X_E$ such that $z \not\in \widehat{e}$, and $G \subseteq F_z$. Now, $R(x,y,z)$, using the definition (\ref{eq:dual-ternary-rel}): if $d' \in D$ and $e' \in E$ are such that $d' \in F_x$ and $h(d',e') \in F_y$, then $e' \in G$ by definition of $G$, so $e' \in F_z$ by the choice of $z$. But we also have $x \in \widehat{d}$ while $z \not\in \widehat{e}$, so by the definition of $\To_{R_h}$ (\ref{eq:imp-of-relation}), we get $y \not\in \widehat{d} \To_{R_h} \widehat{e}$. This concludes the proof of (\ref{eq:imp-dual-to-prove}).

  From (\ref{eq:imp-dual-to-prove}), we conclude in particular that $R_h$
  satisfies the second condition in the definition of compatibility: indeed, if $U$
  is a clopen down-set of $X_D$ and $V$ is a clopen up-set of $X_E$, then there
  exist $d \in D$ and $e \in E$ such that $\widehat{d} = U$ and $X_E \setminus
    \widehat{e} = V$. Now (\ref{eq:imp-dual-to-prove}) implies that $R_h[U, \_, V]$, which is the complement of $\widehat{d} \To_R \widehat{e}$, is clopen.

  For the last condition in the definition of compatibility, note that, for any $y \in X_F$, we may rewrite the definition in (\ref{eq:dual-ternary-rel}) to get \[R_h[\_,y,\_] = \bigcap \{ (X_D \setminus \widehat{d}) \times \widehat{e} \mid (d,e) \in D \times E, h(d,e) \in F_y \}, \] which is clearly closed in the product $X_D \times X_E$.

  Thus, $R_h$ is a compatible relation. It remains to prove that $R_h$ is the only compatible relation for which $\To_R$ is equal to $h$. Let $R \subseteq X_D \times X_F \times X_E$ be any compatible relation and suppose that $\widehat{d} \To_R \widehat{e} = \widehat{h(d,e)}$ for every $d \in D$, $e \in E$. Let $(x,y,z) \in X_D \times X_F \times X_E$ be arbitrary. First, if $R(x,y,z)$, then for any $d \in F_x$ and $e \in E$ such that $h(d,e) \in F_y$, we have $y \in \widehat{d} \To_R \widehat{e}$, so that $e \in F_z$; thus, $R_h(x,y,z)$. Conversely, suppose that $R(x,y,z)$ does not hold. This means that $(x,z)$ is not in the set $C := R[\_,y,\_]$, which is closed and an up-set in the product space $X_D \times X_E^\op$, using that $R$ is compatible. By Proposition~\ref{prop:order-normality-Priestley} applied to $X_D \times X_E^\op$, which is a Priestley space by Exercise~\ref{exe:Priestley-product-space}, there exists a clopen down-set $L \subseteq X_D \times X_E^\op$ such that $(x,z) \in L$ and $L$ is disjoint from $C$. Using Exercise~\ref{exe:Priestley-binary-product}, there then also exist $d \in D$ and $e \in E$ such that $(x,z) \in \widehat{d} \times \widehat{e}^c$, and $\widehat{d} \times \widehat{e}^c$ is disjoint from $C$. It follows that $y$ is in $\widehat{d} \To_R \widehat{e}$: whenever $R(x',y,z')$ for some $x' \in \widehat{d}$, we must have $z' \in \widehat{e}$, since $z' \in \widehat{e}^c$ would contradict that $\widehat{d} \times \widehat{e}^c$ is disjoint from $C = R[\_,y,\_]$. Now, by the assumption that $\widehat{d} \To_R \widehat{e} = \widehat{h(d,e)}$, we get that $y \in \widehat{h(d,e)}$. But then $R_h(x,y,z)$ does not hold, since $d \in F_x$ and $h(d,e) \in F_y$ but $e \not\in F_z$.
\end{proof}
Just as in the case of modal algebras (Definition~\ref{dfn:modal-algebra}), one may define a \emph{distributive lattice with implication-type operator} to be a pair $(D, \To)$, where $D$ is a distributive lattice and $\To \colon D \times D \to D$ is an implication-type operator on $D$. Then, in the same way as was done for modal algebras in Corollary~\ref{cor:modal-object-duality}, it follows easily from Proposition~\ref{prop:implication-operator-objects} that any distributive lattice with implication-type operator $(D,\To)$ is isomorphic to $(\ClD(X), \To_R)$, where $X$ is the Priestley dual space of $D$ and $R$ is the dual relation $R_\To \subseteq X \times X \times X$ of the operator $\To$. Moreover, $(X,R)$ is the unique such object, up to relation-preserving isomorphism of Priestley spaces (that is, order-homeomorphisms $f$ such that both $f$ and $f^{-1}$ preserve the relation). Thus, Proposition~\ref{prop:implication-operator-objects} gives an `object correspondence' between distributive lattices with implication-type operators and pairs $(X,R)$ with $X$ a Priestley space and $R$ a compatible relation.

In order to turn this object correspondence into a full-fledged duality for distributive lattices equipped with implication-type operators, one needs to study when distributive lattice homomorphisms preserve an implication-type operator, and what this says about the dual functions between Priestley spaces. One will then obtain a notion of \emphind{bounded morphism} with respect to ternary relations, analogous to Definition~\ref{dfn:boundedmorphism} in the previous section. Instead of performing this general analysis here, we defer it to the specific cases where we need it. In particular, in Theorem~\ref{thm:esakia-duality} we will establish a full dual equivalence for Heyting algebras, which is the special case where the implication-type operator on the distributive lattice is the residual of the $\wedge$ of the lattice. Further, in Sections~\ref{sec:freeprofmonoid}~and~\ref{sec:EilReittheory} of Chapter~\ref{ch:AutThry}, we will develop quotient-subspace duality for implication-type operators in the special setting of residuation operations on a Boolean algebra, which are in particular implication-type operators.

\ourexercises

\begin{ourexercise}\label{exe:residuals}
Let $L_1,L_2$, and $L$ be lattices and $h \colon L_1 \times L_2 \to L$
a residuated binary function. 
\begin{enumerate}
  \item Show that $h$ is a binary operator of diamond type. \hint{Use Exercise~\ref{exe:adjunctions}.}
  \item Show that its right residual $r$ is an operation of implication type, and that 
           its left residual $l$ is an operation of implication type if we reverse the order 
           of the coordinates. 
 \item Show that if $P=P_1=P_2$ and $h$ is commutative, then the right and left 
          residuals are equal, up to switching the order of the input coordinates. That 
          is, for $p\in P$ and $q\in Q$ we have $r(p,q)=l(q,p)$.
\item Show that if $B$ is a Boolean algebra, then the binary meet operation $\wedge$
         is residuated, and its residual is the logical implication operation $\to$, defined, for $a, b \in B$, by $a\to b:=\neg a\vee b$.
\end{enumerate}                           
\end{ourexercise}
\begin{ourexercise}\label{exe:residuals-Galois}
  Let $L_1,L_2$, and $L$ be lattices and $h \colon L_1 \times L_2 \to L$
  a residuated binary function. Show that for each $a\in L$, the pair of functions $r(-,a)\colon L_1\to L_2$ and $l(a,-)\colon L_2\to L_1$ form a Galois connection.
  \end{ourexercise}
  
\begin{ourexercise}\label{exe:residualscomplete}
Let $L_1, L_2$ and $L$ be complete lattices and $h \colon L_1 \times L_2 \to L$ an order-preserving function. We say that $h$ is a \emphind{complete  operator of diamond type} provided that it preserves all joins in each coordinate, that is, for any $U \subseteq L_1$, $b \in L_2$, $V \subseteq L_2$, $a \in L_1$ we have 
\[ f\big(\bigvee U, b\big) = \bigvee_{u \in U} f(u, b) \quad \text{ and }  \quad f\big(a,\bigvee V\big) = \bigvee_{v \in V} f(a, v)\ .\] 
\begin{enumerate}
  \item Show that $h$ is a complete operator of diamond type if, and only if, $h$ is residuated.
  \item Show that $h$ is residuated if, and only if, $h$ has a right residual that is a complete operator of implication type.
\end{enumerate}
\end{ourexercise}

\begin{ourexercise}\label{exe:difference-finite-meet-binary}
  This exercise shows the difference, for a function $h \colon D \times E \to F$,
  between `operator of implication-type' and `finite-meet-preserving as a function
  from the product lattice $D^\op \times E$ to $F$'.  
  \begin{enumerate}
  \item Give an example of a function $h \colon \mathbf{2} \times \mathbf{2} \to \mathbf{2}$
          which is an operator of implication type, but does not preserve finite meets as
          a map from the product lattice $\mathbf{2}^\op \times \mathbf{2}$ to
          $\mathbf{2}$.  \hint{Boolean implication is an operator of implication
          type.}  \item Give an example of a function $h \colon \mathbf{2} \times
            \mathbf{2} \to \mathbf{2}$ that preserves finite meets as a map from
          $\mathbf{2}^\op \times \mathbf{2}$ to $\mathbf{2}$, but which is not an operator
          of implication type. \hint{The only equations that can fail are the ones concerning $\top$ and $\bot$ in Definition~\ref{dfn:implicationtype}. } \item Suppose that $h
            \colon D \times E \to F$ is finite-meet-preserving as a function from $D^\op
            \times E$ to $F$. Prove that $h$ is an operator of implication type if, and only
          if, $h(\bot,e) = \top$ and $h(d, \top) = \top$ for all $d \in D$, $e \in E$.
  \end{enumerate}
\end{ourexercise}

\begin{ourexercise}
  Verify that the operations defined in Examples~\ref{exa:MV} and
  \ref{exa:residualfromdot} are indeed operators of implication type.
\end{ourexercise}

\begin{ourexercise}\label{exe:discrete-dual-res-family}
This exercise generalizes Example~\ref{exa:residualfromdot}, which showed how to use a binary operation to obtain an implication-type operator. Here we show how to obtain an implication-type operator obtained from any ternary relation. 
  Let $X$ be a
  set and $R\subseteq X^3$ be a ternary relation on $X$. We use $R$ to define
  three binary operations $\cdot, /, \backslash \colon \cP(X)^2 \to \cP(X)$, as
  follows:
  \begin{align*}
    u \cdot t      & := \{z \in X \mid  \exists x \in u, y \in t \text{ such that } R(x,y,z)\},                \\
    u \backslash v & := \{y \in X \mid \forall x \in u, z \in X, \text{ if } R(x,y,z) \text{ then } z \in v\}, \\
    v / t          & := \{x \in X \mid \forall y \in t, z \in X, \text{ if } R(x,y,z) \text{ then } z \in
    v\}.
  \end{align*}
  Note that, when $R$ is the graph of a function $\cdot \colon X^2\to X$, the definitions of $\backslash$ and $\slash$ are the same as in Example~\ref{exa:residualfromdot}.
  Show that these operations form a residuated family, that is, for all $s,t,u\in
    \cP(X)$ we have
  \[ u \cdot t \subseteq v\ \iff\ t\subseteq u\backslash v\ \iff\ u\subseteq v/t.  \]
\nl{$u \cdot v$}{binary operator associated to a ternary relation $R$, equal to $R[u,v,\_]$}{}
\nl{$u \backslash v$}{operator of implication type associated to a ternary relation $R$, left residual of $\cdot$}{}
\nl{$u \slash v$}{operator of implication type associated to a ternary relation $R$, right residual of $\cdot$}{}
\end{ourexercise}

\begin{ourexercise}\label{ex:ToR-is-implication} This exercise leads to the verification of Lemma~\ref{lem:ToR-is-implication} and explains the origin of the first two conditions of the definition of compatible relation of implication type. In particular, it gives a generalization of the discrete duality between posets and down-set lattices described in Exercise~\ref{exe:discrete-duality}. For the definition of complete operator of diamond type (see Exercise~\ref{exe:residualscomplete}).
  \begin{enumerate}
    \item \label{itm:discrete-binary-op-duality}Let $X$, $Y$, and $Z$ be posets. Show that under the discrete duality between posets and down-set lattices, complete operators of diamond type
          \[
            f\colon\cD(X)\times\cD(Y)\to\cD(Z)
          \]
          correspond to ternary relations $R\subseteq X\times Y\times Z$ satisfying for any $x, x' \in X, y, y' \in Y$, and $z, z' \in Z$
          \[
            \big[ \ x \leq x', y \leq y', z' \leq z, \text{ and }R(x,y,z)\ \big] \implies R(x',y',z').
          \]
    \item Show that for $R\subseteq X\times Y\times Z$, $U\in\cP(X)$, and  $W\in\cP(Z)$, we have
          \[
            U \To_R W=Y\setminus R[U,{\_}\,,Z\setminus W]
          \]
          and this operation is the right residual of the complete operator of diamond type defined by 
          $$(U,V)\mapsto R[U,V,\_\,] \text{ for }U\in\cP(X) \text{ and }  V\in\cP(Y)\ .$$ 
      \item Show that the operation $\To_R$ restricts correctly to the down-set lattices when $R$ satisfies the property in (\ref{itm:discrete-binary-op-duality}).
    \item Verify the details of the proof of Lemma~\ref{lem:ToR-is-implication}.
  \end{enumerate}
\end{ourexercise}

\section{Heyting algebras and Esakia duality}\label{sec:esakia-heyting}
In this section we introduce the class of \emphind{Heyting algebras}. Using the fact that Heyting algebras are a subclass of distributive lattices, equipped with a very special operator of implication type, we will derive a duality for them, using Priestley duality and the general results of the previous section.

\begin{definition}
  Let $L$ be a distributive lattice and let $a\in L$. If the set $\{c \in L : a \wedge c=\bot\}$ has a maximum in $L$, we denote it by $a^*$, and we call $a^*$ the \emphind{pseudocomplement} of $a$ with respect to $b$. Further, if, for $a, b \in L$, the set $\{c \in L : a \wedge c \leq b\}$ has a maximum in $L$, we denote it by $a \hto b$, and we call $a \hto b$ the \emphind{relative pseudocomplement} of $a$ with respect to $b$.
\nl{$a \hto b$}{relative pseudocomplement of $a$ with respect to $b$, also known as the Heyting implication or Heyting arrow}{}

  A \emphind{Heyting algebra} $H$ is a distributive lattice such that the relative pseudocomplement $a \hto b$ exists for any $a, b \in H$.

  Let $H$ and $K$ be Heyting algebras. A \emphind{Heyting homomorphism} is a lattice homomorphism $h \colon H \to K$ which moreover preserves relative pseudocomplements, that is, $h(a \hto b) = h(a) \hto h(b)$ for all $a, b \in H$.
\end{definition}

\begin{remark}
Pseudocomplement is a weakening of the notion of complement in a Boolean algebra. Being a Heyting algebra is a notion of `hereditary' pseudocomplementedness, in the sense that a distributive lattice $L$ is a Heyting algebra if, and only if, for each $b\in L$, the lattice ${\uparrow}b$, which is an unbounded sublattice of $L$, but is a bounded lattice in its own right, is pseudocomplemented. Also, any Boolean algebra is a Heyting algebra in which $a\hto c$ is given by $\neg a\vee b$; see Exercise~\ref{exer:pseudocomplement}. Finally, there is an order-dual notion to Heyting algebra, that we call \emphind{co-Heyting algebra}, in which the set $\{c \in L \colon b \leq a \vee c\}$ has a minimum for every $a, b \in A$; see Exercise~\ref{exe:coheyting}.
\end{remark}

Note that a distributive lattice $H$ is a Heyting algebra if, and only if, for any $a \in H$, the function $a \wedge (-) \colon H \to H$ has an upper adjoint. Indeed, a relative pseudocomplement of $a$ with respect to $b$ is an element $a\hto b$ of $H$ such that, for any $c \in H$,
\begin{equation}\label{eq:heyting-adj} a \wedge c \leq b \iff c \leq a \hto b.
\end{equation}
Interestingly, Heyting algebras may also be characterized as those distributive lattices 
for which the embedding map into their Boolean envelope has an upper adjoint (see Exercise~\ref{exe:heyting-S4}).

There is an interesting, and at first potentially confusing, tension in the
definition of Heyting algebras. While we chose to define a Heyting algebra as a
distributive lattice with an additional property, we often view it in practice as a
distributive lattice that comes equipped with the additional \emph{structure} of
a Heyting implication. It is possible to give an
alternative, equivalent definition of Heyting algebras that is purely
equational, see for example \cite[II.1, Example 11]{BurSan2000}.
It is important to remember here that any given distributive
lattice admits \emph{at most} one Heyting implication. This situation is
analogous to the more familiar case of Boolean algebras: a distributive lattice
admits at most one Boolean negation, which we then often regard as additional
structure. However, a difference between the two situations is that, while a
distributive lattice homomorphism between Boolean algebras always preserves
Boolean negation (see Exercise~\ref{exe:distuncomp}), not every distributive lattice
homomorphism between Heyting algebras preserves the Heyting implication
(see Exercise~\ref{exe:morphismsHAdifferent}).

\begin{example}\label{exa:HAs}
  Any Boolean algebra is a Heyting algebra, with the implication definable as $a \to b := \neg a \vee b$; this example was already mentioned in Example~\ref{exa:BoolandHeyt} in the previous section (also see Exercise~\ref{exer:pseudocomplement}). Any frame is a Heyting algebra, because the functions $a \wedge (-)$ in a frame preserve arbitrary joins, and therefore have an upper adjoint by the adjoint functor theorem for complete lattices (see Exercise~\ref{exe:adjointexistsiff}). However, not every homomorphism between frames preserves the Heyting implication (see Exercise~\ref{exe:morphismsHAdifferent}).  The Heyting implication in a frame $F$ can be computed explicitly,  for any $a, b \in F$, as
  \[ a \hto b = \bigvee \{ c \in F \mid a \wedge c \leq b\}.\]
  In particular, the following are examples of Heyting algebras:
  \begin{enumerate}
    \item Any finite distributive lattice is a Heyting algebra.
    \item The open set lattice of any topological space is a Heyting algebra, in which $a \hto b$ is naturally interpreted using the interior operator (see Exercise~\ref{exe:HAinterior}).
          \item\label{itm:HAdownsets} Let $(P, \leq)$ be a preorder. The lattice of down-sets $\Down(P)$ is a complete Heyting algebra, as follows from the preceding example. Given two down-sets $U, V \in \Down(P)$, the relative pseudocomplement of $U$ with respect to $V$ may be calculated as follows (see Exercise~\ref{exe:downsetsHA}):
          \begin{equation}\label{eq:downsetshto}
            U \hto V = \{ p \in P \mid \forall q \in {\downarrow}p, \text{ if } q \in U \text{ then } q \in V\}.
          \end{equation}
          Note that equation (\ref{eq:downsetshto}) can also be written as $U \hto V = P \, {\setminus} \, {\uparrow}(U \setminus V)$.
  \end{enumerate}
\end{example}

Let $H$ be a Heyting algebra. For any $a, b \in H$ and $S, T \subseteq H$ such that $\bigvee S$ and $\bigwedge T$ exist, we have:
\begin{align}
  a \wedge \Big(\bigvee S\Big) & = \bigvee_{s \in S} (a \wedge s), \label{eq:heytingframe}    \\
  b \hto \Big(\bigwedge T\Big) & = \bigwedge_{t \in T} (b \hto t), \label{eq:heytingsecondco} \\
  \Big(\bigvee T\Big) \hto b   & = \bigwedge_{t \in T} (t \hto b). \label{eq:heytingfirstco}
\end{align}
The first two equations follow from equation (\ref{eq:heyting-adj}), because
lower adjoints preserve existing suprema, and upper adjoints preserve existing
infima (see Exercise~\ref{exe:adjunctions}). For a proof of the third equation (see
Exercise~\ref{exe:htofirstcoordinate}). Thus, the operation $\to$ on any Heyting
algebra is an \emph{operator of implication type}\index{operator!of implication type} in the sense of Section~\ref{sec:generalopduality}.

We remark that (\ref{eq:heytingframe}) shows that any \emph{complete} Heyting algebra is a frame. Thus, in light of Example~\ref{exa:HAs}, a lattice is a frame if, and only if, it is a complete Heyting algebra. However, we emphasize that the notion of \emph{morphism} depends on whether we view such a lattice as a frame or as a complete Heyting algebra: frame homomorphisms are required to preserve finite meets and arbitrary joins, whereas Heyting algebra homomorphisms are required to preserve finite meets, \emph{finite} joins, and the relative pseudocomplement. Exercise~\ref{exe:morphismsHAdifferent} guides you towards specific examples showing the difference.

\subsection*{Esakia spaces}
We now derive from the duality for implication-type operators developed in the previous section a duality for Heyting algebras. This duality is originally due to \cite{Esa1974} and is known as Esakia duality in the literature. An English translation of Esakia's original 1985 monograph was published in 2019 \parencite{Esakia2019}.

From the perspective of the previous section, if a distributive lattice $L$ is a Heyting algebra, then the implication-type operator $\hto$ has a dual ternary relation $R_{\hto}$, as defined in Definition~\ref{dfn:dual-relation-of-implication}. Since $\hto$ is definable from the distributive lattice structure of $L$, the relation $R_{\hto}$ ought to also be definable directly from the dual space of $L$. As a first step towards deducing Esakia duality from Priestley's, we now show that this is indeed the case.
\begin{lemma}\label{lem:Rhto-from-order}
  Let $L$ be a Heyting algebra and let $X$ be the dual Priestley space of $L$. Then the ternary relation $R_{\hto}$ dual to the Heyting implication $\hto$ is given by
  \[ R_{\hto} = \{(x,y,z) \in X^3 \mid z \leq x \text{ and } z \leq y \}.\]
\end{lemma}
\begin{proof}
  Recall that the definition of $R_{\hto}$ (Definition~\ref{dfn:dual-relation-of-implication}) says: $R_{\hto}(x,y,z)$ if, for any $a, b \in L$, if $a \in F_x$ and $a \to b \in F_y$, then $b \in F_z$. We show that this is equivalent to the condition that $F_x \subseteq F_z$ and $F_y \subseteq F_z$.

  First suppose that $R_{\hto}(x,y,z)$. Let $a \in F_x$. Note that $\top \leq a \hto a$ because $\top \wedge a \leq a$. Thus, $a \hto a \in F_y$. Hence, $a \in F_z$. Now let $c \in F_y$. Note that $c \leq \top \hto c$ because $c \wedge \top \leq c$. Thus, $\top \hto c \in F_y$, and also $\top \in F_x$. Hence, $c \in F_z$. Conversely, suppose that $F_x \subseteq F_z$ and $F_y \subseteq F_z$. Let $a \in F_x$ and $b \in L$ be such that $a \to b \in F_y$. Then $a \in F_z$ and $a \hto b \in F_z$, so $a \wedge (a \hto b) \in F_z$. Now note that $a \wedge (a \hto b) \leq b$, since $a \hto b \leq a \hto b$. Thus, $b \in F_z$.
\end{proof}
For the remainder of this section, if $X$ is a Priestley space, we write $E$ for the ternary relation defined by
\[ E := \{(x,y,z)\in X^3 \ \mid \ z \leq x \text{ and } z \leq y \}.\]
Lemma~\ref{lem:Rhto-from-order} says that, if $X$ is dual to a Heyting algebra, then $E = R_{\hto}$. In this case, $E$ is a compatible relation in the sense of Definition~\ref{dfn:compatible-implication}. We now reverse the question, and ask, when is the relation $E$, which can be defined on any Priestley space, a compatible relation? We note immediately that $E$ satisfies the first, order-theoretic, condition in Definition~\ref{dfn:compatible-implication}. Also, for any $y \in Y$, we have
\[ E[\_,y,\_] = ({\geq_X}) \cap (X \times {\downarrow} y),\]
which is always a closed subset of $X \times X$, since in any compact ordered space, the (reverse) order is closed, and the downward closure of a point is closed (Proposition~\ref{prop:cos-closed}).
Thus, the relation $E$ is compatible if, and only if, it satisfies the second
condition in Definition~\ref{dfn:compatible-implication}: for any clopen
down-set $U$ of $X$ and any clopen up-set $V$ of $X$, the set $E[U,\_,V]$ must
be clopen. From the definition of $E$, we note that $E[U,\_,V] = {\uparrow}
  ({\downarrow} U \cap V) = {\uparrow} (U \cap V)$ when $U$ is a down-set. With
these preliminary considerations, we are now ready to characterize the Priestley
spaces that are dual to Heyting algebras.
\begin{proposition}\label{prop:esakiaspace}
  Let $L$ be a distributive lattice with dual Priestley space $X$. The following are equivalent:
  \begin{enumerate}
    \item[(i)] $L$ is a Heyting algebra;
    \item[(ii)] for any clopen subset $K$ of $X$, the generated up-set ${\uparrow} K$ is clopen;
    \item[(iii)] for any open subset $U$ of $X$, the generated up-set ${\uparrow} U$ is open.
  \end{enumerate}
\end{proposition}
\begin{proof}
  We show that (i) and (ii) are equivalent; the equivalence of (ii) and (iii) is left as Exercise~\ref{exe:Esakia-equiv}.
  First, if $L$ is a Heyting algebra, then the relation $R_{\hto} = E$ is compatible. By the remarks preceding the proposition, this means that ${\uparrow} (U \cap V)$ is clopen for any clopen up-set $U$ and clopen down-set $V$. Now recall from Lemma~\ref{lem:Priestleybase} that any clopen set of $X$ is a finite union of sets of the form $U \cap V$, with $U$ a clopen up-set and $V$ a clopen-down-set. Thus, it follows that ${\uparrow} K$ is clopen for any clopen set $K$.
  For the converse, suppose that (ii) holds. By the above remarks, the relation $E$ is in particular compatible, and therefore by Lemma~\ref{lem:ToR-is-implication} defines an implication-type operator $\To_E$ on $\ClD(X)$. Unraveling the definition of $\To_E$ in this case, we see that, for any down-sets $U$ and $V$ of $X$,
  \begin{align*}
    U \To_E V & = \{ y \in X \ \mid \ \forall x \in U, z \in X, \text{ if } z \leq x \text{ and } z \leq y \text{ then } z \in V \} \\
              & = \{ y \in X \ \mid \ \forall z \in X, \text{ if } z \leq y \text{ and } z \in U, \text { then } z \in V\}.
  \end{align*}
  As remarked in Example~\ref{exa:HAs}.\ref{itm:HAdownsets}, this is the Heyting implication on the (complete) Heyting algebra $\mathcal{D}(X)$. Since $\To_E$ is well-defined on the bounded sublattice $\ClD(X)$ of $\mathcal{D}(X)$, it now follows from Exercise~\ref{exe:heytingsub} that $\ClD(X)$ is a Heyting algebra. By Priestley duality, $L$ is isomorphic to $\ClD(X)$, and is therefore also a Heyting algebra.
\end{proof}
A Priestley space $X$ that satisfies the equivalent conditions in Proposition~\ref{prop:esakiaspace} is called an \emphind{Esakia space}. Note that the proof of Proposition~\ref{prop:esakiaspace} in particular implies that, if $L$ is a Heyting algebra with dual Esakia space $X$, then for any $a, b \in L$,
\[ \widehat{a \to b} = \widehat{a} \to \widehat{b},\]
that is, the function $a \mapsto \widehat{a}$ is a Heyting homomorphism from $L$ to $\mathcal{D}(X)$.

\begin{proposition}\label{prop:esakiamorphisms}
  Let $L$ and $M$ be Heyting algebras with Priestley dual spaces $X$ and $Y$, respectively. Let $h \colon M \to L$ be a lattice homomorphism with dual Priestley morphism $f \colon X \to Y$; that is, for all $a \in M$, we have $\widehat{h(a)} = f^{-1}(\widehat{a})$. The following are equivalent:
  \begin{enumerate}
    \item[(i)] for all $y \in Y$, $f^{-1}({\uparrow} y) \subseteq {\uparrow} f^{-1}(y)$.
    \item[(ii)] the homomorphism $f^{-1} \colon \Down(Y) \to \Down(X)$ is a Heyting homomorphism;
    \item[(iii)] the homomorphism $h$ is a Heyting homomorphism, that is, preserves $\hto$.
  \end{enumerate}
\end{proposition}

\begin{proof}
  (i) $\Rightarrow$ (ii). Since $f$ is order preserving, $f^{-1}$ sends down-sets to down-sets, and thus it is a (complete) lattice homomorphism from $\Down(Y)$ to $\Down(X)$. We just need to show that $f^{-1}$ preserves the relative pseudocomplement. Let $U, V \in \Down(Y)$. By adjunction, we have $U\cap(U\hto V)\subseteq V$ and thus $f^{-1}(U) \cap f^{-1}(U \to V) \subseteq f^{-1}(V)$. Again by adjunction we obtain $f^{-1}(U \to V) \subseteq f^{-1}(U) \to f^{-1}(V)$.

  For the reverse inclusion, we use the definition of the Heyting implication on $\Down(X)$, see (\ref{eq:downsetshto}). Let $x \in f^{-1}(U) \to f^{-1}(V)$. We need to show that $f(x) \in U \to V$. To this end, let $y \leq f(x)$ and $y \in U$. By (i), we then have $x \in {\uparrow} f^{-1}(y)$, so pick $x' \leq x$ such that $f(x') = y$. Since $y \in U$, we have $x' \in f^{-1}(U)$. Thus, since $x \in f^{-1}(U) \to f^{-1}(V)$ by assumption, we get $x' \in f^{-1}(V)$. Hence, $y = f(x') \in V$, as required.

  (ii) $\Rightarrow$ (iii). Since $\widehat{(-)}$ is a Heyting homomorphism and
  $f^{-1}$ is a Heyting homomorphism, their composition is a Heyting
  homomorphism. By Priestley duality, this composition sends $a$ to
  $\widehat{h(a)}$. It follows that $h$ itself is a Heyting homomorphism, since
  $\widehat{(-)}$ is a lattice isomorphism between $L$ and the clopen down-sets
  of $X$.

  (iii) $\Rightarrow$ (i). Suppose $y\leq f(x)$. By the definition of the Priestley topology on $Y$ and the fact that it is $T_0$, it follows that
\[
\{y\}=\bigcap\{\widehat{b}\cap(\widehat{c})^c\mid y\in\widehat{b}\text{ and }y\not\in\widehat{c}\,\}.
\]
Therefore
\[
f^{-1}(y)=\bigcap\{f^{-1}(\widehat{b})\cap [f^{-1}(\widehat{c})]^c\mid  b,c\in M, y\in\widehat{b},\text{ and }y\not\in\widehat{c}\,\}.
\]
Also, 
\[
{\downarrow}x=\bigcap\{\widehat{a}\mid a\in L, x\in\widehat{a}\},
\]
 so ${\downarrow}x\cap f^{-1}(y)$ is equal to
\[
\bigcap\{\,\widehat{a}\,\cap f^{-1}(\widehat{b})\cap [f^{-1}(\widehat{c})]^c\mid a\in L,\ b,c\in M, x\in\widehat{a}, y\in\widehat{b},\text{ and }y\not\in\widehat{c}\,\}.
\]
Note that the family of sets we take the intersection of is down directed. Now showing that $x\in{\uparrow} f^{-1}(y)$ is equivalent to showing that ${\downarrow}x\cap f^{-1}(y)$ is non-empty, which, by compactness, is equivalent to showing that each set in the family of sets we take the intersection of is non-empty.
  
 To this end, let $a\in L$ with $x\in\widehat{a}$ and $b,c\in M$ with $y\in\widehat{b}$ and $y\not\in\widehat{c}$. Now $y\leq f(x)$, $y\in\widehat{b}$ and $y\not\in\widehat{c}$ implies, by (\ref{eq:downsetshto}), that $f(x)\not\in\widehat{b}\hto\widehat{c}$. That is, by (iii),
\begin{align*}
x \not\in f^{-1}(\widehat{b}\hto\widehat{c})=\widehat{h(b\hto c)}
   =\widehat{h(b)\hto h(c))}=f^{-1}(\widehat{b})\hto f^{-1}(\widehat{c}).
\end{align*}
It follows that $\widehat{a}\not\subseteq f^{-1}(\widehat{b})\hto  f^{-1}(\widehat{c})$. By adjunction we obtain $\widehat{a}\cap f^{-1}(\widehat{b})\not\subseteq f^{-1}(\widehat{c})$ or, equivalently, $\widehat{a}\cap f^{-1}(\widehat{b})\cap [f^{-1}(\widehat{c})]^c\neq\emptyset$ as required.
\end{proof}

Order preserving functions satisfying condition (i) in Proposition~\ref{prop:esakiamorphisms} are called p-morphisms or bounded morphisms in the literature.  Note that the definition is a special case of the definition for Kripke frames given in Definition~\ref{dfn:boundedmorphism}. This is because, being a Heyting algebra homomorphism for the lattices of clopen up-sets is the same as being a $\Box$ homomorphism of the Boolean envelopes of these lattices for the modal operation given by $\leq$ viewed as a Kripke relation (see Exercise~\ref{exe:heyting-S4}). 
 \begin{definition}
  An order-preserving function $f \colon X \to Y$ between posets is called a \emphind{p-morphism} or \emphind{bounded morphism} if, for any $x \in X$ and $y \in Y$, if $y \leq f(x)$, then there exists $x' \leq x$ such that $f(x') = y$.

  We denote by $\cat{Esakia}$ the category of Esakia spaces with continuous p-morphisms.
\end{definition}
Observe that a function $f\colon X\to Y$ between posets satisfies the reverse inclusion of (i) in Proposition~\ref{prop:esakiamorphisms}, that is, $f^{-1}({\uparrow} y) \supseteq {\uparrow} f^{-1}(y)$ for all $y \in Y$ if, and only if, it is order preserving. Therefore, an equivalent definition of $p$-morphism is: a continuous function $f \colon X \to Y$ such that $f^{-1}({\uparrow} y) = {\uparrow} f^{-1}(y)$ for all $y \in Y$.

The following duality theorem for Heyting algebras is now an immediate consequence of Priestley duality in Chapter~\ref{ch:priestley}, Proposition~\ref{prop:esakiaspace}, and Proposition~\ref{prop:esakiamorphisms}.
\begin{theorem}\label{thm:esakia-duality}
  The category $\cat{HA}$ of Heyting algebras with Heyting homomorphisms is dually equivalent to the category $\cat{Esakia}$ of Esakia spaces with continuous p-morphisms.
\end{theorem}
Just as duality for modal algebras yields a semantics for modal logic (see Section~\ref{sec:kripke}), a similar Kripke-style semantics for intuitionistic propositional logic can be obtained from this duality between Heyting algebras and Esakia spaces.

An \emphind{intuitionistic formula} is a term in the signature of Heyting algebras, that is, built from propositional variables, using the operations $\bot$, $\top$, $\to$, $\vee$ and $\wedge$. If $\phi$ is an intuitionistic formula with variables in a set $V$, we define $\models \phi$ to mean: for every function $v \colon V \to L$, where $L$ is a Heyting algebra, we have $\overline{v}(\phi) = \top$. Here, $\overline{v}(\phi)$ is the unique interpretation of a formula $\phi$ extending the given interpretation $v$ of variables.

If $(X, \leq)$ is a partially ordered set, then an \emph{admissible valuation}\index{valuation!admissible, for intuitionistic Kripke models} is a function $v \colon V \to \mathcal{D}(X)$, assigning to every variable a down-set of $X$. We say that an intuitionistic formula $\phi$ with variables in $V$ \emph{holds} at a point $x \in X$ if $x \in \overline{v}(\phi)$. A tuple $(X, \leq, v)$ where $(X,\leq)$ is a poset and $v$ is an admissible valuation is known in the literature as an \emphind{intuitionistic Kripke model}. An intuitionistic formula $\phi$ is \emph{satisfiable} if there exists an intuitionistic Kripke model such that $\overline{v}(\phi)$ is non-empty, and \emph{consistent}\index{consistent!intuitionistic formula} if there exists an interpretation to a Heyting algebra such that $\overline{v}(\phi) \neq \bot$.

\begin{theorem}\label{thm:int-completeness}
  Every consistent intuitionistic formula is satisfiable.
\end{theorem}
\begin{proof}
  Let $\phi$ be a consistent intuitionistic formula and let $v \colon V \to H$ be a valuation of the variables occurring in $\phi$ such that $v(\phi) \neq \bot$. Let $X$ the Esakia space dual to $H$. We define the admissible valuation $v' \colon V \to \mathcal{D}(X)$ by $v'(p) := \widehat{v(p)}$ for every $p \in V$. Since $\widehat{(-)}$ is a homomorphism of Heyting algebras, we have $\overline{v'}(\phi) = \widehat{v(\phi)}$, and since $\widehat{(-)}$ is injective, we get $\overline{v'}(\phi) \neq \emptyset$, as required.
\end{proof}

We show now that, with only a little extra work, we may use the above proof technique to show that every consistent intuitionistic formula is satisfiable in a \emph{finite} intuitionistic Kripke model. This argument proves the so-called \emphind{finite model property} of intuitionistic propositional logic and is due to \cite{McKTar1948}.

\begin{proposition}\label{prop:fmp-algebraic}
  If an intuitionistic formula $\phi$ is consistent, then there exists a valuation $v$ from the variables of $\phi$ into a finite Heyting algebra such that $\overline{v}(\phi) \neq \bot$.
\end{proposition}
\begin{proof}
  Denote by $V$ the set of propositional variables occurring in $\phi$. Let $v_0 \colon V \to H$ be a valuation to a Heyting algebra such that $\overline{v_0}(\phi) \neq \bot$. Let $F$ be the bounded sublattice of $H$ generated by the set of elements $\overline{v_0}(\psi)$, where $\psi$ ranges over the subformulas of $\phi$, that is, the formulas that occur in the construction tree of $\phi$. Note that $F$ is finite since it is a distributive lattice generated by a finite set; thus, $F$ is a Heyting algebra. Let $v \colon V \to F$ be the co-restriction of $v_0$ to $F$. We may now show by induction that for any subformula $\psi$ of $\phi$, $\overline{v}(\psi) = \overline{v_0}(\psi)$. The only non-trivial case is when $\psi = \psi_1 \to \psi_2$ for some formulas $\psi_1$ and $\psi_2$. By the induction hypothesis, $\overline{v}(\psi_i) = \overline{v_0}(\psi_i)$ for $i = 1, 2$, and by definition, $\overline{v}(\psi) = \overline{v}(\psi_1) \hto_F \overline{v}(\psi_2)$ and $\overline{v_0}(\psi) = \overline{v_0}(\psi_1) \hto_H \overline{v_0}(\psi_2)$. The equality $\overline{v}(\psi) = \overline{v_0}(\psi)$ now follows from the general fact that if $a, b \in F$ are such that $a \to_H b \in F$, then $a \hto_F b = a \hto_H b$ (see Exercise~\ref{exe:heytingsub2}). In particular, $\overline{v}(\phi) = \overline{v_0}(\phi) \neq \bot$, as required.
\end{proof}
\begin{corollary}
  Every consistent intuitionistic formula is satisfiable in a finite model.
\end{corollary}
\begin{proof}
  In the proof of Theorem~\ref{thm:int-completeness}, we may take the Heyting algebra $H$ to be finite, by Proposition~\ref{prop:fmp-algebraic}. The poset dual to the finite Heyting algebra, equipped with the admissible valuation given in Theorem~\ref{thm:int-completeness}, gives a finite model in which the formula is satisfiable.
\end{proof}
Esakia duality may be used to deduce many more results on Heyting algebras and Kripke semantics. We refer the interested reader to \cite{GehEsak} for more information, as well as more historical information on Esakia duality.

\ourexercises

\begin{ourexercise}\phantomsection\label{exer:pseudocomplement}
  \begin{enumerate}
    \item Show that in any Boolean algebra, $\neg$ is a pseudocomplement.     
    \item Show that any Heyting algebra is pseudocomplemented with $a^*$ 
    defined by $a\hto\bot$.
    \item Let $L$ be a distributive lattice and $b\in L$. Let $L_b$ be the distributive 
    lattice whose underlying set is ${\uparrow}b$ in $L$ with the induced order from 
    $L$. Show that $L_b$ has the same $\top$ as $L$, that $\wedge$ and $\vee$ are 
    the restrictions of the operations on $L$, but $\bot=b$. 
    \item Show that if $L$ is a Heyting algebra, then, for each $b\in L$, the distributive 
    lattice $L_b$ is pseudocomplemented with the pseudocomplement given by 
    $a\hto b$ for $a\in L_b$.
    \item Conversely, let $L$ be a distributive lattice for which $L_b$ is 
    pseudocomplemented for each $b\in L$. Show that $L$ is a Heyting algebra in 
    which $a\hto b$ is given by the pseudocomplement in $L_b$ of $a\vee b$.
    \item Show that if $L$ is a Heyting algebra, then $(a\hto\bot)\vee b$ is always 
    an operation of implication-type.
    \item Show that any Boolean algebra is in particular a Heyting algebra, in which 
    $a\hto b$ is given by $\neg a\vee b$.
    \item Give an example of a Heyting algebra in which the operation 
    $(a\hto\bot)\vee b$ is not the relative pseudocomplement.
  \end{enumerate}
\end{ourexercise}
    
\begin{ourexercise}\label{exe:HAinterior}
  Let $X$ be a topological space, and $\Omega(X)$ its lattice of open sets. Prove that, for any open sets $U, V \subseteq X$, the relative pseudocomplement, $U \to V$, of $U$ with respect to $V$ is the interior of the set $(X \setminus U) \cup V$.
\end{ourexercise}

\begin{ourexercise}\label{exe:downsetsHA}
  The aim of this exercise is to prove equation (\ref{eq:downsetshto}). Let $(P, \leq)$ be a preorder and let $\tau$ be the dual Alexandrov topology for $\leq$ on $P$ (that is, $\tau$ is the topology of down-sets of $\leq$).
  \begin{enumerate}
    \item Prove that, for any $S \subseteq P$, a point $p \in P$ lies in the
	    interior of $S$ in the topology $\tau$ if, and only if, ${\downarrow} p \subseteq S$.
    \item From (a) and Exercise~\ref{exe:HAinterior}, deduce (\ref{eq:downsetshto}).
  \end{enumerate}
\end{ourexercise}

\begin{ourexercise}\label{exe:non-HA}
  Show that a complete distributive lattice may fail to be a Heyting algebra. \hint{Explain why this is essentially the same question as Exercise~\ref{exe:counterexamples-complete}.\ref{ite:compDLnotframe}.}
\end{ourexercise}

\begin{ourexercise}\label{exe:morphismsHAdifferent}
  Give examples of:
  \begin{enumerate}
    \item a frame (or even complete lattice) homomorphism between complete Heyting algebras which is not a Heyting homomorphism. \hint{It suffices to find a homomorphism between finite distributive lattices which does not preserve $\hto$. You may use finite duality to construct such an example.}

    \item a Heyting homomorphism between frames which is not a frame homomorphism. \hint{Consider $\mathbb N^\infty$, the chain of natural numbers with a top added, and the map to $2$ which sends all natural numbers to $0$ and $\infty$ to $1$.}
  \end{enumerate}
\end{ourexercise}

\begin{ourexercise}\label{exe:htofirstcoordinate}
  The aim of this exercise is to prove the third algebraic property of Heyting algebras, namely, that $\hto$ reverses arbitrary joins in the first coordinate into arbitrary meets. Let $H$ be a Heyting algebra.
  \begin{enumerate}
    \item Prove that, for any $a, b, c \in H$, $c \leq a \hto b \iff a \leq c \hto b$.
    \item Use the previous item to show that, for any $b \in H$, the function $(-) \hto b \colon H \to H^\op$ has an upper adjoint.
    \item Deduce that, for any $b \in H$ and $T \subseteq H$ such that $\bigvee T$ exists, we have $\big(\bigvee T\big) \hto b = \bigwedge_{t \in T} (t \hto b)$.
  \end{enumerate}
\end{ourexercise}

\begin{ourexercise}\label{exe:coheyting}
  A lattice $L$ is a \emphind{co-Heyting algebra} if, for every $a \in L$ the function $a \vee (-)$ has a lower adjoint. We will denote the value of this lower adjoint at $b \in L$ by $a \leftsquigarrow b$ here.
  \begin{enumerate}
    \item Prove that, for any topological space $X$, the closed set lattice is a co-Heyting algebra, with $C \leftsquigarrow D = \overline{C \setminus D}$ for any closed sets $C, D \subseteq X$.
    \item Prove that, for any preorder $(P, \leq)$, the down-set lattice $\Down(P)$ is a co-Heyting algebra.
  \end{enumerate}
\end{ourexercise}

\begin{ourexercise}\label{exe:Esakia-equiv} Show the equivalence of (ii) and (iii) in
  Proposition~\ref{prop:esakiaspace}. \hint{For (iii) $\Rightarrow$ (ii),
  recall Proposition~\ref{prop:cos-closed}. For (ii) $\Rightarrow$ (iii), use the
  fact that ${\uparrow}$ preserves arbitrary unions. } 
  \end{ourexercise}

\begin{ourexercise}\label{exe:heytingsub}
  Let $H$ be a Heyting algebra and let $L$ be a bounded sublattice of $H$ with the property that, for any $a, b \in L$, the implication $a \to_H b$ is in $L$. Prove that $L$ is a Heyting algebra, with implication given by the restriction of $\to_H$ to $L$.
\end{ourexercise}

\begin{ourexercise}\label{exe:heytingsub2}
  Let $L$ be a bounded sublattice of a Heyting algebra $H$. Suppose that $a, b \in L$ and that $a \to_H b$ is also in $L$. Show that $a \to_H b$ is a relative pseudocomplement of $a$ with respect to $b$ in $L$. Conclude in particular that, if $L$ is finite, and $a, b \in L$ are such that $a \to_H b$ is in $L$, then $a \to_L b = a \to_H b$.
\end{ourexercise}

\begin{ourexercise}\label{exe:heyting-S4}
  This exercise outlines the algebraic content of the so-called \emphind{Gödel translation} of intuitionistic logic into the modal logic known as {\bf S4}, and outlines why the bounded morphisms of this section are a special case of those between Kripke frames, in Section~\ref{sec:kripke}. This exercise is based on the results in \cite[Section 4]{GehEsak}, and solutions and further references may be found there.

  Let $L$ be a distributive lattice, and let $e \colon L \to B$ be the \emphind{Boolean envelope} of $L$, as defined in Definition~\ref{def:booleanenvelope}.
  \begin{enumerate}
    \item Prove that $L$ is a Heyting algebra if, and only if, the map $e$ has an upper adjoint, $g$.
    \item Suppose that $L$ is a Heyting algebra. Prove that the composite function $\Box := e \circ g \colon B \to B$ preserves finite meets, and satisfies the {\bf S4} axioms $\Box a \leq \Box \Box a$ and $\Box a \leq a$ for every $a \in B$.
    \item Conversely, if $\Box$ is a finite-meet-preserving function on a Boolean algebra $A$ satisfying the {\bf S4} axioms, prove that the image of $\Box$ is a (bounded) sublattice $L$ of $A$, which admits a relative pseudocomplement given, for $a,b\in L$, by $a\hto b=\Box(\neg a\vee b)$, where the operations on the right-hand side are those of $A$. Further, show that the inclusion $L\hookrightarrow A$ has an upper adjoint.
    \item Extend the object correspondence to a categorical equivalence: the category of Heyting algebras is equivalent to the subcategory of modal algebras consisting of those algebras $(A, \Box)$ that satisfy the {\bf S4} axioms and such that the image of $\Box$ generates $A$.
    \item Use this result to explain why bounded morphisms as defined in this section are a special case of the bounded morphisms of Section~\ref{sec:kripke}.
  \end{enumerate}
\end{ourexercise}

\section{Boolean closure and alternating chains}\label{sec:altchains}
In this section, we give an application of \emph{discrete} duality that we will use in Chapter~\ref{ch:AutThry}. The result we present here has a combinatorial flavor, although it also forms the kernel of a more general topological principle, explored in much more detail in \cite{BGKS2020}. In the exercises of this section, we ask the reader to work out a few elements of this more general theory.

Let $L$ be a sublattice of a Boolean algebra $B$. In accordance with Corollary~\ref{cor:boolenv-hull}, we denote by $L^-$ the Boolean subalgebra of $B$ generated by $L$. As a Boolean algebra, $L^-$ is the Boolean envelope of $L$, and we have an inclusion $L^- \into B$ of Boolean algebras. The general aim of the work that we do a piece of here, is to give a concrete dual characterization of the elements $a \in B$ that belong to $L^-$.

A basic idea, going back to \cite{Hausdorff14}, is that $L^-$ admits a stratification, according to the complexity of describing an element of $L^-$ in terms of elements of $L$. More concretely, if $a \in L^-$ then $a$ can be written as $a = a_1 - (a_2 - (\dotsm (a_{n-1}-a_n) \dotsm))$ for some elements $a_n \leq a_{n-1} \leq \cdots \leq a_1$ of $L$; such a sequence is called a \emphind{difference chain} for $a$. This leads to a well-defined notion of degree of $a$ over $L$ as the minimum length of such chains. However, in general, there may not always be a \emph{least} difference chain for  $a \in L^-$, where we compare difference chains for $a$ by saying that one difference chain $a_n \leq a_{n-1} \leq \cdots \leq a_1$ is less than another difference chain $b_m \leq b_{m-1} \leq \cdots \leq b_1$ if $n\leq m$ and $a_i\leq b_i$ for $i\leq n$. In this section, we concentrate on a simple situation where there always is a least difference chain for elements of $L^-$ over $L$.

In what follows, let $(X,\leq)$ be a poset and write $B = \cP(X)$ and $L = \cU(X)$ for the power set Boolean algebra and the sublattice of up-sets of $X$, respectively. Our aim is to describe $L^-$, the collection of subsets of $X$ that can be written as a Boolean combination of up-sets of $X$. The crucial notion that will allow us to do so is the following.

\begin{definition}\label{dfn:alternating-chain}
Let $(X, \leq)$ be a poset and let $a \in \cP(X)$. An \emphind{alternating chain} for $a$ is a chain $x_1 \leq y_1 \leq x_2 \leq y_2 \cdots \leq x_{n} \leq y_n$ of elements of $X$ such that, for every $1 \leq i \leq n$, $x_i \in a$ and $y_i \not\in a$. We call the number $n$ the \emph{height} of the chain; an alternating chain of height $n$ thus contains $2n$ elements. The \emphind{alternation height} of $a$, $h(a) \in \bN \cup \{\infty\}$, is defined as 
\[ h(a) := \sup \{n \in \bN \mid \text{there exists an alternating chain of height } n \text{ for } a\}.\]
We say an element $a \in \cP(X)$ has \emph{bounded alternation height} if $h(a) < \infty$.
\end{definition}
In particular, the alternation height of a set $a \in \cP(X)$ is by definition $0$ when no alternating chains exist for it; note that this happens if, and only if, $a$ is an up-set. 

The aim of this section is to prove the following theorem.
\begin{theorem}\label{thm:discrete-alt-height}
Let $(X,\leq)$ be a poset. The Boolean subalgebra of $\cP(X)$ generated by $\cU(X)$ consists of the elements of bounded alternation height.
\end{theorem}
Denote by $A$ the set of elements of $\cP(X)$ of bounded alternation height, and $L = \cU(X)$; our aim is to prove that $L^- = A$. The left-to-right inclusion is relatively straightforward, as we show in the following proposition.
\begin{proposition}\label{prop:bddaltheight-is-ba}
Let $(X, \leq)$ be a poset. The set $A$ of elements of $\cP(X)$ that have bounded alternation height is a Boolean subalgebra of $\cP(X)$ that contains all up-sets of $(X,\leq)$. In particular, $\cU(X)^- \subseteq A$.
\end{proposition}
\begin{proof}
As noted above, any up-set has alternation height $0$.
Now suppose that $a \in A$ has alternation height $h < \infty$. Note that if $w$ is an alternating chain for the complement $a^c$ of $a$ of length $n+1$, then removing its first and last element gives an alternating chain for $a$ of length $n$. Thus, the alternation height of $a^c$ is at most $h+1$.

Let $a_1, a_2 \in A$. We claim that the alternation height of $a_1 \cup a_2$  is less than $h := 2\max(h(a_1),h(a_2))+1$. Indeed, suppose towards a contradiction that $x_1 \leq y_1 \leq \cdots \leq x_h \leq y_h$ is an alternating chain for $a_1 \cup a_2$ of height $h$. Then we have $x_n \in a_1 \cup a_2$ for every $1 \leq n \leq h$, so by the pigeon-hole principle, there is $i \in \{1,2\}$ such that for $k \geq \frac{h}{2}$ indices, we have $x_n \in a_i$. Consider the subchain $x_{n_1} \leq y_{n_1} \leq \cdots \leq x_{n_k} \leq y_{n_k}$ where $n_1 < \dotsm < n_k$ are the indices such that $x_{n} \in a_i$. This is clearly an alternating chain for $a_i$, because, for any $1 \leq n \leq h$, we have $y_n \not\in a_1 \cup a_2$, so in particular $y_n \not\in a_i$. However, the height of this alternating chain is $k \geq \frac{h}{2} > h(a_i)$, contradicting the definition of $h(a_i)$.

The `in particular' statement now follows immediately from the fact that $\cU(X)^-$ is the Boolean subalgebra of $\cP(X)$ generated by $\cU(X)$ (Corollary~\ref{cor:boolenv-hull}).
\end{proof}
We now show a somewhat more technically involved lemma, which in particular establishes the inclusion $A \subseteq L^-$, and thus the equality of Boolean algebras $L^- = A$ claimed in Theorem~\ref{thm:discrete-alt-height}.
\begin{lemma}\label{lem:altchain}
  Suppose that $a \in \cP(X)$ has alternation height $h \in \mathbb{N}$. Then there exist $a_0, \dots, a_{h}$ and $b_0, \dots, b_{h}$ in $\cU(X)$ such that 
  \[ a = \bigcup_{n=0}^{h} (a_n \setminus b_n).\]
\end{lemma}
\begin{proof}
  Throughout this proof, we slightly extend the definition of alternating chain, and we will call a sequence $x_1 \leq y_1 \leq \cdots \leq x_n \leq y_n \leq x_{n+1}$ of elements of $X$ with $x_i \in a$ for $1 \leq i \leq n+1$ and $y_i \not\in a$ for $1 \leq i \leq n$ an \emph{alternating chain for} $a$ of height $n + \frac{1}{2}$.

  We inductively define the following sequence of elements of $\cU(X)$:
  \begin{equation*}
  \begin{aligned}[c]
  &a_0 := {\uparrow} a, &b_0 := {\uparrow}(a_0 \cap a^c), &\\
  &a_{n+1} := {\uparrow} (b_n \cap a), &b_{n+1} := {\uparrow}(a_{n+1} \cap a^c) &\text{ for } 0 \leq n \leq h.
  \end{aligned}
  \end{equation*}
  We now first show by induction on $n$ that, for any $n \geq 0$ and any $x \in X$,
  \begin{itemize}
  \item $x \in a_n$ if, and only if, there exists an alternating chain for $a$ of height $n+\frac{1}{2}$ whose last point is $\leq x$, and
  \item $x \in b_n$ if, and only if, there exists an alternating chain for $a$ of height $n + 1$ whose last point is $\leq x$.
\end{itemize}
  For $n = 0$, note that an alternating chain of height $\frac{1}{2}$ is just a single element of $a$, which gives the first item. For the second item, unraveling the definition of $b_0$, we see that $x \in b_0$ if, and only if, there exist $y_0 \not\in a$ and $x_0 \in a$ with $x_0 \leq y_0 \leq x$, as required. Now, let $n \geq 1$ be arbitrary. For any $x \in X$, by definition we have $x \in a_{n}$ if, and only if, there exists $x' \in b_{n-1} \cap a$ with $x' \leq x$. Suppose that $x \in a_n$. By induction, since $x' \in b_{n-1}$, there is then an alternating chain for $a$ of height $n$ whose last point is $\leq x'$. Since $x' \in a$, extending this sequence by adding $x'$ to the end gives an alternating chain for $a$ of height $n + \frac{1}{2}$ and $x' \leq x$. Conversely, if $x_1 \leq y_1 \leq \cdots \leq x_n \leq y_n \leq x_{n+1}$ is an alternating chain for $a$ of height $n + \frac{1}{2}$ with $x_{n+1} \leq x$, then $x_{n+1} \in a$, and the alternating chain $x_1 \leq y_1 \leq \cdots \leq x_n \leq y_n$ shows, by induction, that $x_{n+1} \in b_{n-1}$, so that $x \in {\uparrow}(b_{n-1} \cap a) = a_{n}$, as required. The proof of the second item is very similar.
  
  Let us now write $c_0 := \emptyset$ and, for each $1 \leq n \leq h+1$,
  \[ c_n := \bigcup_{i=0}^{n-1} (a_i \setminus b_i).\] 
  The statement of the lemma is that $a = c_{h+1}$. Towards proving this, we will show that, for each $0 \leq n \leq h+1$, 
  \begin{equation}\label{eq:alt-height-proof-key}
  c_n \subseteq a \subseteq c_n \cup a_n.
  \end{equation}
  Note that the fact that the alternation height of $a$ is $h$ implies, using the characterization above, that $b_h = a_{h+1} = \emptyset$. Thus, (\ref{eq:alt-height-proof-key}) for $k = h+1$ will give the desired equality $a = c_{h+1}$.
  
  For the first inclusion of (\ref{eq:alt-height-proof-key}), it suffices to show that $a_n \setminus b_n \subseteq a$ for every $0 \leq n \leq h$. Indeed, if $x \in a_n \setminus b_n$, pick an alternating chain of height $n+\frac{1}{2}$ for $a$ whose last point is $\leq x$. Since $x \not\in b_n$, adding $x$ to the end of this alternating chain does \emph{not} give an alternating chain of height $n+1$, so we must have $x \in a$.

  For the second inclusion of (\ref{eq:alt-height-proof-key}), we proceed by induction on $n$. When $n = 0$, simply note that indeed $a \subseteq a_0$. For the inductive step, suppose that $x \in a$. We distinguish two cases, according to whether or not $x \in b_{n-1}$. 
  
  Case 1: $x \in b_{n-1}$. Then there exists an alternating chain for $a$ of height $n$ whose last point is $\leq x$. Since $x \in a$, we can extend this sequence to an alternating chain of height $n+\frac{1}{2}$ whose last point is $x$, so $x \in a_n$.

  Case 2: $x \not\in b_{n-1}$. Then, if $x \in a_{n-1}$, we are done because $a_{n-1} \setminus b_{n-1} \subseteq c_n$ by definition. Assume $x \not\in a_{n-1}$. By the inductive hypothesis, we have $a \subseteq c_{n-1} \cup a_{n-1}$, so $x \in c_{n-1}$. Clearly, $c_{n-1} \subseteq c_n$, so we conclude $x \in c_n$.
\end{proof}

\ourexercises

\begin{ourexercise}\label{exer:difchain-least}
This exercise is about the collection of elements $\{a_i,b_i\}_{i\in \N}$ that 
we construct for an element $a\in\cP(X)$ in Lemma~\ref{lem:altchain}.
\begin{enumerate}
\item  Show that $a_0\supseteq b_0\supseteq a_1\supseteq\dots\supseteq b_n$ 
for all $n\in\N$.
\item Show that 
 \[ \bigcup_{i=0}^{n} (a_i \setminus b_i)= 
 a_0\setminus(b_0\setminus(a_1\setminus\dots(a_n\setminus b_n)\dots)).
 \]
 \item For an element $a\in\cP(X)$, we call $a_0\supseteq b_0\supseteq 
 a_1\supseteq\dots\supseteq b_n$ with $a_i,b_i\in\cU(X)$ and 
  \[ a = \bigcup_{n=0}^{h} (a_n \setminus b_n)\]
  a \emph{difference chain} for $a$ and we order difference chains 
  for $a$ by coordinate-wise inclusion (if one is shorter than the other, then
  we consider it extended with empty sets). Show that if $a$ has a 
  difference chain then it has a least such and it is the one we define 
  in the proof of Lemma~\ref{lem:altchain}.
\end{enumerate}
\end{ourexercise}

\begin{ourexercise}\label{exer:co-Heyting}
Suppose $H$ is a co-Heyting algebra, and $X$ its dual space. Consider $X$ 
with the opposite of the Priestley order so that $H$ is isomorphic to the lattice 
of clopen up-sets via $a\mapsto\widehat{a}$. Let $V\subseteq X$ be a clopen in 
$X$. Show that the difference chain for $V$ as defined in 
Exercise~\ref{exer:difchain-least} consists entirely of clopen up-sets of $X$. 
Conclude that each element of $H^-$ has a least difference chain over $H$.
\end{ourexercise}

\begin{ourexercise}\label{exer:DLalternationheight}
Show that for an arbitrary bounded distributive lattice $L$, each element 
$a\in L^-$ is given by a difference chain $a_0>b_0>\dots a_n>b_n$ of 
elements of $L$ such that
\[
a =a_0\setminus(b_0\setminus(a_1\setminus\dots(a_n\setminus b_n)\dots))
= \bigcup_{n=0}^{h} (a_n \setminus b_n),
\]
and that this gives rise to a well-defined notion of alternation height for 
elements $a\in L^-$ over $L$. \hint{If $a\in L^-$, then $a\in F^-$, where 
$F$ is some finite sublattice of $L$. The finite sublattice $F$ is then in particular a co-Heyting 
algebra, and provides a least difference chain for $a$ over $F$, which is 
also a difference chain for $a$ over $L$.}
\end{ourexercise}

\begin{ourexercise}\label{exer:no-least-diffchain}
Give an example of a distributive lattice $L$, and an element $a\in L^-$ to 
show that there is not a least difference chain for $a$ in general.
\end{ourexercise}

\begin{ourexercise}\label{exer:compsat-diffchain}
Show that the closed down-sets of the dual Priestley space of a distributive 
lattice $L$ form a co-Heyting algebra $H$ containing $L$, and that the 
alternation height $a\in L^-$ is the same over either $L$ or $H$. 
\end{ourexercise}

\notessec
\label{fn:canext}
The theory of canonical extensions
 allows one to make more precise the idea, mentioned above 
Proposition~\ref{prop:relation-dual-to-box}, that the duality for operators treated 
in this chapter is an extension of the finite case; see also 
Exercise~\ref{exe:finiterelationduality}. 
In that theory, any distributive lattice $L$ 
embeds into a ``finite-like'' lattice, $L^\delta$, and any finite-meet-preserving 
function $h \colon M \to L$ is shown to lift in a unique way to a completely meet 
preserving function $h^\delta \colon M^\delta \to L^\delta$, which then yields the 
dual relation $R$ from the dual space of $X$ to the dual space of $Y$, exactly 
by the ideas from the finite duality outlined here. This algebraic definition of the 
dual relation was at the core of the work in \cite{JonTar1951, JonTar1952}, also 
see for example \cite{GehJon1994, GehJon2004}, which generalize this theory 
from Boolean algebras to distributive lattices and in addition study generalizations 
of the theory to arbitrary operations that do not need to preserve joins or meets.

A property of a lattice of operators is called \emphind{canonical} if it is preserved under canonical extension. The study of canonicity and its connections to first order definability is a major research theme in modal logic, see for example \cite{GHV03, HV05}. The field of \emph{Sahlqvist theory} provides syntactic  conditions on properties of lattices with operators that guarantee  canonicity, see for example \cite[Chapter~5]{BRV2001}, \cite{GNV2005,ConGhiPal2014, ConradiPalmigiano2020}.

The structures that we call `Kripke Boolean spaces' in Section~\ref{sec:kripke} are typically presented in the modal logic literature in the form of
\emph{descriptive general frames}, which give an isomorphic category. The
latter however avoids having to  speak explicitly about topology, see for example
\cite[Ch.~5]{BRV2001} or \cite[Ch.~8]{ChaZak1997}. We prefer to make the
Boolean topology on general frames explicit, because we think it clarifies the
mathematical content of the general semantics, and in particular the link with
Stone duality.
\chapter{Categorical duality}\label{ch:categories}
In this chapter, we introduce notions from category theory and use them to rephrase the duality theorems that we have seen so far. For our purposes, category theory provides a general language that allows us to talk about duality from a higher level viewpoint. We limit ourselves to explaining only those notions from category theory that we need in this book; this is by no means a complete introduction to category theory. In particular, we do not go into deep connections between the topological duality theory developed in this book and the rich category theory of frames and toposes; for more on this, see for example \cite{Johnstone1986, Mak1993, MM1992, Car2017}. Since we want to keep category-theoretic preliminaries to a minimum,  all the category theory treated in this chapter is entirely standard and well-known, and can be found in a lot more detail in, for example, \cite{Mac1971, Awo2010, Lei2014}.

\section{Definitions and examples of categories}\label{sec:categorydef}
A \emph{category} formalizes a class of mathematical objects, the morphisms between them, and how the morphisms compose. For example, the class of sets and functions, with usual function composition, is a category, but so is the class of sets and binary relations, with relational composition.

A category is a two-sorted structure: it has an \emph{objects} sort and a \emph{morphisms} sort. Every morphism $f$ in a category has an associated \emph{domain} object, $\dom(f)$, and \emph{codomain} object, $\cod(f)$. For every `composable' pair of morphisms $(f,g)$, that is, a pair such that $\dom(g) = \cod(f)$, there is a morphism, $g \circ f$, the \emph{composition} of $g$ after $f$, and this operation of composition is required to be associative. Finally, every object $A$ in a category has an associated \emph{identity morphism} $1_A$, which acts as a neutral element for the composition. We now give the formal definition. Although this definition may look lengthy and abstract at first, the conditions are natural and easily seen to be satisfied in the examples that follow.

\begin{definition}\label{def:category}\index{category}
A \emphind{category} $\cat{C}$ is a tuple $(\ob \cat{C}, \mor \cat{C}, \dom, \cod, \circ, 1_{(-)})$, where $\ob \cat{C}$ and $\mor \cat{C}$ are classes of \emph{objects} and \emph{morphisms}, $\dom$ and $\cod$ are assignments from $\mor \cat{C}$ to $\ob \cat{C}$, $1_{(-)}$ is an assignment from $\ob \cat{C}$ to $\mor \cat{C}$, and $\circ$ is an assignment from $\{(g,f) \ | \  \dom(g) = \cod(f)\}$ to $\mor \cat{C}$, with the following properties:
\begin{enumerate}
	\item for any $f,g \in \mor \cat{C}$ with $\dom(g) = \cod(f)$,  
	\[\dom(g \circ f) = \dom(f) \text{ and } \cod(g \circ f) = \cod(g).\]
	\item for any $f, g, h \in \mor \cat{C}$, if $\dom(g) = \cod(f)$ and $\dom(h) = \cod(g)$, then 
	\[h \circ (g \circ f) = (h \circ g) \circ f.\]
	\item for any $A \in \ob \cat{C}$, 
	\[\dom(1_A) = A = \cod(1_A).\]
	\item for any $f \in \mor \cat{C}$ with $\dom(f) = A$ and $\cod(f) = B$,
	\[ f \circ 1_A = f \text{ and } 1_B \circ f = f.\]
\end{enumerate}
\end{definition}
\nl{$\ob \cat{C}$}{objects of a category $\cat{C}$}{}
\nl{$\mor \cat{C}$}{morphisms of a category $\cat{C}$}{}

While Definition~\ref{def:category} gives the official definition of the structure usually called a category, alternative definitions are possible. In particular, the definition of category does not need to be two-sorted, as the structure of a category is entirely determined by its class of morphisms: the objects can be encoded in the class of morphisms using the identity morphisms; see, for example, \cite[p.~9]{Mac1971} or \cite[Definition~3.53, p.~42]{AdaHerStr1990}.

The structure of categories is rich, due to many derived notions, which we begin to introduce now. Let $\cat{C}$ be a category. For any $A, B \in \ob \cat{C}$, the \emph{Hom-class}\index{Hom-class}\index{Hom-set} from $A$ to $B$ in $\cat{C}$ is the collection of morphisms $f$ with $\dom(f) = A$ and $\cod(f) = B$; this class is denoted by $\Hom_{\cat{C}}(A,B)$, or simply $\cat{C}(A,B)$. We also use the notation $f \colon A \to B$ to mean $f \in \Hom_{\cat{C}}(A,B)$ when the ambient category is clear. The notation $\circ$ and parentheses are often omitted: for example, we write $hgf$ for the morphism $h \circ (g \circ f) = (h \circ g) \circ f$.
\nl{$\Hom_{\cat{C}}(A,B)$}{morphisms from $A$ to $B$ in a category $\cat{C}$}

Note that the components defining a category are not assumed to be sets. This level of generality is necessary, because, for example, the collection of all sets is not a set itself, and we want to be able to consider the category of sets. Unless explicitly mentioned otherwise, all the categories we work with in this book are \emph{locally small}\index{category!locally small}, which means that, for any fixed pair of objects $A$, $B$ in the category, the Hom-class from $A$ to $B$ is a set. A category is called \emph{small}\index{category!small} if the class of morphisms (and, hence, also the class of objects) is a set, and \emph{large}\index{category!large} otherwise.

\begin{example}\label{exa:bigcats}
The following are examples of categories that we have already seen. In each example, we only state what $\ob \cat{C}$ and $\mor \cat{C}$ are; composition and identity are what you expect.
\begin{enumerate}
\item The category $\cat{DL}$ of distributive lattices and their homomorphisms.
\item The category $\TopCat$ of topological spaces and continuous functions.
\item The category $\cat{Priestley}$ of Priestley spaces and continuous order-preserving functions.
\item The category $\cat{BA}$ of Boolean algebras and their homomorphisms.
\item The category $\cat{BoolSp}$ of Boolean spaces and continuous functions between them.
\item The category $\cat{HA}$ of Heyting algebras and their homomorphisms.
\item The category $\cat{Esakia}$ of Esakia spaces and continuous p-morphisms.
\item The categories $\cat{DL}_f$ and $\cat{BA}_f$ of \emph{finite} distributive lattices and \emph{finite} Boolean algebras, respectively, and homomorphisms between them.
\item The category $\cat{Set}$ of sets and functions; the category $\cat{Set}_f$ of finite sets and functions between them.
\item The category $\cat{Pos}$ of posets and order-preserving functions; the category $\cat{Pos}_f$ of finite posets and order-preserving functions between them.
\end{enumerate}
\end{example}
The categories described in the above example correspond to the way we motivated categories above, as `mathematical universes'. In technical terms, these are large categories that are \emph{concrete}\index{category!concrete}, meaning intuitively that each object is a set (with additional structure) and each morphism is a (special kind of) function between the sets underlying the objects. We will be a bit more precise about this concept of a `concrete category' in  Example~\ref{exa:cat-functor-examples}.\ref{itm:forgetfulfunctor} below. 
However, there are other, more abstract examples of categories. 
\begin{example}\label{exa:smallcats} We give two classes of typical examples of categories, \emph{thin} categories and \emph{monoids}. %
\begin{enumerate}\index{category!preorder as a}
\item \label{itm:thin}Let $(P,\leq)$ be a preorder. There is a category $\cat{P}$ with $P$ as the set of objects and ${\leq}$ as the set of morphisms; that is, a morphism is a pair $(p,q) \in P \times P$ with $p \leq q$, $\dom((p,q)) := p$ and $\cod((p,q)) := q$. Exercise~\ref{exe:preordercategory} asks you to define the additional structure that makes $\cat{P}$ into a category. Categories $\cat{P}$ that arise in this way are called \emphind{thin}\index{category!thin}; they are characterized by the property that there is at most one morphism between any pair of objects.
\item\label{ite:monoid} \index{category!monoid as a} Let $(M,\cdot,1)$ be a monoid, that is, a set with an associative operation $\cdot$ such that $1 \cdot m = m = m \cdot 1$ for all $m \in M$. There is a category $\cat{M}$ with a single object, $\ast$, and $M$ as the set of morphisms. The domain and codomain of any morphism is by definition $\ast$, the identity morphism $1_\ast$ on the object $\ast$ is $1$, and composition is defined by $\cdot$. The axioms of a monoid immediately give that $\cat{M}$ is indeed a category.
\end{enumerate}
\end{example}
Both of these examples of categories are in a sense degenerate, but in two different directions: few morphisms (preorders), and one object (monoids). A general category can thus be thought of as either `a multi-arrow preorder' or `a multi-object monoid'.

\medskip

One more central notion in categories, that we introduce now because we will need it throughout this chapter, is that of an \emph{isomorphism}.

\begin{definition}\label{dfn:internaliso}\index{isomorphism}
Let $f \colon A \to B$ be a morphism in a category $\cat{C}$. A \emph{two-sided inverse} to $f$ is a morphism $g \colon B \to A$ in $\cat{C}$ such that $g \circ f = 1_A$ and $f \circ g = 1_B$. The morphism $f$ is called an \emph{isomorphism} if a two-sided inverse exists for it. Two objects $A$ and $B$ are called \emph{isomorphic} if there exists an isomorphism from $A$ to $B$.
\end{definition}
If $f$ is an isomorphism, then the inverse of $f$ is unique, and is denoted by $f^{-1}$. This morphism $f^{-1}$ is also an isomorphism, and the notion of being isomorphic is in fact an equivalence relation on the objects of $\cat{C}$ (see Exercise~\ref{exe:inverseunique})\footnote{Note that, since $\ob \cat{C}$ is not necessarily a set, while `equivalence relations' are usually taken to be subsets, when we say `equivalence relation' in this chapter, more formally we mean: a binary predicate that satisfies reflexivity, symmetry, and transitivity; but we will not need to worry about size issues.}. A \emph{skeleton}\index{category!skeleton}\index{skeleton} of a category $\cat{C}$ is a full subcategory $\cat{S}$ of $\cat{C}$ such that each object is isomorphic to exactly one object of $\cat{S}$ (see Exercise~\ref{exe:skeleton} in Section~\ref{sec:internal}).

\begin{example}
In any of the concrete categories from Example~\ref{exa:bigcats}, `isomorphism' as defined in Definition~\ref{dfn:internaliso} means exactly what one would expect. In particular, an isomorphism between sets is just a bijection, and an isomorphism between Priestley spaces is an order-homeomorphism.

In a thin category $\cat{P}$ as in 
Example~\ref{exa:smallcats}.\ref{itm:thin}, two objects
are isomorphic if, and only if, they are equivalent in the sense of Exercise~\ref{exe:reflection}. In a single-object category $\cat{M}$ as in Example~\ref{exa:smallcats}.\ref{ite:monoid}, an isomorphism corresponds to an \emph{invertible} element of the monoid. In particular, the monoid $M$ is a \emph{group} if, and only if, every morphism of $\cat{M}$ is an isomorphism. More generally, a category $\cat{G}$ is called a \emph{groupoid} if every morphism in $\cat{G}$ is an isomorphism.
\end{example}

\ourexercises
\begin{ourexercise}\label{exe:preordercategory}\index{category!thin}
Let $(P,\leq)$ be a preordered set. Define $\ob \cat{P} := P$ and $\mor \cat{P} := {\leq}$. For any $(p,q) \in \mor \cat{P}$, let $\dom (p,q) := p$ and $\cod (p,q) := q$.
\begin{enumerate}
\item Why does there exist, for every $p \in P$, a unique morphism $1_p$ with domain and codomain $p$?
\item For any $(p,q), (q,r) \in \mor \cat{P}$, define $(q,r) \circ (p,q) := (p,r)$. Why is $(p,r)$ a morphism in $\cat{P}$?
\item Verify that the structure defined above turns $\cat{P}$ into a category in which there is at most one morphism between any two objects. A category of this kind is called \emph{thin}.
\item Conversely, if $\cat{C}$ is a small category in which there is at most one morphism between any two objects, define a preorder $\leq$ on the set $\ob \cat{C}$ such that $\cat{C}$ is the category associated to the preorder $(\ob \cat{C}, \leq)$.
\end{enumerate}
\end{ourexercise}

\begin{ourexercise}\label{exe:categorymonoid}
Let $\cat{C}$ be a small category with one object. Define a monoid $M$ such that $\cat{C}$ is the category associated to the monoid $M$, as in Example~\ref{exa:smallcats}.\ref{ite:monoid}.
\end{ourexercise}

\begin{ourexercise}\label{exe:catofrels}
For any sets $X, Y, Z$ and relations $R \subseteq X \times Y$ and $S \subseteq Y \times Z$, recall that the \emphind{relational composition} of $R$ and $S$ is the relation $R \cdot S \subseteq X \times Z$ defined by 
\[ R \cdot S := \{(x,z) \in X \times Z\mid\exists\, y \in Y \text{ such that } (x,y) \in R \text{ and } (y,z) \in S\}.\]
Verify that the following properly defines a category $\cat{Rel}$: let $\ob \cat{Rel}$ be the class of sets, and for any sets $X$, $Y$, let $\Hom_{\cat{Rel}}(X,Y)$ be the set of relations $R \subseteq X \times Y$. For any set $X$, let $1_X$ be the diagonal relation on $X$. Finally, for $R \in \Hom_{\cat{Rel}}(X,Y)$ and $S \in \Hom_{\cat{Rel}}(Y,Z)$, define $S \circ R := R \cdot S$. 
\end{ourexercise}

\begin{ourexercise}\label{exe:inverseunique}
Let $f \colon A \to B$ be a morphism in a category $\cat{C}$. 
\begin{enumerate}
\item Prove that if $g$ and $g'$ are both two-sided inverses to $f$, then $g = g'$.
\item Suppose that $f$ is an isomorphism. Prove that its two-sided inverse $f^{-1}$ is also an isomorphism.
\item Define a relation ${\sim}$ on $\ob \cat{C}$ by: $A \sim B$ if, and only if, $A$ is isomorphic to $B$. Prove that ${\sim}$ is an equivalence relation.
\end{enumerate}
\end{ourexercise}

\section{Constructions on categories}\label{sec:external}
Categories can be studied from two different perspectives, both of which are important; one is \emph{external}, the other \emph{internal}. %
From the \emph{external} perspective, one regards categories themselves as mathematical objects. From this point of view, one may consider mappings \emph{between} categories and constructions \emph{on} categories. For example, in this section we will describe the opposite category, subcategories of a category, functors between categories, adjunctions between categories and, most important for our purposes, dualities. %

From the \emph{internal} perspective, the structure of a category, viewed as a
universe for performing mathematical operations, allows us to mimick many natural constructions that we know from mathematics \emph{within} a given, fixed category. In the previous section, we already saw in Definition~\ref{dfn:internaliso} that the familiar notion of isomorphism between mathematical objects can be defined internally in any category. We will see in the next section, Section~\ref{sec:internal}, that within many of the categories relevant to us, there is an analogue of `injective' and `surjective' map, of product, disjoint unions, power sets, and more. 

Many interesting results in category theory come from the interplay between the two perspectives. As a case in point, duality theory uses an external construction (a duality) to obtain knowledge about internal constructions (for example, quotients and subobjects) of the categories at hand.

As an aside, we note that many of the constructions that we describe here as
external constructions \emph{on} categories can also be viewed as internal
constructions \emph{within} a category, namely, the category $\cat{Cat}$ of
categories (to be properly defined shortly). This phenomenon, which is both
beautiful and a bit jarring, will mostly be ignored in this chapter. 
We recall that our aim in this chapter is not to give the reader a complete
course in category theory, but to only introduce categories to the extent that a
working duality theorist needs them.

\subsection*{The opposite category}
\nl{$\cat{C}^\op$}{the opposite of a category}{Cop}
For a category $\cat{C}$, the \emphind{opposite category}, $\cat{C}^\op$, is the category with the same objects, morphisms, and identity as the category $\cat{C}$, but the domain and codomain assignments are interchanged, and the order of writing composition is reversed. More formally:
\begin{definition}
Let $\cat{C} = (\ob \cat{C}, \mor \cat{C}, \dom, \cod, \circ, 1_{-})$ be a
category. The \emphind{opposite category}, $\cat{C}^\op$, is defined as $(\ob
\cat{C}, \mor \cat{C}, \cod, \dom, \bullet, 1_{-})$, where, if $f \colon A \to
B$, $g \colon B \to C$ are composable morphisms in $\cat{C}$ with composite $g
\circ f$, then the composition of the corresponding morphisms $g \colon C \to B$
and $f \colon B \to A$ in $\cat{C}^\op$, $f \bullet g$, is defined to be the
morphism $C \to A$ in $\cat{C}^\op$ corresponding to $g \circ f$. Thus, for any
objects $A, B$, we have by definition $\Hom_{\cat{C}^\op}(A,B) = \Hom_{\cat{C}}(B,A)$.
\end{definition}

\subsection*{Subcategories}
If $\cat{C}$ and $\cat{D}$ are categories with $\ob \cat{D} \subseteq \ob \cat{C}$ and $\mor \cat{D} \subseteq \mor \cat{C}$ and the assignments $\dom$, $\cod$, $\circ$, and $1_{(-)}$ of $\cat{D}$ are restrictions of those of $\cat{C}$, then $\cat{D}$ is called a \emphind{subcategory} of $\cat{C}$.

An important special case of a subcategory is that of a \emphind{full subcategory}. A subcategory $\cat{D}$ of $\cat{C}$ is called a \emph{full subcategory} if $\Hom_{\cat{D}}(A,B) = \Hom_{\cat{C}}(A,B)$ for every $A, B \in \ob \cat{D}$. If $D$ is a collection of objects in a category $\cat{C}$, then we can always consider \emph{the full subcategory of $\cat{C}$ whose objects are in $D$}.

For example, the categories $\cat{Set}_f$, $\cat{DL}_f$ and $\cat{BA}_f$ are full subcategories of $\cat{Set}$, $\cat{DL}$ and $\cat{BA}$, respectively. Also, $\cat{BA}$ is a full subcategory of $\cat{DL}$, using the fact that every lattice homomorphism between Boolean algebras is a Boolean homomorphism (see Exercise~\ref{exe:distuncomp}.\ref{itm:lathomBA}). The category $\cat{HA}$ is a \emph{non}-full subcategory of $\cat{DL}$: not every lattice homomorphism between Heyting algebras is a Heyting homomorphism (see Exercise~\ref{exe:morphismsHAdifferent}). A subcategory may only restrict morphisms, while keeping the same class of objects: for example, the category $\cat{Rel}$ introduced in Exercise~\ref{exe:catofrels} has $\cat{Set}$ as a non-full subcategory: when moving from $\cat{Rel}$ to $\cat{Set}$, the objects stay the same, but the morphisms are those relations that are functions. 
The category $\cat{Heyt}_f$ of finite Heyting algebras is an example of a non-full subcategory of $\cat{DL}_f$: every finite distributive lattice is a finite Heyting algebra, so the objects are the same, but not every distributive lattice homomorphism preserves the Heyting implication (see Exercise~\ref{exe:morphismsHAdifferent}).

\subsection*{Functors and natural transformations}
The correct notion of `homomorphism' between categories is that of a \emph{functor}. We have already seen several examples of functors. For example, in Definition~\ref{def:Priestleydualspace}, we gave a functor from $\cat{DL}$ to $\cat{Priestley}^\op$: we associated to any object $L$ of the category $\cat{DL}$ an object $X_L$ of the category $\cat{Priestley}$, and, in Proposition~\ref{prop:priestleymorphisms}, to any morphism $h \colon L \to M$ of the category $\cat{DL}$ a morphism $f \colon X_M \to X_L$ of the category $\cat{Priestley}$, that is, a morphism $f \colon X_L \to X_M$ in $\cat{Priestley}^\op$.
\begin{definition}
Let $\cat{C}$ and $\cat{D}$ be categories. A \emphind{functor} $F$ from $\cat{C}$ to $\cat{D}$ is a pair of assignments, $\ob \cat{C} \to \ob \cat{D}$ and $\mor \cat{C} \to \mor \cat{D}$ such that the following properties hold: 
\begin{enumerate}
  \item for any $f \in \mor \cat{C}$, $\dom(F(f)) = F(\dom(f))$ and $\cod(Ff) = F(\cod(f))$, 
  \item for any $A \in \ob \cat{C}$, $F(1_A) = 1_{FA}$, and 
  \item for any composable pair $(f,g) \in \mor \cat{C}$, $F(g \circ f) = F(g) \circ F(f)$. 
\end{enumerate}
Note a slight abuse of notation in this definition: both the object and the morphism assignment are denoted by $F$, even though they are strictly speaking two separate components of the functor $F$. This notation rarely leads to confusion.

Let $F \colon \cat{C} \to \cat{D}$ be a functor. Then $F$ is called 
\begin{itemize}
  \item \emphind{full} if, for any two objects $A$, $B$ of $\cat{C}$, the assignment $f \mapsto F(f)$ is surjective as a map from $\Hom_{\cat{C}}(A,B)$ to $\Hom_{\cat{D}}(FA,FB)$,
  \item \emphind{faithful} if, for any two objects $A$, $B$ of $\cat{C}$, the assignment $f \mapsto F(f)$ is injective as a map from  $\Hom_{\cat{C}}(A,B)$ to $\Hom_{\cat{D}}(FA,FB)$,
  \item \emphind{essentially surjective} if for every object $B$ in $\cat{D}$, there exists an object $A \in \cat{C}$ such that $FA$ is isomorphic to $B$ in $\cat{D}$.
\end{itemize}
\end{definition}
Functors can be composed in the obvious way, and on any category $\cat{C}$ there is an obvious \emph{identity functor}, denoted $1_{\cat{C}}$. In this way, as mentioned above, the collection of locally small categories and functors between them is itself again a category, \cat{Cat}! We will not pursue this point any further here.

\begin{example}\label{exa:cov-functor-examples} We give some first examples of functors.
  \begin{enumerate}
    \item If $P$ and $Q$ are preorders, viewed as categories as in Example~\ref{exa:smallcats}, then a functor $f \colon P \to Q$ is essentially the same thing as an order-preserving function: given a function on the underlying sets, the only non-trivial requirement for it to extend to a functor between the associated categories $\cat{P}$ and $\cat{Q}$ is that, if $p \leq p'$, then there must be a morphism from $f(p)$ to $f(p')$ in $Q$, in other words, $f(p) \leq f(p')$.
      \item If $M$ and $N$ are monoids, viewed as categories as in Example~\ref{exa:smallcats}, then a functor $f \colon M \to N$ is essentially the same thing as a monoid homomorphism.
        \item \label{itm:forgetfulfunctor} There is a functor $U \colon \cat{DL} \to \cat{Set}$ that sends any distributive lattice $L$ to its underlying set, and any homomorphism $h \colon L \to M$ to itself. The functor $U$ is clearly faithful, but not full. This functor $U$ is called the \emphind{forgetful functor}\index{functor!forgetful}, and in fact such a functor exists in a much wider setting, including any category of `algebraic structures'. While the functor may look rather trivial, it is important for defining \emph{free objects} categorically, as we will see in Example~\ref{exa:adjunctions}.\ref{itm:freeDL} below. Categories that admit a well-behaved forgetful functor to $\cat{Set}$ are sometimes called \emphind{concrete categories}, and enjoy special properties, also see \cite{AdaHerStr1990}.
\end{enumerate}
\end{example}

\begin{example}
  \label{exa:cat-functor-examples}
We recall several examples of functors from a category $\cat{C}$ to a category $\cat{D}^\op$ that we have already seen in this book. Such functors are sometimes called \emph{contravariant}\index{functor!contravariant} functors from $\cat{C}$ to $\cat{D}$, and in this context, an actual functor from $\cat{C}$ to $\cat{D}$ is called a \emph{covariant}\index{functor!covariant} functor. %
\begin{enumerate}
  \item\label{itm:finitedownfunctor} The functor $\Down \colon \Pos_{f} \to (\DL_{f})^\op$ is defined by sending a finite poset to its lattice of down-sets, and an order-preserving function $f \colon P \to Q$ to the lattice homomorphism $\Down(f) \colon \Down(Q) \to \Down(P)$, $D \mapsto f^{-1}(D)$.
\item\label{itm:finitepowerfunctor} The functor $\mathcal{P} \colon \Set_{f} \to (\BA_{f})^\op$ is defined by sending a finite set to its Boolean power set algebra, and a function $f \colon X \to Y$ to the Boolean homomorphism $\mathcal{P}(f) \colon \mathcal{P}(Y) \to \mathcal{P}(X)$, $S \mapsto f^{-1}(S)$. If we view the category $\Set_{f}$ as the full subcategory of $\Pos_{f}$ consisting of the finite posets with the discrete order, then $\mathcal{P}$ is the restriction of $\Down$ to this full subcategory.
\item \label{itm:priestleyfunctors} The functor $\ClD \colon \Priest \to \DL^\op$ is defined by sending a Priestley space $X$ to the lattice $\ClD(X)$ of clopen down-sets of $X$, and a continuous order-preserving function $f \colon X \to Y$ to the lattice homomorphism $\ClD(f) \colon \ClD(Y) \to \ClD(X)$, $D \mapsto f^{-1}(D)$. If we view $\Pos_{f}$ as the full subcategory of $\Priest$ consisting of the finite Priestley spaces (see Exercise~\ref{exe:fullsubcats}), then $\ClD$ restricts to the functor $\Down$ in item \ref{itm:finitedownfunctor} above.
\item The functor $\Clp \colon \BoolSp \to \BA^\op$ is defined by sending a Boolean space $X$ to the Boolean algebra $\Clp(X)$ of clopen subsets of $X$, and a continuous function $f \colon X \to Y$ to the Boolean homomorphism $f^{-1} \colon \Clp(Y) \to \Clp(X)$. If we view $\Set_{f}$ as the full subcategory of $\BoolSp$ consisting of the finite discrete spaces, then $\Clp$ restricts to the functor $\mathcal{P}$ in item \ref{itm:finitepowerfunctor} above.
\item \label{itm:priestley-compatible-functor} Write $\Priest_{R^{\uparrow}}$ for the category whose objects are Priestley spaces and whose morphisms are upward Priestley compatible relations; see Definition~\ref{dfn:compatiblerelation}. Write $\DL_{\wedge}$ for the category whose objects are distributive lattices and whose morphisms are finite-meet-preserving functions. The results in Section~\ref{sec:unaryopduality} show in particular that there is a functor $\ClD \colon \Priest_{R^{\uparrow}} \to \DL_{\wedge}^\op$, which acts on objects as the functor described under item \ref{itm:priestleyfunctors} above, and on morphisms sends a relation $R \subseteq X \times Y$ to the finite-meet-preserving function $\forall_{R^{-1}}$.
\end{enumerate}
\end{example}
In each of the above examples, there is also a functor in the other direction. For example, recall from Theorem~\ref{thm:birkhoffduality} that there is a functor $\mathcal{J} \colon (\DL_f)^{\op} \to \Pos_f$ which is defined by sending a finite distributive lattice $L$ to its poset of join-irreducibles $\cJ(L)$, and a homomorphism $h \colon L \to M$ to the restriction of its lower adjoint, which is an order-preserving function $\cJ(M) \to \cJ(L)$.

How are the functors $\Down$ and $\mathcal{J}$ related to each other? A first guess might be that they form an \emph{isomorphism} between categories, in the sense that they are mutually inverse to each other, but this is not literally the case. (We give the notion of isomorphism of categories at the end of this section.) The composition of the functors $\Down$ and $\mathcal{J}$ is `almost' the identity, but only \emph{up to an isomorphism between objects}: for any lattice $L$, there is an isomorphism $\alpha_L$ between $L$ and $\Down(\cJ(L))$. While the isomorphism $\alpha_L$ has the lattice $L$ as a parameter, its definition is `consistent' as $L$ varies.

To precisely define an \emph{equivalence} of categories, we need to make precise what we mean by `consistent' in the previous sentence. To this end, we now introduce the notion of \emph{natural transformation}. A natural transformation should be thought of as a morphism \emph{between} two functors, and, as we shall see, the family of isomorphisms $\alpha := (\alpha_L)_{L \in \ob \DL_f}$ is an example. The family $\alpha$ is a special kind of natural transformation, because each of its components is an isomorphism; this is not required in the general definition of natural transformation. %

\nl{$\phi \colon F \To G$}{a natural transformation $\phi$ from a functor $F$ to a functor $G$}{}
\begin{definition}
  Let $F$ and $G$ be two functors from a category $\cat{C}$ to a category $\cat{D}$. A \emph{natural transformation} from $F$ to $G$, notation $\phi \colon F \To G$, is an $(\ob \cat{C})$-indexed collection $\phi = (\phi_A)_{A \in \ob \cat{C}}$ of morphisms in $\cat{D}$ that satisfies the following properties:
  \begin{enumerate}
  \item for every $A \in \ob \cat{C}$, $\dom(\phi_A) = FA$ and $\cod(\phi_A) = GA$;
  \item for every morphism $f \colon A \to A'$ in $\cat{C}$, we have $Gf \circ \phi_A = \phi_{A'} \circ Ff$; that is, the following square commutes:
    \begin{center}
    \begin{tikzpicture}
  \matrix (m) [matrix of math nodes,row sep=3em,column sep=3em,minimum width=3em]
  {
     FA & FA' \\
     GA & GA'\\};
  \path[-stealth]
    (m-1-1) edge node [above] {$Ff$} (m-1-2)
    (m-2-1) edge node [above] {$Gf$} (m-2-2)
    (m-1-1) edge node [left] {$\phi_A$} (m-2-1)
    (m-1-2) edge node [right] {$\phi_{A'}$} (m-2-2);
  \end{tikzpicture}
  \end{center}
\end{enumerate}
A \emph{natural isomorphism} is a natural transformation all of whose components are isomorphisms in the category $\cat{D}$.
\end{definition}
Under this definition, the family $\alpha = (\alpha_L)_{L \in \cat{DL}_f}$ defined above is an example of a natural isomorphism $\alpha \colon 1_{\cat{DL}_f} \To \Down \circ \cJ$, where we recall that $1_{\cat{DL}_f}$ denotes the identity functor on the category $\cat{DL}_f$.

Natural transformations are the morphisms in a (large) category of functors, that is, there is a composition of natural transformations and an identity natural transformation $1_F \colon F \To F$, for any functor $F$. Moreover, when $\phi \colon F \To G$ is a natural transformation and $K, L$ are functors that can be post-composed and pre-composed with $F$, then we also have well-defined natural transformations $K\phi \colon KF \To KG$ and $\phi_L \colon FL \To GL$. In this category, a natural transformation is a natural isomorphism if, and only if, it is an isomorphism in the category, that is, if there exists a natural transformation that is its two-sided inverse. For more precise statements (see Exercise~\ref{exe:naturaliso}).

\begin{example}
Let $P$ and $Q$ be preorders. Recall from Example~\ref{exa:cat-functor-examples} that functors from $P$ to $Q$ are (given by) order-preserving functions. Since there is at most one morphism between any two objects of $P$, or of $Q$, there is at most one natural transformation between any two functors in this setting. If $f, g \colon P \to Q$ are two order-preserving functions, viewed as functors, then there exists a natural transformation $\alpha \colon f \to g$ if, and only if, $f$ is pointwise below $g$, that is, $f(p) \leq g(p)$ for every $p \in P$.

For monoids $M$ and $N$, a natural transformation between two functors (that is, homomorphisms) $f, g \colon M \to N$ is given by an element $\alpha \in N$ such that, for every $m \in M$, $\alpha f(m) = g(m) \alpha$.
\end{example}
\subsection*{Equivalences and dualities}
Finally, we come to the categorical notion that is central to this book.
\begin{definition}\label{dfn:equivalence}
Let $\cat{C}$ and $\cat{D}$ be categories. 
  A pair of functors $F \colon \cat{C} \leftrightarrows \cat{D} \colon G$ is called an \emphind{equivalence} between $\cat{C}$ and $\cat{D}$ if there exist natural isomorphisms $\alpha \colon 1_{\cat{C}} \To GF$ and $\beta \colon FG \To 1_{\cat{D}}$.
An equivalence between $\cat{C}$ and $\cat{D}^\op$ is called a \emph{dual
equivalence} or \emph{duality} between $\cat{C}$ and $\cat{D}$. 

The categories
$\cat{C}$ and $\cat{D}$ are called
\emph{equivalent}\index{equivalent!categories} if there exists an equivalence
between them, and $\cat{C}$ and $\cat{D}$ are called \emph{dually equivalent}
or simply \emph{dual} if there exists a dual equivalence between
them.\index{dually equivalent categories}\index{dual categories}
\end{definition}

In this definition, since the notion of natural isomorphism is symmetric, the direction of $\alpha$ and $\beta$ clearly does not matter, and the asymmetric choice, where $\alpha$ has the identity functor as codomain and $\beta$ has the identity functor as domain, may look a bit strange. However, we will see shortly that the notion of equivalence is a special case of the notion of adjunction, where the direction of the natural transformations \emph{does} matter. As is common in the literature, we formulated Definition~\ref{dfn:equivalence} to resemble the definition of adjunction (Definition~\ref{dfn:adjunction}) as much as possible.

The following theorem, whose proof relies on the axiom of choice, is very
useful when proving that a functor is an equivalence, as it allows us to avoid
explicitly defining natural isomorphisms, but instead to just check three
properties for one of the two equivalence functors. One may view this theorem as
a categorified version of the fact that if a function $f \colon C \to D$ between
sets is bijective, then it has a two-sided inverse.
\begin{theorem}\label{thm:charequivalence} 
Let $F \colon \cat{C} \to \cat{D}$ be a functor. The following are equivalent:
\begin{enumerate}
\item[(i)] There exists a functor $G \colon \cat{D} \to \cat{C}$ such that $(F,G)$ is an equivalence between $\cat{C}$ and $\cat{D}$.
\item[(ii)] The functor $F$ is full, faithful and essentially surjective.
\end{enumerate}
\end{theorem}
\begin{proof}
  We leave the implication (i) $\Rightarrow$ (ii) to the reader (see Exercise~\ref{exe:necequivalence}). Conversely, suppose that $F$ is full, faithful and essentially surjective. Since $F$ is essentially surjective, for every object $B$ in $\cat{D}$, pick an object, $GB$, in $\cat{C}$ and an isomorphism $\beta_B \colon F(GB) \to B$. Let $g \colon B \to B'$ be a morphism in $\cat{D}$, and define $g' \colon F(GB) \to F(GB')$ to be the composite map $(\beta_{B'})^{-1} \circ g \circ \beta_B$. Since $F$ is full and faithful, there exists a unique morphism, $Gg \colon GB \to GB'$, such that $F(Gg) = g' = (\beta_{B'})^{-1} \circ g \circ \beta_B$. In other words, applying $\beta_{B'}$ to both sides of this equality, $Gg$ is the unique morphism $GB \to GB'$ such that $\beta_{B'} F(Gg) = g \beta_B$. (We will often omit the symbol $\circ$ in the remainder of this proof to improve readability.)

  We show that the assignment $G \colon \cat{D} \to \cat{C}$ given by $G(B)=GB$ and $G(g)=Gg$ is a functor. For any object $B$ in $\cat{D}$, notice that $\beta_B F(1_{GB}) = \beta_B 1_{F(GB)} =\beta_B= 1_B \beta_B$. Thus, by the uniqueness in the definition of $G$ on morphisms, $G(1_B) = 1_{GB}$. Now let $g_1 \colon B_1 \to B_2$ and $g_2 \colon B_2 \to B_3$ be a pair of composable morphisms in $\cat{D}$. We will prove that $G(g_2 \circ g_1) = G(g_2) \circ G(g_1)$ by showing that the morphism $G(g_2) \circ G(g_1)$ satisfies the defining property of $G(g_2 \circ g_1)$. Indeed, using that $F$ is a functor and the defining properties of $G(g_1)$ and $G(g_2)$, we get
  \[ \beta_{B_3} F(G(g_2) G(g_1)) = \beta_{B_3} FG(g_2) FG(g_1) = g_2 \beta_{B_2} FG(g_1) = g_2 g_1 \beta_{B_1}.\]
  It is immediate from the definition of $G$ that $\beta \colon FG \To 1_{\cat{D}}$ is a natural transformation, and thus a natural isomorphism.

  Finally, we construct a natural isomorphism $\alpha \colon 1_{\cat{C}} \To GF$. For any object $A$ in $\cat{C}$, let $\alpha_A$ be the unique morphism $A \to GF(A)$ such that $F(\alpha_A) = \beta_{FA}^{-1}$. We leave it as a highly instructive exercise to the reader to check that $\alpha$ is a natural isomorphism (see Exercise~\ref{exe:sufequivalence-alpha} for hints).
\end{proof}

\begin{example}\label{exa:cat-duality-examples}
  All the functors in Example~\ref{exa:cat-functor-examples} are part of a dual equivalence. Indeed, consider again the functor $\Down \colon \cat{Pos}_f \to (\cat{DL}_f)^\op$ from Example~\ref{exa:cat-functor-examples}.\ref{itm:finitedownfunctor}. We check explicitly that $\Down$ is full, faithful, and essentially surjective. It was shown in Proposition~\ref{prop:birkhoff} of Chapter~\ref{ch:order} that every finite distributive lattice is isomorphic to some lattice of the form $\Down(P)$, where we can take for $P$ the finite poset of join-prime elements of the lattice. Thus, $\Down$ is essentially surjective. It was shown in Proposition~\ref{prop:finDLmorphisms} that, for arbitrary finite posets $P, Q$, the assignment $f \mapsto \Down(f) := f^{-1}$ is a bijection between the set of order-preserving functions from $P$ to $Q$ and the set of lattice homomorphisms from $\Down(Q)$ to $\Down(P)$. This means exactly that $\Down$ is full and faithful. Therefore, by Theorem~\ref{thm:charequivalence}, $\Down$ is part of a dual equivalence. The functor $\mathcal{J} \colon (\cat{DL}_f)^\op \to \cat{Pos}_f$ gives the (up to natural isomorphism unique) functor in the other direction. To see this, one may either trace the proof of Theorem~\ref{thm:charequivalence} in this specific case and check that $\mathcal{J}$ works as a choice for $G$ in that proof, or directly exhibit the required natural isomorphisms. A natural isomorphism $\alpha \colon 1_{(\cat{DL}_f)^\op} \To \Down \cJ$ is given by the family of maps $\widehat{(-)} \colon L \to \Down(\cJ(L))$, which are isomorphisms according to Proposition~\ref{prop:birkhoff}. 
A natural isomorphism $\beta \colon 1_{\cat{Pos}_f} \To \cJ \Down$ is given by
the family of isomorphisms $\beta_P \colon P \to \cJ (\Down(P))$ that send $p
\in P$ to ${\downarrow} p$, given by
Proposition~\ref{prop:finite-poset-double-dual}. This concludes our detailed
proof of Theorem~\ref{thm:birkhoffduality}. 

You are asked to supply a similar
proof of the Priestley Duality Theorem~\ref{thm:priestleyduality} via
Theorem~\ref{thm:charequivalence} in
Exercise~\ref{exe:priestley-is-duality} below, also see
Theorem~\ref{thm:priestley-duality-long} for a direct proof that does  not use
Theorem~\ref{thm:charequivalence}.

Note that, once the functors have been defined,  these proofs always follow the
same pattern, and have very little to do with the precise categories at hand.
The real work in establishing a (dual) equivalence between two specific
categories is in defining a full, faithful, essentially surjective functor from
one category to another. The rest is, as some category theorists like to call
it, `\emphind{abstract nonsense}'. 
\end{example}

In the next subsection, we will encounter some pairs of functors that do not form equivalences, but that do enjoy a looser bond.

\subsection*{Adjunctions}
An important weakening of the notion of equivalence is that of an \emph{adjunction} between categories. We already encountered adjunctions in the context of preordered sets in the very first section of this book (Definition~\ref{dfn:poset-adjunction-def}). We used adjunctions between posets several times in Section~\ref{sec:quotients-and-subs}, to formulate the correspondence between quotient lattices and subspaces. Adjunctions between categories generalize adjunctions between preorders, if we view a preorder as a category with at most one morphism between two objects (see Example~\ref{exa:smallcats}). We have already seen examples of such categorical adjunctions `in action': both the \emph{Boolean envelope} (Section~\ref{sec:boolenv-duality}) and the \emph{free distributive lattice} (Section~\ref{sec:free-description}) are examples of adjoint constructions; see Example~\ref{exa:adjunctions} below for details. 

\begin{definition}\label{dfn:adjunction}
  Let $\cat{C}$ and $\cat{D}$ be categories, and let $F \colon \cat{C} \leftrightarrows \cat{D} \colon G$ be a pair of functors between them. We say that $(F, G)$ is an \emph{adjunction}\index{adjunction!between categories} if there exist natural transformations $\eta \colon 1_{\cat{C}} \to GF$ and $\epsilon \colon FG \to 1_{\cat{D}}$ satisfying the following two properties:
  \begin{enumerate}
  \item\label{itm:triangle-F} for any object $A$ of $\cat{C}$, $\epsilon_{FA} \circ F(\eta_A) = 1_{FA}$, and
  \item\label{itm:triangle-G} for any object $B$ of $\cat{D}$, $G(\epsilon_B) \circ \eta_{GB} = 1_{GB}$.
  \end{enumerate}
  In this situation, $F$ is called \emph{left} or \emph{lower adjoint} to $G$ and $G$ is called \emph{right} or \emph{upper adjoint} to $F$; a common notation for this is $F \dashv G$. \index{adjoint!between categories} When the natural transformations $\eta$ and $\epsilon$ are explicitly specified, the natural transformation $\eta$ is called the \emphind{unit} of the adjunction, and $\epsilon$ is called the \emphind{co-unit} of the adjunction.
\end{definition}

  Properties \ref{itm:triangle-F} and \ref{itm:triangle-G} in Definition~\ref{dfn:adjunction} are called \emphind{triangle identities}, as they can be expressed by saying that the following two triangles of natural transformations commute:
  \begin{center}
    \begin{tabular}{cc}
  \begin{tikzpicture}
  \matrix (m) [matrix of math nodes,row sep=3em,column sep=3em,minimum width=3em]
  {
     F & FGF \\
       & F\\};
  \path[-stealth]
    (m-1-1) edge node [above] {$F\eta$} (m-1-2)
    (m-1-1) edge node [below] {$1_{F}$} (m-2-2)
    (m-1-2) edge node [right] {$\epsilon_{F}$} (m-2-2);
  \end{tikzpicture}
      &
  \begin{tikzpicture}
      \matrix (m) [matrix of math nodes,row sep=3em,column sep=3em,minimum width=3em]
  {
     G & GFG \\
       & G   \\
   };
  \path[-stealth]
    (m-1-1) edge node [above] {$\eta_G$} (m-1-2)
    (m-1-1) edge node [below] {$1_{G}$} (m-2-2)
    (m-1-2) edge node [right] {$G\epsilon$} (m-2-2);

  \end{tikzpicture}
      \end{tabular}
    \end{center}
    While Definition~\ref{dfn:adjunction} does not make the unit-co-unit pair $(\eta, \epsilon)$ part of the structure of an adjunction, and there may exist distinct choices of unit and co-unit for the same pair of functors, in the examples we consider in this book, there will always be a natural choice for $\eta$ and $\epsilon$, that we will specify when needed. 
An adjoint pair of functors is an \emph{equivalence} exactly if the unit and co-unit can be chosen to be natural isomorphisms, in the sense that, for any equivalence, \emph{there exist} natural isomorphisms that satisfy the triangle identities (see Exercise~\ref{exe:adjoint-equivalence}).\footnote{There is a slight subtlety here: not \emph{every} pair of natural isomorphisms $1_{\cat{C}} \to GF$ and $FG \to 1_{\cat{D}}$ satisfies the triangle identities, but there is always \emph{some} choice that does. An equivalence together with explicit natural isomorphisms that satisfy the triangle identities is sometimes called an \emph{adjoint equivalence} in the literature, to distinguish it from the non-adjoint situation, but for all of the (dual) equivalences in this book there is a straightforward choice of natural isomorphisms that do satisfy the triangle identities, so we need not make this distinction.}

The reader may recall that, for preorders $P$ and $Q$, two order-preserving functions $f \colon P \leftrightarrows Q \colon g$ form an adjunction if, for any $p \in P$ and $q \in Q$, $f(p) \preceq_Q q$ if, and only if, $p \preceq_P g(q)$. An equivalent definition, given in Exercise~\ref{exe:adjunctions}, says that the function $1_P$ is pointwise below the function $gf$, and that $fg$ is pointwise below $1_Q$, which shows why the notion of adjunction between categories generalizes that between preorders. In the setting of categories, we have the following statement, which is closer to the original definition of adjunction between preorders.
\begin{proposition}\label{prop:hom-set-bijections}
If $F \colon \cat{C} \leftrightarrows \cat{D} \colon G$ is an adjunction, then there exists, for any objects $A$ in $\cat{C}$ and $B$ in $\cat{D}$, a bijection between the sets $\Hom_{\cat{D}}(FA, B)$ and $\Hom_{\cat{C}}(A, GB)$.
\end{proposition}
In other words, the cardinality of the set of morphisms from $FA$ to $B$ in $\cat{D}$ is equal to that of the set of morphisms from $A$ to $GB$ in $\cat{C}$. Since, in a preorder, the number of morphisms between two objects is at most one, the reader should now clearly see the correspondence with the preorder definition of adjunction.

We note here that Proposition~\ref{prop:hom-set-bijections} can also lead to an alternative definition of adjunction between categories, by adding the information that the family of bijections $\Hom_{\cat{D}} (FA, B) \to \Hom_{\cat{C}}(A, GB)$ in fact forms a \emph{natural isomorphism} between  $\Hom_{\cat{D}}(F-, -)$ and $\Hom_{\cat{C}}(-, G-)$, when these assignments are viewed as functors from $\cat{C}^\op \times \cat{D}$ to $\cat{Set}$. Proving that this definition is equivalent to our Definition~\ref{dfn:adjunction} is beyond the scope of this book, and we will work with Definition~\ref{dfn:adjunction} as the official definition of adjunction. For a treatment of adjunctions that takes the $\Hom$-set definition as primitive, see, for example, \cite[Ch.~4]{Mac1971}.

 An important property of adjunctions is that they compose: when $F \colon \cat{C} \leftrightarrows \cat{D} \colon G$ and $F' \colon \cat{D} \leftrightarrows \cat{E} \colon G'$ are functors such that $F$ is left adjoint to $G$ and $F'$ is left adjoint to $G'$, then $F'F$ is left adjoint to $GG'$ (see Exercise~\ref{exe:adjunctions-compose}).

\begin{example}\phantomsection\label{exa:adjunctions}
  \begin{enumerate}
  \item \label{itm:freeDL} Let $U \colon \cat{DL} \to \cat{Set}$ denote the \index{forgetful functor}\index{functor!forgetful} forgetful functor of Example~\ref{exa:cat-functor-examples}.\ref{itm:forgetfulfunctor}. The left adjoint to this forgetful functor is the \emphind{free distributive lattice} functor, $F_{\cat{DL}} \colon \cat{Set} \to \cat{DL}$, whose object part was constructed in Section~\ref{sec:free-description}. The universal property allows one to extend this object assignment to a functor which is left adjoint to $U$. You are asked to supply the details of this example in Exercise~\ref{exe:freeDL}.
\item Similarly to the previous item, the \emphind{Boolean envelope} of a distributive lattice (see Section~\ref{sec:boolenv-duality}) is the image of the left adjoint $(-)^- \colon \cat{DL} \to \cat{BA}$ to the forgetful functor $U \colon \cat{BA} \to \cat{DL}$. This item, together with the previous one and the fact that adjunctions compose, gives a more abstract proof of Lemma~\ref{lem:freeDLplusboolenv}: the composition of $(-)^-$ and $F_{\cat{DL}}$ gives a left adjoint $F_{\cat{BA}} \colon \cat{Set} \to \cat{BA}$ to the forgetful functor $U \colon \cat{BA} \to \cat{Set}$.
\item\label{itm:free-join-semi} The \emphind{free join-semilattice} over a poset $P$ was given in Exercise~\ref{exe:fingendown}, as the lattice of finitely generated down-sets of $P$. This can be used to show that the forgetful functor $U \colon \cat{SL} \to \cat{Pos}$, where $\cat{SL}$ is the category of semilattices, has a left adjoint, $\Downfin \colon \cat{Pos} \to \cat{SL}$. 
\item \label{itm:DLplus} Let $\cat{DL}^+$ denote the (non-full) subcategory of distributive lattices that are isomorphic to a lattice of the form $\Down(P)$, for $P$ a poset, with \emph{complete} homomorphisms between them. (The class of distributive lattices that are isomorphic to the down-set lattice of some poset may be characterized as ``completely distributive lattices that have enough completely join-irreducible elements''. Completely distributive lattices will be studied extensively in Section~\ref{sec:dom}.) There is a (non-full) inclusion functor $I \colon \cat{DL}^+ \to \cat{DL}$. The functor $I$ has a left adjoint, $()^\delta$, which may be constructed as follows. For any distributive lattice $L$, let $L^\delta$ be the lattice $\Down(X_L)$, where $X_L$ is the poset underlying the Priestley dual space of $L$. For any distributive lattice homomorphism $h \colon L \to M$, let $h^\delta$ be the complete homomorphism $(f_h)^{-1} \colon L^\delta \to M^\delta$, where $f_h \colon X_M \to X_L$ is the (order-preserving) function dual to $h$. The functor $()^\delta$ is (naturally isomorphic to) the so-called \emphind{canonical extension} functor.

  The unit of the adjunction $I \dashv ()^\delta$ is the function $\eta_L \colon L \to \Down(X_L)$ that sends $a \in L$ to $\hat{a}$. The co-unit of the adjunction is defined, for $C$ an object of $\cat{DL}^+$, as the map $\epsilon_C \colon \Down(X_C) \to C$ which sends a down-set $D \subseteq X_C$ to $\bigvee_{x \in D} \bigwedge F_x$, where the join and meet are taken in the complete lattice $C$.
\item \label{itm:freemonoid} 
Another basic example of a free construction, which
will become relevant in Chapter~\ref{ch:AutThry} of this book, is the
\emphind{free monoid} on an input alphabet. Recall from
Example~\ref{exa:smallcats}.\ref{ite:monoid} that a \emphind{monoid} is a small
category with one object; that is, a \emph{monoid} is a tuple $(M, \cdot, 1)$
where $\cdot$ is an associative binary operation on $M$, that is, $u \cdot (v \cdot
w) = (u \cdot v) \cdot w$ for all $u, v, w \in M$, and $1$ is a neutral element, that is,
$1 \cdot u = u = u \cdot 1$ for all $u \in M$. In this context, a functor is
called a \emphind{monoid homomorphism}, that is, a function $f \colon M \to N$ that
preserves the multiplication and the identity element. We denote the category whose objects are 
monoids and whose morphisms are monoid homomorphisms by \cat{Mon}. Note that this category $\cat{Mon}$ is different from the category associated to a single monoid $M$ that we saw in Example~\ref{exa:smallcats}.\ref{ite:monoid}. 
We have an obvious forgetful functor $U \colon
\cat{Mon} \to \cat{Set}$. This functor has a left adjoint, which we denote
$(-)^* \colon \cat{Set} \to \cat{Mon}$, and which sends a set $X$ to the
collection of finite sequences over $X$, equipped with multiplication given by
concatenation of sequences, for which the empty sequence is an identity
element. Exercise~\ref{exe:freemonoid} asks you to show that this indeed yields
a left adjoint to the forgetful functor.
\nl{$A^*$}{free monoid over a set $A$}{}
\end{enumerate}
\end{example}
The similarity between items \ref{itm:freemonoid} and \ref{itm:freeDL} in Example~\ref{exa:adjunctions} is no coincidence: general results of universal algebra, going back to Birkhoff, show that any so-called \emph{variety} of algebras admits a free construction (see Exercise~\ref{exe:birkhoff}). In fact, \emph{any} adjunction between categories may be understood as `a universal construction', in the following sense.
\begin{theorem}\label{thm:adj-universal-arrow}
Let $G \colon \cat{D} \to \cat{C}$ be a functor. The following are equivalent:
\begin{enumerate}
  \item[(i)] The functor $G$ has a left adjoint,
  \item[(ii)] For every object $A$ of $\cat{C}$, there exist an object $F(A)$ of $\cat{D}$ and a morphism $\eta_A \colon A \to G(F(A))$ of $\cat{C}$ such that, for any object $B$ of $\cat{D}$ and any $\cat{C}$-morphism $f \colon A \to G(B)$, there is a unique $\cat{D}$-morphism $\bar{f} \colon F(A) \to B$ such that $(G\bar{f}) \circ \eta_A = f$.
\end{enumerate}
\end{theorem}
This theorem generalizes the characterization of functors that are part of an equivalence, Theorem~\ref{thm:charequivalence}, in the following sense: if a functor $G$ is essentially surjective, then one may always choose an object $F(A)$ and an \emph{iso}morphism $\eta_A \colon A \to G(F(A))$; if $G$ is then moreover full and faithful, then any such choice will satisfy the universal property stated in (ii) of Theorem~\ref{thm:adj-universal-arrow}, and  will thus have a left adjoint, with which $G$ will then form an equivalence. The proof of Theorem~\ref{thm:adj-universal-arrow} is beyond the scope of this book; see for example \cite[Section 3.1]{Borceux1}.

We end our discussion of adjunctions by connecting it to equivalences in a different useful way: one may always obtain a (dual) equivalence from a (contravariant) adjunction. Indeed,  Theorem~\ref{thm:charequivalence} can be applied to prove that any (contravariant) adjunction restricts to a maximal (dual) equivalence. This associated duality is obtained by restricting the functors to the full subcategories given by the objects for which the natural map between the object and its double dual is an isomorphism. We record this fact here, and leave the precise proof as Exercise~\ref{exe:proofofrestricttoduality}.
\begin{corollary}\label{cor:restricttoduality}
Let $F \colon \cat{C} \leftrightarrows \cat{D} \colon G$ be an adjunction. Denote by $\mathrm{Fix}(\cat{C})$ the full subcategory of $\cat{C}$ consisting of those objects $A$ for which the unit $\eta_A \colon A \to GFA$ is an isomorphism. Similarly, denote by $\mathrm{Fix}(\cat{D})$ the full subcategory of $\cat{D}$ consisting of those objects $B$ for which the co-unit $\epsilon_B \colon FGB \to B$ is an isomorphism. Then $F$ and $G$ restrict to an equivalence between $\mathrm{Fix}(\cat{C})$ and $\mathrm{Fix}(\cat{D})$.
\end{corollary}

\subsection*{Isomorphisms}
Finally, the notion of isomorphism between categories is easy to state, but often too restrictive.

\begin{definition}\label{dfn:iso-between-cats}\index{isomorphism!between categories}
A functor $F \colon \cat{C} \rightarrow \cat{D}$ is called an \emphind{isomorphism} between categories $\cat{C}$ and $\cat{D}$ if there exists a functor $G \colon \cat{D} \rightarrow \cat{C}$ such that $G \circ F$ is equal to the identity functor on $\cat{C}$ and $F \circ G$ is equal to the identity functor on $\cat{D}$. 
A \emphind{dual isomorphism} between $\cat{C}$ and $\cat{D}$ is an isomorphism between $\cat{C}$ and $\cat{D}^\op$. 
Two categories $\cat{C}$ and $\cat{D}$ are called \emphind{isomorphic} if there exists an isomorphism between them, and \emphind{dually isomorphic} if there exists a dual isomorphism between them.
\end{definition}
Note that an isomorphism between locally small categories is, by definition, a functor that is an isomorphism, in the sense of Definition~\ref{dfn:internaliso}, when viewed as a morphism in the (large) category of locally small categories.

\begin{example}\label{exa:topfiniso}
As indicated in Exercise~\ref{exer:TOPfin}, the category $\TopCat_f$ of finite $T_0$ topological spaces is isomorphic to the category $\Pos_f$ of finite partially ordered sets. Indeed, for a finite $T_0$ space $(X,\tau)$, let $F(X,\tau) := (X,\leq_{\tau})$ be the specialization order of $(X,\tau)$, and, conversely, for a finite order $(X,\leq)$, let $G(X,\leq) := (X,\alpha(X))$ be the set $X$ equipped with the Alexandrov topology induced by $\leq$. These object assignments $F$ and $G$ extend to functors by sending any function to itself.
\end{example}

\begin{example}\label{exa:stcomp-kord-iso}
In Section~\ref{sec:comp-ord-sp} in Chapter~\ref{chap:TopOrd}, we studied a bijective correspondence between compact ordered spaces and stably compact spaces: to a compact ordered space $(X, \tau, \leq)$, assign the stably compact space $(X, \tau^{\uparrow})$, and conversely to a stably compact space $(X,\sigma)$, assign the compact ordered space $(X,\sigma^p,\leq_\sigma)$. We saw in Proposition~\ref{prop:compordaspatch} that these assignments are mutually inverse. Compact ordered spaces form a category whose morphisms are the continuous order-preserving maps. If we equip the class of stably compact spaces with the \emph{proper maps} defined in Exercise~\ref{exer:propermaps}, then it will become isomorphic to the category of compact ordered spaces; you will be asked to supply the details in Exercise~\ref{exe:stcomp-kord-iso} below. In the next chapter, Section~\ref{sec:StoneSpaces}, we will study the restriction of this isomorphism of categories to the full subcategory of Priestley spaces; the image on the side of stably compact spaces is the class of \emph{spectral spaces}.
\end{example}

\ourexercises

\begin{ourexercise}\label{exe:fullsubcats}
This exercise guides you through the details of identifying $\Pos_{f}$ with a full subcategory of $\Priest$. 
\begin{enumerate}
\item Prove that the topology on a finite compact Hausdorff space is discrete.
\item\label{itm:iotaobjects} Prove that, if $(X,\leq)$ is a finite poset, then $\iota(X,\leq) := (X,\tau,\leq)$, where $\tau$ is the discrete topology on $X$, is a Priestley space.
\item Prove that a function $f \colon X \to Y$ between finite Priestley spaces is a morphism in $\Priest$ if, and only if, $f$ is order preserving.
\item Conclude that the object assignment $\iota$ in item \ref{itm:iotaobjects} extends to a full and faithful functor $\iota \colon \Pos_{f} \to \Priest$.
\end{enumerate}
\end{ourexercise}

\begin{ourexercise}\label{exe:self-dual}
  A category is called \emph{self-dual} if $\cat{C}$ is equivalent to $\cat{C}^\op$.
  \begin{enumerate}
  \item \label{ite:rel-relop} Verify that the category $\cat{Rel}$ is self-dual, and in fact, \emph{isomorphic} to the category $\cat{Rel}^\op$.

    We now outline a proof that $\cat{Set}$ is not self-dual.
  \item Exhibit an object $1$ in $\cat{Set}$ with the property that for any object $X$ there is exactly one morphism $X \to 1$. (Such an object is called a \emph{terminal} object in a category, and is unique up to isomorphism, as we will see in the next section.)

  \item Show that $\cat{Set}$ has the following property: for any two distinct morphisms $f, g \colon X \to Y$, there exists a morphism $x \colon 1 \to X$ such that $fx \neq fy$. A category with this property is said to be \emph{generated by its terminal object} or \emphind{well-pointed}.

  \item What is the terminal object in $\cat{Set}^\op$?

  \item Prove that $\cat{Set}^\op$ is not generated by its terminal object.

  \item  Conclude that $\cat{Set}$ is not self-dual.
    \end{enumerate}
  \end{ourexercise}

\begin{ourexercise}\label{exe:skeleton}
This exercise outlines a possible definition for a \emphind{skeleton}\index{category!skeleton} of a category. The special case of preorders was already given in Exercise~\ref{exe:reflection}.
Let $\cat{C}$ be a category. Choose a collection of objects $S$ in $\cat{C}$
such that every object in $\cat{C}$ is isomorphic to exactly one object of $S$;
note that this step uses the axiom of choice for the
collection $\ob \cat{C}$, which is in general not a set, but we ignore
set-theoretic issues here. 
Prove that the inclusion functor of the full subcategory $\cat{S}$ on the
collection of objects $S$ is an equivalence.
\end{ourexercise}

\begin{ourexercise}\label{exe:naturaliso}
  Let $F, G, H \colon \cat{C} \to \cat{D}$ be functors, let $\phi \colon F \To G$ and $\psi \colon G \To H$ be natural transformations.
  \begin{enumerate}
  \item Give a definition of the composition $\psi \circ \phi \colon F \To H$, and show that it is again a natural transformation.
    \item Show that there is an identity natural transformation $1_F \colon F \To F$.
    \item\label{ite:map-functor-nat} Given a functor $K \colon \cat{D} \to \cat{E}$, show that $K\phi \colon KF \To KG$, defined at an object $C$ of $\cat{C}$ by $(K\phi)_A := K(\phi_A)$, is a natural transformation.
    \item Given a functor $L \colon \cat{B} \to \cat{C}$, show that  $\phi_L \colon FL \To GL$, defined at an object $B$ of $\cat{B}$ by $(\phi_L)_B := \phi_{LB}$, is a natural transformation.
    \item Prove that, if $\theta \colon G \to F$ is a natural transformation such that $\theta \circ \phi = 1_F$ and $\phi \circ \theta = 1_G$, then $\phi$ and $\theta$ are natural isomorphisms.
      \item Prove that, if $\phi$ is a natural isomorphism, then defining $\theta_A := (\phi_A)^{-1}$ for every $A \in \ob \cat{C}$ yields a natural transformation $\theta$ that is a two-sided inverse for $\phi$.
  \end{enumerate}
\end{ourexercise}

\begin{ourexercise}\label{exe:necequivalence}
Let $F \colon \cat{C} \leftrightarrows \cat{D} \colon G$ be an equivalence of categories. Prove that the functors $F$ and $G$ are full, faithful and essentially surjective.
\end{ourexercise}

\begin{ourexercise}\label{exe:sufequivalence-alpha}
  This exercise asks you to complete the proof of the sufficiency direction of Theorem~\ref{thm:charequivalence}. Consider the family of morphisms $(\alpha_A)_{A \in \ob \cat{C}}$ defined at the end of the proof of that theorem in the text.
  \begin{enumerate}
  \item Prove that $\alpha \colon 1_{\cat{C}} \Rightarrow GF$ is a natural transformation. \hint{For a morphism $f \colon A \to A'$ in $\cat{C}$, show first the equality
    \[\beta_{FA'} \circ F(\alpha_A \circ f) = \beta_{FA'} \circ F(GF(f) \circ \alpha_A),\] and then explain why this is enough to conclude that $\alpha_{A'} f = GF(f) \alpha_A$.}
  \item Prove that $\alpha_A$ is an isomorphism for every object $A$. \hint{Consider a morphism $\gamma_A \colon GFA \to A$ such that $F\gamma_A = \beta_{FA}$.}

  \end{enumerate}
\end{ourexercise}

\begin{ourexercise}\label{exe:priestley-is-duality}
  Use Theorem~\ref{thm:charequivalence} to show that the results from Chapter~\ref{ch:priestley} imply that the functors from Example~\ref{exa:cat-functor-examples}.\ref{itm:priestleyfunctors} form a dual equivalence. In other words, prove Theorem~\ref{thm:priestleyduality} (Priestley duality).

  {\it Note.} We will give more details, and in particular a somewhat more explicit proof of Priestley duality in Section~\ref{sec:duality-categorically} below.
\end{ourexercise}

\begin{ourexercise}\label{exe:adjoint-equivalence}
  Show that every equivalence can be improved to an adjoint equivalence. That is, let $F \colon \cat{C} \leftrightarrows \cat{D} \colon G$ be an equivalence of categories. Prove that there exist natural isomorphisms $\eta \colon 1_{\cat{C}} \to GF$ and $\epsilon \colon FG \to 1_{\cat{D}}$ satisfying the triangle identities. \hint{Consider the proof of Theorem~\ref{thm:charequivalence}.}
\end{ourexercise}

\begin{ourexercise}\label{exe:freeDL}
  This exercise asks you to supply the details omitted in Example~\ref{exa:adjunctions}.\ref{itm:freeDL}, about the free distributive lattice on a set. For any set $X$, denote by $\eta_X \colon X \to UF_{\cat{DL}}(X)$ the embedding of $X$ into the free distributive lattice over $X$.
  \begin{enumerate}
  \item Prove that $F_{\cat{DL}}$ extends to a functor $\mathbf{Set} \to \mathbf{DL}$.
  \item Explicitly define a co-unit $\epsilon \colon F_{\cat{DL}} U \to 1_{\cat{DL}}$.
  \item Prove that $\eta$ and $\epsilon$ are natural transformations satisfying the triangle identities.
    \end{enumerate}
  \end{ourexercise}

\begin{ourexercise}\label{exe:freemonoid}
Repeat Exercise~\ref{exe:freeDL} for the free monoid functor from Example~\ref{exa:adjunctions}.\ref{itm:freemonoid} (see Exercise~\ref{exe:birkhoff} for a general result of which this exercise and the previous one are specific instances).
\end{ourexercise}

\begin{ourexercise}\label{exe:birkhoff} This exercise generalizes Exercise~\ref{exe:freeDL} and \ref{exe:freemonoid} above. Let $\mathcal{V}$ be a variety of finitary algebraic structures with at least one constant symbol in the signature. Denote by $\cat{V}$ the category whose objects are the algebras in $\mathcal{V}$, and whose morphisms are the homomorphisms, and denote by $U \colon \cat{V} \to \cat{Set}$ the forgetful functor. Prove that $U$ has a left adjoint $F_{\cat{V}} \colon \cat{Set} \to \cat{V}$. \hint{This exercise essentially asks to prove a well-known theorem of Birkhoff, namely, that any finitary variety has free algebras. The proof mimicks exactly the algebraic construction of the free distributive lattice given in Section~\ref{sec:free-description}. Also see \cite[Ch.~II]{BurSan2000}. 

The reason for assuming that the signature has a constant symbol is that, if the signature has no constant symbol, then $F_{\cat{V}}(\emptyset)$ has to be empty, and some authors do not allow the empty set as an algebra.}
  \end{ourexercise}

\begin{ourexercise}\label{exe:proof-of-adjunction-equiv-def}
  Prove Proposition~\ref{prop:hom-set-bijections}.
\end{ourexercise}

\begin{ourexercise}\label{exe:adjunctions-compose}
Let $F \colon \cat{C} \leftrightarrows \cat{D} \colon G$ and $F' \colon \cat{D} \leftrightarrows \cat{E} \colon G'$ be adjunctions. Prove that $F'F \colon \cat{C} \leftrightarrows \cat{D} \colon GG'$ is an adjunction.
\end{ourexercise}

\begin{ourexercise}\label{exe:equivalence-is-equivalence}
    Prove that `equivalence' is an equivalence relation on
    categories. That is,
    show that any category $\cat{C}$ is equivalent to itself, that $\cat{C}$
    equivalent to $\cat{D}$ implies $\cat{D}$ equivalent to $\cat{C}$, and that
    if $\cat{C}$ is equivalent to $\cat{D}$ and $\cat{D}$ is equivalent to
    $\cat{E}$, then $\cat{C}$ is equivalent to $\cat{E}$.
\end{ourexercise}
\begin{ourexercise}\label{exe:proofofrestricttoduality}
Show that Corollary~\ref{cor:restricttoduality} indeed follows from Theorem~\ref{thm:charequivalence}.
\end{ourexercise}

\begin{ourexercise}\label{exe:stcomp-kord-iso}
Using the results of Exercise~\ref{exer:propermaps}, prove that the bijective correspondence of Theorem~\ref{thrm:COSpace-StabCompSp} extends to an isomorphism of categories, as outlined in Example~\ref{exa:stcomp-kord-iso}. (Also see Theorem~\ref{thm:Stone-isom-Priestley} in the next chapter.)
\end{ourexercise}

\section{Constructions inside categories}\label{sec:internal}
We now come to the internal perspective on a category. We will consider many types of internal notions and constructions in a category: monomorphisms and epimorphisms; products and coproducts; equalizers and coequalizers; and, most generally, limits and colimits. The reader will notice that these notions come in pairs; why? As we discussed early on in Section~\ref{sec:external}, any notion in a category $\cat{C}$ can also be considered in the so-called \emph{opposite} category $\cat{C}^\op$, in which the direction of the morphisms is reversed. This means that, for any concept of a categorical `gadget', one can also define a `co-gadget' in $\cat{C}$ to be ``a thing in $\cat{C}$ that is a gadget when viewed in $\cat{C}^\op$''. The name for the dual notion is often formed by prefixing `co-' to the original name: product and coproduct, equalizer and coequalizer, limit and colimit.\footnote{A notable exception to this terminological rule is `covariant' vs. `contravariant' functor, defined in the previous section, and also `monomorphism' vs. `epimorphism' introduced here.} %
In this section, we often explicitly introduce both a notion and its dual. While not the most economical, this will allow the reader to see both the formal similarity and the fact that these notions can have very different concretizations in specific examples of categories.

\begin{notation}
Throughout the rest of this chapter, we will often omit the symbol $\circ$ for composition within a category, and we write strings like `$hgf$' as shorthand for `$h \circ g \circ f$'.
\end{notation}

\subsection*{Monomorphisms, epimorphisms and coequalizers}
We already saw one type of special morphism in a category: an isomorphism, defined in  Definition~\ref{dfn:internaliso}. We now introduce two related, but weaker, notions.
\begin{definition}
Let $f \colon A \to B$ be a morphism in a category $\cat{C}$. We say that $f$ is a 
\begin{itemize}
\item \emphind{monomorphism} if, for any $g, h \colon A' \to A$ in $\cat{C}$ such that $fg = fh$, we have $g = h$.
\item \emphind{epimorphism} if, for any $g, h \colon B \to B'$ in $\cat{C}$ such that $gf = hf$, we have $g = h$.
\end{itemize}
\end{definition}
We often use the abbreviations `mono' and `epi', omitting the suffix `-morphism'. 
It is not hard to prove that in any category, an iso is both mono and epi. The converse, however, \emph{fails} in many categories (see Exercises~\ref{exe:monoepiconcrete}~and~\ref{exe:monoepiiso} for examples).

In the category of sets, monomorphisms are just injective functions, and epimorphisms are surjective functions. \emph{However, this correspondence breaks down in most other concrete categories!} In particular, categories of algebras often have non-surjective epimorphisms; for an example (see Exercise~\ref{exe:monoepiconcrete} for the case of distributive lattices). An abstract notion that  corresponds more closely to `surjective', at least in categories of algebras, is that of a \emph{regular epimorphism}, which we introduce now.
\begin{definition}\label{def:coequalizer}
  Let $p_1, p_2 \colon C \to A$ and $q \colon A \to Q$ be morphisms in a category $\cat{C}$. We say that $q$ is a \emphind{coequalizer} of $(p_1,p_2)$ if (i) $qp_1 = qp_2$, and (ii) for any morphism $q' \colon A \to Q'$  in $\cat{C}$ such that $q'p_1 = q'p_2$, there exists a unique morphism $\alpha \colon Q \to Q'$ in $\cat{C}$ such that $\alpha q = q'$.

  A \emphind{regular epimorphism}\index{epimorphism!regular} $q \colon A \to Q$ in $\cat{C}$ is a morphism which is the coequalizer of some pair of morphisms.
\end{definition}
The terminology suggests that any morphism which arises as a coequalizer is in fact an epimorphism; you are asked to prove that this is indeed the case in Exercise~\ref{exe:coequalizers}. %
In the following example, we explicitly compute coequalizers in the category $\cat{Set}$.
\begin{example}\label{exa:coequalizer-set}
Let $p_1, p_2 \colon C \to A$ be functions between sets. Define $\equiv$ to be the smallest equivalence relation on the set $A$ that contains the pair $(p_1(c), p_2(c))$, for every $c \in C$. Let $Q := A/{\equiv}$ be the quotient, and $q \colon A \to Q$ the quotient map. We show that $q$ is a coequalizer of $(p_1, p_2)$. Clearly, $qp_1 = qp_2$, since $\equiv$ contains all pairs $(p_1(c),p_2(c))$. Suppose that $q' \colon A \to Q'$ is any function such that $q'p_1 = q'p_2$. Note that the equivalence relation $\ker(q') = \{(a,b) \in A^2 \mid q'(a) = q'(b)\}$ contains every pair $(p_1(c), p_2(c))$, and therefore must contain the equivalence relation $\equiv$, since this is the smallest such. Therefore, there is a well-defined and unique factorization $\alpha \colon Q \to Q'$ such that $\alpha q = q'$.
\end{example}
Exercise~\ref{exe:regularepiquotient} outlines a proof that in the category of
lattices, the regular epimorphisms coincide with the surjective homomorphisms.
Dually, we have notions of \emphind{equalizer} and \emphind{regular
monomorphism}, that we will discuss in more detail below. In
the category of lattices, a lattice homomorphism is a monomorphism if, and only if, it is
injective.  Not every monomorphism is regular in the category of distributive
lattices (see Exercise~\ref{exe:monoepiconcrete}). 

The definition of coequalizer is an example of a \emph{universal property}: it states the existence of a `minimal solution' for a configuration of morphisms in a category; in this case, a coequalizer $q$ of $(p_1,p_2)$ `solves for $x$' the equation $xp_1 = xp_2$, and does so in a `minimal' way, in the sense that any other solution $q'$ \emph{factors through} $q$. A coequalizer is a special case of a \emph{colimit} in a category, as we will explain below. 
Any concept defined by a universal property is \emph{unique up to unique isomorphism}. Concretely, for coequalizers this means that if $q_1$ and $q_2$ are both coequalizers of $(p_1,p_2)$, then the unique maps between the codomains of $q_1$ and $q_2$ which are guaranteed to exist by the definition must be inverse to each other (see Exercise~\ref{exe:coequalizers-unique}).

We will see several more examples of concepts defined by universal properties in this section. In fact, universal properties, limits and adjunctions are all closely related concepts; the interested reader may find, after digesting the definition of limit, an indication of the connection in Exercise~\ref{exe:limit-adjunction}.

\subsection*{Products and coproducts}
\begin{figure}
\begin{center}
\begin{tikzpicture}
  \matrix (m) [matrix of math nodes,row sep=3em,column sep=3em,minimum width=3em]
  { & A & \\
     & & \\
     & D_1 \times D_2 & \\
     D_1 & & D_2\\};
  \path[-stealth]
    (m-1-2) edge[dashed] node [right] {$\langle a_1, a_2 \rangle$} (m-3-2)
    (m-1-2) edge[bend right=30] node [above,xshift=-5pt] {$a_1$} (m-4-1)
    (m-1-2) edge[bend left=30] node [above,xshift=5pt] {$a_2$} (m-4-3)
    (m-3-2) edge node [right,xshift=5pt] {$\pi_1$} (m-4-1)
    (m-3-2) edge node [left,xshift=-5pt] {$\pi_2$} (m-4-3);
\end{tikzpicture}
\end{center}
\caption{The universal property of binary product in a category.}\label{fig:prod-diagram}
\end{figure}
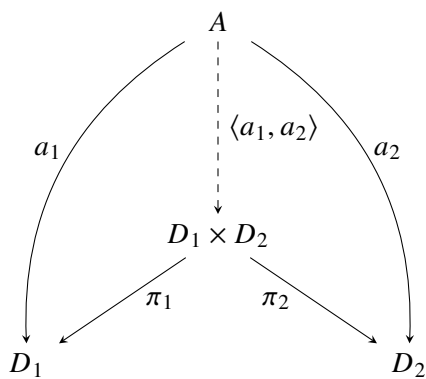

As we saw in Example~\ref{exa:coequalizer-set}, coequalizers are an appropriate categorical generalization of ``quotient map''. The Cartesian product and disjoint union of sets admit a similar generalization, namely as product and coproduct, respectively.
\begin{definition}\label{def:prod}
Let $I$ be a set, and $(D_i)_{i \in I}$ an $I$-indexed family of objects in a category $\cat{C}$. %
A \emph{product}\index{product!in categories} of the family $(D_i)_{i \in I}$ is an object $L$ together with an $I$-indexed family of morphisms $(\pi_i \colon L \to D_i)_{i \in I}$ such that, for any object $A$ and any $I$-indexed family of morphisms $(a_i \colon A \to D_i)_{i \in I}$, there exists a unique morphism $\alpha \colon A \to L$ such that $\pi_i \circ \alpha = a_i$ for every $i \in I$.

If a family has a product, we often denote it by $\prod_{i \in I}D_i$, and, if $I = \{1, 2\}$, simply by $D_1 \times D_2$. 
In the case $I = \emptyset$, a product of the unique $\emptyset$-indexed family is called a \emphind{terminal object}, and is often denoted by $1$. 
\end{definition}
Note that an object $T$ in a category $\cat{C}$ is terminal if, and only if, for every object $A$ in $\cat{C}$, there is a unique morphism from $A$ to $T$. 
When $a_1 \colon A \to D_1$ and $a_2 \colon A \to D_2$ are morphisms, the unique
morphism $A \to D_1 \times D_2$ is often denoted $\langle a_1, a_2 \rangle$.
The universal property can then be expressed in a diagram, as in Figure~\ref{fig:prod-diagram}.
Also, when $f_1 \colon A_1 \to D_1$ and $f_2 \colon A_2 \to D_2$ are morphisms,
then there is a unique morphism $f \colon A_1 \times A_2 \to D_1 \times D_2$
defined by $\langle f_1 \circ \pi_1, f_2 \circ \pi_2 \rangle$. This morphism $f$
is denoted $f_1 \times f_2$.

It is easy to verify that the Cartesian product of sets gives a product in $\cat{Set}$. The terminal object in $\cat{Set}$ is the one-element set.
Products exist in any category of finitary algebraic structures and are given by Cartesian product in those categories; in particular, the categorical product in categories of lattices is just given as in $\cat{Set}$ (see Exercise~\ref{exe:cart-product-cat}). The product in Priestley spaces is also given as in the category of ordered topological spaces, because the product of a collection of Priestley spaces is a Priestley space (see Exercise~\ref{exe:Priestley-product-cat}).

If a category has finite products, then this allows us to `internalize'
algebraic structure in the category; this perspective was popularized by Lawvere, giving rise
to a notion of Lawvere theory. While we do not need to develop this theory here, we give
an example that will be relevant in Chapter~\ref{ch:AutThry}, and gives the
flavor of the idea.
\begin{example}\label{exa:topmon}
	A \emph{monoid} internal to a category $\cat{C}$\index{monoid!internal
	to a category}\index{topological monoid} is an object $M$ of $\cat{C}$,
	together with two morphisms $I \colon 1 \to M$ and $\mu \colon M \times M
	\to M$ such that the following three diagrams commute:
	\begin{center}
\begin{tikzcd}
	1 \times M \arrow{r}{I \times \id_M} \arrow{rd}[swap]{\pi_2}  & M \times M
	\arrow{d}{\mu} \\
					     & M
\end{tikzcd}
\begin{tikzcd}
	M \times 1 \arrow{r}{\id_M \times I} \arrow{rd}[swap]{\pi_1}  & M \times M
	\arrow{d}{\mu} \\
					     & M
\end{tikzcd}

\begin{tikzcd}
	(M \times M) \times M \arrow{rr}{\mu \times \id_M} \arrow{d}[swap]{\alpha} & & M
	\times M \arrow{d}{\mu} \\
	M \times (M \times M) \arrow{r}{\id_M \times \mu} & M \times M
	\arrow{r}{\mu}
									  & M,
\end{tikzcd}
\end{center}
where, in the last diagram, the morphism $\alpha$ is the \emph{associator} defined by $\langle f_1, f_2
\rangle$, where $f_1 \colon (M \times M) \times M \to M$ is the composition of
the two projections on the first coordinate, and $f_2 \colon (M \times M) \times
M \to M \times M$ is $\pi_2 \times \id_{M}$.
These diagrams correspond to the monoid axioms $1 \cdot a = a$, $a \cdot 1 = a$,
and $(a \cdot b) \cdot c = a \cdot (b \cdot c)$, respectively.

For example, a \emph{topological monoid}, which we will encounter again in
Definition~\ref{def:topmon} is a monoid internal to $\cat{Top}$, and a
\emph{Priestley monoid} is a monoid internal to $\cat{Priestley}$. For these
concrete categories, it is easy to see that, for instance, a topological monoid
$(M, I, \mu)$ is the same thing as a monoid on the underlying sets for which
$\mu$ is continuous ($I$ always is because $1$ is a one-point space). An example
of a different flavor is a \emph{monad}, which can be seen as a monoid internal to the
category of endofunctors; however, properly defining monads in this way requires replacing in the above example the product $\times$ by the composition of endofunctors, known as \emph{monoidal structure} on a category. 
\end{example}

We now give the dual definition of \emph{coproduct}. 
\begin{definition}\label{def:coproduct}
  Let $I$ be a set, and $(D_i)_{i \in I}$ an $I$-indexed family of objects in a category $\cat{C}$. %
A \emph{coproduct}\index{coproduct!in categories} of the family $(D_i)_{i \in I}$ is an object $C$ together with an $I$-indexed family of morphisms $(\iota_i \colon D_i \to C)_{i \in I}$ such that, for any object $B$ and any $I$-indexed family of morphisms $(b_i \colon D_i \to B)_{i \in I}$, there exists a unique morphism $\beta \colon C \to B$ such that $\beta \circ \iota_i = b_i$ for every $i \in I$.

If a family has a coproduct, we often denote it by $\sum_{i \in I}D_i$, and, if $I = \{1, 2\}$, simply by $D_1 + D_2$. In the case $I = \emptyset$, a coproduct of the unique $\emptyset$-indexed family is called an \emphind{initial object}, and is often denoted by $0$. %
  \end{definition}
\begin{remark}
  We occasionally used the notation $A + B$ earlier in this book for the \emphind{symmetric difference} of two sets $A$ and $B$, usually as subsets of a common set $X$, as for example in Example~\ref{exa:remainderofbetaX}. Note that this symmetric difference $A + B$ is \emph{not} in general a coproduct of the objects $A$ and $B$ in the category  $\cat{Set}$, but it is when the sets $A$ and $B$ are disjoint.
\end{remark}
An initial object in a category is an object $0$ such that there exists a unique morphism from $0$ to any object of the category.
Coproducts in $\cat{Set}$ are given by disjoint unions and the initial object is the empty set. However, coproducts are often more complicated to compute in other categories of algebras, and spaces. For any set $V$, $\sum_{v \in V} \mathbf{3}$ is the free distributive lattice over the set of generators $V$, since $\mathbf{3}$ is the free distributive lattice on a single generator (see Exercise~\ref{exe:freeiscoproduct}). Finite coproducts in Priestley spaces are disjoint unions, but infinite coproducts involve a compactification step (see Exercise~\ref{exe:priestley-coproduct}). 

\begin{example}\label{exa:prod-as-meet}
Let $\cat{P}$ be a preorder, viewed as a category, and let $(p_i)_{i \in I}$ be a family of objects in $\cat{P}$. A product of $(p_i)_{i \in I}$ is exactly an \emph{infimum}, and  a coproduct  is exactly a \emph{supremum}. A terminal object in the category $\cat{P}$ is a greatest element, and an initial object is a least element. Note that the convention of drawing a product \emph{above} the objects (as in Figure~\ref{fig:prod-diagram}) clashes with the convention of Hasse diagrams (as for example in Figure~\ref{fig:diamond}), in which an object $q$ of $\cat{P}$ is depicted above an object $p$ if $p \leq q$.
\end{example}

As we have seen in Chapter 1, the value of infima and suprema depend on the ambient poset they are taken in. Thus, the example above illustrates clearly that products and coproducts depend not only on the objects in question, but also on the ambient category the product or coproduct is taken in. 
 
\subsection*{Limits and colimits}
In our definition of \emph{limit} and \emph{colimit} below (Definition~\ref{def:limits}), the following formal notion of \emph{diagram} plays a central role. We have been drawing diagrams informally throughout this book, but a diagram can itself be seen as a categorical object. We first give the formal definition.
\begin{definition}\label{dfn:cat-diagram}
  Let $\cat{C}$ be a category and $\cat{I}$ a small category. An \emph{$\cat{I}$-shaped diagram in $\cat{C}$} is a functor $D \colon \cat{I} \to \cat{C}$.
  \end{definition}

\begin{figure}[htp]
\begin{center}
\begin{tikzpicture}
  \matrix (m) [matrix of math nodes,row sep=4em,column sep=4em,minimum width=4em]
  {
     D1 & D2 \\
     D3 & D4\\};
  \path[-stealth]
    (m-1-1) edge node [above] {$Dx_{12}$} (m-1-2)
    (m-2-1) edge node [below] {$Dx_{34}$} (m-2-2)
    (m-1-1) edge node [above right] {$Dx_{14}$} (m-2-2)
    (m-1-1) edge node [left] {$Dx_{13}$} (m-2-1)
    (m-1-2) edge node [right] {$Dx_{24}$} (m-2-2);
  \end{tikzpicture}
\end{center}
\caption{An example of a commutative diagram, formally.} 
\label{fig:square-diag}
\end{figure}
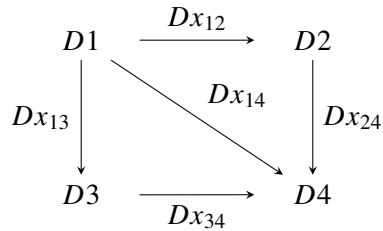

As a motivating example for this formal definition, which can look surprising at first, let us see how we would describe the following diagram (Figure~\ref{fig:square-diag}) according to Definition~\ref{dfn:cat-diagram}. 

Let $\cat{I}$ be the category defined by $\ob \cat{I}  := \{1, 2, 3, 4\}$, with the following five non-identity morphisms: $x_{12} \colon 1 \to 2$, $x_{13} \colon 1 \to 3$, $x_{24} \colon 2 \to 4$, $x_{34} \colon 3 \to 4$, and $x_{14}$, which is by definition equal to both $x_{24} \circ x_{12}$ and $x_{34} \circ x_{13}$. This indeed defines a category $\cat{I}$; associativity is trivial because there is no way to compose three non-identity morphism in this category. Now note that a functor $D \colon \cat{I} \to \cat{C}$ is essentially the same thing as a commutative square in $\cat{C}$ as depicted in Figure~\ref{fig:square-diag}, where we note that the value of $Dx_{14}$ is determined by functoriality of $D$ and the equalities $x_{24} \circ x_{12} = x_{14} = x_{34} \circ x_{13}$ holding in $\cat{I}$. Indeed, we have $Dx_{14} = D(x_{24} \circ x_{12}) = Dx_{24} \circ Dx_{12}$, and similarly $Dx_{14} = Dx_{34} \circ Dx_{13}$, showing that the diagram commutes. 

We are now ready to define limits and colimits. Two instances of this very general concept that can be helpful to keep in mind while reading this definition are the following:
\begin{itemize}
\item Let $I$ be a set, and let $\cat{I}$ be the small category with
$\ob{\cat{I}} = I$ and with only identity morphisms, no morphisms between
distinct objects -- this is called the \emphind{discrete category} on $I$. An
$\cat{I}$-shaped diagram in $\cat{C}$ is just given by an $I$-indexed family of objects in $\cat{C}$. A \emph{product}\index{product!special case of limit} as defined in
Definition~\ref{def:prod} above is a limit of an $\cat{I}$-shaped diagram, and a
coproduct is a \emph{colimit}\index{coproduct!special case of colimit} of an
$\cat{I}$-shaped diagram.  
\item Let $\cat{J}$ be the small category with two objects, $1$ and $2$, and two
morphisms from $1$ to $2$. A $\cat{J}$-shaped diagram in a category $\cat{C}$ is
a pair of objects in $\cat{C}$ with a parallel pair of morphisms between them. A
colimit of a $\cat{J}$-shaped diagram is a coequalizer as defined in
Definition~\ref{def:coequalizer}, and a limit of such a diagram is called an 
\emphind{equalizer}.
Spelling out the definition of the latter, an \emph{equalizer} of a pair of parallel arrows $f_1, f_2 \colon A \rightrightarrows B$ in a category $\cat{C}$ is an object $E$ and a morphism $e \colon E \to A$ such that $f_1e = f_2e$ and, for any morphism $a \colon C \to A$ such that $f_1a = f_2a$, there exists a unique $\alpha \colon C \to E$ such that $a = e\alpha$.
\end{itemize}
We now give the general definition of limit and colimit of a diagram in a category.
\begin{definition}\label{def:limits}
  Let $\cat{C}$ be a category, let $\cat{I}$ be a small category, and let $D \colon \cat{I} \to \cat{C}$ be an $\cat{I}$-shaped diagram in $\cat{C}$.

  A \emph{cone above} \index{cone} the diagram $D$ is an object $A$ together with an $\ob \cat{I}$-indexed family of morphisms in $\cat{C}$,  $(a_i \colon A \to Di)_{i \in \ob \cat{I}}$, such that, for any morphism $x \colon i \to j$ in $\cat{I}$, we have $(Dx) \circ a_i = a_j$.
  A \emph{limit} \index{limit!in categories} of the diagram $D$ is a universal cone above $D$. That is, a limit of $D$ is an object $L$ with an $\ob \cat{I}$-indexed family of morphisms, $(\pi_i \colon L \to Di)_{i \in \ob \cat{I}}$, such that, for any cone $(a_i \colon A \to Di)_{i \in \ob \cat{I}}$ above $D$, there exists a unique morphism $\alpha \colon A \to L$ in $\cat{C}$ such that $\pi_i \alpha = a_i$ for every $i \in \ob \cat{I}$.  

  Similarly, a \emph{cone below}, or also \emph{co-cone under},\index{co-cone} $D$ is an object $B$ together with morphisms $b_i \colon Di \to B$ for every $i \in \ob \cat{I}$ such that, for any $x \colon i \to j$ in $\cat{I}$, $b_j \circ Dx = b_i$. A \emph{colimit} \index{colimit} is a universal cone below $D$, that is, an object $C$ together with morphisms $\iota_i \colon Di \to C$ for every $i \in \ob \cat{I}$ such that, for any cone $(b_i \colon Di \to B)_{i \in \ob \cat{I}}$ below $D$, there exists a unique $\beta \colon C \to B$ such that $\beta \iota_i = b_i$ for every $i \in \ob \cat{I}$.
\end{definition}
We note here that a cone above a diagram $D$ can be viewed more conceptually as a natural transformation from a constant functor $\Delta$ to $D$, and a cone below $D$ as a natural transformation from $D$ to a constant functor (see Exercise~\ref{exe:cones-constant} for more details). Using this observation, one may see a limit of $D$ as a terminal object in a category of cones over $D$. 

As we allowed the shape of the indexing category of the diagram to be \emph{any}
small category in Definition~\ref{dfn:cat-diagram}, the reader may reasonably
worry that limits and colimits could have incredibly complicated shapes. While
this is true, the relevant shapes that we need to consider in practice are
rather limited, thanks to a few bits of general category theory that we will develop in Proposition~\ref{prop:complete-category} below.
In particular, we will explain why it is
enough for our purposes to understand the following three types of limits/colimits:
\begin{itemize}
\item product and coproduct;
\item equalizer and coequalizer;
\item projective limit and directed colimit.
\end{itemize}
\subsection*{Projective limits and directed colimits}
The definitions of a projective limit and directed colimit
use the notion of directed poset: recall that a poset $P$ is \emph{directed} (sometimes called \emph{up-directed})
if it is non-empty, and any pair of elements $p, q$ has a common upper bound,
and \emph{down-directed} if $P^\op$ is up-directed (see
p.~\pageref{def:directed}).  
\begin{definition}
A diagram $D \colon \cat{I} \to \cat{C}$  is \emph{directed} if $\cat{I}$ is the category associated to an up-directed partial order, and \emph{projective} if $\cat{I}$ is the category associated to a down-directed partial order. 
A \emph{projective limit} is a limit of a projective diagram; a \emph{directed colimit} is a colimit of a directed diagram.\footnote{In algebraic literature, what we call `directed colimit'
is sometimes called `direct limit' or `inductive limit'. We avoid this
terminology because we prefer to emphasize `co-' in the name of this concept, to
signal the fact that this generalizes a `co-product'. On the other hand, in
category theory, what we call `projective limit' is sometimes called
`down-directed limit', or `codirected limit'.} 
\end{definition}
Some sources are more liberal in their definitions of directed colimit and
projective limit, by allowing the shape of the diagram, \cat{I}, to be a
directed \emph{pre}order rather than partial order. However, this does not give
a more general notion of limit, up to isomorphism (see
Exercise~\ref{exe:preorder-or-poset-colimit}). A seemingly even more general definition is that of a colimit with respect to a \emph{filtered} category, that is, a non-empty category in which any pair of objects admits a morphism to a common object, and any parallel pair of morphisms has a cone under it. Surprisingly, colimits with respect to filtered categories are no more general than  colimits with respect to directed posets; the proof of this fact is beyond the scope of this book, see for example \cite[Theorem 1.5]{AdaRos1994}.

In a concrete category, projective limits and directed colimits of \emph{finite}
objects are often particularly interesting to look at, as they tend to inherit
some of the properties of finite objects. As one example, in categories of
algebraic structures, \emph{every} algebra is a directed colimit of finitely
presented algebras. This is essentially the defining property of what is known
as an $\aleph_0$-accessible category in the literature, see for example
\cite{AdaRos1994}. 

In the particular case of distributive lattices, the situation simplifies further, since any finitely presented distributive lattice is in fact finite, as follows immediately from the fact that free distributive lattices are finite. Thus, every distributive lattice is a directed colimit of \emph{finite} distributive lattices. We now sketch a direct proof of the latter fact, which should give an impression of the basics of the much more general theory. 

\begin{example}\label{exa:DLasIndDLfin} Let $L$ be a distributive lattice. Let
	$\Sub_f(L)$ denote the poset of finite sublattices of $L$, which is
	up-directed, because if $F$ and $G$ are finite sublattices of $L$, then
	they are both contained in the sublattice generated by $F \cup G$, which
	is still finite. Let $D : \Sub_f(L) \to \cat{DL}$ be the diagram defined
	on objects by $D(F) := F$, and sending any sublattice inclusion $F
	\subseteq G$ to the homomorphism $D(F \subseteq G) := i_{F,G} \colon F
	\into G$. Note that a co-cone under this diagram $D$ is given by the
	inclusion maps $\iota_F \colon F \to L$.  Indeed, this is a colimit: if
	$(j_F \colon F \to M)_{F \in \Sub_f(L)}$ is another co-cone under this
	diagram, then a unique $\alpha \colon L \to M$ may be defined by sending
	$a \in L$ to the element $j_{L_a}(a)$, where $L_a$ denotes the (at most
	three-element) sublattice of $L$ generated by $a$. To prove that
	$\alpha$ is indeed a homomorphism, one may use the fact that any pair of
	elements generates a finite sublattice of $L$, and that the $j_F$ form a
	co-cone. 

Generalizing this example slightly, let $P$ be a directed subset of $\Sub_f(L)$.
One may then prove that the \emph{union} of the finite sublattices in $P$ is a
sublattice $L'$ of $L$, and that $L'$ is the directed colimit of the restricted
diagram $D|_P$.  The directed colimit of $P$ is (isomorphic to)
$L$ itself if, and only if, every $a \in L$ lies in some sublattice in the
collection $P$. This characterization of directed colimits will be used in
Section~\ref{sec:bifinite} to give a concrete definition of \emph{bifinite}
distributive lattices.
\end{example}

For an example of the importance of projective limits of finite objects, it is
interesting to move to the topological side of Priestley duality. Indeed,
Priestley spaces may be characterized as the \emphind{profinite posets}, in a
sense that we make precise now in the following example.
\begin{example}\label{exa:priestley-as-profinite posets}
In the category $\cat{TopOrd}$ of ordered topological spaces with continuous order-preserving functions, we have a full
subcategory isomorphic to the category of finite posets. Indeed, for every
finite poset $P$, we have an ordered topological space, that we still denote by
$P$, by equipping $P$ with the discrete topology. With this topology, $P$ is a
finite Priestley space, and any order-preserving function between finite posets is
continuous for the discrete topology, as we saw in Exercise~\ref{exe:fullsubcats} above. We denote this category of finite posets with the discrete topology by $\cat{TopOrd}_{f,d}$. We now use Example~\ref{exa:DLasIndDLfin} and Priestley
duality to show that every Priestley space is a projective limit of finite
Priestley spaces. 

Let $X$ be a Priestley space and let $L$ be its dual
distributive lattice.  Since $L$ is the colimit of the diagram $D \colon
\Sub_f(L) \to \cat{DL}$ given in Example~\ref{exa:DLasIndDLfin}, it follows by a
purely formal argument, since $\cat{DL}^\op$ is equivalent to $\cat{Priestley}$,
that $X$ is the limit of the diagram $D' \colon (\Sub_f(L))^\op \to
\cat{Priestley}$ that sends a finite sublattice $F$ of $L$ to the Priestley dual
space of $F$. We may give a more direct description of the diagram $D'$ using the
sublattice-quotient duality of Section~\ref{sec:quotients-and-subs}: the domain
of $D'$ is (isomorphic to) the poset of compatible preorders $\preceq$ on $X$
for which $X/{\preceq}$ is finite, and it sends such a pre-order $\preceq$ to
the finite Priestley space $X/{\preceq}$.  Since finite Priestley spaces always
have discrete topology, an equivalent way to describe $D'$ is as the diagram of
all continuous surjective maps $X \onto P$ where $P$ is a finite poset, and $X$
is the limit of this diagram in the category $\cat{TopOrd}$. Also, since any limit of
finite posets is a Priestley space, we conclude that Priestley spaces exactly coincide
with the projective limits of finite posets in $\cat{TopOrd}$. 

As a special case,
Boolean spaces exactly coincide with the projective limits of finite sets with the
discrete topology in $\cat{Top}$, and are therefore sometimes called \emphind{profinite sets}.
The ``profinite'' point of view on Stone-Priestley duality will play an important role
in Chapter~\ref{ch:AutThry}, where it is linked to the theory of \emph{profinite monoids}.
The categorical definition of profinite monoids is analogous to the ones in this example: 
the category of finite monoids is isomorphic to a full subcategory of
the category $\cat{TopMon}$ topological monoids (Example~\ref{exa:topmon}) by
equipping each finite monoid with the discrete topology.
A \emphind{profinite monoid} is then an object of $\cat{TopMon}$ that is a
projective limit of finite monoids. Further see
Definition~\ref{def:topmon}.
\end{example}

\begin{remark}\label{rem:priestley-as-procompletion}
We remark that the arguments given in Example~\ref{exa:priestley-as-profinite posets} can be turned around to give a
categorical proof of the Priestley duality theorem itself from Birkhoff duality
for \emph{finite} distributive lattices, as we now outline, without giving all
the details. Example~\ref{exa:DLasIndDLfin} shows, in categorical terms, that
the category $\cat{DL}$ may be obtained by taking directed colimits of objects
from the category $\cat{DL}_f$. Therefore, its (formal) dual category
$\cat{DL}^\op$ can be described by taking projective limits of the category
$(\cat{DL}_f)^\op$. The latter category is equivalent to the category of finite
posets, by Birkhoff duality, and this is isomorphic to the to the full
subcategory $\cat{TopOrd}_{f,d}$ of $\cat{TopOrd}$ on the finite posets with the
discrete topology. It follows that $\cat{DL}^\op$ is equivalent to the closure
of $\cat{TopOrd}_{f,d}$ under projective limits. This gives a `formal' version of Priestley duality. 

Note that we have not used any choice
principles to deduce this duality, only Birkhoff duality and abstract category theory. 
In this point of view on Priestley duality, the choice principle needs to
be used only in the very last step, namely, to prove that the closure under
projective limits of $\cat{TopOrd}_{f,d}$
is indeed the same as the category of Priestley spaces. This requires in
particular showing that any product of finite spaces is compact, which is a
`finite' Tychonoff theorem, and is in fact equivalent over choice-free set
theory to the Stone prime ideal theorem. 
\end{remark}

To understand the importance of products, equalizers, and projective limits for understanding general limits, we now 
first look at a special, \emph{finite} case, where just binary products and equalizers will suffice. A generalization of this result, also involving projective limits, will be proved after that, in Proposition~\ref{prop:complete-category}. 

Let $\cat{C}$ be a category,
and consider two morphisms in $\cat{C}$, $f_1 \colon A_1 \to B$ and $f_2 \colon
A_2 \to B$, with a common codomain $B$, as in the following diagram.
%
\begin{center}
\begin{tikzpicture}
  \matrix (m) [matrix of math nodes,row sep=4em,column sep=4em,minimum width=4em]
  {
     & A_1 \\
     A_2 & B\\};
  \path[-stealth]
    (m-2-1) edge node [below] {$f_2$} (m-2-2)
    (m-1-2) edge node [right] {$f_1$} (m-2-2);
  \end{tikzpicture}
\end{center}
%
Such a pair of morphisms is called a \emphind{co-span}
in category theory, and a limit of such a diagram is called a
\emphind{pullback}; the dual notions are called \emphind{span} and
\emphind{pushout}, respectively. 

Now assume that binary products and equalizers exist in
$\cat{C}$. We show that the limit (pullback) of any such co-span can be constructed
in $\cat{C}$, as follows. First let $P = A_1 \times A_2$ be the product of $A_1$
and $A_2$ in $\cat{C}$, and consider two morphisms $g_1, g_2 \colon P
\rightrightarrows B$ defined by $g_i := f_i \pi_i$ for $i = 1, 2$, where $\pi_i
\colon P \to A_i$ is the projection morphism. Now let $e \colon E \to P$ be the
equalizer of $g_1$ and $g_2$. We may now prove that the pair of maps $p_1, p_2
\colon E \to A_i$, defined by $p_i := \pi_i e$ for $i = 1, 2$, is a limit cone for
the co-span. Indeed, using the above definitions, we have that
\[ f_1 p_1 = f_1 \pi_1 e = g_1 e = g_2 e = f_2 \pi_2 e = f_2 p_2,\]
so that the $p_i$ form a cone over the diagram. Furthermore, if $a_1 \colon A \to A_1$ and $a_2 \colon A \to A_2$ are another cone over the same diagram, then, writing $a$ for the morphism $\langle a_1, a_2 \rangle \colon A \to A_1 \times A_2$, we get
\[ g_1 a = f_1 \pi_1 a = f_1 a_1 = f_2 a_2 = f_2 \pi_2 a = g_2 a,\]
using the definitions and the fact that the $a_i$ form a cone. Thus, there exists a unique $\alpha \colon A \to E$ such that $e \alpha = a$. It follows that $a_i = \pi_i a = \pi_i e \alpha = p_i \alpha$ for $i = 1, 2$. Also, such $\alpha$ is unique: if $\alpha'$ is such that $a_1 = p_1 \alpha'$ and $a_2 = p_2\alpha'$, then $\pi_i e \alpha' = p_i \alpha' = a_i$ for $i = 1, 2$, so that $e \alpha' = e \alpha$ by the universal property of the product $P$. Therefore, $\alpha' = \alpha$ since equalizers are monomorphisms (see Exercise~\ref{exe:coequalizers}). 

A generalization of the above construction allows one to show that, if all equalizers and all set-sized products exist in a category, then the limit of any diagram exists. Here, recall that the shape of a diagram is by definition a \emph{small} category. We prove this fact now, and further connect it to the construction of \emph{projective} limits in a category.
\begin{proposition}
  \label{prop:complete-category}
Let $\cat{C}$ be a category. The following are equivalent:
\begin{enumerate}
  \item[(i)] Any diagram (of small shape) has a limit in $\cat{C}$. 
  \item[(ii)] Any finite diagram has a limit in $\cat{C}$, and any projective diagram has a limit in $\cat{C}$.
  \item[(iii)] All equalizers and small products exist in $\cat{C}$.
\end{enumerate}
\end{proposition}
\begin{proof}
It is trivial that the first item implies the second. 

For the implication (ii)
$\Rightarrow$ (iii), since an equalizer is the limit of a finite diagram, it
suffices to prove that any set-sized product exists in $\cat{C}$. Let $(A_i)_{i
\in I}$ be a family of objects of $\cat{C}$ indexed by a set $I$. We will show
that a product of the $A_i$ can be constructed as a limit of a diagram, $B$, of
shape $\mathcal{P}_f(I)^\op$, the poset of finite subsets of $I$, ordered by
reverse inclusion -- this is clearly a down-directed poset, since it is a
meet-semilattice, so the limit of $B$ will exist by assumption. For each $F
\in \mathcal{P}_f(I)$, let $B(F)$ be a product of the finite set of objects
$\{A_i \mid i \in F\}$, which exists in $\cat{C}$ by assumption, and, for $i
\in F$, denote by $\pi^F_i \colon B(F) \to A_i$ the projection of this cone on
$A_i$. In the special case where $F$ is a singleton, we choose $B(\{i\}) = A_i$, which is the product of the single object $A_i$. 
If $F \subseteq G$, then we have a cone with summit $B(G)$ over $\{A_i
\mid i \in F\}$, for each $i \in F$, we have $i \in G$, so we have $\pi^G_i
\colon B(G) \to A_i$. Let $\beta_{F,G} \colon B(G) \to B(F)$ be the unique
factorization of this cone, so that for every $i \in F$, $\pi^F_i \beta_{F,G} =
\pi^G_i$, and define $B(G \supseteq F) := \beta_{F,G}$. Now let $L$ be a limit
of the diagram $B$ with $\lambda_F \colon L \to B(F)$ the legs of the limit
cone, for each finite $F \subseteq I$. In particular, for every $i \in I$, write
$\pi_i$ for the leg $\lambda_{\{i\}} \colon L \to B(\{i\}) = A_i$. We claim that
$\pi_i \colon L \to A_i$ gives a product of the $A_i$. Indeed, if $f_i \colon A
\to A_i$ is any family of morphisms, then we may construct, for each finite
subset $F \subseteq I$, a unique morphism $f_F \colon A \to B(F)$ such that
$\pi^F_i f_F = f_i$ for each $i \in F$, using the universal property of the
product $B(F)$. The morphisms $f_F$ give a cone over the diagram $B$, because
for any morphism $G \supseteq F$ in $\mathcal{P}_f(I)^\op$, we have $\pi^F_i
\beta_{F,G} f_G = \pi^G_i f_G = f_i = \pi^F_i f_F$, so that $\beta_{F,G} f_G =
f_F$, using the uniqueness part of the universal property of the product $B(F)$.
The unique factorization $\alpha \colon A \to L$ satisfies $\pi_i
\alpha = f_i$, as required.

Finally, we sketch a proof that (iii) implies (i); we emphasize that this is an
infinite generalization of the construction of pullbacks from binary products
and equalizers, described just above this proposition. Let $D \colon \cat{I} \to
\cat{C}$ be a diagram indexed by a small category $\cat{I}$. Observe that giving
a cone $a_i \colon A \to Di$ over the diagram is the same thing as giving a
single morphism $a \colon A \to \prod_{i \in \ob(I)} Di$ with the property that,
for all morphisms $f \colon i \to j$ in $\cat{I}$, $A(f) \circ a_i = a_j$. The
latter property, in turn, can be more synthetically expressed as follows. Write
$P$ for the product $\prod_{i \in \ob \cat{I}} Di$, with projection maps $\pi_i
\colon P \to Di$ for every object $i$ of $\cat{I}$, and write $M$ for the
product $\prod_{f \in \mor \cat{I}} D(\cod(f))$, with projection maps $\rho_f
\colon M \to D(\cod(f))$ for every morphism $f$ of $\cat{I}$.  Now the crucial
observation is that $a \colon A \to P$ is a cone exactly if $d \circ a = e \circ
a$, where $d, e \colon P \rightrightarrows Q$ are two morphisms of $\cat{C}$,
defined, respectively, by requiring that for every morphism $f$ of $\cat{I}$:
\[ \rho_f \circ d := \pi_{\cod(f)}, \quad \rho_f \circ  e := D(f) \circ \pi_{\dom(f)}.\]
From this observation, it follows that the equalizer of the parallel pair of
maps $d, e$ is a limit cone of the diagram $D$. For more information, we refer to, for instance, 
\cite[Theorem~12.3]{AdaHerStr1990} or \cite[Proposition~5.1.26]{Lei2014}.  
\end{proof}
A category is called \emphind{complete} if it satisfies the equivalent
properties in Proposition~\ref{prop:complete-category}, and
\emphind{co-complete} if $\cat{C}^\op$ is complete. Observe that the dual of
Proposition~\ref{prop:complete-category} says that a category has all small
colimits if, and only if, all co-equalizers and set-sized coproducts exist, which
happens if, and only if, all finite
colimits and directed colimits exist.  Note also that, if we specialize to the
case of a thin category, then the equivalence of (i) and (ii) in
Proposition~\ref{prop:complete-category} (and its dual version) reduce to a
familiar characterization of complete lattices: a poset is a complete lattice if
and only if all finite sets and all up-directed sets have a supremum. This fact
was already stated in Exercise~\ref{exe:directed} and also plays an important
role in the background in domain theory (Chapter~\ref{chap:DomThry}).

\ourexercises
\begin{ourexercise}\label{exe:monoepiconcrete}
Let $f \colon L \to M$ be a homomorphism between distributive lattices.
\begin{enumerate}
\item Prove that $f$ is a monomorphism if, and only if, $f$ is injective.
\item Prove that if $f$ is surjective, then $f$ is an epimorphism.
\item Consider an order-embedding $e$ of the three-element chain $\mathbf{3}$ into the four-element Boolean algebra, $\mathbf{2} \times \mathbf{2}$, which is clearly not surjective. Prove that $e$ is an epimorphism in the category of distributive lattices.
\item Conclude that in $\cat{DL}$, there exist morphisms which are both mono and epi, but not iso.
\item Also prove that $e$ is a non-regular monomorphism in $\cat{DL}$.
\end{enumerate}
\end{ourexercise}
\begin{ourexercise}
\label{exe:spacesepi}
This exercise examines epimorphisms in categories of topological spaces; note in particular that the notion of ``epimorphism'' is \emph{not} stable under taking full subcategories.
\begin{enumerate}
  \item Prove that in the category $\cat{Top}$, any epimorphism is surjective.
  \item Prove that the epimorphisms in the category of Hausdorff topological spaces are exactly the continuous functions whose image is dense in the codomain. 
  \item Deduce from the previous item that in $\cat{Priestley}$, the epimorphisms are exactly the surjective continuous functions. 
\end{enumerate}
\end{ourexercise}
\begin{ourexercise}
\label{exe:booleanepi}
This exercise characterizes epimorphisms in the categories $\cat{DL}$ and $\cat{BA}$.
\begin{enumerate}
\item Prove that in $\cat{BA}$, the epimorphisms are exactly the surjective homomorphisms. \hint{One direction always holds. For the other, use Stone duality and the fact that the dual of a surjective homomorphism is injective.}
\item Let $L$ be a distributive lattice and let $e_L \colon L \to L^-$ be its Boolean envelope (Definition~\ref{def:booleanenvelope}). Prove that $e_L$ is an epimorphism in $\cat{DL}$.
\item Prove that in $\cat{DL}$, a homomorphism $f \colon L \to M$ is an epimorphism if, and only if, the unique extension $\overline{e_M \circ f} \colon L^- \to M^-$ of $e_M \circ f \colon L \to M^-$ is surjective.
\end{enumerate}
\end{ourexercise}

\begin{ourexercise}\label{exe:monoepiiso}
Let $f \colon A \to B$ be a morphism in a category $\cat{C}$.
\begin{enumerate}
\item Prove that, if $f$ is an isomorphism, then $f$ is both a monomorphism and an epimorphism.
\item Prove that, if $f$ has a left inverse (that is, there exists $g \colon B \to A$ such that $gf = 1_A$), then $f$ is a monomorphism.
\item Formulate the dual statement of (b) and conclude that it is true.
\item Consider the category with one object (that is, a monoid), a countable set of morphisms $\{f_n \ : \ n \in \mathbb{N}\}$, composition given by $f_n \circ f_m := f_{n + m}$, and identity $f_0$. Prove that, for every $n \geq 1$, $f_n$ is mono and epi, but not iso.
\end{enumerate}
\end{ourexercise}

\begin{ourexercise}\label{exe:coequalizers}
Prove that any equalizer is a monomorphism, and any coequalizer is an epimorphism.\index{coequalizer!is an epimorphism}\index{equalizer!is a monomorphism}
\end{ourexercise}
\begin{ourexercise}\label{exe:coequalizers-unique}
Suppose that $q_1 \colon A \to Q_1$ and $q_2 \colon A \to Q_2$ are both coequalizers for the same pair of maps $p_1, p_2 \colon C \to A$.\
\begin{enumerate}
\item Show that there exist morphisms $\alpha_1 \colon Q_1 \to Q_2$ and $\alpha_2 \colon Q_2 \to Q_1$ such that $\alpha_1q_1 = q_2$ and $\alpha_2 q_2 = q_1$.
\item Prove that $\alpha_2\alpha_1 = \id_{Q_1}$ and $\alpha_1\alpha_2 = \id_{Q_2}$. \hint{Show first that $\alpha_2\alpha_1$ is a factorization of $q_1$ through itself, that is, that $\alpha_2\alpha_1q_1 = q_1$, and then use the uniqueness of $\alpha$ in the definition of coequalizers.}
\end{enumerate}

\end{ourexercise}

\begin{ourexercise}\label{exe:regularepiquotient}
This exercise outlines a proof that the regular epimorphisms in  $\cat{Lat}$ are exactly the surjective homomorphisms. To readers with more knowledge of universal algebra, it should be clear that the same argument applies in any category whose objects are the algebras in a variety and whose morphisms are the homomorphisms.
\begin{enumerate}
\item Prove that an equivalence relation $\theta$ on a lattice $A$ is a congruence if, and only if, $\theta$ is a (bounded) sublattice of the product algebra $A \times A$.
\item Let $p_1, p_2 \colon C \to A$ be a pair of lattice homomorphisms. Define $\theta$ to be the smallest congruence on $A$ which contains $(p_1(c),p_2(c))$ for every $c 
\in C$. Prove that the quotient $q \colon A \onto A/{\theta}$ is a coequalizer of $(p_1,p_2)$. Conclude that regular epimorphisms are surjective.
\item Let $\theta$ be a congruence on a lattice $L$. Prove that the quotient homomorphism $q \colon A \onto A/{\theta}$ is the coequalizer of the pair of maps $p_1, p_2 \colon \theta \to A$, where $p_i$ is the projection onto the $i^\mathrm{th}$ coordinate.
\item Conclude that a lattice homomorphism is a regular epimorphism if, and only if, it is surjective.
\end{enumerate}
\end{ourexercise}

\begin{ourexercise}\label{exe:cones-constant}
  Let $\cat{C}$ be a category and let $\cat{I}$ be a small category.
  \begin{enumerate}
  \item Show that, for any object $C \in \cat{C}$, the functor $\Delta_C$, defined by $\Delta_Ci = C$ for every object and $\Delta_Cx = 1_C$ for every morphism, is a diagram, called the $\cat{I}$-shaped \emphind{constant diagram}\index{diagram!constant} at $C$.
  \item Let $D$ be a diagram. Show that a cone above $D$ is essentially the same thing as a natural transformation $p \colon \Delta_C \Rightarrow D$, for some object $C$ of $\cat{C}$.
  \item Reformulate the definition of \emph{limit} purely in terms of natural transformations, without ever referring to specific objects and morphisms in $\cat{I}$ and $\cat{C}$.
  \item Dualize this exercise for \emph{cones below} and \emph{colimits}.
  \end{enumerate}
\end{ourexercise}

\begin{ourexercise}\label{exe:cart-product-cat}
Prove that the Cartesian product of posets/lattices/distributive lattices is a product in the categorical sense of Definition~\ref{def:prod}.
\end{ourexercise}
\begin{ourexercise}\label{exe:Priestley-product-cat}
Recall from Exercise~\ref{exe:Priestley-product-space} that the Cartesian product of a family of Priestley spaces, defined by
equipping the product set with the product topology and the point-wise partial order, is again a Priestley space. Conclude that products exist in the category of Priestley spaces.
\end{ourexercise}

\begin{ourexercise}\label{exe:priestley-coproduct}
Let $(X_i)_{i \in I}$ be a collection of Priestley spaces. For each $i \in I$, write $L_i$ for the lattice of clopen down-sets of $X_i$. Prove that the Priestley dual space $X$ of the product $\prod_{i \in I} L_i$ is the coproduct $\sum_{i \in I} X_i$ in Priestley spaces. Show that there is an injective function from the disjoint union of the sets $X_i$ to $X$, which is in general not surjective.
\end{ourexercise}
\begin{ourexercise}\label{exe:left-ad-pres-colimits}
Let $F \colon \cat{C} \to \cat{D}$ be a functor, and suppose that $F$ has a right adjoint. Prove that $F$ preserves any colimit that exists in $\cat{C}$, that is, if $D \colon \cat{I} \to \cat{C}$ is a diagram in $\cat{C}$ and $C$ is a colimit of $D$ with a co-cone given by $\lambda_i \colon D(i) \to C$ for each object $i$ of $\cat{I}$, then the image of this co-cone under $F$, $F\lambda_i \colon FD(i) \to FC$, gives a colimit of the diagram $F \circ D \colon \cat{I} \to \cat{D}$. 
\end{ourexercise}
\begin{ourexercise}\label{exe:freeiscoproduct}
  This exercise elaborates on the comment made in Example~\ref{exa:twogen-logic} to the effect that the two element poset $(2,\preceq)$ appearing in the concrete incarnation of $F_{\cat{DL}}(V)$ derived in Section~\ref{sec:free-description} is the poset dual to $\bf 3$, the three element lattice.
  \begin{enumerate}
  \item Prove that $\bf 3$ is a free distributive lattice over $\{p\}$, by exhibiting a universal arrow.
  \item Prove that for any set $V$ we have $F_{\cat{DL}}(V)=\sum_{v\in V} F_{\cat{DL}}(\{v\})$.
  \item Conclude, via duality, that the Priestley dual of $F_{\cat{DL}}(V)$ is the product space $(2,\preceq,\delta)^V$, where $\delta$ is the discrete topology on $2$, see Proposition~\ref{prop:freeDL-duality}.
    \end{enumerate}
  \end{ourexercise}

\begin{ourexercise}\label{exe:limit-adjunction}
  Let $\cat{C}$ be a category and let $\cat{I}$ be a small category. Denote by $[\cat{I}, \cat{C}]$ the category of $\cat{I}$-shaped diagrams in $\cat{C}$, with natural transformations between them.
  \begin{enumerate}
  \item Show that there is a functor $\Delta \colon \cat{C} \to [\cat{I}, \cat{C}]$ which sends any object $C$ of $\cat{C}$ to the constant diagram $\Delta_C$ (see Exercise~\ref{exe:cones-constant}), and a morphism $f \colon C \to C'$ to the natural transformation $\Delta C \to \Delta C'$ that is $f$ at every component.
    
  \item Suppose that $\Delta$ has a left adjoint, $L \colon [\cat{I}, \cat{C}] \to \cat{C}$. Prove that, for any diagram $D \in [\cat{I}, \cat{C}]$, $L(D)$ is a colimit of $D$.
  \item Similarly, prove that if $R$ is a right adjoint for $\Delta$, then $R(D)$ is a limit of $D$, for any $\cat{I}$-shaped diagram $D$ in $\cat{C}$.
  \end{enumerate}
\end{ourexercise}
\begin{ourexercise}\label{exe:preorder-or-poset-colimit}
Let $I$ be a preordered set that is directed, that is, any finite subset of $I$ has an upper bound. Let $\pi \colon I \to J$ be the poset reflection of $I$ (see Exercise~\ref{exe:reflection}). Let $D \colon \cat{I} \to \cat{C}$ be a diagram.
\begin{enumerate}
  \item Prove that $J$ is directed.
  \item Prove that there is a diagram $D' \colon \cat{J} \to \cat{C}$ such that $D'(\pi(i))$ is isomorphic in $\cat{C}$ to $D(i)$ for every $i \in I$.
  \item Suppose that $D$ has a colimit in $\cat{C}$. Prove that the same object gives a colimit of $D'$ in $\cat{C}$.
\end{enumerate}
\end{ourexercise}

\section{Priestley duality categorically}\label{sec:duality-categorically}
In this section we apply the categorical point of view introduced in this chapter to Priestley duality. While the general language of category theory is useful for expressing our duality theorem, we will also require some knowledge that seems to be rather specific to Priestley duality.
In particular, in this section we further elaborate on what the correspondences between monomorphisms and epimorphisms on either side of the duality concretely look like in the case of Priestley duality. This yields connections between quotients and sublattices, and between subspaces and quotient lattices, which are crucial to some of the applications in the later chapters of this book. Further, we show how Priestley duality can be used to easily describe free distributive lattices.

\subsection*{Priestley's dual 2-equivalence}
Recall from Chapter~\ref{ch:priestley} that, to any Priestley space $X$, we associated a distributive lattice $\ClD(X)$ of clopen down-sets of $X$, and, to any continuous order-preserving map $f \colon X \to Y$, the distributive lattice homomorphism $f^{-1} \colon \ClD(Y) \to \ClD(X)$. In the reverse direction, if $L$ is a distributive lattice, we associate to it the Priestley space $X_L$ (Definition~\ref{def:Priestleydualspace}), and if $h \colon M \to L$ is a distributive lattice homomorphism, then Proposition~\ref{prop:priestleymorphisms} showed that there is a unique continuous map $f \colon X_L \to X_M$ with the property that $\widehat{h(a)} = f^{-1}(\widehat{a})$ for every $a \in M$. In this section, we will write $X(h)$ for this unique map $f$.

With the categorical language developed in this chapter, we can now properly state and prove Priestley's duality theorem. 

\begin{theorem}\label{thm:priestley-duality-long}
The assignments $L \mapsto X_L$ and $(h \colon M \to L) \mapsto (X(h) \colon X_L \to X_M)$ constitute a contravariant functor, $X \colon \DL \to \Priest$, and the pair $(\ClD, X)$ is a contravariant equivalence of categories. 
\end{theorem}
\begin{proof}
  To see that $X$ is a functor, let $k \colon N \to M$ and $h \colon M \to L$ be morphisms in $\DL$. Note that $X(k) \circ X(h)$ is a Priestley morphism for which it holds that, for any $a \in N$, $(X(k) \circ X(h))^{-1}(\widehat{a}) = X(h)^{-1}(X(k)^{-1}(\widehat{a})) = X(h)^{-1}(\widehat{k(a)}) = \widehat{h(k(a))}$. By the uniqueness of $X(h \circ k)$, we must have $X(h \circ k) = X(k) \circ X(h)$. The fact that $X(1_L) = 1_{X_L}$ for any distributive lattice $L$ can be proved similarly.

  For any distributive lattice $L$, let $\alpha_L \colon L \to \ClD(X_L)$ be the homomorphism defined by $a \in L \mapsto \widehat{a} \in \ClD(X_L)$. For any Priestley space $X$, we recall the definition of an order-homeomorphism $\beta_X \colon X \to X_{\ClD(L)}$ from Proposition~\ref{prop:priestley-double-dual}; note that, compared to the notation of Definition~\ref{dfn:equivalence}, we here use the notation $\beta$ for the inverse map. For any $x \in X$, $\beta(x)$ is defined as the point associated to the  prime filter $\cB(x) = \{K \in \ClD(L) \mid x \in K\}$, the clopen-down-set neighborhoods of $x$. To see that $\alpha \colon 1_{\DL} \to \ClD \circ X$ is a natural transformation, let $h \colon M \to L$ be a homomorphism between distributive lattices. Then, by definition of $X$, we have, for any $a \in M$, that $X(h)^{-1}(\widehat{a}) = \widehat{h(a)}$, showing the naturality of $\alpha$. Finally, to see that $\beta \colon 1_{\Priest} \to X \circ \ClD$ is a natural transformation, let $f \colon X \to Y$ be a continuous order preserving function. We want to show that the following diagram commutes:
\begin{center}
  \begin{tikzcd}
	X \arrow{r}{f} \arrow{d}[swap]{\beta_X}  & Y \arrow{d}{\beta_Y} \\
	X_{\ClD(X)} \arrow{r}{X(f^{-1})} & X_{\ClD(Y)}
  \end{tikzcd}
\end{center}
  Let $x \in X$ be arbitrary. Let us write $F_1$ and $F_2$ for the prime filters of $\ClD(Y)$ that correspond to the two points $X(f^{-1})(\beta_X(x))$ and $\beta_Y(f(x))$, respectively. We need to show that $F_1 = F_2$. Recall from Exercise~\ref{exe:morphisms-concrete} that, for any homomorphism $h \colon M  \to L$, the prime filter corresponding to $X(h)(x)$ is $h^{-1}(F_x)$. Thus, for any $K \in \ClD(Y)$, we have 
  \[ K \in F_1 \iff f^{-1}(K) \in \beta_X(x) \iff x \in f^{-1}(K) \iff f(x) \in K \] 
  which is equivalent to $K \in F_2$.
\end{proof}

We now also prove that this duality has a further property, which, in
categorical terms, makes it an equivalence \emph{enriched over partially ordered sets}. The general definition of
enriched equivalence is beyond the scope of this book. 
In our categories $\DL$ and $\Priest$, the enriched structure
amounts to the fact that the sets of morphisms are in fact \emph{posets}. 
That
is, if $L$ and $M$ are distributive lattices, then the set $\Hom_{\DL}(L,M)$ is
partially ordered by the point-wise order: for a parallel pair of homomorphisms $h, h' \colon L \to M$, we say $h \leq h'$ if $h(a) \leq_M h'(a)$ for every $a \in L$. Similarly, for any Priestley spaces $X$ and $Y$, $\Hom_{\Priest}(X,Y)$ is also ordered point-wise. The fact that Priestley duality is a poset-enriched now boils down to the fact that the equivalence functors are anti-order-isomorphisms between the $\Hom$-posets, as we show now. 
\begin{proposition}\label{prop:priestley-2-eq}
  For any distributive lattices $L$, $M$, the function $$X \colon \Hom_\DL(M, L) \to \Hom_\Priest(X_L,X_M)$$ is an anti-isomorphism of partial orders.

  That is, for any parallel pair of homomorphisms $h, h' \colon M \to L$, we have $h \leq h'$ if, and only if, $X(h') \leq X(h)$, where the symbol $\leq$ denotes the pointwise order between lattice homomorphisms and between Priestley morphisms, respectively. 
\end{proposition}
We note that a similar result to this proposition, in the context of duality for relations and finite-meet-preserving functions, was already shown in the proof of Proposition~\ref{prop:relation-dual-to-box}; compare in particular equation (\ref{eq:2-eq-box-relations}) in that proof.
\begin{proof}
Recall (see Exercise~\ref{exe:morphisms-concrete}) that $X(h)$ sends any point $x \in X_L$ to the point $X(h)(x)$ with associated prime filter $F_{X(h)(x)} := h^{-1}(F_x)$. Suppose that $h \leq h'$ pointwise and let $x \in X_L$. For any $a \in M$, $h(a) \in F_x$ implies $h'(a) \in F_x$, since $F_x$ is an up-set. Hence, $F_{X(h)(x)} \subseteq F_{X(h')(x)}$, which means by definition of the order on Priestley spaces that $X(h')(x) \leq X(h)(x)$. For the converse, suppose that $h \nleq h'$ in the pointwise order. Pick $a \in M$ such that $h(a) \nleq h(a')$. By the distributive prime filter-ideal theorem (Theorem~\ref{thm:DPF}), pick a prime filter $F_x$ in $L$ such that $h(a) \in F_x$ and $h'(a) \not\in F_x$. For the corresponding point $x \in X_L$, we have $h^{-1}(F_x) \not\subseteq (h')^{-1}(F_x)$. Thus, $X(h')(x) \nleq X(h)(x)$.
\end{proof}

We note in particular  (see Exercise~\ref{exe:adjunctionpreservedbyduality}) 
that it follows from Proposition~\ref{prop:priestley-2-eq} that, for a pair of 
homomorphisms between distributive lattices $f \colon L \leftrightarrows 
M \colon g$, $(f,g)$ is an adjoint pair if, and only if, the pair $X(f) \colon X_M 
\leftrightarrows X_L \colon X(g)$ is an adjoint pair between the underlying 
posets of the dual Priestley spaces, where $X(f)$ is the lower adjoint 
and $X(g)$ is the upper adjoint.
Essentially this fact (phrased in the setting of spectral spaces) 
will be crucial in Section~\ref{sec:bifinite} of Chapter~\ref{chap:DomThry}.

\subsection*{A categorical perspective on subs and quotients}
\begin{theorem}\label{thm:subquotient-duality-categorically}
Let $L$ be a distributive lattice and $X_L$ its dual space. The complete 
lattice of sublattices of $L$ is anti-isomorphic to the complete lattice of 
quotient spaces of $X_L$, and the complete lattice of quotients of $L$ is 
anti-isomorphic to the complete lattice of closed subspaces of $X_L$.
\end{theorem}
\begin{proof}
By Proposition~\ref{prop:surj-inj} , Priestley duality restricts to a dual equivalence 
between injective morphism into $L$ and surjective morphisms from $X_L$, and also 
between surjective homomorphisms from $L$ and embeddings in $X_L$. These lattices 
are concretely realized as stated in the theorem. To see that these are order reversing 
bijections, see Proposition~\ref{prop:sublattice-quotientspace} and 
Proposition~\ref{prop:quotientlattice-subspace}.  
\end{proof}

\subsection*{A categorical perspective on operators}
We revisit the duality for unary operators described in Section~\ref{sec:unaryopduality}, giving in particular the proof of the generalized duality Theorem~\ref{thm:unaryboxduality}, for which we now have the required terminology in place. You are invited to fill in some small details in Exercise~\ref{exe:compatible-relation-duality-details}.
\begin{proof}[Proof of Theorem~\ref{thm:unaryboxduality}]\phantomsection\label{unaryboxproof}
Let us write $\DL_m$ for the category of distributive lattices with finite-meet-preserving functions, and $\Priest_{R^{\uparrow}}$ for the category of Priestley spaces with upward-compatible relations; the latter is indeed a category under relational composition, with the identity morphism on a Priestley space $X$ given by the upward Priestley compatible relation ${\geq} \subseteq X \times X$. We define a functor $P_{\Box} \colon \Priest_{R^{\uparrow}} \to \DL_m$ by sending an object $X$ to the lattice of clopen down-sets of $X$, and sending an upward Priestley compatible relation $R \subseteq X \times Y$ to the finite-meet-preserving function $P_{\Box}(R) := \forall_{R^{-1}}[-] \colon \ClD(Y) \to \ClD(X)$. It is straightforward to verify that $P_\Box$ is a well-defined functor. Moreover, Proposition~\ref{prop:relation-dual-to-box} shows that $P_\Box$ is full and faithful. Finally, $P_\Box$ is essentially surjective because any distributive lattice $L$ is isomorphic (also in $\DL_m$) to $\ClD(X)$ where $X$ is its dual Priestley space. It follows from Theorem~\ref{thm:charequivalence} that $P_{\Box}$ is part of a dual equivalence.
\end{proof}
We also show how a duality theorem for finite-\emph{join}-preserving functions can be deduced from this in a purely abstract way.

\begin{theorem}\label{thm:unarydiamondduality}
  The category of distributive lattices with fi\-nite-join-pre\-ser\-ving functions is dually equivalent to the category of Priestley spaces with downward Priestley compatible relations.
  \end{theorem}
\begin{proof}
The category $\DL_j$ of distributive lattices with finite-join-preserving functions is isomorphic to the category of distributive lattices with finite-meet-preserving functions, via the isomorphism which sends a distributive lattice $L$ to its opposite $L^\op$, and a morphism to itself. Similarly, the category of Priestley spaces with downward Priestley compatible relations is isomorphic to the category of Priestley spaces with upward Priestley compatible relations, via the isomorphism which sends a Priestley space $X$ to the space with the same topology, and the opposite order. The claimed dual equivalence is obtained as the conjugate of the dual equivalence of Theorem~\ref{thm:unaryboxduality} under these isomorphisms.
\end{proof}
Since the proof of Theorem~\ref{thm:unarydiamondduality} is very abstract, we also give a more concrete description of the action of the dual equivalence functor, $P_\Diamond \colon \DL_j \to \Priest_{R^{\downarrow}}$, on morphisms. Fix two Priestley spaces $X$ and $Y$. The bijection $$P_\Diamond(X,Y) \colon \Priest_{R^{\downarrow}}(X,Y) \to \DL_j(\ClD(Y),\ClD(X))$$ may be defined as a composition of four bijections: 

\begin{align} \label{eq:bijection-chain}
	\Priest_{R^{\downarrow}}(X,Y) &\stackrel{\id}{\longrightarrow} \Priest_{R^{\uparrow}}(X^\op,Y^\op) \\
	&\stackrel{\ref{thm:unaryboxduality}}{\longrightarrow} \DL_m(\ClU(Y), \ClU(X)) \nonumber \\
	&\stackrel{\id}{\longrightarrow} \DL_j(\ClU(Y)^\op, \ClU(X)^\op) \nonumber \\
	&\stackrel{c}{\longrightarrow} \DL_j(\ClD(Y),\ClD(X)). \nonumber
\end{align}
Here, the first and third bijections are identity functions, using the obvious
facts that a downward compatible relation from $X$ to $Y$ is the same as an
upward Priestley compatible relation from the opposite Priestley space $X^\op$
to $Y^\op$, and that a finite-meet-preserving function from a lattice $L$ to a
lattice $M$ is the same as a finite-join-preserving function from $L^\op$ to
$M^\op$. The second bijection is given by the dual equivalence of
Theorem~\ref{thm:unaryboxduality}, modulo the fact that a clopen down-set of
the Priestley space $X^\op$ is the same thing as a clopen up-set of $X$, and
similarly for $Y$. The fourth and last bijection is ``conjugation by
complementation'': $c$ sends a finite-join-preserving function $h \colon \ClU(Y)^\op \to
\ClU(X)^\op$ to the finite-join-preserving function $c(h) \colon \ClD(Y) \to \ClD(X)$,
defined by $c(h)(K) := X \setminus h(Y \setminus K)$. Unraveling the
definitions, the bijection $P_{\Diamond}(X,Y)$ can be seen to send $R$ to
$\exists_{R^{-1}}$, as defined in (\ref{eq:diamondR}) at the end of Section~\ref{sec:unaryopduality}.

\begin{example}\label{exa:modal-alg-functor}
  This example gives a functorial point of view on unary modal operators on a distributive lattice, which we studied in Section~\ref{sec:unaryopduality}. Recall that a \emphind{meet-semilattice} is a poset $(M, \leq)$ in which all finite sets have an infimum; equivalently, a meet-semilattice is a commutative monoid $(M, \wedge, \top)$ in which $a \wedge a = a$ for all $a \in M$. 
  Denote by $\cat{SL}$ the category of meet-semilattices with functions preserving all finite meets.
  For any meet-semilattice $M$, we make a set of symbols $\Box M=\{\Box a\mid a\in M\}$ and define the distributive lattice 
  \[
  F_\Box(M) := F_{\mathbf{DL}} (\Box(M))/\theta_M,
  \]
  where we recall from Section~\ref{sec:free-description} that $F_{\mathbf{DL}}(V)$ denotes the \emphind{free distributive lattice} over a set $V$, and we define $\theta_M$ to be the congruence generated by the set of pairs of the form 
  \[(\Box\big(\bigwedge G\big),\bigwedge\{\Box a\mid a\in G\})\ , \]
  where $G$ ranges over the finite subsets of $M$.  
  Since this equivalence relation $\theta_M$ exactly says that `box preserves finite meets', the lattice $F_\Box(M)$ can be thought of as the lattice generated by the act of `freely adding one layer of a box operation'. The function $a \mapsto [\Box a]_{\theta_M}$ is an order-embedding that preserves finite meets (see Exercise~\ref{exer:boxfunctor}).
  
  The constructor $F_\Box$ may be extended to a functor from the category of
  meet-semilattices to the category of distributive lattices, in such a way
  that, for any finite-meet-preserving function $h\colon M\to M'$  of
  meet-semilattices, we have $F_\Box h([\Box a]_{\theta_M}) = [\Box
  h(a)]_{\theta_{M'}}$ (see Exercise~\ref{exer:boxfunctor}). Pre-composing
  $F_\Box$ with the forgetful functor from distributive lattices to
  meet-semilattices, we obtain a functor $T_\Box$ from $\DL$ to itself. This
  functor $T_\Box$ is closely related to the category of positive $\Box$-modal
  algebras,
  where a \emphind{positive $\Box$-modal algebra}\index{modal algebra!positive} is defined as a pair $(L,\Box)$, with $L$ a distributive lattice and $\Box \colon L \to L$ a finite-meet-preserving function. 
Indeed, there is a one-to-one correspondence between $\Box$ operations on a distributive lattice $L$ and so-called $T_\Box$-algebra structures based on $L$ (see Exercise~\ref{exer:boxfunctor}). 
  
  The functor $T_\Box \colon \DL \to \DL$ is one of the program constructors in Abramsky's \emph{program logic}, which we will study in Chapter~\ref{chap:DomThry}. It is dual to the \emphind{Vietoris space} construction of topology, and is known in domain theory as the \emphind{Smyth powerdomain}, see also Definition~\ref{def:uppervietoris} in the next chapter. For a survey of dualities for modal logics in a general categorical setting, and the related development of \emph{coalgebraic modal logic}, see \cite{Ven2007}.
  \end{example}
  \ourexercises

\begin{ourexercise}\label{exe:adjunctionpreservedbyduality}
  Prove that, for any pair of homomorphisms between distributive lattices $f \colon L \leftrightarrows M \colon g$, $(f,g)$ is an adjoint pair if, and only if, the pair $X(f) \colon X_M \leftrightarrows X_L \colon X(g)$ is an adjoint pair between the underlying posets of the dual Priestley spaces, where $X(f)$ is the lower adjoint and $X(g)$ is the upper adjoint. \hint{Combine Proposition~\ref{prop:priestley-2-eq} and Exercise~\ref{exe:adjunctions} of Chapter~\ref{ch:order}.}
\end{ourexercise}

\begin{ourexercise}\label{exe:compatible-relation-duality-details}
	This exercise asks you to fill in a few details of the proof of Theorem~\ref{thm:unaryboxduality}.
	\begin{enumerate}
	\item Prove that the $\Priest_{R^{\uparrow}}$, the category of Priestley spaces with upward Priestley compatible relations, is indeed a category.
	\item Prove that the functor $P_\Box$ defined in the proof of Theorem~\ref{thm:unaryboxduality} is indeed a well-defined functor.
	\end{enumerate}
\end{ourexercise}

\begin{ourexercise}\label{exer:boxfunctor}
  For a meet-semilattice $M$, consider the lattice $F_\Box(M)$ defined in Example~\ref{exa:modal-alg-functor}. For every $a \in M$, define $e(m) := [\Box m]_{\theta_{M}}$.
  \begin{enumerate}
  \item Prove that, for any distributive lattice $L$ and any fi\-nite-meet-pre\-ser\-ving function $f \colon M \to L$, there exists a unique $\DL$ homomorphism $\bar{f} \colon F_\Box(M) \to L$ such that $\bar{f} \circ e = f$.
  \item Prove that $e$ is finite-meet-preserving, injective, and that the image of $e$ generates $F_\Box(M)$.
  \item  Prove that the object assignment $F_\Box$ extends to a functor $\cat{SL} \to \DL$, which is left adjoint to the forgetful functor $U \colon \DL \to \cat{SL}$ that sends any distributive lattice  $(L, \top,\bot,\wedge,\vee)$ to its \emphind{meet-semilattice reduct}, $(L,\top,\wedge)$, and any lattice homomorphism to itself, now viewed as a meet-semilattice homomorphism.
  
  {\it Note.} It is a general fact in category theory that ``universal arrows yield left adjoints'', see for example \cite[Section 3.1]{Borceux1} for the general theory, which you only prove in a specific example here.
  \end{enumerate}
  For a functor\footnote{There is a related, but different, notion of \emph{algebra for a monad}, which we do not treat in this book (at least not explicitly).} $T$ from $\cat{C}$ to itself, an \emph{algebra for the functor $T$}, or simply \emph{$T$-algebra}, \index{algebra for an endofunctor} is defined as a pair $(A,h)$ where $A$ is an object in $\cat{C}$ and $h\colon T(A)\to A$ is a morphism in $\cat{C}$. A \emph{morphism of $T$-algebras}, $g\colon(A,h)\to (A',h')$, is defined to be a morphism $g\colon A\to A'$ in $\cat{C}$ so that $g \circ h = h' \circ T(g)$, that is, the following diagram commutes:
  \[
  \begin{tikzpicture}
    \node (A) {$T(A)$};
    \node (B) [node distance=3cm, right of=A] { $T(A')$};
    \node (C) [node distance=2cm, below of=A]{$A$};
    \node (D) [node distance=3cm, right of=C] {$A'$};
    \draw[->] (A) to node [above,midway] {$T(g)$} (B);
    \draw[->] (C) to node [above,midway]{$g$} (D);
     \draw[->] (A) to node [left,midway] {$h$} (C);
    \draw[->] (B) to node [right,midway]{$h'$} (D);
  \end{tikzpicture}
  \]
  Recall that a \emphind{positive $\Box$-modal algebra} is a pair $(A,\Box)$ where $A$ is a distributive lattice and $\Box \colon A \to A$ preserves finite meets.
  Write $T_\Box$ for the composite functor $\DL \stackrel{U}{\to} \cat{SL}
  \stackrel{F_\Box}{\to} \DL$. The rest of this exercise makes precise the idea that $T_\Box$ \emph{freely adds one layer} of the unary modal $\Box$ operator. For further details, see \cite{Abr1988,Ghi1995,BeKu2007}. The analogous construction for a binary implication-type operator  is discussed and applied in Section~\ref{sec:funcspace}.
  \begin{enumerate}[resume]
  \item Show that the category of $T_\Box$-algebras is
  isomorphic to the category of positive $\Box$-modal algebras.
  \end{enumerate}
We will now construct the free positive $\Box$-modal algebra over a distributive lattice by iterated application of the functor $F:=1+T_\Box$. 
Here, $1$ is the identity functor, and $+$ is the coproduct computed in $\DL$; more explictly, for any distributive lattice $L$, $F(L)$ is the lattice $L + T_\Box(L)$, and if $h \colon L \to M$ is a lattice homomorphism then $F(h) \colon L + T_\Box(L) \to M + T_\Box(M)$ is given by $h + T_\Box(h)$.
   
Let $L$ be a distributive lattice, and consider the countable chain of distributive lattices 
  \begin{equation*}
  L \stackrel{e_0}{\hookrightarrow} L+T_{\Box}(L)\stackrel{e_1}{\hookrightarrow} L+T_{\Box}\big(L+T_{\Box}(L)\big) \stackrel{e_2}{\hookrightarrow} \ \dotsm,
  \end{equation*}
  where the first embedding $e_0$ is the inclusion of $L$ in the coproduct
  $L+(-)$, and the subsequent embeddings are defined by $e_{n+1}:=e_0+T_{\Box}(e_n)$.
  \begin{enumerate}[resume]
  \item Show that the colimit of the chain is the free positive $\Box$-modal algebra over $L$.

  \item Show that, when viewing the chain as a tower of sublattices of the free positive $\Box$-modal algebra over $L$, for each $n \geq 0$, the $n$th lattice consists precisely of the elements that can be described by a term in which the maximum nesting depth of the operation $\Box$ is less than or equal to $n$.
  \end{enumerate}
  \end{ourexercise}


\chapter{Omega-point duality}\label{chap:Omega-Pt}
In this chapter we consider a Stone-type duality that applies to spaces more
general than the ones dual to distributive lattices. We begin
(Section~\ref{sec:StoneSpaces}) by giving Stone's original 
duality~\parencite{Sto1937/38} for bounded distributive lattices with a 
class of (unordered) topological
spaces that we call \emph{spectral} spaces in this book.
In particular, we show in Theorem~\ref{thm:Stone-isom-Priestley} that the
category of spectral spaces is isomorphic to the category of Priestley spaces, 
establishing exactly how  Stone and Priestley duality for bounded distributive lattices are two presentations of the same mathematical result.
This recasting of Priestley duality uses the correspondence between compact ordered spaces and stably compact spaces already discussed in Section~\ref{sec:comp-ord-sp}. Then we introduce
an adjoint pair of functors (Section~\ref{sec:Omega-Pt-adjunction}) between the
category of topological spaces and a certain category of complete lattices and we show that this adjunction restricts to a duality between so-called sober
spaces and spatial frames (Section~\ref{sec:Omega-Pt-duality}). Throughout these sections, we in particular pay attention to how all these dualities fit together.
In the final section, Section~\ref{sec:funcspace}, as an application of the duality between spectral spaces and distributive lattices, we study two classical constructions on spectral spaces: powerdomains and function spaces. This leads us to introduce the notion of \emph{preserving joins at primes}, that will play a role in both of the application chapters~\ref{chap:DomThry} and \ref{ch:AutThry}. 

\section{Spectral spaces and Stone duality}\label{sec:StoneSpaces}
In the Priestley duality of Chapter~\ref{ch:priestley}, we obtain certain compact ordered spaces as duals of distributive lattices. However, in the original formulation by Stone, the dual spaces of distributive lattices are certain stably compact spaces (Definition~\ref{dfn:stab-comp-sp}). In fact the relationship between Stone's dual spaces and Priestley spaces is given by restricting the bijective correspondence between compact ordered spaces and stably compact spaces established in Theorem~\ref{thrm:COSpace-StabCompSp} of Chapter~\ref{chap:TopOrd}. The precise relationship between distributive lattices, Priestley spaces, and spectral spaces will be given in the commutative triangle of equivalences in Figure~\ref{fig:dl-priest-spec-diagram} below.
The main technical ingredient needed for establishing these relationships is the following.

\begin{theorem}\label{thrm:Priestley-Stone}
Let $(X,\tau)$ be a stably compact space. The following statements are equivalent:
\begin{enumerate}
\item[(i)] $(X,\tau^p,\leq_\tau)$ is a Priestley space;
\item[(ii)] $(X,\tau^p,\geq_\tau)$ is a Priestley space;
\item[(iii)] $(X,\tau)$ has a base of compact-opens.
\item[(iv)] $(X,\tau^\partial)$ has a base of compact-opens.
\end{enumerate}
\end{theorem}

\begin{proof}
Since the definition of a Priestley space is self-dual with respect to order-duality, it follows that (i) and (ii) are equivalent. Since, under the correspondence of Theorem~\ref{thrm:COSpace-StabCompSp}, $(X,\tau)$ is the stably compact space corresponding to $(X,\tau^p,\leq_\tau)$ and $(X,\tau^\partial)$ is the stably compact space corresponding to $(X,\tau^p,\geq_\tau)$, it suffices to show that (iii) is equivalent to (i) and/or (ii) in order to show that all four statements are equivalent.

So suppose $(X,\tau)$ is a stably compact space for which $(X,\tau^p,\leq_\tau)$ is a Priestley space and let $x\in X$ and $U\in\tau$ with $x\in U$. We will construct a compact-open set $V$ such that $x \in V \subseteq U$. Note that $C:=U^c$ is a closed down-set which does not contain $x$. Hence, for each $y\in C$ we have $x\nleq y$ and thus there is a clopen up-set $V_y$ of the Priestley space $(X,\tau^p,\leq_\tau)$ with $x\in V_y$ and $y\in V_y^c$. It follows that $\{V_y^c\}_{y\in C}$ is an open cover of $C$ in $(X,\tau^p,\leq_\tau)$. Since $C$ is closed and thus compact, it follows that there is a finite subset $F\subseteq C$ so that $\{V_y^c\}_{y\in F}$ is a cover of $C$. Consequently, the set
$V :=\bigcap\{V_y\mid y\in F\}$ 
is a clopen up-set of $(X,\tau^p,\leq_\tau)$ with $x\in V\subseteq U$. Finally, since $V$ is an open up-set in $(X,\tau^p,\leq_\tau)$, it is open in the topology $(\tau^p)^{\uparrow}$, which is equal to $\tau$ by Proposition~\ref{prop:compordaspatch}. Also, since $V$ is closed in $(X,\tau^p,\leq_\tau)$, it is compact in $(X,\tau^p,\leq_\tau)$, and thus also in the smaller topology of $(X,\tau)$, since the compact topologies are a down-set in $\TopLat(X)$ (see Exercise~\ref{ex:down-and-up-in-Top}). We have thus shown that $(X,\tau)$ possesses a base of compact-opens.

For the converse, suppose $(X,\tau)$ is a stably compact space with a base of compact-opens, and let $x,y\in X$ with $x\nleq_\tau y$. By definition of the specialization order, it follows that there is an open $U\in\tau$ with $x\in U$ and $y\not\in U$. Also, since $(X,\tau)$ has a base of compact-opens, it follows that there is a compact-open $V\subseteq X$ with $x\in V\subseteq U$. The fact that $V$ is open in $(X,\tau)$ implies that it is an open up-set in $(X,\tau^p,\leq_\tau)$. Furthermore, since $V$ is compact and open, and thus in particular saturated, it is closed in $(X,\tau^\partial)$, and thus also in $(X,\tau^p,\leq_\tau)$. Altogether, we have that $V$ is a clopen up-set in $(X,\tau^p,\leq_\tau)$ and that $x\in V$ and $y\not\in V$, as required.%
\end{proof}

We recall from Chapter~\ref{ch:priestley} that the morphisms of Priestley spaces are the continuous and order-preserving maps, see also Example~\ref{exa:bigcats}. Thus, in order to obtain a (non-full) subcategory of $\TopCat$ that is dual to the category of bounded distributive lattices, we have to consider the corresponding stably compact spaces with proper maps (see Exercise~\ref{exer:propermaps}). Finally, if $X$ and $Y$ are stably compact spaces satisfying the equivalent statements of Theorem~\ref{thrm:Priestley-Stone}, one can simplify the definition of proper maps: a map $f\colon X\to Y$ is a proper continuous map if, and only if, the inverse image of a compact-open is compact-open  (see Exercise~\ref{exer:Stone-maps}). In this context, proper maps are also known as \emph{spectral maps}\index{spectral map}.

The above considerations naturally suggest the following definition of the
category of \emph{spectral spaces}. A shorter but more \emph{ad hoc} definition can be found in Exercise~\ref{exer:Stone-spaces}; alternative definitions in the literature often also explicitly include the statement that a spectral space is \emphind{sober} (see Definition~\ref{def:sober} below), but we will derive this as a consequence of our definition (see Exercise~\ref{exer:sober-vs-wellfiltered} below). The spaces we call spectral spaces are also known as \emphind{coherent spaces} or (non-Hausdorff) \emphind{Stone spaces} in less recent literature, but these terms have also been used with other meanings, so we will avoid using them.
\begin{definition}\label{def:StoneSpace}
A \emphind{spectral space} is a stably compact space satisfying the equivalent statements of Theorem~\ref{thrm:Priestley-Stone}. A \emphind{spectral map} is a continuous function between spectral spaces such that the inverse image of any compact-open set is compact (and open). We denote by $\Spec$ the category of spectral spaces with spectral maps.
\end{definition}

\begin{definition}\label{def:StoneFunctor}
Let $L$ be a distributive lattice.  The \emph{Stone dual space of $L$} or \emph{spectral space of $L$}\index{Stone dual space of a lattice}\index{spectral space of a lattice}, denoted $\St(L)$, is a topological space $(X,\tau)$, where $X$ comes with bijections $F_{(\ )}\colon X\to\PrFilt(L)$, $I_{(\ )}\colon X\to\PrIdl(L)$, $h_{(\ )}\colon X\to\Hom_{\mathbf{DL}}(L,\btwo)$ so that for all $x\in X$ and all $a\in L$ we have
\[
a\in F_x\quad\iff\quad a\not\in I_x \quad\iff\quad h_x(a)=\top,
\]
and the topology $\tau$ is generated by the sets
\[
\eta(a)=\{x\in X\mid a\in F_x\}=\{x\in X\mid a\not\in I_x\}=\{x\in X\mid h_x(a)=\top\}
\]
where $a$ ranges over the elements of $L$. This object assignment extends to a contravariant functor
\[
\St\colon\DL\to\Spec
\]
where a homomorphism $h\colon L \to M$ is sent to the function  $f\colon X_M\to X_L$ given by
requiring that $F_{f(x)}=h^{-1}(F_x)$, for every $x \in X_M$.
\end{definition}

\nl{$\KO(X)$}{the set of compact-open subsets of a topological space $X$}{}
Given a topological space $X$ we denote by ${\KO}(X)$ the collection of compact-open subsets of $X$. Note that, if $(X,\tau)$ is a spectral space with specialization order $\leq_\tau$, then ${\KO}(X)$ ordered by inclusion is the distributive lattice which corresponds under Priestley duality to the Priestley space $(X,\tau^p,\geq_\tau)$, where we recall that $\tau^p$ denotes the \emphind{patch topology}, that is, the smallest topology containing both $\tau$ and $\tau^\partial$. Indeed, the clopen down-sets of $(X,\tau^p, \geq_\tau)$ are the clopen up-sets of $(X,\tau^p, \leq_\tau)$, which are exactly the compact-open sets of $(X,\tau)$ by Proposition~\ref{prop:with-without-order-scs}.\ref{itm:koisclup}. By definition, any proper continuous map $f \colon X \to Y$ yields a distributive lattice homomorphism $f^{-1} \colon {\KO}(Y) \to {\KO}(X)$. In summary, we have the following commutative triangle of categorical (dual) equivalences, where the equivalence between $\Priest$ and $\Spec$ is in fact an isomorphism, see Theorem~\ref{thm:Stone-isom-Priestley} below.

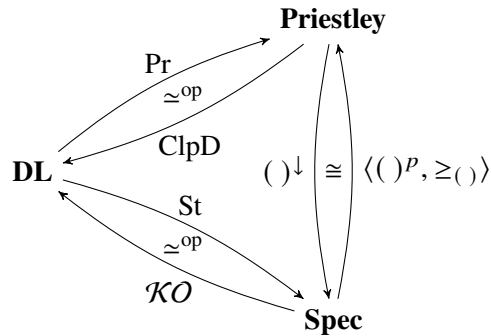
\begin{figure}[htp]
\begin{center}
\begin{tikzpicture}
  \node(DL) at (0,0) {$\DL$};
  \node(Pr) at (4,2) {$\Priest$};
  \node(Sp) at (4,-2) {$\Spec$};
  
  \node at (2,1) {$\simeq^\op$};
  \node at (2,-1) {$\simeq^\op$};
  \node at (4, 0) {$\cong$};

  \draw[->,>=stealth'] (DL) to[bend left=11] node[above] {$\mathrm{Pr}$} (Pr);
  \draw[->,>=stealth'] (Pr) to[bend left=11] node[below=.1] {$\ClD$} (DL);

  \draw[->,>=stealth'] (DL) to[bend left=11] node[above] {$\mathrm{St}$} (Sp);
  \draw[->,>=stealth'] (Sp) to[bend left=11] node[below=.1] {$\KO$} (DL);

  \draw[->,>=stealth'] (Pr) to[bend right=11] node[left] {$(\ )^{\downarrow}$} (Sp);
  \draw[->,>=stealth'] (Sp) to[bend right=11] node[right] {$\langle(\ )^{p},\geq_{(\ )}\rangle$} (Pr);
\end{tikzpicture}
\caption{A commutative triangle of equivalences. The equivalence between $\Priest$ and $\Spec$ is an isomorphism of categories, see Theorem~\ref{thm:Stone-isom-Priestley}. The dual equivalence between $\DL$ and $\Priest$ is the one described in Theorem~\ref{thm:priestley-duality-long}. The functors  between $\DL$ and $\Spec$ have been defined in this section; it is a dual equivalence because it is the composition of a dual equivalence with an isomorphism.}
\label{fig:dl-priest-spec-diagram}
\end{center}
\end{figure}

\begin{theorem}\label{thm:Stone-isom-Priestley}
  The categories {\Spec} and {\Priest} are isomorphic. On objects the isomorphisms are given, respectively, by the restriction of the patch construction for stably compact spaces, equipped with the reverse of the specialization order:
  \begin{align*}
    \langle(\ )^{p},\geq_{(\ )}\rangle \ \colon \ {\Spec}\  & \to {\Priest}\\
             (X,\tau)\  & \mapsto (X,\tau\vee\tau^\partial, \geq_\tau)
  \end{align*}
  and by intersecting with the dual Alexandrov topology of the order:
  \begin{align*}
  (\ )^{\downarrow}\colon {\Priest} & \to {\Spec}\\
             (X,\rho,\leq) & \mapsto (X,\rho\cap \Down(X,\leq)).
  \end{align*}
  On maps, the isomorphism of categories is simply the identity.
  \end{theorem}
  
  \begin{proof}
  By Theorem~\ref{thrm:Priestley-Stone} and the definition of spectral spaces, it follows that
  all spectral spaces are stably compact spaces and that the correspondence of Theorem~\ref{thrm:COSpace-StabCompSp}, applied with a reversal of the orderings, %
   between stably compact spaces and compact ordered spaces provides a bijective correspondence between spectral spaces and Priestley spaces. Furthermore, Exercise~\ref{exer:propermaps} shows that this bijective correspondence extends to an isomorphism of categories between the category {\Priest} and that of spectral spaces with proper maps. Finally, by Exercise~\ref{exer:Stone-maps}, proper maps between spectral spaces are precisely those for which the inverse image of a compact-open is compact-open, which in turn, by Definition~\ref{def:StoneSpace}, are the morphisms of the category {\Spec}.
   \end{proof}

\ourexercises

\begin{ourexercise}\label{exer:Stone-maps}
Show that a map between spectral spaces is proper if, and only if, the inverse image of a compact-open is compact-open. Further show that  a map between Boolean spaces is proper if, and only if, it is continuous.
\end{ourexercise}

\begin{ourexercise}\label{exer:non-proper}
  Show, by giving an example, that there may exist continuous maps $f\colon X\to Y$ between spectral spaces $X$ and $Y$ which are not proper.
  \end{ourexercise}

\begin{ourexercise}\label{exer:Stone-spaces}
Show that a topological space $X$ is a spectral space if, and only if, it is well filtered and ${\KO}(X)$ is a bounded sublattice of $\cP(X)$ and a base for the topology of $X$.
\end{ourexercise}

\begin{ourexercise}\label{exer:Stone-duality}
Using Priestley duality, Theorem~\ref{thrm:COSpace-StabCompSp} of Chapter~\ref{chap:TopOrd}, and Theorem~\ref{thrm:Priestley-Stone} above, show that
\begin{enumerate}
\item The assignment $X\mapsto {\KO}(X)$ may be extended to a contravariant functor from $\Spec$ to $\DL$ by sending a proper map $f\colon X\to Y$ to the map
\[
\KO(f)\colon\KO(Y)\to\KO(X), \ U\mapsto f^{-1}(U).
\]
 \item The pair of functors $\St$ and $\KO$ yield a duality between the categories $\Spec$ and $\DL$.
\end{enumerate}
\end{ourexercise}

\begin{ourexercise}\label{exer:co-compactdualofStone}
  Show that if $(X,\tau)$ is a spectral space, then the co-compact dual of $\tau$ is generated by the complements of the compact-opens of $(X,\tau)$. That is
  \[
  \tau^\partial=\langle U^c\mid U\in\KO(X,\tau)\rangle.
  \]
  \end{ourexercise}

\section{The Omega-point adjunction}\label{sec:Omega-Pt-adjunction}

The collection of compact-open subsets of a spectral space is closed under finite intersections and finite unions and is thus a distributive lattice. In fact, Stone duality tells us, among other things, that \emph{any} distributive lattice occurs as  $\KO(X)$ for an, up to homeomorphism, unique spectral space $X$.

Topological spaces that are not spectral spaces will in general not have a base of compact-opens, but they will have lots of bases that are sublattices of the power set. Using these bases, equipped with certain \emph{proximity} relations, leads to an alternative duality for stably compact spaces \parencite{Smy92,JS96}.
If a compact space has a base of compact-opens, then, by compactness, this is necessarily the smallest base which is closed under finite unions. Topological spaces that are not spectral spaces  will in general also not have a smallest base. To obtain a canonical base beyond spectral spaces, we are forced to take the base of all opens. 
But the collection of open subsets of a topological space is closed under arbitrary unions and is thus a complete lattice in the inclusion order. As a consequence, if we want to extend Stone duality to a duality between certain lattices and spaces beyond those with bases of compact-opens, we are forced to consider a category of complete lattices.

\begin{definition}\label{def:frame}
A \emphind{frame} is a complete lattice in which finite meets distribute over infinite joins, that is, it satisfies the Join Infinite Distributive law (see (\ref{eq:JID})). A map between frames is said to be a \emphind{frame homomorphism} provided it preserves finite meets and arbitrary joins. We denote by $\Frame$ the category of frames with frame homomorphisms.
\end{definition}

\begin{example}\label{exa:frame-of-opens}
Given a topological space $X$, the collection $\Omega(X)$ of all opens is closed under binary intersections and arbitrary unions. Since, in the power set lattice $\cP(X)$, (arbitrary) intersections distribute over arbitrary unions, $\Omega(X)$ is a frame. Also, by the very definition of continuity, a continuous map $f\colon X\to Y$ yields a map
\[
\Omega(f)\colon\Omega(Y)\to \Omega(X)
\]
 given by $U\mapsto f^{-1}[U]$. Since the inverse image map preserves (arbitrary) intersections and arbitrary unions, the map $\Omega(f)$ is a frame homomorphism.
\end{example}

\begin{proposition}
$\Omega\colon\TopCat\to\Frame$ is a contravariant functor.
\end{proposition}

\begin{proof}
It is straightforward to verify that $\Omega(\id_X)=\id_{\Omega(X)}$ for any topological space $X$, and that $\Omega(f\circ g)=\Omega(g)\circ\Omega(f)$ whenever $g\colon X\to Y$ and $f\colon Y\to Z$ are continuous maps.
\end{proof}

If we want a way to recapture a space from its frame of open sets, first, we need a way of distinguishing points using open sets. As soon as a space is $T_0$ this is possible via neighborhood filters. Note that a similar definition, but restricted to clopen down-sets, was already used when we proved Priestley duality (Proposition~\ref{prop:priestley-double-dual}).

\nl{$\cN(x)$}{the neighborhood filter of a point $x$ in a topological space}{}
\begin{definition}\label{def:neighborhood-filter}
Let $X$ be a topological space and $x\in X$. The \emphind{neighborhood filter} of $x$ is given by
\[
\cN(x)=\{ U\in\Omega(X)\mid x\in U\}.
\]
\end{definition}

\begin{definition}\label{def:compl-prime-filter}
Let $L$ be a frame. A subset $F\subseteq L$ is a \emphind{completely prime filter} of $L$ if it is a proper filter satisfying, for all $S\subseteq L$,
\[
\bigvee S\in F \quad\implies\quad S\cap F\neq\emptyset.
\]
We denote by $\CompPrFilt(L)$ the set of completely prime filters of $L$.
\nl{$\CompPrFilt(L)$}{the set of completely prime filters of a frame $L$}{}
\end{definition}

The proof of the following proposition is left as Exercise~\ref{exer:proveNx}.

\begin{proposition}\label{prop:Nx}
Let $X$ be a topological space and $x\in X$. Then $\cN(x)$ is a completely prime filter of $\Omega(X)$.  
Furthermore, $X$ is $T_0$ if, and only if, the assignment $x\mapsto\cN(x)$ is injective.
\end{proposition}

Recall from Chapter~\ref{ch:order} that, for any distributive lattice $L$, $\cM(L)$ denotes the set of (finitely) meet-prime elements of $L$.

\begin{proposition}\label{prop:frame-points}
Let $L$ be a frame and $F\subseteq L$. The following statements are equivalent:
\begin{enumerate}[label=(\roman*)]
\item the set $F$ is a completely prime filter of $L$;
\item the complement $L\setminus F$ of $F$ is a principal down-set ${\downarrow}m$, for some $m \in \cM(L)$;
\item the characteristic function $\chi_F\colon L\to \btwo$, where $\chi_F(a)=1$ if, and only if, $a\in F$, is a frame homomorphism.
\end{enumerate}
\end{proposition}

\begin{proof}
First suppose $F$ is a completely prime filter of $L$. Define $m:=\bigvee (L\setminus F)$. Clearly, $L\setminus F\subseteq{\downarrow}m$. Since $F$ is completely prime and $(L\setminus F)\cap F=\emptyset$, it follows that $m\not\in F$. Therefore, since $L\setminus F$ is a down-set, $
{\downarrow}m\subseteq L\setminus F$, 
and thus ${\downarrow}m=L\setminus F$. Finally, for any finite $M\subseteq L$, if  $\bigwedge M\leq m$, then $\bigwedge M\not\in F$, because $m \not\in F$ and $F$ is an up-set. If $M\subseteq F$, since $M$ is finite and $F$ is a filter, then $\bigwedge M\in F$. So there is $m'\in M$ with $m'\not\in F$ or equivalently $m'\leq m$. Thus, $m \in \cM(L)$, and we have shown that (i) implies (ii).

Now suppose that $L\setminus F={\downarrow}m$ for some $m \in \cM(L)$. To see that $\chi_F$ preserves arbitrary joins, note that, for any $A \subseteq L$, $\chi_F(\bigvee A) = 1$ if, and only if, $\bigvee A \nleq m$, if, and only if, $a \nleq m$ for some $a \in A$, which is equivalent to $\bigvee \chi_F[A] = 1$. %
That $\chi_F$ preserves finite meets follows directly from the definition of meet-prime: for any finite $A \subseteq L$, $\chi_F(\bigwedge A) = 1$ if, and only if, $\bigwedge A \nleq m$, which, since $m \in \cM(L)$, happens if, and only if, $a \nleq m$ for every $a \in A$, and this is in turn equivalent to $\bigwedge \chi_F[A] = 1$.

For the last implication, suppose that $\chi_F\colon L\to \btwo$ is a frame homomorphism. Since $\chi_F$ preserves finite meets, it follows that $F=\chi_F^{-1}(1)$ is a filter. If $A\subseteq L$ with $\bigvee A\in F$, then $\chi_F(\bigvee A)=1$ and thus $\chi_F(a)=1$ for some $a\in A$. That is, $A\cap F\neq\emptyset$ and thus $F$ is completely prime.
\end{proof}

We are now ready to give the definition of the $\Pt$ functor which takes us from frames back to topological spaces. By Proposition~\ref{prop:frame-points}, the \emph{space of points} of a frame may be seen as based on the set of homomorphisms into the frame $\btwo$, on the set of completely prime filters, or on the set of meet-primes. As we did for the Priestley dual space in Chapter~\ref{ch:priestley}, we give a `neutral' description here, leaving it as an exercise (Exercise~\ref{exer:pt-functor}) to check the details of the equivalence of the two descriptions, since the proofs are very similar to those already given in Section~\ref{sec:topologize} for Priestley duality. Note, however, that we only obtain a contravariant adjunction here; we will show how it restricts to a duality in the next section.

\nl{$\Pt(L)$}{the space of points of a frame $L$}{}
\begin{definition}\label{def:Pt-functor}
Given a frame $L$, we denote by $\Pt(L)$ a topological space which is determined up to homeomorphism by the fact that it comes with three bijections $F\colon\Pt(L)\to \CompPrFilt(L)$, $m\colon\Pt(L)\to \cM(L)$, and $f_{(-)}\colon\Pt(L)\to\Hom_\Frame(L,\btwo)$ satisfying
\[
\forall x\in\Pt(L)\quad\forall a\in L\quad (\ a\in F(x)\ \iff\  a\nleq m(x)  \ \iff\ f_x(a)=1\ ),
\]
and the opens of $\Pt(L)$ are the sets $\widehat{a}$ for $a\in L$, defined, for $x\in\Pt(X)$, by 
\[
x\in\widehat{a}\ \iff\ a\in F(x)\ \iff\ a\nleq m(x)\ \iff\ f_x(a)=1.%
\]
The assignment $L\mapsto\Pt(L)$ extends to a contravariant functor
\[
\Pt\colon\Frame\to\TopCat
\]
by sending a frame homomorphism $h\colon L\to M$ to the continuous map $\Pt(h)\colon\Pt(M)\to\Pt(L)$ given by
\begin{align*}
\Pt(h)(x)=y\quad &\iff\quad h^{-1}( F(y))= F(x)\\
                           &\iff\quad  h^{-1}( {\downarrow}m(y))= {\downarrow}m(x)\\
                           &\iff\quad f_y=f_x\circ h.
\end{align*}
\end{definition}

\begin{theorem}\label{thrm:OmegaPoint-adjunction}
The functors $\Pt\colon\Frame\rightleftarrows\TopCat\colon\Omega$ form a contravariant adjunction with corresponding natural transformations given, for $L$ a frame, by
\begin{align*}
\eta_L\colon L &\to \Omega(\Pt(L))\\
                     a &\mapsto \{x\in\Pt(L)\mid a\in  F(x)\}
\end{align*}
and, for $X$ a topological space,
\begin{align*}
\varepsilon_X\colon & X\to \Pt(\Omega(X))%
\end{align*}
sends a point $x$ of $X$ to the element of $\Pt(\Omega(X))$ which corresponds to the neighborhood filter of $x$.
\end{theorem}

\begin{proof}
In order to check that the pair $(\Pt,\Omega)$ forms a contravariant adjunction it suffices to check that:
\begin{enumerate}
\item The assignments $\eta$ and $\epsilon$ define natural transformations, that is, for each frame homomorphism $h\colon L\to M$ and each continuous map $f\colon X\to Y$, the following diagrams commute:
\[
\begin{tikzpicture}
\node(PY) at (0,1.5) {$L$};
\node(PX) at (3,1.5) {$\Omega(\Pt(L))$};
\node(OY) at (0,0) {$M$};
\node(OX) at (3,0) {$\Omega(\Pt(M))$};

\draw[->,>=stealth'] (PY)to node[above] {$\eta_L$} (PX);
\draw[->,>=stealth'] (OY) to node[below] {$\eta_M$} (OX);
\draw[->,>=stealth'] (PY) to node[left] {$h$}(OY);
\draw[->,>=stealth'] (PX) to node[right]  {$\Omega(\Pt(h))$}(OX);
\end{tikzpicture}
\qquad
\begin{tikzpicture}
\node(PY) at (0,1.5) {$X$};
\node(PX) at (3,1.5) {$\Pt(\Omega(X))$};
\node(OY) at (0,0) {$Y$};
\node(OX) at (3,0) {$\Pt(\Omega(Y))$};

\draw[->,>=stealth'] (PY)to node[above] {$\varepsilon_X$} (PX);
\draw[->,>=stealth'] (OY) to node[below] {$\varepsilon_Y$} (OX);
\draw[->,>=stealth'] (PY) to node[left] {$f$} (OY);
\draw[->,>=stealth'] (PX) to node[right]  {$\Pt(\Omega(f))$}(OX);
\end{tikzpicture}
\]
\item The triangle identities hold, that is, for each frame $L$ and each space $X$, the following diagrams commute:
\[
\begin{tikzpicture}
\node(PY) at (0,1.5) {$\Pt(L)$};
\node(PX) at (3,1.5) {$\Pt(\Omega(\Pt(L)))$};
\node(OX) at (3,0) {$\Pt(L)$};

\draw[->,>=stealth'] (PY)to node[above] {$\varepsilon_{\Pt(L)}$} (PX);
\draw[->,>=stealth'] (PY) to node[below left] {$1$}(OX);
\draw[->,>=stealth'] (PX) to node[right]  {$\Pt(\eta_L)$}(OX);
\end{tikzpicture}
\qquad
\begin{tikzpicture}
\node(PY) at (0,1.5) {$\Omega(X)$};
\node(PX) at (3,1.5) {$\Omega(\Pt(\Omega(X)))$};
\node(OX) at (3,0) {$\Omega(X)$};

\draw[->,>=stealth'] (PY)to node[above] {$\eta_{\Omega(X)}$} (PX);
\draw[->,>=stealth'] (PY) to node[below left] {$1$} (OX);
\draw[->,>=stealth'] (PX) to node[right]  {$\Omega(\varepsilon_X)$}(OX);
\end{tikzpicture}
\]
\end{enumerate}

We leave these verifications as an exercise for the reader.
\end{proof}

\ourexercises
\begin{ourexercise}
Give an example of a topological space $X$ such that:
\begin{enumerate}
\item $\Omega(X)$ is not closed under arbitrary intersections;
\item $\Omega(X)$ does not satisfy the Meet Infinite Distributive law (that is, the order dual of JID);
\item There is a frame homomorphism from $\Omega(X)$ into $\btwo$ which does not preserve arbitrary meets.
\end{enumerate}
\end{ourexercise}

\begin{ourexercise}\label{exer:infimum-in-frame}
Let $X$ be a topological space and let $L = \Omega(X)$. Prove that, for any $S \subseteq L$, the infimum of $S$ is the interior of the set $\bigcap_{U \in S} U$, that is, $\bigwedge S = \intr(\bigcap S)$.
\end{ourexercise}

\begin{ourexercise}
  \label{exer:proveNx}
Prove Proposition~\ref{prop:Nx}.
\end{ourexercise}

\begin{ourexercise}\label{exer:pt-functor}
Show that $\Pt(L)$ is well-defined and that it is a topological space. Also show that the bijection $f_{(\ )}\colon\Pt(L)\to\Hom_\Frame(L,\btwo)$ embeds $\Pt(L)$ as a closed subspace of the space $\btwo^L$, where $\btwo$ has the Sierpinski topology (see Exercise~\ref{ex:Sierpinski}). Further, show that $\Pt(h)$ is well-defined and that it is a continuous map. Finally show that $\Pt$ is a contravariant functor as stated in Definition~\ref{def:Pt-functor}.
\end{ourexercise}

\begin{ourexercise}
Complete the proof of Theorem~\ref{thrm:OmegaPoint-adjunction}.
\end{ourexercise}

\section{The Omega-point duality}\label{sec:Omega-Pt-duality}

Recall from Corollary~\ref{cor:restricttoduality} that any contravariant adjunction, such as the one given by the functors $\Pt$ and $\Omega$ in the previous section, restricts to a maximal duality. This associated duality is obtained by restricting the functors to the full subcategories given by the objects for which the component of the natural transformation from the identity to the composition of the functors in the appropriate order is an isomorphism. 

Accordingly we want to characterize those frames $L$ for which $\eta_L\colon L\to\Omega(\Pt(L))$ is an isomorphism and those topological spaces $X$ for which $\varepsilon_X\colon X\to\Pt(\Omega(X))$ is a homeomorphism.

\begin{proposition}\label{prop:spatial}
Let $L$ be a frame. The following statements are equivalent:
\begin{enumerate}[label=(\roman*)]
\item $\eta_L\colon L\to\Omega(\Pt(L))$ is an isomorphism;
\item $\eta_L$ is injective;
\item For all $a,b\in L$, if $b\nleq a$ then there is $x\in\Pt(L)$ with $b\in F(x)$ and $a\not\in F(x)$;
\item For all $a\in L$, \quad $a=\bigwedge(\cM(L)\cap{\uparrow} a).$
\end{enumerate}
\end{proposition}

\begin{proof}
By definition of $\Pt(L)$, the opens of this space are precisely the sets in the image of $\eta_L$ and thus $\eta_L$ is always surjective. It follows that items (i) and (ii) are equivalent since a frame homomorphism is an isomorphism if, and only if, it is both injective and surjective.

Since $\eta_L$ is always a frame homomorphism, if it is injective, then it is an order embedding (see Exercise~\ref{exe:injlatthom}). On the other hand, if it is  an order embedding, then it is clearly injective. So (ii) is equivalent to saying that $\eta_L$ is  an order embedding. Further, again because $\eta_L$ is always a frame homomorphism, it is always order preserving. So (ii) is equivalent to $\eta_L$ being order reflecting. But (iii) is exactly the contrapositive of the statement that $\eta_L$ is order reflecting.

We now use Proposition~\ref{prop:frame-points} to prove the equivalence of (iii) and (iv).

In order to prove that (iii) implies (iv), we reason by contraposition:  suppose that (iv) fails, we show that (iii) also fails. Let $a\in L$ be such that (iv) fails. Denote by $b$ the element $\bigwedge (\cM(L) \cap {\uparrow}a)$. We must then have $b \nleq a$, since $a \leq b$ always holds. Now let $x \in \Pt(L)$ be arbitrary such that $b \in F(x)$. To show that (iii) fails, we need to show that $a \in F(x)$. By Proposition~\ref{prop:frame-points}, pick $m \in \cM(L)$ such that $F(x) = ({\downarrow}m)^c$. Then $b \nleq m$, so in particular $a \nleq m$, by the very definition of $b$. Thus, $a \in F(x)$ as required. 

For the converse, suppose  (iv) holds and let $a,b\in L$ with $b\nleq a$. Since  $a=\bigwedge(\cM(L)\cap{\uparrow} a)$, we have $b\nleq\bigwedge(\cM(L)\cap{\uparrow} a)$ and thus there is $m\in\cM(L)$ with $a\leq m$ but $b\nleq m$. Again by Proposition~\ref{prop:frame-points}, $F=({\downarrow}m)^c$ is a completely prime filter of $L$, and thus there is an $x\in\Pt(L)$ with $F(x)=F$. It follows that $b\in F(x)$ and $a\not\in F(x)$ and thus (iii) holds.
\end{proof}

\begin{remark}
As we saw in Section~\ref{sec:finDLduality}, in a finite distributive lattice, each element is a finite join of join-irreducible elements and a finite meet of meet-irreducible elements. This is not the case for frames in general (see Exercise~\ref{exer:meet-irred}).
\end{remark}

\begin{definition}\label{def:spatial}
A frame is said to be \emphind{spatial} provided it satisfies the equivalent conditions of Proposition~\ref{prop:spatial}. We denote by {\SpFrame} the full subcategory of {\Frame} given by the spatial frames.
\end{definition}

As we have already observed, the natural morphism $\varepsilon_X\colon X\to\Pt(\Omega(X))$ is injective if, and only if, $X$ is a $T_0$ space, see Proposition~\ref{prop:Nx}.

\begin{proposition}\label{prop:sober}
Let $X$ be a $T_0$ space. The following statements are equivalent:
\begin{enumerate}[label=(\roman*)]
\item $\varepsilon_X\colon X\to\Pt(\Omega(X))$ is a homeomorphism;
\item $\varepsilon_X$ is surjective;
\item for all $y\in\Pt(\Omega(X))$ there exists an $x\in X$ so that  $ F(y)=\cN(x)$;
\item the join-irreducible elements of the lattice $\cC(X)$ of closed sets in $X$ are precisely the closures of points, $\overline{\{x\}}={\downarrow} x$, for $x\in X$.
\end{enumerate}
\end{proposition}

\begin{proof}
Recall that the opens of $\Pt(\Omega(X))$ are, by definition, the sets of the form $\widehat{U}=\{y\in\Pt(L)\mid U\in F(y)\}$, where $U \in \Omega(X)$. We have, for any $x \in X$ and $U \in \Omega(X)$,
\begin{align*}
\varepsilon_X(x)\in \widehat{U} & \iff U\in F(\varepsilon_X(x))=\cN(x)\\
						& \iff x\in U.
\end{align*}
Thus $\varepsilon_X$ is always an open continuous embedding of topological spaces and therefore it is clear that (i) is equivalent to (ii). Also, by definition of $\varepsilon_X$, an $x\in X$ is sent to the unique element $y\in\Pt(\Omega(X))$ with the property that $F(y)=\cN(x)$. Thus (iii) is simply spelling out the statement that $\varepsilon_X$ is surjective and therefore (ii) and (iii) are equivalent.

To prove that (iii) implies (iv), suppose that (iii) holds and let $C\in\cC(X)$ be join irreducible. Then, since the lattice of closed sets of $X$ is order dual to the lattice of open sets of $X$, we have that $U=C^c$ is a meet-irreducible element of the frame $\Omega(X)$. By Proposition~\ref{prop:frame-points}, it follows that the set
\[
F_C=({\downarrow}U)^c=\{V\in\Omega(X)\mid V\not\subseteq U\}=\{V\in\Omega(X)\mid C\cap V\neq \emptyset\}
\]
is a completely prime filter of $\Omega(X)$. Now let $y\in\Pt(\Omega(X))$ be the point of $\Omega(X)$ corresponding to $F_C$, that is, the point with $F(y)=F_C$. Since (iii) holds, there is $x\in X$ with $F_C=\cN(x)$. That is, for all $V\in\Omega(X)$ we have
\[
C\cap V\neq \emptyset \quad\iff\quad x\in V.
\]
We use this to show that $C=\overline{\{x\}}$. To this end let $D$ be any closed subset of $X$ and let $V=D^c$ be the complementary open set. Then we have the following string of equivalences
\[
C\subseteq D\iff C\cap V=\emptyset \iff x\not\in V \iff x\in D.
\]
Thus we obtain the desired conclusion
\[
C=\bigcap\{D\in\cC(X)\mid C\subseteq D\}=\bigcap\{D\in\cC(X)\mid x\in D\}=\overline{\{x\}}
\]
and every join-irreducible closed set is the closure of a singleton. The fact that closures of singletons always are join irreducible in the lattice of closed sets is left for the reader as Exercise~\ref{ex:singleton-joinirr}. 

To prove that (iv) implies (iii), let $y\in\Pt(\Omega(X))$ and let $F=F(y)$ be the corresponding completely prime filter of $\Omega(X)$. By Proposition~\ref{prop:frame-points}, the complement of $F$ is the down-set of some open set $U$ which is meet irreducible in the lattice $\Omega(X)$. It follows that $C=U^c$ is a closed subset of $X$ which is join irreducible in the lattice $\cC(X)$ of closed sets in $X$. By (iv), there exists $x\in X$ with $C=\overline{\{x\}}$. Now let $V\in\Omega(X)$, then we have
\[
V\not\in F \iff V\subseteq U\iff C\cap V=\emptyset \iff \overline{\{x\}}\cap V=\emptyset \iff x\not\in V
\]
and thus $V\in F$ if, and only if, $x\in V$ so that $F=\cN(x)$ as required.
\end{proof}

\begin{definition}\label{def:sober}
A topological space is said to be \emphind{sober} provided it is a $T_0$ space and  satisfies the equivalent conditions of Proposition~\ref{prop:sober}. We denote by {\SOB} the full subcategory of {\TopCat} given by the sober topological spaces.
\end{definition}

It now follows, using Corollary~\ref{cor:restricttoduality}, that the functors $\Omega$ and $\Pt$, properly restricted and co-restricted, yield a duality between the category of spatial frames and the category of sober spaces.

\begin{theorem}\label{thrm:Omega-Point-duality}
The $\Omega$-$\Pt$ adjunction cuts down to a duality between the category of spatial frames with frame homomorphisms and the category of sober topological spaces with continuous maps.
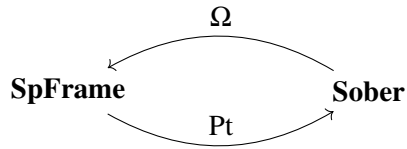
\[
\begin{tikzpicture}
  \node (A) {{\SpFrame}};
  \node (B) [node distance=4cm, right of=A] { {\SOB}};
  \draw[->,bend right] (A) to node [above,midway] {$\Pt$} (B);
  \draw[->, bend right] (B) to node [above,midway]{$\Omega$} (A);
\end{tikzpicture}
\]
\end{theorem}
We end this section with a summary of the relationship between the $\Omega$-$\Pt$, the Stone, and the Priestley dualities. Recall that we already showed how Stone and Priestley duality for distributive lattices relate in Figure~\ref{fig:dl-priest-spec-diagram}. We now add the $\Omega$-$\Pt$ duality into the mix.

Since the $\Omega$-$\Pt$ duality is for unordered topological spaces rather than ordered spaces, we compare it to Stone duality rather than to Priestley duality. As we have seen in Exercise~\ref{exer:sober-vs-wellfiltered}, every spectral space is sober, thus {\Spec} is a subcategory of {\SOB}. Note that it is not a full subcategory, since not every continuous map between spectral spaces is proper (see Exercise~\ref{exer:non-proper}). Also, Exercise~\ref{exer:join-approx} outlines a generalization of Stone duality for the category $\Spec_c$ of spectral spaces with continuous, rather than spectral functions -- the dual category of distributive lattices then has certain relations as its morphisms. The category $\Spec_c$ can be used to give a categorical description of the category of stably compact spaces with continuous functions: it is the so-called \emph{Karoubi envelope} or \emph{splitting by idempotents}. The topological fact underlying this theorem is that every stably compact space is a continuous retract of a spectral space, see for example \cite[Thm.~VII.4.6]{Johnstone1986} and \cite[Sec.~2.2]{GoolPhD} for more information.

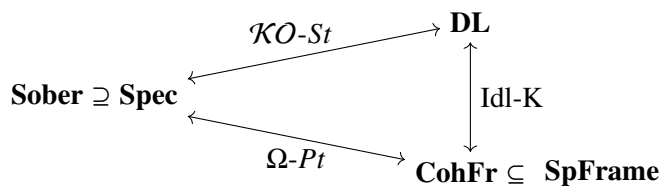
\begin{figure}[htp]
  \begin{center}
\begin{tikzpicture}
  \node (A) {$\mathbf{Sober}\supseteq \Spec$};
  \node (B) [node distance=5cm, yshift=1cm, right of=A] {${\DL}$};
  \node (C) [node distance=5cm, yshift=-1cm, right of=A] {$\mathbf{CohFr}\subseteq $};
  \node (C') [node distance=1.75cm, right of=C]{$\mathbf{SpFrame}$};
  \draw[<->] (B) to node [right,midway]{$\Idl$-$\Kel$} (C);
  \draw[<->] (A) to node [above,midway] {${\KO}$-$St$\hspace*{.6cm}} (B);
  \draw[<->] (A) to node [below,midway]{$\Omega$-$Pt$} (C);
\end{tikzpicture}
\end{center}
\caption{Comparing the $\Omega$-$\Pt$ and Stone duality on objects.}
\label{fig:stone-omega-pt-objects}
\end{figure}

The $\Omega$-$Pt$ duality is not directly a generalization of Stone duality as, under the $\Omega$-$Pt$ duality, a spectral space is sent to its entire open set lattice rather than just to the lattice of compact-opens. However, the category $\DL$ embeds into the category $\Frame$ via the \emph{ideal completion} $A\mapsto \Idl(A)$, which freely adds directed joins to $A$. A \emphind{compact} element $k$ of a frame $F$ is one such that, for every directed subset $S$ of $F$, we have $k\leq\bigvee S$ implies $k\leq s$ for some $s\,{\in}\, S$, also see Section~\ref{sec:dom} for more on compact elements in the more general setting of directedly complete posets. In particular,  we will encounter in the next chapter a more general version of this construction, $P \mapsto \Idl(P)$, which turns an arbitrary poset into a so-called \emphind{algebraic domain}.  For an arbitrary frame $F$, we denote by $\Kel(F)$ the set of compact elements of $F$, which always forms a join-subsemilattice of $F$. In $\Idl(A)$, the compact elements are the principal ideals, thus we have $A\cong \Kel(\Idl(A))$. Finally, calling \emph{coherent frames}\footnote{Coherent frames are also known as \emph{arithmetic frames} in the literature.} those frames whose compact elements form a sublattice which generates the frame by directed joins, we obtain the following diagram (Figure~\ref{fig:stone-omega-pt-objects}) which illustrates how to move back and forth between the Stone duality and the $\Omega$-$\Pt$ duality.

\noindent This diagram is not the whole story as it does not specify what
happens with morphisms. Stone duality acts on maps in $\Spec$, which is not a
full subcategory of the category $\SOB$ of sober spaces with continuous maps.
The $\Omega$-$\Pt$ duality, on the other hand, works on the full subcategory
given by the spectral spaces. To get the maps dual to the proper maps we need to
restrict the frame homomorphisms between coherent frames to those that carry
compact elements to compact elements. An alternative way of seeing this is
bitopological: a function between lattices is a homomorphism if, and only if, it
is a homomorphism between the order duals of the two lattices. In fact lattice
homomorphisms correspond precisely to those continuous maps between spectral spaces that are also continuous with respect to the co-compact dual spectral topologies. Working with spectral spaces equipped with \emph{both} these topologies and the  corresponding notion of frames, one gets a natural bitopological description of lattice homomorphisms, see \cite{Pi94,JuMo06,Bezetal2010}. 
Alternatively to restricting the morphisms on the spectral spaces and coherent frames to fit those of bounded distributive lattices, we can weaken the notion of morphism on ${\DL}$ to correspond to lattice homomorphisms $A\,{\to}\, Idl(B)$, which in turn may be seen as certain `join-approximable' relations from $A$ to $B$ (see Exercise~\ref{exer:join-approx} below).

\ourexercises

\begin{ourexercise}\label{exer:meet-irred}
Give an example of an infinite frame in which each element is a finite meet of meet-irreducible elements and one in which it is not the case.
\end{ourexercise}

\begin{ourexercise}\label{exer:non-spatial}
This exercise contains material beyond the scope of this book. It is given mainly as an indication of a means of getting one's hands on non-spatial frames.
\begin{enumerate}
\item Show that non-atomic complete Boolean algebras are examples of non-spatial frames. To this end you may proceed as follows.
\begin{enumerate}[label=(\arabic*)]
\item Every complete Boolean algebra is a frame, that is, it satisfies (JID);
\item A complete Boolean Algebra is a spatial frame if, and only if,  it is \emphind{atomic}, that is, for all $b,c\in B$, if $b\nleq c$ then there is an atom $a$ of $B$ with $a\leq b$ but $a\nleq c$.
\end{enumerate}
\item Show that there exist non-atomic complete Boolean algebras. 
\hint{We outline two rather distinct ways of proving this. (1) Show that there exist non-atomic Boolean algebras (for example, the Lindenbaum-Tarski algebra of Classical Propositional Logic on a countable set of primitive propositional variables, see Section~\ref{sec:free-description}, has no atoms at all), and show that the MacNeille completion \parencite{Mac37} of a Boolean algebra is atomic if, and only if, the original Boolean algebra is atomic; or (2) Show that the regular open subsets of a Hausdorff space form a complete Boolean algebra, and that this algebra has no atoms if the original space has no isolated points. Recall that an open subset in a topological space is called \emph{regular} if it is equal to the interior of its closure.}
\end{enumerate}
\end{ourexercise}

\begin{ourexercise}\label{exer:non-sober}
  Show that, if $X$ is an infinite set equipped with the topology generated by the cofinite subsets, then $X$ is $T_0$, but not sober.
\end{ourexercise}

\begin{ourexercise}\label{exer:sobriety}
Show that all Hausdorff spaces are sober but that some $T_1$ spaces are not.
\end{ourexercise}
\begin{ourexercise}\label{ex:singleton-joinirr}
  Let $X$ be a topological space.
  Prove that the closure of a singleton set is join irreducible in the lattice of closed subsets of $X$.
\end{ourexercise}
\begin{ourexercise}\label{exer:sober-vs-wellfiltered}
Show that if $X$ is a locally compact space, then $X$ is well filtered if, and only if, it is sober. Use this in combination with Exercise~\ref{exer:Stone-spaces} to show that a $T_0$ topological space is a spectral space if, and only if, the following two properties hold:
\begin{enumerate}
\item The collection of compact-open subsets of $X$ forms a bounded sublattice of $\cP(X)$ and a base for the topology of $X$, and 
\item $X$ is sober.
\end{enumerate}
\end{ourexercise}

\begin{ourexercise}\label{exer:sober-and-order}
Let $L$ be a frame.
\begin{enumerate}
\item Show that for $x,y\in \Pt(L)$ we have $x\leq y$ in the specialization order if, and only if, $f_x\leq f_y$ if, and only if, $F(x)\leq F(y)$.
\item Show that binary joins and meets need not exist in the specialization order in a sober space.
\item\label{itm:pt-is-dcpo} Show that the poset of completely prime filters of a frame is closed under directed union. Conclude that any sober space is a dcpo in its specialization order. 

\item Show that each open of $\Pt(L)$ is Scott open with respect to the specialization order of $\Pt(L)$.
\end{enumerate}
\emph{Note.} Item (\ref{itm:pt-is-dcpo}) of this exercise has a counterpart for distributive lattices and spectral spaces (see Exercise~\ref{exe:unionoffilters}). Also, in Proposition~\ref{prop:sober-vs-dcpo} we give a direct, elementary proof of the results in this exercise that does not use $\Omega$-$\Pt$ duality.
\end{ourexercise}

\begin{ourexercise}\label{exer:join-approx}
Let $L$ and $M$ be distributive lattices. A relation $R\subseteq L\times M$ is called \emph{join-approximable} (see, for example, \cite[Definition~7.2.24]{AbJu94}) provided that, for any $a, a' \in L$ and $b, b' \in M$, the following four properties hold: 
\begin{itemize} 
  \item if $a' \geq a {R} b \geq b'$, then $a' {R} b'$; 
  \item if $aRb$ and $aRb'$ then $a{R}b \vee b'$; 
  \item if $a{R}b$ and $a'{R}b$ then $a \wedge a'{R} b$; 
  \item if $a \vee a' {R} b$ then there exist $c, c' \in K$ such that $a{R}c$, $a'{R}c'$, and $b \leq c \vee c'$.
\end{itemize}

Denote by $X$ and $Y$ the spectral spaces dual to $L$ and $M$, respectively.
\begin{enumerate}
\item Show that there is a one-to-one correspondence between continuous functions $f\colon Y\to X$ and {\DL} homomorphisms $h\colon L\to \Idl(M)$. \hint{Use the fact that $\Omega(Y)$ is isomorphic to $\Idl(M)$.}
\item Show that there is a one-to-one correspondence between {\DL} homomorphisms $h\colon L\to \Idl(K)$ and join-approximable relations $R\subseteq L\times K$.
\item Conclude that the category whose objects are distributive lattices and whose morphisms are join-approximable relations is dually equivalent to the category of spectral spaces with continuous functions between them. 
\end{enumerate}
\emph{Note.} A related result on algebraic domains will be stated in Exercise~\ref{exer:cat-algebraic-domain}.
\end{ourexercise}

\section{Duality for spaces of relations and functions}\label{sec:funcspace}
As an application of the duality between distributive lattices and spectral spaces, in this section, we will develop a general duality theory for spaces of relations and functions, which we will apply to domains in Chapter~\ref{chap:DomThry}. We begin in the setting of general spectral spaces, and we will start with a class of constructors, on the lattice side, which are central in logic applications as already considered in Chapter~\ref{ch:methods}, namely that of freely adding one layer of a modal box, and of an implication-type operator. We will see in particular that the dual of a layer of implication-type operator gives a `binary relation space' (Theorem~\ref{thm:contfunc-upviet-spectral}). In order to obtain from this binary relation space a function space, we identify (Definition~\ref{def:jpp}) a property that we call ``preserving joins at primes'', which we will encounter again in Chapters~\ref{chap:DomThry}~and~\ref{ch:AutThry}. We prove (Corollary~\ref{cor:spectral-funcspace}) that, for certain lattices, this property gives a lattice-theoretic description of the function space construction, which will allow us to prove in the next chapter that bifinite domains are closed under that construction.

\subsubsection*{Freely adding a layer of unary modal operator}
An important idea from categorical algebra is that the addition of algebraic structure to an object $A$ in a category $\cat{C}$ can sometimes be understood as a homomorphism $TA \to A$, where $T$ is an endofunctor on $\cat{C}$. In particular, Exercise~\ref{exer:boxfunctor} spells out how a positive modal $\Box$-algebra $(A, \Box)$ can alternatively be specified by giving a distributive lattice $A$ and a lattice homomorphism $F_{\Box}(A) \to A$. Here, we can think of $F_{\Box}$ as ``freely adding one layer of unary $\Box$ operator'' to a given distributive lattice. Thus, an important step in analyzing positive modal $\Box$-algebras is to analyze this `signature' functor $F_\Box$ itself. In this subsection, we will identify the topological dual of this functor, and, in the next subsection, we carry out a similar program for the functor $F_{\to}$, which ``freely adds a layer of implication-type operator'' to a distributive lattice.

Before treating general implication-type operators, we first consider the case of a unary modal operator, which, as we will see, corresponds to the filter space construction. We also show how this construction specializes to the \emphind{upper Vietoris space}, which is also known as the \emphind{Smyth powerdomain} in the context of the spectral domains that we discuss in the next chapter.

We recall the construction of the free distributive lattice over a meet-semilattice (also see Example~\ref{exa:modal-alg-functor} and Exercise~\ref{exer:boxfunctor} for more details). Let $M$ be a meet-semilattice. The distributive lattice $F_{\Box}(M)$ is characterized uniquely by the universal property that any finite-meet-preserving function $M \to L$, with $L$ a distributive lattice, lifts uniquely to a $\DL$ homomorphism $F_{\Box}(M) \to L$.  Recall from Example~\ref{exa:modal-alg-functor} that we may construct $F_\Box(M)$ as the quotient of the free distributive lattice $F_{\DL}(\Box M)$ over the set of `formal boxes' $\{\Box a \ \mid \ a \in M\}$ under the congruence $\theta_M$ generated by the set of pairs of the form 
\begin{equation}\label{eq:meet-pres-scheme}
a_G \approx b_G \, \text{ where } a_G := \Box\Big(\bigwedge G\Big) \text{ and } b_G := \bigwedge \{\Box a : a \in G\}, 
\end{equation}
and where $G$ ranges over the finite subsets of $M$. We follow here the same convention as in Section~\ref{sec:quotients-and-subs}, that the pairs involved in generating a congruence are denoted by $a \approx b$ instead of $(a,b)$. We now give a method for calculating the Stone-Priestley dual space of this lattice $F_{\Box}(M)$, which we will subsequently generalize to the setting of implication-type operators.

Denote by $X$ the Priestley dual space of $F_\Box(M)$ and by $Y$ the Priestley dual space of $F_{\DL}(\Box M)$. Recall from Proposition~\ref{prop:freeDL-duality} that we may regard the space $Y$ as the \emph{ordered generalized Cantor space}\index{Cantor space!ordered} over $\Box M$, as given in Definition~\ref{dfn:ordered-general-cantor}; that is, the underlying set of $Y$ is $2^{\Box M}$, the partial order $\preceq_Y$ on $Y$ is the pointwise order induced by the order on $2$ in which $1 \leq 0$, that is, for $y, y' \in Y$, 
\[ y \preceq_Y y' \iffdef \text{ for all } m \in M, \text{ if } y'(\Box m) = 1 \text{ then } y(\Box m) = 1 \ , \] 
and the Priestley topology on $Y$ is generated by the subbase of sets of the form $\sem{m \mapsto i} = \{y \in 2^{\Box M} \mid y(\Box m) = i\}$, where $m \in M$ and $i \in \{0,1\}$. Since $F_\Box(M)$ is the quotient of $F_{\DL}(\Box M)$ by $\theta_M$, by Theorem~\ref{thrm:quotientlattice-subspace}, the Priestley space $X$ is order-homeomorphic to the closed subspace of $Y$ consisting of the points of $Y$ that `respect' all the pairs generating $\theta_M$; using the notation of Theorem~\ref{thrm:quotientlattice-subspace} and (\ref{eq:meet-pres-scheme}) above, $X$ is order-homeomorphic to $\sem{\{a_G \approx b_G \ \mid \ G \subseteq M \text{ finite } \} }$. Now note that, for any finite subset $G \subseteq M$, a point $y \in 2^{\Box M}$ respects $(a_G, b_G)$ if, and only if,
\[y(\Box\Big(\bigwedge G\Big)) = \bigwedge_{a \in G} y(\Box a).\]
In other words, the points $y \in Y$ that respect all pairs of the form (\ref{eq:meet-pres-scheme}) are exactly the points $y \in 2^{\Box M}$ that preserve finite meets when viewed as functions $M \to 2$, 
by precomposing with the bijection $a\mapsto \Box a$. Finally, note that the meet-preserving functions $M \to 2$ correspond exactly to the \emph{filters} of $M$, by considering the inverse image of $1$. Let us write $\phi$ for the induced bijection $X \to \Filt(M)$, which can be defined concretely, for $x \in X$, by
$$\phi(x) := \{m \in M \ \mid \ x(\Box m) = 1\}.$$
Under this bijection $\phi$, the partial order $\preceq_X$ on the space $X$, inherited from the order $\preceq_Y$ on $Y$, corresponds to the \emph{reverse} inclusion on $\Filt(M)$, that is, for any $x, x' \in X$,
$x \preceq_X x'$ if, and only if, $\phi(x) \supseteq \phi(x').$
Since the set $\Box M$ generates $F_\Box(M)$ as a lattice, the lattice of compact-opens of the spectral topology $\tau^{\downarrow}$ on $X$ is also generated by the sets of the form $\widehat{\Box a}$, as $a$ ranges over the elements of $M$. Note that, for any $a \in M$, the direct image under $\phi$ of the set $\widehat{\Box a}$ is the set
\nl{$\tilde{a}$}{the subset of the filter space defined by an element $a$}{}
$$\widetilde{a} := \{ F \in \Filt(M) \ \mid \ a \in F\}.$$
Note that the collection of sets $\widetilde{a}$, as $a$ ranges over the elements of $M$, is closed under finite intersections, since $\Filt(M) = \widetilde{\top}$ and $\widetilde{a} \cap \widetilde{b} = \widetilde{a \wedge b}$ for any $a, b \in M$. Denote by $\rho$ the topology generated by the base $\{\widetilde{a} \mid a \in M\}$. It then follows from the above calculations that the function $\phi$, viewed now as a map between topological spaces $(X,\tau^{\downarrow}) \to (\Filt(M),\rho)$, is a homeomorphism. Therefore, $\rho$ is a spectral topology,  and $F_\Box(M)$ is isomorphic to the lattice of compact-opens of $\rho$, via the unique homomorphism extending the function which sends $\Box a$ to $\tilde{a}$, for $a \in M$. In particular, each set of the form $\widetilde{a}$ is compact-open in $\rho$. Note also that any compact-open set in $(\Filt(M), \rho)$ is a finite union of sets from this base, since any open set in $\rho$ is a union of sets from the base.

We summarize our findings in the following proposition.
\begin{proposition}\label{prop:filter-box}
Let $(M, \wedge,\top)$ be a meet-semilattice. The spectral space dual to $F_\Box(M)$ is homeomorphic to the set of filters $\Filt(M)$ of $M$, equipped with the topology generated by the base $\{\widetilde{a} \mid a \in M\}$. Moreover, the compact-open sets of $\Filt(M)$ are exactly the finite unions of sets in this base.
\end{proposition}

To finish this first subsection, suppose that $M$, in addition to being a meet-semilattice, is itself also a distributive lattice. One may then directly describe the space of filters $\Filt(M)$ in terms of the spectral space dual to $M$, using the following general construction.
\begin{definition}\label{def:uppervietoris}
Let $X$ be a topological space. The \emphind{upper Vietoris space} on $X$ is the topological space  $\mathcal{V}^{\uparrow}(X)$ whose points are the compact-saturated subsets of $X$, and whose topology is generated by the base consisting of the sets
  \[
  \Box U:=\{K\in\mathcal{V}^{\uparrow}(X)\mid K\subseteq U\},\ \text{ for } U \in \tau.
  \]
\end{definition}
Note that, for $K, K' \in \mathcal{V}^{\uparrow}$, we have $K \leq K'$ in the specialization order if, and only if, $K' \subseteq K$.

In the context of general non-Hausdorff topological spaces, the name `Vietoris hyperspace' has been associated to various spaces whose points are the subsets of a given topological space, with a certain topology induced from that space. Often, the points of such a hyperspace are the \emph{closed} sets, and indeed, in our definition here, the points of $\cV^{\uparrow}(X)$ are in fact the closed subsets of the space $X^\partial$, the co-compact dual of $X$.  
\begin{proposition}\label{prop:uppervietoris-filter}
Let $M$ be a distributive lattice, and let $X$ be the spectral space of $M$. The topological space $\Filt(M)$ is homeomorphic to the upper Vietoris space $\mathcal{V}^{\uparrow}(X)$.
\end{proposition}
\begin{proof}
Let $f \colon \mathcal{V}^{\uparrow}(X) \to \Filt(M)$ be the function that sends a compact-saturated set $K$ of $X$ to the filter $\{a \in M \ \mid \ K \subseteq \widehat{a}\}$. This function $f$ has a two-sided inverse, $g \colon \Filt(M) \to \mathcal{V}^{\uparrow}(X)$, which sends a filter $F$ of $M$ to the set $\bigcap \{ \widehat{a} \ \mid \ a \in F\}$, which is indeed compact and saturated in $X$. 

Now, to see that $g$ is a two-sided inverse to $f$, let $F$ be a filter of $M$ and $b \in M$. Then, using well-filteredness of the spectral space $X$, $g(F) \subseteq \widehat{b}$ if, and only if, there exists $a \in F$ such that $\widehat{a} \subseteq \widehat{b}$, which is clearly equivalent to saying that $b \in F$. Thus, $fg(F) = F$. Conversely, for any compact-saturated set $K$ of $X$, we have $gf(K) = \bigcap \{ \widehat{a} \ \mid \ K \subseteq \widehat{a} \} = K$, as can be seen from the fact that $X^\partial$ is a spectral space in which the sets $\widehat{a}^c$ form a base for the open sets. 

To see that $f$ is continuous, let $a \in M$ be arbitrary. For any $K \in
\cV^{\uparrow}(X)$, we have $f(K) \in \widetilde{a}$ if, and only if,  $K
\subseteq \widehat{a}$ if, and only if, $K \in \Box \widehat{a}$, so $f^{-1}(\widetilde{a}) = \Box \widehat{a}$, which is open in $\cV^{\uparrow}(X)$.  

Finally, to see that $g$ is continuous, let $\Box U$ be a basic open set of $\cV^{\uparrow}(X)$, where $U$ is an arbitrary open in $X$. Let $F$ be a filter of $M$ and suppose that $g(F) \in \Box U$. We show that there exists an open set around $F$ which is mapped entirely inside $\Box U$ by $g$. By definition of $g$ and of $\Box U$, we have $g(F) = \bigcap \{\widehat{a} \ \mid \ a \in F\} \subseteq U$, and since $X$ is a spectral space, we may write $U$ as the union $\bigcup \{ \widehat{b} \ : \ b \in M, \widehat{b} \subseteq U\}$. Using the well-filteredness of $X$, pick some $a \in F$ such that $\widehat{a} \subseteq U$. Then $F \in \widetilde{a}$, and for any $F' \in \widetilde{a}$, we have $g(F') \subseteq \widehat{b} \subseteq U$, so $g(F') \in \Box U$.
\end{proof}
Note in particular that it follows from the above results that spectral spaces are closed under the upper Vietoris construction. Indeed, for a spectral space $X$, writing $M := \KO(X)$, the space $\cV^{\uparrow}(X)$ is homeomorphic to $\Filt(M)$ by Proposition~\ref{prop:uppervietoris-filter}, which is the dual space of $F_{\Box}(M)$, and is therefore a spectral space.

\begin{corollary}
Spectral spaces are closed under the upper Vietoris construction.
\end{corollary}

\subsection*{Freely adding a layer of implication-type operator}
The above calculation of the dual space of $F_\Box(M)$ suggests the following general methodology, which we will follow here to construct function spaces via their dual lattices. Suppose that a lattice $K$ is described by a set of generators $V$ and a set of equalities $R$ between distributive lattice terms in the generators, and that we want to compute its Priestley dual space $(X,\leq,\tau)$. One starts from a free distributive lattice over $V$, whose dual space is $2^V$. Then every equality $(a,b) \in R$ between generators yields a quotient of $F_{\DL}(V)$ and thus a closed subspace of $2^V$. Intersecting all these subspaces, as $(a,b)$ ranges over the equalities in $R$, yields the dual space $X$ of the distributive lattice in question. This is an application of  the quotient-lattice--subspace duality of Section~\ref{sec:quotients-and-subs}, where algebraic equations yield subspaces. 

An interesting point of this method is that one may prove properties of this dual space by examining the corresponding properties of the lattice; for example, the space $(X,\tau^{\downarrow})$ will \emph{automatically} be a spectral space, by virtue of being the dual space of a distributive lattice. This method will be used in our study of domain theory in Section~\ref{sec:DTLF}, and has also been used in the study of free algebras in varieties of lattice ordered algebras, see for example \cite{Ghilardi1992, BezhanishviliGehrke11, Ghilardi2010, CouvGo12}.

We now apply the same methodology to implication-type operators. 
Let $L$ and $M$ be distributive lattices and $X$ and $Y$ their respective dual spaces.  We will think of the  free distributive lattice over the set $L\times M$ as generated by `formal implications' 
\[
a\to b
\]
 for $a\in L$ and $b\in M$. For this reason we will here denote elements of $L\times M$, when we consider them as generators of the free distributive lattice, by $a\to b$ rather than the usual $(a,b)$. We aim to calculate (see Theorem~\ref{thm:contfunc-upviet-spectral} below) the dual space of the quotient $F_\to(L,M)$ of $F_{\DL}(L \times M)$ given by the congruence generated by the following two sets of equalities, which we also refer to as (equational) `schemes' in what follows: 
\begin{align}\label{eq:meet-snd-co}
  a\to\bigwedge G \approx \bigwedge\{a\to b\mid b\in G\}, \text{ for each $a\in L$ and finite $G\subseteq M$, }
\end{align}
\begin{align}\label{eq:meet-fst-co}
\bigvee F\to b \approx \bigwedge\{a\to b\mid a\in F\}, \text{ for each finite $F\subseteq L$ and $b\in M$}.
\end{align}

\begin{remark}\label{rem:link-with-ch4}
Notice that the schemes (\ref{eq:meet-snd-co}) and (\ref{eq:meet-fst-co}) together precisely say that $\to$, viewed as an operation from $L\times M$ to the algebra we are building is an implication-type operator, in the sense of Definition~\ref{dfn:implicationtype} in Section~\ref{sec:generalopduality}. Note also that the implication $\to$ that we build here will \emph{not} in general be a Heyting implication, that is, it will not be adjoint to a meet operation $\wedge$.

Also, as in the unary case, $F_\to$ corresponds to freely adding a layer of implication-type operator in the sense that the coproduct $L+F_\to(L,L)$ is the sublattice of depth $\leq 1$ terms in the free distributive lattice with a  implication-type operator over $L$ (see Exercise~\ref{exer:boxfunctor}).
\end{remark}
Again using Proposition~\ref{prop:freeDL-duality}, as we did in the case of a unary operator above, we identify the points of the dual space of the free distributive lattice on $L \times M$ with subsets of $L \times M$. Concretely, any prime filter $F$ of $F_{\DL}(L \times M)$ gives a subset $S_F = F \cap (L \times M)$ of $L \times M$, and conversely, if $S$ is any subset of $L\times M$, then the filter $F_S$ of $F_{\DL}(L \times M)$ generated by $S$ is prime. Indeed, by Exercise~\ref{exe:generatefilter}, we have 
\[
F_S=\{u \in F_{\DL}(L \times M) \ \mid \ \bigwedge T \leq u \text{ for some finite }  T \subseteq S\}, 
\]
and this filter is prime: for any $T \subseteq S$ finite, $\bigwedge T$ is join prime in ${F}_{\DL}(L\times M)$. That is, in disjunctive normal form, any $u\in {F}_{\DL}(L\times M)$ is of the form $u=\bigvee_{i=1}^n u_i$, where each $u_i$ is the meet of a finite set of generators $U_i\subseteq L \times M$ and we have
\[
\bigwedge T \leq\bigvee_{i=1}^n \bigwedge U_i  \quad\iff\quad T \supseteq U_i \text{ for some }i\text{ with } 1\leq i\leq n.
\]

Under this correspondence, if $S \subseteq L \times M$ is any subset and $t, u\in F_{\DL}(L \times M)$ is a pair of lattice terms in the set of generators $\{ a \to b : a \in L, b \in M\}$, then we say \emph{$S$ satisfies the equation $t \approx u$} when the prime filter $F_S$ is in the subspace $\sem{t \approx u}$ defined by the equality $t \approx u$, in the sense of the quotient-lattice--subspace duality of Theorem~\ref{thrm:sublattices-inequations}. In other words, saying that $S$ satisfies $t \approx u$ simply means that $t \in F_S$ if, and only if, $u \in F_S$. The dual space of the congruence generated by a relation $R$ on $F_{\DL}(L \times M)$ may then be identified with the collection of subsets $S \subseteq L \times M$ that satisfy all the equalities in $R$.

In particular, in the following lemma we calculate the subspaces of  $2^{L\times M}$ corresponding to the schemes (\ref{eq:meet-snd-co}) and (\ref{eq:meet-fst-co}). For any $S\subseteq L\times M$, $a\in L$, and $b\in M$, we define $S_a := \{b\in M\mid a\to b\in S\}$ and $S^b:=\{a\in L\mid a\to b\in S\}$.
\begin{lemma}\label{lem:freeoperatordual}
A subset $S\subseteq L\times M$ satisfies (\ref{eq:meet-snd-co}) if, and only if,, for each $a\in L$, $S_a$ is a filter of $M$ and $S\subseteq L\times M$ satisfies schema (\ref{eq:meet-fst-co}) if, and only if, for each $b\in M$, $S^b$ is an ideal of $L$.
\end{lemma}

\begin{proof}
We prove just the first statement as the second follows by symmetry and order
duality. Since $a \to \bigwedge G$ is a generator, $a\to\bigwedge G\in F_S$ if
and only if $a\to\bigwedge G\in S$, while $\bigwedge\{a\to b\mid b\in G\}\in
F_S$ if, and only if, $\{a\to b\mid b\in G\}\subseteq S$. That is, $S$ satisfies all equalities in (\ref{eq:meet-snd-co}) if, and only if, for each $a\in L$, the set $S_a$ is closed under finite meets. That is,  if, and only if, $S_a$ is a filter of $M$.
\end{proof}

We will establish, in Theorem~\ref{thm:contfunc-upviet-spectral} below, a first connection between the implication-type operator and the function space construction. To do so, we need the following general definition of a topology on a function space.

\nl{$[X,Y]$}{the space of continuous functions from $X$ to $Y$}{}
\begin{definition}\label{def:funcspace}
  Given two topological spaces $Z_1$ and $Z_2$, we denote the set of continuous functions from $Z_1$ to $Z_2$ by $[Z_1,Z_2]$. The \emphind{compact-open topology}\footnote{Note that the adjective compact-open in `compact-open topology', which is common in the literature, refers to the fact that the subbase is given by sets of functions determined by a compact and a (different) open set. This contrasts with the use of `compact-open set' throughout this book, meaning a single set that is compact and open, also common in the literature.} on $[Z_1,Z_2]$ is defined to be the topology generated by the subbase consisting of the sets
\nl{$K \Rightarrow U$}{the set of functions $f$ for which the direct image of $K$ is contained in $U$; a generic open set in the subbase for the topology on $[X,Y]$}{}
\[
K\Rightarrow U :=\{f\in[Z_1,Z_2]\mid f[K]\subseteq U\} \ \text{ for }  K\in\KS(Z_1) \ \text{ and } U\in\Omega(Z_2),
\]
where we recall that $\KS(Z_1)$ denotes the collection of compact-saturated subsets of $Z_1$ and $\Omega(Z_2)$ denotes the collection of open subsets of $Z_2$.
\end{definition}
See Exercise~\ref{ex:compact-open topology} for more detail and some basic properties used in the proofs below.
\begin{lemma}\label{lem:comp-open-base}
Let $L$ be a distributive lattice with dual spectral space $X$, and let $M$ be a meet-semilattice. Then the compact-open topology on $[X,\Filt(M)]$ is generated by the subbase consisting of the sets $\,\wa\to\simb$ for $a\in L$ and $b\in M$.
\end{lemma}

\begin{proof}
By Proposition~\ref{prop:filter-box}, the space $\Filt(M)$ with the topology generated by the sets $\simb$ for $b\in M$ is a spectral space, as it is the dual of  $F_\Box(M)$. Thus it follows, by Exercise~\ref{ex:compact-open topology}, that the sets $\wa\to V$ for $a\in L$ and $V$ compact-open in $\Filt(M)$ form a subbase for the topology on $[X,\Filt(M)]$. The compact-opens of $\Filt(M)$ are finite unions of sets of the form $\simb$ for $b\in M$. Thus we must show that such sets are in the topology generated by the $\wa\to\simb$ with $a\in L$ and $b\in M$. To this end, let $a \in L$ and $G \subseteq M$ be finite. Define $C := \{c \in L^G \ \mid \ a \leq \bigvee_{g \in G} c_g \}$. We will prove that
\[
\wa\to\bigcup_{g\in G}\simg=\bigcup\Big\{\bigcap_{g\in G}(\widehat{c_g}\to \simg) \ \mid \  c \in C\Big\}.
\]
For the right-to-left inclusion, let $c \in C$ be arbitrary and suppose that $f\in\bigcap_{g\in G}(\widehat{c_g}\to\simg)$. For any $x \in \widehat{a}$, since $a \leq \bigvee_{g \in G} c_g$, there is $g\in G$ with $x\in\widehat{c_g}$, and since $f\in\widehat{c_g}\to\simg$ it follows that $f(x)\in\simg$. Thus $f\in\wa\to\bigcup_{g\in G}\simg$.

For the converse, let $f\in\wa\to\bigcup_{g\in G}\simg$. For every $x\in\wa$, pick $g_x\in G$ with $f(x)\in\widetilde{g_x}$. Since $f$ is continuous, $f^{-1}(\widetilde{g_x})$ is open. Since $X$ is a spectral space, pick $a_x\in L$ with $x\in\wa_x\subseteq f^{-1}(\widetilde{g_x})$. Then $\wa \subseteq \bigcup_{x\in\wa} \wa_x$, and thus, by compactness of $\wa$, there is a finite $F\subseteq\wa$ with $\wa\subseteq\bigcup_{x\in F}\wa_x.$

Now, for each $g\in G$, define $c_g := \bigvee\{a_x\mid x\in F\text{ and } g_x=g\}$. Then, since each $g_x\in G$, 
\[
a \leq\bigvee_{x\in F}a_x = \bigvee_{g \in G}\big(\bigvee\{a_x \mid  x\in F\text{ and } g_x=g\}\big)=\bigvee_{g \in G} c_g.
\] 
That is, $c=(c_g)_{g\in G} \in C$. Also, for each $x\in\wa$, we have $\wa_x\subseteq f^{-1}(\widetilde{g_x})$, or equivalently, $f[\wa_x]\subseteq \widetilde{g_x}$, so $f\in\wa_x\to\simg_x$. Therefore we obtain
\begin{align*}
f\in\bigcap_{x\in F}(\wa_x\to\simg_x) &=\bigcap_{g\in G}\big(\bigcap\{\wa_x\to\simg\mid\ x\in F\text{ and }g_x=g\} \big)\\
                                                           &=\bigcap_{g\in G}\big(\big[\bigcup\{\wa_x\mid x\in F\text{ and }g_x=g\} \big]\to\simg\big)
							=\bigcap_{g\in G}(\widehat{c}_g\to\simg), 
\end{align*}
where the second equality holds since $\to$, with the second argument fixed, sends finite unions in the first argument to finite intersections (see Exercise~\ref{ex:compact-open topology}.\ref{itm:alg-props-arrow}.\ref{itm:arrow-of-imptype1}). Note that this is part of saying that $\to$ is an operation of implication-type.
\end{proof}
Using our results above, we will now give a concrete description of the dual space of $F_{\to}(L,M)$ for arbitrary distributive lattices. In the following theorem, we will show that, viewed on the dual spectral spaces $X$ and $Y$, the construction $F_{\to}$ corresponds to a space of upward Priestley compatible relations from $X$ to $Y$, also see Exercise~\ref{exe:upw-comp-upper-viet}.

\begin{theorem}\label{thm:contfunc-upviet-spectral}
Let $L$ and $M$ be distributive lattices and let $X$ and $Y$ be their respective dual spectral spaces. The spectral space dual to the distributive lattice ${F}_\to(L,M):={F}_{\DL}(L\times M)/\theta$, where $\theta$ is the congruence generated by the schemes (\ref{eq:meet-snd-co}) and (\ref{eq:meet-fst-co}) together, is homeomorphic to the space $[X,\cV^{\uparrow}(Y)]$ of continuous functions from $X$ to the upper Vietoris space of $Y$ in the compact-open topology.
\end{theorem}

\begin{proof}
By Proposition~\ref{prop:uppervietoris-filter}, the space $\mathcal{V}^{\uparrow}(Y)$ is homeomorphic to $\Filt(M)$. It is therefore equivalent  to prove that the dual space of $F_{\to}(L,M)$ is homeomorphic to the space $[X, \Filt(M)]$ with the compact-open topology.
Also, by Lemma~\ref{lem:freeoperatordual} and the considerations preceding it, the dual space of ${F}_\to(L,M)$ is the subspace $Z$ of $2^{L\times M}$ consisting of those subsets $S$ of $L \times M$ for which $S_a$ is a filter of $M$ for each $a\in L$ and $S^b$ is an ideal of $L$ for each $b\in M$. 

To prove the theorem, we will define a mutually inverse pair of continuous maps $\phi \colon Z \leftrightarrows [X, \Filt(M)] \colon \psi$.
For $S \in Z$, define a function $\phi(S)$ by 
\[
\phi(S)\colon X\longrightarrow \Filt(M), \quad x\mapsto \bigcup \{S_a\mid a\in F_x\}.
\]
We need to show that $\phi(S)$ is a well-defined continuous function from $X$ to $\Filt(M)$. To see that the function $\phi(S)$ is well-defined, that is, that $\phi(S)(x)$ is a filter for every $x \in X$, it suffices to show that the collection of filters $\{ S_a \mid a \in F_x\}$ is directed, by Exercise~\ref{exe:unionoffilters}. To this end, note that the assignment $a\mapsto S_a$ is order reversing; indeed, if $b \in S_{a'}$, then $a' \in S^b$, so $a \in S^b$ since $S^b$ is a down-set, so that $b \in S_a$. Thus, since $F_x$ is down-directed, the collection $\{S_a \mid a \in F_x\}$ is up-directed, as required.

We now show that $\phi(S)$ is continuous. Let $b \in M$. For every $x \in X$, we
have $x \in \phi(S)^{-1}(\simb)$ if, and only if, there exists $a \in F_x$ such that $a \to b \in S$; in a formula:
\begin{equation}\label{eq:phi-inv-simb}
\phi(S)^{-1}(\simb) = \bigcup_{c \in S_b} \widehat{c}.
\end{equation}
In particular, $\phi(S)^{-1}(\simb)$ is open for every $b \in M$, so $\phi(S)$ is a continuous function from $X$ to $\Filt(M)$.

For the inverse function $\psi$,  let $f\colon X\to\Filt(M)$ be continuous. Define 
\[
\psi(f):=\{a\to b \in L \times M \mid \wa\subseteq f^{-1}(\simb)\}.
\]
We first show that $\psi(f) \in Z$. Note that, for any $b \in M$,
\[
\psi(f)^b=\{a\in L\mid \wa\subseteq f^{-1}(\simb)\},
\]
which is an ideal because $\widehat{(\ )}$ preserves finite joins. Similarly, for any $a \in L$, 
\[
\psi(f)_a=\{b\in M\mid f[\wa]\subseteq \simb\},
\]
and since $\widetilde{(\ )}$ preserves finite meets, this is a filter. So $\psi(f) \in Z$. 

To show that $\phi$ and $\psi$ are mutually inverse, first note that, for any $b \in M$, $f \in [X, \Filt(M)]$ and $x \in X$, we have
\[ b \in \phi(\psi(f))(x) \iff \exists a \in L , x \in \widehat{a} \subseteq f^{-1}(\simb) \iff b \in f(x),\]
where we use in the second equivalence that $f$ is continuous. Thus, $\phi(\psi(f)) = f$. Also, for any $S \in Z$, $a \in L$, and $b \in M$, we have
\[ a \to b \in \psi(\phi(S)) \iff \widehat{a} \subseteq \phi(S)^{-1}(\simb) = \bigcup_{c \in S^b} \widehat{c},\]
recalling (\ref{eq:phi-inv-simb}). By compactness of $X$ and the fact that $\widehat{(\ )}$ is an order embedding, the latter inclusion is equivalent to: there exists a finite subset $C$ of $S^b$ such that $a \leq \bigvee C$. Since $S^b$ is an ideal, this is in turn equivalent to $a \in S^b$, that is, $a \to b \in S$, as required.

It remains to prove that the bijection $\phi$ and $\psi$ between $Z$ and $[X,\Filt(M)]$  is a homeomorphism. Observe that a basic open of $Z$, which is of the form
\[
\widehat{a\to b}=\{S\in Z\mid a\to b\in S\},
\]
 is sent by $\phi$ to the set
\[
\wa\to\simb=\{f\in [X,\Filt(M)]\mid \widehat{a}\subseteq f^{-1}(\simb)\}.
\]
By Lemma~\ref{lem:comp-open-base} these sets form a base for $[X,\Filt(M)]$ with the compact-open topology, and thus both $\phi$ and $\psi$ are continuous.
\end{proof}
We may recover the upper Vietoris space construction itself as a special case of the `relation space' $[X,\cV^{\uparrow}(Y)]$ used in Theorem~\ref{thm:contfunc-upviet-spectral}. Indeed, notice that $\cV^{\uparrow}(X) \cong [X^\partial,\mathbb{S}] \cong [X^\partial,\cV^{\uparrow}(1)]$, where $2^{\downarrow}$ is the Sierpinski space and $1$ is the one element space (see Exercise~\ref{exe:upperviet-spectral}).

\begin{remark}
Let $X$ and $Y$ be spectral spaces with dual lattices $L$ and $M$, respectively. Then we have shown that $[X,\cV^{\uparrow}(Y)]$ is again a spectral space. Even though $\cV^{\uparrow}(Y)$ is dual to the lattice $F_\Box(M)$, obtained by freely adding a layer of unary dual operator to $M$, notice that the compact-opens are not given by the elements of $L$ and $F_\Box(M)$, but rather by the elements of $L$ and $M$, see Example~\ref{ex:[X,Y]notspectral}.
\end{remark}

\begin{remark}
Notice that  in Theorem~\ref{thm:contfunc-upviet-spectral} we are considering
only  operators of implication type but, with some order flips, this result may
be transposed to other types of operators. For example, consider the
construction that takes distributive lattices $L$ and $M$ and produces a
distributive lattice $F_{\bullet}(L,M)$ by freely adding a binary operation
$\bullet$ that preserves finite meets in both coordinates; that is, keeping the analogue for $\bullet$ of scheme (\ref{eq:meet-snd-co}), and replacing (\ref{eq:meet-fst-co}) by $(\bigwedge F) \bullet b \approx \bigwedge \{ a \bullet b \ \mid \ a \in F\}$ for each finite $F \subseteq L$ and $b \in M$. The dual space of this construction is then $[X^\partial,\cV^{\uparrow}(Y)]$, where $X^\partial$ denotes $X$ equipped with the co-compact dual of the topology of $X$.
\end{remark}

\subsection*{Preserving joins at primes}
Let $X$ and $Y$ be spectral spaces with dual distributive lattices $L$ and $M$, respectively.
We are interested in the space of continuous functions $[X,Y]$, which can be regarded as a (generally non-spectral) subspace of the space of compatible relations $[X,\cV^{\uparrow}(Y)]$, as we explain now (see Exercise~\ref{ex:XinVietorisX} for more details). 

Denote by $\eta \colon Y \to \cV^{\uparrow}(Y)$ the embedding of $Y$ in $\cV^{\uparrow}(Y)$ given by $y \mapsto {\uparrow} y$. Then we have an injective function $[X,Y] \into [X,\cV^{\uparrow}(Y)]$ by sending $f \in [X,Y]$ to $\eta \circ f$, and a base for the topology on $[X,Y]$, now viewed as a subspace of $[X, \cV^{\uparrow}(Y)]$, is given by the sets $\wa \Rightarrow \wb = (\wa \Rightarrow \simb) \cap [X,Y]$, for $a \in L$ and $b \in M$. Even though $Y$ is a Priestley-closed subspace of $\cV^{\uparrow}(Y)$, it is not in general the case that $[X,Y]$ is a Priestley-closed subspace of $[X,\cV^{\uparrow}(Y)]$, reflecting the fact that $[X,Y]$ is not in general a spectral space; we give an explicit example of this occurrence in Example~\ref{ex:[X,Y]notspectral} below.
One would need to move to frames, sober spaces, and geometric theories to describe $[X,Y]$ as the dual of a quotient. 
However, we have a finitary approximation of the subspace $[X,Y]$, that we will give in Theorem~\ref{thrm:preserves-joins-at-primes}. Here, recall from Section~\ref{sec:quotients-and-subs} that if $Z$ is a closed subspace of a Priestley space $P$, then the corresponding congruence $\theta$ on the dual lattice $A$ of $P$ is given by
\[
\theta = \{ (a,b) \in A^2 \ \mid \ \widehat{a} \cap Z = \widehat{b} \cap Z\},
\]
so that, for any $a, b \in A$, we have
\begin{equation}\label{eq:dual-order-cong} 
  [a]_{\theta} \leq [b]_{\theta} \text{ if, and only if, for every } z \in Z, \text{ if } z \in \widehat{a} \text{ then } z \in \widehat{b}.
\end{equation}
\begin{theorem}\label{thrm:preserves-joins-at-primes}
Let $L$ and $M$ be distributive lattices, $X,Y$ their respective dual spaces, and $Z$ a closed subspace of $[X, \cV^{\uparrow}(Y)]$, viewed as a Priestley space rather than a spectral space. Denote by $\theta$ the congruence on $F_{\to}(L,M)$ corresponding to $Z$. The following are equivalent:
\begin{enumerate}[label=(\roman*)] 
  \item $Z$ is a subspace of $[X,Y]$,
  \item for every $x\in X$, $a\in F_x$, and finite subset $G\subseteq M$, there is $c\in F_x$  such that
\begin{equation}\label{eq:jppmodtheta}
[a\to \bigvee G]_\theta\ \leq\ [\bigvee\{c\to b\mid b\in G\}]_\theta.
\end{equation}
\end{enumerate}
\end{theorem}

\begin{proof}
 In light of Proposition~\ref{prop:uppervietoris-filter}, we will work with $\Filt(M)$ instead of $\cV^{\uparrow}(Y)$, and we may consider $Z$ as a subspace of $[X, \Filt(M)]$. As explained above, under the identification $\cV^{\uparrow}(Y) \cong \Filt(M)$, the subspace $Y$ of $\cV^{\uparrow}(Y)$ corresponds to the subspace $\PrFilt(M)$ of \emph{prime} filters of $M$, so $[X,Y]$ corresponds to $[X,\PrFilt(M)]$.

First suppose that (i) holds, that is, $Z \subseteq [X, \PrFilt(M)]$. Let  $x\in X$, $a\in F_x$, and $G\subseteq M$ finite. We need to show that there exists $c \in F_x$ such that (\ref{eq:jppmodtheta}) holds. We will first show that, for any $f \in Z$ with $f\in\widehat{a\to \bigvee G}$, there exist $c_f \in F_x$ and $b_f \in G$ such that $f \in \widehat{c_f \to b_f}$. To see this, let $f \in Z$  be arbitrary and suppose that $f \in \widehat{a \to (\bigvee G)}$. Then, by definition, $f[\wa]\subseteq\widetilde{\bigvee G}$ and thus, as $x\in\wa$, we have $\bigvee G\in f(x)$. Now, since $f \in Z \subseteq [X, \PrFilt(M)]$, the filter $f(x)$ is prime, so we may pick $b_f\in G$ with $b_f\in f(x)$, or equivalently, $f(x)\in\simb$. Since $f$ is continuous, pick $c_f\in L$ with $x\in\widehat{c_f}$ and $f[\widehat{c_f}]\subseteq \widetilde{b_f}$, or equivalently, $f\in\widehat{c_f}\to\widetilde{b_f}=\widehat{c_f \to b_f}$. Now the sets $\widehat{c_f \to b_f}$, for $f$ ranging over the Priestley-closed set $Z \cap \widehat{a \to (\bigvee G)}$, are a cover of this set. By compactness, pick a finite subcover, indexed by $f_1, \dots, f_n$, say, and define $c := \bigwedge_{i=1}^n c_{f_i}$. Then $c \in F_x$ since each $c_{f_i}$ is in $F_x$. Let us show that for this $c$, (\ref{eq:jppmodtheta}) holds. Writing $d := \bigvee \{c \to b \ \mid \ b \in G\}$, we want to show that $[a \to \bigvee G]_{\theta} \leq [d]_{\theta}$. For any $f \in Z$, if $f \in \widehat{a \to \bigvee G}$, then $f \in \widehat{c_{f_i} \to b_{f_i}}$ for some $1 \leq i \leq n$. It now follows that $f \in \widehat{d}$, since 
\[ c_{f_i} \to b_{f_i} \leq c \to b_{f_i} \leq d\]
where we use first that $\to$ is order preserving in the first coordinate, and then that $b_{f_i} \in G$. Thus, we have shown that for any $f \in Z$, if $f \in \widehat{a \to \bigvee G}$, then $f \in \widehat{d}$. Using (\ref{eq:dual-order-cong}), we thus conclude (\ref{eq:jppmodtheta}).

Conversely, suppose (ii) holds.  
Let $f\in[X,\Filt(M)]$ be in $Z$ and $x\in X$. If $G\subseteq M$ is finite and $\bigvee G\in f(x)$, then, as $f$ is continuous, there is $a\in L$ with $f[\wa]\subseteq\widetilde{\bigvee G}$, or equivalently, $f\in\wa\to\widetilde{\bigvee G}=\widehat{a\to(\bigvee G)}$. Now, using (\ref{eq:jppmodtheta}), it follows that there are $c\in F_x$ and  $b\in G$ with  $f\in\widehat{c\to b}$. Thus $x\in\widehat{c}$ and $f[\widehat{c}]\subseteq \widetilde{b}$ and thus $b\in f(x)$. That is, we have shown that $f(x)$ is a prime filter and thus that $f\in[X,\PrFilt(M)]\cong[X,Y]$ as required.
\end{proof}
\begin{definition}\label{def:jpp}
  Let $K, L$, and $M$ be distributive lattices, and let ${\to} \colon L \times M \to K$ be an implication-type operator. We say that $\to$ \emph{preserves joins at primes} if, for every prime filter $F$ of $L$, $a \in F$, and finite subset $G$ of $M$, there exists $c \in F$ such that
  \[ a \to \bigvee G \leq \bigvee \{c \to b \ \mid \ b \in G\}.\]
  For a congruence $\theta$ on $F_{\to}(L,M)$, we also say that $\theta$ \emph{makes }$\to$ \emph{preserve joins at primes}, or $\to$ \emph{preserves joins at primes modulo $\theta$} if the equivalent properties in Theorem~\ref{thrm:preserves-joins-at-primes} hold.\index{preserve joins at primes}\index{preserve joins at primes!modulo a congruence}\index{join preserving!at primes} 
\end{definition}
\begin{remark}\label{rem:iii}
To explain the above terminology, note that the property of preserving joins at primes modulo $\theta$ is equivalent to the property that, for each $x \in X$, the following operation preserves finite joins:
\begin{align*}
x\to(-)\colon M &\to \Idl({\mathbb F}_\to(L, M)/\theta)\\ 
                            b &\mapsto \langle[a\to b]_\theta\mid a\in x \rangle_{\rm Idl}.
\end{align*}  
In a lattice with enough join-primes, there is actually a largest congruence
that makes $\to$ preserve joins at primes, see
Corollary~\ref{cor:spectral-funcspace} below. This congruence is used crucially
in Theorem~\ref{thm:Bif-CCC} in Chapter~\ref{chap:DomThry}.

The property of preserving joins at primes is closely related to being determined by finite quotients. See Theorem~\ref{thm:bintopalg-dual-to-jpp} in Chapter~\ref{ch:AutThry} and \cite[Theorem~3.18]{Geh16} for another occurrence of this phenomenon in the setting of topological algebras on Boolean spaces. For a study of this notion via canonical extensions, see \cite{FussnerPalmigiano2019}.
\end{remark}
Recall that an element $p$ in a lattice $L$ is said to be \emphind{join prime} provided that, for any finite $F \subseteq L$, $p\leq\bigvee F$ implies that there exists $a\in F$ with $p\leq a$. Further, we say that $L$ \emph{has enough join-primes} provided every element of $L$ is the join of a finite set of join-prime elements. Also recall from Chapter~\ref{ch:order} that we denote by $\cJ(L)$ the poset of join-prime elements of $L$, with order inherited from $L$.
The following special case of Theorem~\ref{thrm:preserves-joins-at-primes}, where one of the two lattices is assumed to have enough join-primes, is central to the treatment of the function space operator in Domain Theory in Logical Form \parencite{Abr}, as we will also see in Section~\ref{sec:DTLF}.

\begin{corollary}\label{cor:spectral-funcspace}
Let $L$ and $M$ be distributive lattices and suppose $L$ has enough join-primes. Further, let $X$ and $Y$ be the dual spaces of $L$ and $M$, respectively. The space $[X,Y]$ of continuous functions from $X$ to $Y$ in the compact-open topology is dual to the congruence $\theta_{\jpp}$ on $F_{\to}(L,M)$ that is generated by the schema, for every $p \in \cJ(L)$ and $G \subseteq M$ finite, 
\begin{equation}\label{eq:jpp-scheme}
 p\to\bigvee G \approx \bigvee\{p\to b\mid b\in G\}.
\end{equation}
 
\end{corollary}

\begin{proof}
We first note that it suffices to prove that a congruence $\theta$ of ${F}_{\to}(L,M)$ makes $\to$ preserve joins at primes if, and only if, $\theta_{\jpp}$ is contained in $\theta$. Indeed, by Theorem~\ref{thrm:preserves-joins-at-primes}, this claim implies that $\theta_{\jpp}$ is the minimum congruence whose dual closed set is contained in $[X,Y]$, and thus the dual $Z_{\jpp}$ of ${F}_{\to}(L,M)/\theta_{\jpp}$ is the maximum Priestley closed subspace of $[X,Y]$. If $Z_{\jpp}$ were a proper subset of $[X,Y]$, then there would be $f \in [X,Y] \setminus Z_{\jpp}$, and $Z_{\jpp} \cup \{f\}$ would still be Priestley closed, since singletons are closed in the Priestley topology. Thus $Z_{\jpp}$ must be the full subspace $[X,Y]$.

To prove that $\theta_{\jpp}$ is indeed the minimum congruence that makes $\to$ preserve joins at primes, suppose first that $\theta$ makes $\to$ preserve joins at primes, and let $p\in \cJ(L)$ and $G\subseteq M$ finite. Instantiating the condition (\ref{eq:jppmodtheta}) for $F_x={\uparrow}p$ and $a=p$, pick an element $c\in{\uparrow}p$ so that 
\begin{align*}
[p\to\bigvee G]_{\theta}&\leq[\bigvee\{c\to b\mid b\in G\}]_{\theta}.\end{align*}
Now $p\leq c$ implies that $c\to b\leq c\to b$ in ${F}_{\to}(L, M)$, so the latter join is at most $[\bigvee \{p \to b \ \mid b \in G\}]_{\theta}$. We conclude that $$[p \to \bigvee G]_{\theta} \leq [\bigvee \{p \to b \ \mid \ b \in G\}]_{\theta}.$$ The other inequality holds because $\to$ is order preserving in its second coordinate. Thus, $\theta_{\jpp}$ is contained in $\theta$.

Conversely, suppose $\theta$ is a congruence that contains $\theta_{\jpp}$ and let $x\in X$, $a\in F_x$, and $G\subseteq M$ be finite. Since $L$ has enough join-primes, there is a finite set $F\subseteq \cJ(L)$ so that $a=\bigvee F$. Also, since $a\in F_x$ and $F_x$ is a prime filter, there is $p\in F$ with $p\in F_x$. Further, as the function $\to$ from $L\times M$ to ${F}_{\to}(L\times M)$ is an implication-type operator, we have  
\[
a\to(\bigvee G)=\bigwedge\{q\to(\bigvee G)\mid q\in F\}\leq p\to(\bigvee G).
\]
It follows that 
\begin{align*}
[a\to(\bigvee G)]_{\theta}&\leq[p\to(\bigvee G)]_{\theta}=[\bigvee\{p\to b\mid b\in G\}]_{\theta}.
\end{align*}
That is, taking $c=p$ shows that $\theta$ makes $\to$ preserve joins at primes. 
\end{proof}

We finish this subsection by giving the promised example that $[X,Y]$ itself is not always a Priestley closed subspace of $[X, \cV^{\uparrow}(Y)]$.

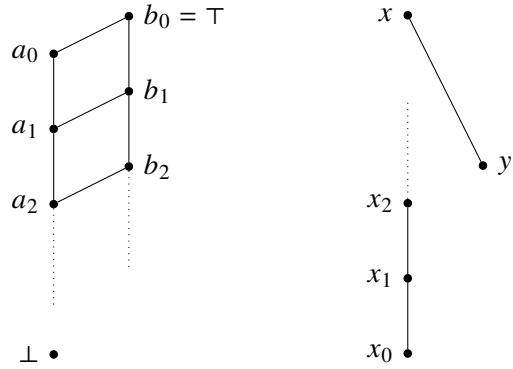
\begin{figure}[!tbp]
  \begin{center}
  \begin{tikzpicture}
  \node[draw,circle,inner sep=1pt,fill,label=right:{$b_0=\top$}] (1) at (0,0) {};
  \node[draw,circle,inner sep=1pt,fill,label=left:{$a_0$}] (a0) at (-1,-.5) {};
  \node[draw,circle,inner sep=1pt,fill,label=left:{$a_1$}] (a1) at (-1,-1.5) {};
  \node[draw,circle,inner sep=1pt,fill,label=left:{$a_2$}] (a2) at (-1,-2.5) {};
  \node[] (x) at (-1,-4) {};
  \draw[dotted] (a2) -- (x);
  \node[draw,circle,inner sep=1pt,fill,label=left:{$\bot$}] (0) at (-1,-4.5) {};
  \node[draw,circle,inner sep=1pt,fill,label=right:{$b_1$}] (b1) at (0,-1) {};
  \node[draw,circle,inner sep=1pt,fill,label=right:{$b_2$}] (b2) at (0,-2) {};
  \node[] (y) at (0,-3.5) {};
  \draw[dotted] (b2) -- (y);
  \draw[thin] (1) -- (a0);
  \draw[thin] (1) -- (b1);
  \draw[thin] (b1) -- (a1);
  \draw[thin] (b1) -- (b2);
  \draw[thin] (b2) -- (a2);
  \draw[thin] (a0) -- (a1);
  \draw[thin] (a1) -- (a2);
  \end{tikzpicture}
  \qquad\qquad
  \begin{tikzpicture}
  \node[draw,circle,inner sep=1pt,fill,label=right:{$y$}] (y) at (0,-2) {};
  \node[draw,circle,inner sep=1pt,fill,label=left:{$x_2$}] (x0) at (-1,-2.5) {};
  \node[draw,circle,inner sep=1pt,fill,label=left:{$x_1$}] (x1) at (-1,-3.5) {};
  \node[draw,circle,inner sep=1pt,fill,label=left:{$x_0$}] (x2) at (-1,-4.5) {};
  \node[draw,circle,inner sep=1pt,fill,label=left:{$x$}] (x) at (-1,0) {};
  \node[] (x') at (-1,-1) {};
  \draw[dotted] (x0) -- (x');
  \draw[thin] (x2) -- (x1);
  \draw[thin] (x1) -- (x0);
  \draw[thin] (y) -- (x);
  \end{tikzpicture}
  \end{center}
  \caption{The lattice $L$ from  Example~\ref{ex:[X,Y]notspectral} and its dual space $X$.}
  \label{fig:nonprincipal prime}
  \end{figure}
  
  \begin{example}\label{ex:[X,Y]notspectral}
  Let $L$ be $\bot\oplus(\bN^{op}\times 2)$ and $X$ its dual space, both as depicted in Figure~\ref{fig:nonprincipal prime}.
  Here $x_i={\uparrow}a_i$ is a principal prime filter for each $i\in\bN$, while $x=L\setminus \{\bot\}$ and $y=\{b_i\mid i\in\bN\}$. Further, we let $Y$ be the finite spectral space depicted  in Figure~\ref{fig:Y}. 
  Note that $[X,Y]$ is not compact since we have the following infinite cover, that clearly cannot have a finite subcover:
  \begin{align*}
  [X,Y]=(X\to{\uparrow}y_1)\cup(X\to{\uparrow}y_2)&\cup\bigcup_{i\in\bN}\big[(\widehat{a}_0\to{\uparrow}y_1)\cap(\widehat{b}_i\to{\uparrow}y_2) \big]\\
  &\cup\bigcup_{i\in\bN}\big[(\widehat{a}_0\to{\uparrow}y_2)\cap(\widehat{b}_i\to{\uparrow}y_1) \big].
  \end{align*}
  This equality is saying that a continuous function $f$ from $X$ to $Y$ either misses $y_1$ or it misses $y_2$ and otherwise $x_0$ gets sent to $y_1$ or $y_2$ and $y$ gets sent to the other. In these last cases only finitely many of the $x_i$ take the same value as $x_0$ and this means that there exists $n \in \bN$ such that the direct image of $\wa_n$ under $f$ is $\{y_0\}$.
  Compare also to the proof of Lemma~\ref{lem:comp-open-base}.
  \end{example}
  
  \begin{figure}[!tbp]
  \begin{center}
  \begin{tikzpicture}
  \node[draw,circle,inner sep=1pt,fill,label=right:{$y_2$}] (y2) at (1,-2) {};
  \node[draw,circle,inner sep=1pt,fill,label=left:{$y_1$}] (y1) at (-1,-2) {};
  \node[draw,circle,inner sep=1pt,fill,label=left:{$y_0$}] (y0) at (0,0) {};
  \draw[thin] (y0) -- (y1);
  \draw[thin] (y0) -- (y2);
  \end{tikzpicture}
  \end{center}
  \caption{The space $Y$ from Example~\ref{ex:[X,Y]notspectral}.}
  \label{fig:Y}
  \end{figure}
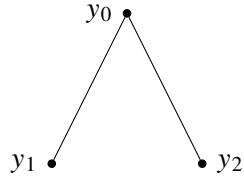
  
\ourexercises

  \begin{ourexercise}\label{exe:upw-comp-upper-viet}
    Let $M$ be a distributive lattice with dual Priestley space $(X, \rho, \leq)$ and let $(\mathcal{V}^{\uparrow}(X),\tau)$ be the upper Vietoris space of the spectral space $(X, \rho^{\uparrow})$. This exercise makes a link between the calculation of the dual space of the lattice $F_{\Box}(M)$ and the duality for unary operators $\Box$ of Section~\ref{sec:unaryopduality}.
    \begin{enumerate}
  \item Prove that $\tau^\partial$, the co-compact dual topology of $\tau$, is generated by the base consisting of finite unions of sets of the form $(\Box a)^c$, for $a \in M$.
  \item Conclude that the Priestley dual space of $F_\Box(M)$ is order-homeomorphic to $(\mathcal{V}^{\uparrow}(X), \tau^p, \leq)$, where $\tau^p$ is the patch topology $\tau \vee \tau^\partial$, and $\leq$ is the \emph{inclusion} order on closed down-sets.
  \item Explain how the result from the preceding item, together with Theorem~\ref{thm:unaryboxduality}, show that upward Priestley compatible relations $R \subseteq X \times Y$ are in bijection with continuous order-preserving functions $f \colon X \to \mathcal{V}(Y)$.
    \end{enumerate}
  \end{ourexercise}
  
  \begin{ourexercise}\phantomsection\label{exe:upperviet-spectral}
  \begin{enumerate}
  \item Let $1$ be the one element space. Show that $\mathcal V^{\uparrow}(1)$ is homeomorphic to the Sierpinski space $\mathbb{S}$.
  \item Let $X$ be a spectral space. Show that $\mathcal V^{\uparrow}(X)$  is homeomorphic to $[X^\partial, \mathbb{S}]$ equipped with the compact-open topology.
  
  \end{enumerate}
  \end{ourexercise}
  
    \begin{ourexercise}\label{ex:compact-open topology}
  Let $X$ and $Y$ be topological spaces. Recall that $[X,Y]$ denotes the set of continuous functions from $X$ to $Y$.
  \begin{enumerate}
  \item Give an example to show that the generating family
  \[
  \{K\to U\mid K\in\KS(X), U\in\Omega(Y)\}
  \]
  may not be closed under intersection or union.
  \item \label{itm:alg-props-arrow} Let $x\in X$, $K,K_1,K_2\in\KS(X)$, $U,U_1,U_2\in\Omega(Y)$, $\cU\subseteq\Omega(Y)$ a directed family, and $\cC\subseteq\KS(X)$ a filtering family. Prove:
  \begin{enumerate}[label=(\arabic*)]
  \item \label{itm:arrow-of-imptype1} $(K_1\to U)\cap(K_2\to U)=(K_1\cup K_2)\to U$;
  \item $(K\to U_1)\cap(K\to U_2)=K\to(U_1\cap U_2)$;
  \item $({\uparrow}x\to U_1)\cup({\uparrow}x\to U_2)={\uparrow}x\to(U_1\cup U_2)$;
  \item $K\to(\bigcup\cU)=\bigcup\{K\to W\mid W\in\cU\}$;
  \item If $X$ is a well-filtered space then $(\bigcap\cC)\to U=\bigcup\{K\to U\mid K\in\cC\}$.
  \end{enumerate}
  \item \label{spec-order-pointwise} Suppose that $X$ is a locally compact space, that is, for any open set $U$ and $x \in U$, there exists a compact set $K$ such that $x \in K \subseteq U$. Prove that the specialization order on the space $[X,Y]$ coincides with the pointwise ordering.
  \item \label{compact-open-domain} Suppose $\cG\subseteq\KS(X)$ generates $\KS(X)$ as a join semilattice and that $\cB$ is a base for $Y$ which is closed under finite unions. Show that
  \[
  \{K\to U \mid K\in\cG, U\in\cB\} 
  \]
  is a subbase for the compact-open topology on $[X,Y]$.

  \end{enumerate}
  \end{ourexercise}
  \begin{ourexercise}\label{ex:XinVietorisX}
  Let $M$ be a distributive lattice and $Y$ its Stone dual space. For every $y \in Y$, define $\eta(y) := {\uparrow} y$, which is a compact-saturated subset of $Y$, and thus an element of $\cV^{\uparrow}(Y)$.
  \begin{enumerate}
    \item Prove that $\eta$ is a spectral embedding.
    \item Show that the map $\eta$ is dual to the quotient of $F_\Box(M)$ under the congruence given by the scheme $\Box(a\vee b)\approx\Box a\vee\Box b$ for $a,b\in M$.
    \item Prove that, under the isomorphism of $\cV^{\uparrow}(Y) \cong \Filt(M)$ of Proposition~\ref{prop:uppervietoris-filter}, the image of $\eta$ is mapped to the set of \emph{prime} filters of $Y$.
    \item Let $L$ be a distributive lattice with Stone dual space $X$. Using the previous items and Lemma~\ref{lem:comp-open-base}, prove that a base for the (not necessarily spectral!) subspace $[X,Y]$ of $[X,\cV^{\uparrow}(Y)]$ is given by the sets $$\wa \To \wb := \{ f \in [X,Y] \ \mid \ f[\wa] \subseteq \wb\}, \text{ for } a \in L, b \in M.$$
  \end{enumerate}
\textit{Note.} In the setting of Priestley spaces, the more delicate issue of a two-sided version of the Vietoris functor was studied in \cite{Palmigiano2004}. Restricted to Boolean spaces, the Vietoris functor was also extensively studied in the context of modal logic \parencite{KupKurVen2004}. This analysis was recently extended to dualities for the category of compact Hausdorff spaces in \cite{BezCarMor2022}.
  \end{ourexercise}
  
\notessec
{The name `Stone space' has also been used for the more restricted class of spaces that we call \emph{Boolean} spaces in this book. We prefer to avoid the name `Stone space' for a class of spaces, to avoid confusion. We do use the terminology `the Stone dual space of a lattice' to refer to the spectral space associated to a distributive lattice through Stone's duality.  Our choice of the terminology `spectral spaces' follows in particular the recent monograph \cite{DicSchTre2019}, the first chapters of which we recommend as useful complementary reading to the material in this chapter.}
\chapter{Domain theory}\label{chap:DomThry}
In this chapter we develop some of the order-topological theory of dcpo's and domains and provide the key duality theoretic elements that were used in \cite{Abr} to solve so-called domain equations.

In denotational semantics, one seeks a category of mathematical objects, so called \emph{denotational types}, whose elements will model programs while the objects themselves model the `types' of the programs. In addition, program constructors should be modelled by functors, so-called \emph{type constructors}. In particular, if $X$ and $Y$ are program types, then we want to be able to form the type which, given a program of type $X$ as input, outputs a program of type $Y$. In the setting of a concrete category, this will mean that we want the set $\Hom(X,Y)$ to be endowed with structure in a natural way that makes it an object of the category whenever $X$ and $Y$ are; a formalization of this idea leads to the categorical definition of a \emphind{Cartesian closed category}.\footnote{We do not need to formally define this notion for our purposes in this chapter, the intuition of `Cartesian closed' meaning `has internal $\Hom$-objects' suffices. We refer to, for example, p. 46 in \cite{MM1992} for the precise definition.}  In this setting, expressions involving the type constructors correspond to program specifications in terms of program constructors and type equations of the form $X\cong F(X)$ correspond to recursive specifications, also known as inductive types. The sought-after category should therefore be closed under a comprehensive set of type constructions, adequate for the needs of semantics of programming languages, as well as under solutions of type equations involving the type constructors. Finally, given the algorithmic nature of computing, it should also admit a reasonable notion of \emph{effective presentability}.

\cite{ScSt71} proposed to look for such a category within the category of dcpos with functions that are continuous with respect to the \emphind{Scott topology}, which we already briefly encountered back in Section~\ref{sec:TopOrd} and will study in more detail starting in Section~\ref{sec:dom} of this chapter. Indeed, Scott originally considered the equation 
\[
X\cong [X,X],
\]
which corresponds to the pure $\lambda$-calculus given by the program constructors of functional abstraction and composition. Scott's solution for this equation looks within the category of dcpos with Scott continuous functions and seeks a dcpo $X$ which is isomorphic to the dcpo of its Scott continuous endomorphisms. 

While the category of dcpos is Cartesian closed, the category as a whole fails to admit a reasonable notion of effective presentation and the general constructions get fairly wild. \emph{Domains}\index{domain}, which we study in Section~\ref{sec:dom}, were introduced precisely as a subcollection of dcpo's with a notion of relatively finite presentability. Further restricting to \emph{algebraic domains}\index{algebraic domain} gives actual finite presentability. Indeed, Scott's solution to $X\cong [X,X]$ is an algebraic domain. However, neither categories of domains or algebraic domains are Cartesian closed, so they do not provide quite the right setting either. A further property of Scott's solution is that it is profinite and thus falls within Stone duality.

Much further work by a large number of researchers confirmed that categorical methods, topology and in particular topological duality are central to the theory, see \cite{ScSt71,Plo76,SmyPlo82,Smyth1983,LW91}.  Abramsky went a step further in his view on duality theoretic methods in this area. Rather than seeing Stone duality and its variants merely as useful technical tools for denotational semantics, he put Stone duality front and center stage. The centrality of duality and the general principles of the theory of domain equations, already put forward in \cite{Abramsky87}, come to their full and clear 
expression in \cite{Abr}, for which he won the IEEE 
Logic in Computer Science Test-of-Time Award in 2007. 

Abramsky casts Stone duality as the mechanism which links programs to their models. Thus Abramsky defines a \emph{program logic}, in which denotational types correspond to theories, and the ensuing Lindenbaum-Tarski algebras of the theories are bounded distributive lattices, whose dual spaces yield the domains as types. The constructors involved in domain equations have duals under Stone duality which are constructors in the program logic, and solutions are obtained as duals of the solutions of the corresponding equation on the lattice side. For this purpose, one needs a Cartesian closed category of domains that are spectral spaces in their Scott topology, and several such had in the meantime been discovered. We will study such \emph{spectral domains}\index{spectral domain}\index{domain!spectral} in Section~\ref{sec:dom-Stone} below, and the more general class of \emph{spectral dcpo's}\index{dcpo!spectral}. While the latter, larger class, is not directly used in domain theory in logical form, its characterization was only previously announced in a conference abstract \parencite{Er2009}, and for completeness' sake we give a proof of it in Section~\ref{sec:dom-Stone}, see Theorem~\ref{thrm:Erne2009}. A reader who wants to get quickly to the domain theory in logical form material may safely skip that part of Section~\ref{sec:dom-Stone}.

In \cite{Abramsky87}, duality is restricted to the so-called Scott domains. These
domains are fairly simple and are closed under many constructors, including
function space, sums, products, and upper and lower powerdomains, but they are
not closed under the convex powerdomain (corresponding to the two-sided Vietoris
construction in topology). In \cite{Abr}, it is shown that his program goes
through for the larger category of \emph{bifinite domains}\index{domain!bifinite}, 
although the mathematics, and especially the duality theory, is much more involved. 
See also the book \cite{Zhang91} where many of the technical aspects of this theory 
were developed concurrently, albeit without the central focus on Stone duality. 

We introduce bifinite domains and study them from a duality-theoretic perspective in Section~\ref{sec:bifinite} below. This category, which was originally introduced by Plotkin, has further closure properties. In particular it is closed under the convex powerdomain construction. As shown in \cite{Smyth83a}, this is optimal when restricted to countably based algebraic domains with least element. The PhD dissertation  \cite{Jung89} completely described all Cartesian closed categories of algebraic domains and, shortly thereafter, Jung introduced two new categories, of $L$- and $FS$-domains, respectively, and showed these to be \emph{the} maximal Cartesian closed categories of continuous domains \parencite{Jung90}. 

This however is far from the end of the story as researchers in the area turned towards the \emph{probabilistic powerdomain}, which forces one into the setting of continuous domains. In order to include the unit interval, which is central in probability theory, one has to go beyond the profinite setting. This leads to the problem of finding a Cartesian closed category of continuous domains which is closed under the probabilistic powerdomain. See \cite{JuTi98} for more on this problem and \cite{GL19} which surveys the state of the art at that time. The notes at the end of this chapter give some more bibliographic details on specific points.

We finish this extended introduction by giving a short outline of the chapter,
in particular giving a road-map for readers who want to get to Section~\ref{sec:DTLF}, in some sense the heart of the chapter, as quickly as possible.

Section~\ref{sec:dom} introduces the notion of a \emph{domain}, which may be seen either as a special kind of poset or as a special kind of sober topological space. In Section~\ref{sec:dom-Stone}, we study the intersection of the class of domains with the class of spectral spaces, and show that, under Stone duality, the domains in this intersection correspond to a very natural class of distributive lattices (Theorem~\ref{thrm:dual-spectral-domains}).  We also establish that these domains are exactly the ones satisfying the properties in the ``2/3 SFP Theorem'' of Plotkin, see \cite[Proposition 4.2.17]{AbJu94}. 

Each of these first two sections also contains an excursion to a purely duality-theoretic result; namely, in the first section, we show that the $\Omega$-$\mathrm{pt}$ duality of Chapter~\ref{chap:Omega-Pt} restricts to the so-called \emph{Hoffmann-Lawson duality} between completely distributive lattices and domains viewed as spaces (Theorem~\ref{thrm:Domain-CDL-duality}); in the second section, in Theorem~\ref{thrm:Erne2009} we characterize the dcpo's that are spectral spaces in their Scott topology, a theorem first announced in \cite{Er2009}. These two results are not directly important for the domain-theoretic applications in the later two sections of the chapter and the rather intricate proofs may be skipped by a reader who wants to get to the duality-theoretic view on bifinite domains and domain equations as fast as possible. 

In Section~\ref{sec:bifinite}, we introduce and study \emph{bifinite} domains, which form the category of domains in which we will constructively solve domain equations, and which moreover are interesting objects for a duality theorist because of their self-dual nature.  We finish the chapter in Section~\ref{sec:DTLF} by showing how the duality-theoretic analysis of the function space construction done in the previous chapter in Section~\ref{sec:funcspace} specializes in the case of bifinite domains, and how this allows one to reconstruct solutions to Scott's equation $X \cong [X, X]$.

\section{Domains and Hoffmann-Lawson duality}\label{sec:dom}

In a state space of computations, we may think of order as `further specification'. That is, two points $p$ and $q$ satisfy $p\leq q$ provided $q$ is a further specified state of the computation than the state $p$. A bottom element may then be thought of as the state of no specification at all, and often domain theorists will consider posets with bottom (called pointed posets) as a convenient set-up. A top element, on the other hand, would further specify all computations in the space. Since one usually would want to consider states leading toward different computations, it is most common in domain theory to consider posets without a top. When a top is added, it is often thought of as the inconsistent, overspecified, state. 
Further, since a complex program which computes its output by finite approximation may be seen as the directed join of the finite approximations, the posets considered in domain theory are often assumed to be closed under suprema of directed sets. The final ingredient is that any point should be obtainable by directed supremum from `finite' or at least `relatively finite' points below it. The pertinent mathematical notions abstracting these ideas are as follows.

 \begin{definition}\label{def:dcpo}
 A poset $P$ is a \emphind{dcpo} (\emphind{directedly complete partial order}) provided every directed subset of $P$ has a supremum in $P$. 
 Let $P$ and $Q$ be dcpos, and $f\colon P\to Q$ a map. Then $f$ is said to \emphind{preserve directed joins} provided $f$ is order preserving and, for every directed set $D$, we have
 \[ f\Big(\bigvee D\Big) = \bigvee f[D]. \]
Note that since $f$ is order preserving, $f[D]$ is directed whenever $D$ is
(see Exercise~\ref{exer:dcpoScott}) and thus, both suprema in this equation
exist. We denote by $\dcpo$ the non-full subcategory of $\POS$ consisting of
dcpos with as morphisms those order-preserving maps which preserve directed
joins. 
\end{definition}

\begin{remark}\label{rem:dcpo-morphism}
Here, we have included in the property of  $f$ preserving directed joins that $f$ is order preserving as this is the most natural setting in which to consider this concept. This reflects the fact that we consider $\dcpo$ as a subcategory of $\POS$. To avoid confusion, every time we state that a function preserves directed joins, we will make sure that it is also clear that this means that it must be order preserving.
\end{remark}
 
  While this category is given entirely in order theoretic terms, as we will see, it is actually isomorphic to the full subcategory of $\TopCat$ consisting of (the underlying sets of) dcpos equipped with the Scott topology of the dcpo. One significant and interesting feature of dcpos is that viewing them as posets, they are first-order structures. This allows access to a class of topological spaces (non first-order) based on first-order structures. However, being closed under directed joins is not a first-order property.
 Recall from Section~\ref{sec:TopOrd} that the \emphind{Scott topology} on a poset $(P,\leq)$ is the topology $\sigma(P,\leq)$ consisting of those sets $U\subseteq P$ which are up-sets in the order and which are inaccessible by directed joins. That is, an up-set $U\subseteq P$ is \emphind{Scott open} provided whenever $D\subseteq P$ is directed and $\bigvee D\in U$ then $D \cap U\neq\emptyset$. Equivalently, we have that a set $C\subseteq P$ is  \emphind{Scott closed} if, and only if, $C$ is a down-set, and whenever $D\subseteq P$ is directed and $D \subseteq C$ then $\bigvee D\in C$. That is, a down-set $C$ is Scott closed if, and only if, $C$ is closed under directed joins. The specialization order of the Scott topology on a poset $P$ is equal to the original order on $P$ (see Exercise~\ref{exer:dcpoScott}.\ref{itm:specord-scott-is-orig}).
 A function $f\colon P\to Q$ between dcpos is said to be \emphind{Scott continuous} provided it is continuous with respect to the Scott topologies on $P$ and $Q$.

 \begin{remark}\label{rem:order-yoga}
     In this chapter, as in the modern literature on domain theory, open sets are \emph{up-sets}. As also noted in Remark~\ref{rem:order-choice}, this clashes somewhat with the choice in Priestley duality theory, which we adhered to up until the previous chapter in this book, that represents a distributive lattice $L$ as the clopen \emph{down-sets} of its Priestley dual space $X$. We already encountered a similar clash in Chapter~\ref{chap:Omega-Pt}: recall from Theorem~\ref{thm:Stone-isom-Priestley} that the lattice of compact-opens of a spectral space $(X,\rho)$ is isomorphic to the the lattice of clopen down-sets of $(X, \rho^p, \leq)$, where $\leq$ is the opposite of the specialization order $\leq_{\rho}$; see also Figure~\ref{fig:dl-priest-spec-diagram} in the previous chapter. Note that another way of saying this is that if $L$ is a distributive lattice with Priestley space $(X,\tau,\leq)$, then the associated spectral space $(X,\tau^{\uparrow})$ of open up-sets has $L^\op$ as its lattice of compact-open sets.

  Using up-sets or down-sets to represent lattice elements is ultimately an arbitrary choice, and there are valid arguments for either choice. Being flexible about this type of `order yoga' is a somewhat cumbersome, but  necessary, part of every duality theorist's life. As a consequence, certain order-theoretic arguments in this chapter may sometimes look `upside down' from those in the earlier chapters of this book. We made the choice to adhere to the conventions from domain theory in this chapter and we will take care to warn the reader throughout the chapter when certain orders are the opposite of the orders used when we discussed Priestley duality.

  As a general remark, for later use in this chapter, we note already that, since Priestley duality is order reversing on morphisms, see Proposition~\ref{prop:priestley-2-eq}, the duality we consider here is order {\it preserving} on morphisms. That is, if $f,g \colon X \rightrightarrows Y$ are two spectral maps between spectral spaces, and $K$ and $L$ are the lattices of compact-opens of $X$ and $Y$, respectively, with $f^*, g^* \colon L \rightrightarrows K$ the dual lattice homomorphisms, then $f \leq g$ in the pointwise order on spectral maps with respect to the specialization order on $Y$, if, and only if, $f^* \leq g^*$ in the pointwise order on lattice homomorphisms with respect to the inclusion order on $K$.
\end{remark}

 \begin{proposition}\label{prop:Scott-dcpo}
  Let $P$ and $Q$ be dcpos, and $f\colon P\to Q$ a map. The following conditions are equivalent:
 \begin{enumerate}[label=(\roman*)]
 \item $f$ is order preserving and preserves directed joins;
 \item $f$ is Scott continuous.
 \end{enumerate}
 \end{proposition}

 \begin{proof}
 To show that (i) implies (ii), suppose $f$ is order preserving and preserves directed joins and let $U\subseteq Q$ be Scott open. First, since $f$ is order preserving, $f^{-1}(U)$ is an up-set. Now suppose $D\subseteq P$ is directed and $\bigvee D\in f^{-1}(U)$. Then $f[D]$ is directed in $Q$ by Exercise~\ref{exer:dcpoScott}.\ref{itm:op-pres-dir} and, as $f$ preserves directed joins, it follows that $\bigvee f[D]=f(\bigvee D)\in U$. Since $U$ is Scott open, it follows that there is $d\in D$ with $f(d)\in U$ and thus $d\in f^{-1}(U)$ and we have proved that $f^{-1}(U)$ is Scott open.

 For the converse, suppose $f$ is Scott continuous.  By Exercise~\ref{exer:cont-implies-op} in Chapter~\ref{chap:TopOrd}, any continuous function is order preserving with respect to the specialization order. Now, the specialization order of a dcpo in the Scott topology is simply the order of the dcpo (see Exercise~\ref{exer:dcpoScott}.\ref{itm:specord-scott-is-orig}). So $f$ is order preserving for the original orders on $P$ and $Q$, respectively.

 Now suppose $D\subseteq P$ is directed. Since $f$ is order preserving, $f[D]$ is directed and thus $\bigvee f[D]$ exists in $Q$. Also since $f$ is order preserving, it follows that $\bigvee f[D]\leq f(\bigvee D)$. Now consider $C={\downarrow}(\bigvee f[D])$. Since it is a principal down-set, it is closed in the Scott topology. By  (ii), $f^{-1}(C)$ is closed in the Scott topology on $P$. Also note that by the definition of $C$ and as $f$ is order preserving we have $D\subseteq f^{-1}(C)$. Therefore $\bigvee D\in f^{-1}(C)$ or, equivalently, $f(\bigvee D)\in C$. That is, $f(\bigvee D)\leq \bigvee f[D]$ as required.
 \end{proof}

 \begin{corollary}\label{cor:dcpoisotop}
 The category $\dcpo$ is isomorphic to a full subcategory of the category 
 $\TopCat$ of topological spaces with continuous maps.
\end{corollary}

Domain theorists often implicitly and harmlessly switch between the two 
perspectives that are provided by Corollary~\ref{cor:dcpoisotop}. However,
this practice, combined with the practice of suppressing the structure in the denotation 
of a mathematical structure may lead to complications: there are topological
spaces $X$ which are dcpo's in their specialization order, even though
the original topology on $X$ is not the Scott topology. To be able to state 
our results in this section, it will sometimes be convenient to keep the distinction 
between the two perspectives clear; to do so, when 
$(P, \leq_P)$ is a dcpo, we refer to $(P, \sigma(P,\leq_P))$ as \emph{the space of the 
dcpo $P$}. \index{dcpo!space of a}

We also note that, despite the isomorphism of categories in 
Corollary~\ref{cor:dcpoisotop}, unexpected things can
happen in switching perspectives if one is not careful; for example, the Scott
topology on the Cartesian product of two dcpo's with the pointwise ordering need
not be equal to the product topology of the Scott topologies on the two dcpo's,
see \cite[Exercise~II-4.26]{Getc2003}. From a categorical point of view, this
means that, if $X$ and $Y$ are dcpo's, then $X\times Y$ is the product of $X$ 
and $Y$ in the ambient category of posets but it is not necessarily the product 
in the ambient category of topological spaces. 
When we restrict to the full subcategory of $\dcpo$ on the objects that are 
continuous, also known as \emph{domains}, see Definition~\ref{def:cmp-wb-cont-alg} 
below, this problem does not occur and finite products of Scott topologies equal the 
Scott topology of the topological products, see \cite[Corollary~II-4.14]{Getc2003}. 

 The notions of dcpo and of Scott continuity also arise naturally from a purely 
 topological point of view, or more specifically, from the point of view of the 
 $\Omega$-$\Pt$ duality, as follows.

 \begin{proposition}\label{prop:sober-vs-dcpo}
 If $(X,\tau)$ is a sober topological space, then $(X,\leq_\tau)$ is a dcpo and $\tau\subseteq\sigma((X,\leq_\tau))$.
\end{proposition}

Thus, any sober topological space $X$ comes with a naturally associated dcpo, that we call \emph{the dcpo of the sober space $X$}. The above proposition then says that the Scott topology of the dcpo of a sober space is always finer than (that is, has at least as many open sets as)  the original sober topology itself. \index{sober space!dcpo of a}

 Before giving a proof of this proposition, we note that a proof via $\Omega$-$\Pt$ duality was outlined for spaces of the form $\Pt(L)$ where $L$ is any frame in Exercise~\ref{exer:sober-and-order}. Essentially, it is a direct consequence of the fact that the specialization order of  $\Pt(L)$ is the inclusion order on completely prime filters and that these are closed under directed unions. By the $\Omega$-$\Pt$ duality, the sober spaces are, up to homeomorphism, precisely the spaces of the form $\Pt(L)$. Thus this proves the proposition. Here  we give a direct proof, not invoking the duality.

 \begin{proof}[Proof of Proposition~\ref{prop:sober-vs-dcpo}]
 Suppose $(X,\tau)$ is sober. It suffices to prove the following:

 {\bf Claim.} For any closed $F \subseteq X$ and for any $D \subseteq F$ directed in the specialization order $\leq_\tau$, the supremum $\bigvee D$ exists and belongs to $F$.

 Indeed, this claim implies the proposition, for the following reasons. First of all, if $D$ is an arbitrary directed subset of $X$, then choosing $F=X$ we just get that $\bigvee D$ exists. Further, since a subset of a poset is Scott closed if, and only if, it is a down-set and is closed under directed joins, we will in fact have shown that every $\tau$-closed set (which is necessarily a down-set in the specialization order) is Scott closed.
 Thus $\tau$ is contained in the Scott topology as claimed.

 We now prove the claim. Let $F$ be closed and $D \subseteq F$ directed in $\leq_\tau$. Define %
 \[
 \cF=\{U\in\tau\mid U\cap D\neq\emptyset\}.
 \]
 It is easy to check that $\cF$ is a completely prime filter; we only show that $\mathcal{F}$ is closed under binary intersection and leave the other parts to the reader. If $U,V\in\cF$, then there are $p,q\in D$ with $p\in U$ and $q\in V$. Since $D$ is directed, there is some $r\in D$ which is above both $p$ and $q$. Now since open sets are always up-sets in the specialization order, it follows that $r\in U\cap V$ and thus $U\cap V\in\cF$ thus showing that $\cF$ is closed under binary intersection.

 Now since $\cF$ is a completely prime filter and $X$ is sober, it follows from the definition of sobriety (Definition~\ref{def:sober}) that there exists $x\in X$ so that $\cF=\cN(x)$. Since every open neighborhood of $x$ intersects $D\subseteq F$ and thus $F$ and since $F$ is closed, it follows that $x\in F$.
 
 Finally, we show that $\bigvee D=x$. Let $y\in D$. We need to show that $y\leq_\tau x$. Let $U\in\tau$ with $y\in U$. Then $D\cap U\neq\emptyset$ and thus $U\in\cF=\cN(x)$. That is, $x\in U$ and thus we have shown that every open containing $y$ contains $x$. That is, $y\leq_\tau x$ as required. On the other hand, suppose $z\in X$ is an upper bound for $D$. We show that $x\leq_\tau z$. Let $U\in\tau$ with $x\in U$, then $U\in \cN(x)=\cF$ and thus $D\cap U\neq\emptyset$. Let $y\in D\cap U$. Since $z$ is an upper bound of $D$, $y\leq_\tau z$ and thus $z\in U$. That is, $x\leq_\tau z$ as required.
 \end{proof}

 \begin{remark}
 The space of a dcpo does not need to be sober, see \cite{Joh81}. In fact, understanding the frames of Scott open sets of dcpos seems a difficult problem, see \cite{HJX}, in which it is shown that dcpos are not determined up to isomorphism by their closed set lattices (which are of course isomorphic to the order dual of the frames of Scott open sets). This will not be a problem for us here, because we will soon restrict to a class of dcpos all of whose spaces are sober, see Proposition~\ref{prop:dom-is-sober} below.
 \end{remark}

 \begin{example}
 Let $X$ be any set. Then $(X,=)$ is a dcpo and
 \[
 \sigma(X,=)=\alpha(X,=)=\delta(X,=)=\cP(X).
 \]
 That is, relative to the trivial order, the Scott topology is equal to the Alexandrov topology and these are equal to the discrete topology on $X$.

 In particular, this example shows that the space of the dcpo of any $T_1$ topological space is discrete. Thus, the inclusion in Proposition~\ref{prop:sober-vs-dcpo} may very well be strict.
 \end{example}
 \begin{example}
 The unit interval in its usual order $([0,1],\leq)$ is a complete lattice, so in particular a dcpo. Its Scott topology is generated by the half-open intervals, that is, the sets of the form $(a,1]=\{x\in[0,1]\mid a<x\}$, for $a \in [0,1]$, and is thus equal to the upper topology; the resulting topological space is in fact a stably compact space (see~Section~\ref{sec:comp-ord-sp}). The reader may verify that the corresponding compact ordered space carries the usual compact Hausdorff topology of the unit interval inherited from the usual topology on the real line and the order is the usual order inherited from the reals.
 \end{example}

The following notions are of fundamental importance in computer science applications of dcpos.
\nl{$\Kel(X)$}{the set of compact elements of a dcpo $X$}{}
 \begin{definition} \label{def:cmp-wb-cont-alg}
 An element $k$ in a dcpo $X$  is said to be a \emphind{compact element} provided for all directed $D\subseteq X$ with $k\leq\bigvee D$ there is $d\in D$ with $k\leq d$. We will denote the set of compact elements of $X$ by $\Kel(X)$. Compact elements are sometimes called \emph{finite elements}\index{finite element} in the literature.

\nl{$x \wayb y$}{the way-below relation between elements $x$ and $y$ of a dcpo}{}
\nl{$\waydown y$}{the set of elements way below $y$}{}
\nl{$\wayup y$}{the set of elements way above $y$}{}
 Let $x,y\in X$. We say that $x$ is \emphind{way below} $y$ and write $x\wayb y$ provided that, for all directed $D\subseteq X$, if $y\leq\bigvee D$, then there is $d\in D$ with $x\leq d$. Further we denote by $\waydown y$ the set of all elements of $X$ that are way below $y$. That is,
 \[
   \waydown y=\{x\in X\mid x\wayb y\},
 \]
 and $\wayup y$ is defined similarly.

 We call a dcpo $X$ {\it a}  \emphind{domain}, also known as a \emphind{continuous dcpo},  provided  each element of $X$ is the directed join of the elements way below it. More explicitly, for a dcpo $X$ to be a domain, for any $x \in X$, the set $\waydown x$ must be directed, and its supremum must be $x$.  We call $X$ an \emphind{algebraic dcpo} or an \emphind{algebraic domain} provided  each element of $X$ is the directed join of the \emph{compact} elements below it.
 \end{definition}

As the nomenclature ``way below'' suggests, $x \wayb y$ implies $x \leq y$ for any elements $x, y$ of a dcpo $X$. Also, an element $x$ in  a dcpo $X$ is compact if, and only if, $x \wayb x$ (see Exercise~\ref{exer:way-below}).
We note that, in the definition of a domain, it suffices to assume that each element $x$ of the dcpo is the directed join of some directed subset of the elements way below $x$; every element of the domain will then in fact be equal to the directed join of the set $\waydown x$ (see Exercise~\ref{exer:weak-cont}). This fact can be quite useful in proofs, as we will see, for example, in the proof of Proposition~\ref{prop:CD-implies-cont}.

We highlight an alternative characterization of \emph{algebraic
domains}\index{domain!algebraic} that will be particularly important later in
this chapter. For a poset $P$, we call an \emphind{order ideal} of $P$ a
down-set that is directed; note that this definition generalizes the notion of
\emph{ideal} for a lattice $L$, see also the remarks following
Theorem~\ref{thrm:Omega-Point-duality}. We denote by $\Idl(P)$ the collection of
ideals of $P$, ordered by inclusion. Now, $\Idl(P)$ is always an algebraic
domain, and an algebraic domain $X$ is always isomorphic to $\Idl(\Kel(X))$;
thus, a domain $X$ is algebraic if, and only if, it is isomorphic to $\Idl(\Kel(X))$. Exercise~\ref{exer:algebraic-domain} asks you to prove this equivalence, via some other equivalent characterizations of algebraic domains.
  \begin{example}
For any set $X$, all elements of the dcpo $(X,=)$  are compact and thus it is an algebraic domain.
 \end{example}
 \begin{example}
  If $(P, \leq)$ is any \emph{finite} poset, then it is an algebraic domain;
  indeed, any directed subset of $P$ is also finite, and therefore contains a
  maximum element. From this, it follows that $x \wayb y$ if, and only if, $x \leq y$, so all elements are compact.
 \end{example}
 \begin{example}
In $([0,1],\leq)$, we have $x\wayb y$ if, and only if, $x=0$ or $x<y$. It follows that $0$ is the only compact element is and that the unit interval is a domain which is not algebraic.
 \end{example}
  \begin{example}
Let $X$ be a set. The partial order $(\cP(X),\subseteq)$ is an algebraic domain: we have $x\wayb y$ if, and only if, $x \subseteq y$ and $x$ is finite (see Exercise~\ref{exe:domainexamples}). Thus, all finite subsets are compact elements of $\cP(X)$, and each subset is the directed union of its finite subsets.
 \end{example}
 \begin{example}
  Let $(X,\tau)$ be a topological space. An open subset $K \subseteq X$ is compact in the topological sense if, and only if, $K$ is a compact element of the frame $\Omega X$. That is, $\Kel(\Omega(X))$ is equal to the set of compact-open subsets of $X$, which was denoted by $\KO(X)$ in Chapter~\ref{chap:Omega-Pt}.
 \end{example}
  \begin{example}\label{exa:partialfunctiondomain}
Let  ${\it Part}(X)$ denote the poset of partial functions on a set $X$ with the order given by $f\leq g$ if, and only if, $g$ extends $f$. That is, ${\dom}(f)\subseteq\dom(g)$ and $f(x)=g(x)$ for all $x\in\dom(f)$. Then again the finite partial functions (that is, those with finite domain) are the compact elements, and ${\it Part}(X)$ is an algebraic domain which is not a lattice (see Exercise~\ref{exe:domainexamples}). It does however have the property that all principal down-sets are complete lattices; in fact, they are complete and atomic Boolean algebras. Note that the total functions on $X$ are the maximal elements of ${\it Part}(X)$.
\end{example}

\begin{definition}
A binary relation $R$ on a set $X$ is called \emphind{interpolating} if $R \subseteq R \cdot R$.
\end{definition}
Note that a relation $R$ is \emph{idempotent}\index{idempotent relation}, i.e,
$R \cdot R = R$ if, and only if, $R$ is transitive and interpolating.
\begin{lemma}\label{lem:contdcpo}
 Let $X$ be a domain. Then $\wayb$ is transitive, interpolating and, for each $x\in X$, the set $\wayup x$ is Scott open. Furthermore, $U\subseteq X$ is Scott open if, and only if,
 \[
 U=\bigcup_{x\in U}\wayup x.
 \]
\end{lemma}

\begin{proof}
Note that $\wayb$ is transitive on any dcpo (see Exercise~\ref{exer:way-below}). To see that $\wayb$ is interpolating, let $x, y \in X$ be such that $x \wayb y$. Since $X$ is continuous, we have 
\begin{align*}
x \wayb y=\bigvee\waydown y&=\bigvee\{z \mid z \in X, z\wayb y\}\\
                                                   &=\bigvee\Big\{\bigvee \waydown z \mid z \in X, z\wayb y\Big\}\\
                                                    &=\bigvee Z,
\end{align*}
where $Z := \bigcup\{ \waydown z \mid z \in X, z\wayb y\}$. 
Since $\waydown y$ is directed, the collection $\{\waydown z \mid z \in X, z\wayb y\}$ is a directed family of sets in the inclusion order. Combining this with the fact that each $\waydown z$ for $z\wayb y$ is directed, it follows that the set $Z$ is directed. Thus, there exist $z, z'\in X$ with $x\leq z'\wayb z\wayb y$. It follows that $x(\wayb\cdot\wayb) y$ (see Exercise~\ref{exer:way-below}).

To show that $\wayup x$ is Scott open, let $D\subseteq X$ be directed with $x\wayb\bigvee D$. Then there is $x'\in X$ with $x\wayb x'\wayb\bigvee D$. Since $x'\wayb\bigvee D$, there is $d\in D$ with $x'\leq d$. Now since $x\wayb x'$ it follows that $x\wayb d$ and thus $\wayup x \cap D \neq \emptyset$, showing that  $\wayup x$ is Scott open.
Consequently, any set $U$ satisfying $U=\bigcup_{x\in U}\wayup x$ is also Scott open.

Finally, let $U\subseteq X$ be Scott open. Clearly, since $U$ is an up-set, in particular $\bigcup_{x\in U}\wayup x\subseteq U$. For the reverse inclusion, let $y\in U$. As $X$ is a domain, $y=\bigvee\waydown y$. Since $U$ is Scott open it follows that there is $x\wayb y$ with $x\in U$. Thus $y\in \bigcup_{x\in U}\wayup x$ as required.
\end{proof}

 \begin{proposition}\label{prop:dom-is-sober}
 Any domain is sober in its Scott topology.
 \end{proposition}

 \begin{proof}
 Let $X$ be a domain and suppose $F\subseteq X$ is a join-irreducible closed set for the Scott topology. We need to show that $F = {\downarrow} x$ for some $x \in X$.   Consider the set
 \begin{align*}
  \waydown F &=\{y\in X\mid \exists x\in F\ y\wayb x\}\\
                                      &=\{y\in X\mid \wayup y\cap F\neq\emptyset\}.
 \end{align*}
 We will first show that $\waydown F$ is directed, so that it has a supremum, $x$, and we will then show that $F = {\downarrow} x$. Note first that $F$, being join-irreducible, is non-empty. Let $x\in F$, then, as $X$ is a domain,  $\waydown x$ is directed and in particular non-empty. Thus $\waydown F$ is not empty. Let $y_1,y_2\in\waydown F$ be arbitrary. Define $F_i=(\wayup y_i)^c$ for $i=1,2$. Then, since $y_i \in \waydown F$, there is $x_i\in F$ with $y_i\wayb x_i$ and thus $F\not\subseteq F_i$ for both $i=1,2$.  Now since $F_1$ and $F_2$ are closed in the Scott topology by Lemma~\ref{lem:contdcpo} and since $F$ is a join-irreducible closed set, it follows that
 \[
 F\not\subseteq F_1\cup F_2=(\wayup y_1\cap\wayup y_2)^c.
 \]
 Thus, pick $z\in F\cap\wayup y_1\cap\wayup y_2$. Since $\wayup y_i$ is open for both $i=1,2$, it follows that $\wayup y_1\cap\wayup y_2$ is open. Therefore, since $\bigvee\waydown z=z\in\wayup y_1\cap\wayup y_2$, there exists $z'\wayb z$ with $z'\in\wayup y_1\cap\wayup y_2$. Thus, in particular, $y_1\leq z'$ and $y_2\leq z'$. Also, since $z'\wayb z\in F$ it follows that $z'\in\waydown F$ and we have shown that $\waydown F$ is directed. Let $x=\bigvee\waydown F$. We show that $F={\downarrow} x$. First, since $\waydown F\subseteq {\downarrow} F \subseteq F$ and since $F$ is closed, it follows that $x=\bigvee\waydown F\in F$. Since $F$ is a down-set, we get the containment ${\downarrow}x\subseteq F$. Conversely, for any $x'\in F$, we have $\waydown x'\subseteq \waydown F$ and thus
 \[
 x'=\bigvee \waydown x'\leq \bigvee\waydown F=x.
 \]
 That is, $F\subseteq{\downarrow} x$ as required.
 \end{proof}

\subsection*{Hoffmann-Lawson duality}
 As a consequence of Corollary~\ref{cor:dcpoisotop} and
 Proposition~\ref{prop:dom-is-sober}, the category $\CONT$ of domains with order-preserving maps that preserve
 directed joins (or equivalently, Scott continuous functions) is
 (isomorphic to) a full subcategory of the category $\SOB$ of sober topological
 spaces.  Thus it is natural to ask which category of spatial frames is dual in
 the $\Omega$-$\Pt$ duality to the category  $\CONT$. The appropriate frames
 are the \emph{completely distributive} ones and the resulting duality is the
 \emphind{Hoffmann-Lawson duality}, as we will prove now.
As noted in the introduction to this chapter, while this is
an interesting excursion to a classical result in domain theory, it is not
directly needed for the applications in Section~\ref{sec:DTLF}, and can be
skipped by readers wanting to get to those applications as quickly as possible.

  Since all suprema exist in a frame, any frame is in fact a complete lattice (see Exercise~\ref{exe:complattsuff} and Exercise~\ref{exer:infimum-in-frame}). This means that we can consider the notion of \emph{complete distributivity} for frames.

 \begin{definition}
 A complete lattice $L$ is said to be \emphind{completely distributive} provided, for any family $\{A_i\}_{i\in I}$ of subsets of $L$, we have
 \[
 \bigwedge_{i\in I}\bigvee A_i =\bigvee\biggl\{\bigwedge \im(\Phi)\mid \Phi\colon I\to L \text{ such that }\Phi(i)\in A_i \text{ for each }i\in I\biggr\}.
 \]
 For a family $\{A_i\}_{i\in I}$ of sets, the functions $\Phi\colon I\to \bigcup_{i\in I}A_i$ such that $\Phi(i)\in A_i$ for each $i\in I$ are called \emphind{choice functions} on $\{A_i\}_{i\in I}$.
 \end{definition}
We will show that $\Omega$-$\Pt$ duality further restricts to a duality between domains and completely distributive lattices.
\begin{theorem}\label{thrm:Domain-CDL-duality}
The $\Omega$-$\Pt$ duality between spatial frames and sober spaces cuts down to a duality between the category $\CDFrame$ of completely distributive complete lattices with frame homomorphisms and the category of domains, $\CONT$.
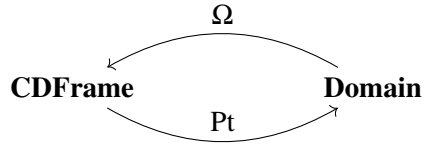
\[
\begin{tikzpicture}
  \node (A) {{\CDFrame}};
  \node (B) [node distance=4cm, right of=A] { {\CONT}};
  \draw[->,bend right] (A) to node [above,midway] {$\Pt$} (B);
  \draw[->, bend right] (B) to node [above,midway]{$\Omega$} (A);
\end{tikzpicture}
\]
\end{theorem}

The remainder of this section is dedicated to proving Theorem~\ref{thrm:Domain-CDL-duality}. We first give an outline of the proof.

 \begin{proof}[Outline of proof of Theorem~\ref{thrm:Domain-CDL-duality}]
We have that:
\begin{enumerate}
\item If $X$ is a domain, then $X$ is sober (Proposition~\ref{prop:dom-is-sober});
\item If $X$ is a domain, then $\Omega(X)$ is completely distributive (Proposition~\ref{prop:cont-implies-compldist});
\item If $L$ is completely distributive, then $L$ is spatial (Corollary~\ref{cor:CD-spatial});
\item  If $L$ is completely distributive, then $\Pt(L)$ is a domain and its topology is the Scott topology (Proposition~\ref{prop:CD-implies-cont} and Lemma~\ref{lem:CompDist-implies-PointScott});
\end{enumerate}
This is precisely what is needed to show that the $\Omega$-$\Pt$ duality cuts down to a duality between domains and completely distributive lattices.
 \end{proof}

Towards proving Proposition~\ref{prop:cont-implies-compldist}, we begin by studying in slightly more detail the class of completely distributive lattices. First, we give a simpler description of complete distributivity.

 \begin{proposition}\label{prop:infdist}
 A complete lattice $L$ is completely distributive if, and only if, for all families $(D_i)_{i\in I}$ of down-sets of $L$, we have
  \[
 \bigwedge_{i\in I}\bigvee D_i =\bigvee\Big(\bigcap_{i\in I} D_i\Big)
 \]
 \end{proposition}

 \begin{proof}
 First note that for any family of sets $(A_i)_{i \in I}$, $ \bigwedge_{i\in I}\bigvee A_i = \bigwedge_{i\in I}\bigvee {\downarrow} A_i $ and, for any choice function $\Phi\colon I\to \bigcup_{i\in I}{\downarrow} A_i$, there is a choice function $\Phi'\colon I\to \bigcup_{i\in I} A_i$ with $\Phi(i) \leq\Phi'(i)$ for every $i \in I$, and thus with $\bigwedge\im(\Phi)\leq\bigwedge\im(\Phi')$. It follows that it suffices to consider families of down-sets in the definition of complete distributivity.

 Now the proposition follows if we can show that, for any family $\{D_i\}_{i \in I}$ of down-sets of $L$, we have
 \begin{equation}\label{eq:intersection-choice-functions}
   \bigcap_{i\in I} D_i=\Big\{\bigwedge \im(\Phi)\mid \Phi \text{ is a choice function for }\{D_i\}_{i\in I}\Big\}.
 \end{equation}
 Let $a\in \bigcap_{i\in I} D_i$. Then the constant function $\Phi_a$ on $I$ given by $\Phi(i)=a$ for all $i\in I$ is a choice function for $\{D_i\}_{i\in I}$ and $\bigwedge \im(\Phi_a)=a$. This proves the left-to-right containment in (\ref{eq:intersection-choice-functions}).

 On the other hand, if $\Phi$ is any choice function for $\{D_i\}_{i\in I}$, then for any $i \in I$, the element $\bigwedge\im(\Phi)$ is below $\Phi(i) \in D_i$, and is thus in $\bigcap_{i\in I} D_i$.
 \end{proof}

 In order to understand completely distributive lattices, we introduce a strengthened variant of the way below relation, which is a relativized version of complete join-primeness in the same way that the way below relation is a relativized version of compactness for elements of a dcpo.


 \begin{definition}\label{def:<<<}
 Let $L$ be a complete lattice and $a,b\in L$. We write $b\wwb a$, or equivalently $b\in\threeheaddownarrow a$, provided that, for any subset $S$ of $L$, if $a \leq \bigvee S$, then there exists $s \in S$ such that $b \leq s$.
 \end{definition}
An element $a$ in a complete lattice is called \emphind{completely join prime} if $a \wwb a$. Note that $a \wwb b$ clearly implies $a \wayb b$. 
The following theorem is due to Raney \parencite{Raney53}.
\begin{theorem}[Raney's Theorem]\label{thrm:Raney53} 
A complete lattice $L$ is completely distributive if, and only if,
\begin{equation}\label{eq:join-of-relprimebelow}
\text{ for every } a \in L, \qquad\qquad a=\bigvee \threeheaddownarrow a.
\end{equation}
\end{theorem}
  \begin{proof}
 Suppose (\ref{eq:join-of-relprimebelow}) holds and let $\{D_i\}_{i\in I}$ be a collection of down-sets of $L$. We define 
   \[
    d := \bigwedge_{i\in I}\bigvee D_i, \quad d' := \bigvee\big(\bigcap_{i\in I} D_i\big).
 \]
 We want to show that $d = d'$. In fact, for each $i \in I$, since $\bigcap_{i\in I} D_i\subseteq D_i$ we have $d' \leq \bigvee D_i$, so $d' \leq d$. We now show the reverse inequality, $d \leq d'$. To this end, let $b\wwb\bigwedge_{i\in I}\bigvee D_i $. Then $b\wwb\bigvee D_i$ for each $i\in I$. By the definition of $\wwb$ there is, for each $i\in I$, an element $d_i\in D_i$ with $b\leq d_i$. Since each $D_i$ is a down-set, it follows that $b\in D_i$ for each $i\in I$ and thus $b\in \bigcap_{i\in I} D_i$. That is,
\[
\threeheaddownarrow d \subseteq  \bigcap_{i\in I} D_i
\]
and thus by (\ref{eq:join-of-relprimebelow}) we have the desired inequality
$d = \bigvee {\threeheaddownarrow d} \leq d'.$

For the converse implication, suppose $L$ is completely distributive. Note that, for any $a\in L$, we have
\[
a=\bigwedge\{\bigvee S\mid S\in\cD(L)\text{ and }a\leq\bigvee S\}.
\]
This is because $a$ is clearly a lower bound of the collection we are taking the infimum of, and $S := {\downarrow a}$ is a down-set of $L$ with $a=\bigvee S$. Now applying complete distributivity we obtain
\[
a=\bigvee(\bigcap\{S\mid S\in\cD(L)\text{ and }a\leq\bigvee S\}).
\]
Finally observe that
\[
\bigcap\{S\mid S\in\cD(L)\text{ and }a\leq\bigvee S\}=\threeheaddownarrow a.\qedhere
\]
 \end{proof}

 \begin{corollary}
 Any completely distributive lattice is continuous.
 \end{corollary}

 This is a consequence of the relation between $\wayb$ and $\wwb$ (see Exercise~\ref{exer:<<<}).

  \begin{corollary}\label{cor:wwb-int}
 The relation $\wwb$ is interpolating on a completely distributive lattice.
 \end{corollary}

 The proof of Corollary~\ref{cor:wwb-int} is a simpler version of the corresponding fact for $\wayb$ on a domain (again see Exercise~\ref{exer:<<<}). For a more substantial consequence of Raney's result, we prove one direction of Hoffmann-Lawson duality.

 \begin{proposition}\label{prop:cont-implies-compldist}
Let $X$ be a domain. The frame $\sigma(X)$ of Scott open subsets of $X$ is completely distributive.
 \end{proposition}

 \begin{proof}
 Let $U\in\sigma(X)$. Recall Lemma~\ref{lem:contdcpo},
  in which we showed that, if $X$ is a domain then,
 for each $x\in X$, the set $\wayup x$ is Scott
 open and $U=\bigcup_{x\in U}\wayup x$. Thus, by Theorem~\ref{thrm:Raney53}, we
 may conclude that $\sigma(X)$ is completely distributive if we can show that
 $\wayup x\wwb U$ for each $x\in U$. To this end, note that if $x\in U$ and
 $U=\bigcup_{i\in I} U_i$, then there is $i\in I$ with $x\in U_i$ and thus $\wayup x
 \subseteq {\uparrow} x\subseteq U_i$ and indeed $\wayup x\wwb U$.
 \end{proof}

To prove that, conversely, any completely distributive lattice $L$ is isomorphic to one of the form $\sigma(X)$ for some domain $X$, we need a few lemmas.

 \begin{lemma}\label{lem:principal-sep}
 Let $L$ be a frame and let $F, G$ be completely prime filters in $L$. If there is $a\in L$ with $F\subseteq {\uparrow} a\subseteq G$, then $F\wayb G$.
 \end{lemma}
 \begin{proof}
If $\{F_i\}_{i\in I}$ is a directed family of completely prime filters and $G=\bigcup_{i\in I} F_i$ then $a\in F_i$ for some $i\in I$ and thus $F\subseteq {\uparrow} a\subseteq F_i$.
 \end{proof}

 \begin{lemma}\label{lem:wwb-CDLat}
 Let $L$ be a completely distributive lattice and $a,b\in L$, then $a\wwb b$ if
 and only if there is $F\in\Pt(L)$ with ${\uparrow} b\subseteq F\subseteq{\uparrow} a$.
 \end{lemma}
 \begin{proof}
Suppose $a\wwb b$. Since $\wwb$ is interpolating by Corollary~\ref{cor:wwb-int},
there is a sequence $\{b_n\}_{n\in\bN}$ with $a\wwb b_{n+1}\wwb b_n\wwb b$ for
all $n\in\bN$. One may easily verify that $F=\bigcup_{n\in \bN}{\uparrow} b_n$ is a
completely prime filter (see Exercise~\ref{exer:<<<-chain}) and that the element
of $\Pt(L)$ given by $F$ satisfies the required property.

For the converse, notice that if ${\uparrow} b\subseteq F\subseteq{\uparrow} a$
and $a\leq\bigvee S$, then there is $s\in S$ with $s\in F$ and thus $b\leq s$. That
is, $a\wwb b$.
  \end{proof}

  \begin{corollary}\label{cor:CD-spatial}
    Any completely distributive lattice is a spatial frame.
\end{corollary}

\begin{proof}
Let $L$ be a completely distributive lattice. If $a,b\in L$ with $a\nleq b$ then, by Raney's Theorem, there is $c\in L$ with
$c\wwb a$ but {\it not} $c \wwb b$. By Lemma~\ref{lem:wwb-CDLat}
pick a completely prime filter $F$ of $L$ with ${\uparrow}a\subseteq F\subseteq {\uparrow}c$ but ${\uparrow}b\not\subseteq F$. That is $a\in F$ and $b\not\in F$ and thus $L$ is spatial.
\end{proof}

Recall from Definition~\ref{def:Pt-functor} that we adopt the `neutral space' notation for $\Pt(L)$, analogously to what we did in earlier chapters for the Priestley dual space of a distributive lattice. That is, we consider the set underlying $\Pt(L)$ as a fresh set of `names' that is in bijection with the set of completely prime filters of $L$. We denote elements of $\Pt(L)$ by $x, y, z, \dots$, and the corresponding completely prime filters respectively by $F_x, F_y, F_z$, et cetera.
\begin{proposition}\label{prop:CD-implies-cont}
 If $L$ is a completely distributive lattice then the dcpo $\Pt(L)$ is a domain.
 \end{proposition}
 \begin{proof}
Let $y\in\Pt(L)$ and $b\in F_y$. Since ${\bigvee}{\threeheaddownarrow b}=b\in F_y$, there is $a\in F_y$ with $a\wwb b$.
By Lemma~\ref{lem:wwb-CDLat}, it follows that there is $x_b\in\Pt(L)$ with
${\uparrow} b\subseteq F_{x_b}\subseteq{\uparrow} a$. Therefore $F_{x_b}\subseteq{\uparrow} a\subseteq F_y$
and, by Lemma~\ref{lem:principal-sep}, we have $x_b\wayb y$. Now notice that
\begin{align*}
F_y=\bigcup\{{\uparrow}b\mid b\in F_y\}&\subseteq\bigcup\{F_{x_b}\mid b\in F_y\}\\
&\subseteq\bigcup\{F_x\mid \exists a\in L\quad F_x\subseteq{\uparrow}a\subseteq F_y\}\subseteq F_y.
\end{align*}
Consequently
\[
y=\bigvee\{x\in\Pt(L)\mid \exists a\in L\quad F_x\subseteq{\uparrow}a\subseteq F_y\}.
\]
By Lemma~\ref{lem:principal-sep}, it follows that $y$ is the join of a subfamily of $\waydown y$.
If we can show that this collection is directed, we can conclude that $\Pt(L)$ is a domain by Exercise~\ref{exer:weak-cont}.
To this end suppose $F_x\subseteq{\uparrow}a\subseteq F_y$ and $F_z\subseteq{\uparrow}b\subseteq F_y$.
Since both $a$ and $b$ belong to $F_y$, we have $c=a\wedge b\in F_y$. Now since $c={\bigvee}{\threeheaddownarrow c}$,
there is $d\wwb c$ with $d\in F_y$. By Lemma~\ref{lem:wwb-CDLat}, it follows that there is $s\in\Pt(L)$ with
${\uparrow} c\subseteq F_s\subseteq{\uparrow} d$. Since $a$ and $b$ are both above $c$, it follows that
$F_x\subseteq{\uparrow}a\subseteq{\uparrow} c\subseteq F_s$ and
$F_z\subseteq{\uparrow}b\subseteq{\uparrow} c\subseteq F_s$.  Finally, since $d\in F_y$, we have  $F_s\subseteq{\uparrow} d\subseteq F_y$ and thus $x,z\wayb s \wayb y$. We conclude that $\{x\in\Pt(L)\mid \exists a\in L\quad F_x\subseteq{\uparrow}a\subseteq F_y\}$ is indeed directed.
 \end{proof}

\begin{corollary}\label{cor:CD-wayb}
If $L$ is a completely distributive lattice and $x,y\in\Pt(L)$, then $x\wayb y$ if, and only if,
there exists an $a\in L$ with $F_x\subseteq{\uparrow} a\subseteq F_y$.
\end{corollary}

\begin{proof}
By Lemma~\ref{lem:principal-sep} the `if' part is always true. For the converse, note that,
in the proof of Proposition~\ref{prop:CD-implies-cont}, we proved that if $y\in\Pt(L)$ then
\[
y=\bigvee\{z\in\Pt(L)\mid \exists a\in L\quad F_z\subseteq{\uparrow}a\subseteq F_y\}
\]
and that this join is directed. Thus, if $x\wayb y$, then there is $z\in\Pt(L)$ and $a\in L$
so that $F_z\subseteq{\uparrow}a\subseteq F_y$ and $x\leq z$, but then $F_x\subseteq{\uparrow}a\subseteq F_y$.
\end{proof}

\begin{lemma}\label{lem:CompDist-implies-PointScott}
If $L$ is a completely distributive lattice then the topology on $\Pt(L)$ is equal to the Scott topology of its specialization order.
\end{lemma}

 \begin{proof}
Since any space in the image of the $\Pt$ functor is sober, it follows by
 Proposition~\ref{prop:sober-vs-dcpo} that the topology of $\Pt(L)$ is
 contained in the Scott topology.

For the converse, let $U\subseteq\Pt(L)$ be Scott open. Recall from Definition~\ref{def:Pt-functor} that the topology of $\Pt(L)$ consists of the sets $\widehat{a}=\{x\in\Pt(L)\mid a\in F_x\}$  for $a\in L$. Let
\[
a=\bigvee\Big\{\bigwedge F_x\mid x\in U\Big\}\ .
\]
We claim that $U=\widehat{a}$. First, if $a\in F_y$, then as $F_y$ is completely prime, there is $x\in U$ with $\bigwedge F_x\in F_y$. It follows that $F_x\subseteq F_y$
or equivalently that $x\leq y$ and thus $y\in U$. That is, $\widehat{a}\subseteq U$.
On the other hand, if $x\in U$ then, since $\Pt(L)$ is a domain, there is $y\in U$ with
$y\wayb x$. By Corollary~\ref{cor:CD-wayb}, there is $b\in L$ with $F_y\subseteq{\uparrow}b\subseteq F_x$. It follows that $b\leq \bigwedge F_y\leq a$. Finally, since $b\in F_x$ also $a\in F_x$ and $x\in\widehat{a}$.
 \end{proof}

This concludes the last piece of the proof of  Theorem~\ref{thrm:Domain-CDL-duality}.

 \ourexercises
  \begin{ourexercise}\label{exer:dcpoScott}
 Let $P$ and $Q$ be dcpos, and $f\colon P\to Q$ a map.
 \begin{enumerate}
 \item\label{itm:specord-scott-is-orig} Show that the specialization order of the topological space $(P,\sigma(P,\leq))$ is the original order $\leq$ on $P$;
\item\label{itm:op-pres-dir} Show that if $f$ is order preserving and $D\subseteq P$ is directed, then so is $f(D)$;
 \item\label{itm:dcpo-in-top} Show that
 \begin{align*}
 S\colon\dcpo  & \to\TopCat\\
 (P,\leq) &\mapsto (P,\sigma(P,\leq))\\
     f       & \mapsto f
     \end{align*}
 is a functor whose image is a full subcategory of $\TopCat$. Further show that if $S(P)$ and $S(Q)$ are homeomorphic as topological spaces, then $P$ and $Q$ are isomorphic as posets.
 \end{enumerate}
 \end{ourexercise}

 \begin{ourexercise}
 Let $(X,\leq)$ be a partially ordered set.
 \begin{enumerate}
 \item Show that the Scott topology on the poset $\cD=\Down(X,\leq)$ is equal to its upper topology;
 \item Show that these topologies are generated by the principal up-sets $U_x=\{V\in\cD\mid x\in V\}$;
 \item Show that the resulting space is stably compact;
 \item Show that the associated compact ordered space is a Priestley space.
 \end{enumerate}
 \end{ourexercise}

  \begin{ourexercise}\label{exer:way-below}
 Let $(X,\leq)$ be a dcpo and $x\in X$.
 \begin{enumerate}
 \item Show that $x$ is way below itself if, and only if, $x$ is a compact element of $X$;
 \item Show that the way below relation is contained in the order relation;
 \item Show that $\leq\cdot\wayb\cdot\leq\,=\,\wayb$.
 \item Show that the set $\waydown x$ is closed under any existing binary joins. That is, if
 $y,z\in\waydown x$ and $y\vee z$ exists in $X$, then $y\vee z\in\waydown x$.
 \end{enumerate}
 \end{ourexercise}

\begin{ourexercise}\label{exe:domainexamples}
  Let $X$ be a set.
  \begin{enumerate}
  \item Prove that in the dcpo $(\cP(X), \subseteq)$, $x \wayb y$ if, and only if, $x \subseteq y$ and $x$ is finite.
  \item Draw the Hasse diagram of the partial order ${\it Part}(X)$ in case $X$ has two elements.
  \item Characterize the way below relation in the dcpo $({\it Part}(X), \leq)$ of Example~\ref{exa:partialfunctiondomain}.
  \item Show that $({\it Part}(X), \leq)$ is not a lattice when $X$ has at least two elements.
  \item Show that, for any $f \in {\it Part}(X)$, the sub-poset ${\downarrow} f$ is isomorphic to $(\cP(\dom(f)), \subseteq)$.
  \end{enumerate}
\end{ourexercise}

  \begin{ourexercise}\label{exer:weak-cont}
  Prove that if $X$ is a dcpo in which every element is the join of some directed set of elements way below it, then $X$ is a domain.
   \end{ourexercise}

  \begin{ourexercise}\label{exer:compldist}
 Prove each of the following statements:
 \begin{enumerate}
  \item A complete lattice $L$ is completely distributive if, and only if, its order dual is completely distributive;
  \item There are frames which are not completely distributive;
  \item There are completely distributive frames with no completely join-irreducible or completely meet-irreducible elements;
 \item A complete Boolean algebra is completely distributive if, and only if, it is atomic (see Exercise~\ref{exer:non-spatial}).
 \end{enumerate}
\end{ourexercise}

\begin{ourexercise}\label{exer:compldist-jirr}
  An element $j$ of a complete lattice $L$ is called \emphind{completely join-irreducible} if, for any $S \subseteq L$, if $j = \bigvee S$, then $j \in S$.
  \begin{enumerate}
    \item Show that a completely join-prime element is always completely join irreducible.
    \item Show that, if $L$ is a frame, then any completely join-irreducible element is completely join prime.
  \end{enumerate}
\end{ourexercise}
\begin{ourexercise}\label{exer:frames-compdist}
  This exercise compares the complete distributivity law to the strictly weaker frame distributivity law.
  \begin{enumerate}
  \item Prove that a complete lattice $L$ is a frame if, and only if, the  complete distributive law holds in $L$ for finite collections, that is, if for any \emph{finite} index set $I$ and any collection $\{A_i\}_{i \in I}$ of subsets of $L$, 
 \[
 \bigwedge_{i\in I}\bigvee A_i =\bigvee\Big\{\bigwedge \im(\Phi)\mid \Phi\colon I\to L \text{ a choice function for } \{A_i\}_{i \in I}\Big\}.
\]
\item Formulate a `directed distributive law' (DDL) such that a complete lattice is completely distributive if, and only if, it is a frame that satisfies (DDL). 
  \end{enumerate}
\end{ourexercise}

    \begin{ourexercise}\label{exer:<<<}
  Let $L$ be a completely distributive lattice and $a\in L$.
   \begin{enumerate}
  \item Show that any finite join of a subset of $\threeheaddownarrow a$ is in $\waydown a$. \hint{Use Exercise~\ref{exer:way-below}.d.}
  \item Show that $L$ is continuous.
  \item Show that $\wwb$ is interpolating. \hint{The proof is similar to the corresponding proof for $\wayb$ given in Lemma~\ref{lem:contdcpo}.}
  \end{enumerate}
   \end{ourexercise}

\begin{ourexercise}\label{exer:<<<-chain}
      Let $L$ be a frame and $S\subseteq L$ which is filtering with respect to $\wwb$. That is, if $a,b\in S$ then there is $c\in S$ with $c\wwb a$ and $c\wwb b$. Show that $F={\uparrow} S$ is a completely prime filter of $L$.
       \end{ourexercise}

\begin{ourexercise}\label{exer:algebraic-domain}
Let $X$ be a domain. Recall that $\Kel(X)$ denotes the poset of compact elements of $X$, and for a poset $P$, $\Idl(P)$ denotes the collection of order ideals of $P$, that is, down-sets that are up-directed. Show that the following conditions on $X$ are equivalent:
\begin{enumerate}[label=(\roman*)]
\item $X$ is algebraic;
\item $X\cong\Idl(\Kel(X))$;
\item The frame of opens of $X$ is isomorphic to $\Up(\Kel(X))$;
\item The frame of opens of $X$ is isomorphic to $\Up(P)$ for some poset $P$;
\item $X$ is isomorphic to $\Idl(P)$ for some poset $P$.
\end{enumerate}
\end{ourexercise}

\begin{ourexercise}\label{exer:cat-algebraic-domain}
Let $P$ and $Q$ be posets. A relation $R\subseteq P\times Q$ is called
\emphind{approximable} provided $\geq_P{\cdot}\, R\,{\cdot}\geq_Q\,=R$ and
$R[p]=\{q\in G\mid pRq\}$ is directed for each $p\in P$. Show that the category
{\Alg} of algebraic domains is equivalent to the category ${\POS_{\rm approx}}$ of
posets with approximable relations via the functors which send an algebraic
domain to its poset of compact elements and a poset to the algebraic domain of
its order ideals. 

\emph{Note.} Exercise~\ref{exer:join-approx} is related to this exercise in the case where the Scott topology on the algebraic domain is spectral, a situation that we will study in Section~\ref{sec:dom-Stone}. The notion of bases and approximable relations between them allows one to generalize the equivalence between ${\POS_{\rm approx}}$ and algebraic domains in this exercise to continuous domains in general \parencite[Section~2.2.6]{AbJu94}.
\end{ourexercise}

\begin{ourexercise}\label{exer:algebraic-domain-COs}
Let $X$ be an algebraic domain. Show that $U\subseteq X$ is compact-open in $X$ if, and only if, $U={\uparrow} F$  for some finite $F\subseteq\Kel(X)$. Conclude that the compact-open subsets of an algebraic domain form a base for the Scott topology on $X$ which is closed under finite unions.
\end{ourexercise}

\section{Dcpos and domains that are spectral}\label{sec:dom-Stone}

It is particularly interesting to know which dcpos and domains are such that their associated space is spectral. The topological spaces that are simultaneously spectral \emph{and} Scott topologies of their specialization order have been characterized by Marcel Ern\'e. While these spaces originally went by the name \emphind{hyperspectral spaces}, we will call them \emphind{spectral dcpo}'s here, in light of our convention that a dcpo is always equipped with its Scott topology. The first result of this section, Theorem~\ref{thrm:Erne2009}, characterizes spectral dcpo's as the coherent sober spaces which have a base of open finitely generated up-sets. While this result is not immediately needed for domain theory in logical form, we believe it answers a very natural question about the relationship between spectral spaces and dcpo's. We therefore include a full proof of it here, as far as we are aware for the first time in writing, since the statement was so far only available in a conference abstract \parencite{Er2009}. Again, a reader who wants to get to the applications of duality to domain theory in logical form as quickly as possible may just read the relevant definitions and statement of Theorem~\ref{thrm:Erne2009}. 

An important class of spaces that is directly relevant to the rest of this chapter and that we begin to study on p.~\pageref{subsec:specdomain} is obtained by restricting the class of spectral dcpo's to those that are also \emph{domains}; we refer to these objects as \emph{spectral domains}\index{spectral domain} here. These spaces have very nice descriptions both as spaces and as posets, which we will study in detail in the second part of this section. We will show in Theorem~\ref{thrm:spectral-algebraic-domain} that, in topological terms, these domains are precisely the spectral spaces whose compact-opens are finite unions of union-irreducible compact-opens, and in order theoretic terms they are precisely the completions under directed joins of so-called finitely \MUB-complete posets, see Corollary~\ref{cor:2/3sfp} below.

In Section~\ref{sec:bifinite}, we will introduce a further subcategory of the category of spectral domains, namely the \emph{bifinite domains}. These are central to Abramsky's Domain Theory in Logical form as treated in \cite{Abr}. 
Figure~\ref{fig:dcpo-and-spaces} gives an overview of the relevant classes of spaces and their relationships.

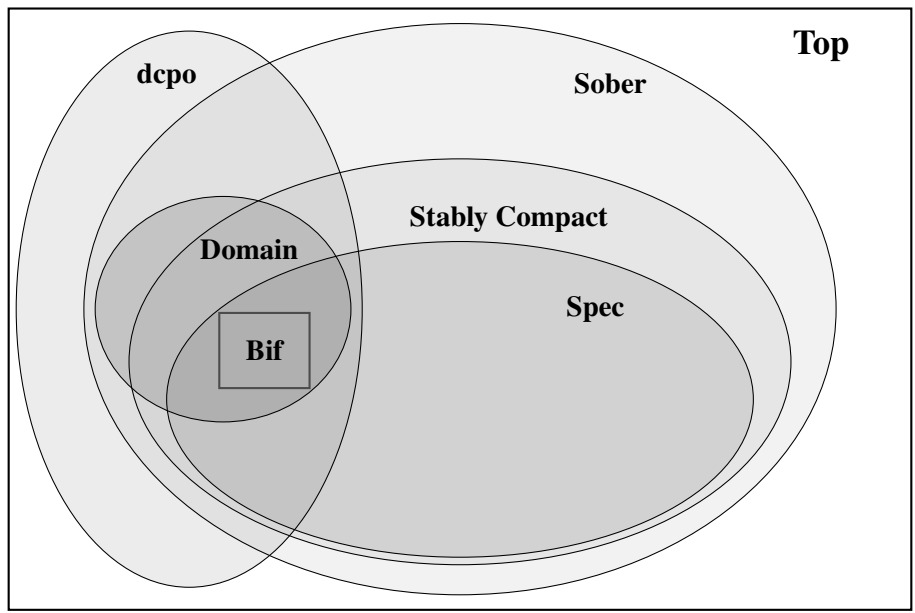
\begin{figure}
 \def\Sober{(90:0cm) ellipse (5cm and 3.8cm)}
 \def\StaCom{(90:-.7cm) ellipse (4.4cm and 2.7cm)}
  \def\Dcpo{(180:3.6cm) ellipse (2.3cm and 3.7cm)}
  \def\Spect{(90:-1.2cm) ellipse (3.9cm and 2.1cm)}
  \def\Domain{(180:3.15cm) ellipse (1.7cm and 1.5cm)}
\begin{center}
 \begin{tikzpicture}
\draw[thick] (-6,-4) rectangle (6,4);
    \draw (4.8,3.5) node[scale=1.2] {$\TopCat$};

\draw[thick] (-3.2,-1.05) rectangle (-2,-0.05);

      \begin{scope}[fill opacity=0.1]
        \fill[gray] \Sober;
      \end{scope}

      \begin{scope}[fill opacity=0.15]
        \fill[gray] \Dcpo;
      \end{scope}

      \begin{scope}[fill opacity=0.2]
        \fill[gray] \Spect;
      \end{scope}

     \begin{scope}[fill opacity=0.1]
       \fill[gray] \StaCom;
      \end{scope}

      \begin{scope}[fill opacity=0.3]
        \fill[gray] \Domain;
      \end{scope}

      \draw \Sober;
      \draw\StaCom;
      \draw \Dcpo;
      \draw \Spect;
      \draw \Domain;

      {
        
        \node[text=black] at ( 1.8,0) {\textbf{Spec}};
        \node[text=black] at (-3.9, 3.1) {\textbf{dcpo}};
        \node[text=black] at (2, 3) {\textbf{Sober}};
        \node[align=left] at (-2.8, .8) {\textbf{Domain}};
          }

       { 
        \node[align=left] at (  -2.6, -0.55) {{\bf Bif}};
         \node[align=left] at (  0.65, 1.2) {{\bf Stably Compact}};
        }
    \end{tikzpicture}
    \caption{The various spaces of interest: Within the sober spaces, the stably compact spaces allow a duality between open and compact-saturated sets, and {\Spec} is a special case of this (see Chapter~\ref{chap:Omega-Pt}).
      Within the dcpos, {\CONT} is a particularly important class of spaces allowing `relatively-finite' approximation. We will identify both the spaces in the intersection of {\dcpo} and {\Spec}, namely \emph{spectral dcpos}, and the spaces in the intersection of {\CONT} and {\Spec}, namely \emph{spectral domains}.
      Bifinite domains form a Cartesian closed category, {\bf Bif}, of spaces in the intersection of {\CONT} and {\Spec} that is the setting of Domain Theory in Logical Form as elaborated by Abramsky.}
  \label{fig:dcpo-and-spaces}
  \end{center}
\end{figure}

\subsection*{Spectral dcpos}
Our goal in this subsection is to characterize spectral dcpo's, defined as follows.
\begin{definition}\label{def:hypspec}
  A topological space $(X, \tau)$ is called a \emphind{spectral dcpo} if the topology $\tau$ is spectral, and $\tau$ is equal to the Scott topology of its specialization order.
  \end{definition}

The characterization will make crucial use of the following notion of a finitely generated up-set.
\begin{definition}\label{def:hypcomp}
  Let $X$ be a poset. A \emphind{finitely generated up-set} is an up-set $T \subseteq X$ such that $T = {\uparrow} F$ for some finite subset $F \subseteq X$.

  In a topological space $X$, a finitely generated up-set in the specialization order of $X$ is sometimes called a \emphind{hypercompact} saturated subset or a \emphind{tooth} of $X$, but we will also use the term finitely generated up-set for this concept, where it is understood that we mean finitely generated up-set in the specialization order.
\end{definition}

A finitely generated up-set in a topological space $X$ is clearly compact, since only finitely many opens are needed to cover a finite
set and, since opens are up-sets, they cover the up-set of any set they contain. Importantly for us, in Scott topologies, the converse also holds.

\begin{lemma}\label{lem:Scottopencompact=hypercompact}
Let $X$ be a poset and $V\subseteq X$ a Scott open subset. Then $V$ is compact in the Scott topology on $X$ if, and only
if, $V$ is a finitely generated up-set.
\end{lemma}

\begin{proof}
We only need to show the necessity, as sufficiency holds in any topological space. Any compact-saturated subset of a topological space is the up-set of its minimal points
(see Exercise~\ref{ex:K=upminK}). We show that, for each $x\in\min(V)$, the set
\[
U_x=(V-\min(V))\cup\{x\}
\]
is Scott open in $X$. To this end, let $D$ be a directed subset of $X$ with $\bigvee D\in U_x$.
Then  $\bigvee D\in V$ and thus there is $y\in D\cap V$. If $y\geq\bigvee D$ then
$y=\bigvee D\in U_x$. Otherwise, there is $z\in D$ with $z\nleq y$. Now, since $D$ is directed,
there is $s\in D$ with $y,z\leq s$. In fact, since $z\nleq y$, $y<s$. Thus
$s\in D\cap (V-\min(V))\subseteq U_x$ as required.

It follows that $\{U_x\mid x\in\min(V)\}$ is an open cover of $V$. By compactness it has a finite
subcover. However, this implies that $\min(V)$ is finite and thus $V$ is a finitely generated up-set.
\end{proof}

We note the following immediate corollary, which will be important in the next subsection when we characterize spectral domains.
\begin{corollary}\label{cor:comp-open-teeth}
  Let $X$ be a spectral dcpo.
A subset of $X$ is compact-open if, and only if, it is an finitely generated open up-set. In particular, a finite intersection of finitely generated open up-sets in $X$ is a finitely generated open up-set.
\end{corollary}
\begin{proof}
The first statement holds by Lemma~\ref{lem:Scottopencompact=hypercompact}. The in particular statement follows because the compact-opens in a spectral space are closed under finite intersections.
\end{proof}

Recall from Section~\ref{sec:comp-ord-sp} that a space is called \emphind{coherent} provided it is compact, and the intersection of any two compact-saturated subsets is again compact. Also recall from Definition~\ref{def:StoneSpace} that a topological space is a spectral space provided it is coherent and sober and has a base of compact-open subsets. Our goal in this subsection is to prove the following characterization of spectral dcpo's, first announced in \cite{Er2009}. 
\begin{theorem}\label{thrm:Erne2009}
Let $(X, \tau)$ be a topological space. Then $X$ is a spectral dcpo if, and only if, the topology $\tau$ is coherent, sober, and has a base of finitely generated open up-sets.
\end{theorem}

The proof of Theorem~\ref{thrm:Erne2009} will be a relatively straightforward combination of things we have already seen, with the notable exception of one crucial step, namely Proposition~\ref{prop:Erne-crucial} below. We first prove a simple lemma on spaces having a base of finitely generated open up-sets.

\begin{lemma}\label{lem:compact-open-are-teeth}
Let $X$ be a topological space. If $X$ has a base of finitely generated open up-sets, then every compact-open in $X$ is a finitely generated up-set.
\end{lemma}

\begin{proof}
Let $U\subseteq X$ be compact-open. Since it is open, it is the union of all the
open finitely generated up-sets it contains, and as $U$ is compact, it is the union of finitely many
such finitely generated open up-sets, ${\uparrow}F_1,\dots,{\uparrow}F_n$. But generating up-sets
is union preserving and thus
\[
U={\uparrow}F_1\cup\ \dotsm \ \cup{\uparrow}F_n={\uparrow}(F_1\cup\ \dotsm \ \cup F_n).\qedhere
\]
\end{proof}

We also require \emphind{Rudin's Lemma}, which is often useful in topology. Like the Alexander Subbase Theorem, it requires a non-constructive principle. Our statement and proof is based on \cite[Lemma III-3.3]{Getc80}.

\begin{lemma}[Rudin's Lemma]\label{lem:rudin}
Let $X$ be a poset and $\cF$ a collection of finite non-empty subsets of $X$ with the property that the collection
\[
\{{\uparrow}F\mid F\in\cF\}
\]
is filtering. Then there is a directed set $D\subseteq\bigcup\cF$ so that $D\cap F\neq\emptyset$
for each $F\in\cF$.
\end{lemma}

\begin{proof}
Consider the collection $\cS$ consisting of all subsets $E\subseteq \bigcup\cF$ with the following two properties
\begin{enumerate}
\item $E\cap F\neq\emptyset$  for all $F\in\cF$;
\item For all $F,G\in\cF$, if $G\subseteq{\uparrow} F$ then $E\cap G\subseteq{\uparrow}(E\cap F)$.
\end{enumerate}
Notice that $\bigcup\cF\in\cS$ and thus $\cS$ is non-empty. We show that any chain $\cC$ in
$\cS$ has a lower bound in $\cS$. To this end, let $D=\bigcap\cC$. We
claim that $D$ has the required properties. First we show that $D\in\cS$. Fix $F\in\cF$ then, as
$\cC$ is a chain, so is $\{E\cap F\mid E\in\cC\}$. Also, since $F$ is finite, it is a finite chain, and
by (a), it consists of non-empty subsets of $F$. It follows that it has a minimum
$\emptyset\neq E_0\cap F\subseteq F$, where $E_0\in\cC$. It follows that
\[
D\cap F=\bigcap\{E\cap F\mid E\in\cC\}=E_0\cap F\neq\emptyset.
\]
Now consider $F,G\in\cF$ with $G\subseteq{\uparrow} F$. As argued above, there are
elements $E_1,E_2\in\cC$ so that $D\cap F=E_1\cap F$ and $D\cap G=E_2\cap G$.
Now since $\cC$ is a chain $E_1$ and $E_2$ are comparable, and we let $E$ be the
smaller of the two. Then, by minimality, we have
\[
D\cap F=E_1\cap F=E\cap F\quad\text{ and }\quad D\cap G=E_2\cap G=E\cap G.
\]
Finally, since $E\in\cS$, we have
$D\cap G=E\cap G\subseteq{\uparrow}(E\cap F)={\uparrow}(D\cap F)$ as required.
It follows, by Zorn's Lemma, that $\cS$ contains minimal elements (that is, maximal elements
with respect to the order given by reverse inclusion).

Before concluding, we make the following observation. Let $E\in\cS$ and $x\in X$ with the
following property
\begin{equation}\label{eq:star}
\forall F\in\cF\quad E\cap F\not\subseteq {\uparrow}x.
\end{equation}
We claim that $E'=E\cap({\uparrow}x)^c$ is again in $\cS$. To this end, first note that (\ref{eq:star})
is equivalent to $E\cap F\cap({\uparrow}x)^c\neq\emptyset$ for all $F\in\cF$, so that (a) holds for
$E'$. For (b), first notice that for any $S\subseteq X$ we have
${\uparrow}S\cap({\uparrow}x)^c\subseteq{\uparrow}[S\cap({\uparrow}x)^c]$. Now if
$F,G\in\cF$ with $G\subseteq{\uparrow} F$ then we have
\begin{align*}
E'\cap G=E\cap G\cap({\uparrow}x)^c&\subseteq{\uparrow}(E\cap F)\cap({\uparrow}x)^c\\
							&\subseteq{\uparrow}[E\cap F\cap({\uparrow}x)^c]={\uparrow}(E'\cap F)
\end{align*}

Now let $D\in\cS$ be minimal, then removing elements from $D$ results in a set not in $\cS$,
thus the negation of  (\ref{eq:star}) holds for each $x\in D$. That is,
\begin{equation}%
\forall x\in D\ \exists F\in\cF\quad D\cap F\subseteq {\uparrow}x.
\end{equation}
We show that this implies that $D$ is directed: Let $x,y\in D$. Pick $F,G\in\cF$ with
$D\cap F\subseteq {\uparrow}x$ and $D\cap G\subseteq {\uparrow}y$. Since
$\{{\uparrow}H\mid H\in\cF\}$ is filtering, there is $H\in\cF$ with $H\subseteq{\uparrow}F,{\uparrow}G$.
Finally, as $D\in\cS$, there is $z\in D\cap H$. It follows that $z\in D$ and
\[
z\in H\subseteq{\uparrow}F,{\uparrow}G\subseteq{\uparrow}x,{\uparrow}y
\]
and thus $D$ is directed.
\end{proof}

\begin{proposition}\label{prop:Erne-crucial}
  Let $X$ be a coherent sober space with a base of finitely generated open up-sets. Then any Scott open set in the specialization order of $X$ is open in $X$.
  \end{proposition}
  
\begin{proof}
Let $U\subseteq X$ be Scott open. We want to show that for each
$x\in U$, there is a finitely generated open up-set ${\uparrow}F$ with $x\in {\uparrow}F\subseteq U$.
By contraposition, assume that $x\in X$ is such that for any finite $F \subseteq X$, if ${\uparrow}F$ is open and contains $x$, then $F \not\subseteq U$. We show that $x \not\in U$.
Consider the collection
\[
\cF:=\{F-U\mid F\text{ is finite and } x\in {\uparrow}F\text{ is open}\}.
\]
We show that $\cG=\{{\uparrow}G\mid G\in\cF\}$ is filtering. For this purpose let $F_1,F_2$
be finite subsets of $X$ with $x\in{\uparrow}F_i$ and ${\uparrow} F_i$ open for $i=1,2$. Since ${\uparrow} F_1$ and ${\uparrow} F_2$ are compact-open in the coherent space, so is ${\uparrow} F_1 \cap {\uparrow} F_2$. By Lemma~\ref{lem:compact-open-are-teeth}, there exists a finite subset $F \subseteq X$ such that ${\uparrow} F_1 \cap {\uparrow} F_2 = {\uparrow} F$.
It follows from this equality that $x\in{\uparrow}F$ is open. We show that $F-U
\subseteq{\uparrow}(F_i-U)$ for $i=1,2$. Since ${\uparrow}F\subseteq{\uparrow}F_i$ we have $F-U\subseteq{\uparrow}F_i$. Let $x\in F-U$ and let $y\in F_i$ with $y\leq x$. If $y\in
U$ then $x\in U$ since $U$ is an up-set. So $y\in F_i-U$ and thus $F-U\subseteq{\uparrow}
(F_i-U)$ for $i=1,2$.

Therefore, by Rudin's Lemma, there is a directed subset
\[
D\subseteq \bigcup\cF\subseteq X-U
\]
 with $D\cap(F-U)\neq\emptyset$ for each finite subset $F$ of $X$ with $x\in{\uparrow}
F$ open. Now, on the one hand, since $U$ is Scott open, it follows that $\bigvee D\not\in U$,
and on the other hand, since
\[
\bigvee D\in\bigcap\{{\uparrow}F\mid F\text{ is finite and } x\in {\uparrow}F\text{ is open}\}={\uparrow}x,
\]
it follows that $x\leq\bigvee D$ and thus $x\not\in U$ as required.
\end{proof}

With Proposition~\ref{prop:Erne-crucial} in hand, we can now prove Theorem~\ref{thrm:Erne2009}.
\begin{proof}[Proof of Theorem~\ref{thrm:Erne2009}]
  First, assuming $X$ is a spectral dcpo, the compact-open sets form a base by definition of spectral spaces, but the compact Scott open sets are the finitely generated open up-sets by Lemma~\ref{lem:Scottopencompact=hypercompact}; thus, $X$ has a base of finitely generated open up-sets.

  Now assume that $X$ is coherent, sober, and has a base of finitely generated open up-sets. The specialization order of a sober space is a dcpo by Proposition~\ref{prop:sober-vs-dcpo}, and the topology is spectral because finitely generated up-sets are always compact in the Scott topology. Also, by Proposition~\ref{prop:sober-vs-dcpo}, $\tau$ is contained in the Scott topology of $\leq_\tau$. But $\leq_\tau$ is also contained in $\tau$, by Proposition~\ref{prop:Erne-crucial}.
\end{proof}

\subsection*{Spectral domains}\label{subsec:specdomain}
Having characterized the spectral dcpos, we now restrict our attention to \emph{domains} whose associated space is spectral.
\begin{definition}
Let $X$ be a domain. We call $X$ a \emphind{spectral domain} if the Scott topology on $X$ is spectral.
\end{definition}

If $X$ is a spectral domain, then $X$ is a spectral dcpo, and $X$ is compact, so by Corollary~\ref{cor:comp-open-teeth},  $X$ can only have a finite number of minimal elements. In any poset, it is clear that the union of two finitely generated up-sets is again a finitely generated up-set. On the other hand, it is not true in general that the finitely generated open up-sets of a domain are closed under binary intersection (see Exercise~\ref{ex:teeth}), while this must be true in a spectral dcpos, and thus in particular in spectral domains.

For a general spectral dcpo, it is hard to understand in order-theoretic terms what it means for the finitely generated open up-sets to be closed under binary intersections, as we do not know which finite subsets $F\subseteq X$ generate opens. However, if $X$ is not just a dcpo but a domain, then the description of the compact-open subsets simplifies substantially, as we show now.

\begin{proposition}\label{prop:spectral-domain-is-algebraic}
Let $X$ be a spectral domain. Then $U\subseteq X$ is compact-open if, and only if, $U={\uparrow}F$ for some finite set $F\subseteq\Kel(X)$ of compact elements of $X$. In particular,
$X$ is an algebraic domain.
\end{proposition}

\begin{proof}
For any $x \in \Kel(X)$, we have ${\uparrow} x = \wayup x$, which is a Scott open set by Lemma~\ref{lem:contdcpo}. Thus, for any subset $F \subseteq \Kel(X)$, ${\uparrow} F = \bigcup \{{\uparrow}x \mid x \in F\}$ is open. Since any finitely generated up-set is compact, it follows that ${\uparrow} F$ is compact-open when $F \subseteq \Kel(X)$ is finite.

For the converse, suppose that $U\subseteq X$ is compact-open. Then, by Corollary~\ref{cor:comp-open-teeth}, $U$ is a finitely generated open up-set. Write $F$ for the finite set of minimal elements of $U$. We show that $F \subseteq \Kel(X)$. Let $x \in F$. By Lemma~\ref{lem:contdcpo}, since $U$ is open, there is $y\in U$ with $y\wayb x$ and thus, in particular, $y\leq x$. Since $x$ is minimal in $U$, it follows that $y=x$, so that $x \wayb x$, that is, $x$ is a compact element of $X$.

For the `in particular' statement, let $y \in X$ be arbitrary. We need to prove that $y = \bigvee ({\waydown} y \cap \Kel(X))$. Since $X$ is a domain, we have that $y = \bigvee {\waydown} y$, so it suffices to prove that, for any $y' \wayb y$, there exists a compact element $x \wayb y$ with $y' \leq x$. Let $y' \wayb y$ be arbitrary. Since $\wayb$ is interpolating by Lemma~\ref{lem:contdcpo}, pick $z \in X$ with $y' \wayb z \wayb y$. Since $z \in {\wayup} y'$, which is open, and since $X$ is a spectral space, there is a compact-open set $U$ such that $z \in U \subseteq {\wayup} y'$. Since $U$ is the upward closure of finitely many compact elements, there exists a compact element $x \in \Kel(X) \cap U$ such that $x \leq z$. Since $x \in U$, we have $y' \wayb x$, so in particular $y' \leq x$. Also, $x \leq z \wayb y$, so $x \wayb y$, as required.
\end{proof}

Recall that and element $p$ in a lattice $L$ is said to be \emphind{join prime} provided that, for any finite $F \subseteq L$, $p\leq\bigvee F$ implies that there exists $a\in F$ with $p\leq a$. Further, we say that $L$ \emph{has enough join-primes} provided every element of $L$ is the join of a finite set of join-prime elements.

Also recall that, when $X$ is a spectral space, ``the distributive lattice dual to $X$'' is by definition the lattice $\KO(X)$ of compact-open sets of $X$, also see Figure~\ref{fig:dl-priest-spec-diagram}.

\begin{corollary}\label{cor:dual-spectral-domain}
Let $X$ be a spectral domain. Then the distributive lattice $L$ dual to $X$ is isomorphic to the poset of finitely generated up-sets of $\Kel(X)$ with the set theoretic operations of intersection and union. In particular, $L$ has enough join-primes.
\end{corollary}
\begin{proof}
Consider the function from the poset of finitely generated up-sets of $\Kel(X)$ to $\KO(X)$ that sends $T \subseteq \Kel(X)$ to the up-set ${\uparrow} T$ of $X$. By Proposition~\ref{prop:spectral-domain-is-algebraic}, this is a well-defined surjective map, and it is clearly an order-embedding. The poset of finitely generated up-sets of $\Kel(X)$ clearly has enough join-primes, as each ${\uparrow}x$ with $x\in\Kel(X)$ is join prime, and finitely generated up-sets of $\Kel(X)$ are finite unions of such.
\end{proof}
Corollary~\ref{cor:dual-spectral-domain} shows that, for a spectral domain $X$, the poset $\Kel(X)$ completely determines the lattice $L$ dual to $X$, and thus also $X$ itself. In particular, the poset $\Kel(X)$ is isomorphic to $\cJ(L)^\op$, as can be seen from  Exercise~\ref{exe:fingendown} (also see Exercise~\ref{exe:join-primes-spectral-domain} below).
We now show that the spectral domains are in fact precisely the duals of the distributive lattices with enough join-primes.
Even though this theorem was first proved through the work of domain theorists in the 1980s and 90s, culminating in the duality-theoretic formulation given here, a closely related result coming from a model- and set-theoretic tradition is~\cite[Proposition~2.8]{BPZ2007}.

\begin{theorem}\label{thrm:dual-spectral-domains}
Stone's duality between spectral spaces and distributive lattices restricts to a duality between spectral domains and  distributive lattices with enough join-primes.
\end{theorem}

\begin{proof}
Corollary~\ref{cor:dual-spectral-domain} establishes that the distributive lattice dual to a spectral domain has enough join-primes. It remains to prove that the spaces of the form $\St(L)$, where $L$ is a  distributive lattice with enough join-primes, are precisely the spectral domains. 
Let $L$ be a  distributive lattice with enough join-primes and $X$ its Stone dual. Then $X$ is clearly spectral. We show that it is a domain in its specialization order. Each join-prime $p$ of $L$ gives rise to a principal prime filter $x={\uparrow}p$ and any such element $x$ is a compact element of $X$. Now let $x,y\in X$ with $x\nleq y$, then there is $a\in L$ with $a\in x$ but $a\not\in y$. Since $a$ is a finite join of join-primes and $x$ is a prime filter, there is $p\in \cJ(L)$ with $p\leq a$ and $p\in x$. Since $a\not\in y$ it follows that  $p\not\in y$. That is, the ${\uparrow}p\leq x$ and ${\uparrow}p\nleq y$ and thus each element of $X$ is a join of compact elements. It remains to show that this join is directed. Let $p,q\in \cJ(L)$ and $x\in X$ with $p,q\in x$. Then $p\wedge q\in x$ and thus, as above for $a$, there is $r\in \cJ(L)$ with $r\leq p\wedge q$ and $r\in x$. It follows that ${\uparrow}p,{\uparrow}q\leq{\uparrow}r\leq x$ in $X$ and thus we have shown that $X$ is a spectral domain.
\end{proof}
The duality in Theorem~\ref{thrm:dual-spectral-domains} is centered around the poset $P=\Kel(X)$ associated to the spectral domain $X$. An interesting feature of spectral domains is that there is a \emph{significant overlap} between the space $X$ and its dual lattice $L$, since this subposet $P=\Kel(X)$ of $X$ is isomorphic to the opposite of the poset $\cJ(L)$ of join-prime elements of $L$. It follows that this poset uniquely determines both $X$ and $L$, since, in this case, $X\cong\Idl(P)$ and $L\cong\Downfin(P^\op)\cong \Upfin(P)$, the free join-semilattice on $P^\op$, which, for this very special type of poset $P$, is not just closed under finite unions but also under finite intersections. 
We end this section by giving an order-theoretic characterization of the domains that are spectral in their Scott topology; that is, in Theorem~\ref{thrm:spectral-algebraic-domain} below, we will characterize the spectrality of a domain $X$ entirely in terms of the poset $P = \Kel(X)$. The characterizations of bifinite domains among spectral domains that we will give in the next section also make heavy use of this poset $P$, which lives on both sides of the duality. %

\nl{$\UB(M)$}{the set of upper bounds of a subset $M$ of a poset}{}
\nl{$\MUB(M)$}{the set of minimal upper bounds of a subset $M$ of a poset}{}
\begin{definition}\label{def:MUB}
Let $X$ be a poset, $x\in X$, and $M\subseteq X$. Then we will write $M\leq x$ provided $x$ is a common upper bound of all the elements of $M$. That is, $x$ belongs to the set 
$\UB(M):=\{x\in X\mid \forall m\in M\ m\leq x\}.$
Further, we denote by $\MUB(M)$ the set of all \emphind{minimal upper bounds} of $M$. That is,
\[
\MUB(M):=\min(\UB(M)).
\]
\end{definition}

We are now ready to characterize, in order-theoretic terms, the spectral domains among the algebraic domains.
\begin{theorem}\label{thrm:spectral-algebraic-domain}
Let $X$ be an algebraic domain. Then the following conditions on $X$ are equivalent:
\begin{enumerate}[label=(\roman*)]
\item the Scott topology on $X$ is spectral;
\item for all finite subsets $F\subseteq\Kel(X)$, the following two properties hold:
\begin{enumerate}
\item $\MUB(F)$ is finite and contained in $\Kel(X)$;
\item $\UB(F)={\uparrow}\MUB(F)$.
\end{enumerate}
\item the finitely generated Scott-open up-sets of $X$ are closed under finite intersections.
\end{enumerate}
\end{theorem}

\begin{proof}
First suppose $X$ is a spectral domain and let $F \subseteq \Kel(X)$ be finite. Then ${\uparrow} x=\wayup x$ is compact-open in $X$ for each $x\in F$ and thus the finite intersection,
\[
\UB(F) = \bigcap_{x\in F}{\uparrow} x
\]
is also compact-open.  Thus, by Proposition~\ref{prop:spectral-domain-is-algebraic}, there is a finite antichain $G\subseteq\Kel(X)$ such that $\UB(F)={\uparrow}G$. It follows that
\[
\MUB(F)=\min(\UB(F))=\min({\uparrow}G)=G.
\]
That is, $G = \MUB(F)$ is a finite subset of $\Kel(X)$ and $\UB(F)={\uparrow}\MUB(F)$, as required.

Now suppose that $X$ is an algebraic domain satisfying the two properties in (ii). Recall from Exercise~\ref{exer:algebraic-domain-COs} that, in any algebraic domain, the finitely generated open up-sets are precisely the sets ${\uparrow}G$ with $G\subseteq\Kel(X)$ finite. Thus, by  distributivity, it suffices to show that finite intersections of sets of the form ${\uparrow} x$, with $x \in \Kel(X)$, are finitely generated open up-sets. Let $F\subseteq\Kel(X)$ be finite and $U := \bigcap_{x \in F} {\uparrow} x$. Then, using (ii)(b),
\[
U=\UB(F)={\uparrow}\MUB(F)
\]
and by (ii)(a), $\MUB(F)$ is a finite subset of $\Kel(X)$, so that $U$ is a finitely generated open up-set.

Finally, to see that (iii) implies (i), recall that a topological space $X$ is spectral if, and only if, it is sober and the collection of compact-open sets is a bounded sublattice of $\cP(X)$ which is a base for the topology, see Definition~\ref{def:StoneSpace} and  Exercises~\ref{exer:Stone-spaces} and \ref{exer:sober-vs-wellfiltered}. Any algebraic domain is sober, and by Exercise~\ref{exer:algebraic-domain-COs}, the finitely generated open up-sets form a base closed under finite unions, and any compact-open is a finitely generated open up-set. Now, since finitely generated open up-sets are closed under finite intersections by (iii), it follows that the finite intersection of compact-opens is a finitely generated open up-set, and thus in particular compact again. Thus, the Scott topology on $X$ is spectral.
\end{proof}

Now combining Proposition~\ref{prop:spectral-domain-is-algebraic} and Theorem~\ref{thrm:spectral-algebraic-domain}, we obtain the following order-theoretic characterization of spectral domains, which have also been called ``2/3 bifinite'' or ``2/3 SFP''\index{bifinite!2/3} domains in the literature (see for example \cite[Proposition~4.2.17]{AbJu94}), because they satisfy two of the three properties that characterize bifinite domains, see Corollary~\ref{cor:mub-bifinite} in the next section.

\begin{corollary}\label{cor:2/3sfp}
Let $X$ be a domain. Then the Scott topology on $X$ is spectral if, and only if, $X$ is algebraic, and for all finite subsets $F\subseteq\Kel(X)$ the following two properties hold:
\begin{enumerate}
\item[(i)] $\MUB(F)$ is finite and contained in $\Kel(X)$;
\item[(ii)] $\UB(F)={\uparrow}\MUB(F)$.
\end{enumerate}
\end{corollary}

\begin{remark}\label{rem:MUB-comp}
Posets $P$ satisfying property (ii) for every finite subset $F \subseteq P$ are called \emph{\MUB-complete}\index{MUB-complete} in \cite{AbJu94}, and this property is sometimes also referred to as `property m'. Combining Corollary~\ref{cor:2/3sfp} with Exercise~\ref{exer:algebraic-domain}, it follows that spectral domains are, up to isomorphism, the posets of the form $\Idl(P)$, for $P$ a \MUB-complete poset in which $\MUB(F)$ is finite for any $F \subseteq P$; we will call such posets \emph{finitely \MUB-complete}\index{MUB-complete!finitely}\index{MUB-complete}. By Corollary~\ref{cor:dual-spectral-domain} the dual bounded distributive lattice is $\Upfin(P)$ (see also Exercise~\ref{ex:dual-spectral-domain}).
\end{remark}
\begin{definition}\label{def:finite-mub-comp}
Let $P$ be a poset. We say $P$ is \emph{finitely \MUB-complete} if, for any finite subset $F \subseteq P$, $\MUB(F)$ is finite, and $\UB(F) = {\uparrow} \MUB(F)$.
\end{definition}
We note that a poset $P$ is finitely \MUB-complete if, and only if, the set $\Upfin(P)$ is a sublattice of $\cP(P)$.
To summarize the results in this section: a spectral domain $X$ is always the domain of ideals of a finitely \MUB-complete poset $P$, whose lattice of finitely generated up-sets is the dual distributive lattice of $X$.

\ourexercises

\begin{ourexercise}\label{ex:K=upminK}
Let $X$ be any $T_0$ space, and $K\subseteq X$ a compact-saturated subset of $X$. Then
$K={\uparrow}\min(K)$. \hint{Given $x\in K$, use Zorn's Lemma on the collection of closed sets
${\downarrow} y$ for $x\geq y\in K$ with reverse inclusion order to show that every element of $K$
is above a minimal one}.
\end{ourexercise}

\begin{ourexercise}\label{ex:teeth}
Let $X$ be a poset and $M,N\subseteq X$.
\begin{enumerate}
\item Show that ${\uparrow}M\cup{\uparrow}N={\uparrow}(M\cup N)$ and conclude that the finitely generated open up-sets of any topological space are closed under binary unions;
\item Show that ${\uparrow} M\cap {\uparrow} N=\bigcup\{\UB(\{m,n\})\mid m\in M\text{ and }n\in N\}$;
\item \label{itm:exa-not-coh} Give an example of a topological space, and finitely generated open up-sets $U$ and $V$ so that $U\cap V$ is not a finitely generated up-set;
\item Give a domain whose associated space gives an example as in (ref{itm:exa-not-coh});
\item Give examples showing that, even if $\UB(M)$ is non-empty, $\MUB(M)$ may be empty, and even if $\MUB(M)$ is non-empty, we may not have $\UB(M)={\uparrow}\MUB(M)$;
\item Show that if $X=\Idl(P)$ and $F\subseteq P$ is finite, then $\MUB(F)\subseteq P$.
\end{enumerate}
\end{ourexercise}

\begin{ourexercise}\label{ex:fin-MUB-comp}
  Prove that a poset $P$ is finitely \MUB-complete if, and only if, any finite intersection of finitely generated up-sets of $P$ is finitely generated. Conclude that this happens if, and only if, $\Upfin(P)$ is a bounded sublattice of $\cP(P)$.
\end{ourexercise}

\begin{ourexercise}\label{ex:dual-spectral-domain}
Let $X$ be a spectral domain. Show that the following three are isomorphic distributive lattices.
\begin{enumerate}
\item $\KO(X)$, the distributive lattice dual to $X$;
\item $\Upfin(\Kel(X))$, the lattice of finitely generated up-sets\index{finitely generated up-set} of the poset $\Kel(X)$, with the inclusion order;
\item The poset reflection\index{poset reflection} of the set $\cPfin(\Kel(X))$ of finite subsets of $\Kel(X)$, equipped with the preorder given by
\[
F\preceq G \ \iff \  \forall y\in G\ \exists x\in F\ (x\leq y).
\]
(The poset reflection was defined in Exercise~\ref{exe:reflection}.)
\end{enumerate}
\end{ourexercise}
\begin{ourexercise}\label{exe:join-primes-spectral-domain}
  Let $X$ be a spectral domain and let $L$ be its lattice of compact open subsets. Prove that $\cJ(L)^\op$ is isomorphic to $\Kel(X)$.  \hint{You can describe an isomorphism directly: send $p \in \Kel(X)$ to the join-prime element ${\uparrow}p$ of $L$.}
\end{ourexercise}

\section{Bifinite domains}\label{sec:bifinite}
Bifinite domains%
\footnote{\emph{SFP domains} were introduced in \cite{Plo76} as certain profinite posets, along with the powerdomain construction, under which they are closed. The name `SFP' is an acronym for `sequences of finite posets'. SFP domains are required to be $\omega$-algebraic, that is, the set of compact elements is countable, and pointed, that is, to have a least element. Bifinite domains may be viewed as the generalization of SFP domains where the $\omega$-algebraicity and pointedness are dropped. The minimal upper bounds point of view that we discuss here was developed by \cite{Smyth1983} and in the thesis \cite{Gunter1985}. Smyth showed that SFP is the largest Cartesian closed category of pointed $\omega$-algebraic domains.} 
are mathematical structures obtained as limits of finite `embedding projection pairs', which we will introduce below. Within domains, one can show that these are actually bi-limits, that is, simultaneously limits of the projections and colimits of the embeddings in the finite `embedding projection pairs'. 
Limits of finite posets clearly make profinite posets, that is, Priestley spaces (recall Example~\ref{exa:priestley-as-profinite posets}), pertinent to the subject, but it is only over time that the duality theoretic point of view came into focus, fully expressed in \cite{Abr}. Here we introduce bifinite domains as a subclass of spectral spaces, and we consider them from a duality theoretic point of view. In Corollary~\ref{cor:mub-bifinite} we give the characterization of bifinite domains via $\MUB$-completeness and in Theorem~\ref{thrm:dual-bifinite} we give a dual characterization in terms of a condition of `conjunctive closure' on distributive lattices, bringing us full circle, showing that these are precisely the structures studied in \cite{Abr}. 

A central notion in the definition of bifiniteness is that of a (finite) \emphind{embedding projection pair}, which we introduce now. This definition will be applied \emph{both} to distributive lattices \emph{and} to spectral domains. In the following definition, the category $\cC$ can be thought of as either $\cat{DL}$ or $\cat{Spec}$, or in fact any category equipped with a faithful functor to $\cat{Poset}$.
\begin{definition}\label{def:embedding-retraction-pair}
Let $\cC$ be a concrete category in which each object is equipped with a partial order and each morphism is order preserving. We say that a pair of morphisms \leftright{C}{D}{f}{g} of $\cC$ is an \emphind{embedding projection pair} (\epp) provided $(f,g)$ is an adjoint pair, $f$ is injective, and $g$ is surjective. Here, $f$ is called the \emph{embedding} and $g$ is called the \emph{projection} of the pair.
  Further, such an {\epp} is said to be \emph{finite} if $C$ is finite, and we call it an {\epp} `of' $D$.
  \end{definition}
Recall from Exercise~\ref{exe:adjunctioninjectivesurjective} that, for an adjoint pair between posets, the left adjoint is injective if, and only if, the right adjoint is surjective; so, for an adjunction to be an {\epp}, it suffices to check one of the two conditions (also see Exercise~\ref{ex:dual-adj}).


We start by characterizing the finite \epp's for distributive lattices and spectral spaces in terms of substructures, see Propositions~\ref{prop:fin-space-epp}~and~\ref{prop:lat-e-r-p} below.

\begin{notation}
  In what follows, we often consider pairs of functions between spectral spaces and distributive lattices. If $X$ and $Y$ are spectral spaces, then we use the notation \leftright{X}{Y}{e}{p} for a pair of spectral maps between these spaces. We then call the \emph{dual pair of homomorphisms} the pair of functions \leftright{K}{L}{i}{h}, where $K$ is the lattice dual to $X$, $L$ is the lattice dual to $Y$, $i$ is the homomorphism dual to $p$, and $h$ is the homomorphism dual to $e$.
\end{notation}

\begin{lemma}\label{lem:dual-adj-pairs}
Let \leftright{X}{Y}{e}{p} be a pair of spectral maps between spectral spaces $X$ and $Y$ and let \leftright{K}{L}{i}{h} be the dual pair of homomorphisms. %
Then $e$ is lower adjoint to $p$ if, and only if, $i$ is lower adjoint to $h$.
\end{lemma}

\begin{proof}
The function $e$ is lower adjoint to $p$ if, and only if, $e\circ p\leq \id_Y$ and $\id_X\leq p\circ e$ (see Exercise~\ref{exe:adjunctions}). 
The dual of the composite map $e\circ p$ is $i\circ h$ and the dual of $\id_Y$ is $\id_L$, so, using Remark~\ref{rem:order-yoga}, $e\circ p\leq \id_Y$ is equivalent to $i\circ h\leq \id_L$. Similarly $\id_X\leq p\circ e$ is equivalent to $\id_K\leq h\circ i$. Thus $e$ is lower adjoint to $p$ if, and only if, $i\circ h\leq\id_L$ and $\id_K\leq h\circ i$. This, in turn, is equivalent to saying that $i$ is lower adjoint to $h$. %
\end{proof}

We have the following corollary of Lemma~\ref{lem:dual-adj-pairs} and Exercise~\ref{ex:dual-adj}.

\begin{corollary}\label{cor:e-r-p-duality}
The Stone dual of an embedding projection pair on either side of the duality is an embedding projection pair on the other side, and the dual of a finite {\epp} is finite. 
\end{corollary}
Note that the dual of the embedding on either side is the projection on the other. The following propositions identify the nature of finite \epp's on either side of Stone-Priestley duality.

Recall that a map $f$ between topological spaces is said to be an \emphind{embedding} provided it is injective and the inverse function $f^{-1}\colon\im(f)\to X$ is also continuous. If $Y$ is a spectral space and $X \subseteq Y$ is a subspace of $Y$ which is itself a spectral space in the subspace topology, then we call the inclusion map $e \colon X \into Y$ a \emphind{spectral subspace embedding}. The following proposition allows us to see \epp's on a spectral space $Y$ as certain spectral subspaces $X$ of $Y$.
\begin{proposition}\label{prop:fin-space-epp}
Let $Y$ be a spectral space and $e\colon X\into Y$ a spectral subspace embedding. Then the following two conditions are equivalent:
\begin{enumerate}[label=(\roman*)]
\item the function $e$ has an upper adjoint $p$, and $p$ is spectral,
\item the subspace $X$ has the following two properties:
\begin{enumerate}[label=(\arabic*)]
\item for all $y\in Y$, the down-set ${\downarrow}y\cap X$ is  principal;
\item for all $U\subseteq X$ compact-open, the up-set ${\uparrow}U$ is open in $Y$. If, in addition, $X$ is finite, then it suffices that ${\uparrow}x$ is open in $Y$ for all $x\in X$.
\end{enumerate}
\end{enumerate}
\end{proposition}
\begin{proof}
  Note first that, essentially by definition, $e$ has an upper adjoint if, and only if, for every $y \in Y$, ${\downarrow} y \cap X$ has a maximum. Thus, (ii)(1) is equivalent to the existence of an upper adjoint. Now suppose the upper adjoint exists, and consider the map $p\colon Y\to X, y\mapsto\max({\downarrow}y\cap X)$. We show that {(ii)(2)} is equivalent to $p$ being a spectral map.

Notice that, for any up-set $U$ of $X$ we have $p^{-1}(U)={\uparrow}U$. Thus $p$ being spectral is the statement that ${\uparrow}U$ is compact-open whenever $U\subseteq X$ is compact-open. Also note that $U$ compact in X implies U compact in $Y$, which, in turn, implies that ${\uparrow}U$ is compact in $Y$. Thus we just need to know that ${\uparrow}U$ is open in $Y$. If, in addition, $X$ is finite, then every compact-open of $X$ is a finite union of principal up-sets ${\uparrow}_Xx$ and thus it suffices to consider these.
\end{proof}

\begin{proposition}\label{prop:lat-e-r-p}
Let $L$ be a bounded distributive lattice and $i\colon K\to L$ a bounded sublattice inclusion. Then the following conditions are equivalent:
\begin{enumerate}[label=(\roman*)]
\item the homomorphism $i$ has an upper adjoint that is also a homomorphism;
\item the sublattice $K$ has the following two properties:
\begin{enumerate}
\item[(1)] For all $b\in L$, ${\downarrow}b\cap K$ is a principal down-set;
\item[(2)] For any prime filter $F$ of $K$ the filter ${\uparrow}F$ in $L$ is again prime.
\end{enumerate}
\end{enumerate}
\end{proposition}
\begin{proof}
Again, the existence of an upper adjoint is equivalent to (ii)(1) since the upper adjoint of the embedding must be given by $h(b)=\max({\downarrow}b\cap K)$. We show that (ii)(2) is equivalent to this upper adjoint being a bounded lattice homomorphism.

Consider the map $h\colon L\to K, b\mapsto\max({\downarrow}b\cap K)$. Since $h$ is the upper adjoint of the inclusion, it preserves all existing meets. Also clearly $h(\bot)=\bot$. We show that (ii)(2) is equivalent to $h$ preserving binary joins. For this purpose, assume $h$ preserves binary joins and let $F$ be a prime filter of $K$. Then as $F$ is proper, so is ${\uparrow}F$. If $b_1\vee b_2\in{\uparrow}F$ in $L$, then there is $a\in F$ with $a\leq h(b_1\vee b_2)=h(b_1)\vee h(b_2)$ in $K$. Now since $F$ is  a prime filter, it follows that $h(b_1)\in F$ or $h(b_2)\in F$ and, since $h(b_i)\leq b_i$, it follows that $b_1\in{\uparrow}F$ or $b_2\in{\uparrow}F$ as required. Conversely, suppose that (ii)(2) holds. We want to show that $h$ preserves binary joins. Since $h$ is order preserving we have $h(b_1)\vee h(b_2)\leq h(b_1\vee b_2)$. Since any filter is the intersection of the prime filters containing it (see Exercise~\ref{exe:intersection-primefilters}), it suffices to show that any prime filter of $K$ containing $h(b_1\vee b_2)$ must contain $h(b_1)\vee h(b_2)$. To this end, let $F$ be a prime filter of $K$ with $h(b_1\vee b_2)\in F$.  Since $h(b_1\vee b_2)\leq b_1\vee b_2$, we have $b_1\vee b_2\in {\uparrow}F$ and by  (ii)(2), it follows that $b_i\in{\uparrow}F$ for $i=1$ or $2$. Now $b_i\in{\uparrow}F$ implies there is $a\in F$ with $a\leq b_i$, and as $a\in K$, it follows that $a\leq h(b_i)$ and thus $h(b_i)\in F$. That is,  $h(b_1)\vee h(b_2)\in F$ as required.
\end{proof}

Just as for spectral spaces, Proposition~\ref{prop:lat-e-r-p} allows us to see \epp's on a distributive lattice $L$ as certain sublattices $K$ of $L$, namely those satisfying the two conditions of Proposition~\ref{prop:lat-e-r-p}(ii). Interestingly, Corollary~\ref{cor:e-r-p-duality} then yields a duality between certain subspaces and certain sublattices, rather than the usual matching in duality of subs and quotients. This is of course because these subs are also quotients, but viewing this as a duality between sublattices and subspaces is interesting relative to a phenomenon we will meet in our second application in Chapter~\ref{ch:AutThry}, where certain residuated Boolean \emph{algebras} are dual to certain profinite \emph{algebras}, namely profinite monoids.

We will be particularly interested in finite \epp's. Since any inclusion of a finite (join semi-)lattice in a lattice has an upper adjoint, and since any prime filter in a finite lattice is principal, generated by a join-prime element, we obtain from Proposition~\ref{prop:lat-e-r-p} the following much simpler description of finite \epp's in a distributive lattice. 

\begin{corollary}\label{cor:fin-lat-e-r-p}
Let $L$ be a distributive lattice and $K\subseteq L$ a finite bounded sublattice, with $i \colon K \into L$ the inclusion map. Then the following conditions are equivalent:
\begin{enumerate}[label=(\roman*)]
\item there exists a bounded lattice homomorphism $h\colon L\to K$ such that $(i,h)$ is an {\epp} in the category $\DL$; 
\item $\cJ(K)\subseteq \cJ(L)$.
\end{enumerate}
\end{corollary}
\begin{definition}\label{def:finite-epp-sublattice}
Let $L$ be a distributive lattice. We say that a finite bounded sublattice $K$ of $L$ is a \emph{finite \epp-sublattice} if it satisfies the equivalent conditions of Corollary~\ref{cor:fin-lat-e-r-p}.
\end{definition}

Note that, if $K$ is a finite \epp-sublattice of $L$, then $K$ has `enough join-primes', in the sense that every element of $K$ is a finite join of join-prime elements (see Lemma~\ref{lem:finiteperfect}), and every join-prime of $K$ is also a join-prime of $L$. However, an infinite lattice $L$ need not have any join-prime elements at all (see Example~\ref{exa:nojoinirr}), and such a lattice does not have any finite \epp. Thus having `enough' finite \epp's is special. For $L$ to be bifinite, in addition, we require that the finite \epp-sublattices of $L$ form a directed diagram. Here, recall that ``directed'' means that for any two finite sublattices $K_1$ and $K_2$ of $L$ which give rise to \epp's, there is a finite sublattice $K_0$ of $L$ which gives rise to an {\epp} and contains both $K_1$ and $K_2$.

\begin{definition}\label{def:spectral-bif}
Let $X$ be a spectral space, $L$ its dual lattice. We say that $X$ and $L$ are \emphind{bifinite} provided either and then both of the following two equivalent conditions are satisfied:
\begin{enumerate}[label=(\roman*)]
\item $X$ is the projective limit in $\Spec$ of the projections of its finite \epp's;
\item $L$ is the directed colimit in $\DL$ of the embeddings of its finite \epp's.
\end{enumerate}
\end{definition}
In light of Definition~\ref{def:finite-epp-sublattice} and the fact that directed colimit of a collection of sublattices of a distributive lattice can be computed as a union (see Example~\ref{exa:DLasIndDLfin}), an equivalent definition for a distributive lattice $L$ to be bifinite is that the finite \epp-sublattices of $L$ form a directed diagram whose union is equal to $L$.
\begin{remark}
The definition we give here is not identical to any of the standard ones. The definition is usually given for algebraic domains, and not for spectral spaces, but as we will observe shortly, any bifinite spectral space is an algebraic domain. That is, the spaces we call bifinite are exactly the same objects as the bifinite domains of the domain theory literature.
\end{remark}

\begin{remark}\label{rem:bilimits}
Given a bifinite lattice or spectral space, we have a directed diagram of \epp's for that structure, and thus we can form two different diagrams: the diagram consisting of the embedding parts of the diagram of finite \epp's and the diagram consisting of the projection parts of the diagram of finite \epp's for the given structure. By definition, the colimit of the diagram of embeddings in a bifinite lattice yields the bifinite lattice itself. Dually the limit of the diagram of projections of the finite \epp's in a bifinite spectral space yields the bifinite space. On the other hand, if we take the limit in $\DL$ of the diagram of projections in a bifinite lattice, in general, we will get a bigger lattice, which has the colimit as a sublattice, and similarly, if we take the colimit in $\Spec$ of the diagram of embeddings in a bifinite spectral space we will get a bigger space, which has the limit as a quotient (see Exercise~\ref{exer:epp-bilimits}).

However, if we embed $\DL$ in $\Frame$ via the ideal completion, as in Figure~\ref{fig:dl-priest-spec-diagram}, and take the full subcategory of $\Frame$ it generates, then we get the category $\DL_{j-\rm approx}$ of distributive lattices with join-approximable relations (see Exercise~\ref{exer:join-approx}) and in this category we have a morphism from the limit in $\DL$ of the diagram of projections in a bifinite lattice to the bifinite lattice $L$ itself. This makes it possible to show that we have coincidence of the colimit of the injections and limit of the projections, thus showing that in this category $L$ is so-called \emphind{bifinite}, that is, simultaneously the colimit and the limit of its diagram of finite \epp's (for more details see Exercise~\ref{exer:epp-bilimits}). Finally, by duality, the same is true if we consider bifinite spectral spaces with arbitrary continuous functions rather than just spectral maps, see also \cite[Subsection~3.3.2]{AbJu94}.
\end{remark}

\begin{corollary}\label{cor:bif-implies-dom}
Let $X$ be a spectral space and $L$ its dual lattice. If $X$ and $L$ are bifinite then $X$ is a spectral domain, or, equivalently, $L$ has enough join-primes.
\end{corollary}
\begin{proof}
If a distributive lattice $L$ is bifinite, then every element $a\in L$ is contained in a finite \epp-sublattice $K$. In $K$, $a$ is a finite join of join-prime elements, and, by Corollary~\ref{cor:fin-lat-e-r-p}, these join-prime elements are also join-prime in $L$. Thus every element of $L$ is a finite join of join-prime elements. Now by Theorem~\ref{thrm:dual-spectral-domains} this is equivalent to $X$ being a spectral domain.
\end{proof}

We will now give an order theoretic characterization of bifinite spectral spaces. By Corollary~\ref{cor:bif-implies-dom}, if $X$ is bifinite with dual lattice $L$, then $X$ is a spectral domain. That is, by Remark~\ref{rem:MUB-comp}, the poset $P=\Kel(X)$, which is isomorphic to $\cJ(L)^{\op}$, is finitely \MUB-complete. We further have $X\cong\Idl(P)$ and $L\cong\Upfin(P)$ (see Exercise~\ref{ex:dual-spectral-domain}). In Proposition~\ref{prop:bifinite-MUB-complete}, we will characterize bifinite domains in terms of order theoretic properties of $P$.

\begin{definition}\label{def:MUB-closure}
Let $X$ be a poset and $F\subseteq X$. We say that $F$ is \emph{\MUB-closed}\index{MUB-closed} provided, for all $G\subseteq F$, we have $\MUB(G)\subseteq F$. The \emph{\MUB-closure}\index{MUB-closure} of $F\subseteq X$ is the least $F'\subseteq X$ such that $F\subseteq F'$ and $F'$ is \MUB-closed.
\end{definition}

\begin{proposition}\label{prop:bifinite-MUB-complete}
Let $P$ be a finitely \MUB-complete poset, and let $X:=\Idl(P)$ be the corresponding spectral domain. Then  a finite subspace $F\subseteq X$ gives rise to an {\epp} in the category $\Spec$ if, and only if, $F\subseteq P$ and $F$ is \MUB-closed. Furthermore, the spectral domain $X$ is bifinite if, and only if, the \MUB-closure of any finite subset of $P$ is finite.
\end{proposition}

\begin{proof}
A finite subspace $F\subseteq X$ gives rise to an {\epp} if, and only if,  it satisfies the equivalent conditions of Proposition~\ref{prop:fin-space-epp}. The finite case in condition (ii)(2) of that proposition is equivalent to $F\subseteq P$. We show that (ii)(1) is equivalent to $F$ being \MUB-closed. %

 Let $F\subseteq P$ be finite and suppose that $\max({\downarrow}x\cap F)$ exists for each $x\in X$. Let $G\subseteq F$ and let $p\in\MUB(G)$ and consider the set ${\downarrow}p\cap F$. By our assumption, there is $p'=\max({\downarrow}p\cap F)$. Also, clearly $G\subseteq{\downarrow}p\cap F$, so $p'\in\UB(G)$. Since $P$ is \MUB-complete, there is $p''\in\MUB(G)$ with $p''\leq p'$. Thus we have $p''\leq p'\leq p$ with both $p'',p\in\MUB(G)$. It follows that $p''=p'=p$ and thus $p\in F$ since $p'\in F$. That is, $F$ is \MUB-closed.

For the converse, again let $F\subseteq P$ be finite, suppose $F$ is \MUB-closed, and let $x\in X$. Further set $G={\downarrow}x\cap F$. Since $x\in\UB(G)$, there is a $p\in\MUB(G)\subseteq F$ with $p\leq x$. But then $p\in G$ and thus $p=\max(G)=\max({\downarrow}x\cap F)$ as required.

Finally, $X$ is bifinite if, and only if, every finite subset $F$ of $P$ is
contained in a finite subset $F'$ which gives rise to an \epp, which we showed
happens if, and only if, $F' \subseteq P$ and $F'$ is \MUB-closed. Thus, $X$ is bifinite if, and only if, the \MUB-closure of any finite subset of $P$ is finite.
\end{proof}

\begin{corollary}\label{cor:mub-bifinite}
Let $X$ be an algebraic domain. Write $P:=\Kel(X)$, so that $X\cong\Idl(P)$. Then $X$ is {bifinite}\index{bifinite} if, and only if, for all finite subsets $F$ of $P$, the following three properties hold:
\begin{enumerate}
\item $\MUB(F)$ is finite,
\item $\UB(F)={\uparrow}\MUB(F)$, and 
\item the \MUB-closure of $F$ is finite.
\end{enumerate}
\end{corollary}

A poset $P$ is called a \emphind{Plotkin order} if it satisfies the three
properties in Corollary~\ref{cor:mub-bifinite} for all finite subsets $F$ of
$P$. Thus Corollary~\ref{cor:mub-bifinite} shows that an algebraic domain $X$ is
bifinite if, and only if, the poset $P = \Kel(X)$ is a Plotkin order.
Note that the properties {(a)} and {(b)} in Corollary~\ref{cor:mub-bifinite} are equivalent to saying that the corresponding domain $X=\Idl(P)$ is spectral, see Corollary~\ref{cor:2/3sfp}.

We finish the subsection with a characterization of the lattices dual to bifinite domains. For this purpose we define the notion of conjunctive closure.

\begin{definition}\label{def:conj-clos}
Let $L$ be a lattice and $S\subseteq L$. We say that $S$ is \emphind{conjunctively closed} provided, for each finite subset $F\subseteq S$, there is a finite $G\subseteq S$ such that
\[
\bigwedge F=\bigvee G.
\]
\end{definition}

\begin{theorem}\label{thrm:dual-bifinite}
Stone's duality between $\Spec$ and $\DL$ restricts to a duality between bifinite domains and bounded distributive lattices $L$ with enough join-primes satisfying the following property: for every finite $S\subseteq \cJ(L)$, there exists $S'\subseteq \cJ(L)$ finite, such that $S' \supseteq S$ and $S'$ is conjunctively closed.
\end{theorem}

\begin{proof}
By Theorem~\ref{thrm:dual-spectral-domains} spectral domains are dual to bounded distributive lattices $L$ with enough join-primes, so we just need to verify that the additional property in the theorem corresponds dually to condition { (iii)} in the definition of bifinite domains (see Exercise~\ref{exer:dual-bif-dom}).
\end{proof}

\begin{example}\label{exa:prime-bifinite}
  We show that $(\bN, |)$, the divisibility lattice of $\bN$, is a bifinite lattice.   Recall that we computed the Priestley dual space of $(\bN, |)$ in Example~\ref{exa:primenumbers-topology}. To prove that the lattice $L := (\bN, |)$ is bifinite, let $n \in \bN$ be arbitrary. Denote by $D_n$ the subset of natural numbers that are $0$ or divisors of $n$, that is, the down-set of $n$ in $L$, with the top element $0$ added. This is a finite bounded sublattice of $L$; we apply Corollary~\ref{cor:fin-lat-e-r-p} to show that it is an \epp-sublattice. Let $j$ be join prime in $D_n$. Then $j | n$, and we show that $j$ must in fact be a positive power of a prime number: since $1$ is the bottom element of $D_n$, $j \neq 1$, and we may write $j$ as the join of $p_1^{k_1}, \dots, p_r^{k_r}$, for some $r \geq 1$, prime numbers $p_1, \dots, p_r$ and $k_1, \dots, k_r \geq 1$. Then each $p_i^{k_i}$ divides $j$, and hence, since $j \in D_n$, $p_i^{k_i}$ is also an element of $D_n$. Moreover, by the assumption that $j$ is join prime and the fact that each $p_i^{k_i}$ is incomparable to $p_j^{k_j}$, we must have $r = 1$. Thus, $j$ is also join prime in $L$. 
  Recall that we gave an example of a quotient space of the dual of $(\bN, |)$ in Example~\ref{exa:primenumbers-maps}.  For any $n \in \bN$, the finite \epp-sublattice $D_n$ defined here dually yields a different, finite quotient of the Priestley dual space of $(\bN, |)$.
  \end{example}
  
\ourexercises

\begin{ourexercise}\label{ex:dual-adj}
This exercise concerns Lemma~\ref{lem:dual-adj-pairs}.
\begin{enumerate}
\item Prove that (ii) implies (i);
\item Show that all of the following are equivalent:
\begin{enumerate}
\item[(i)] $e$ is injective;
\item[(ii)] $p$ is surjective;
\item[(iii)] $i$  is injective;
\item[(iv)] $h$ is surjective,
\end{enumerate}
and that if these hold, then $e$ and $i$ are embeddings.
\end{enumerate}
\end{ourexercise}

\begin{ourexercise}\label{exer:dual-bif-dom}
  Let $X$ be a spectral domain, $P := \Kel(X)$ and $L$ the lattice dual to $X$. Show that the following are equivalent:
  \begin{enumerate}[label=(\roman*)]
    \item for every finite $F \subseteq P$, the \MUB-closure of $F$ is finite;
    \item for every finite $S \subseteq \cJ(L)$, there exists finite $S' \subseteq \cJ(L)$ such that $S' \supseteq S$, and $S'$ is conjunctively closed.
  \end{enumerate}
\hint{Use the isomorphism between $\cJ(L)^\op$ and $P$, and the fact that, for a finite collection $F \subseteq P$, if $G = \MUB(F)$, then $$\UB(F) = \bigcap_{p \in F} {\uparrow} p = \bigcup_{q \in G} {\uparrow} q.$$}
\end{ourexercise}

\begin{ourexercise}\label{exer:epp-bilimits}
Let $\cC$ be a concrete category in which each object is equipped with a partial order and each morphism is order preserving. Let 
\[ \{\leftright{C}{D}{f}{g}\} \]  
be a directed diagram of \epp's in $\cC$. 
Further let $C_f$ be the colimit in $\cC$ of the embeddings of the diagram and $C_g$ the limit in $\cC$ of the projections of the diagram.
\begin{enumerate}
    \item Show that there is a unique morphism $C_f\to C_g$ in $\cC$.
    \item Consider now the particular case $\cC = \DL$. Let $L$ be a bifinite lattice and let 
    \[ \{\leftright{L}{K}{i}{h}\}  \] 
    be the associated diagram of finite \epp's. Prove that the colimit $L_i$ of the embedding parts of the diagram is 
    isomorphic to $L$, and that the homomorphism $L \to L_h$, obtained from the previous item by composing with this isomorphism, is an embedding. 
    
    Conclude also that, by Stone duality, if $X$ is a bifinite spectral space and \[ \{\leftright{X}{Y}{e}{p}\} \] 
    is the associated diagram of finite \epp's, then the limit $X_p$ of the projection parts is homeomorphic to $X$, and the spectral map $X_e\to X$, obtained from the previous item by composing with this homeomorphism, is a quotient map. 
    \item Prove that, if $L$ is a bifinite lattice and $X$ is the dual spectral space, then with the notations of the previous item, $L_h$ is dual to $X_e$ and the morphism  $L\to L_h$ is dual to the spectral map $X_e\to X$.
    \item Now consider the situation where $L=\Upfin(P)$, where $P$ is a Plotkin order. Show that $L_h\cong\Up(P)$ and conclude that $L$ and $L_h$ are not necessarily isomorphic.
    \item Still under the assumption that $L=\Upfin(P)$, where $P$ is a Plotkin order, show that the relation 
    \[
    R=\{(U,V)\in\Up(P)\times\Upfin(P)\mid U\supseteq V\}
    \]
    is a join-approximable relation from the lattice $L_h \cong \Up(P)$ to the lattice $L = \Upfin(P)$, and conclude that $L=\Upfin(P)$ is the bilimit of the diagram of \epp's for $L$ in the category of distributive lattices with join-approximable relations. 
  \end{enumerate}
\end{ourexercise}

\begin{ourexercise}\label{exe:prime-bifinite}
This exercise is about the lattice $\mathbb{N}$ with the divisibility order. We already studied this lattice as a running example for Priestley duality in Chapter~\ref{ch:priestley}, see Examples~\ref{exa:prime-divisibility}, \ref{exa:primenumbers-topology}, and \ref{exa:primenumbers-maps}. Here 
we consider the spectral space $X$ dual to $(\mathbb{N}, \mid)$, and its specialization order. Recall from Remark~\ref{rem:order-yoga} that the specialization order on $X$ is the \emph{opposite} of the order on the Priestley space dual to $(\mathbb{N}, \mid)$, depicted in Figure~\ref{fig:prime-naturalnumbers}.

\begin{enumerate}
\item Show that for any $x \in X$, the principal up-set ${\uparrow} x$ is either finite or equal to $X$.
\item Show that the Scott topology on the spectral space $X$ coincides with the Alexandroff topology.
\item Deduce that $X$ is an algebraic domain in which all elements are compact, and that the poset $\Kel(X) = X$ is finitely MUB complete.
\item Identify those finite subsets of $X$ which correspond to finite $\epp$'s, and show that $X$ is bifinite.
\end{enumerate}
\end{ourexercise} 

\section{Domain theory in logical form}\label{sec:DTLF}
In this section, we analyze the Domain Theory in Logical Form (DTLF) program in the setting of bifinite domains as in \cite{Abr}, emphasizing the duality theoretic aspects. We will illustrate DTLF by considering duality for the function space construction and a domain equation based on this constructor. 

\subsection*{Bifinite domains are closed under function space}
In this subsection we use the duality developed for function spaces in
Section~\ref{sec:funcspace} to show that the category $\bf Bif$ of bifinite
domains with Scott continuous functions is closed under the function space
construction. Recall that, for topological spaces $X$, $Y$, the set $[X,Y]$ of
continuous functions from $X$ to $Y$ is equipped with the \emph{compact-open
topology}: a subbase is given by the sets of functions $K \Rightarrow U := \{f \in [X,Y] \
\mid \ f[K] \subseteq U\}$, for $K \subseteq X$ compact and $U \subseteq Y$
open, see Definition~\ref{def:funcspace}. We aim to prove here
(Theorem~\ref{thm:Bif-CCC} below) that this space $[X,Y]$ is bifinite if both
$X$ and $Y$ are. For this purpose, we will make crucial use of the fact that
this space is  given by the smallest congruence that is \emphind{join preserving
at primes}, the property that we identified in Definition~\ref{def:jpp} above. 
Recall that, if $X$ and $Y$ are bifinite, then in particular they are spectral domains by Corollary~\ref{cor:bif-implies-dom} and thus, by Corollary~\ref{cor:spectral-funcspace}, $[X,Y]$ is a spectral space. We will need to prove that $[X,Y]$ can still be approximated by finite \epp's, if both $X$ and $Y$ have this property. To this end, we first prove a lemma about lifting finite \epp's to function spaces. 
\begin{lemma}\label{lem:ep-pairs}
Let $X, Y$ be two spectral spaces such that $[X,Y]$ is a spectral space, and let $X', Y'$ be finite posets, equipped with the Alexandrov topology. If $\leftright{X'}{X}{e_X}{p_X}$ and $\leftright{Y'}{Y}{e_Y}{p_Y}$ are  \epp's, then the functions
\begin{align*}
  e\colon [X',Y']&\to [X,Y]\ \colon p\\
             f'\quad&\mapsto\ e_Yf'p_X\\
             p_Yfe_X&\mapsfrom \quad f.
  \end{align*}
are spectral, and form an {\epp} between $[X',Y']$ and $[X,Y]$. Moreover, for any up-sets $K' \subseteq X'$, $U' \subseteq Y'$, $p^{-1}(K' \Rightarrow U') = p_X^{-1}(K') \Rightarrow p_Y^{-1}(U')$.
\end{lemma}

\begin{proof}
  Let $\leftright{X'}{X}{e_X}{p_X}$ and $\leftright{Y'}{Y}{e_Y}{p_Y}$ be finite \epp's. We define 
\begin{align*}
e\colon [X',Y']&\to [X,Y]\ \colon p\\
           f'\quad&\mapsto\ e_Yf'p_X\\
           p_Yfe_X&\mapsfrom \quad f.
\end{align*}
Note that, for any $f' \in [X',Y']$, we have 
\[pe(f')=p_Ye_Yf'p_Xe_X=\id_{Y'}f'\id_{X'}=f',\] showing that 
$p \circ e=\id_{[X',Y']}$. 
Also, for any $f \in [X,Y]$, we have 
\[ep(f)=e_Yp_Yfe_Xp_X\leq\id_{Y}f\id_{X}=f,\] showing that 
$e\circ p\leq\id_{[X,Y]}$. Thus, $(e,p)$ is an adjoint pair of monotone functions between the underlying posets, $e$ is injective, and $p$ is surjective. It remains to show that $e$ and $p$ are spectral maps.

Thus, in order to prove that $e$ is continuous, and hence spectral because $[X',Y']$ is finite, it suffices to show that $e^{-1}(K \Rightarrow U)$ is an up-set for every $K \subseteq X$ compact and $U \subseteq Y$ open. Recall that the specialization order on a function space coincides with the pointwise order when the domain is locally compact (see Exercise~\ref{ex:compact-open topology}.\ref{spec-order-pointwise}), and so certainly in this case. Suppose that $e(f') \in K \Rightarrow U$ for some $f' \in [X',Y']$, and that $f' \leq g$ for some $g \in [X',Y']$. We show that $e(g) = e_Y g p_X \in K \Rightarrow U$. Let $x \in K$ be arbitrary. Note that $f' p_X[K] \subseteq e_Y^{-1}(U)$ by assumption, and that $f'p_X(x) \leq gp_X(x)$. Thus, since $e_Y^{-1}(U)$ is an up-set by continuity of $e_Y$, it follows that $gp_X(x) \in e_Y^{-1}(U)$, that is, $e_Y g p_X(x) \in U$, as required.

Note that the fact that $p$ is spectral will follow from the ``moreover''
statement, because both $p_X$ and $p_Y$ are spectral, and the sets of the form
$K' \Rightarrow U'$, with $K'$ and $U'$ up-sets of $X'$ and $Y'$, respectively,
form a base for the topology on $[X',Y']$. Thus, to prove the moreover
statement, let $K' \subseteq X'$ and $U' \subseteq Y'$ be up-sets and let $f \in
[X,Y]$ be arbitrary. We will show that $p(f) \in K' \Rightarrow U'$ if and only
if $f \in p_X^{-1}(K') \Rightarrow p_Y^{-1}(U')$. For the left-to-right direction, suppose that $p(f) \in K' \Rightarrow U'$ and let $x \in X$ be such that $p_X(x) \in K'$. Then, using the adjunction between $e_X$ and $p_X$,
\[ p_Yf(x) \geq p_Y fe_Xp_X(x) = p(f)(p_X(x)),\] 
and the latter is an element of $U'$ by assumption, so that $p_Y f(x) \in U'$ since $U'$ is an up-set. 
Conversely, suppose that $f \in p_X^{-1}(K') \Rightarrow p_Y^{-1}(U')$ and let $x' \in K'$ be arbitrary. Using that $p_Xe_X(x') = x'$, we get $f e_X(x') \in p_Y^{-1}(U')$, which shows that $p(f)(x') \in U'$, as required.
\end{proof}

\begin{theorem}\label{thm:Bif-CCC}
Let $X$ and $Y$ be bifinite domains. Then $[X,Y]$ is a bifinite domain.
\end{theorem}
\begin{proof}
Denote by $L$ and $M$ the bifinite distributive lattices of compact-opens of $X$ and $Y$, respectively. Note that $L$ in particular has enough join-primes, by Corollary~\ref{cor:bif-implies-dom}. Thus, Corollary~\ref{cor:spectral-funcspace} applies, and we obtain that the space $[X,Y]$ is spectral, and its distributive lattice of compact-opens is (up to isomorphism) the lattice $K := F_{\to}(L,M)/{\theta_{\jpp}}$. By duality, it is thus equivalent to show that this lattice $K$ is bifinite, that is, that every element of $K$ lies in some finite \epp-sublattice. Since the elements of the form ${a} \Rightarrow {b}$, with $a \in L$ and $b \in M$, generate the lattice $K$, it suffices to show the result for these elements. Let $a \in L$ and $b \in M$. Since $L$ and $M$ are bifinite, pick  finite \epp-sublattices $L'$ and $M'$ of $L$ and $M$, respectively, with dual \epp's $\leftright{X'}{X}{e_X}{p_X}$ and $\leftright{Y'}{Y}{e_Y}{p_Y}$. Combining Lemma~\ref{lem:ep-pairs} with Lemma~\ref{lem:dual-adj-pairs}, the lattice of compact-opens of $[X',Y']$ is isomorphic to a finite \epp-sublattice $K'$ of $K$, with embedding given by $p^{-1}$ and projection given by $e^{-1}$. This sublattice $K'$ in particular contains the element ${a} \Rightarrow b$ because $p^{-1}(a \Rightarrow b)$ is equal to $p_X^{-1}(a) \Rightarrow p_Y^{-1}(b)$. Thus, $K$ is a bifinite lattice. Finally, since $X$ and $Y$ are bifinite, they are spectral domains and thus carry their Scott topologies. It follows in particular that $[X,Y]$ is the set of Scott continuous functions from $X$ to $Y$. Finally, since $[X,Y]$ is bifinite, it is also a spectral domain by Corollary~\ref{cor:bif-implies-dom}, so that the spectral topology on $[X,Y]$ coincides with the Scott topology of its specialization order.
\end{proof}

\subsection*{Variations on a domain equation}
 
In the above proof of Theorem~\ref{thm:Bif-CCC}, we proved a fact about the type
constructor of function spaces of domains by transferring it to a dual
construction on the corresponding distributive lattices. This is the fundamental
idea of Domain Theory in Logical Form \parencite{Abramsky87, Abr}: all of the domain or
type constructors of interest (sums, products, various forms of power domains,
functions spaces, etc.) are dual to endofunctors on the category of bifinite
lattices that automatically preserve directed colimits of sequences since they
are given `algebraically', that is, freely by generators and relations. It
follows that any composition $F$ of these endofunctors also preserves directed
colimits. This allows one in particular to build fixed points, that is, solutions
of equations of the form $X\cong F(X)$, by iterating $F$, starting from a
morphism $A\to F(A)$ (that is, a co-algebra for the functor $F$). A canonical such fixed point is the \emphind{least fixed point} of $F$, which can be constructed as the colimit of the following sequence, starting from the free distributive lattice on the empty set, $\mathbf{2}$:
 \begin{equation} \label{eq:iterative-sequence}
 \mathbf{2}\stackrel{e_0}{\longrightarrow}F(\mathbf{2})\stackrel{e_1}{\longrightarrow}F^2(\mathbf{2})\stackrel{e_2}{\longrightarrow}\ \dotsm.
\end{equation}
Here, $F$ is any directed-colimit-preserving functor on $\DL$, $e_0$ is the unique morphism from $\mathbf{2}$ to $F(\mathbf{2})$, given by the fact that $\mathbf{2}$ is initial in $\DL$, and $e_{n+1}:=F(e_n)$. If in addition the maps of the sequence consist of \epp's then we can conclude that the fixed point obtained as the colimit of (\ref{eq:iterative-sequence}) is again bifinite.
 
 We now illustrate these ideas by considering the most classical example of a domain equation, namely $X\cong[X,X]$. One may first look for the least solution of this equation, but, as classical as this case is, it is also an anomaly, in that the least solution is not so interesting, since it is the one-element space. On the dual side, this is reflected by the fact that the iterative sequence (\ref{eq:iterative-sequence}), applied in the case where $F$ is the functor sending a lattice $L$ to $F_{\to}(L,L)/{\theta_\jpp}$, never `gets off the ground', that is, the very first map $e_0 \colon \mathbf{2} \to F_{\to}(\mathbf{2},\mathbf{2})/{\theta_{\jpp}}$ is already an isomorphism. This follows from the observation that $[1,1] \cong 1$ and Corollary~\ref{cor:spectral-funcspace}, but it is also not hard to give an elementary algebraic proof directly from the definition of the lattice $F_{\to}(\mathbf{2},\mathbf{2})/{\theta_\jpp}$: notice first that already in $F_{\to}(\mathbf{2},\mathbf{2})$, we have the equalities 
 \[
 0\to 0=1\ ,\quad 0\to 1=1\ ,\quad 1\to 1=1
 \] 
 because an implication-type operator by definition sends any pair with $\bot$ in the first coordinate or $\top$ in the second coordinate to $\top$; indeed, these are the `empty set' instances of the schemes (\ref{eq:meet-snd-co}) and (\ref{eq:meet-fst-co}) that define $F_{\to}(L,M)$. Moreover, since $\theta_\jpp$ makes $\to$ preserve joins at primes, and $1$ is a join prime element in $\mathbf{2}$, we must have $1\to 0=0$ in $F_{\to}(\mathbf{2},\mathbf{2})/{\theta_{\jpp}}$. Thus, all four generators of $F_{\to}(\mathbf{2},\mathbf{2})/{\theta_{\jpp}}$ are equal to $0$ or $1$.

For the rest of this section, let $F_{\End}$ denote the functor on bifinite
lattices that sends $L$ to $F_{\to}(L,L)/{\theta_\jpp}$, which is the dual
lattice of the endomorphism type constructor $X \mapsto [X,X]$ on bifinite
domains, by Corollary~\ref{cor:spectral-funcspace}. The previous paragraph shows
that the least fixpoint of $F_{\End}$ is $\mathbf{2}$, so we can not directly
use the construction in (\ref{eq:iterative-sequence}). There are, essentially,
two ways around this problem. One is to change the functor with which we induct;
the other is to start the induction `higher up'. Among the amended functors, a
minimal and very natural choice is to consider $T(X)=\bot\oplus[X,X]$, usually
denoted $[X,X]_\bot$, which takes the function space and adds a new bottom.
Clearly the dual is $L \mapsto F_{\End}(L) \oplus \bf 1$, which adds a new top.
This functor is natural from the programming languages point of view and has
been studied extensively, because $X\cong[X,X]_\bot$ is the domain equation
corresponding to what is known as the ``lazy lambda calculus'', in which lambda
terms are identified with their so-called weak head normal form
\parencite{Abr90,AbrOng93}. In this case the canonical solution, that is, the least fixed point of the sequence (\ref{eq:iterative-sequence}), is non-trivial but its function space, is a retract of $X$, and not isomorphic to $X$. This construction is the subject of Example~\ref{ex:lazylambda} below.
 
 From the theoretical point of view, it is of course also important to show that there \emph{exist} non-trivial, on-the-nose, solutions to the equation  $X \cong [X,X]$. The way to obtain this is to stick with the original functor $F_\End$ as the functor we iterate, but to start the iterative process in (\ref{eq:iterative-sequence}) with a bigger lattice than $\mathbf{2}$. This is the approach taken by Scott at the very start of domain theory, see for example \cite{Scott72, Scott80}. This however inherently involves an ad hoc choice of a starting point for the iterative sequence. In \cite{SmyPlo82}, this approach is cast in a general scheme based on solving the domain equation in a comma category. The solution of $X \cong [X,X]$ starting from Sierpinski space, or, on the lattice side, from the $3$-element chain, is the subject of Example~\ref{ex:Scott}. 
 
 It is also worth mentioning the treatment in the book  \cite{LamSco86} on categorical logic, which solves a three-way equation $X \cong [X,X] \cong X \times X$. More generally, modelling the lambda calculus and associated programming languages is a vast field in theoretical computer science which goes far beyond the confines of duality theory, see for example \cite{Ba84-2, BarMan2022}, and also see \cite{MaSa08} for a more recent survey on universal algebra in lambda theories.
 
 \begin{example}\label{ex:Scott}
We want a solution to the domain equation $X\cong[X,X]$, starting the iterative process from the Sierpinski space $\mathbb S$. In dual form, this means we want a solution to the equation $L\cong F_{\End}(L)$ within bifinite lattices, starting from the three element chain, which we will denote by $\three$. As we have explained above, in order to get started, we need an embedding $\three\to F_{\End}(\three)$ or, dually speaking, a projection $[\mathbb S,\mathbb S]\to \mathbb S$. Also, as we shall see, it will be important that these morphisms are part of embedding-projection-pairs. So let's first see whether this is possible. As the calculation is easier on the space side, we look at it there. First, $[\mathbb S,\mathbb S]$ is the poset of all order-preserving endofunctions on $\mathbb S$. There are $2^2=4$ functions from a two-element set to itself. In the case of a two element chain, just one is not order preserving. That is, $[\mathbb S,\mathbb S]=\{\underline{0},\id_{\mathbb S},\underline{1}\}$, where  $\underline{0}$ is the constantly $0$ function, $\id_{\mathbb S}$ is the identity function, and $\underline{1}$ is the constantly $1$ function. Since the projection map $[\mathbb S, \mathbb S]\to\mathbb S$ must be surjective, $\underline{0}$ and $\underline{1}$ must be sent to $0$ and $1$, respectively. There are now two choices for where to send $\id_{\mathbb S}$: for $k \in \{0,1\}$, let $p^{(k)}$ denote the function $[\mathbb S, \mathbb S] \onto \mathbb S$ that sends $\underline{0}$ to $0$, $\underline{1}$ to $1$, and $\id_{\mathbb{S}}$ to $k$, see Figure~\ref{fig:Scotts-epps}. 

\begin{figure}[htp]
  \begin{center}
    \begin{tikzpicture}
      \polab{(0,0)}{$\underline{0}$}{left}
      \polab{(0,1)}{$\id_{\mathbb S}$}{left}
      \polab{(0,2)}{$\underline{1}$}{left}
      \draw (0,0) -- (0,2);

      \polab{(2,0)}{$0$}{right}
      \polab{(2,2)}{$1$}{right}
      \draw (2,0) -- (2,2);
      
      \draw[->] (0.2,0.1) to[bend left=15] (1.8,0.1);
      \draw[->] (0.2,1) to[bend left=15] (1.8,0.2);
      \draw[->] (0.2,2.1) to[bend left=15] (1.8,2.1);

      \draw[dotted,->] (1.8,-0.1) to[bend left=15] (0.2,-0.1);
      \draw[dotted, ->] (1.8,1.9) to[bend left=15] (0.2,1.9);

      \node at (0, -1.5) {$[\mathbb S, \mathbb S]$};
      \node at (2, -1.5) {$\mathbb S$};

      \draw[->] (0.5, -1.3) to[bend left=15] node[above] {$p^{(0)}$} (1.7, -1.3);
      \draw[->,dotted] (1.7, -1.7) to[bend left=15] node[below] {$e^{(0)}$} (0.5, -1.7);

      \polab{(4,0)}{$\underline{0}$}{left}
      \polab{(4,1)}{$\id_{\mathbb S}$}{left}
      \polab{(4,2)}{$\underline{1}$}{left}
      \draw (4,0) -- (4,2);

      \polab{(6,0)}{$0$}{right}
      \polab{(6,2)}{$1$}{right}
      \draw (6,0) -- (6,2);
      
      \draw[->] (4.2,0.1) to[bend left=15] (5.8,0.1);
      \draw[->] (4.2,1) to[bend right=15] (5.8,1.8);
      \draw[->] (4.2,2.1) to[bend left=15] (5.8,2.1);

      \draw[dotted,->] (5.8,-0.1) to[bend left=15] (4.2,-0.1);
      \draw[dotted, ->] (5.8,1.9) to[bend right=15] (4.2,1.1);
      \node at (4, -1.5) {$[\mathbb S, \mathbb S]$};
      \node at (6, -1.5) {$\mathbb S$};

      \draw[->] (4.5, -1.3) to[bend left=15] node[above] {$p^{(1)}$} (5.7, -1.3);
      \draw[->,dotted] (5.7, -1.7) to[bend left=15] node[below] {$e^{(1)}$} (4.5, -1.7);

    \end{tikzpicture}
    \end{center}
\caption{The two embedding-projection pairs in Example~\ref{ex:Scott}. The solid lines represent the projections $p^{(0)}$ and $p^{(1)}$, the dotted lines represent the embeddings $e^{(0)}$ and $e^{(1)}$.}
\label{fig:Scotts-epps}
\end{figure}
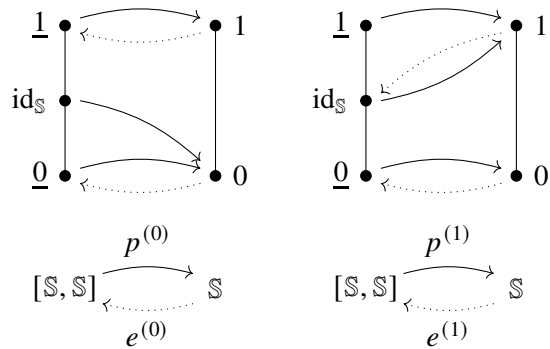

Both $p^{(0)}$ and $p^{(1)}$ have lower adjoints, $e^{(0)}$ and $e^{(1)}$, respectively. Note that $e^{(1)}$ sends $0$ to $\underline{0}$ and $1$ to $\id_{\mathbb S}$, while $e^{(0)}$ sends each element of $\mathbb S$ to the corresponding constant function. Note also that for each $k \in \{0,1\}$, the pair $(e^{(k)}, p^{(k)})$ is an \epp, since all order-preserving maps between finite sets are spectral maps. We now choose to work with $(e^{(0)}, p^{(0)})$, because the function $e^{(0)}$ is the more natural embedding of $\mathbb S$ into $[\mathbb S, \mathbb S]$. Using Lemma~\ref{lem:ep-pairs}, we obtain  by induction a sequence of spaces
\[
Y_0:=\mathbb S\qquad\text{ and } \qquad Y_{n+1}=[Y_n,Y_n] \text{ for all } n \geq 0,
\]
and a sequence of \epp's 
\[
e_0:=e^{(0)}, p_0:=p^{(0)} \qquad\text{ and } \qquad      e_{n+1}:=e_n\circ\underline{\ \ }\circ p_n, \ p_{n+1}:=p_n\circ\underline{\ \ }\circ e_n
\]
and thus we have a sequence of finite \epp's between domains
\[
\mathbb S\ \leftrightarrows \ [\mathbb S,\mathbb S]\ \leftrightarrows \ [[\mathbb S,\mathbb S],[\mathbb S,\mathbb S]]\ \leftrightarrows \ \dotsm
\]
and dually a sequence of finite \epp's between finite distributive lattices
\[
\three\ \leftrightarrows \ F_{\End}(\three) \ \leftrightarrows \ F_{\End}(F_{\End}(\three)) \ \leftrightarrows \dotsm
\]
It remains to show that the inverse limit of the spaces, or equivalently, the colimit of the lattices is bifinite. This is easier to do on the lattice side, as we show now.

Let $L$ be the colimit of a sequence or \epp's $e_n\colon L_n\leftrightarrows L_{n+1}\colon p_n$. As in the proof of Theorem~\ref{thm:Bif-CCC}, we use the fact that a directed colimit of inclusions is, up to isomorphism, just given by union, so that each $L_n$ is an \epp-sublattice of $L$. The embedding $e^n\colon L_n\to L$ is the one given by the colimit, and the projection $p^n\colon L\to L_n$ is given for $a\in L_m$ by $p_n\circ\ \dotsm\ \circ p_{m-1}(a)$ if $m>n$, by the identity if $m=n$, and by $e_{n-1}\circ\ \dotsm\ \circ e_m$ if $m< n$. We leave it as an exercise to check that this is a well-defined embedding projection pair (see Exercise~\ref{exe:colimit-e-r-ps}).
\end{example}

 \begin{example}\label{ex:lazylambda}
 We want to construct the least solution to the domain equation $X\cong[X,X]_\bot$, or dually, to the equation $L\cong F_{\to}(L,L)/{\theta_\jpp}\oplus \bf 1$ within bifinite distributive lattices. The development is precisely the same as in the above example once we have an {\epp} to initialize the process. Again, it is simplest to argue on the space side. The terminal bifinite space is the one-element space, which we will denote by $1$, and $[1,1]_\bot$ is the two-element ordered space $\bot<1$. The difference here, which makes the process canonical, is that we have only one choice for $p\colon [1,1]_\bot \to 1$. Its lower adjoint is $e\colon 1\to [1,1]_\bot, 1\mapsto \bot$. Once we have this initial {\epp}, the arguments are identical to the ones in Example~\ref{ex:Scott}. We leave the details as an exercise.
 See \cite[Section~4]{AbrOng93} for further details and the computational intuition behind the equation $X\cong[X,X]_\bot$.
\end{example}

\ourexercises

\begin{ourexercise}\label{exe:colimit-e-r-ps}
Denote by $\mathbf{Bif}_{\epp}$ the category of bifinite lattices with as morphisms the \epp's. Let $U \colon \mathbf{Bif}_{\epp} \to \DL$ be the functor that is the identity on objects, and sends a morphism $(e,p)$ to $e$. Let $D \colon I \to \mathbf{Bif}_{\epp}$ be a directed diagram. Prove that the colimit of the diagram $U \circ D$ is again a bifinite lattice.
\end{ourexercise}

\begin{ourexercise}\label{exe:details}
Fill in the details of the construction in Examples~\ref{ex:Scott} and~\ref{ex:lazylambda}.
\end{ourexercise}

\begin{ourexercise}\label{exer:two-omega}
  Let $2^{\leq \omega}$ be the set of pairs $(I, f)$, where $I$ is a down-set of
  the ordinal $\omega$ and $f$ is a function $I \to 2$. We put the strict prefix
  ordering on $2^{\leq \omega}$, that is, $(I,f) \leq (J,g)$ if, and only if, $I \subseteq J$ and $g|_I = f$. Show that $(2^{\leq \omega}, \leq)$ is a bifinite domain. 
\end{ourexercise}

\notessec
{The treatment of domain theory in this chapter focuses on this material as an application of Stone-Priestley duality and is based on this book's authors' interpretation of some of the material in \cite{Abr}. For more insight into denotational semantics and domain equations, we refer the interested reader to the specialized literature; we recommend in particular \cite{AbJu94}. Much more on $\lambda$-calculus and its models can be found for example in \cite{Ba84-2, BarMan2022}.}

{\cite{Abramsky87} 
was the first account of Domain Theory in Logical Form in
the setting of Scott domains based on the author's Ph.D. dissertation. 
Abramsky also already
announced there the conceptual straightforwardness of extending to the bifinite
setting modulo some technical intricacies to do with the identification of
primes, and that this extension was being worked on by Glynn Winskel's student
Guo-Qiang Zhang at the time. Indeed, the Ph.D. dissertation
\cite{UCAM-CL-TR-185}, later published as \cite{Zhang91}, takes among other
\cite{Abramsky87} as starting point and treats the extension to the bifinite
setting using the prevalent point of view in domain theory at the time.
Concurrently, \cite{Abr} also resolved the required
intricacies based on an improvement of earlier work on bifinite domain theory in the dissertation \cite{Gunter1985}.}

{The Hoffmann-Lawson duality was established independently by \cite{hoffmann81} and 
\cite{lawson79}. Note that this duality is not to be confused with the duality between distributive continuous lattices and locally quasicompact sober space developed around the same time in \cite{KHHofJDLaw78}; the one for completely distributive lattices involved R.~E. Hoffmann while the one for distributive continuous lattices involved K.~H. Hofmann, and both dualities have the same J.~D.~Lawson as their second-named author.}

{The results on completely distributive lattices presented in this Chapter are due to \cite{Raney53}, which introduced the relation $\wwb$ many years before the introduction of the relation $\wayb$ in \cite{Scott72}. Thus it is somewhat backwards to introduce $\wwb$ as a more specialized version of $\wayb$ as, given the chronology, $\wayb$ was in fact a generalization of the already existing $\wwb$.}

\chapter{Automata theory}\label{ch:AutThry}

In this chapter we discuss applications of duality theory to automata theory.
The theory of regular languages and automata is an original computer science
topic with a rich theory and a wide and still growing range of practical and
theoretical applications \parencite{HandbookI,HandbookII,Almetal2020}. The
availability of sophisticated mathematical tools from algebra and topology is
one of the main strengths of the classical framework. Finite algebras were
introduced into the theory early on by Myhill, and Rabin and Scott, and their
power was established by the effective characterization of star
free languages by means of their syntactic monoids
\parencite{Schutzenberger56,Sch65}. Syntactic monoids provide an abstract and
canonical notion of so-called \emph{recognition} for regular languages.
Eilenberg's Theorem \parencite{Eilenberg76} supplies a general framework in which to
apply the strategy of Sch\"utzenberger's result. 

The success of the algebraic method was greatly augmented by the introduction
in the 1980s of profinite topology and in particular profinite monoids, see
\cite{Almeida1995, Weil02,Almeida05,Pin2009}. A most powerful combination in this setting is that of Eilenberg's and
Reiterman's theorems. Reiterman's Theorem is a generalization of Birkhoff's
Variety Theorem from universal algebra. It states that pseudovarieties of
finite algebras are precisely the ones given by \emph{profinite identities}
\parencite{Reit82}.  Thus Eilenberg-Reiterman theory allows for the equational
description of certain classes of regular languages. This
can lead to decidable criteria for membership in the corresponding classes,
in case a finite equational basis can be found, or, at least, an infinite basis
that can be effectively tested on a finite algebra \parencite{AlmeidaWeil98}.

In this chapter, we will see the rudiments of this theory from a duality point
of view. Thus, recognition by monoids, and in particular syntactic monoids, will
be defined duality-theoretically. We will see that the dual space of the
Boolean algebra of regular languages over a given alphabet with certain
additional operations is the so-called \emph{free profinite monoid}. We will then 
show that the subalgebra-quotient duality of Section~\ref{sec:quotients-and-subs},
when applied to this pair, leads to the profinite equational
method of Eilenberg-Reiterman theory. We will not cover the Eilenberg-Reiterman
theory in full generality, but will see a particular example of it and we provide
references to the recent research literature throughout the chapter 
for readers who want to get the fully general picture.
More bibliographic references can also be found in the notes section at the end
of this chapter.

In Section~\ref{sec:syntmon}, after introducing basic notions from automata theory,
we show how the classical construction of a syntactic monoid for a regular
language can be understood via duality. In particular, we will see how this
is a special case of the general duality for implication operators developed
earlier in Section~\ref{sec:generalopduality}. The results in the first section
use \emph{discrete} duality in the case of free monoids and finite monoids,
from which we then also derive that the collection of all regular languages is a
Boolean algebra. Section~\ref{sec:freeprofmonoid} shows how to
extend these ideas to \emph{profinite} monoids, and establishes as
its main result that the \emph{free} profinite monoid is the dual of the Boolean
algebra of regular languages, enriched with some additional operations. We show
also how this result fits in a more general duality theorem for Boolean spaces
equipped with a binary operation.
Sections~\ref{sec:EilReittheory}~and~\ref{sec:openmult} show how this allows one
to apply the methods of duality theory in the setting of regular languages,
focusing in particular on characterizations of classes of languages via
profinite equations. A more detailed outline of the rest of the chapter is given at the
beginning of Section~\ref{sec:freeprofmonoid}, after we introduce some of the
basic notions that we will need.
\section{The syntactic monoid as a dual space}\label{sec:syntmon}
We started our introduction to duality theory back in Chapter~\ref{ch:order}
with the discrete duality for finite distributive lattices, see
Section~\ref{sec:finDLduality}. The finite case came first chronologically and,
while it is simpler than the full topological theory, it contains many of the
central ideas of the subject in embryonic form. A foundational result of
algebraic automata theory is the fact that one may associate, with each finite
state automaton, a finite monoid with a universal property. Our first goal in
this chapter is to show that this foundational result is a
case of discrete duality for complete and atomic Boolean algebras with
additional operations. This will then lead to the applications of duality theory
to the theory of regular languages that we look at later in the chapter. 
\subsection*{Automata, quotienting operations, and residuation}

We start with automata, and show how they naturally lead us to consider dual
Boolean algebras equipped with quotienting operations.
An automaton is a very simple and weak model of computation; 
it is a machine that reads once over a finite sequence of input values
and only in the order in which it is entered, and only produces one output bit at the end of its computation. Readers who have previously encountered notions of computability will notice that automata are very far from being Turing complete  (see for example Exercise~\ref{exer:RegA-not-complete}.\ref{itm:exp-not-reg}).  
We will typically 
denote the input alphabet by $A$ and will assume throughout that it is 
\emph{finite}. In fact, \emph{profinite} input alphabets are also natural in this setting, but 
we shall not need them in what we discuss in this chapter.

Recall that a \emphind{monoid} is a triple $(M, \cdot, 1)$, where $\cdot$ is an associative
binary operation on $M$ and $1$ is a \emph{neutral}\index{neutral element} or
\emphind{identity element} for $\cdot$, that is, $1 \cdot m = m = m \cdot 1$ for
all $m \in M$.  We will often omit
$\cdot$ in the notation, writing simply $mn$ instead of $m \cdot n$, as is
common practice in algebra. When multiple monoids $M$ and $N$ are in play, we  sometimes use the notations $1_M$ and $1_N$ to distinguish their neutral elements.
\nl{$1_M$}{neutral element of a monoid $M$}{oneM}
Note that, for an
input alphabet $A$, the set of all finite sequences over $A$ with binary
concatenation is the \emph{free monoid} generated by $A$, see
Example~\ref{exa:adjunctions}.\ref{itm:freemonoid} and
Exercise~\ref{exe:freemonoid}.\index{free monoid} Recall that we
denote this algebra by $A^*$. We suppress the commas in sequences, so that we
think of an element $w\in A^*$ as a \emphind{finite $A$-word}, often simply
called a \emph{word}. That is, a word is a possibly empty string 
$w = a_1 \ \dotsm\ a_n$, where $n \geq 0$ is a number called the \emph{length} of $w$, and $a_i \in A$ for $1 \leq i \leq n$.
\begin{notation}
  If $w \in A^*$ is a finite word, 
  as a notational convention, whenever we write
  $w = a_1 \ \dotsm \ a_m$, we mean that the $a_i$ are elements of the alphabet $A$,
  and thus that $m$ is equal to the length of $w$.
\end{notation}

We now define a notion of \emph{automaton}, also known as
\emph{non-deterministic finite state automaton} in the literature. 
\begin{definition}\label{def:FSA}
  A (finite state) \emph{automaton} \index{automaton} is a structure
  $\cA=(Q,A,\delta,I,F)$ where
  \begin{enumerate}
      \item \label{itm:statesdefFSA} $Q$ is a finite set, whose elements we call \emphind{states};
    \item $A$ is a finite set, called the \emph{(input) alphabet};
    \item \label{itm:transdefFSA} $\delta$ is a subset of $Q\times A\times Q$,
        called the \emph{transition relation};\index{transition relation}
    \item $I\subseteq Q$ is a set,  whose elements we call \emph{initial
        states};\index{initial state}
    \item $F\subseteq Q$ is a set, whose elements we call  \emph{final} or
        \emph{accepting} states.\index{final state}\index{accepting state}
  \end{enumerate}
  We will call a tuple of the form $\cM=(Q,A,\delta)$ satisfying
  (\ref{itm:statesdefFSA})--(\ref{itm:transdefFSA}), without the initial and finite
  states specified, a \emph{(finite state) machine}\index{machine}. Given an automaton
  $\cA=(Q,A,\delta,I,F)$, we call the associated machine $\cM=(Q,A,\delta)$ the \emphind{underlying
    machine} of $\cA$.

  Let $w = a_1 \ \dotsm\ a_n \in A^*$ with each $a_i\in A$. An automaton $\cA$
  \emphind{accepts} $w$ if there exists a sequence of states $q_0, \dots, q_n
  \in Q$ such that $q_0 \in I$, $q_n \in F$, and $(q_{i-1}, a_i, q_i) \in
  \delta$ for every $1 \leq i \leq n$; if no such sequence exists, $\cA$
  \emphind{rejects} $w$. In particular, note that $\cA$ accepts the empty word
  $\epsilon$ if, and only if, $I \cap F \neq \emptyset$. We denote by $L(\cA)$
  the set of words accepted by $\cA$, which is also known as the
  \emphind{language recognized by $\cA$}. A subset $L \subseteq A^*$ is called
  a \emphind{regular language} if there exists an automaton $\cA$ such that $L
  = L(\cA)$.%
  \footnote{The name \emph{regular language} originates from a
  different approach, via \emph{regular expressions}, that is, expressions,
  built up from single letters by finite union, finite concatenation product,
  and the \emph{Kleene star}, which yields the submonoid of $A^*$ generated by
  a subset of $A^*$. The concept we have defined is more aptly called
  \emph{recognizable language}. However, it is a non-trivial theorem, proved by
  Kleene at the very beginning of automata theory, that a language is
recognizable if, and only if, it is regular \parencite{Kleene1956}, and the
adjective regular is more commonly used, even if one focuses on recognition as
we do here.} 
\end{definition} 
The above definition allows for non-determinism:
when in state $q$, upon reading a letter $a \in A$, the automaton has the
possibility to move to any of the states in the set $\delta[q,a,\_]$. When this
set is always non-empty, the automaton is called \emph{complete}\index{complete
automaton}, and when this set has cardinality at most $1$, and the set of
initial states has cardinality $1$, the automaton is
called \emph{deterministic}\index{deterministic automaton}.\footnote{In some
references, the term `deterministic' includes `complete' by definition. Also
note that we can safely ignore the case where the set of initial states $I$ is empty.} 
Any
regular language of finite words can be recognized by a complete deterministic
automaton, as will also follow from the characterizations we give in
Theorem~\ref{thrm:reglanguage} below. We note, however, that determinism does
not come for free: Exercise~\ref{exe:non-det-example} gives an example of a
non-deterministic automaton which recognizes a language with $n$ states, but
any deterministic automaton recognizing the same language requires $2^n$
states.
The connection with duality that we expose here is natural for the general
notion of finite state machine that we defined above, and does
not require determinism a priori.

When $\cM = (Q,A,\delta)$ is a machine, 
for each $a \in A$, we will write 
\[\delta_a := \{(q,q') \in Q \times Q \ \mid \ (q,a,q') \in
\delta\},\] 
\nl{$\delta_x$}{reachability relation of a machine with transition relation $\delta$ when reading a word $x$}{}
which is a binary relation on $Q$. We extend this definition to words as follows.
For
any non-empty $w = a_1 \ \dotsm\ a_m \in A^*$, let $\delta_x$ be the relational
composition $\delta_{a_1} \cdot \ \cdots \ \cdot \delta_{a_m}$, and also define
$\delta_\epsilon$ to be the identity relation on $Q$. A different way to say
this is that $w \mapsto \delta_w$ is the unique monoid homomorphism from the
free monoid $A^*$ to the monoid of binary relations on $Q$ extending the
assignment $a \mapsto \delta_a$, see Exercise~\ref{exe:deltahom}. The relation $\delta_w$ 
is the \emphind{reachability relation} of the machine $\cM$ when reading a word $w \in A^*$. 
A word $w$
is then accepted by an automaton $\cA = (Q,A,\delta,I,F)$ based on $\cM$ if,
and only if, there exist $q_0 \in I$ and $q_f \in F$ such that $(q_0,q_f) \in
\delta_w$.
\begin{example}
  Consider the deterministic finite-state automaton
  \[
    \cA=(Q, A, \delta,\{1\},\{2\}),
  \]
  where $Q := \{1,2,3\}$, $A :=\{a,b\}$, and $\delta := \{(1,a,2),(2,b,3),(3,a,2)\}.$
  We will often depict such an automaton by a graph as follows:
  
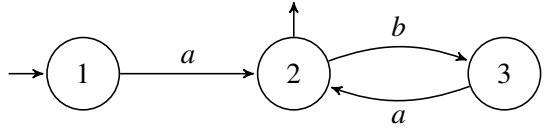
\begin{center}
    \begin{tikzpicture}[->,>=stealth',shorten >=1pt,auto,node distance=2.8cm,
        semithick]
      \tikzstyle{every state}=[fill=white,draw=black,text=black]

      \node[state] (A)                {$1$};
      \node[state] (B)  [right of=A]  {$2$};
      \node[state] (C)  [right of=B]  {$3$};

      \initarr{A}
      \node (FINALB) [above=0.5cm of B] {};
      \path (B) edge (FINALB);
      \path (A) edge              node  {$a$}  (B)
      (B) edge[bend left=20]  node  {$b$}  (C)
      (C) edge[bend left=20]  node  {$a$}  (B);

    \end{tikzpicture}
  \end{center}
  Here the nodes of the graph are the states, and the directed edges are labeled by letters and
  correspond to the elements of the transition relation. The initial states are indicated by an
  entering arrow while the final states are indicated by an arrow out. From the definition of
  acceptance, note that a word $w$ is accepted by this automaton $\cA$ if, and only if, $w$ begins
  with an $a$, followed by an arbitrary number of repetitions of the string $ba$. Thus, we have our
  first example of a regular language,
  \[ L(\cA) = \{a(ba)^n \mid n \geq 0\} = \{a, aba, ababa, abababa, \dots\}.\]
  The same transition relation $\delta$ on this set of states $Q$ may be used to recognize other
  regular languages; for example, the reader may verify that the automaton
  \[ \cA' := (Q, A, \delta, \{2\}, \{2, 3\})\]
  recognizes the language $L(\cA') = \{(ba)^n, (ba)^n b \mid n \geq 0\}$. Other regular languages
  can be obtained by choosing other sets of initial and final states. Crucially, however, if we fix
  any finite state machine $\mathcal{M} = (Q,A,\delta)$, then there are only finitely many choices for a pair of subsets $(I,F) \in
    \mathcal{P}(Q)^2$. Thus, only finitely many regular languages can be recognized by the automata
  based on a fixed machine $\mathcal{M}$. This essentially trivial observation will be important in
  what follows.
\end{example}
We will build on the idea of the above example to give an algebraic characterization of regular languages as
certain special elements of the Boolean algebra $\cP(A^*)$, enriched with additional operations.

\nl{$L(\cA)$}{language recognized by an automaton $\cA$}{}
\nl{$L(I,F)$}{language recognized by a finite state machine $\cM$ when taking $I$ as the set of initial states and $F$ as the set of final states}{}
\begin{definition}
  Let $\cA = (Q, A, \delta, I_0, F_0)$ be an automaton and let $L = L(\cA)$ and
  $\cM=(Q, A, \delta)$ the underlying machine of $\cA$. For any pair of subsets
  $(I, F) \in \cP(Q)^2$, denote by $L(I,F)$ the language recognized by the
  automaton $(Q, A, \delta, I, F)$; in particular, $L(\cA) := L(I_0, F_0)$. 
\nl{$\cR(\cM)$}{set of languages recognizable by a finite state machine $\cM$}{}
\nl{$\cL(\cM)$}{lattice of languages generated by $\cR(\cM)$}{}
\nl{$\cB(\cM)$}{Boolean algebra of languages generated by $\cR(\cM)$}{}
  We write
  \[\cR(\cM) := \{L(I,F) \ \mid \ (I,F) \in \cP(Q)^2\ \} \subseteq \cP(A^*)\]
  for the set of languages which can be recognized by automata obtained from
  $\cA$ by changing its initial and final states. We denote by $\cL(\cM)$ the
  \emph{sublattice} generated by $\cR(\cM)$ and by $\cB(\cM)$ the \emph{Boolean
  subalgebra} generated by $\cR(\cM)$. 
\end{definition}

\begin{remark}\label{cB(cA)finite}
  Let $\cM$ be a finite state machine. Note that the set $\cR(\cM)$ has at most
  $2^{2|Q|}$ elements and is thus in particular \emph{finite}. Therefore, $\cL(\cM)$ and $\cB(\cM)$ are both finite, since finitely generated distributive lattices and Boolean algebras are finite (see Exercise~\ref{exe:locallyfinite}). It follows that both  $\cL(\cM)$ and $\cB(\cM)$ are, trivially, complete sublattices of $\cP(A^*)$. 
\end{remark}

The idea that we will develop now is that the finite lattice $\cL(\cM)$ and the
finite Boolean algebra $\cB(\cM)$ can be
viewed as \emph{dual} algebraic structures  
associated to the machine $\cM$. 
The dynamic behavior of a machine $\cM$ is captured by the action of the
alphabet on the set of states. Taking into account the direction-reversing
nature of duality, it stands to reason that the appropriate structure to
consider on the dual collection of languages is given by 
\emph{quotienting} operations, defined, for any language $L \subseteq A^*$ and
letter
$a \in A$, by 
\nl{$x^{-1}L$}{left quotient of a language $L$ by the word $x$}{}
\nl{$Ly^{-1}$}{right quotient of a language $L$ by the word $y$}{}
\[ a^{-1}L := \{w \in A^* \ | \ aw \in L\}, \quad La^{-1} := \{w \in A^* \ | \
wa \in L\}\ ,\] 
or more generally, for any words $x, y \in A^*$, 
\[ x^{-1}L := \{w \in A^* \ | \ xw \in L\}, \quad Ly^{-1} := \{w \in A^* \ | \ wy
\in L\}\ .\]
We call the operations $L \mapsto x^{-1}L$ and $L \mapsto Ly^{-1}$ the left
and right \emphind{quotienting operations} on the Boolean algebra
of languages $\cP(A^*)$. Note that these operations are Boolean algebra
endomorphisms on $\cP(A^*)$; indeed, each of the operations computes the
inverse image of a function from $A^*$ to $A^*$.  These quotienting operations
may thus be seen
as a \emph{biaction} of the free monoid $A^*$ on the Boolean algebra
$\cP(\cA^*)$; we will further explore this
in Definition~\ref{def:biaction} below. 
For now, we say that a collection of languages
$\mathcal{L} \subseteq \cP(\cA^*)$ is \emphind{closed under quotienting} if,
for every $x,y \in A^*$ and $L \in \mathcal{L}$, we have that $x^{-1}L \in
\mathcal{L}$ and $Ly^{-1} \in \mathcal{L}$. 

A key insight is that the collection $\cR(\cM)$ of languages recognized by
a machine $\cM$ is closed under quotienting, as we will prove
now. Recall that, when $R$ is a relation from a set $X$ to $Y$, we write $R[U]$
for the \emphind{relational direct image} of $U$ under $R$, that is, $R[U] = \{y \in Y \
\mid \ uRy \text{ for some } u \in U\}$, and similarly $R^{-1}[V]$ for the
\emphind{relational inverse image}. We apply this notation here to the reachability
relations $\delta_x$ for $x \in A^*$ induced by a finite state machine.
\begin{lemma}\label{lem:autfinite}
  Let $\cM = (Q,A,\delta)$ be a machine. For any $x, y \in A^*$ and $(I,F) \in \cP(Q)^2$, we have
  \[ x^{-1} L(I,F) = L(\delta_x[I], F) \quad \text{ and }\quad L(I,F)y^{-1} = L(I, \delta_y^{-1}[F]).\]
  In particular, $\cR(\cM)$ is closed under quotienting.
\end{lemma}
\begin{proof}
Note that
$(q, q') \in \delta_x$ if, and only if, there exists a path from $q$ to $q'$
labeled by the letters of the word $x$. Using this characterization, 
for any $x, w \in A^*$ and $(I,F) \in \cP(Q)^2$, we now have that 
\[ w \in x^{-1}L(I,F) \text{ if, and only if, } xw \in L(I,F)\ , \] 
if, and only if, there exist
states $q_0, q_1, q_2 \in Q$ such that
\[q_0 \in I, (q_0,q_1) \in \delta_x, (q_1,q_2) \in \delta_w  \text{ and } q_2 \in
F\ .\] 
This last condition is clearly equivalent to the condition that the word $w$ is
in the language $L(\delta_x[I],F)$. This
proves the first equality. The proof of the second equality is similar.
\end{proof}
We have thus associated, with any finite state machine $\cM$, the set
$\cR(\cM)$ which, as Lemma~\ref{lem:autfinite} shows, is closed under the
quotienting operations.
We now show that it follows from this that the lattice $\cL(\cM)$ and the Boolean algebra
$\cB(\cM)$ are also closed under quotienting, which means that they can be considered
as algebraic structures in and of themselves, independently from the machine $\cM$.

\begin{proposition}\label{prop:LandB-closed-under-quotienting}
  Let $\cM=(Q, A, \delta)$ be a machine. Then $\cL(\cM)$ and
  $\cB(\cM)$ are closed under quotienting.
\end{proposition}
\begin{proof}
    Let $x \in A^*$. We will prove that $x^{-1}L \in \cB(\cM)$ for every $L \in
    \cB(\cM)$; the proofs of the other properties are essentially the same.
    Let $\cB_x$ denote
    the inverse image of $\cB(\cM)$ under the quotienting operation by $x$, that is, 
\[ \mathcal{B}_x := \{L \in \cP(A^*) \ | \ x^{-1}L \in \cB(\cM) \}\ . \]
    Then, since the quotienting operation is an endomorphism, $\cB_x$ is a Boolean subalgebra of
    $\cP(\cA^*)$, and it contains $\cR(\cM)$ by Lemma~\ref{lem:autfinite}.
    Thus, $\cB_x$ contains $\cB(\cM)$, which is what we needed to prove.  
\end{proof}
Combining Remark~\ref{cB(cA)finite} and
Proposition~\ref{prop:LandB-closed-under-quotienting}, we now aim to  
deduce stronger closure properties for $\cL(\cM)$ and $\cB(\cM)$, namely, they
are not just closed under quotienting by single words, but by any \emph{set} of
words. In order to express this closure property, which also plays a central role in the
more general theory of varieties of regular languages alluded to in the
introduction of this chapter, we need to introduce \emph{residuation} operations on
$\cP(M)$ for an arbitrary monoid $M$.

 For any monoid $M$, the Boolean algebra $\cP(M)$ carries a \emph{complex multiplication}, $\cdot$, defined by
\[ u \cdot t := \{\,mn \mid  m \in u, n \in t\,\}.\]
Since this operation $\cdot$ is a complete binary operator on a complete Boolean algebra, it has a left and a right residual, which we denote by $\backslash$ and $\slash$, respectively, see Example~\ref{exa:residualfromdot} and Exercise~\ref{exe:discrete-dual-res-family}. These are the two unique operations, $\backslash, \slash \colon \cP(M)^2 \to \cP(M)$, such that, for any $t,u,v \in \cP(M)$,
\[ u \cdot t \subseteq v\ \iff\ t\subseteq u\backslash v\ \iff\ u\subseteq v/t.  \]
More explicitly, in terms of our duality for implication-type operators of
Section~\ref{sec:generalopduality}, $\backslash$ and $\slash$ are the
implication-type operators obtained from the relation $R$ given by the monoid
multiplication on $M$ by using the first and the third, or, respectively, the first and the second coordinates as inputs:
\begin{align*}
  u \backslash v & = M-( R[u,\_\ ,M-v])=\{m \in M \ \mid \ km \in v \text{ for all } k \in u\},        \\
  v \slash t     & = M-( R[\_\ ,t,M-v])= \{m \in M \ \mid \ m\ell \in v \text{ for all } \ell \in t\}.
\end{align*}
In fact, as we will see below, if we want to witness a monoid operation $\cdot$ on a set $M$ on its discrete dual Boolean algebra $\mathcal{P}(M)$, then the residual operations $\backslash$ and $\slash$ are more important than the complex multiplication. This is not so surprising since duality should reverse the direction of a map.\label{pg:residual-defs}
Note that, starting with the following definition, we distinguish between
residuat\emph{ed} algebras, where the multiplication $\cdot$ is present as part
of the structure, and residuat\emph{ion}  algebras, where the multiplication
$\cdot$ is absent. While the difference is small in this section, we will  later
consider residuation algebras in more generality, see
Definition~\ref{def:residuation-alg}.
\begin{definition}\label{dfn:dual-residuated-ba}
  Let $M$ be a monoid. We call the Boolean algebra with three additional operations $(\cP(M), \cdot, \backslash, \slash)$ the \emph{(discrete) dual residuated Boolean algebra}\index{residuated Boolean algebra!dual to a monoid} of $M$. The reduct $(\cP(M), \backslash, \slash)$, a Boolean algebra with two additional operations, will be called the \emph{dual residuation algebra} of $M$\index{residuation algebra!dual to a monoid}. %

  A \emphind{residuation ideal} is a sublattice $\cL$ of $\cP(M)$ such that,
for any $r \in \cP(M)$ and $s \in \cL$, both $r \backslash s \in \cL$ and $s
\slash r \in \cL$. A residuation ideal is said to be \emph{Boolean} provided
the underlying sublattice is a Boolean subalgebra and it is said to be
\emph{complete} provided the underlying sublattice is a complete sublattice of
$\cP(M)$. 
\end{definition}
Residuation ideals are closed under intersections, so that for any collection
$S \subseteq \cP(M)$, there is a smallest residuation ideal containing $S$,
which we call the residuation ideal \emph{generated by}\index{residuation
ideal!generated} $S$. Note also that the collections of Boolean, complete, and
complete Boolean residuation ideals are also all closed under intersections and
thus there are also corresponding notions of generated complete and/or Boolean
residuation ideals (see Exercise~\ref{exe:res-ideals-intersection} and
Definition~\ref{dfn:res-ideal-gen}).

\emph{Quotienting operations} may now be defined for a general monoid $M$ by
restricting the residual
operations to the atoms of $\cP(M)$, as follows. %
\begin{definition}\label{def:actbyquot}
  Let $M$ be a monoid, $m\in M$, and  $s\in\cP(M)$. We write
  \[ m^{-1}s := \{m\}\backslash s, \quad s m^{-1} := s \slash \{m\},\]
  and we call $m^{-1}s$ and $sm^{-1}$ the \emph{left and right quotients} of $s$ by $m$, respectively. We will call operations $s \mapsto m^{-1} s$ and $s \mapsto s m^{-1}$ \emphind{quotienting operations} on the Boolean algebra $\cP(M)$.

  A subset $S \subseteq \cP(M)$ is \emphind{closed under quotienting} if, for
  every $s \in S$ and $m \in M$, $m^{-1}s \in S$ and $sm^{-1} \in S$. The
  \emphind{closure under quotienting} of a subset $S \subseteq \cP(M)$ is defined
  as the smallest subset $\cQ(S)$ that contains $S$ and is closed under quotienting; concretely, 
  \[ \mathcal{Q}(S) = \{ m^{-1} s n^{-1} \ \mid \ m, n \in M, s \in S\}.\]
  \nl{$\cQ(S)$}{closure under quotienting of a subset $S$ of the residuated Boolean algebra dual to a monoid}{}
\end{definition}
Note in particular that the operations $L
\mapsto x^{-1}L$ and $L \mapsto Ly^{-1}$, for $x, y \in A^*$ and $L \subseteq
A^*$, which we introduced in Lemma~\ref{lem:autfinite} above, are the
quotienting operations on $\cP(A^*)$ given by the monoid $A^*$. Being closed
under quotienting is clearly a consequence of being a residuation ideal; we now
prove that the converse holds in the complete setting.
\begin{proposition}\label{prop:residl-quotienting}
	Let $M$ be a monoid and let $\cL$ be a complete sublattice of $\cP(M)$. Then $\cL$ is a
	residuation ideal if, and only if, $\cL$ is closed under quotienting.
\end{proposition}
\begin{proof}
	Suppose $\cL$ is a complete sublattice of $\cP(M)$ that is closed under
	quotienting. Let $s \in \cL$ and $r \in \cP(M)$. Since $r = \bigcup_{m
	\in r} \{ m \}$, we have
	\[ r \backslash s = \big( \bigcup_{m \in r} \{m \} \big) \backslash s =
	\bigcap_{m \in r} m^{-1} s, \] 
	and the latter is in $\cL$ because $\cL$ is closed under quotienting and
	complete. In the same way, $s / r \in \cL$.
\end{proof}
Summarizing what we know about the lattice $\cL(\cM)$ and the Boolean algebra
$\cB(\cM)$ associated to a finite state machine $\cM$, we have the following.
\begin{proposition}\label{prop:residuationideal}
  Let $\cM=(Q, A, \delta)$ be any finite state machine. Then $\cL(\cM)$ and
  $\cB(\cM)$ are finite residuation ideals of $\cP(A^*)$. 
\end{proposition}
\begin{proof}
By Proposition~\ref{prop:LandB-closed-under-quotienting}, 
	$\cL(\cM)$ and $\cB(\cM)$ are closed under quotienting. 
	 By Remark~\ref{cB(cA)finite}, $\cL(\cM)$ and $\cB(\cM)$ are finite, hence
    complete, sublattices of $\cP(A^*)$. Proposition~\ref{prop:residl-quotienting}
	now gives the statement. 
\end{proof}

Proposition~\ref{prop:residuationideal} will allow us to express regularity of a
language $L$ in a way that is independent from any particular machine
recognizing $L$, only using the residuation structure on $\cP(A^*)$.  
\begin{definition}\label{dfn:res-ideal-gen}
  For any language $L \in \cP(A^*)$, we write $\cL(L)$ for the residuation ideal of $\cP(A^*)$ generated
by $\{L\}$ and $\cB(L)$ for the Boolean residuation ideal of $\cP(A^*)$ generated by $\{L\}$.
\end{definition}
\nl{$\cL(L)$}{the residuation ideal generated by a language $L$}{}
\nl{$\cB(L)$}{the Boolean residuation ideal generated by a language $L$}{}

\begin{corollary}\label{cor:reg-implies-finite}
  Let $L \in \cP(A^*)$ be regular. Then $\cL(L)$ and $\cB(L)$ are finite and thus complete residuation ideals of $\cP(A^*)$.
  Moreover, both are effectively computable from any automaton that recognizes $L$.
\end{corollary}
\begin{proof}
Let $\cA$ be any automaton that recognizes $L$, and write $\cM$ for its underlying machine. Then $L \in
\cL(\cM)$, so by
Proposition~\ref{prop:residuationideal}, we must have $\cL(L) \subseteq
\cL(\cM)$,  and the latter is finite. In particular, $\cL(L)$ can be computed effectively
by closing $\{L\}$ under quotienting and lattice operations of $\cL(\cM)$, which are given by Lemma~\ref{lem:autfinite}.
 The argument for $\cB(L)$ is the same. 
\end{proof}
We will see in Theorem~\ref{thrm:reglanguage} that the conditions in
Corollary~\ref{cor:reg-implies-finite} are in fact enough to conclude that $L$
is regular. For now, we see that $\cL(L)$ and $\cB(L)$ are finite algebraic
structures canonically associated with a regular language $L$. 
We will show in the next subsection that $\cB(L)$ is in fact dual to the
syntactic monoid of the regular language $L$.

\subsection*{Discrete duality for monoid quotients and residuation ideals}
We will now need to use some universal algebraic notions for monoids. 
The general definitions of universal algebra apply to monoids just as they
apply to lattices or Boolean algebras; we briefly recall the basics of
these definitions in the case of monoids and refer to, for example, \cite{BurSan2000,Wechler1992}, for more information.

A function $f \colon M \to N$ between monoids is called a \emph{homomorphism}\index{homomorphism!of monoids} if $f(xy) = f(x)f(y)$ for all $x, y \in M$, and $f(1_M) = 1_{N}$.
A \emphind{congruence}\index{monoid!congruence} on a monoid $M$ is an
equivalence relation $\theta$ such that, for any $m_1, m_1', m_2, m_2' \in M$,
if $m_i \ \theta \ m_i'$ for $i = 1,2$, then $m_1m_2 \ \theta \ m_1' m_2'$.
Note that it is equivalent to say that $\theta$ is invariant under
multiplication on both sides, that is, if $m, m' \in M$ with $m \ \theta \ m'$
and $x \in M$, then $xm \ \theta \ xm'$ and $mx \ \theta \ m'x$.
The \emph{quotient}\index{monoid!quotient}\index{quotient!of a monoid} of $M$ by the congruence $\theta$ is the monoid based on the quotient set $M/{\theta}$ uniquely defined by the requirement that $[m_1]_{\theta} [m_2]_{\theta} = [m_1m_2]_{\theta}$ for every $m_1, m_2 \in M$; its neutral element is $[1_M]_{\theta}$. For any homomorphism $f \colon M \to N$ between monoids, the \emphind{kernel} of $f$, $$\ker(f) := \{(m,m') \mid f(m) = f(m')\},$$ is a congruence on $M$, and $f$ can be factored as the surjective homomorphism $M \onto M/{\ker(f)}$ that sends $m$ to $[m]_{\ker(f)}$ followed by an injective homomorphism $M/{\ker{f}} \into N$ well-defined by requiring that $[m]_{\ker(f)}$ is sent to $f(m)$, for every $m \in M$. This construction is known as the \emph{first isomorphism theorem} for monoids.\index{first isomorphism theorem} \index{monoid!first isomorphism theorem}
In particular, for any surjective monoid homomorphism $M \onto N$, $N$ is isomorphic to $M/{\ker(f)}$, so that surjective homomorphisms with domain $M$ may be identified with congruences on $M$. As in the case of distributive lattices (see Exercise~\ref{exe:monoepiconcrete}), surjective homomorphisms in the category of monoids are epimorphisms, but the converse is not true (see Exercise~\ref{ex:nonsurjective-epi-mon}).

A \emph{semigroup}
is a pair $(S, \cdot)$ where $\cdot$ is an associative binary operation on $M$.
Thus, ``semigroup = monoid minus neutral element''.\footnote{We will not be
concerned much with semigroups in this chapter, although most of the results
that we state and prove for monoids also hold for semigroups.} %
A \emphind{subsemigroup} of a monoid $M$ is a subset of $M$ that is 
closed under  multiplication. A subsemigroup is a \emphind{submonoid} 
if it moreover contains $1$. Notice that a subsemigroup of a monoid $M$ 
may be a monoid with a different identity element than $M$ (see Exercise~\ref{exe:monoidbasic}). 
If $M$ is a monoid and $S$ is a subset, there exist a smallest subsemigroup 
and a smallest submonoid of $M$ that contain $S$, these are called the 
\emph{subsemigroup generated by $S$}\index{generated subsemigroup} and 
the \emph{submonoid generated by $S$}\index{generated submonoid}.

A cornerstone result in the algebraic theory of regular languages establishes a
close connection between languages recognized by finite automata and
finite quotients of free monoids. Our aim in the rest of this section is
to show that this connection is an instance of discrete duality. 

We saw in Section~\ref{sec:quotients-and-subs},
Theorem~\ref{thrm:BA-sub-discrete}, that quotients of a \emph{set} correspond
to complete Boolean subalgebras of its power set. We are now in a situation
where the set $M$ additionally carries the structure of a monoid, and we will
show in Theorem~\ref{thrm:monoid-quotient-duality} that the monoid quotients of
$M$ correspond dually to complete Boolean residuation ideals of $\cP(M)$. The
first step towards our proof of Theorem~\ref{thrm:monoid-quotient-duality} will
be to view monoid quotients as biactions; see Proposition~\ref{prop:biactions}
below. We will then use this theorem to characterize regular languages. We
first give the necessary algebraic definitions.
\begin{definition}\label{def:biaction}
  Let $M$ be a monoid and $S$ a set. A \emphind{left action}\index{action!left} of $M$ on $S$ is a function $\lambda \colon M \times S \to S$ such that, for every $m, m' \in M$ and $x \in S$, we have
  \[\lambda(m'm,x) = \lambda(m',\lambda(m,x)) \text{ and } \lambda(1_M,x) = x.\]
  A left action may also be viewed as a homomorphism from $M$ to the monoid of functions $S \to S$ under composition, by sending $m \in M$ to the function $\lambda_m \colon S \to S$ defined by $\lambda_m(x) := \lambda(m,x)$; the left action is called \emph{faithful} if $m \mapsto \lambda_m$ is injective. If the left action is clear from the context, we write $m \cdot x$ for $\lambda(m,x)$, noting that the action law then gives that $(m'm) \cdot x = m' \cdot (m \cdot x)$. If  $M$ has left actions on two sets $S$ and $S'$, then a function $h \colon S \to S'$ is a \emph{morphism} if, for any $m \in M, x \in S$, we have $m \cdot_S h(x) = h(m \cdot_{S'} x)$. The notion of a \emphind{right action}\index{action!right} of $M$ on $S$ is defined symmetrically.

  A \emphind{biaction}\index{action!two-sided} of $M$ on $S$ is a pair $(\lambda, \rho)$, where $\lambda$ is a left action, $\rho$ is a right action, and, for any $m, m' \in M$ and $x \in S$, we have
  \[\rho(\lambda(m,x),m') = \lambda(m,\rho(x,m')).\]
  This axiom implies that, when the pair $(\lambda, \rho)$ is clear from the context, we can unambiguously write $m \cdot x \cdot m'$ for the result of acting on $x$ by $m$ on the left and by $m'$ on the right, with the actions being performed in whichever order. A biaction can thus be described alternatively as a function $\alpha \colon M \times S \times M \to S$ for which $\alpha(-,-,1_M)$ is a left action, $\alpha(1_M,-,-)$ is a right action, and $\alpha(m,x,m') = \alpha(1_M,\alpha(m,x,1_M),m') = \alpha(m,\alpha(1_M,x,m'),1_M)$ for every $m, m' \in M$ and $x \in S$. A \emphind{biaction morphism} is a function that is an action morphism for both the induced left and right actions.
\end{definition}
\begin{example}\label{exa:Cayley}
  The \emphind{Cayley representation} of a monoid $M$ is given by letting it act on its own underlying set: for any $m, x \in M$, define $\lambda(m,x) := mx$. This left action is faithful, since $\lambda_m(1_M) = \lambda_{m'}(1_M)$ implies $m = m'$. Similarly, the monoid $M$ acts on itself on the right, by setting $\rho(x,m) := xm$ for any $x, m \in M$. The associativity axiom of the monoid says that $(\lambda, \rho)$ is a biaction, and we call this the \emphind{Cayley biaction} of $M$ on itself.

  Applying discrete duality to the Cayley functions $\lambda_m$ and $\rho_m$, for each $m \in M$, we obtain two complete homomorphisms $\lambda_m^{-1}, \rho_m^{-1} \colon \cP(M) \to \cP(M)$. These two functions then give a biaction of $M$ on the set $\cP(M)$, concretely given by:
  \begin{align*}
    M\times\cP(M)\times M & \to \cP(M),                                 \\
    (m,s,n)               & \mapsto \{ x \in M \mid mxn \in s \}.
  \end{align*}
  This biaction of $M$ on $\cP(M)$ corresponds to looking at residuation with respect to \emph{singleton sets}. As long as we stick to the discrete duality, it is equivalent whether one considers residuation by general sets, or the biaction given by  residuation with respect to singleton sets.
  This fact is exploited in the development below, in particular in the proof of Theorem~\ref{thrm:monoid-quotient-duality}.
\end{example}
We now show that biactions, taken through discrete duality, give biactions.
\begin{proposition}\label{prop:biaction-dual}
  Let $M$ be a monoid and suppose $M$ has a biaction on $S$. Then $M$ has a biaction on $\cP(S)$ given, for $m_1, m_2 \in M$ and $u \in \cP(S)$, by  
  \[ m_1^{-1} \ u \ m_2^{-1} := \{ x \in S \ \mid \ m_1 \cdot x \cdot m_2 \in u
  \},\]
  and the map $u \mapsto m_1^{-1} u m_2^{-1}$ is a complete Boolean algebra
  homomorphism.
  The induced left action is $u \mapsto u m_2^{-1}$ and the induced right action is
  $u \mapsto m_1^{-1}u$. If $M$ also has a biaction on $T$ and $h \colon S \to
  T$ is a biaction morphism, then $h^{-1} \colon \cP(T) \to \cP(S)$ is a
  biaction morphism for the dual biactions. In particular, for any $m_1, m_2 \in
  M$, the quotienting operation $\cP(M) \to \cP(M)$ given by $u \mapsto
  m_1^{-1} \ u \ m_2^{-1}$ is a complete Boolean
  algebra homomorphism. 
\end{proposition}
\begin{proof}
  The fact that $(m,x)\mapsto m\cdot x$ is a \emph{left} action of $M$ on $S$ corresponds via duality to the fact that $(u,m)\mapsto m^{-1}u$ is a \emph{right} action of $M$ on $\cP(S)$, and symmetrically on the other side. Note that for this reason we write the dual of the right action on the left and the dual of the left action on the right. The compatibility of the two actions on $S$ is dual to the compatibility of the two actions on $\cP(S)$, and likewise the dual $h^{-1}$ of $h$ is a biaction morphism (see Exercise~\ref{exe:biaction-duality} for details). The `in particular' statement is the application of the proposition to the Cayley biaction of $M$ on itself. 
\end{proof}
We note a consequence of Proposition~\ref{prop:biaction-dual}, which generalizes the
argument given in Proposition~\ref{prop:LandB-closed-under-quotienting}.

\begin{proposition}\label{prop:sublat-quotienting}
	Let $M$ be a monoid and $S \subseteq \cP(M)$. The (complete) sublattice generated
	by $\cQ(S)$ and the (complete) Boolean subalgebra generated by $\cQ(S)$
	are all closed under quotienting. In particular, for any regular language $L$,  $\cL(L)$ is the sublattice of $\cP(A^*)$ generated by $\cQ(\{L\})$ and $\cB(L)$ is the Boolean subalgebra of $\cP(A^*)$ generated by $\cQ(\{L\})$.
\end{proposition}

\begin{proof} Denote by $\cL$ the sublattice generated by $\cQ(S)$.  Let $m_1,
	m_2 \in M$ be arbitrary. Since the function $h \colon u \mapsto m_1^{-1}
	u m_2^{-1}$ is a lattice homomorphism, the inverse image $h^{-1}(\cL)$
	is a sublattice of $\cP(M)$. Also note that $\cQ(S) \subseteq
	h^{-1}(\cQ(S)) \subseteq h^{-1}(\cL)$, and thus $\cL \subseteq
h^{-1}(\cL)$. We conclude that, for any $u \in \cL$, $m_1^{-1} u m_2^{-1} \in
\cL$. Thus $\cL$ is closed under quotienting. The proofs for the generated complete
sublattice, the generated
Boolean subalgebra, and the generated complete Boolean subalgebra are analogous.
The in particular statement now follows by the proof of Corollary~\ref{cor:reg-implies-finite}.
\end{proof}

In the following proposition, we will consider a monoid $M$, a set $S$, and a
surjective function $h \colon M \onto S$ from $M$ onto a set $S$. Note that, since
$h$ is surjective, there is at most one binary operation $\cdot_S$ on $S$ with
respect to which $h$ becomes a homomorphism, namely, whenever $s_1, s_2 \in S$,
we can pick $m_i \in M$ such that $s_i = h(m_i)$, and then $s_1 \cdot_S s_2$
must be $h(m_1 m_2)$. This operation $\cdot_S$ is well-defined if, and only if,
the kernel of $h$ is a monoid congruence. We will now show that this is also
equivalent to having a biaction of $M$ on $S$ that is respected by the function
$h$.
\begin{proposition}\label{prop:biactions}
Let $M$ be a monoid, $S$ a set, and let $h \colon M \onto S$ be a surjective
function. The following are equivalent:
  \begin{enumerate}[label=(\roman*)]
    \item there exists a well-defined monoid structure on $S$ such that $h$ is a homomorphism;
    \item there exists a biaction of $M$ on $S$ such that  $h$ is a morphism
	    from the Cayley biaction on $M$ to the biaction  on $S$.
  \end{enumerate}
\end{proposition}
\begin{proof}
Assume (i) holds. Define the biaction
of $M$ on $S$ by
\begin{equation}\label{eq:alphah}
  \alpha_h \colon M \times S \times M \to S, \quad \alpha(m_1,s,m_2) := h(m_1) \cdot_S s \cdot_S h(m_2).
\end{equation}
 We first show that $(m,s) \mapsto \alpha_h(m,s,1_M)$ is a left action. Indeed,
 for any $m, m' \in M$ and $s \in S$ we have
  \[ \alpha_h(m'm,s,1_M) = h(m')h(m)s = \alpha_h(m',\alpha_h(m,s,1_M),1_M),\]
  and also
  \[ \alpha_h(1_M,s,1_M) = h(1_M)sh(1_M) = s.\]
  Symmetrically, $(s,m) \mapsto \alpha_h(1_M,s,m)$ is a right action. By a similar calculation, $\alpha_h$ is a biaction.
  Finally, $h$ is a morphism of biactions, since, for the left actions, for any $m, m' \in M$, we have
  \[ h(m \cdot_M m') = h(mm') = h(m) h(m') = \alpha_h(m, h(m'),1_M),\]
  and similarly for the right action.

Now assume (ii) holds.   Then, for any $m_1, m_2 \in M$, we have
  \[ h(m_1m_2) = h(\lambda_{m_1}(m_2)) = \alpha(m_1,h(m_2),1_M), \] 
and similarly $h(m_1m_2) = \alpha(1_M,h(m_1),m_2).$
  From this, it follows that if $(m_1,m_1') \in \ker(h)$ and $(m_2,m_2') \in
  \ker(h)$, then we must have 
  \[ h(m_1m_2) = \alpha(m_1,h(m_2),1_M) = h(m_1m_2') = \cdots = h(m_1'm_2').\]
  Thus, $\ker(h)$ is a monoid congruence, from which (i) follows by the remarks
  preceding the proposition.
\end{proof}

Propositions~\ref{prop:biactions}~and~\ref{prop:biaction-dual} together imply
that surjective monoid homomorphisms $M \onto N$ are dual to maps $\cP(N) \into
\cP(M)$ that preserve the
quotienting operations on the power set algebras of $M$ and $N$. We now
characterize the images of such maps as the complete Boolean residuation ideals
of $\cP(M)$.
\begin{theorem}\label{thrm:monoid-quotient-duality}
  Let $M$ be a monoid and $(\cP(M), \cdot, \backslash, \slash)$ its discrete dual residuated Boolean algebra.
  Let $h\colon M\twoheadrightarrow N$ be a surjective function from $M$ to a monoid $N$, and denote by $\cB$ the image of $h^{-1} \colon \cP(N) \into \cP(M)$. Then
  $h$ is a monoid homomorphism if, and only if, $\cB$ is a residuation ideal of $\cP(M)$. %
  In particular, under discrete duality, monoid quotients of $M$ correspond to complete Boolean residuation ideals of $\cP(M)$.
\end{theorem}

\begin{proof}
  First suppose $h$ is a monoid homomorphism. By
  Proposition~\ref{prop:biactions}, $h$ induces a biaction of $M$ on $N$, and by
  Proposition~\ref{prop:biaction-dual}, $h^{-1}$ is a morphism for the dual
  biaction of $M$ on $\cP(N)$, and thus the image of $h^{-1}$ is closed
  under quotienting. By Proposition~\ref{prop:residl-quotienting}, it is
  therefore a residuation ideal.

  For the converse, suppose that the image $\cB$ of $h^{-1}$ is a residuation ideal. Then in
  particular we get a biaction of $M$ on $\cP(N)$, as follows. For any $m, m' \in
  M$ and $u \in \cP(N)$, let $m \cdot u \cdot m'$ be the unique element of $\cP(N)$ such that
\[ h^{-1}(m \cdot u \cdot m') = m^{-1} h^{-1}(u) (m')^{-1} \ .\] 
Note that the latter is indeed in $\cB$ because it is a
  residuation ideal. The dual of this biaction of $M$ on $\cP(N)$ now gives a biaction of $M$ on $N$ for
  which $h$ is a biaction morphism (see Exercise~\ref{exe:biaction-duality}), so Proposition~\ref{prop:biactions} implies
  that $h$ is a monoid homomorphism. 
\end{proof}
In summary, Theorem~\ref{thrm:monoid-quotient-duality} shows that the anti-isomorphism of Theorem~\ref{thrm:BA-sub-discrete}, between equivalence relations on a set $M$ and complete Boolean subalgebras of $\cP(M)$, restricts to an anti-isomorphism between monoid congruences on $M$ and complete Boolean residuation ideals of $\cP(M)$. For any surjective monoid homomorphism $h \colon M \onto N$, we define the \emph{corresponding residuation ideal}\index{residuation ideal!corresponding to a monoid quotient} of $\cP(M)$  to be
\[
  \cB:=\{h^{-1}(u)\mid u\subseteq N\}.
\]
Conversely, when $\cB$ is a complete Boolean residuation ideal of $\cP(M)$, with $M$ a monoid, we have a unique congruence $\theta$ on $M$ for which $\cB$ is the corresponding residuation ideal; we call $\theta$ the \emph{congruence dual to} $\cB$. Recall from the proof of Theorem~\ref{thrm:BA-sub-discrete} that this congruence $\theta$ may be defined explicitly, for $x, y \in M$, by
\[ x \theta y \iff \text{for all } a \in \cB, x \in a \text{ if and only
if } y \in a.\]

The following definition, namely that of recognition from automata theory, codifies a slight generalization of the duality relationship of Theorem~\ref{thrm:monoid-quotient-duality}, in that $h$ need not be surjective.
\begin{definition}
  Let $M$ be a monoid, $s\in\cP(M)$, and $h\colon M\rightarrow N$ a monoid homomorphism. We say that $s$ is \emphind{recognized} by $h$ provided there is $u\subseteq N$ such that $s=h^{-1}(u)$; in this case, we call $h$ a \emphind{recognizer} for $s$. That is, $s$ is recognized by $h$ if, and only if, it belongs to the residuation ideal corresponding to the surjective homomorphism $M\twoheadrightarrow \im(h)$. Further, we say that a monoid $N$ \emphind{recognizes} $s$ provided there is a monoid homomorphism $h\colon M\rightarrow N$ which recognizes $s$.
\end{definition}
When $h \colon M \rightarrow N$ is a monoid homomorphism, we call the surjective function $h \colon M \onto \im(h)$, which is
still a homomorphism, now onto the submonoid $\im(h)$ of $N$, the \emphind{surjective co-restriction} of the homomorphism $h$.

It follows from Theorem~\ref{thrm:monoid-quotient-duality}, combined with Theorem~\ref{thrm:BA-sub-discrete}, that we always have a least recognizer for any $s\in\cP(M)$, namely the dual of the complete Boolean residuation ideal generated by $s$.

\begin{corollary}\label{cor:syntactic}
  Let $M$ be a monoid and let $s\in\cP(M)$. Then the complete Boolean residuation ideal generated by $s$ corresponds to the monoid quotient
  \[
    h_s\colon M \twoheadrightarrow M/{\equiv_s},
  \]
  where $m\equiv_s m'$ if, and only if,
  \[
	  \text{for all } x,y\in M \qquad (xmy\in s \iff  xm'y\in s).
  \]
  This is the least recognizer of $s$, in the sense that it factors through the surjective co-restriction of any other recognizer.
\end{corollary}

The first statement of Corollary~\ref{cor:syntactic} follows from the fact 
that the complete Boolean residuation ideal generated by $s$ is equal to the smallest complete Boolean
subalgebra containing $\cQ(\{s\})$, by
Propositions~\ref{prop:residl-quotienting}~and~\ref{prop:sublat-quotienting}.
We leave the details of the rest of the Corollary as
Exercise~\ref{exe:syntactic}. 

\begin{definition}\label{def:syntactic}
  Let $M$ be a monoid, $s\in\cP(M)$. Then the congruence $\equiv_s$, defined by $m\equiv_s m'$ if, and only if,
  \[
	  \text{for all } x,y\in M, \qquad xmy\in s \ \iff\ xm'y\in s
  \]
  is called the  \emphind{syntactic congruence} of $s$. The quotient map
  \[
    h_s\colon M \twoheadrightarrow M/{\equiv_s}
  \]
  is called the  \emphind{syntactic homomorphism} of $s$, and $M_s=M/\equiv_s$ is called the  \emphind{syntactic monoid} of $s$.
\end{definition}
\nl{$\equiv_s$}{syntactic congruence of a subset $s$ of a monoid $M$}{}
\nl{$M_s$}{syntactic monoid of a subset $s$ of a monoid $M$}{}
We are now ready to prove the main theorem of the section. {The equivalence of (i), (iv) and (v) in Theorem~\ref{thrm:reglanguage} is a classical result in the automata theory literature, but our point of view here, via duality for residuation ideals, is more recent~\parencite{GGP2008}.}
\index{regular language!has finite syntactic monoid}
\index{regular language!generates finite residuation ideal}
\begin{theorem}\label{thrm:reglanguage}
  Let $L \in \cP(A^*)$. The following conditions are equivalent:
  \begin{enumerate}[label=(\roman*)]
    \item the language $L$ is regular;
    \item the residuation ideal $\cL(L)$ is finite;
    \item the Boolean residuation ideal $\cB(L)$ is finite;
    \item the syntactic monoid $M_L$ of $L$, is finite;
    \item the language $L$ is recognized by a finite monoid.
  \end{enumerate}
Moreover, if $L$ is regular then the syntactic monoid is effectively computable
from any automaton recognizing $L$.
\end{theorem}

\begin{proof}
  If $L$ is regular, then, by Corollary~\ref{cor:reg-implies-finite}, $\cL(L)$
  is finite. If $\cL(L)$ is finite, then $\cB(L)$ is finite as it is the Boolean
  algebra generated by $\cL(L)$. If  $\cB(L)$ is finite, then its dual set is
  finite, and this is the syntactic monoid $M_L$, by
  Corollary~\ref{cor:syntactic}. Note that this dual of $\cB(L)$ is indeed
  effectively computable, since $\cB(L)$ can be computed from any automaton recognizing $L$ 
  by Corollary~\ref{cor:reg-implies-finite}, and its dual is the set of atoms of $\cB(L)$ equipped with the multiplication
  dual to the quotienting operations. If $M_L$ is finite, then as $L$ is recognized by $M_L$, it is recognized by a finite monoid. Finally, if $L$ is recognized by a finite monoid $N$, then there is a monoid homomorphism  $h\colon A^*\rightarrow N$ and $u\subseteq N$ so that $L=h^{-1}(u)$. We define an automaton $\cA = (N, A, \delta, \{1_N\}, u)$ by setting
  \[
    (n,a,n')\in\delta\quad\iff\quad nh(a)=n'.
  \]
  Then, for $x=a_1\ \dotsm\ a_n\in A^*$, we have that $\cA$ accepts $x$ if, and only if, $h(x)=h(a_1)\ \dotsm\ h(a_n)=1_Nh(a_1)\ \dotsm\ h(a_n)\in u$ if, and only if, $x\in L$.
\end{proof}
\nl{$\RegA$}{Boolean algebra of regular languages over a finite alphabet $A$}{}
\begin{definition}\label{def:RegA}
  Let $A$ be a finite set. We denote by $\RegA$ the subset of $\cP(A^*)$ consisting of the regular languages over $A$.
\end{definition}
Using the characterization of regular languages in 
Theorem~\ref{thrm:reglanguage}, we can now show that the set $\RegA$ of all regular languages in fact carries a lot of algebraic structure. {Apart from its interest in computer science, the infinite Boolean algebra $\RegA$ is of interest in constructive mathematics as it is a non-trivial countably infinite subalgebra of a power set algebra which exists constructively.}
\begin{proposition}\label{prop:RegBoolResidl}
  The set $\RegA$ over $A$ is a Boolean residuation ideal of $\cP(A^*)$.
\end{proposition}
\begin{proof}
  Let $L$ be a regular language over $A$. We have seen, in Theorem~\ref{thrm:reglanguage}, that $L$ is recognized by a finite monoid quotient $h\colon A^*\twoheadrightarrow M$. That is, $L=h^{-1}(u)$ where $u:=h[L] \subseteq M$ is the direct image of $L$ under $h$. It follows that the complementary language $L^c=A^*\setminus L$ is also recognized by $h$, via the set $u^c=M\setminus u$. Thus the set of regular languages is closed under complementation. To see that it is closed under intersection, notice that if $L$ is recognized by  $h\colon A^*\twoheadrightarrow M$ and $K$ is recognized by  $g\colon A^*\twoheadrightarrow N$ via $v\subseteq N$, then $L\cap K$ is recognized by the product map
  \[
    \gen{h,g}\colon A^* \to M\times N, \quad w\mapsto (h(w),g(w))
  \]
  via the subset $u \times v$.

It remains to show that $\RegA$ is a residuation ideal. For any $L \in \RegA$, the residuation ideal $\cL(L)$ is finite by Theorem~\ref{thrm:reglanguage}. Thus, for any language $K \in \cP(A^*)$, the languages $K \backslash L$ and $L/K$ are in the finite residuation ideal $\cL(L)$, so that the residuation ideals $\cL(K \backslash L)$ and $\cL(L/K)$ are contained in $\cL(L)$ and thus also finite. By Theorem~\ref{thrm:reglanguage}, $K \backslash L$ and $L/K$ are regular, as required.
\end{proof}

Notice that even though we took $h$ and $g$ in the proof of Proposition~\ref{prop:RegBoolResidl} to be surjective, the monoid homomorphism $\gen{h,g}$ is not necessarily surjective.  Of course we could have gotten a surjective monoid morphism by looking at the surjective reflection of the product map. This is one of the reasons why it is more natural to relax the notion of recognition to monoid morphisms rather than insisting on surjections. Also note that, by taking products of homomorphisms, as in the proof of Proposition~\ref{prop:RegBoolResidl}, we obtain the following corollary.
\begin{corollary}\label{cor:reg-directed}
  If $L_1,\dots,L_n$ are regular languages over $A$, then there exists a finite monoid quotient $h\colon A^*\twoheadrightarrow M$ that recognizes all the languages in the Boolean subalgebra generated by $L_1,\dots,L_n$.
\end{corollary}

In the next section, we will will apply \emph{topological} duality to the Boolean algebra with operators $(\RegA, \slash, \backslash)$.

\ourexercises
\begin{ourexercise}\label{exe:singletonlanguages}
  Let $A$ be a finite alphabet and $u\in A^*$. Show that the singleton $\{u\}$ is a regular language by exhibiting a finite state automaton which recognizes this singleton language. Conclude that the Boolean algebra of all finite and cofinite subsets of $A^*$ is contained in $\RegA$.
\end{ourexercise}
\begin{ourexercise}\label{exe:non-det-example}
    Let $A$ be a finite alphabet. For $w \in A^*$, write $c(w)$ for the set of
    letters that occur in $w$. Consider the language
    \[ L = \{ w \in A^* \ \mid \ c(w) \subsetneq A \}.\] 
    \begin{enumerate}
        \item Give a finite automaton with $|A|$ states that recognizes
            $L$. \hint{In each state, verify that a particular letter does not
            appear in the word; make all states initial.} 
        \item Show that $c \colon A^* \to \mathcal{P}(A)$ is a monoid homomorphism, where the monoid operation on $\mathcal{P}(A)$ is $\cup$, with neutral element $\emptyset$.
        \item Prove that $c$ is, up to isomorphism, the syntactic homomorphism of $L$.
        \item Conclude that $L$ can not be recognized by a monoid that has fewer than $2^{|A|}$ elements. 
    \end{enumerate}
    \emph{Note.} The syntactic monoid of a regular language $L$ is in general isomorphic to the monoid of transitions of the so-called \emphind{minimal automaton} of $L$, see for example \cite[Ch.~1, Sec.~4.3]{HandbookI}.
\end{ourexercise}
\begin{ourexercise}\label{exer:RegA-not-complete}
Let $A$ be a one element alphabet.
\begin{enumerate}
\item  Show that $A^*$ is isomorphic to the monoid $(\bN,+)$. 
\end{enumerate}
{In the rest of this exercise, we will identify the monoid $A^*$ with $(\bN, +)$.}
\begin{enumerate}[resume]
\item Show that the finite quotients of $(\bN,+)$ are given by the congruences $\theta_{N,q}$, where $N,q\in\bN$ and 
\[
\theta_{N,q}:=\{(n,n+kq) \ \mid \ N \leq n\text{ and } k\in\bN\}\ .
\]
\item Conclude that if $M$ is a finite monoid that is generated by a single element, then $M$ is isomorphic to $\bN/{\theta_{N,q}}$ for some $N, q \in \bN$. 
\end{enumerate}
\emph{Note.} An important fact in the theory of finite monoids, closely related to the preceding two items, is that every finite semigroup generated by a single element contains a unique idempotent; see the proof of Proposition~\ref{prop:omega}, which contains hints for this exercise.
\begin{enumerate}[resume]
\item Show that the collection of regular languages in $\bN$ is generated as a Boolean algebra by the singleton languages and the languages of the form 
\[ q\N+r := \{qn + r \ \mid \ n \in \bN \} \quad \text{ for } q,r\in N. \] 
\item \label{itm:exp-not-reg} Show that the language $\{2^n\mid n\in\bN\}$ is not regular.
\item Show that $\RegA$ is not complete. 
\hint{Use Exercise~\ref{exe:counterexamples-complete}.}
\end{enumerate}
\end{ourexercise}

\begin{ourexercise}\label{exe:even}
  Let $A$ be a finite set and $\mathbb Z_2$ be the monoid of integers modulo $2$ under addition. Consider the finite state automaton
  \[\mathcal A=(\mathbb Z_2, A,\delta, \{0\}, \{0\}), \text{where } \delta=\{(0,a,1),(1,a,0)\mid a\in A\}.\]
  Show that $\mathcal A$ recognizes the language consisting of all words in $A^*$ of even length.
\end{ourexercise}

\begin{ourexercise}\label{exe:monoidbasic}
  Prove that the identity element of a monoid is unique. Give an example of a monoid $M$ with a subsemigroup $N$ that is not a submonoid of $M$, but such that $N$ is a monoid in its own right.
\end{ourexercise}
\begin{ourexercise}\label{exe:res-ideals-intersection}
  Let $M$ be a monoid. Suppose that $\cJ$ is a collection of residuation ideals of $\cP(M)$. Prove that $\bigcap \cJ$ is again a residuation ideal. Prove that if, moreover, the residuation ideals in $\cJ$ are all Boolean and/or complete, then the same is true for $\bigcap \cJ$.
\end{ourexercise}
\begin{ourexercise}\label{exe:biaction-duality}
  This exercise fills in the details of Proposition~\ref{prop:biaction-dual}. Let $M$ be a monoid, $S$ a set, and let $\lambda \colon M \times S \to S$ and $\rho \colon S \times M \to S$ be functions. For any $m \in M$, write $\lambda_m(s) := \lambda(m,s)$ and $\rho_m(s) := \rho(s,m)$. Let $r \colon \cP(S) \times M \to \cP(S)$ and $\ell \colon M \times \cP(S) \to \cP(S)$ be the functions defined, for any $m \in M$ and $u \in \cP(S)$, by 
  \[ r(u,m) := \lambda_m^{-1}(u), \quad \ell(m,u) := \rho_m^{-1}(u).\]
  \begin{enumerate}
    \item Prove that $r$ is a right action of $M$ on $\cP(S)$ if, and only if, $\lambda$ is a left action of $M$ on $S$. Conclude that also $\ell$ is a left action if, and only if, $\rho$ is a right action.
    \item Prove that $(\lambda,\rho)$ is a biaction of $M$ on $S$ if, and only if, $(\ell,r)$ is a biaction of $M$ on $\cP(S)$.
    \item Prove that, if $(\lambda',\rho')$ is a biaction of $M$ on $S'$ with dual biaction $(\ell',r')$ on $\cP(S')$, then $h \colon S \to S'$ is a biaction morphism if, and only if, $h^{-1} \colon \cP(S') \to \cP(S)$ is a biaction morphism.
  \end{enumerate}
\end{ourexercise}

\begin{ourexercise}\label{ex:nonsurjective-epi-mon}
  Show that the inclusion $(\mathbb{N},+) \into (\mathbb{Z},+)$ is an epimorphism in the category of monoids and homomorphisms, although this is clearly not a surjective function.
\end{ourexercise}

\begin{ourexercise}\label{exe:deltahom}
  Prove that the function $\delta$ defined before Lemma~\ref{lem:autfinite} is a homomorphism of monoids.
\end{ourexercise}

\begin{ourexercise}\label{exe:syntactic}
  Let $M$ be a monoid and let $s\in\cP(M)$.
  \begin{enumerate}
    \item Show that the dual of the complete Boolean residuation ideal generated by $s$ is given by
          \[
            h_s\colon M \twoheadrightarrow M/\equiv_s
          \]
          where $m\equiv_s m'$ if, and only if,
          \[
            \forall x,y\in M \qquad (xmy\in s \ \iff\ xm'y\in s).
          \]
    \item Further show that this is the least recognizer of $s$ in the sense that if $h\colon M \rightarrow N$ is a homomorphism of monoids, then there is a unique map $\tilde{h}\colon \im(h) \twoheadrightarrow M/\equiv_s$ such that $\tilde{h}\circ h'=h_s$, where $h'\colon M\twoheadrightarrow\im(h)$ is the surjective co-restriction of $h$ to its image.
  \end{enumerate}
\end{ourexercise}

\begin{ourexercise}\label{exe:hom-res}
  A consequence of (the proof of) Theorem~\ref{thrm:reglanguage}, which gets a bit lost between the lines because the proof is construed as a cycle around five conditions, is the following fact which we ask you to prove directly.

  Let $h\colon M\to N$ be a homomorphism between monoids and let $u\subseteq N$ and $s\subseteq M$.

  \begin{enumerate}
    \item Show directly from the definitions that
          \[
            s\backslash h^{-1}(u)=h^{-1}(h[s]\backslash u).
          \]
    \item Conclude that, if $h$ is surjective, then for any $p \subseteq N$, we have 
          \[ h^{-1}(p) \backslash h^{-1}(u) = h^{-1}(p \backslash u).\]
    \item Show that, if $N$ is finite, then there exist $m_1,\dots, m_k\in s$ with
          \[
            s\backslash h^{-1}(u)=\bigcap_{i=1}^k m_i^{-1} h^{-1}(u)
          \]
  \end{enumerate}
\end{ourexercise}

\section{Regular languages and free profinite monoids}\label{sec:freeprofmonoid}
In Section~\ref{sec:syntmon} above, we considered regular languages one at a
time, and we saw that the dual of the residuation ideal generated by such a
language is finite -- essentially by virtue of the simple fact that any finite
state machine can only be equipped with finitely many choices of pairs $(I,F)$,
of initial and final states. We placed this result within discrete duality, and,
as a consequence, we were able to see that regular languages over an alphabet
$A$ are precisely those subsets of $\cP(A^*)$ that are recognized by finite
monoids, and that the collection of regular languages is a Boolean residuation ideal of $\cP(A^*)$. 
Now we can give an idea of the shape of the rest of this chapter.

In the current section, we want to consider the set of all regular languages
together as a whole. Since the Boolean residuation ideal of regular languages 
over $A$ is not complete (see Exercise~\ref{exer:RegA-not-complete}), we 
have to switch from the discrete duality to Stone's topological duality. We
will be rewarded by seeing that its dual is a very natural object from
topological algebra, namely, the \emph{free profinite monoid} over $A$.

Next, in Section~\ref{sec:EilReittheory}, we apply the
subalgebra--quotient-space duality to the pair: (regular languages over $A$, the
free profinite monoid over $A$), focusing in on subalgebras of regular languages that are residuation
ideals. This subject culminates in the pairing of so-called pseudovarieties
of regular languages and relatively free profinite monoids that is known in
automata theory via the combination of \emph{Eilenberg's Theorem} and
\emph{Reiterman's Theorem}. We will not give the general theory, but we will
illustrate with the example of so-called piecewise testable languages. In
particular, we will use the duality to give a proof of \emphind{Simon's theorem}, which gives a
decidable characterization of the piecewise testable languages among the regular
ones via profinite monoid equations.

Finally, in Section~\ref{sec:openmult}, we again apply the
subalgebra--quotient-space duality to the pair: (regular languages over $A$, the
free profinite monoid over $A$), but this time we focus on subalgebras of regular languages
that are closed
under concatenation. This points in the direction of
\emph{categorical logic} and \emph{hyperdoctrines} rather than
automata theory and is beyond the scope of this book. We just give a few
elementary observations in the short final section.

\subsection*{Duality for the Boolean residuation algebra of regular languages}

Our main focus in this subsection will be to prove the following theorem; the terms `Boolean residuation algebra' and `free profinite monoid' will be defined below.

\begin{theorem}\label{thrm:dualRegA}
  Let $A$ be a finite set. Then $\RegA$ is a Boolean residuation ideal of $\cP(A^*)$, and the dual space of the Boolean residuation algebra $(\RegA,\backslash,/)$ is the free profinite monoid over $A$.
\end{theorem}

In order to regard $\RegA$ as an algebraic structure in its own right, without necessarily having to refer to its representation inside the algebra $\cP(A^*)$, we now introduce the following notion of Boolean residuation algebra, a Boolean algebra equipped with two implication-type operators, that are linked to each other via a Galois property.\footnote{In this book, for simplicity we only consider the case where the operators are binary, although a more general version exists in the literature, see \cite{Geh16}.}
\begin{definition}\label{def:residuation-alg}
  A \emphind{Boolean residuation algebra} is a tuple $(B, \backslash, /)$, where $B$ is a Boolean algebra, and $\backslash, /:B^2\to B$ are binary operations with the following properties:
  \begin{enumerate}
    \item the operation $\backslash$ preserves finite meets in the second coordinate, that is, $a\backslash \top = \top$ and $a\backslash (b_1\wedge b_2)=(a\backslash b_1)\wedge (a\backslash b_2)$ for all $a, b_1, b_2 \in B$,
    \item the operation $/$ preserves finite meets in the first coordinate, that is, $\top/b=\top$ and $ (a_1\wedge a_2)/b=(a_1/b)\wedge (a_2/b)$ for all $a_1, a_2, b \in B$,
    \item the two operations $\backslash$ and $/$ are linked by the following \emph{Galois property}: for all $a, b, c \in B$,
    $b\leq a\backslash c$ if, and only if, $a\leq c/b$.
  \end{enumerate}
\end{definition}
Note that, under the assumption of (c), conditions (a) and (b) in Definition~\ref{def:residuation-alg} are equivalent.  The Galois property (c) also implies that $\backslash$ and $\slash$ are implication-type operators, in the sense of Definition~\ref{dfn:implicationtype}, that is, the following equations automatically hold:
  \begin{itemize}
    \item $\bot \backslash b = \top$ and $(a_1\vee a_2)\backslash b=(a_1\backslash b)\wedge (a_2\backslash b)$ for all $a_1, a_2, b \in B$, and
    \item $a/\bot = \top$ and $a/(b_1\vee b_2)=(a/b_1)\wedge (a/b_2)$ for all $a, b_1, b_2 \in B$.
  \end{itemize}

Now, for the following corollary to Proposition~\ref{prop:RegBoolResidl} of the previous section, simply note that, since $\RegA$ is closed under $\backslash$ and $\slash$ with arbitrary denominators from $\cP(A^*)$, it is in particular a residuation algebra in its own right.
\begin{corollary}\label{cor:RegABooResAlg}
  Let $A$ be a finite set. Then $\RegA$, the set of all regular languages over $A$, is a Boolean residuation algebra.\end{corollary}
\begin{remark}\label{rem:RegExp}
  Recall that $\cP(A^*)$ is a residuated Boolean algebra, that is, a Boolean algebra equipped with a monoid operation which preserves join that is residuated (see Exercise~\ref{exe:residualscomplete}). In fact, $\RegA$ is not only a residuat\emph{ion} algebra, but even a residuat\emph{ed} Boolean subalgebra of  $\cP(A^*)$. That is, in addition to being closed under the residuation operations, it is also closed under concatenation product of languages. However, $\RegA$ being closed under concatenation is a different phenomenon than what we are looking at here. In fact, for \emph{any} monoid $M$, the set of all the subsets of $M$ recognized by finite monoids is a residuation ideal in $\cP(M)$ and therefore, in particular,  closed under the residual operations, but this Boolean subalgebra is in general \emph{not} closed under the concatenation product of $\cP(M)$ (see Exercise~\ref{exe:RecM}). We will study closure under concatenation product in Section~\ref{sec:openmult}.
\end{remark}

Now we are ready to consider the dual space of the Boolean residuation algebra $\RegA$. We will prove that this dual space is equipped with a monoid operation such that it is, up to isomorphism, the \emphind{free profinite monoid} on $A$.  We now give the necessary definitions.
\begin{definition}\label{def:topmon}
A \emphind{topological monoid} is a tuple $(M, \tau, \cdot)$, where $(M,\tau)$ is a topological space, $(M, \cdot)$ is a monoid, and the function $\cdot \colon M^2 \to M$ is continuous. 
A \emphind{discrete monoid} is a topological monoid with the discrete topology. 
A topological monoid is \emph{profinite}\index{profinite monoid} if it is a projective limit of finite discrete monoids in the category of topological monoids. 
\end{definition}
Unless mentioned otherwise, a finite monoid is always equipped with the discrete topology. 
\nl{$\hat{f}$}{the unique continuous homomorphic extension to the free profinite monoid of a function $f$ from generators to a finite monoid}{}
\begin{definition}\label{dfn:freeprofmonoid}
  Let $A$ be a finite set. A \emphind{free profinite monoid} over $A$ is a profinite monoid $X$ together with a function $\eta \colon A \to X$, such that, for any finite monoid $M$ and function $f \colon A \to M$, there exists a unique continuous homomorphism $\widehat{f} \colon X \to M$ such that $\widehat{f} \circ \eta = f$, as in the following diagram:
  \[
      \begin{tikzcd}[column sep=large,row sep=large,arrows={-Stealth}]
      A \arrow[r,"\eta"]\arrow[dr,"f"']& X\arrow[d,"\widehat{f}",dashed] \\%
      & M
    \end{tikzcd}
  \]
  \end{definition}
In Definition~\ref{dfn:freeprofmonoid}, the continuity of $\widehat{f}$ is with respect to the discrete topology on $M$, and means exactly that $\widehat{f}^{-1}(m)$ is clopen for every $m \in M$.

If $(X,\eta)$ and $(X',\eta')$ are free profinite monoids over $A$, then there is a unique homeomorphism $\phi$ between $X$ and $X'$ such that $\phi \circ \eta = \eta'$, see Proposition~\ref{prop:profinite-char} and Exercise~\ref{exe:freeprofunique}.  The fact that free profinite monoids  exist can be showed 
with a construction in the same vein as the one for profinite ordered sets in Remark~\ref{rem:priestley-as-procompletion}, by taking the limit of a diagram of (discrete) finite monoids in the category of topological monoids. However, our proof that the free profinite monoid exists here will be by exhibiting it as the dual space of the Boolean residuation algebra $(\RegA,\backslash,/)$.

For the rest of this subsection, let $A$ be a finite set and let $X$ denote the dual space of the Boolean algebra $\RegA$. Towards proving that $X$ is the free profinite monoid over $A$, note first that, since $\RegA$ is a Boolean subalgebra of $\cP(A^*)$, its dual space $X$ is a topological quotient of the dual space of $\cP(A^*)$, which is
the Stone-{\v C}ech compactification $\beta A^*$ of the set $A^*$. In particular, we have a natural map $A^*\to X$ given by sending $u\in A^*$ to the point of $X$ corresponding to the ultrafilter consisting of all regular languages that contain the word $u$. Since the singleton $\{u\}$ is regular, this map is injective and the image consists entirely of isolated points (see Exercises~\ref{exe:singletonlanguages}~and~\ref{exe:singletons-dual}). In other words, the map $A^*\to X$ embeds the discrete space $A^*$ in $X$. We will henceforth identify $A^*$ with its image in $X$ and consider $A^*$ as contained in $X$.
\nl{$\bar{h}$}{the unique continuous extension to the free profinite monoid of a homomorphism $h$ from the free monoid to a finite monoid}{}
\begin{remark}\label{rem:reglangdual}
  In what follows, we will use the following consequence of Stone duality, which is the `topological half' of the statement that $X$ is the free profinite monoid over $A$ (see Exercise~\ref{exe:reglangdual} for a proof). Let $h\colon A^*\twoheadrightarrow M$ be a surjective homomorphism to a finite monoid $M$. Viewing the finite set $M$ as a discrete topological space, the function $h$ has a \emph{unique continuous extension} $\bar{h}\colon X\twoheadrightarrow M$, that is, $\bar{h}^{-1}(m)$ is clopen for every $m \in M$ and $\bar{h}(w) = h(w)$ for all $w \in A^*$. Further note that, if $\ell\in\RegA$ is recognized by $h$ via $P\subseteq M$, that is, if $\ell=h^{-1}(P)$, then the corresponding subset $\widehat{\ell}$ of $\widehat{A^*}$ is equal to $\widehat{h^{-1}(P)}=\bar{h}^{-1}(P)$.
\end{remark}
To properly formulate and prove the `monoid half' of the statement that $X$ is the free profinite monoid over $A$, we first need to construct a monoid operation on $X$. This monoid operation comes from the dual of the residuation operations $\backslash$ and $/$ of the Boolean residuation algebra $\RegA$. Recall that in Definition~\ref{dfn:dual-relation-of-implication} we defined a ternary relation dual to any implication-type operator between distributive lattices. Instantiating (\ref{eq:dual-ternary-rel}) in that definition for the specific implication-type operator $\backslash \colon \RegA \times \RegA \to \RegA$, we obtain the ternary relation $R_\backslash$ on $X$ defined by 
\begin{equation} \label{eq:dual-of-backslash}
R_\backslash(x,y,z) \iff \text{for all } \ell, k \in \RegA, \text{if } x \in \widehat{\ell} \text{ and } y \in \widehat{\ell \backslash k}, \text{ then } z \in \widehat{k}.
\end{equation}

\begin{lemma}\label{lem:R-funct1}
  Let $x,y,z\in X$. Then $R_{\backslash}(x,y,z)$ if, and only if, for any finite monoid quotient $h\colon A^*\twoheadrightarrow M$, we have $\bar{h}(z) = \bar{h}(x)\bar{h}(y)$.
\end{lemma}
\begin{proof}
  For the left to right direction, suppose $R_{\backslash}(x,y,z)$ and let $h\colon A^*\twoheadrightarrow M$ be a finite monoid quotient with unique continuous extension $\bar{h} \colon X \to M$. Write $p := \bar{h}(x)$ and $q := \bar{h}(y)$; we need to show that $\bar{h}(z) = pq$. Consider the regular languages $\ell := h^{-1}(p)$ and $k := h^{-1}(pq)$.
  Using Remark~\ref{rem:reglangdual}, since 
  $\widehat{\ell}=\bar{h}^{-1}(p)$, we have $x\in \widehat{\ell}$. Also,
  \[
    \ell\backslash k=h^{-1}(p^{-1}\{pq\}) 
  \]
  by Exercise~\ref{exe:hom-res}, using that $h$ is a monoid homomorphism. Since obviously 
  $q\in p^{-1}\{pq\}$, we get 
  $y\in \bar{h}^{-1}(p^{-1}\{pq\})=\widehat{\ell \backslash k}$, using Remark~\ref{rem:reglangdual} again. 
  Now, because $R_{\backslash}(x,y,z)$ by assumption, we must have $z\in \widehat{k} = \bar{h}^{-1}(pq)$, as required.

  For the right to left direction, suppose $\bar{h}(z) = \bar{h}(x) \bar{h}(y)$ for all finite monoid
  quotients $h\colon A^*\twoheadrightarrow M$. To show $R_{\backslash}(x,y,z)$, suppose that $\ell, k\in\RegA$ are such that 
  $x\in\widehat{\ell}$ and $y\in\widehat{\ell\backslash k}$; we need to prove $z \in \widehat{k}$. By Corollary~\ref{cor:reg-directed}, pick
  a finite monoid quotient $h\colon A^*\twoheadrightarrow M$ which recognizes both languages $\ell$ and $k$, that is, pick $p,q\subseteq M$ such that $\ell=h^{-1}(p)$ and $k=h^{-1}(q)$. By Remark~\ref{rem:reglangdual}, we have $\widehat{\ell}=\bar{h}^{-1}(p)$ and
  $\widehat{k}=\bar{h}^{-1}(q)$. Since $x\in\widehat{\ell}$, we have $\bar{h}(x)\in p$. Also, using Exercise~\ref{exe:hom-res},
  \[
    \ell\backslash k=h^{-1}(p)\backslash h^{-1}(q)=h^{-1}(p\backslash q),
  \]
  so that $\widehat{\ell\backslash k}=\bar{h}^{-1}(p\backslash q)$ and thus $\bar{h}(y)\in p\backslash q$. It now follows from the definition of $p \backslash q$ that $\bar{h}(x)\bar{h}(y)\in q$. Now since $\bar{h}(z) = \bar{h}(x) \bar{h}(y)$ by hypothesis, we get $z\in\bar{h}^{-1}(q) = \widehat{k}$, as required.
\end{proof}

\begin{lemma}\label{lem:R-funct2}
Let $x, y \in X$. The set of regular languages
\begin{align*} 
\mu_{x,y} := \{\ell \in \RegA \ \mid \ &\text{there is a finite monoid quotient } h \colon A^* \twoheadrightarrow M \text{ which } \\
&\text{  recognizes $\ell$ and such that } h^{-1}(\bar{h}(x) \bar{h}(y)) \subseteq \ell\}
\end{align*}
is an ultrafilter of $\RegA$. Moreover, the point $z \in X$ for which $F_z = \mu_{x,y}$ is the unique point in the set $R_{\backslash}[x,y,\_]$.
\end{lemma}
\begin{proof}
  Note first that, for any $\ell \in \RegA$, if $h \colon A^* \twoheadrightarrow M$ is a finite monoid quotient which recognizes $\ell$, then, for any $m \in M$, exactly one of the following two properties holds:
  \[ h^{-1}(m) \subseteq \ell \quad \text{ or } \quad h^{-1}(m) \subseteq A^* \setminus \ell \ .\]
  Indeed, at most one of the two can hold because $h^{-1}(m)$ is non-empty, since $h$ is surjective. Also, since $h$ recognizes $\ell$, we have $\ell = h^{-1}(u)$ for some $u \subseteq M$, so that at least one of the two must hold: the first in case $m \in u$ and the second in case $m \not\in u$.

  It follows that, for every regular language $\ell$, since there exists some finite monoid quotient $h \colon A^* \onto M$ recognizing $\ell$, we have that either $\ell$ or $A^* \setminus \ell$ is in $\mu_{x,y}$, applying the above argument to $m := \bar{h}(x) \bar{h}(y)$. It is also clear that $\emptyset \not\in \mu_{x,y}$. It remains to prove that $\mu_{x,y}$ is a filter.
%
  To this end, let $\ell_1, \ell_2$ be regular languages. If $\ell_1$ and $\ell_2$ are both in $\mu_{x,y}$, then, for $i = 1, 2$, pick finite monoid quotients $h_i\colon A^*\twoheadrightarrow M_i$, recognizing $\ell_i$ and such that $h_i^{-1}(\overline{h_i}(x)\overline{h_i}(y)) \subseteq \ell_i$. 
  Define the homomorphism $h\colon A^*\twoheadrightarrow M$ as the co-restriction of the homomorphism $h_1\times h_2 \colon A^* \to M_1 \times M_2$ to its image. Then $h$ is a finite monoid quotient of $A^*$ that recognizes $\ell_1 \cap \ell_2$, as we saw in the proof of Proposition~\ref{prop:RegBoolResidl}. Moreover, for $i = 1, 2$, we have
  \[h^{-1}(\overline{h}(x)\overline{h}(y))\subseteq h_i^{-1}(\overline{h_i}(x)\overline{h_i}(y))\ . \] 
  Combining this with the assumption that $h_i^{-1}(\overline{h_i}(x)\overline{h_i}(y)) \subseteq \ell_i$, we see that
  $h^{-1}(\overline{h}(x)\overline{h}(y)) \subseteq \ell_1 \cap \ell_2$, so that $\ell_1 \cap \ell_2 \in \mu_{x,y}$. To see that $\mu_{x,y}$ is an up-set, suppose that $\ell_1 \in \mu_{x,y}$ and that $\ell_1 \subseteq \ell_2$. Pick a finite monoid quotient $h_1 \colon A^* \onto M_1$ recognizing $\ell_1$ such that $h^{-1}(\overline{h_1}(x)\overline{h_1}(y)) \subseteq \ell_1$ and pick some finite monoid quotient $h_2 \colon A^* \onto M_2$ recognizing $\ell_2$. Defining $h \colon A^* \onto M$ in the same way as before, this morphism still recognizes $\ell_2$, and  $h^{-1}(\overline{h}(x)\overline{h}(y))$ intersects non-trivially with $\ell_1$, so it also intersects non-trivially with $\ell_2$. 

  For the moreover statement, denote by $z$ the unique point of $X$ such that $F_z = \mu_{x,y}$. We use Lemma~\ref{lem:R-funct1} to show that $R_{\backslash}(x,y,\_) = \{z\}$. To see that $R_{\backslash}(x,y,z)$, let $h \colon A^* \onto M$ be a finite monoid quotient. Then the language $\ell := h^{-1}(\bar{h}(x)\bar{h}(y))$ is regular and, clearly, $\ell \in \mu_{x,y} = F_z$, which means that $z \in \widehat{\ell} = \bar{h}^{-1}(\bar{h}(x)\bar{h}(y))$. We conclude that $R_\backslash(x,y,z)$ by Lemma~\ref{lem:R-funct1}. Conversely, if $R_{\backslash}(x,y,z')$ for any point $z' \in X$, then let $\ell \in \mu_{x,y}$ be arbitrary. Pick a finite monoid quotient $h \colon A^* \onto M$ recognizing $\ell$ with $h^{-1}(\bar{h}(x)\bar{h}(y)) \subseteq \ell$. Write $\ell' := h^{-1}(\bar{h}(x)\bar{h}(y))$.  By Lemma~\ref{lem:R-funct1}, $\bar{h}(z') = \bar{h}(x)\bar{h}(y)$, so $z' \in \widehat{\ell'} \subseteq \widehat{\ell}$. Thus, $\mu_{x,y} \subseteq F_{z'}$, and the two ultrafilters must be equal, so $z' = z$.
\end{proof}

We summarize Lemma~\ref{lem:R-funct1} and Lemma~\ref{lem:R-funct2} in the following corollary.

\begin{corollary}\label{cor:R-funct}
  The ternary relation $R_{\backslash}$ dual to $\backslash$ is a total binary operation $\star\colon X\times X \to X$, and the ultrafilter corresponding to $x\star y$ is $\mu_{x,y}$. 
 Moreover, for any $x, y \in X$, $x \star y$ is the unique element of $X$ such that, for every finite monoid quotient $h\colon A^*\twoheadrightarrow M$, we have $\bar{h}(x\star y)=\bar{h}(x)\bar{h}(y)$.
\end{corollary}
Of course, we also have the ternary relation $R_{/}$ dual to the operator $/$, which is of `reverse' implication-type, in the sense that it is monotone in the first and antitone in the second coordinate. Concretely, we may define, analogously to (\ref{eq:dual-of-backslash}),
\begin{equation}\label{eq:dual-of-slash}
R_{/}(x,y,z)  \iff \text{for all } \ell, k \in \RegA, \text{ if } x \in \widehat{\ell} \text{ and } y \in \widehat{k / \ell}, \text{ then } z \in \widehat{k}.
\end{equation}
However, this ternary relation $R_{/}$ does not add any additional structure to the dual space, as it is closely related to $R_{\backslash}$ above. Indeed, the Galois property in the definition of the Boolean residuation algebra $\RegA$ (see Definition~\ref{def:residuation-alg}) expresses the fact that the relations $R_{/}$ and $R_{\backslash}$ are related via
\[ R_{\backslash}(x,y,z) \iff R_{/}(y,x,z) \] 
for any $x, y, z \in X$, so that $R_{/}$ essentially contains the same information as $R_{\backslash}$. This may be proved in a similar way to the case of relations dual to an adjunction considered in Exercise~\ref{exe:converse-adjunction} (see Exercise~\ref{exe:slash-backslash-converse}).

Finally, we show that the operation $\star$ makes $X$ into a topological monoid and that it satisfies the universal property of the free profinite monoid over $A$. Recall that we view $A^*$ as a discrete subspace of $X$, and we denote the empty word by $\epsilon$.

\begin{lemma}\label{lem:Xfreeprof}
The triple $(X,\star,\epsilon)$ is the free profinite monoid over $A$.
\end{lemma}

\begin{proof}
  First we show that $(X,\star)$ is a monoid with neutral element the empty word $\epsilon \in A^* \subseteq X$. Let $x,y,z\in X$. In order to prove $(x \star y) \star z = x \star (y \star z)$, we show that the points have the same clopen neighborhoods. Let $\ell\in\RegA$ be arbitrary, and suppose that $(x\star y)\star z\in\widehat{\ell}$; we will show that $x \star (y \star z) \in \widehat{\ell}$. Pick $h\colon A^*\twoheadrightarrow M$ a finite monoid quotient which recognizes $\ell$ via $P$, so that $\widehat{\ell}={\bar{h}}^{-1}(P)$. Repeatedly using the defining property of $\star$ in Corollary~\ref{cor:R-funct}, we have 
  \[ \bar{h}((x \star y) \star z) = (\bar{h}(x)\bar{h}(y))\bar{h}(z), \text{ and } \bar{h}(x \star (y \star z)) = \bar{h}(x)(\bar{h}(y)\bar{h}(z)).\]
  Now, since $M$ is a monoid, these two elements of $M$ equal, and thus in particular we also have $x \star (y \star z) \in \bar{h}^{-1}(P) = \widehat{\ell}$.
  Similar reasoning 
  shows that $\epsilon$ is a neutral element for $\star$ (see Exercise~\ref{exe:epsilon-freeprof}).

  Next, we need to show that the function $\star$ is continuous. Suppose $x,y\in X$, $\ell\in\RegA$ and $x\star y\in\widehat{\ell}$. Again, let $h\colon A^*\twoheadrightarrow M$ be a finite quotient which recognizes $\ell$ via $P$, then $\widehat{\ell}=f^{-1}(P)$ where $f\colon X\twoheadrightarrow M$ is the continuous extension of $h$. Let $\ell_1=h^{-1}(\bar{h}(x))$ and $\ell_2=h^{-1}(\bar{h}(y))$. Then $\ell_1,\ell_2\in\RegA$ and, for any $x'\in\widehat{\ell}_1$ and $y'\in\widehat{\ell}_2$, we have
  \[\bar{h}(x'\star y')=\bar{h}(x') \bar{h}(y')=\bar{h}(x)\bar{h}(y)=\bar{h}(x\star y)\in P.\]
  That is, $x'\star y'\in\widehat{\ell}$ as required.

  Finally, we show that $(X,\star,\epsilon)$ has the universal property of the free profinite monoid over $A$ with respect to the function $\eta$ which sends each letter $a \in A$ to the one-letter word $a$ in $X$. To this end, let $M$ be a finite monoid and let $f\colon A\to M$ be any function. Then, since $A^*$ is the free monoid over $A$, we obtain a unique homomorphic extension $f^*\colon A^* \to M$. 
  Since $\im(f^*)$ is closed in the discrete space $M$, any continuous extension of $f^*$ will map onto $\im(f^*)$, so it suffices to consider the surjective co-restriction $h\colon A^*\onto\im(f^*)$ of $f^*$ to its image. As we have been using extensively throughout, see Remark~\ref{rem:reglangdual}, $h$ extends uniquely to a continuous map $\bar{h}\colon X\onto\im(f^*)$, and, by Corollary~\ref{cor:R-funct}, $\bar{h}$ is also a monoid morphism with respect to $\star$. In a diagram,
  
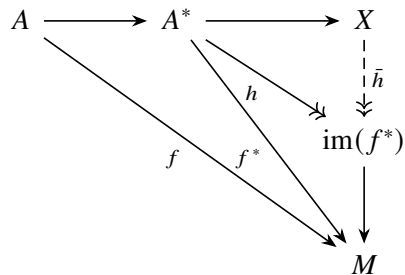
\[
    \begin{tikzcd}[column sep=large,row sep=large,arrows={-Stealth}]
      A\arrow[r]\arrow[ddrr,"f"']&A^* \arrow[r]\arrow[ddr,"f^*"']\arrow[dr,->>,"h"']& X\arrow[d,"\bar{h}",->>,dashed] \\%
      &&\im(f^*)\arrow[d]\\
      && M
    \end{tikzcd}
  \]
  Each of the inner triangles of the diagram commutes and thus so does the outer triangular diagram. Finally, the outer vertical arrow $X \to M$ is the unique such continuous homomorphism, since $A$ generates $A^*$, $A^*$ is dense in $X$, and $\bar{h}$ is both a monoid homomorphism and continuous (see Exercise~\ref{exe:uniqueness-freeprof}).
\end{proof}

This completes the proof of Theorem~\ref{thrm:dualRegA}. We chose to give a very
concrete proof of this theorem as this is maybe more hands-on and tangible for
the novice and because such a proof is not readily available in the literature.
Usually, in research papers, a more conceptual proof is given: $\RegA$ is
clearly the filtered colimit of its finite residuation ideals, so, by duality,
its dual space is the co-filtered, or projective, limit of the duals of
these finite residuation ideals, which, by
Theorem~\ref{thrm:reglanguage}, are the finite monoid quotients of $A^*$. It
remains to check that the projective limit  of the finite monoids in Boolean
spaces equipped with a profinite monoid structure is indeed the multiplication
operation on the profinite monoid. This yields a quicker, simpler, but maybe
somewhat magical proof of Theorem~\ref{thrm:dualRegA}. See for example
\cite[Theorem~4.4]{Geh16} for further details of this alternate proof.

Now that we know that a free profinite monoid over $A$ exists and is unique up to isomorphism, we can speak of  \emph{the} free profinite monoid over $A$, and we will denote it by $\wh{A^*}$.
\nl{$\hat{A^*}$}{the free profinite monoid over a finite set $A$}{}

We finish this subsection by briefly discussing the more general notion of
profinite monoid, of which the free profinite monoid will be our main example.
Since these results are not central to the rest of the chapter, we leave the
proofs of the statements as extended exercises.

\nl{$\ClpCon(M)$}{the set of clopen congruences on a topological monoid $M$}{}
Let us first give a more concrete characterization of profinite monoids, that is often
convenient to work with. When $M$ is a topological monoid, we will denote by
$\ClpCon(M)$ the set of \emph{clopen congruences} on $M$, that is,
congruences $\theta$ such that $\theta$ is clopen as a subset of $M \times M$. Recall from Exercise~\ref{exe:finitesubs-of-BA} that for a congruence $\theta$
to be clopen, it is equivalent to say that 
$M/{\theta}$ is
finite and $M \onto M/{\theta}$ is continuous for the discrete topology on
$M/{\theta}$. Note further that, for a congruence $\theta$ such that $M/{\theta}$ is finite, the continuity condition
is equivalent to requiring that each equivalence class $[m]_{\theta}$ is clopen.
\index{congruence!clopen}
\begin{proposition}\label{prop:profinite-char}
A topological monoid $M$ is profinite if, and only if, the space underlying $M$
is compact, and, for every $x, y \in M$, if $x \neq y$, then there exists a
continuous homomorphism $f \colon M \to F$, with $F$ a finite discrete monoid,
such that $f(x) \neq f(y)$. 

If $M$ is profinite, then $M$ is isomorphic to the closed submonoid of
$\prod_{\theta \in \ClpCon(M)} M/{\theta}$ consisting of the tuples
$(x_{\theta})_{\theta \in \ClpCon(M)}$ such that, whenever $\theta \subseteq
\theta'$ and $x \in M$ is such that $[x]_{\theta'} = x_{\theta'}$, then
$[x]_{\theta} = x_{\theta}$. 
\end{proposition}
\begin{proof}
	See Exercise~\ref{exe:profinite-char}.
\end{proof}
\begin{remark}\label{rem:hunter}
Recall that Boolean topological spaces are the same thing as profinite
\emph{sets}; see Example~\ref{exa:priestley-as-profinite posets} in Chapter~\ref{ch:categories}.
It is in particular easy to see from Proposition~\ref{prop:profinite-char} that
the topology of a profinite monoid must be Boolean, since it is a closed
submonoid of a product of finite spaces. It is less immediate that the converse
holds true, in the following sense: if $\cdot$ is a continuous associative multiplication on a
Boolean space $M$ with neutral element $1$, then $(M,\cdot,1)$ is a profinite monoid. Towards showing the
condition in the first part of Proposition~\ref{prop:profinite-char}, if $x \neq
y$, there is, by zero-dimensionality of the space, a clopen set separating $x$ from $y$, but it
remains to show that this clopen set can be realized as a union of equivalence classes for
some $\theta \in \ClpCon(M)$. This is true in the case of monoids, where it is known
as `Hunter's lemma', but it does not
generally hold true for other algebraic structures based on Boolean spaces; see
\cite{AGK21} for a recent analysis of the question of what makes a topological
algebra based on a Boolean space profinite.
\end{remark}
The characterization of Proposition~\ref{prop:profinite-char} in particular
allows us to define an element of a profinite monoid by giving its value modulo
$\theta$, for each $\theta \in \ClpCon(M)$. The reader will note that this is very
similar to the way we characterized the multiplication $\star$ on $X$ in
Corollary~\ref{cor:R-funct} above. We now give one more instance of such a
definition `via the finite quotients', which we will use in the next section.
An element $e$ of a monoid is called \emphind{idempotent} if $e \cdot e = e$. 
\begin{proposition}\label{prop:omega}
Let $M$ be a profinite monoid. For any $x \in M$, there is a unique idempotent
element in the closure of the set $\{x^n \ \mid \ n \geq 1\}$.
\end{proposition}
\begin{proof}
First note that, for a finite monoid $M$, there is always a unique idempotent
element in the set $\{x^n \ \mid \ n \geq 1\}$, for every $x \in M$. Indeed, by
the pigeon-hole principle, pick $n, p \geq 1$ such that $x^n = x^{n+p}$. 
Let $r$ be the remainder of the division of $-n$ by $p$, that is, choose the unique $0 \leq r < p$ such that $n+r$ is a multiple of $p$. Then $x^{n+r}$ is
idempotent, as shown by the computation:
\[ x^{n+r} x^{n+r} =
x^{n+r+qp} = x^{n+qp} x^r = x^{n+r}.\]
If, for any $m \geq 1$, $x^m$ is also idempotent, then $x^{m} = (x^m)^{n+r} = (x^{n+r})^m = x^{n+r}$,
so this idempotent is unique. This establishes the proposition for the special
case of finite (discrete) monoids.

Now, if $M$ is a profinite monoid and $x \in M$, let us write $P := \{x^n \ \mid
\ n \geq 1\}$. If $y$ is an idempotent element in the closure of $P$, then for
any continuous homomorphism $f \colon M \to F$, with $F$ a finite discrete
monoid, $f(y)$ must be idempotent, and since $y$ is in the closure of $P$, it
must be in the closed set $f^{-1}(f[P])$, so that $f(y)$ is
a power of $f(x)$. By Proposition~\ref{prop:profinite-char} and the fact that
$f(x)$ has a unique idempotent power for every $f$, there can be at
most one such element $y$. Moreover, such an element exists, because if $\theta
\subseteq \theta'$ and $[x]_{\theta'}^m$ is idempotent, then $[x]_{\theta}^m$ is
idempotent, as well, so, using the second part of
Proposition~\ref{prop:profinite-char} we have the element 
$x^\omega$ of $M$ defined by the condition that, for each $\theta \in \ClpCon(M)$, 
$[x^\omega]_{\theta}$ is the unique idempotent power of $[x]_{\theta}$ in
$M/{\theta}$. 
\end{proof}
 \begin{definition}\label{def:omega}
	Let $M$ be a profinite monoid. For $x \in M$, we denote by
	$x^{\omega}$ the unique idempotent element in the closure of $\{x^n
	\ \mid \ n \geq 1\}$.
\end{definition}
Similarly, there exists, for any $x \in M$, an element $x^{\omega-1}$ in the closure of $\{x^n \ \mid \ n \geq 1 \}$ which is
uniquely determined by the condition that $x^{\omega-1} \cdot x = x\cdot x^{\omega-1} = x^{\omega}$ (see Exercise~\ref{exe:omegaminus1}).

\subsection*{A general duality for Boolean residuation algebras}
Above, we have seen that Stone duality on objects  restricts to a correspondence between finitely generated free profinite monoids and the Boolean residuation algebras of the form $\RegA$, for $A$ a finite alphabet. In order to show that this is a more general phenomenon, we will show here that
\begin{enumerate}
  \item \label{itm:bin-top-duals} all binary topological algebras whose
	  underlying space is Boolean are dual spaces of Boolean residuation algebras;
  \item \label{itm:duals-jpp} the residuation algebras dual to binary topological algebras are precisely those which preserve joins at primes.
\end{enumerate}
Item (\ref{itm:duals-jpp}) is particularly interesting as it makes a link to Domain Theory in Logical Form as treated in Sections~\ref{sec:funcspace}~and~\ref{sec:DTLF}, where the concept of preserving joins at primes is also central, see Corollary~\ref{cor:spectral-funcspace} and Theorem~\ref{thm:Bif-CCC}.

By a \emphind{binary topological algebra}, we mean a pair $(X,f)$, where $X$ is
a topological space and $f \colon X^2 \to X$ is a continuous binary operation on
$X$, with no equational axioms assumed on $f$. A set equipped with a binary
operation satisfying no axioms has also been called \emphind{magma} in the
algebraic and categorical literature; from that perspective, we here consider
magmas internal to the category $\TopCat$, and in particular in the full
subcategory $\BoolSp$. The fact that we restrict ourselves to \emph{binary}
operations is just for simplicity of notation and a more general theorem for
arbitrary arity can be proved, see \cite[Section~3]{Geh16}. 
We will use the term `binary topological algebra on a Boolean space' for one
whose underlying space is Boolean; these structures could also be called
`Boolean-topological magmas'.

Recall from Section~\ref{sec:generalopduality} that implication-type operators on a lattice are dual to \emph{compatible} ternary relations on its dual space. Thus, to show (\ref{itm:bin-top-duals}), we start by showing that a binary operation on a Boolean space is continuous if, and only if, its graph is a compatible ternary relation in the sense of Definition~\ref{dfn:compatible-implication}. By the \emph{graph}\index{graph of a binary operation} of a binary operation $f \colon X^2 \to X$, we here mean the ternary relation 
\[ \{ (x,y,f(x,y)) \ \mid \ x, y \in X \}. \]
Since we are in the special case of Boolean spaces, we will here specialize the definition of compatibility for a ternary relation to that setting. Note that a relation $R \subseteq X^3$ is then
\emph{compatible}\index{compatible!relation of implication type on a Boolean
space}\index{implication type!compatible relation of} if for any clopen subsets $U, V$ of $X$, $R[U, \_, V]$ is clopen, and for any $y \in X$, the set $R[\_,y,\_]$ is closed. 

\begin{proposition}\label{prop:cont-compatible}
  Let $X$ be a Boolean space, $f\colon X^2\to X$ a binary operation on $X$, and $R\subseteq X^3$ its graph. Then $f$ is continuous if, and only if, its graph is a compatible relation on $X$.
\end{proposition}

\begin{proof}
Throughout this proof, for any $y \in X$, we denote by $i_y \colon X \to X^2$ the injective function $x \mapsto (x,y)$, which is continuous by Exercise~\ref{exer:easy-product}, and we write $f_y := f \circ i_y$, that is, $f_y \colon X \to X$ is the function defined by $f_y(x) = f(x,y)$ for every $x \in X$.
  
  Suppose first $f$ is continuous. Then $f_y$ is also continuous, and for any $y \in X$, the set $R[\_,y,\_]$ is the graph of the function $f_y \colon X \to X$. This graph is closed because $X$ is Hausdorff (see Exercise~\ref{exer:Hausdorff}).
  Let $U,V\subseteq X$ be clopen. First notice that
  \[
    R[U,\_,V]=\pi_2[(U\times X)\cap f^{-1}(V)],
  \]
  where we denote by $\pi_2\colon X^2\to X$ the projection on the second coordinate. Further, since $X$ is compact, the projection $\pi_2$ is a closed mapping by Proposition~\ref{prop:proj-along-comp}, and thus $R[U,\_,V]$ is closed. Finally, by Exercise~\ref{exer:easy-product}, $\pi_2$ is always an open mapping, and thus $R[U,\_,V]$ is clopen as required.

  For the converse, suppose the graph $R$ is compatible. We show first that $f_y$ is continuous for every $y \in X$. For any $S\subseteq X$, we have
  \[
    f_y^{-1}(S)=\pi_1[R[\_,y,\_]\cap(X\times S)],
  \]
  that is, the inverse image of $S$ under $f_y$ is the same as the direct image under $\pi$ of the set $R[\_,y,\_] \cap (X \times S)$. Again, since $\pi_1$ is a closed mapping and $R[\_,y,\_]$ is closed, this shows that $f_y^{-1}(S)$ is closed whenever $S$ is closed. 
  Now let $V\subseteq X$ be clopen; we show that $f^{-1}(V)$ is open. Let $(x,y)\in f^{-1}(V)$ be arbitrary. Then $x\in f_y^{-1}(V)$, which is open, so pick a clopen $U \subseteq X$ satisfying $x\in U\subseteq f_y^{-1}(V)$. Consider the set 
  \[U \backslash V := R[U,\_,V^c]^c = \{ w \in X \ \mid \ \text{for all } u \in U, f(u,w) \in V\},\] 
  which is clopen because $R$ is compatible. Note that, since $U \subseteq f_y^{-1}(V)$, we have $y \in U \backslash V$, and clearly $U \times (U \backslash V) \subseteq f^{-1}(V)$. Thus, $U \times (U \backslash V)$ is an open neighborhood around $(x,y)$ contained in $f^{-1}(V)$, and we conclude that $f^{-1}(V)$ is open, as required.
\end{proof}

Notice that by symmetry of assumptions, exactly the same proof with the role of the first two coordinates of the relation $R$ switched would work just as well. Thus, if $f$ is a continuous operation on $X$, then its dual Boolean algebra is a residuation algebra. We get the following corollary.

\begin{corollary}\label{cor:booresalg} Every binary topological algebra on a Boolean space is the dual space of a Boolean residuation algebra.
\end{corollary}
For later use in this section, we recall one direction of the object duality of Corollary~\ref{cor:booresalg} in more concrete terms. 
\begin{notation}
  Throughout the rest of this section, in order to simplify notation, when $B$ is the Boolean algebra of clopens of a Boolean space $X$, we will tacitly identify ultrafilters of $B$ with points of $X$; that is, if $U \in B$ and $x \in X$ then $U$ belongs to the ultrafilter corresponding to $x$ if, and only if, $x \in U$. 
\end{notation} 
Let $(X,f)$ be a binary topological algebra on a Boolean space, let $R$ be the graph of $f$, and let $B$ be the Boolean algebra of clopens dual to $X$. 
The algebra dual to $(X,R)$ is then the Boolean residuation algebra $(B,\backslash,/)$ where, if $U,V\in B$, then
  \[
    U\backslash V=R[U,\_,V^c]^c=\{y\in X\mid \text{ for all } u \in U, \ f(u,y)\in V\}.
  \]
  Here, we have simply applied the definition of implication-type operator associated with a ternary relation (\ref{eq:imp-of-relation}) in the specific case where the relation is the graph of a continuous binary function. The definition of the other residual $/$ is similar but with the role of the first two coordinates switched. When $U = \{x\}$, we will use $x^{-1}V$ as a notation for $\{x\} \backslash V$.

We now proceed to characterize exactly \emph{which} Boolean residuation algebras can occur as the dual of a binary topological algebra on a Boolean space. This is where the notion of preserving joins at primes re-appears.

Let $(B, \backslash, \slash)$ be a Boolean residuation algebra, as defined in Definition~\ref{def:residuation-alg}. Also recall from Definition~\ref{def:jpp} what it means for an implication-type operator to preserve joins at primes. That definition was given for an implication-type operator $\to \colon L \times M \to K$, which is antitone in the first and monotone in the second coordinate, and thus directly specializes to the residual $\backslash$ on a Boolean residuation algebra: the operation $\backslash$ \emph{preserves joins at primes} if, for every ultrafilter $F$ of $B$, $a \in F$, and finite subset $G$ of $B$, there exists $a' \in F$ such that 
\[ a \backslash \Big(\bigvee G \Big) \leq \bigvee_{g \in G} (a' \backslash g).\] 
The other residual, $/$, is of `reverse' implication-type, in that it is monotone in the first and antitone in the second coordinate. We say that the operation $/$ \emph{preserves joins at primes} if for every ultrafilter $F$ of $B$, $a \in F$, and finite subset $G$ of $B$, there exists $a' \in F$ such that
  \[ \Big(\bigvee G\Big)/a \leq \bigvee_{g \in G} (g/a').\]
\begin{definition}\label{dfn:jpp-resalg}
  Let $(B,\backslash,\slash)$ be a Boolean residuation algebra and let $X$ be its dual. We say that the algebra $B$ \emph{preserves joins at primes}\index{preserve join at primes!for a residuation algebra} provided that both operations $\backslash$ and $/$ preserve joins at primes.
\end{definition}
It will in fact follow from the proof of Theorem~\ref{thm:bintopalg-dual-to-jpp} below that the one residual $/$ preserves joins at primes if, and only if, the other residual $\backslash$ preserves joins at primes, so it would suffice in Definition~\ref{dfn:jpp-resalg} to assume only one of the two.

Before proving Theorem~\ref{thm:bintopalg-dual-to-jpp}, we isolate one lemma that contains an important step.
\begin{lemma}\label{lem:esakia-for-backslash}
  Let $(X,f)$ be a binary topological algebra on a Boolean space $X$, and let $(B, \backslash, \slash)$ be the dual Boolean residuation algebra of clopen subsets of $X$. For any $x \in X$ and $V \in B$, there exists $U \in B$ such that $x \in U$ and $x^{-1} V = U \backslash V$.
\end{lemma}
\begin{proof}
  Let $V \in B$ and $x \in X$. We first show that  the collection 
  \[ \mathcal{C} := \{ K \backslash V \ \mid \ K \in B \text{ and } x \in K \} \] 
  is a cover of the set $x^{-1}V$. Indeed, if $y \in x^{-1}V$, then $f(x,y) \in V$, so by continuity of the function $f_y : x \mapsto f(x,y)$, since $x \in f_y^{-1}(V)$, there is a clopen set $K$ around $x$ such that $K \subseteq f_y^{-1}(V)$, which means that $y \in K \backslash V$. Now, since all sets in $\mathcal{C}$ are clopen, and $x^{-1}V$ is closed (even clopen) and hence compact, pick a finite set $\mathcal{K} \subseteq B$ such that $x \in K$ for every $K \in \mathcal{K}$, and $x^{-1} V \subseteq \bigcup_{K \in \mathcal{K}} K \backslash V$. Define $U := \bigcap \mathcal{K}$. Then clearly $x \in U$, so $U \backslash V \subseteq x^{-1} V$, and, for every $K \in \mathcal{K}$, we have that $U \subseteq K$, so $K \backslash V \subseteq U \backslash V$. Thus, 
  \[ x^{-1} V \subseteq \bigcup_{K \in \mathcal{K}} K \backslash V \subseteq U \backslash V \subseteq x^{-1} V.\] 
  In particular, $x^{-1} V = U \backslash V$, as required.
\end{proof}
\begin{theorem}\label{thm:bintopalg-dual-to-jpp}
  The dual algebras of binary topological algebras on Boolean spaces are precisely the Boolean residuation algebras preserving joins at primes.
\end{theorem}

\begin{proof}
  First suppose that $(X,f)$ is a binary topological algebra and $(B,\backslash,/)$ is its dual Boolean residuation algebra. We need to show that $(B,\backslash,/)$ preserves joins at primes. We only show this for $\backslash$, the argument for $/$ is symmetric. By induction, it suffices to treat just the cases where $G$ is empty or $G$ contains two elements; when $G$ contains one element there is nothing to do. 
  If $G$ is empty, then $\bigvee G = \emptyset$. Let $U \in B$ and suppose that $u \in U$. We have 
  \[
    U \backslash {\emptyset} \subseteq \{y\in X\mid f(u,y)\in \emptyset\}= \emptyset,
  \]
  so we may take $U' := U$.

  Now suppose that $G = \{V_1, V_2\}$ and let $U \in B$ and $u \in U$. Note first that 
  \[ U \backslash (V_1 \cup V_2) \subseteq u^{-1}(V_1 \cup V_2) = u^{-1}V_1 \cup u^{-1}V_2,\] 
  using that $V \mapsto u^{-1} V$ is a Boolean algebra homomorphism.
  By Lemma~\ref{lem:esakia-for-backslash}, pick $U_1, U_2 \in B$ such that $u \in U_i$ and $u^{-1}V_i = U_i \backslash V_i$ for $i = 1,2$. Now define $U' := U_1 \cap U_2$, which still contains $u$. Plugging this into the previous inclusion, we get
  \[ U \backslash (V_1 \cup V_2) \subseteq (U_1 \backslash V_1) \cup (U_2 \backslash V_2) \subseteq (U' \backslash V_1) \cup (U' \backslash V_2),\]
  using for the second inclusion that $\backslash$ is antitone in its first coordinate.
 
  For the converse, suppose that $(B,\backslash,/)$ is a Boolean residuation algebra such that $\backslash$ preserves joins at primes, and let $(X,R)$ be its dual. We will show that $R$ is the graph of a binary function $f \colon X^2 \to X$, that is, that $R[x,y,\_]$ is a singleton for all $x, y \in X$. Let $x, y \in X$ be arbitrary. Note that it follows from the definition of $R$ as the ternary relation dual to the implication-type operator $\backslash$ (Definition~\ref{dfn:dual-relation-of-implication}) that,  for any $z \in X$, we have $R(x,y,z)$ if, and only if, $z \in \bigcap F$, where
  \[ F := \{ V \in B \ \mid \ \text{ there exists } U \in B \text{ such that } x \in U \text{ and } y \in U \backslash V \}.\]
  We show that $F$ is an ultrafilter, from which it follows that $\bigcap F$ is a singleton, since $X$ is a Boolean space (also see Exercise~\ref{exe:intersection-ultrafilter}). Using that, for any $U \in B$, the function $V \mapsto U \backslash V$ preserves finite meets, one may prove that $F$ is a filter because it is a directed union of filters (see Exercise~\ref{exer:filter-from-backslash}). To get that $F$ is an ultrafilter, we use the assumption that $\backslash$ preserves joins at primes. First, we show that $F$ is proper. Indeed, for any $U \in B$ such that $x \in U$, the definition of preserving joins at primes applied in the case $G = \emptyset$ gives that $U \backslash \emptyset = \emptyset$, so in particular $y \not\in U \backslash \emptyset$, so that $\emptyset$ is not in $F$. Now let $V \in B$ be arbitrary. Applying the assumption that $\backslash$ preserves joins at primes to the element $\top \in B$ and $G = \{V, V^c\}$, since $x \in \top$ and $\top \backslash (V \cup V^c) = \top \backslash \top = \top$, we can pick $U \in B$ with $x \in U$ and $(U \backslash V) \cup (U \backslash V^c) = \top$.  Then either $y \in U \backslash V$ or $y \in U \backslash V^c$, so one of $V$ and $V^c$ is in $F$, as required.

  Finally, let $f \colon X^2 \to X$ be the function defined, for $x,y \in X$, by taking $f(x,y)$ the unique element in $R[x,y,\_]$. Since the relation $R$ is compatible, it follows from Proposition~\ref{prop:cont-compatible} that $f$ is continuous. 
\end{proof}

\ourexercises
\begin{ourexercise}
  Prove the remarks directly following Definition~\ref{def:residuation-alg}.
\end{ourexercise}
\begin{ourexercise}\label{exe:singletons-dual}
  Let $W$ be a set and let $B$ be a Boolean subalgebra of $\cP(W)$. Denote the dual space of $B$ by $X$, which is a topological quotient of $\beta X$. Consider the composition $i \colon W \to \beta W \onto X$. Assume further that $B$ contains $\{w\}$ for every $w \in W$. Show that $i$ is injective and its image is a dense subspace of $X$, and that the subspace topology on $i[W]$ is discrete.
\end{ourexercise}
\begin{ourexercise}\label{exe:slash-backslash-converse}
  Let $(B, \backslash, \slash)$ be a Boolean residuation algebra with dual space $X$. Consider the relations $R_{\backslash}$ and $R_{\slash}$ defined by (\ref{eq:dual-of-backslash}) and (\ref{eq:dual-of-slash}), respectively. Prove that, for any $x,y,z \in X$,
  \[ R_{\backslash}(x,y,z) \iff R_{\slash}(y,x,z).\]
\end{ourexercise}
\begin{ourexercise}\label{exe:epsilon-freeprof}
  Prove that $\epsilon$ is a neutral element for the operation $\star$ on $X$, following a similar proof to the first paragraph of the proof of Lemma~\ref{lem:Xfreeprof}.
\end{ourexercise}
\begin{ourexercise}\label{exe:uniqueness-freeprof}
  This exercise gives some more details to show the uniqueness of the extension constructed in the proof of Lemma~\ref{lem:Xfreeprof}.
  \begin{enumerate}
    \item Let $X$ and $Y$ be Hausdorff topological spaces and $D$ a dense subspace of $X$. Prove that if $h, h' \colon X \to Y$ are continuous and $h|_D = h'_D$, then $h = h'$.
    \item Recall from Exercises~\ref{exe:singletonlanguages}~and~\ref{exe:singletons-dual} that in particular $A^*$ is dense in the dual space $X$ of $\RegA$. Use this to conclude that, for any function $f \colon A \to M$ with $M$ a finite monoid, there can be at most one continuous homomorphism $X \to M$ extending $f$.
  \end{enumerate}
\end{ourexercise}
\begin{ourexercise}\label{exe:cts-mon-hom}
	Let $X$ and $Y$ be Hausdorff topological monoids, and suppose that $D$
	is a dense submonoid of $X$. Let $f \colon X \to Y$ be a continuous
	function such that $f(u)f(v) = f(uv)$ for every $u, v \in D$.
	Prove that $f$ is a homomorphism $X \to Y$.
\end{ourexercise}
\begin{ourexercise}\label{exe:profinite-char}
	Let $M$ be a topological monoid. As in
	Proposition~\ref{prop:profinite-char}, we write $\ClpCon(M)$ for the
	collection of clopen congruences $\theta$ on $M$. We say that
	\emph{the clopen congruences separate points} if, for every $x, y
	\in M$, if $x \neq y$ then there is $\theta \in \ClpCon(M)$ such that $(x,y)
	\not\in \theta$.
\begin{enumerate}
	\item Prove that, if $M$ is the projective limit of a diagram of finite
		discrete monoids, then $M$ embeds as a closed submonoid of
		a product of finite monoids, and is therefore in particular
		compact and the congruences in $\ClpCon(M)$ separate points.
	\item Now assume $M$ is compact and that the clopen congruences separate
		points. Show that the sets of the form
		\[ K_{\theta, x_0} := \{ x \in M \mid x \theta x_0 \} \] 
where $\theta \in \ClpCon(M)$ and $x_0 \in M$, form a base of clopen sets for the
topology on $M$.
\item Let us write $Q$ for the diagram of shape $\ClpCon(M)$ that sends $\theta$ to
	$M/{\theta}$, and that sends an inclusion $\theta \subseteq \theta'$ to
	the unique factorization of $M \onto M/{\theta}$ through $M \onto
	M/{\theta'}$. Using the previous item, show that if $M$ is compact and the congruences
	in $\ClpCon(M)$ separate points, then the collection of maps $(M \onto
	M/{\theta})_{\theta \in \ClpCon(M)}$ is a limiting cone over the diagram
	$Q$.
\end{enumerate}
\end{ourexercise}
\begin{ourexercise}\label{exe:freeprofunique}
 Let $A$ be a set and suppose that $\eta \colon A \to X$ and $\eta' \colon A \to X'$ are both free profinite monoids over $A$. This exercise shows that there is a unique function $\phi \colon X \to X'$ which is both a homeomorphism and a monoid isomorphism and satisfies $\phi \circ \eta = \eta'$.  We use the results proved in Exercise~\ref{exe:profinite-char}.
  \begin{enumerate}
    \item Let $\theta \in \ClpCon(X')$ and write $F := X'/{\theta}$. Consider the function $f \colon A \to F$ defined by sending $a \in A$ to $[\eta'(a)]_{\theta}$. By the universal property of $X$, extend $f$ uniquely to a continuous homomorphism $\widehat{f} \colon X \to F$. Prove that if $\phi \colon X \to X'$ is a continuous homomorphism such that $\phi \circ \eta = \eta'$, then we must have $[\phi(x)]_\theta = \widehat{f}(x)$ for all $x \in X$. \hint{Use the uniqueness of the extension $\wh{f}$.}
    \item Show that there exists a continuous homomorphism $\phi \colon X \to \prod_{\theta \in \ClpCon(X')} X'/{\theta}$.
    \item Show that $\phi$ is injective.
    \item Show that the image of $\phi$ is equal to the image of the  embedding $X' \hookrightarrow \prod_{\theta \in \ClpCon(X')} X'/{\theta}$.
  \end{enumerate}
\end{ourexercise}
\begin{ourexercise}\label{exe:omegaminus1}
Let $M$ be a profinite monoid and $x\in M$.
\begin{enumerate}
	\item Show that the closure of $\{x^n \ \mid \ n \in \mathbb{N}_{\geq 1}\}$ is a commutative    
	subsemigroup of $M$.
	\item Analogously to Proposition~\ref{prop:omega}, prove that there exists a unique element
	$y$ in the closure of $\{x^n \ \mid \ n \in \mathbb{N}_{\geq 1}\}$ such that $yx = xy = x^{\omega}$.
\end{enumerate}
\end{ourexercise}
\begin{ourexercise}\label{exe:RecM}
  Let $M$ be a monoid. As in the case of $A^*$, we say that a subset $L$ of $M$ is \emph{recognizable} if it is recognized by some homomorphism $h \colon M \to N$ with $N$ a finite monoid. We denote the set of recognizable subsets of $M$ by ${\rm Rec}_M$. Show that the results for $A^*$ go through at this level of generality and thus, in particular, that
  \begin{enumerate}
    \item ${\rm Rec}_M$ is a Boolean residuation ideal in $\cP(M)$;
    \item The dual space of the Boolean residuation algebra $({\rm
	    Rec}_M,\backslash,/)$ is the \emph{profinite completion of $M$}.
	    That is, the dual space admits a continuous monoid multiplication,
	    and the resulting monoid $\widehat{M}$ is the projective limit
	    of the diagram of finite monoid quotients of $M$.  The profinite
	    completion of $M$ may also be identified as the image of $M$ under
	    the left adjoint from monoids to topological monoids based on a
	    Boolean space, that is, there is  
	    a monoid homomorphism
	    $\eta \colon M \to \widehat{M}$ such that, for any topological
	    monoid $N$ on a Boolean space and
	    any homomorphism $h \colon M \to N$, there exists a
	    unique continuous homomorphism $\widehat{h} \colon \widehat{M} \to
	    N$ such that $\widehat{h} \circ \eta = h$. 
Compare this with
	    Remark~\ref{rem:priestley-as-procompletion} and
	    Exercise~\ref{exe:limit-adjunction}; further see
	    \cite[Section~4]{Geh16}.
  \end{enumerate}
   Note that ${\rm Rec}_M$ need not be closed under concatenation product in $\cP(M)$, as Exercise~\ref{exe:RecMnotclosedunderproduct} shows.
\end{ourexercise}

\begin{ourexercise}\label{exe:RecMnotclosedunderproduct}
  This exercise uses the definition of recognizable set in an arbitrary monoid as given in Exercise~\ref{exe:RecM}; it is based on \cite[Ex.~III.12.9]{Eilenberg74}, also see \cite[Ex.~2.6]{Sak2009}.
  \begin{enumerate}
    \item Show that $\{0\}$ is not recognizable in the monoid $(\bZ, +)$.
    \item Let $M$ be the structure $(\bZ \cup \{e, x\}, \cdot)$, where $\cdot$ is defined by $k\ell = k + \ell$ for $k,l \in \bZ$, $xk = kx = k$ for all $k \in \bZ$, $x^2 = 0$, and $em = me = m$ for all $m \in M$. Show that $M$ is a monoid.
    \item Let an equivalence relation defined by $m \equiv n$ if, and only if $m = n$ or $m$ and $n$ are both in $\bZ$. Show that ${\equiv}$ is a congruence.
    \item Prove that $L := \{x\}$ is recognizable in $M$ but $LL$ is not.
  \end{enumerate}
\end{ourexercise}
\begin{ourexercise}\label{exe:reglangdual}
  Let  $h\colon A^*\twoheadrightarrow  M$ be a finite monoid quotient and consider the inverse image function $h^{-1} \colon \cP(M) \to \RegA$. Let $X$ denote
  the dual space of $\RegA$.
  \begin{enumerate}
    \item Explain why $h^{-1}(P)$ is regular for every $P \subseteq M$.
    \item By applying Stone duality to the Boolean algebra homomorphism $h^{-1}$, show that $h$ extends uniquely to a continuous function $\bar{h}\colon X\twoheadrightarrow M$.
    \item Show that, for any $w \in A^*$, $\bar{h}(w) = h(w)$, that is, the diagram below commutes.
          
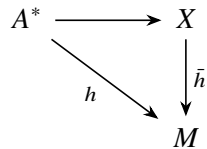
\[
            \begin{tikzcd}[column sep=large,row sep=large,arrows={-Stealth}]
              A^* \arrow[r]\arrow[dr,"h"']& X\arrow[d,"\bar{h}"] \\%
              & M
            \end{tikzcd}
          \]
    \item Show that if $\ell\in\RegA$ is recognized by $h$ via $P\subseteq M$, then
          \[
            \widehat{\ell}=\widehat{h^{-1}(P)}=\bar{h}^{-1}(P).
          \]
  \end{enumerate}
\end{ourexercise}

\begin{ourexercise} \label{exer:filter-from-backslash}
  Let $h \colon L \times L \to L$ be a binary operation on a distributive lattice that is antitone in the first coordinate and preserves finite meets in the second coordinate. Prove that, if $F_1, F_2 \subseteq L$ are filters, then the set
  \[ F := \{ b \in L \ \mid \ \text{ there exists } a \in F_1 \text{ such that } h(a, b) \in F_2 \} \]
  is a filter in $L$. Use this to conclude that the set $F$ defined at the end of the proof of Theorem~\ref{thm:bintopalg-dual-to-jpp} is a filter.
\end{ourexercise}

\section{Equations, subalgebras, and profinite monoids}\label{sec:EilReittheory}
As follows from what we have seen in Sections~\ref{sec:quotients-and-subs}~and~\ref{sec:freeprofmonoid}, given a finite alphabet $A$, the set of all Boolean quotient spaces of the Boolean space underlying the free profinite monoid over $A$ is in one-to-one correspondence with the set of Boolean subalgebras of $\RegA$. What is more is that these are all given by sets of what we called equations, which are really just pairs of elements of the dual space, see Definition~\ref{dfn:spatial-equations} and Corollary~\ref{cor:subalgebras-equations}. The interest of these equations is that, in order to describe a quotient space-subalgebra pair, one can use sets of equations that do not satisfy any special properties -- in contrast with
compatible equivalence relations (Definition~\ref{def:compatible} and Corollary~\ref{cor:subalgebras-equations}), which are more difficult to understand. We recall in particular Example~\ref{exa:equationssubalg}, which showed that in the presence of additional structure, by the use of equations, we may be able to obtain much smaller and simpler sets that characterize the subalgebras than the equivalence relations obtained from the full compatible preorders. This methodology has been greatly exploited in automata and formal language theory, and is in fact behind a great number of \emph{decidability results} in the area. In this section, we want to give you a glimpse of the use of equations and inequations in language theory from a duality-theoretic perspective.

\subsection*{Profinite monoid inequations and equations}
We first recall the definitions connecting sublattices and inequations as given in Section~\ref{sec:quotients-and-subs} (Definition~\ref{dfn:spatial-inequations}), 
specialized to our use case here.
For any element $L \in \RegA$,  we define a binary relation $\preceq_L$ on $\widehat{A^*}$ by
\[ x \preceq_L y \iff y \in \widehat{L} \text{ implies } x \in \widehat{L}.\]
This is a preorder, and $\preceq_L$ is clopen as a subset of $\widehat{A^*}$.
When $\mathcal{L} \subseteq \RegA$, we also define 
\[\preceq_{\mathcal{L}} {:=} \bigcap_{L \in \mathcal{L}} \preceq_L\ .\]
\begin{remark}\label{rem:synt-preorder}
  In the literature on automata and finite monoids, a \emphind{syntactic preorder} is commonly associated with a language $L \subseteq A^*$. The classical definition~\parencite[Section~2]{Schutzenberger56}, still in use today~\parencite[Section~1.4.4]{HandbookI}, is the following. For any $L \subseteq A^*$ and $x, y \in A^*$, define the syntactic preorder $\preceq_L^s$ of $L$ by: $x \preceq_L^s y$ if, and only if, for every $u,v \in A^*$, if $x \in u^{-1}Lv^{-1}$ then $y \in u^{-1}Lv^{-1}$. This means that $x \preceq_L^s y$ if, and only if, $y \preceq_{\mathcal{Q}(L)} x$ according to our definitions, where $\cQ(L)$ denotes the closure under quotienting of $\{L\}$; that is, $\preceq_L^s {=} (\preceq_{\cQ(L)})^\op$. The original language $L$ is then \emph{up}ward closed in the preorder $\preceq_L^s$, while it is \emph{down}ward closed in the preorder $\preceq_L$ that we obtain from duality theory.

  As we have seen throughout this book (see Remarks~\ref{rem:order-choice}~and~\ref{rem:order-yoga}), when applying duality theory to a different field, some conventions are bound to clash with the existing ones in the field, and one needs to choose whether to adapt or co-exist. We here chose to keep the same definition of preorder associated to a sublattice as we did earlier in the book, which is thus the opposite of the definition of syntactic preorder. As it happens, this choice fits well with some of the existing literature on recognition by ordered monoids, in which the definition of syntactic preorder was also reversed, specifically some of the papers that use ordered monoids to analyze the class of piecewise testable languages~\parencite{StrThe88,Pin1995variety,HenckellPin98}. Since piecewise testable languages will also be our focus later in this section, this gives another reason to adhere to the convention that we choose here.
\end{remark}

For any $x, y \in\wh{A^*}$, we define 
\[
  \sem{x\preceq y} := \{L \in \RegA \ \mid \ \text{ if } y \in
  \widehat{L}, \text{ then } x \in \widehat{L} \},
\]
and note that this is a sublattice of $\RegA$. As in Definition~\ref{dfn:spatial-inequations}, we also use the notation $L \models x \preceq y$ to mean that $L \in \sem{x \preceq y}$. 
If $E$ is a set of pairs of $\wh{A^*}$, we also define 
\[ 
\sem{E} := \bigcap_{(x,y) \in E} \sem{x \preceq y}.
\] 
The content of Proposition~\ref{prop:sublattice-quotientspace}, in this setting, is that we have
an adjunction 
\begin{equation}\label{eq:sub-quotient-adjunction} \preceq_{-} \colon \mathcal{P}(\RegA) \leftrightarrows 
\mathcal{P}(\wh{A^*} \times \wh{A^*})^\op \colon \sem{-} \end{equation}
whose fixed points on the left are sublattices of $\RegA$, and whose fixed points on the
right are compatible preorders on $\wh{A^*}$; see
Definition~\ref{def:compatible}. 
Recall that, for a subset
$S$ of $\RegA$, $\cQ(S)$ denotes the closure under quotienting of $S$. 
When $E$ is a relation on $\wh{A^*}$, let us denote by $m(E)$ the
\emph{bi-action invariant closure} of $E$, that is, the relation defined by
\[ m(E) := \{ (uxv, uyv) \ \mid \ (x,y) \in E, u, v \in \wh{A^*} \}.\] 
We now prove that the operations of closure under quotienting and closure under
the bi-action are dual to each other, in the following
sense.
\begin{proposition}\label{prop:mult-equation}
	For any set of pairs $E \subseteq \wh{A^*} \times \wh{A^*}$ and any subset $S
	\subseteq \RegA$, we have 
	\[ S \subseteq \sem{m(E)} \iff \cQ(S) \subseteq \sem{E}.\]
	In particular, if $\sem{E}$ is closed under quotienting, then $\sem{E} =
	\sem{m(E)}$.
\end{proposition}
\begin{proof}
	We first show that, for any $L \in \RegA$ and $x, y \in \wh{A^*}$, if
	$\cQ(L) \subseteq \sem{x \preceq y}$, then $L \models uxv \preceq uyv$
	for every $u,v \in \wh{A^*}$. Indeed, if $\cQ(L) \subseteq \sem{x
	\preceq y}$, then at least for every $u,v \in A^*$, we have $u^{-1}Lv^{-1} \models x \preceq y$, which means that $L \models uxv \preceq uyv$.
	To show that this extends to all of $\wh{A^*}$, consider the
	function $f \colon \wh{A^*} \times \wh{A^*} \to \wh{A^*} \times
	\wh{A^*}$ defined by $f(u,v) := (uxv, uyv)$. This function is continuous
	because the multiplication on $\wh{A^*}$ is continuous. Denote by $C$ the
	inverse image of the set $(\wh{L}^c \times \wh{A^*}) \cup (\wh{A^*}
	\times \wh{L})$; this set is closed, and the assumption that $\cQ(L) \subseteq \sem{x
	\preceq y}$ gives that $A^* \times A^*$ is contained in $C$. Thus, since
	$A^* \times A^*$ is dense in $\wh{A^*} \times \wh{A^*}$, we get that $C
	= \wh{A^*} \times \wh{A^*}$, meaning that for any $(u,v) \in\wh{A^*}
	\times \wh{A^*}$, we have $L \models uxv \preceq uyv$,
	that is, $L \in \sem{m(\{x \preceq y\})}$. 

	This establishes the proposition in case $S$ and $E$ are singletons. For
	the general case, note that $S \subseteq \sem{m(E)}$ if, and only if, $L
	\in \sem{m(\{x \preceq y\})}$ for every $L \in S$ and every $(x,y) \in
	E$, and $\cQ(S) \subseteq \sem{E}$ if, and only if, $\cQ(L) \subseteq
	\sem{(x,y)}$ for every $L \in S$ and every $(x,y) \in E$.

	For the in particular statement, if $\sem{E}$ is closed under
	quotienting, then taking $S := \sem{E}$ in the equivalence gives $\sem{E} \subseteq
	\sem{m(E)}$, and the other inclusion is clear because $E \subseteq
	m(E)$ and the map $\sem{-}$ is antitone.
\end{proof}
\nl{$L \models x \stackrel{*}{\preceq} y$}{a regular language $L$ satisfies a profinite monoid inequation given by the pair $(x,y)$ of elements of the free profinite monoid}{}
\nl{$\sem{x \stackrel{*}{\preceq} y}$}{the set of regular languages $L$ which satisfy a profinite monoid inequation given by the pair $(x,y)$ of elements of the free profinite monoid}{}
For $x,y \in \wh{A^*}$, let us write $L \models x \stackrel{*}{\preceq} y$ 
if, and only if, $L \in
\sem{m(\{x \preceq y\})}$, that is, $L \models uxv \preceq uyv$ for all $u, v
\in \wh{A^*}$. We say in this case that $L$ \emphind{satisfies the profinite monoid
inequation} $x \stackrel{*}{\preceq} y$, and we write 
\[ \sem{x \stackrel{*}{\preceq} y} := \{L \in \RegA \ \mid \ L \models x \stackrel{*}{\preceq} y \} \ .\]
Similarly, we write $L \models x
\stackrel{*}{=} y$ if, and only if, $L \models uxv \preceq uyv$ and $L \models
uyv \preceq uxv$ for all $u, v \in \wh{A^*}$. In this case, we say that $L$
\emphind{satisfies the profinite monoid equation} $x \stackrel{*}{=} y$.
\nl{$L \models x \stackrel{*}{=} y$}{a regular language $L$ satisfies a profinite monoid equation given by the pair $(x,y)$ of elements of the free profinite monoid}

It follows in particular from Proposition~\ref{prop:mult-equation} that, in the
adjunction (\ref{eq:sub-quotient-adjunction}), the sublattices on the left that
are closed under quotienting correspond exactly to the compatible preorders on
the right that are \emphind{monoid-compatible}, in the sense of the following
definition.  Note that this notion generalizes the definition of congruence, which is just a
monoid-compatible equivalence relation.

\begin{definition}\label{dfn:monoid-compatible}
	A binary relation $R$ on a monoid $M$ is \emphind{monoid-compatible} if,
	for every $u,v, x \in M$, if $u{R}v$, then $ux {R} vx$ and $xu {R} xv$.
\end{definition}

We now show that profinite monoid equations for a regular language
may be `tested' on the 
syntactic monoid. A similar result holds for 
profinite monoid inequations and ordered monoids 
(see Remark~\ref{rem:synt-preorder}), but we do not fully develop that theory  here.

Recall that,  when $h \colon A^* \to M$ is a monoid homomorphism,
with $M$ a finite monoid, we denote by
$\bar{h} \colon \wh{A^*} \to M$ its unique continuous extension, which is also
a homomorphism.
\begin{theorem}\label{thm:profinite-equation-synt-mon}
  Let $L \in \RegA$ and $(x,y) \in \wh{A^*} \times \wh{A^*}$.
  The following are equivalent:
  \begin{enumerate}[label=(\roman*)]
    \item the language $L$ satisfies the profinite monoid equation $x \stackrel{*}{=} y$;
    \item the (finite) Boolean residuation ideal $\cB(L)$ is contained in $\sem{x \approx y}$;
    \item for the syntactic homomorphism $h_L \colon A^* \onto M_L$ of $L$, we have $\overline{h_L}(x) = \overline{h_L}(y)$;
    \item there exists a homomorphism $h \colon A^* \to M$, with $M$ a finite monoid, which recognizes $L$ and such that $\bar{h}(x) = \bar{h}(y)$;
  \end{enumerate}
\end{theorem}
\begin{proof}
  By definition, (i) means that $L \in \sem{m(\{x \preceq y, y \preceq x\})}$.
  By Proposition~\ref{prop:mult-equation}, this is equivalent to $\cQ(\{L\}) \subseteq \sem{x \approx y}$.
  Since $\sem{x \approx y}$ is a Boolean subalgebra of $\RegA$, and $\cB(L)$ is the Boolean subalgebra of $\cP(A^*)$ generated
  by $\cQ(\{L\})$ (Proposition~\ref{prop:sublat-quotienting}), this gives that (i) is equivalent to (ii). 

  For (ii) $\implies$ (iii), recall that the languages recognized by $M_L$ are precisely those in $\cB(L)$ (Corollary~\ref{cor:syntactic}). Thus in particular
  the language $K := h_L^{-1}(\{\overline{h_L}(x)\})$ is in $\cB(L)$, and $x \in \widehat{K}$, so (ii) gives that $y \in \widehat{K}$. Using Exercise~\ref{exe:reglangdual}, we have that $\widehat{K} = \overline{h_L}^{-1}(\{\overline{h_L}(x)\})$, so in particular $\overline{h_L}(x) = \overline{h_L}(y)$, as required.
 
  (iii) $\implies$ (iv) is trivial, since $h_L$ recognizes $L$.

  For (iv) $\implies$ (i), note that we have, for any $u, v \in \wh{A^*}$,
  \[ \bar{h}(uxv) = \bar{h}(u) \bar{h}(x) \bar{h}(v) = \bar{h}(u) \bar{h}(y) \bar{h}(v) = \bar{h}(uyv), \] 
  using that $\bar{h}$ is a homomorphism (Corollary~\ref{cor:R-funct}).
  Using  Exercise~\ref{exe:reglangdual} again, $\wh{L} = \bar{h}^{-1}(P)$, where $P$ is a subset of $M$ 
  such that $L = h^{-1}(P)$. 
  Hence, for any $u, v \in \wh{A^*}$, we have
  \[ xuv \in \wh{L} \iff \bar{h}(uxv) \in P \iff \bar{h}(uyv) \in P \iff xyv \in \wh{L} \ , \] 
  as required.
\end{proof}

A major question in the theory of automata is whether a certain subclass has a decidable \emphind{membership problem} within the class of regular languages. More precisely, given an automaton, is there an effective procedure that allows us to decide whether or not the language recognized by the automaton is in the given class. A powerful method, coming from the algebraic approach to language theory, is to give an \emph{equational} criterion on the syntactic monoid of a regular language $L$ for belonging to the given class. 

Indeed, Proposition~\ref{prop:mult-equation} shows that profinite monoid quotients of $\wh{A^*}$ correspond dually to the Boolean subalgebras of $\RegA$ that are closed under quotienting. 
This is the kernel of Eilenberg's Theorem \parencite{Eilenberg76} seen from a duality theoretic point of view. Eilenberg considered \emph{pseudovarieties of regular languages}\index{pseudovariety of regular languages}, which are families $V = (V_A)_{A \in \mathbf{Set}_f}$ of Boolean subalgebras of $\RegA$ closed under quotienting, satisfying the further property that the Boolean algebras $V_A$ are stable under inverse image by homomorphisms between free monoids. An extension of Proposition~\ref{prop:mult-equation} then shows that the class of syntactic monoids of languages in $V$ is definable by profinite equations that are moreover invariant under substitutions. Such classes of finite monoids are called \emph{pseudovarieties of finite monoids}\index{pseudovariety of finite monoids}, and, by a finite version of Birkhoff's theorem, are exactly those classes of finite monoids closed under homomorphic images, submonoids, and finite products.
Eilenberg's Theorem can thus be viewed a duality between pseudovarieties of regular languages and pseudovarieties of finite monoids; for further details, see \cite{Geh16}.

Rather than exposing this general theory any further here, we illustrate the method by working out a particular case, namely that of piecewise testable languages. In this case, a decidable criterion for membership was first obtained in \cite{Simon1975}, who showed that the piecewise testable languages are exactly the languages recognized by finite \emph{$\cJ$-trivial}\index{J-trivial@$\cJ$-trivial} monoids. Our aim in the rest of this section is to obtain a profinite equational characterization of piecewise testable languages by duality-theoretic means and we will explain how this shows that the class has a decidable membership problem (Corollary~\ref{cor:Simon}).

\subsection*{Positively piecewise testable languages}
In order to define the \emphind{positively piecewise testable} and \emphind{piecewise testable} languages, we need the \emph{subword ordering}.  To
this end, it is convenient to think of a finite word $w \in A^*$ as a function
$|w| \to A$, where $|w|$ is identified with the totally ordered set $1 < 2 <
\cdots < |w|$, and the function sends $i$ to the $i^\mathrm{th}$ letter of
$w$. In this context, we often also write $w(i)$ for the $i^\mathrm{th}$ letter
of $w$. We say that $u \in A^*$ is a (scattered) \emphind{subword} of $w$ if
there exists an order embedding $\phi$ from $|u|$ to $|w|$ such that $w(\phi(i))
= u(i)$ for every $i \in |u|$, that is, if $u$ is a subsequence of $w$. For
example, the word $abbc$ is a subword of $acabbac$, but not of $acaba$, and not
of $abc$ either. We will write $u \subword w$ if $u$ is a subword of $w$, and
we note that $\subword$ defines a partial order on $A^*$.  
\nl{$\subword$}{the subword partial order}{}
\begin{definition} Let $A$ be a finite alphabet and $L \subseteq A^*$ a
	language. A language is called \emphind{positively piecewise testable}\footnote{What we call a positively piecewise testable language is also known as a \emphind{shuffle ideal}~\parencite{Pin1995variety}, or as a language in the (marked) \emphind{polynomial closure} of $\{\emptyset, A^*\}$~\parencite{Pin2011theme}.}
	if it is a finitely generated up-set in the partial order $(A^*,\subword)$, and
	\emphind{piecewise testable} if it is a Boolean combination of
	positively piecewise testable languages.  We denote by $\PTplus_A$ the
lattice of positively piecewise testable languages and by $\PT_A$ the Boolean
algebra of piecewise testable languages.  
\end{definition}
\nl{$\PTplus$}{the collection of positively piecewise testable languages (in a fixed finite alphabet)}{}
\nl{$\PT$}{the collection of piecewise testable languages (in a fixed finite alphabet)}{}
For the rest of this section, we fix a finite alphabet $A$ and write $\PTplus$ and
$\PT$ for $\PTplus_A$ and $\PT_A$, respectively. We also write $\Reg$ for
$\RegA$. In formal language theory, a language class is often seen as a fibered
collection of lattices, or Boolean algebras, where the alphabet may vary and morphisms between alphabets induce
morphisms between the lattices. We do not need this structure in this section,
and we only remark here that it is related to the fibrational approach via hyperdoctrines that we
will point to at the end of Section~\ref{sec:openmult}.

Any positively piecewise testable language is regular, since it is easy to check
with a finite automaton whether any of a finite number of subwords appear in a
word (see Exercise~\ref{exe:subword-regular}). It follows that any piecewise
testable language is regular,  since $\Reg$ is a Boolean algebra
(Proposition~\ref{prop:RegBoolResidl}).

Towards the equational characterization of piecewise testable languages, we
first show how the collection $\PTplus$ may be identified with an \emph{in}equational
property on the free profinite monoid. 
To this end, first note that the collection $\PTplus$ of finitely generated 
up-sets of $(A^*, \subword)$
is a sublattice of $\cP(A^*)$. By Exercise~\ref{ex:fin-MUB-comp} this 
statement is equivalent to the statement that
$(A^*, \subword)$ is finitely \MUB-complete in the sense of
Definition~\ref{def:finite-mub-comp}, that is, for any finite set $F \subseteq
A^*$, there is a finite set $G \subseteq A^*$ such that ${\uparrow} G =
\bigcap_{w \in F} {\uparrow} w$. The latter is indeed the case, but in fact something much
stronger is true. Even though $(A^*, \subword)$ is not bifinite, it is what is
known as a \emphind{well-quasi-order}, that is, \emph{every} up-set of $(A^*,
\subword)$ is finitely generated. Indeed, this was shown early on by Higman.
For completeness, we include a proof. Note that while the following proof uses
the Axiom of Choice, one can also give a constructive proof, for example by induction
on the size of the alphabet, see \cite{MurRus1990}.

\begin{theorem}[\cite{Higman1952}, Theorem 4.4]\label{thm:higman}
For any finite alphabet $A$, the poset $(A^*, \subword)$ is a well-quasi-order.
\end{theorem}

\begin{proof}
Let $\cS=\{U\in\Up(A^*,\subword)\mid \min(U)\text{ is infinite}\}$. If $\cS$ is empty, 
then we are done. Suppose $\cS\neq\emptyset$. Let $\cC$ be a chain in $\cS$. 
Clearly $\bigcup\cC$ is an up-set. Notice that every up-set is the up-set of its 
minimal elements since the down-set of any word is finite. Suppose $\bigcup\cC=
{\uparrow}F$ with $F$ finite. Then, since $\cC$ is a chain, there is a single $V\in\cC$ 
with $F\subseteq V$. But then $V\subseteq\bigcup\cC={\uparrow}F\subseteq V$, 
which is a contradiction since $V\in\cS$. So $\bigcup\cC\in\cS$ and by Zorn's 
Lemma (Lemma~\ref{lem:zorn}), there is a maximal element $U\in\cS$. We 
now have 
\[
\min(U)=\bigcup\{\min(U)\cap aA^*\mid a\in A\}
\]
and thus $U_a=\min(U)\cap aA^*=a V_a$ is infinite for some $a\in A$. Now let 
$W= U\cup{\uparrow} V_a$. Then $V_a$ is infinite and each element in $V_a$ 
is minimal in $W$. To see this, notice that if $w\in V_a$, then $aw\in\min(U)$ 
and thus $w\not\in U$. Further, if $v\in V_a$ and $v\subword w$, then 
$av\subword aw$ and $av,aw\in\min(U)$, so $v=w$. Thus $W\in\cS$ and 
$V_a\not\subseteq U$ so that $U\subsetneq W$, which contradicts the 
maximality of $U$.
\end{proof}

Since $\PTplus$ is the lattice of up-sets of a well-quasi-order, we
note that we are in a special case of the spectral domains considered in
Section~\ref{sec:dom-Stone}. Thus, by Corollary~\ref{cor:2/3sfp} and
Exercise~\ref{ex:dual-spectral-domain}, the spectral space dual to $\PTplus$
is isomorphic to the poset $\Idl(A^*)$ in the Scott topology.
Since $\PTplus$ is a sublattice of the Boolean algebra of regular languages,
its dual space $\Idl(A^*,\subword)$ can also be characterized by inequations 
on the dual space $\wh{A^*}$ of $\Reg$. 

Note in particular that $\PTplus$ is closed under quotienting. Indeed, for any
words $u, v \in A^*$ and $L \in \PTplus$, we have $u^{-1} L v^{-1} \in
\PTplus$: if $uwv \in L$ and $w \subword w'$, then $uwv \subword uw'v$, so
$uw'v \in L$ as well. Thus, in order to obtain a characterization of $\PTplus$
via profinite inequations, Proposition~\ref{prop:mult-equation} implies that 
we may look for a characterization via profinite \emph{monoid} inequations.  

\begin{proposition}\label{prop:PTpluseq}
The lattice $\PTplus$ contains exactly the regular languages that satisfy the
profinite monoid inequation $x \stackrel{*}{\preceq} \epsilon$, for every $x \in
\wh{A^*}$.
\end{proposition}

\begin{proof}
Let us first show that any language $L$ in $\PTplus$ satisfies the profinite
monoid equation. Since $\PTplus$ is closed under quotienting, it suffices to
show that $\PTplus \subseteq \sem{x \preceq \epsilon}$, by
Proposition~\ref{prop:mult-equation}. But if $L$ is in $\PTplus$, then $\epsilon
\in L$ clearly implies $L = A^*$, since $\epsilon \subword w$ for every $w \in
A^*$, so that $\widehat{L} = \wh{A^*}$, and thus $x \in
\wh{L}$ for every $x \in \wh{A^*}$.

Conversely, suppose that $L \models x \stackrel{*}{\preceq} \epsilon$ for every
$x \in \wh{A^*}$. Suppose that $w \subword w'$ and $w \in L$. By induction on
the difference in length between $w'$ and $w$, it suffices to treat the case
where $w = uv$ and $w' = uav$ for some $u, v \in A^*$ and $a \in A$. Now, since
$L \models a
\stackrel{*}{\preceq} \epsilon$, we also get $L \models uav \preceq uv$.
Therefore, since $uv \in L$, we have $uav \in L$.
\end{proof}

Propositions~\ref{prop:mult-equation}~and~\ref{prop:PTpluseq} now yield
the following.

\begin{theorem}\label{thm:ptplus-dual}
	The Priestley dual space of $\PTplus$ is homeomorphic to the quotient
	space $\wh{A^*}/{\preceq}$, where $\preceq$ is the smallest compatible 
	preorder containing the profinite monoid inequations $x \stackrel{*}{\preceq} 
	\epsilon$ for every $x \in \wh{A^*}$.
\end{theorem}
Note that this result in particular gives an equational way of looking at the
spectral domain $\Idl(A^*)$, as follows. As noted above, the spectral space dual to
$\PTplus$ can also be described as $\Idl(A^*)$ equipped with the Scott
topology. Thus, using Theorem~\ref{thm:Stone-isom-Priestley},
we also get the following.
\begin{corollary}\label{cor:priestley-ptplus}
The Priestley dual space of $\PTplus$ is homeomorphic to
$(\Idl(A^*), \sigma^p, \supseteq)$, where $\sigma^p := \sigma \vee \sigma^\partial$ and
$\sigma$ is the Scott topology on $(\Idl(A^*),\subword)$.
\end{corollary}
We further have a canonical map $\pi \colon \wh{A^*} \onto \Idl(A^*)$, dual to
the inclusion $\PTplus \into \Reg$, given by
continuously extending the inclusion $A^* \into \Idl(A^*)$ that sends $w \in
A^*$ to the principal ideal ${\uparrow} w$. Concretely, this map sends $u \in
\wh{A^*}$ to the ideal $\{w \in A^* \ \mid \ u \in \wh{{\uparrow}w}\}$ (see
Exercise~\ref{exe:ptplus-idl}). 
Theorem~\ref{thm:ptplus-dual} implies that $\pi(u) \leq \pi(v)$ if, and only
if, $u \preceq_{\PTplus} v$. 
Defining an appropriate monoid multiplication on $\Idl(A^*)$, the map $\pi$
becomes a continuous monoid homomorphism (again see 
Exercise~\ref{exe:ptplus-idl} for more details).

\subsection*{Piecewise testable languages}
Now, moving towards a first characterization of the dual space of $\PT$, 
recall from Proposition~\ref{prop:boolenv} and Corollary~\ref{cor:boolenv-hull} 
that, since $\PT$ is the Boolean envelope of $\PTplus$, the Stone dual space 
of $\PT$ is simply given by forgetting the order of the Priestley dual space of 
$\PTplus$.  We thus also get the following result.
\begin{corollary}\label{cor:pt-dual}
	The dual space of $\PT$ is 
	homeomorphic to $\wh{A^*}/{(\preceq \cap \succeq)}$, where $\preceq$ is
	the preorder of Theorem~\ref{thm:ptplus-dual}.
\end{corollary}
\begin{remark}
  It was recently proved that the dual space of the Boolean algebra of 
  piecewise testable languages over an alphabet $A$ with $|A| = n$ 
  is homeomorphic to the ordinal 
  $\omega^{\omega^{n-1}} + 1$ with the interval topology 
  \parencite[Theorem~11]{Pouzet2023}, based on earlier, more general 
  results by \cite{BPZ2007}. 
\end{remark}
While Corolary~\ref{cor:priestley-ptplus} tells us something about what the dual space
of $\PTplus$ looks like, Corollary~\ref{cor:pt-dual} and Theorem~\ref{thm:ptplus-dual} tells us 
something about the (in)equations satisfied by $\PTplus$ and $\PT$, and the latter 
type of information is what can lead to a decidable characterization.
In order to obtain such a decidable equational characterization for piecewise testable
languages, we need a better description of the congruence ${\preceq \cap
\succeq}$ of Corollary~\ref{cor:pt-dual}. Let us introduce a notation for this
congruence.
\begin{definition}
	For $A$ a finite alphabet, we write $\approx_{\PT}$ for the
	congruence on $\wh{A^*}$ dual to the subalgebra $\PT \leq \Reg$.  
\end{definition}

We first identify two properties of
this congruence that will be important for this description. These properties
take the form of \emph{quasi-identities}\index{quasi-identity}, which is the term commonly used in
universal algebra for an expression of the form 
\[ \text{ if } \big(s_1 = t_1 \text{ and} \ \dotsm \ \text{and } s_n = t_n\big), \text{ then } (s_0 = t_0) \ , \] 
where the $s_i$ and $t_i$ are terms.
\begin{lemma}\label{lem:PT-quasieq}
	Let $\preceq$ be a monoid-compatible preorder on a monoid $M$ 
	such that $x\preceq 1$ for every $x \in M$, and denote its associated 
	equivalence relation by $\approx$.  Then, for any $x, y, u \in M$ we have 
	\begin{enumerate}
		\item if $uxy \approx u$, then $ux \approx u$, and
		\item if $yxu \approx u$, then $xu \approx u$.
	\end{enumerate}
\end{lemma}
\begin{proof}
	Let $x, y, u \in M$. Note that, since $\preceq$ is monoid-compatible 
	and $x \preceq 1$ for every $x\in M$, we also get $ux\preceq u\cdot 1=u$ 
	and $uxy \preceq ux \cdot 1 = ux$. Now, for the first item, suppose that 
	$uxy \approx u$. Then 
\[ u \approx uxy \preceq ux \preceq u,\]
	so $ux \approx u$. The proof of the second item is symmetric.
\end{proof}

Using Theorem~\ref{thm:ptplus-dual}, 
the congruence $\approx_\PT$ dual to $\PT$ satisfies both of the properties 
of Lemma~\ref{lem:PT-quasieq}. We now aim to show that these properties are
equivalent to profinite monoid equations, and that the topological monoid
congruence $\approx_\PT$ is in fact \emph{generated} by those equations. 
This will yield a
proof of Simon's theorem and, as we will see, the decidability of the membership 
problem for piecewise testable languages.

Let us first show how the properties of Lemma~\ref{lem:PT-quasieq} can be
written as profinite monoid equations. Here, the operations $\omega$ and $\omega-1$ that 
we introduced in Definition~\ref{def:omega} and immediately after (see also 
Exercise~\ref{exe:omegaminus1}) are important, as they allows us to construct
profinite terms.
\begin{lemma}\label{lem:PT-equations}
Let $\approx$ be a closed congruence on a profinite monoid $M$.  
\begin{enumerate}
	\item The following are equivalent:
\begin{enumerate}[label=(\roman*)]
	\item For any $u, x, y \in M$, if $uxy \approx u$ then $ux \approx u$.
	\item For any $x, y \in M$, $(xy)^\omega x \approx (xy)^\omega$.
\end{enumerate}
\item The following are equivalent:
\begin{enumerate}[label=(\roman*)]
	\item For any $u, x, y \in M$, if $yxu \approx u$ then $xu \approx u$.
	\item For any $x, y \in M$, $x(yx)^\omega  \approx (yx)^\omega$.
\end{enumerate}
\end{enumerate}
Moreover, if any of these properties hold, then $z^{\omega+1} \approx z^{\omega}$ 
for any $z\in M$.
\end{lemma}
\begin{proof}
	We only prove that (i) and (ii) in (a) are equivalent, and imply that
	$z^{\omega+1} \approx z^{\omega}$, the proof for (b) is symmetric.
	Suppose (i) holds. We first show that
	$z^{\omega +1}  \approx z^\omega$ for any $z \in X$. Indeed, 
	applying (i) with $u = z^\omega$, $x = z$, and $y =
	z^{\omega-1}$, we have $uxy = z^\omega z^\omega = u$, so that
	$z^{\omega+1} \approx z^\omega$. 
	In particular, $(xy)^\omega xy \approx (xy)^\omega$, so that (i) now
	gives $(xy)^\omega x \approx (xy)^\omega$.
	
	Conversely, suppose (ii) holds and suppose that $uxy \approx u$. Then,
	for any $n$, we have $u \approx u(xy)^n$, so, since $\approx$ is closed,
	$u \approx u(xy)^\omega$. Thus,
	\[ ux \approx u(xy)^\omega x \approx u(xy)^\omega \approx u.\qedhere\]
	\end{proof}

\begin{remark}\label{rem:green}
The two profinite equations (a)(ii) and (b)(ii) are known in the literature as the profinite equations for
\emph{$\cR$-trivial}\index{R-trivial@$\cR$-trivial} and
\emph{$\cL$-trivial}\index{L-trivial@$\cL$-trivial} monoids,
respectively, and it is known that together they characterize the class of 
\emph{$\cJ$-trivial}\index{J-trivial@$\cJ$-trivial} finite monoids. The equation $z^\omega \approx
z^{\omega+1}$ that appeared in Lemma~\ref{lem:PT-equations} characterizes the
class of \emph{aperiodic} or
\emph{$\cH$-trivial}\index{aperiodic}\index{H-trivial@$\cH$-trivial} monoids. The letters
$\cL$, $\cR$, $\cJ$ and $\cH$ refer to \emphind{Green's equivalence relations} on 
monoids, and the adjective `trivial' asserts that these relations are equal to
the diagonal in the profinite monoid under consideration. We will not need to consider
these relations in any more detail here, and our development is independent of
any of the facts mentioned in this remark. 
We refer the reader who wants to 
know more about Green's relations to Exercise~\ref{exe:Green}, and 
standard texts on monoid theory, such as 
for example \cite{Eilenberg74, Eilenberg76,Almeida1995,RhodesSteinberg2008, JEP-MPRI}.
\end{remark}

We here make a slightly non-standard definition of $\cJ$-trivial profinite monoid 
which however coincides with the usual definition, see Exercise~\ref{exe:Green}.
\begin{definition}\label{def:Jtrivial}
Denote by $\approx_{\cJ}$ the smallest profinite monoid congruence on $\wh{A^*}$ 
that contains, for every $x,y\in\wh{A^*}$, the pairs 
\[ ((xy)^\omega x, (xy)^\omega) \text{ and } (x(yx)^\omega, (yx)^\omega) \ . \]
Denote the residuation ideal of $\RegA$ dual to $\approx_{\cJ}$ by
$\JT$. We call a profinite monoid $M$ $\cJ$-trivial provided for all $x,y\in M$, we have  
$(xy)^\omega x=(xy)^\omega$ and $x(yx)^\omega=(yx)^\omega$, and we call the 
elements of $\JT$ the \emph{$\cJ$-trivial languages}\index{J-trivial@$\cJ$-trivial!language}
 of $\Reg$.
\end{definition}
By Theorem~\ref{thm:profinite-equation-synt-mon}, a regular language $L$ is $\cJ$-trivial if, and only if, its syntactic monoid is $\cJ$-trivial if, and only if, it is recognized by some $\cJ$-trivial monoid.

Our final aim in this section is to prove that $\PT=\JT$ or equivalently that 
${\approx_{\cJ}}={\approx_{\PT}}$. The inclusion $\PT \subseteq \JT$ actually 
follows immediately from Lemmas~\ref{lem:PT-quasieq}~and~\ref{lem:PT-equations}:

\begin{corollary}\label{cor:PTsoundness}
For any $x, y \in 
\wh{A^*}$, we have
\[ 
(xy)^\omega x \approx_{\PT} (xy)^\omega, \text{ and } x(yx)^\omega
\approx_{\PT} (yx)^\omega.
\]
That is,  ${\approx_{\cJ}} \subseteq {\approx_{\PT}}$ and dually $\PT\subseteq\JT$.
\end{corollary}

It now remains to show that $\JT\subseteq\PT$. 
Since $\PT$ is defined as the Boolean algebra
generated by the lattice $\PTplus = \cU(A^*,\subword)$, we are in the
situation of Section~\ref{sec:altchains}.  Theorem~\ref{thm:discrete-alt-height}
therefore gives that a language $L$ is in $\PT$ if, and only if, it has
bounded alternation height in the subword order $\subword$. We will prove that
any language $L$ that satisfies the profinite monoid equations defining $\approx_\cJ$ (and thus the
quasi-identities of Lemma~\ref{lem:PT-quasieq}) has bounded alternation height.

For this final step, we need to develop a small amount of finite monoid
theory, which gives just a flavor of this rich field.
Before giving the formal details of the proof, we first give an intuition.
Let $w = a_1\ \dotsm\ a_k$ be a finite word and let $\underline{(\ )}\colon A^*\to M$ be a 
homomorphism to a finite monoid. We will analyze how the value $\underline{w}$ in $M$ 
can be computed `from left to right'; a similar analysis applies by computing the value 
`from right to left'. Reading $w$ from left to right, the value $\underline{w}$ can be 
computed, starting from $1_M$,  by first computing $\underline{a_1}$, then 
$\underline{a_1a_2}$, then $\underline{a_1a_2a_3}$, and so on. In this
computation, which takes $k$ steps, only certain steps will change the value;
positions where this happens will be called `unstable' below. By definition,
only the unstable positions of a word contribute to its value, and we may thus
\emph{reduce} a word $w$ to obtain a subword $r(w)$ that only contains its
(right) unstable positions. The idea of the rest of the proof is that, if $M$ is a finite monoid 
that satisfies the quasi-identities, then the number of reduced words is
bounded. This will give a bound on the alternation height as a subset of the poset 
$(A^*,\subword)$ of any regular language $L$ that is recognized by $M$.

Now, more formally, we make the following definitions. A reader familiar with
finite monoid theory may recognize the similarity of the reduction notion that we
introduce here with the Karnofsky-Rhodes expansion, see for example
\cite[Section~2]{RSS2022}.

\begin{definition}
Let $M$ be a finite monoid and $\underline{(-)}
\colon A^* \to M$ a homomorphism; we write $\underline{w}$ for the image of $w \in A^*$ 
under the homomorphism. For any finite word $w = a_1 \ \dotsm\ a_k \in A^*$ with each 
$a_i\in A$ and for each $0 \leq i \leq k$, we write $w_i$ for the length $i$ prefix of $w$, that 
is, $w_0 := \epsilon$ and $w_i := w_{i-1}a_i$ for all $1 \leq i \leq k$. We will say a position 
$i \in \{1,\dots,k\}$ is (right) \emph{stable} if $\underline{w_i} = \underline{w_{i-1}}$, and
(right) \emph{unstable} otherwise. 

Let $t_1 <\cdots < t_\ell$ be an enumeration of the unstable positions of $w$. We 
define $r(w) := a_{t_1}\ \dotsm\ a_{t_\ell}$, the subword of $w$ on the unstable positions. 
A simple induction on the number of unstable positions in $w$ shows that the words $w$ 
and $r(w)$ have the same value, that is, $\underline{w} = \underline{r(w)}$. We call $r(w)$ 
the \emph{right reduction} of $w$, and we call a word in $A^*$ \emph{right reduced} if all 
of its positions are unstable, that is, if $w=r(w)$.

The notions of \emph{left (un)stable}, \emph{left reduction}, and \emph{left
reduced} are defined analogously, using the sequence of suffixes of $w$. We
write $\ell(w)$ for the left reduction of $w$.
\end{definition}

 We begin by showing that the quasi-identity identified in Lemma~\ref{lem:PT-quasieq}.a 
 implies that, if $M$ is finite, then there are finitely many right reduced words for any 
 homomorphism from $A^*$ to $M$. We only state and prove the direction 
that we need in the proof, however, the converse is also true (see 
Exercise~\ref{exe:R-triv-reduced}).

\begin{lemma}\label{lem:finitereduced}
Let $\underline{(-)}\colon A^* \to M$ be a homomorphism to a finite $\cJ$-trivial monoid $M$. 
Then there are finitely many right reduced words for the homomorphism $\underline{(-)}$.
\end{lemma}
\begin{proof}
  We show that, for any right reduced word $w = a_1 \ \dotsm\ a_k$ in $A^*$ of length $k$, 
  and for any $1 \leq i < j\leq k$, the prefixes $w_i$ and $w_j$ have distinct values. This 
  will then imply in particular that a right reduced word can only have length at most $|M|$, 
  and thus there are finitely many right reduced words. 

Let $w = a_1 \ \dotsm\ a_k \in A^*$ be a right reduced word  of length $k$, and suppose, 
towards a contradiction, that there exist $1 \leq i < j \leq k$ such that $\underline{w_i} = 
\underline{w_j}$. This means that 
	\[
	 \underline{a_1 \cdots a_i a_{i+1} \cdots a_j} = \underline{a_1\cdots a_i}.
	 \]
Thus, applying the quasi-identity of Lemma~\ref{lem:PT-equations}.a(i) with 
$u := a_1 \cdots a_i$, $x := a_{i+1}$, $y := a_{i+2} \cdots a_j$, we see that 
$\underline{uxy} = \underline{u}$, so $\underline{ux} = \underline{u}$. But this means 
that $\underline{w_i} = \underline{w_{i+1}}$, so that $i+1$ is a right stable point, 
contradicting the assumption that $w$ is right reduced.
\end{proof}

\begin{lemma}\label{lem:fallpoints}
	Let $\underline{(-)} \colon A^* \to M$ be a homomorphism to a finite monoid $M$, 
	$w \in A^*$, $a \in A$ and suppose that the word $r(w) a$ is a subword of $w$. 
	Then $\underline{wa} = \underline{w}$.
\end{lemma}
\begin{proof}
Let $w \in A^*$, $a \in A$ and suppose that the word $r(w) a$ is a subword of $w$.
If $r(w)=\varepsilon$, then $1_M=\underline{b}=\underline{w}$ for all $b\in A$ occurring 
in $w$ and thus, in particular, $\underline{wa}=\underline{w}\cdot\underline{a}=1_M
=\underline{w}$.
Now suppose $w = a_1\ \dotsm\ a_k$ with each $a_i\in A$ and let $t$ be an unstable 
position in $w$. Further, let $t'$ be the last unstable position preceding $t$. If $t$ is the first 
unstable position in $w$, we take $t'=0$. We claim that the letter $a_{t}$ cannot appear 
 anywhere strictly between positions $t'$ and $t$. To see this, notice that, since for every 
 $t' < j < t$, $j$ is stable, we have $\underline{w_{t'}}=\underline{w_{j}}=\underline{w_{j-1}a_j}
 =\underline{w_{t'}a_j}$, while $\underline{w_{t}}=\underline{w_{t-1}a_t}=\underline{w_{t'}a_t}
 \neq\underline{w_{t'}}$. 
These remarks imply that, if $\phi$ is an embedding of $r(w) = a_{t_1} \dots a_{t_\ell}$
as a subword of $w$, then for every $1 \leq i \leq \ell$, $\phi(i) \geq t_i$. 

Finally, let $\psi$ an embedding of $r(w)a$ as a subword of $w$, then $\psi$ must send 
the last position to some position $j$ with $j > t_\ell$. Then, since all positions after $t_{\ell}$ 
are stable, we get in particular that $\underline{w_{j-1}}=\underline{w_j} =\underline{w} $, so 
\[
\underline{wa} = \underline{w_{j-1}a} = \underline{w_j} = \underline{w}.\qedhere
\]
\end{proof}
Applying the exact same reasoning when starting to read $w$ from right to left, for the 
\emph{left} reduction $\ell(w)$ of $w$, we have, analogously to Lemma~\ref{lem:fallpoints}, 
that $a \ell(w) \subword w$ implies $\underline{aw}=\underline{w}$, and, analogously to 
Lemma~\ref{lem:finitereduced}, that the quasi-identity 
$\underline{yxu}=\underline{u}\To\underline{xu}=\underline{u}$ implies that there are finitely 
many left reduced words.

\begin{proposition}\label{prop:equat-implies-PT}
Let $M$ be a finite $\cJ$-trivial monoid, $\underline{(-)} \colon A^* \to M$ a homomorphism,
and $L$ a regular language recognized by $\underline{(-)}$. Then $L$ is piecewise testable.
\end{proposition}
\begin{proof}
	As explained above, by Theorem~\ref{thm:discrete-alt-height}, it suffices to prove that $L$ 
	has bounded alternation height in the poset $(A^*, \subword)$.	The images of the functions 
	$r \colon A^* \to A^*$ and $\ell \colon A^* \to A^*$ with respect to $\underline{(-)}$ are finite by
	Lemma~\ref{lem:finitereduced} and its symmetric `left' version. Therefore, the set
	\[ F := \{ uav \ \mid \ u \in \im(r), a \in A, v \in \im(\ell) \}\]
	is also finite. For any word $w \in A^*$, denote by $f(w)$ the set of subwords of $w$ that are in 
	$F$. Clearly, if $w \subword w'$ then $f(w)\subseteq f(w')$. 

	{\bf Claim.} Suppose $w \subword w'$. If $\underline{w}\neq\underline{w'}$, then $f(w)\subsetneq f(w')$. 

	{\it Proof of Claim.}
	Let $w \subword w'$. Then $f(w) \subseteq f(w')$. We prove the claim by contraposition.
	To this end, suppose that $f(w) = f(w')$. Since $w'$ can be obtained from $w$ by inserting a finite 
	number of letters, we can assume $w = uv$ and $w'=uav$ for some $u,v \in A^*$ and $a \in A$; 
	if we prove for this special case that $\underline{w} = \underline{w'}$, then the contrapositive of the
	claim follows by an easy induction on the length difference between $w'$ and $w$.

         Now, since $uav = w'$, we have $r(u)a\ell(v) \subword w'$, so $r(u) a \ell(v) \in f(w')$. Since 
         $f(w)= f(w')$, we also have that $r(u) a \ell(v) \in f(w)$, so $r(u) a \ell(v)$ is a subword of $w = uv$. 
         This implies that either $r(u)a\subword u$ or $a \ell(v) \subword v$, according to where the embedding
	sends the `middle' position $a$. Let us assume $r(u) a \subword u$, the other case is symmetric. 
	Lemma~\ref{lem:fallpoints} gives that $\underline{ua} = \underline{u}$, and thus 
	\[
	\underline{w} =\underline{u}\underline{v} = \underline{uav} = \underline{w'},
	\]
	and the claim is proved.

	Finally, let $x_1 \subword y_1 \cdots \subword x_n \subword y_n$ be alternating chain
	for $L$, and let $P$ be the image of $L$ by $\underline{(-)}$. Since $L$ is recognized by
	$\underline{(-)}$, we have $w\in L$ if, and only if, $\underline{w}\in P$. Since 
	$x_1 \subword y_1 \cdots \subword x_n \subword y_n$ alternates in and out of $L$, no 
	two consecutive elements of the chain can have the same image by $\underline{(-)}$.
	Therefore, by the claim, we have  
	\[
	f(x_1)\subsetneq f(y_1) \cdots \subsetneq f(x_n)\subsetneq f(y_n).
	\]
	It follows that $n\leq |F|+1$.
\end{proof}
Putting the above results together, we get that the profinite equations defining
$\approx_{\cJ}$ characterize the congruence $\approx_{\PT}$ dual to $\PT$.
\begin{theorem}\label{thm:Simon-dual}
	The Boolean residuation algebra $\PT$ of piecewise testable languages is
	dual to the smallest profinite monoid congruence $\approx$ on $\wh{A^*}$
	that contains $((xy)^\omega x, (xy)^\omega)$ and $(x(yx)^\omega, (yx)^\omega)$ 
	for every $x, y \in \wh{A^*}$, that is,
  \[ \PT = \sem{(xy)^\omega x \stackrel{*}{=} (xy)^\omega, \ x(yx)^\omega \stackrel{*}{=} (yx)^\omega} \ . \]
\end{theorem}

\begin{proof}
	Recall that $\approx_{\cJ}$ is by definition the smallest profinite
	monoid congruence on $\wh{A^*}$ that contains the equations, and
	$\approx_{\PT}$ is the congruence dual to $\PT$; we show they are equal. 
	By Corollary~\ref{cor:PTsoundness}, $\approx_{\cJ}$ is contained in
	$\approx_{\PT}$. On the other hand, if $L$ is a regular language
	satisfying the equations defining $\approx_{\cJ}$, then
	Proposition~\ref{prop:equat-implies-PT} together with
	Lemma~\ref{lem:PT-equations} imply that $L$ is piecewise testable. Thus,
	the Boolean algebra associated with $\approx_{\cJ}$ is contained in $\PT$,
	so that $\approx_{\PT}$ is contained in $\approx_{\cJ}$.
\end{proof}
It follows in particular that the membership problem for piecewise testable
languages is decidable. Indeed, by Theorem~\ref{thrm:reglanguage}, the
syntactic monoid of a regular language $L$ is computable from any automaton
recognizing it, and Theorem~\ref{thm:profinite-equation-synt-mon}, applied to
the profinite monoid equations defining $\approx_{\cJ}$, shows that these equations
may be tested on the syntactic monoid. We conclude the following.
\begin{corollary}\label{cor:Simon}
A regular language $L$ is piecewise testable if, and only if, its syntactic monoid
$M_L$ satisfies the equations $(xy)^\omega x = (xy)^\omega$ and $x(yx)^\omega =
(yx)^\omega$ for every $x, y \in M_L$. In particular, it is decidable for a
given regular language whether it is piecewise testable.
\end{corollary}

\ourexercises

\begin{ourexercise}\label{exe:subword-regular}
	Let $A$ be a finite alphabet and $w \in A^*$. Construct a finite
	automaton that recognizes the language ${\uparrow} w := \{v \in A^* \ |
	\ w \subword v\}$. Deduce that any finitely generated up-set (and thus,
	by Theorem~\ref{thm:higman}, any up-set) of $(A^*, \subword)$ is
	regular. 
\end{ourexercise}

\begin{ourexercise}\label{exe:Green}
  Let $M$ be a monoid. For any $m, n \in M$, we define 
  \[ m \preceq_{\cL} n \iffdef \text{ there exists } x \in M \text{ such that } xn = m \ , \] 
  \[ m \preceq_{\cR} n \iffdef \text{ there exists } y \in M \text{ such that } ny = m \ , \] 
  \[ m \preceq_{\cJ} n \iffdef \text{ there exist } x, y \in M \text{ such that } xny = m \ . \] 
\begin{enumerate}
  \item Prove that $\preceq_{\cL}, \preceq_{\cR}, \preceq_{\cJ}$ are preorders on $M$, and that ${\preceq_{\cL}} \cup {\preceq_{\cR}} \subseteq {\preceq_{\cJ}}$.
\end{enumerate}
The equivalence relations associated to these preorders are denoted
$\cL, \cR, \cJ$, respectively. 
In general, a monoid is called \emph{$\cL$-trivial},
\emph{$\cR$-trivial}, or \emph{$\cJ$-trivial} if the corresponding
equivalence relation is equal to the diagonal (also known as `trivial') relation.
We show in this exercise that this definition of $\cJ$-trivial 
coincides with the one given
in the main text for profinite monoids, and also that $\cJ$-trivial is equivalent to the conjunction of $\cL$-trivial and $\cR$-trivial for 
profinite monoids.
\index{Green's equivalence relations}

\begin{enumerate}[resume]
  \item Prove that if $\cJ$ is the trivial relation, then $\cL$ and $\cR$ are both the trivial relation.
\end{enumerate}
Now suppose that $M$ is a profinite monoid. 
\begin{enumerate}[resume]
  \item Prove that, for any $x, y \in M$, $(xy)^{\omega} = x (yx)^{\omega-1} y$.
  \item Prove that, for 
  any $x, y \in M$, 
  \[ (xy)^\omega x {\mathcal{R}} (xy)^\omega \text{ and }
  x(yx)^\omega {\mathcal{L}} (yx)^\omega . \]  
  Conclude that if $\mathcal{R}$ is the trivial relation then the equation $(xy)^\omega x = (xy)^\omega$ holds in $M$, and if $\mathcal{L}$ is the trivial relation, then the equation $x (yx)^\omega = (yx)^\omega$ holds in $M$.
  \item Conversely, prove that if both $(xy)^\omega x = (xy)^\omega$ and $x(yx)^\omega = (yx)^\omega$ for all $x,y \in M$, then $\mathcal{J}$ is the trivial relation.
  
  \hint{Show first that, if $mxy = m$ in $M$, then $m(xy)^\omega = m$. Then generalize this idea and use the given equations.}
  \item Conclude that, in any profinite monoid $M$, the following are equivalent:
  \begin{enumerate}[label=(\roman*)]
  \item the equivalence relation $\mathcal{J}$ is trivial;
  \item the equivalence relations $\mathcal{L}$ and $\mathcal{R}$ are both trivial;
  \item the equations $(xy)^\omega x = (xy)^\omega$ and $x(yx)^\omega = (yx)^\omega$ hold for all $x,y \in M$.
  \end{enumerate}
\end{enumerate}
\end{ourexercise}

\begin{ourexercise}\label{exe:ptplus-idl}
	Recall that $\Idl(A^*, \subword)$ is a spectral domain in the Scott
	topology, whose compact-opens are the sets of the form $K_F := \{I \in
	\Idl(A^*) \ \mid \ I \cap F \neq \emptyset\}$, 
	with $F$ a finite subset of $A^*$. Also recall that a function 
	$\pi \colon \wh{A^*} \to \Idl(A^*)$ was defined by $\pi(u) 
	:= \{w \in A^* \mid u \in \wh{{\uparrow}w}\}$.
\begin{enumerate}
	\item Prove that $\pi$ is well-defined, and that in particular for any
		$u \in A^*$, $\pi(u) = {\downarrow}u$, the principal down-set of
		$u$. 
	\item Prove that $\pi$ is continuous with respect to the patch topology on
		$\Idl(A^*)$. \hint{It suffices to prove that $\pi^{-1}(K_F)$ is
		clopen in $\wh{A^*}$, for every finite subset $F$ of $A^*$}.
	\item Prove that there is a well-defined continuous monoid
		multiplication on $\Idl(A^*)$ given, for $I, J
		\in \Idl(A^*)$, by 
		\[ I \cdot J := \{ ww' \ \mid \ w \in I, w' \in J\}.\]

	\item Prove that $\pi(uv) = \pi(u)\cdot\pi(v)$ for any $u, v \in A^*$.
			\item Conclude, using Exercise~\ref{exe:cts-mon-hom}, that $\pi$ is a
		continuous monoid homomorphism $\wh{A^*} \to \Idl(A^*)$.
\end{enumerate}
\end{ourexercise}

\begin{ourexercise}\label{exe:R-triv-reduced}
	Let $\underline{(-)} \colon A^* \to M$ be a monoid homomorphism and 
	suppose that there exist $u, x, y\in A^*$ such that $\underline{uxy} = 
	\underline{u}$ but $\underline{ux} \neq \underline{u}$. Prove that there 
	exist infinitely many right reduced words for the homomorphism $A^*\to M$.
\end{ourexercise}

\begin{ourexercise}
  Prove that the complement of the language $L$ of Exercise~\ref{exe:non-det-example} 
  is positively piecewise testable.
\end{ourexercise}

\begin{ourexercise}\label{exe:majority}
In this exercise we let $A = \{a,b\}$ and consider the language 
\[ L := \{ w \in A^* \ \mid \ |w|_a > |w|_b \}, \]
where, for $w \in A^*$ and $x \in A$, we write $|w|_x$ for the number of 
occurrences of the letter $x$ in the word $w$. The language $L$ is often 
called `majority'  in the literature, since it contains the set of words in 
which the majority of letters is $a$. For any $k \in \bZ$, define the language
\[ L_k := \{w \in A^* \ \mid \ |w|_a - |w|_b > k \}. \] 
\begin{enumerate}
  \item Show that the closure under quotienting $\cQ(L)$ of $L$ is equal to 
  $\{L_k \ \mid \ k \in \bZ\}$. Conclude in particular that $L$ is not regular.
  \item Show that $\bZ$ is, up to isomorphism, the discrete dual of the 
  complete Boolean algebra generated by $\cQ(L)$ and that the monoid 
  morphism $\{a,b\}^*\to\bZ, w\mapsto |w|_a - |w|_b$ is the discrete dual 
  of the inclusion of this algebra in the powerset of $\{a,b\}^*$. Verify that 
   the kernel of this monoid morphism is the syntactic congruence, 
   $\equiv_L$, and thus that the morphism is, up to isomorphism, the 
   syntactic morphism of $L$. 
  \item Deduce that the quotienting Boolean algebra $\cB(L)$ generated 
  by $L$ is equal to the Boolean algebra generated by the languages 
  $L_k$. Show that it is not a complete Boolean algebra.
  \item Prove that $\cB(L)$ is isomorphic to the Boolean subalgebra $M$ 
  of $\cP(\bZ)$, considered in Example~\ref{exa:equationssubalg}, 
  consisting of the subsets $S$ of $\bZ$ such that both $S \cap \bZ^+$ 
  is finite or co-finite and $S \cap \bZ^-$ is finite or co-finite. Here, as in 
  Example~\ref{exa:equationssubalg}, $\bZ^+$ denotes the set of 
  positive integers and $\bZ^-$ denotes the set of negative integers.
 \item Write $R_+$ for the graph of the addition operation on $\bZ$, 
 that is, $R_+ = \{(x,y,z) \in \bZ^3 \ \mid \ x + y = z\}$, and denote by 
 $\overline{R_+}$ the closure of $R_+$ in the space $(\bZ^{+\infty}_{-\infty})^3$. 
 Show that, for any $u \in \bZ^{+\infty}_{-\infty}$, the elements 
 $(-\infty,+\infty,u)$ and $(+\infty, -\infty, u)$ are in $\overline{R_+}$; that is, 
 $\overline{R_+}$ contains the set $\{(-\infty,+\infty),(+\infty,-\infty)\}
 \times \bZ^{+\infty}_{-\infty}$. Conclude that there is no continuous 
 binary function on $\bZ^{+\infty}_{-\infty}$ that extends the addition 
 on $\bZ$.
 \item Using Example~\ref{exa:equationssubalg}, show that $\cB(L)$ is given by 
 the following set of spatial equations on $\beta A^*$:
 $$\{\mu x\approx \mu, x\mu\approx \mu \mid \mu \in \beta A^*, x\in A\}$$
 Here, $\mu x$ and $x\mu$ denote the left and right action of $A^*$ on $\beta A^*$ 
 induced by the biaction dual to the quotienting biaction of $A^*$ on $\cP(A^*)$.
\end{enumerate}
\end{ourexercise}

\section{Open multiplication}\label{sec:openmult}

What better to end with than something open? This book is just a small sampler
of the things you can do with duality -- even within our two chosen
applications of domain theory and automata theory, the theory goes much 
beyond what we have touched on here. 
In the previous sections, we saw that the multiplication of a profinite monoid
can be seen as dual to residuation structure on the dual Boolean algebra, or as
dual to a biaction structure by quotienting operations. In this final section,
we consider a third possibility for seeing the multiplication of a profinite
monoid as the dual of a binary operation on the dual Boolean algebra, which
applies only in the special case where multiplication is an open map.

Theorem~\ref{thm:bintopalg-dual-to-jpp} establishes a connection between the
 \emph{residual} operations $(\backslash, \slash)$ on a Boolean algebra and 
 continuous binary operations $f$ on its dual Boolean space. One may wonder 
 what it takes for the \emph{forward} image map given by the binary operation to
be the dual of an operation. Using Definition~\ref{dfn:compatiblerelation}, we
obtain the following requirements on $f$: 
\begin{enumerate}
  \item For all $x\in X$ the inverse image $f^{-1}(x)$ is closed;
  \item The forward image of a pair of clopens is clopen.
\end{enumerate}
Without assuming that $f$ is also continuous, these are not  very natural
conditions on a map between topological spaces. However, we do obtain 
the following useful corollary.

\begin{corollary}\label{cor:openmult}
  Let $(X,f)$ be a binary topological algebra based on a Boolean space, and
  let $B$ be the dual Boolean algebra of clopen sets. The operation $f$ is an open
  mapping if, and only if, $B$ is closed under the forward operation
  \[
	(U, V) \mapsto U \bullet V := \{f(u,v) \ \mid \ u \in U, v \in V \},
  \]
  that is, if, and only if, the complex multiplication of two clopen sets of $X$
  is again clopen.  In this case, the graph of $f$ is the relational dual to the
  operation $\bullet$ on $B$.  
\end{corollary}
\begin{proof}
  Since all continuous maps from compact spaces to Hausdorff spaces are closed
  mappings, it follows that $B$ is closed under the complex operation $\bullet$ if, 
  and only if, $f$ is an open map. The conditions required for the graph $R$ of $f$ 
  to be the dual of this operation $\bullet$ are that $f^{-1}(x)=R[\_,\_,x]$ is closed 
  for each $x\in X$ and that $U \bullet V=R[U,V,\_]$ is clopen whenever both $U$ 
  and $V$ are clopen. For $f$ continuous and $X$ Hausdorff, the first condition 
  always holds. The second condition holds if $f$ is an open mapping. 
\end{proof}

Note that, in the special case of languages and the free profinite monoid, a 
subtlety arises from the mixing of discrete and topological duality that happens 
there. More precisely, if $B$ is a Boolean algebra of regular languages, then 
it is possible that $B$ fails to be closed under the multiplication in $\cP(A^*)$, 
even though its natural embedding in the powerset of the dual space $X$ of 
$B$ is, a Boolean subalgebra closed under complex multiplication. That is, 
the Boolean subalgebra $\{\wh{L} \ \mid \ L \in B\}$ of $\cP(X)$, is closed 
under complex multiplication, see Example~\ref{exa:grouplang-not-closed} 
below.

\begin{example}\label{exa:RegAclosed}
  One of the, if not the, most famous theorems about regular languages 
  is Kleene's Theorem \parencite{Kleene1956} which asserts that $\RegA$ 
  is closed under concatenation product in $\cP(A^*)$. This fact also has a 
  very natural duality proof in categorical logic \parencite{JM21}.
\end{example}

\begin{example}\label{exa:grouplang-not-closed} A free profinite group, or in
	fact any topological group, has open multiplication. Indeed, for any
	element $g$ of a topological group $G$, the multiplication $x \mapsto x
	\cdot g$ is a homeomorphism of $G$ with inverse $x \mapsto x g^{-1}$,
	and it is thus in particular an open map. Therefore, for any $U, V
	\subseteq G$ with $U$ open, $U \cdot V = \bigcup_{g \in V} Ug$ is open.
However, the product of group languages is not
	necessarily a group language. Indeed, let $L=(a^2)^*a$, then $L$
	is recognized by the finite group $\mathbb Z_2$, but it is not difficult
	to see that $L^2=(a^2)^+$ can not be recognized by a finite group (see
	Exercise~\ref{ex:group-not-concat-closed}).  We refer the reader to
	\cite{RibZal2010} for an extensive survey of the theory of profinite
	groups. 
\end{example}

We do have the following sufficient condition for an algebra of regular 
languages to be closed under concatenation product, in the situation of 
Corollary~\ref{cor:openmult} above.
\begin{proposition}
	Let $f \colon \wh{A^*} \onto X$ be a profinite monoid quotient, 
	let $B$ be the Boolean subalgebra of $\RegA$ dual to $X$, 
	$\widehat{(\ )}$ the Stone embedding of $B$ in $\cP(X)$, 
	and assume that the multiplication on $X$ is an open mapping. Then
	for $L,K\in B$, the element $M\in B$ corresponding to the clopen 
	$\widehat{L}\widehat{K}$ is the least $N\in B$ for which $LK\subseteq N$
	holds in $\cP(A^*)$.
	If in addition $\{u\}$ is in $B$ for every $u\in A^*$, then $B$ is 
	closed under concatenation product in $\cP(A^*)$.
\end{proposition}
\begin{proof}
  Since the multiplication is open in $X$, for every $L,K\in B$, there is an $M\in B$ 
  with $\widehat{L}\widehat{K}=\widehat{M}$. We first show that, even if $LK\not\in B$, 
  this entails that $M$ is the least element of $B$ above $LK$ in $\cP(A^*)$. Let  $N\in B$. 
  Then we have
  \begin{align*}
  LK\subseteq N &\iff L\subseteq N/K\\
                          &\iff \widehat{L}\subseteq \widehat{N/K}=\widehat{N}/\widehat{K}\\
                          &\iff \widehat{M} = \widehat{L}\widehat{K}\subseteq \widehat{N}\\
                          &\iff M\subseteq N.
  \end{align*}
 Finally, suppose that $\{u\}\in B$ for every $u\in A^*$, then also $\{u\}^c\in B$,
 and thus $u\not\in LK$ implies $LK \subseteq \{u\}^c$, which implies $M\subseteq\{u\}^c$,
 and thus $u \not\in M$. By contraposition, we have $M\subseteq LK$  and thus $LK=M$.
\end{proof}
Openness of the multiplication has been studied in language theory, but it also 
appears in categorical first-order logic. Classical
research in categorical first-order logic is concerned with extending Stone
duality to a setting where the Boolean algebras are equipped with quantifiers, a
topic that is also still under development today, see for example \cite{Lur2019,
GM22}. A fundamental structure in that field is that of a \emphind{Boolean
hyperdoctrine}, which is a functor $P$ valued in Boolean algebras which is
required to satisfy two algebraic conditions on the morphisms, called the
Beck-Chevalley and Frobenius conditions, which have a natural meaning through
duality theory.

In one recent result in this direction~\parencite{JM21}, openness is combined with
\emphind{equidivisibility}, a notion which has been considered in language
theory \parencite{AC09}.  A semigroup $S$ is called \textit{equidivisible} if, for
all elements $u_1, u_2, v_1, v_2 \in S$ such that $u_1u_2 = v_1v_2$, there
exists an element $k \in S$ such that either $u_1k = v_1$ and $u_2 = kv_2$, or
$v_1k = u_1$ and $v_2 = ku_2$. 
Relative to duality theory, this condition on the
binary operation $\cdot$ to be equidivisible  is precisely the `back' condition
of bounded morphism (Definition~\ref{dfn:boundedmorphism}).   In
\cite[Theorem~4.3]{JM21}, it is shown that a profinite semigroup $S$ is dual to a
first-order theory enriching the theory of bounded linear orders if, and only
if, the multiplication on $S$ is open and equidivisible. To obtain this result,
it is shown that, under Stone duality,  the openness of the multiplication is
dual to the Beck-Chevalley condition, and the equidivisibility corresponds to
the Frobenius condition. This theorem thus establishes yet another connection
between language theory and topological algebra, mediated by extended Stone
duality.

\ourexercises
\begin{ourexercise}\label{ex:group-not-concat-closed}
	Let $A = \{a\}$ and $L \subseteq A^*$ the language of words of odd
	length.
	\begin{enumerate}
	\item Show that the syntactic monoid of $L$ is the two-element group
		$\mathbb{Z}_2$, so that $L$ is a group language.
	\item Show that the syntactic monoid of  $L \cdot L$, the language 
	of words of non-zero even length, is the three-element monoid 
	$(\mathbb{Z}_2)^I$, that is, the monoid obtained from
		$\mathbb{Z}_2$ by adjoining a new neutral element.
	\item Deduce that $L \cdot L$ can not be recognized by a finite group.
	\end{enumerate}
\end{ourexercise}

\notessec
{The fact that the Stone dual of the Boolean algebra of regular
languages over an alphabet is dual to the topological space underlying the
free profinite monoid is well-known but played a small role
historically in the area. An early exception is the article
\cite{Pippenger97}, see also \cite{Almeida1995}, and more
recently \cite{RhodesSteinberg2008}, where the connection between profinite
semigroups and Boolean rings with coalgebraic structure is exploited, and
\cite{GGP2008,GGP2010}, where the connection with modern topological
duality as it is applied in logic was first developed. The point of these
two last approaches is that it is not only the underlying space that can be
studied via duality but the entire profinite semigroup or monoid structure.
The duality theoretic underpinnings of the work in \cite{GGP2008,GGP2010} were
worked out in \cite{Geh16}. For other key publications see
\cite{Pin2009,BrPi2009,GKP2016,Pin2017}. }

{Starting with the definitions on page~\pageref{pg:residual-defs} and 
throughout this chapter, we fundamentally view the dual of an additional 
$n$-ary operation on a lattice or a Boolean algebra as an $n+1$-ary 
relation on the dual space. Such a relation may also be seen as so-called 
\emph{co-algebraic structure}. Co-algebraic structure is, by definition, 
categorically dual to algebraic structure. For this reason, if we start, 
instead, from algebraic structure on the spatial side, then we get 
co-algebraic structure on the dual lattices or Boolean algebras. This 
formal duality point of view, for monoids, is exploited in 
\cite[Section 8.4]{RhodesSteinberg2008}. What we show in 
Section~\ref{sec:syntmon} is that, if we restrict ourselves to surjective 
morphisms, then, in the setting of discrete duality, monoid structure 
on a set is dual to residuation structure on the dual Boolean algebras. 
This provides a kind of duality between algebras and algebras. This 
can be lifted to non-surjective morphisms and to general Stone and 
Priestley duality as well, see \cite{Geh16}.
} %

{The definition of syntactic congruence (Definition~\ref{def:syntactic}) 
is classical in the theory of regular languages and can be traced back 
to the works \cite{Schutzenberger56,Myhill1957,Nerode1958}. In 
Section~\ref{sec:EilReittheory}, we encountered the notion of syntactic 
\emph{preorder}, which generalizes syntactic congruences to an ordered 
setting, see Remark~\ref{rem:synt-preorder}. This notion is already 
present in the work of \cite{Schutzenberger56} and was pioneered by 
\cite{Pin1995variety}, long before the connection with Priestley duality 
was realized.}

{We are grateful to Jean-Éric Pin and Howard Straubing for many 
enlightening discussions on the proof of Simon's theorem discussed in 
Section~\ref{sec:EilReittheory}. This proof 
will be the topic of a forthcoming joint paper of this book's authors with
Jérémie Marquès, who made significant contributions to the proof in this
section, in particular Lemmas~\ref{lem:finitereduced}~and~\ref{lem:fallpoints}.  
Similar ideas, although not formulated using profinite methods, are present in 
Stern's proof of Simon's theorem~\parencite{Stern85}.
Many proofs of Simon's theorem have been published since the original, 
see in particular also \cite{Almeida91,Higgins97,HenckellPin98,Klima11}.}
\nocite{*}
\cleardoublepage
\phantomsection
\markboth{Notation}{Notation}
\addcontentsline{toc}{chapter}{Notation}
{
\printindex[notation]}

\cleardoublepage

\phantomsection
\markboth{Index}{Index}
\addcontentsline{toc}{chapter}{Index}
{\small
\printindex}

\cleardoublepage
\phantomsection
\markboth{Bibliography}{Bibliography}
\addcontentsline{toc}{chapter}{Bibliography}

\printbibliography
\end{document}